\documentclass[11pt]{book}   
\usepackage{latexsym,amssymb,amsmath,amsfonts,epsf}
\usepackage{tabularx}
\usepackage{theorem}
\usepackage{rotating}
\usepackage{makeidx} 
\makeindex 
\renewcommand{\baselinestretch}{1.5}
\newtheorem{teo}{Theorem}[section]
\newtheorem{definition}[teo]{Definition}
\newtheorem{lemma}[teo]{Lemma}
\newtheorem{corollary}[teo]{Corollary}
\newtheorem{proposition}[teo]{Proposition}
\newtheorem{theorem}[teo]{Theorem}

{\theoremstyle{plain} \theorembodyfont{\rmfamily} 
{\theoremstyle{plain} \theorembodyfont{\rmfamily} 
\newenvironment{proof}{\begin{trivlist}
    \item[\hskip\labelsep{\bf Proof.}]}{$\hfill\Box$\end{trivlist}}

\newtheorem{_rem}[teo]{Remark}
\newtheorem{_eje}[teo]{Example}
\newtheorem{_ass}[teo]{Assumptions}

\newenvironment{rem}{\def\@begintheorem##1##2{\trivlist%
 \item[\hskip\labelsep{\sffamily\bfseries ##2\ ##1}]}\begin{_rem}}{$\hfill\spadesuit$ \end{_rem}}
\newenvironment{eje}{\def\@begintheorem##1##2{\trivlist%
 \item[\hskip\labelsep{\sffamily\bfseries ##2\ ##1}]}\begin{_eje}}{$\hfill\spadesuit$\end{_eje}}
\newenvironment{ass}{\def\@begintheorem##1##2{\trivlist%
 \item[\hskip\labelsep{\sffamily\bfseries ##2\ ##1}]}\begin{_ass}}{$\hfill\spadesuit$\end{_ass}}
\makeatother

\DeclareMathOperator{\op}{Op}

\DeclareMathOperator{\Sp}{Sp}

\newcommand{\one}{\mathbf{1}}
\newcommand{\bp}{\mathbf{p}}
\newcommand{\bq}{\mathbf{q}}

\newcommand{\ip}[2]{ \langle #1 , #2 \rangle }
\newcommand{\iip}[2]{ ( #1 , #2 ) }
\newcommand{\smp}[3]{ \sigma_{#3}( #1 , #2 ) }
\newcommand{\norm}[1]{ \Vert #1 \Vert }
\newcommand{\alg}[1]{ \text{alg}(\{  #1  \})}
\newcommand{\salg}[1]{ \text{$*$-alg}\{ #1 \} }

\newcommand{\osalg}[1]{ C^*( \{ #1 \}) }
\newcommand{\np}{\norm{\mathbf{p}}}
\newcommand{\nps}{{\np}^2}

\newcommand{\sbr}[2]{ [ #1, #2 ]_{sb} }
\newcommand{\rga}{\cA_g^r}
\newcommand{\rgs}{{\cH}^{gr}}
\newcommand{\sn}[1]{\text{span}\{ #1 \}}
\newcommand{\sqb}[1]{[\{ #1 \}]}

\newcommand{\fM}{\mathfrak{M}}

\newcommand{\fD}{\mathfrak{D}}

\newcommand{\fX}{\mathfrak{X}}
\newcommand{\fB}{\mathfrak{B}}
\newcommand{\vL}{\varLambda}

\newcommand{\fF}{\mathfrak{F}}
\newcommand{\fH}{\mathfrak{H}}
\newcommand{\fS}{\mathfrak{S}}
\newcommand{\fT}{\mathfrak{T}}

\newcommand{\cA}{\mathcal{A}}
\newcommand{\cG}{\mathcal{G}}
\newcommand{\cS}{\mathcal{S}}
\newcommand{\cSet}{\widehat{\mathcal{Q}}}
\newcommand{\cSetin}{\mathcal{Q}}
\newcommand{\cSe}{\mathcal{S}^{ext}}
\newcommand{\cH}{\mathcal{H}}

\newcommand{\cL}{\mathcal{L}}
\newcommand{\cP}{\mathcal{P}}
\newcommand{\cB}{\mathcal{B}}

\newcommand{\cX}{X}
\newcommand{\cXs}{\cX_s}
\newcommand{\cXsp}{\cX_{s'}}
\newcommand{\cXspp}{\cX_{s''}}
\newcommand{\Qs}{Q_s}
\newcommand{\Qsp}{Q_{s'}}
\newcommand{\Qspp}{Q_{s''}}
\newcommand{\cR}{\mathcal{R}}
\newcommand{\cD}{\mathcal{D}}
\newcommand{\cE}{\mathcal{E}}
\newcommand{\cC}{\mathcal{C}}
\newcommand{\cF}{\mathcal{F}}
\newcommand{\cO}{\mathcal{O}}
\newcommand{\cN}{\mathcal{N}}
\newcommand{\cT}{\mathcal{T}}
\newcommand{\cU}{\mathcal{U}}
\newcommand{\cM}{\mathcal{M}}

\newcommand{\vp}{\varphi}
\newcommand{\la}{\lambda}
\newcommand{\bla}{ \boldsymbol{\lambda}}
\newcommand{\bmu}{\boldsymbol{\mu}}

\newcommand{\bZ}{\mathbb{Z}}

\newcommand{\drhb}{\tilde{\drb}}
\newcommand{\drb}{\delta}
\newcommand{\rd}{\mathrm{d}}

\newcommand{\tu}{{\gh_{s}}}
\newcommand{\tp}{\tilde{ph}}
\newcommand{\fDp}{\tilde{\fD}}
\newcommand{\fXp}{\tilde{\fX}}

\newcommand{\mH}{Y}
\newcommand{\bch}{Q}

\newcommand{\fL}{\fH^{\mH}}
\newcommand{\cQ}{\overline{Q}}

\newcommand{\LL}{1}
\newcommand{\JJ}{2}
\newcommand{\fDL}{\fD_{\LL}}
\newcommand{\fDJ}{\fD_{\JJ}}
\newcommand{\fXL}{\fX_{\LL}}
\newcommand{\fXJ}{\fX_{\JJ}}
\newcommand{\fHL}{\fH_{\LL}}
\newcommand{\fHJ}{\fH_{\JJ}}
\newcommand{\fLL}{\fL_{\LL}}
\newcommand{\fLJ}{\fL_{\JJ}}
\newcommand{\cHL}{\cH_{\LL}}
\newcommand{\cHJ}{\cH_{\JJ}}

\newcommand{\fst}{\frac{1}{\sqrt{2}}}
\newcommand{\FC}{(2(2\pi)^3))^{-\frac{1}{2}}}
\newcommand{\FCC}{(2(2\pi)^3)^{-1}}

\newcommand{\Reg}{\fS_{r}}
\newcommand{\SReg}{\fS_{sr}}
\newcommand{\sfu}{\mathcal{SF}}
\newcommand{\cAe}{\cA_{ext}}
\newcommand{\ncss}{auxiliary symplectic space }
\newcommand{\ncsf}{auxiliary symplectic form }
\newcommand{\HT}{HT \,}
\newcommand{\KOB}{KO Abelian BRST\,}
\newcommand{\KO}{KO \,}
\newcommand{\MM}{m}
\newcommand{\WW}{W}

\def\im{{\rm Im}}
\def\re{{\rm Re}}
\newcommand{\gh}{\eta}
\newcommand{\cgh}{\rho}

\newcommand{\pfit}{ \smallskip \noindent}

\def\XP#1!{\renewcommand{\baselinestretch}{.7}\marginpar{
{\footnotesize #1}\hfil}\renewcommand{\baselinestretch}{1.5}}

\def\XB{\marginpar{
{\footnotesize\bf Change~starts-----}\lower 11pt\hbox{\mathsurround=0pt$
\!\!\displaystyle{
\Bigg\downarrow}$\mathsurround=3pt}}}

\def\XE{\marginpar{{\footnotesize\bf Change~ends-----}\raise 
10pt\hbox{\mathsurround=0pt$
\!\!\displaystyle{
\Bigg\downarrow}$\mathsurround=3pt}}}

\def\@begintheorem#1#2{\trivlist%
 \item[\hskip \labelsep{\sffamily\bfseries #2\ #1}]\itshape}

\newenvironment{beweis}{{\em Proof:}}{\hfill $\rule{2mm}{2mm}$
\vspace{3mm}}

\DeclareMathAlphabet{\Ma}{U}{msa}{m}{n}
\DeclareMathAlphabet{\Mb}{U}{msb}{m}{n}
\DeclareMathAlphabet{\Meuf}{U}{euf}{m}{n}

\def\got#1{\Meuf{#1}}

\DeclareSymbolFont{ASMa}{U}{msa}{m}{n}
\DeclareSymbolFont{ASMb}{U}{msb}{m}{n}
\DeclareMathSymbol{\hrist}{\mathord}{ASMa}{"16}
\DeclareMathSymbol{\varkappa}{\mathalpha}{ASMb}{"7B}
\DeclareMathSymbol{\CrPr}{\mathord}{ASMb}{"6F}

  \def\al #1.{{\mathcal{#1}}}
  \def\ot #1.{{\got{#1}}}
 
\def\CCRX{\overline{\Delta( \fX,\,\sigma)}}
  \def\ccr #1,#2.{\overline{\Delta(#1,\,#2)}}

  \def\b #1.{{\bf #1}}
  \def\cross#1.{\mathrel{\mathop{\times}\limits_{#1}}}
                  
  \def\C{\Mb{C}}
  \def\N{\Mb{N}}
  
  \def\R{\Mb{R}}
  \def\Z{\Mb{Z}}
 \def\un{{\one}}

\def\f #1,#2.{\mathsurround=0pt \hbox{${#1\over #2}$}\mathsurround=5pt}

  \def\wt{\widetilde}

  \def\cross #1.{\mathrel{\raise 3pt\hbox{$\mathop\times\limits_{#1}$}}}
  \def\ol #1.{\overline{#1}}
\def\b #1.{{\bf #1}}
\def\slim{\mathop{\hbox{\rm s-lim}}}
\def\ker{{\rm Ker}\,}
\def\aut{{\rm Aut}\,}
\def\Aut{{\rm Aut}\,}
\def\dom{{\rm Dom}\,}
\def\sp{{\rm sp}}
\def\ran{{\rm Ran}\,}
\def\rep{{\rm Rep}\,}

\def\rlf{{R(\lambda,f)}}
\def\rsl{\mathord{\al R.(\fX,\sigma)}}
\def\s #1.{_{\smash{\lower2pt\hbox{\mathsurround=0pt $\scriptstyle #1$}}\mathsurround=5pt}}
\def\set #1,#2.{\left\{\,#1\;\bigm|\;#2\,\right\}}
\def\maprightu #1;{\smash{\mathop{\longrightarrow}\limits^{#1}}}
\def\maprightd #1;{\smash{\mathop{\longrightarrow}\limits_{#1}}}
\def\maprightt #1,#2.{\mathrel{\smash{\mathop{\longrightarrow}\limits_{#1}^{#2}}}}
\def\chop{\hfill\break}
\def\j{\phi}

\def\ie{\textit{i.e.\ }}
\def\eg{\textit{e.g.\ }}

\def\margin #1.{\marginpar{#1}}

\def\WD{{\cal O}}
\def\Ob{{\cal P}}

\title{}\author{}\date{}

\begin{document}

\frontmatter

\chapter*{\center \Huge{The Mathematical Structure of the Quantum BRST Constraint Method}} 
\thispagestyle{empty}
\vspace{-1cm}
\begin{center} \Large{.....} \\ \vspace{2cm}
\normalsize{A thesis presented to} \\ \vspace{0.5cm} \large{The University of 
New South Wales} \\ \vspace{0.5cm} \normalsize{in fulfillment of the thesis 
requirement} \\ \normalsize{for the degree  of} \\ \vspace{0.5cm} 
\large{Doctor of Philosophy} \\ \vspace{0.5cm} \normalsize{by}

\vspace{2cm}

{\scshape \LARGE{Patrick Costello}}

\vspace{2cm}

 30/03/08

\end{center}


\newpage
\thispagestyle{empty}
\newpage

\chapter{Abstract}

This thesis describes mathematical structures of the quantum BRST constraint method. Ultimately, the quantum BRST structures are formulated in a $C^*$-algebraic context, leading to comparison of the quantum BRST and the Dirac constraint method in a mathematically consistent framework.

Rigorous models are constructed for the heuristic examples of BRST for quantum electromagnetism (BRST-QEM) and Hamiltonian BRST with a finite number of constraints. This facilitates comparison between the results produced by the BRST method, and the results of the $T$-procedure of Grundling and Hurst for the quantum Dirac constraint method. 

The different constraint methods are shown not to be equivalent for the examples of Hamiltonian BRST with a finite number of constraints that close, and a BRST-QEM model constructed using the Resolvent Algebra of Buchholz and Grundling with covariant test function space. Moreover, this leads to the following three consequences: 

The quantum BRST method, and quantum Dirac method of constraints, are not equivalent in general.  

Examples of quantum Hamiltonian BRST  can be constructed to show that the BRST method does not remove the ghosts in the BRST physical algebra.  This occurs since quantum Hamiltonian BRST selects multiple copies of the physical state space selected by the Dirac algorithm, and the ghosts are not removed from the BRST-physical state space. Extra selection criteria are required to select the correct physical space, which do not gaurantee correspondence between the Dirac and BRST physical algebras.

Conversely, the BRST physical algebra and Dirac physical algebra coincide when QEM is encoded in the auxiliary Resolvent Algebra. This is a rigorous example of Lagrangian BRST, hence quantum Lagrangian and quantum Hamiltonian BRST are not equivalent constraint methods.
\textbf{}
\newpage

\thispagestyle{empty}

\newpage

\chapter{Acknowledgments}
I would like to say my thanks and express sincere appreciation to my supervisor Dr. Hendrik Grundling for all his time, effort and patience. I would like to express my gratitude to Professor Klaus Fredenhagen for his insights into the correct statement of the BRST charge. Thanks also to Dr. Ben Warhurst for many stimulating discussions along the way. Finally, I am very grateful to the School of Mathematics and Statistics at the University of New South Wales and the University for their support.

\newpage
\thispagestyle{empty}
\newpage

\tableofcontents

\newpage

\thispagestyle{empty}

\newpage

\mainmatter

\chapter{Introduction}
This thesis describes the structures of the quantum BRST constraint method in a mathematically consistent manner, ultimately formulating it in a $C^*$-algebraic context. The BRST quantum constraint method is widely used in modern theoretical physics, finding application in many areas such as renormalization of gauge theories \cite{BRS1976,KuOj79,Wein2005II}, constraint theory in classical and quantum Hamiltonian constraint systems \cite{HenTei92,Bat1987,vHol2005,Mc1991} and constraint theory in string theory \cite{KaOg1982,GrSchWit1987,Gr2006}.  As a general quantum constraint method it suffers from the problems that it comes in a multitude of varieties, is commonly defined in the context of a specific model rather than in an independent algorithmic manner and often requires the addition of model-dependent constraints. Although there has been analysis of quantum BRST in varying degrees of mathematical rigour, as discussed below, there is still a lack of unification of the various results as well as a lack of comparison to alternative quantum constraint methods in mathematically precise manner. The following develops a mathematical framework to address these issues. First this framework is developed in the setting of operators acting on a Krein space, which is used to analyse common BRST examples and to discuss problematic issues related to BRST in the literature. The framework is then extended to an abstract $C^{*}$-algebraic context where is used to compare quantum BRST to the Dirac quantum constraint method in a mathematically consistent setting, and to further investigate the formulation of an abstract quantum BRST algorithm.

We start with a brief historical survey of BRST. As the BRST related literature is vast, the account given here covers the part which is central to the development of most branches of BRST. BRST theory officially began with the work by Becchi, Rouet and Stora (BRS) \cite{BRS1976,BRS1975}, and separately by Tyutin (T) \cite{Ty1975}, in relation to renormalization of gauge theories using path integral methods. In \cite{BRS1976} these authors discovered that it was possible to define a superderivation $\drb$ (or $\mathbf{s}$ in their notation) on the fields involved (gauge, ghost, etc) such that $\drb^2=0$; that the Lagrangians were $\drb$ invariant; and that $\drb$ produced the Slavnov-Taylor identities which correspond to gauge invariance in perturbation theory. In this way these authors were able to prove the renormalizability and unitarity of the $S$-matrix in gauge theory models that satisfy the Slavnov-Taylor identities, which they demonstrate explicitly for the example of the $SU(2)$-Higgs-Kibble model.  

Following this, Kugo and Ojima formulated BRST as a constraint theory, published in a series of papers \cite{KuOj1978I,KuOj1978II,KuOJ1979II,KuHa1979,KuOj79} including the major work \cite{KuOj79} (referred to here as K\&O \cite{KuOj79}) that exhibits many of the general structures associated to BRST today. They work in an operator formalism, where fields are operator valued distributions acting on some inner product space $\cH$, and the construction $\drb$ for non-abelian gauge theories is interpreted as `replacing the gauge parameter by ghost parameter' in the gauge transformations. An  explicit formula for the BRST charge $Q$ is given such that it generates $\drb$, is such that $Q^2=0$, and hermicity assignments for the ghost fields are assumed which make $Q$ hermitian. The latter two conditions are incompatible with a Hilbert inner product and thus forces the inner product on $\cH$ to  be indefinite. As $Q^2=0$ we have $\ran Q \subset \ker Q$ and K\&O then go on to assume formally the the physical subspace is $\ker Q /\ran Q$ and to investigate consequent structures. They find that the physical subspace has a positive definite inner product using an argument based on a specific structure associated to the gauge theory set-up called the `quartet mechanism' applied to non-interacting theories and invoking asymptotic completeness (K\&O \cite{KuOj79} p46).  Under the assumption that the observables should be be those operators that factor to $\ker Q /\ran Q$, they prove that the \emph{local observables}  (as defined on K\&O \cite{KuOj79} p46) are the same as $\ker \drb$ factored to $ \ker Q /\ran Q$, which is used to give a proof colour confinement, K\&O \cite{KuOj79} Theorem 5.11 p69. 

The results in K\&O \cite{KuOj79}, while undoubtedly of great importance, rely on formal arguments as well as special features associated to gauge theories. It is natural to consider how these results can be extended to other theories with constraints and once done, to see if the inner product on $\ker Q /\ran Q$ is positive definite, if the operators which factor to $\ker Q /\ran Q$ always are in $\ker \drb$, how to make the structures rigorous, what the relation is to other constraint methods such as Dirac. As the K\&O approach is based upon modifying a symmetry of the Lagrangian, we refer to similar approaches as the \emph{Lagrangian approach to BRST}. Following K\&O \cite{KuOj79} the Lagrangian BRST structures are further analysed and developed by many authors such as \cite{KuUe1980,AbNak2002,AbNakOj1997,AbNak1988,AbNak1996,Nak1979,Nish1996,Nish1984,KaOg1982,Schf2001,NakOj1990,St1989}. Positivity of the physical space for gauge and string models is studied in \cite{Nak1979,St1989,AbNakOj1997,AbNak1996} where it is noted that positivity  is not guaranteed to hold in general \cite{St1989,Nak1979}, and that perturbative arguments may not be applicable to gauge theories in \cite{St1989,AbNak1996}. Kugo  and Uehara \cite{KuUe1980} investigate the observable algebra for Yang-Mills type gauge theories and find that a further restriction to ghost number zero BRST observables is needed to guarantee correspondence to the expected algebra (\cite{KuUe1980} Theorem 1 p1398). Also noteworthy is \cite{KaOg1982}, where the authors apply the Lagrangian BRST approach to bosonic string theory and find that 2-nilpotence of the BRST charge requires the critical dimension condition of $D=26$ and that the physical subspace has positive definite inner product via a `quartet mechanism'-like argument.

Concurrently with the development of the Lagrangian BRST approach, Fradkin et al. \cite{FraFradk1977,FraVil1975,BatFra1983,Bat1987}, published a series of papers with a general method for proving unitarity of the $S$-matrix for arbitrary gauge fixing conditions for degenerate Hamiltonian systems. This method extends the original system by ghost variables and uses the bracket structure of the constraints to construct an  2-nilpotent operator $Q$ ($\Omega$ in their notation) to generate a superderivation $\drb$. It is argued that we can extend the Hamiltonian $H_0$ of the original unextended system to $H_{\Psi}=H_0+\drb(\Psi)$, where $\Psi$ is an arbitrary function of the extended variables, and that the $S$-matrix derived from the extended system is independent of $\Psi$. In this way we achieve independence of the $S$-matrix from the `gauge fixing condition', in the form of the function $\Psi$. Although originally based on the path integral formalism, this method was set in an operator formalism in \cite{BatFra1983} where the basic objects $Q$ and $\drb$ had strikingly similar properties to the K\&O construction, but where the the method now applies to a general Hamiltonian system with constraints, and that the construction of $Q$ depended only on the commutation relations of the constraints. We refer to the above approach as the \emph{BFV (Batalin-Fradkin-Vilkovisky) approach to BRST}.

The BFV and Lagrangian BRST approach to constraints were brought together in the comprehensive paper by Henneaux \cite{Hen1985}, in which the BFV approach is reviewed and the major step of formulating it for classical Hamiltonian constraint theories is taken. Quantum theory is then considered in both operator and path integral formalisms where it is shown how the BRST charge for Yang-Mills as in K\&O \cite{KuOj79} can be arrived at using the BFV approach, that the extended system acts on a space of Berezin superfunctions, and a heuristic argument for equivalence between the Quantum Dirac constraint method and Quantum BRST method is given. The BFV formalism as described in this paper, and its subsequent development we will refer to as the \emph{Hamiltonian approach to BRST}.

Since \cite{Hen1985}, Henneaux and many others have significantly extended the analysis of Hamiltonian BRST, both in the classical and quantum case. An exhaustive list of these works is not given here, but can be found in H\&T \cite{HenTei92} which is the standard reference for Hamiltonian BRST theory in the physics literature. The first half of this book summarises the major structures now associated to classical BRST. The ghosts and conjugate ghosts that appear in classical BRST are shown to arise naturally from the geometric structures of Koszul-Tate and longitudinal differentials associated to the constraint surface and gauge transformations. Using homological perturbation theory, they show that the BRST superderivation $\drb$ arises naturally as the sum of these differentials plus higher order terms and satisfies $\drb^2=0$. This produces the BRST cohomology, and it is shown that the algebra of observables of the original constraint theory is isomorphic to the BRST cohomology at ghost number zero. The generator $Q$ of $\drb$ is canonically constructed as in H\&T \cite{HenTei92} and the analysis is extended to cover the reducible case, \ie where the constraint functions are linearly dependent. We will not treat classical BRST further in this thesis as it has already been put on a firm mathematical foundation, e.g. \cite{FiHenStTei1989,Sch1994,Sta1998,BarBraHen00}; we refer the reader to these papers and their bibliographies for details. 

With classical BRST being well defined, most authors construct Hamiltonian Quantum BRST structures by analogy to classical BRST with particular emphasis on the cohomological aspects. For example, the BRST observables are usually taken to be those $\drb$-invariant operators of ghost number zero rather than operators that factor to $\ker Q/ \ran Q$. While Quantum Hamiltonian BRST is constructed in several heuristic formalisms, such as the  path integral, we are only concerned here with the operator formalisms as this is most easily made rigorous. We study the operator formalism for Hamiltonian BRST given in H\&T \cite{HenTei92} Chapter 13 as this has the basic structures of most versions of Hamiltonian BRST. As in the Lagrangian case, there are many questions to ask of the structure of Hamiltonian Quantum BRST. Fundamental to a probabilistic interpretation is the positivity of the inner product on the physical subspace. 

Positivity of the inner product on the physical subspace is sometimes assumed to be true, or that it will be only for ghost number zero states in the physical subspace, but it is noted in several places, e.g. H\&T \cite{HenTei92} p311, Grigore and Scharf \cite{GrSc2003} p644, that extra selection conditions beyond this are still needed for this to hold. Henneaux \cite{Hen1988} (section 8) gives a discussion of this complication for the case of the string model, along with a discussion of other difficulties related to the superfunction representations commonly used in Hamiltonian BRST, e.g. indefiniteness of the physical inner product and fractional ghost number. Solutions are proposed to both problems. These problems are also noted by other authors e.g. Landsman and Linden \cite{LanLin1992} p425, McMullan and Paterson \cite{McPat1989II} p489 Fuster and van Holten \cite{vHo2006} p7, p9, and essentially the same solution is proposed as in Henneaux \cite{Hen1988}. A useful method for calculating the  physical states which utilizes the \emph{dsp}-decomposition and the BRST Laplacian operator is noted in several places, and treated most systematically for the Hamiltonian framework in van Holten \cite{vHol2005} Section 2.7, and for the Lagrangian framework in Scharf \cite{Schf2001} p21. 

As well as issues related to the state space for both Lagrangian and Hamiltonian BRST, calculation of the quantum  BRST observable algebra is not analysed in detail at the heuristic level. In the Lagrangian picture, K\&O \cite{KuOj79} give Proposition 5.9 p68 which shows heuristically that BRST invariant observables only act on $\ker Q/ \ran Q$ by their restriction to the subspace with no gauge particle creators, Scharf \cite{Schf2001}  p127 proves unitarity of the physical $S$-matrix using a similar argument and the \emph{dsp}-decomposition. In the Hamiltonian picture, H\&T \cite{HenTei92} Section 14.2.1b, gives a structure theorem that relates the \emph{dsp}-decomposition to the operators in $\ker \drb$, but does not compare this to the Dirac method of selection of observables. These results all rely on essentially the same idea, expressed below in Theorem \eqref{pr:krdel} and Theorem \eqref{th:kerdrand}. In \cite{Hen1993I,Hen1993II} a heuristic proof is given to show that  restricting to ghost number zero operators in $\ker \drb$ for the physical observables removes the ghost part of the algebra for Quantum Yang-Mills theory. The proof does not use the structure theorem in H\&T \cite{HenTei92}.  However in Henneaux and Teitelboim \cite{HenTei1987} it is shown that $\drb$ must be restricted to the original  observables (without ghosts) to obtain the correct results for the example of quantization of generalised magnetic monopoles. 

Mathematically consistent treatment of BRST structures is needed to deal with the above issues. As BRST theory has existed for 30 years, work in this direction has of course been done. General BRST structures has been examined by Horuzhy et al. in a series of papers \cite{HoVo89,AzKh89,HoVo92,HoVo292,HoVo93,HoVo97,VoKh2000}, in particular in Horuzhy and Voronin \cite{HoVo89} the BRST charge is analysed $Q$ as a possibly unbounded Hermitian operator acting on a Krein space, and various results including the \emph{dsp}-decomposition are obtained. Ghost number operators are studied in Azizov and Khoruzhii \cite{AzKh89} and conditions for the existence of these operators are given. The natural Lie superalgebra generated by $Q$, its Hilbert-adjoint $Q^*$, ghost number operator $G$ and the BRST-Laplacian $\Delta=\{Q^*,Q\}$ is $l(1,1)$ and its representation theory is carefully analysed in Horuzhy and Voronin \cite{HoVo93} where a general description of how these representations decompose is given. A concrete example of BRST for the Schwinger model is discussed Horuzhy and Voronin \cite{HoVo292} in the Lagrangian approach, where the BRST charge is constructed along the lines of `replacing the gauge parameter by ghost parameter'. A problem is found in applying this approach to the LS solution (as defined in \cite{HoVo292}), but it works well for the CF solution \cite{CapFer1981} with the result that the correct physical space is selected. It should be noted that the BRST charge for the Schwinger models has the same form as that for BRST-QEM when written in terms of creators and annihilators which is why we get the correct physical result. The BRST physical algebra is not discussed in this example, however the selection of the physical algebra is considered in \cite{HoVo92} where it is shown that for theories with structures similar to QEM, the operator cohomology is what we expect. The proof of this fact relies on using an `operator \emph{dsp}-decomposition' and is entirely algebraic. It is conjectured at the end of \cite{HoVo92} that the spatial \emph{dsp}-decomposition can be used to simplify the operator cohomology calculation. We do this here using a structure theorem similar to the one mentioned in H\&T \cite{HenTei92}, but extended to cover the infinite dimensional case. Further papers \cite{HoVo97,VoKh2000} are mainly concerned with studying $l(1,1)$.  

In contrast to the direct study of the general structures of Quantum BRST, a second avenue to investigate Quantum BRST rigorously is to apply a well-defined quantization scheme to classical BRST. As classical BRST is mathematically sound and gives results equivalent to Dirac, it is natural to see if this will carry over to the quantum theory, and this is studied in \cite{BorHerWal2000,DuvElTuy1990,KosSte1987}. Duval, Elhadad and Tuynman \cite{DuvElTuy1990} found that by using Geometric Quantization as well as the BRST procedure gave they obtained the same problems with the indefinite metric of the physical subspace, fractional ghost number, etc, mentioned above with respect to Hamiltonian BRST (\cite{DuvElTuy1990} p543), but could obtain the desired results by modifications made on a case by case by basis. This seems a natural approach to take, however due to obstruction results in quantization \cite{GotHenTuy1996,Got2000}, and the fact that there are quantum constraint systems with no classical analogue,  direct analysis of quantum BRST structures is still of primary importance.

Quantum BRST ideas have been applied to the rigorous formulations of perturbation theory. 
The idea of calculating interacting quantum field theories  (QFT's) via perturbations of free theories is fundamental to all approaches to interacting QFT, and the well known divergences associated to the process are the bane of rigorous formulations of QFT. 
The approach taken in \cite{Schf2001,DuSh1999,DuFred1999,DuFred1998,Holl2008} is to define interacting objects via formal power series in such a way that each term satisfies certain properties we expect from our final theory, such as gauge invariance. Scharf \cite{Schf2001} does this using the Epstein-Glaser \cite{EpGl1973} approach to constructing the $S$-matrix where the terms in the $S$-matrix are smeared time ordered operator valued distributions constructed iteratively without the use of a ultraviolet cut-off function by the use of causality. BRST structures are used to define perturbative gauge invariance (PGI). 
For given models, the $S$-matrix is defined as the formal power series of the appropriate smeared time ordered free matter, gauge, and ghost fields with the $n$-th order term in the series denoted by $(1/n!)T_n(x_1,\dots,x_n)$. Both $Q$ and $\drb$ are constructed using the Lagrangian method of `replacing the gauge parameter by ghost parameter'. PGI is then defined to be that $\drb(T_n(x_1,\dots,x_n))$ is a divergence of particular form (cf. \cite{Schf2001} equation (3.1.18) p99) and is justified as it implies more traditional expressions of PGI, \eg the Ward-Takahashi identities in the case of Quantum Electrodynamics (QED)  (cf. \cite{Schf2001} p100). PGI as defined in \cite{Schf2001} for a general theory also ensures that smearing the $S$-matrix terms $T_n(x_1,\dots,x_n)$ over the Schwartz functions $g\in \cS(\R^4,\R^n)$ and taking the limit $g \to \one$ gives that $\drb(T_n(g))\to 0$ (cf. \cite{Schf2001}  equation (3.1.25) p 100). That is, the terms in the $S$-matrix perturbation series are such that $T_n \in \ker \drb$ in the adiabatic limit. Scharf makes the point that PGI defined using $T_n$ and $\drb$ this way makes sense even when the $S$-matrix does not exist in the adiabatic limit, \eg for the massless field case (\cite{Schf2001} p101). Hence BRST structures, such as the superderivation $\drb$, are fundamental to the meaning of gauge invariance in this approach. Similar ideas are used in \cite{DuFred1998,DuFred1999} where observables, physical state space, and BRST charge $Q_{int}$ for interacting theories are constructed via formal power series derived from the free theory. These ideas are also used in \cite{Holl2008} which provides an extensive analysis of renormalization of Yang-Mills in curved space time. It is assumed throughout these approaches that the physical interpretations require that the physical subspace, $\ker Q/\ran Q$, has positive definite inner product, and also that the physical algebra is defined to be $\ker \drb / \ran \drb$. 

This concludes our brief historical review of the current BRST literature  relevant to this thesis. Some other important works not mentioned so far are e.g. \cite{ZJ1993,Wein2005II,NakOj1990}. For a more extensive guide to the BRST literature, the interested reader can consult consult the bibliographies of \cite{HenTei92,NakOj1990, vHol2005, BarBraHen00,Bat1987}.

From the preceeding discussion it is clear that the quantum BRST constraint method needs a mathematically rigorous treatment for several reasons. We need to investigate positivity of the physical subspace in the Lagrangian and Hamiltonian frameworks. We need to investigate the selection of the algebra of BRST observables. We need to interpret Lagrangian BRST structures in a well-defined mathematical framework, and give an account of free theories such as QEM as in Scharf \cite{Schf2001} in this framework. We need to be able to compare the different approaches to quantum BRST in a mathematically precise way, as well as compare quantum BRST to the quantum Dirac method. 

This thesis pursues these aims and is orgainized broadly in two parts. Chapter \ref{ch:heu}, Chapter \ref{ch:GenStruct} and Chapter \ref{ch:BRSTQEM} make up the first part, and are aimed at developing clear consistent mathematical structures in a concrete representation. Initially, the standard quantum BRST examples are defined formally as found in the literature, then the model independent structures of quantum BRST common to the standard BRST examples in the literature are defined. The standard quantum BRST examples are then developed rigorously in light of these frameworks, facilitating discussion of problematic issues related to BRST in the literature from a consistent viewpoint. 

Much insight can be gained from analysing QFT structures in a representation-free setting, and the second part is aimed at analysing quantum BRST in such a setting. Chapter \ref{ch:CsBRST} encodes the structures of the previous chapters in a $C^*$-algebraic context, and gives an abstract formulation of the general structures of the quantum BRST constraint method. 
 This is a non-trivial mathematical problem for the case of BRST-QEM as most of the basic objects are defined as unbounded operators acting on a Krein-Fock space. The formulation of quantum BRST in $C^*$-algebraic structures allows a comaparison to the quantum Dirac constraint algorithm from a general and rigorous standpoint.

This is not to say that a $C^*$-algebraic framework is the only way to analyse quantum BRST rigorously, as evidenced by the many other approaches to rigorous quantum BRST already mentioned in the introduction. The $C^*$-algebraic viewpoint, however, gives a methodology to draw together and analyse in a consistent mathematical framework coming from the different varieties of quantum BRST.

A summary of the contents of the Chapters is as follows:
 
Chapter \ref{ch:heu} gives an account of the standard heuristic BRST structures such as the BRST charge $Q$, BRST superderivation $\drb$, and ghost structure. The Lagrangian and Hamiltonian BRST approaches will be discussed by giving an example of each. The purpose of this chapter is to motivate the rigorous developments in later chapters and to connect with standard examples as found in the literature. We will see that at the heuristic level, BRST for Quantum Electromagnetism (QEM) selects the correct physical subspace, while for Hamiltonian BRST with a finite number of bounded constraints, BRST has problems with indefiniteness of the physical inner product and that extra selection criteria are required to specify the correct physical space.  

Chapter \ref{ch:GenStruct} makes the heuristic structures rigorous in terms of operators acting on a Krein space. There are several technical difficulties involved in this as many of the basic objects involve algebras of unbounded operators acting on Krein spaces. We first give an account of the ghost algebra as the CAR algebra \cite{BraRob21981} smeared over a Krein test function space with a certain simple structure motivated by the BRST-QEM example. We connect this with the Berezin superfunction representations that are also commonly used for the Ghost Algebra. The natural spatial decomposition associated to the BRST charge $Q$, the \emph{dsp}-decomposition, is given as well as its connection to the physical operator cohomology. Using this we give a simple example of where Hamiltonian BRST will not coincide with the Dirac constraint method at the spatial, or algebra level. 

Chapter \ref{ch:BRSTQEM} is devoted to make the heuristic BRST-QEM example rigorous as well as other Bose-Fock theories with similar one-particle test function space. We find that in this rigorous formulation, the BRST-physical state space and BRST-physical algebra are what we would expect using other constraint methods such as the Gupta-Bleuler method. Other examples fit easily into this framework. Finally, we combine Hamiltonian BRST with a finite number of bounded abelian constraints with abelian Bose-Fock BRST in a natural way which yields a BRST charge $Q$ that selects the correct physical subspace. Hence this yields a BRST-algorithm for a general finite number of bounded abelian constraints that selects the correct physical subspace without the need of imposing any further constraint conditions. This is an interesting result as it is  not the case for the usual Hamiltonian BRST algorithm (cf. MCPS problem as discussed in Subsection \ref{sbs:hhamexphsp} below). To author's knowledge, this algorithm, using the combination of Hamiltonian and Bose-Fock BRST, is an original construction.

Chapter \ref{ch:CsBRST} is then devoted to encoding the previous structures and examples in a $C^*$-algebraic framework. There are major obstacles in doing this such as that the BRST charge $Q$ is unbounded and that the superderivation $\drb$ maps bounded operators to unbounded operators in physical examples such as BRST-QEM. To deal with these we first develop $C^*$-BRST theory with bounded structures to see what features we need to encode. Unbounded superderivations with similar properties to the BRST superderivation have been studied in relation to supersymmetry by Buchholz and Grundling \cite{HendrikBuch2006} using the Resolvent Algebra for the canonical commutation relations developed in Buchholz and Grundling \cite{HendrikBuch2007}. Following a similar procedure we encode BRST for Abelian Bosonic theories in a $C^*$-algebraic way and analyse the resulting structures. We find that when this is done in a way most natural with respect to Poincar{\'e} covariance, the resulting BRST constraint procedure does not remove the ghost part of the algebras. This is a surprise given that the ghosts are removed in BRST-QEM example in Chapter \ref{ch:BRSTQEM}. The reason for this is that we have taken `resolvents' of {Krein}-symmetric fields rather than {Hilbert}-symmetric fields  (Appendix \ref{ap:IIP}) when constructing the corresponding Resolvent Algebra, and so the structures involved do not directly correspond to those in the QEM example. However we do find that Dirac constraining ($T$-procedure Appendix \ref{app:Tp}) produces the expected results. We modify the Resolvent Algebra so that it corresponds exactly with the bosonic example in Chapter \ref{ch:BRSTQEM} and find that BRST-constraining now gives the correct results for strongly regular states (Theorem \ref{th:brstpalgra2}), but that now there are problems in defining the Poincar{\'e} transformations at the algebraic level. 

Finally, using the bounded BRST case and QEM examples as a guide, we give a general formulation of what a $C^*$-algebraic BRST theory looks like,  conditions for states to give physical representations, how the physical algebra is defined at the algebra level, and how the examples developed fit into the $C^*$-framework.

We will use the following conventions in this thesis:
\begin{itemize}
\item The metric tensor for Lorentz space is $g_{\mu \nu}=diag(-1,1,1,1)$.
\item The mantle of the positive light cone is:
\[
C_{+}:=\{ p \in \mathbb{R}^4\,|\, p_{\mu}p^{\mu}=0, \,p_0>0\},
\]
\item Inner products $\ip{\cdot}{\cdot}$ are anti-linear in the first argument and linear in the second.
\end{itemize}

\chapter{Heuristics}\label{ch:heu}
The aim of this chapter is to an outline of the heuristic structures which are to be made rigorous. As this chapter is heuristic, calculations will be formal in nature, with technicalities such as domains of definition and boundedness issues postponed to the following chapters. We will also relegate long proofs to an appendix so as not to interupt the flow of the text. In following the chapters where we are concerned with the technical details, all proofs will be given with the statement of the fact. We first give an account of the features common to the main quantum BRST approaches followed by examples of the Lagarangian and Hamiltonian BRST. The first example given is heuristic BRST quantum electromagnetism (BRST-QEM) done in the Lagragnian approach as found in Scharf \cite{Schf2001}, the second being the example of the Hamiltonian set-up with a finite number of constraints which `close' as given in \cite{HenTei92,vHol2005,LanLin1992,McPat1989I,McPat1989II,Bat1987}. Finally we will discuss the structures that need to be made rigorous, and issues that need to be resolved in the following chapters. 

\section{Set-up}\label{sec:heursetup}
First, BRST begins as an extension procedure, where a  quantum system with constraints is extended by `ghost variables', a description of which is as follows. 

We start with an involutive field algebra $\cA_0$, with involution $\dag$, of possibly unbounded operators which describes the physical observables of the system. This acts on a inner product space $\cH_0$, with (possibly indefinite) inner product $\iip{\cdot}{\cdot}_0$, where $\dag$ is the adjoint with respect to $\iip{\cdot}{\cdot}_0$. We then extend $\cA_0$ by tensoring on a `ghost algebra', 
\[
\cA_g:= \alg{\gh_j, \, \cgh_j \,| j \in I }
\]
where $I$ is some index set determined by the constraints, and the variables have the following anticommutation relations,
\begin{equation}\label{eq:hghstrel}
\{ \gh_j, \cgh_k \}= \delta_{jk} \one, \qquad \{ \gh_j, \gh_k \}=\{ \cgh_j, \cgh_k \}=0,
\end{equation}
We refer to the $\gh_j$ as \emph{ghost variables} and $\cgh_j$ as \emph{conjugate ghost variables}. The ghosts are chosen so that they are associated with the `unphysical degrees of freedom' of $\cA_0$ in some manner, such as associating a ghost to each linearly independent Dirac constraint.

Furthermore, we assume that $\cA_g$ acts as operators on a vector space $\cH_g$, called the \emph{ghost space}, with indefinite inner product $\iip{\cdot}{\cdot}_g$. The involution with respect to this inner product is also denoted $\dag$ and we assume the following hermicity assignments,
\[
\gh_j^{\dag}=\gh_j, \qquad \cgh_j^{\dag}=\cgh_j.
\]
The total extended BRST algebra is $\cA=\cA_0 \otimes \cA_g$ and it acts on $\cH=\cH_0 \otimes \cH_g$, with the usual tensor inner product $\iip{\cdot}{\cdot}$. The adjoint with respect to $\iip{\cdot}{\cdot}$ is given on the elementary tensors as $(A\otimes B)^{\dag}:=A^{\dag} \otimes B^{\dag}$ and extended accordingly.

\begin{rem}\label{rm:heugh}
\begin{itemize}
\item[(i)]Usually the conjugate ghost variables are taken to be anti-$\dag$ hermitian (\cite{HenTei92} p190). That is $\tilde{\cgh}^{\dag}=-\tilde{\cgh}$ and $\{ \gh_j, \cgh_k \}= \delta_{jk} i\one$. By redefining $\cgh:= i \tilde{\cgh}$ we get the above equivalent algebra.
\item[(ii)] The equations \eqref{eq:hghstrel} and hermicity assignments of the ghosts imply that $\iip{\cdot}{\cdot}$ must be indefinite. To see this first note that $\gh_j^2= (1/2) \{ \gh_j, \gh_j \}=0$ and so $\gh_j \cgh_j \gh_j= \gh_j( \one- \gh_j \cgh_j)=\gh_j$. Similarly, $\cgh_j \gh_j \cgh_j= \cgh_j$, and therefore if we let $A=i(\gh_j \cgh_j- \cgh_j \gh_j)$ we have that $A^{\dag}=A$ and $A^2= -(\gh_j \cgh_j- \cgh_j \gh_j)^2=-(\gh_j \cgh_j\gh_j \cgh_j+\cgh_j \gh_j\cgh_j \gh_j)= -\{\gh_j, \cgh_j\}=-\one$. So for $\psi \in \cH_g$ we have $\iip{A\psi}{A\psi}=\iip{\psi}{A^2\psi}=-\iip{\psi}{\psi}$, and so $\cH_g$ cannot be definite or semi-definite.
\item[(iii)]
$\cH$ and $\cH_g$ are almost always assumed to be a Krein space with a Hilbert-topology (see Appendix \eqref{ap:IIP}). An explicit construction of the ghost algebra and $\cH_g$ will be given in chapter \eqref{sec:ghost}, where $\cA_g$ is a representation of the CAR algebra over a particular test function space. We will also show the connection to the usual representation on super function space as given in \cite{HenTei92} p495.
\end{itemize}
\end{rem}

The full algebra $\cA$ is then given a ghost grading in the following manner. Suppose that $M \in \cA$ is a monomial of elements in $\cA_0$, ghosts, and conjugate ghosts. Then $M$ is given ghost number $gh(M)$= number of $\gh$ - number of $\cgh$, and we define $\cG_n=\sn{ M\in \cA \,|\, gh(M)=n}\subset \cA$, the finite linear combinations of such $M$. We then get that $\cA= \sn{ A\in \cG_n\,|\, n \in \mathbb{Z}}$.  This grading is referred to as the \emph{ghost grading}, and the above is well defined as a consequence of the relations \eqref{eq:hghstrel}, which will be seen more explicitly in subsection \eqref{sbs:ghgrad}. 

We also have a $\bZ_2$ grading of $\cA$ given by $\cA^{+}:= \sn{A\in\cG_{2n}\,|\,n \in \mathbb{Z}}$ and $\cA^{-}:= \sn{A\in\cG_{2n+1}\,|\,n \in \mathbb{Z}}$. $\cA^{+}$ and $\cA^{-}$ are referred to as the \emph{even} and \emph{odd} subalgebras respectively, and it is not hard to see that the map $\gamma: \cA^{+} \cup \cA^{-} \to \cA$ given by:
\[
\gamma(A_0):=A_0, \, A_0 \in  \cA^{+}, \qquad \gamma(A_1):=-A_1, \, A_1 \in  \cA^{-},
\]
extends to an automorphism on $\cA$ such that $\gamma^2= \iota$. We call this the $\bZ_2$-grading automorphism, and is essential for constructing the BRST superderivation below.

There is also commonly assumed to be a \emph{ghost number} operator, $G\in \op(\cH_g)$. This is an operator such that, $[\one \otimes G,A]=nA$ for $A \in \cG_n$, and with our ghost hermicity assignment, $G^{\dag}=-G$. The existence of this will be discussed in subsection \eqref{sbs:ghgrad}

Furthermore $\cH_g$ is assumed to have a ghost grading, which is given by the eigenspaces of $G$, \ie $\cH_{gn}:=\{ \psi\in \cH_g \,|\, G\psi=n \psi\}$. Using $G^{\dag}=-G$, we see that $\iip{\psi_n}{\psi_m}_g=0$ for $\psi_n \in \cH_{g,n}$, $\psi_m \in \cH_{g,m}$, $n \neq -m$. In particular, we see that $\cH_{gn}$ are neutral subspaces for $n \neq 0$. A good discussion of ghost number structure and operators is given in \cite{AzKh89}, and we will discuss it in Subsection \ref{sbs:ghgrad}.

\section{Charge}\label{sec:heurcha}

After extending $\cA_0$ by the ghost variables $\cA_g$, an element $Q\in \cA$ is constructed from the constraints or gauge transformations in such a way that $Q^2=0$, $Q^\dag=Q$, and $Q\in \cG_1$. $Q$ is one of the central objects in heuristic BRST theory and is called the {\emph{BRST Charge}}. The purpose of $Q$ is that it is used as a constraint condition to determine the physical subspace in $\cH$ by:
\[
\cH^{BRST}_{c}:= \ker Q,
\]
Various arguments are used to derive formulas for $Q$. We will come back to these formulas and reasoning in the following examples, but first we list basic facts about structures derived from such $Q$:

\pagebreak

\begin{enumerate}
\item $Q^2=0$ implies that $\ran Q \subset \ker Q$. $\cH$ is almost always a Krein space, for which we assume the Hilbert-topology (see Appendix \eqref{ap:IIP}) on $\cH$.  $Q$ is usually closed with respect to this hence $\ker Q$ is closed (\cite{Con1985} Proposition X 1.13 p310) and $\overline{\ran Q} \subset \ker Q$.
\item $Q^{\dag}=Q$ implies that $\overline {\ran Q} [\perp] \ker Q$, where $[\perp]$ denotes $\iip{\cdot}{\cdot}$-orthogonality. In particular, as $\overline{\ran Q} \subset \ker Q$ we have that $\overline{\ran Q}$ is a neutral subspace. 
\item Define an indefinite inner product on $\ker Q/ \overline{\ran Q}$ as follows. Let $\hat{\psi}:= \psi + \overline{\ran Q}$ and $\hat{\vp}:= \vp + \overline{\ran Q}$, then
\begin{equation}\label{eq:hiipkQ}
\iip{\hat{\psi}}{\hat{\vp}}_p:=\iip{\psi}{\vp}.
\end{equation}
As $\overline{\ran Q}$ is neutral this is independent of the choice of representatives $\psi, \vp$. Therefore $\ker Q/ \overline{\ran Q}$ will be positive only when $\ker Q$ is a positive semidefinite subspace of $\cH$. 
\item Let $J$ be the fundamental symmetry for $\cH$ (see Appendix \eqref{ap:IIP}). As $Q^{\dag}=Q$, we get from lemma \eqref{lm:khad} that $Q^*=JQ^{\dag}J=JQJ$. It follows from this, the Krein structure of $\cH$, and $Q^2=0$ that $\cH$ has the decomposition,
\begin{align}\label{eq:hdsp}
\cH = & \; \cH_d \oplus \cH_s \oplus \cH_p,
\end{align}
where the direct sum $\oplus$ is with respect to the Hilbert inner product, and $\cH_s:=\ker Q \cap \ker Q^*$, $\cH_d:=\overline{ \ran Q}$, $\cH_p:=\overline {\ran Q^*}$, $\ker Q= \cH_d \oplus \cH_s$. This extremely useful result is in the literature in many places \cite{HenTei92} p328, \cite{vHol2005} p44, \cite{HoVo89} p682, \cite{Schf2001} p22, and following \cite{Nish1984} it will be referred to as the $dsp$-decomposition. Note that is also commonly referred to as the Hodge decomposition (\cite{vHol2005} p44, \cite{HenTei92} p328). The proof of this is given in section \eqref{sec:dsprig}, which is also in \cite{HoVo89}, and covers the case where $Q$ can be unbounded.
\item From the $dsp$-decomposition it follows that $\overline{ \ran Q}$ is the isotropic subspace of $\ker Q=\cH_s\oplus \cH_p$. Hence $\ker Q / \overline{\ran Q}$ is a Krein space with the indefinite inner product \eqref{eq:hiipkQ} (see Appendix \eqref{ap:IIP}). If $\ker Q$ is  positive semidefinite, then $\overline{ \ran Q}$ is its neutral subspace (see Appendix \eqref{ap:IIP}), hence the indefinite inner product \eqref{eq:hiipkQ} is positive definite.
\item Note that $\cH_s=\ker Q \cap \ker Q^*= \ker (Q^*Q + QQ^*)= \ker \Delta$ by positivity, where $\Delta:=\{Q, Q^*\}$ 
\end{enumerate}

From the above remarks, we see that $\ker Q$ is degenerate with repect to $\iip{\cdot}{\cdot}$ with isotropic part $\overline{\ran Q}$, and so the final physical state space is defined to be,
\[
\cH^{BRST}_{phys}:= \ker Q / \overline{\ran Q}.
\]
If $\cH^{BRST}_{phys}$ is to be an acceptable physical space, we would like it to have positive definite inner product, and as discussed above, this will be true only when $\ker Q$ is a positive semidefinite subspace. As we shall see, this is not always true in simple examples such as the example below for abelian Hamiltonian BRST (Remark \eqref{rm:indhpssp}).

Now $\ker Q$ and $\cH^{BRST}_{phys}$ are usually not easy to calculate in examples, however from the $dsp$-decomposition we see that,
\[
\ker Q= \ker \Delta \oplus \overline{\ran Q}.
\]
So there is a natural identification of $\cH^{BRST}_{phys}$ with $\ker \Delta$, and in many examples it turns out that $\ker \Delta$ is straightforward to calculate in contrast to $\ker Q$. 
\begin{rem}
The natural identification above motivates some authors to say that $\ker \Delta$ and $\cH^{BRST}_{phys}$ are equivalent definitions of the physical subspace (\cite{Schf2001} p22). However in examples,  there are natural $\beta \in \aut(\cA)$ that $\dag$-automorphisms but not $*$-automorphisms and are such that $\beta(Q)=Q$ but $\beta(Q^*)\neq Q^*$, e.g. the relativistic boosts for BRST-QEM.

Hence structures associated to $\alpha$ naturally factor to $\ker Q /\overline{\ran Q}$ but not $\ker \Delta$. Also it is not true in general that elements $A \in \cA$ that factor to $\hat{A} \in \op(\cH^{BRST}_{phys})$ are the same as the elements $A \in \cA$ that preserve $\ker \Delta$, and so the different choices of the physical subspace may have different physical algebras associated to them. For both these reasons we do not use the definition of  $\ker \Delta$ as the physical subspace.
\end{rem}

\section{Superderivation}\label{sbs:hedsperv}
Another fundamental feature of the BRST machinery is the BRST superderivation (see Appendix \eqref{ap:SS} for superalgebra). If $\gamma$ is the $\bZ_2$-grading automorphism on $\cA$, then this is defined as,
\[
\drb(A)=\sbr{Q}{A}:=QA-\gamma(A)Q, \qquad A \in \cA.
\]
and we have that,
\[
\drb(AB)=\drb(A)B+ \gamma(A)\drb(B).
\]

\pagebreak

Basic facts about $\drb$:
\begin{enumerate}
\item It follows from $Q^2=0$, and the super-Jacobi identity (see Appendix \eqref{ap:SS}), that $\drb^2=0$ and we have that $\ran \drb \subset \ker \drb$. 
\item As $Q^{\dag}=Q$ we have that $\ker \drb$ and $\ran \drb$ are both $\dag$-algebras, but as $Q^* \neq Q$  they are not $*$-algebras.
\item $Q \in \cA^{-}$ implies that $\gamma(Q)=-Q$ and hence $\gamma \circ \drb \circ \gamma= -\drb$. Therefore $\gamma(\ran \drb)=\ran \drb$ and $\gamma (\ker \drb) = \ker \drb$.
\item For $B=\drb(C) \in \ran \drb$, $A \in \ker \drb$ we have that by the preceding remark that $\gamma(A) \in \ker \drb$ and so $AB= \gamma^2(A)\drb(C)=\drb(\gamma(A)C)\in \ran \drb$. Similarly $BA \in \ran \drb$ and so we see that $\ran \drb \unlhd \ker \drb$, where $\unlhd$ denotes that $\ran \drb$ is an ideal of $\ker \drb$. To check if this is a proper containment, it turns out that we have to check whether $\cH_s$ is trivial or not. This non-trivial result is Corollary \eqref{cr:oneinker} below.
\item It is easy to check that for $A \in \ker \drb$ we have $A \ker Q \subset \ker Q$ and $A \ran Q \subset \ran Q$, and for $B \in \ran \drb$, $B \ker Q \subset \ran Q$. Therefore we have that an $A \in \ker \drb$ factors through to an operator on $\cH^{BRST}_{phys}$ and $A \in \ran \drb$ factors trivially. 
\end{enumerate}

Now as we are taking $\cH^{BRST}_{phys}= \ker Q / \overline{\ran Q}$ as our physical subspace, we see that by the above comments $\ker \drb$ will lift to the physical subpace and $\ran \drb$ will lift trivially. As we also have that $\ran \drb \unlhd \ker \drb$, we define the heuristic physical algebra as,
\[
\cP^{BRST}:= \ker \drb / \ran \drb.
\]

As we are still dealing with heuristics, we will not concern ourselves with topologies on $\ker \drb$ and $\ran \drb$. Also notice that as $\ker \drb$ and $\ran \drb$ are not $*$-algebras, it is not clear that $\cP^{BRST}$ will be a $*$-algebra under the usual quotient inner product. We will discuss this further in Subsection~\ref{sbs:AbcuQ}, where it turns out that under certain assumptions $\cP^{BRST}$ has a $C^*$-involution coming from the Hilbert inner product on $\ker Q / \overline{\ran Q}$.

As in the state space case, $\cP^{BRST}$ turns out to be difficult to calculate, the only complete explicit calculation in an example the author has found is in \cite{HoVo92}. This is done in a very direct manner using some specific properties of the algebras involved in the example, and not using the \emph{dsp}-decomposition. It turns out however that the \emph{dsp}-decomposition can be used to give conditions on when $A\in \ran \drb$ for a general BRST set-up, with $\cH_s$ playing a pivotal role (see Corollary \eqref{cr:oneinker}). This provides a more economical proof of the result in \cite{HoVo92}, and can be used more generally to show facts such as $\one \in \ran \drb$ iff $\cH_s = \{0\}$ (theorem \eqref{pr:krdel}).

\begin{rem}
The definition of the physical algebra is not uniform, when it is defined at all. The above cohomological definition is motivated by the classical BRST theory, as discussed in the introduction. The above definition coincides with the Hamiltonian approach, H\&T \cite{HenTei92} p301, \cite{McPat1989II} p490, and some authors using the Lagrangian approach \cite{HoVo92} p1319.
 
Originally a more Dirac-like approach was assumed by K\&O \cite{KuOj79} p59, which they proved to be equivalent to the above definition for \emph{local} observables in the field theories they were considering (Chapter 5.1-5.2 in particular proposition 5.4 and proposition 5.9). We will return to this issue in Chapter \eqref{sbs:alalg}.
\end{rem}

With the basic objects of BRST defined in a heuristic fashion, we want to investigate the correspondence between the results using $Q$ and $\drb$ and those using the traditional Dirac method. With this in mind we give two examples which  correspond to the two main approaches to BRST, Lagrangian BRST and Hamiltonian BRST. These examples will be analysed in a mathematically precise framework in the following chapters. At this stage we will only investigate the state space claim, leaving the analysis of the algebra until we have developed enough rigourous theory for a meaningful discussion in Section \eqref{sbs:AbcuQ}.

The first example is heuristic BRST quantum electromagnetism (BRST-QEM) as presented in Scharf \cite{Schf2001} Chapter 1. This is an example of Lagrangian BRST as described in Kugo \& Ojima \cite{KuOj79} p19, with abelian gauge groups, and in this case we find that we get the expected state space.

The second example is that of Hamiltonian BRST with Lie algebra of constraints. Hamilitonian BRST also used to cover the case where the constraints do not `close', \ie form a Lie algebra, but the example given is sufficient to highlight the important features of Hamiltonian BRST, some of which are problematic, without adding unnecessary complication. We will see that for this example the BRST state space is bigger than the Dirac state space, and has indefinite inner product, and that even adding the extra condition of restricting to the ghost zero subspace does not fix this problem. Therefore the approaches in these two examples are not equivalent. 

Following the examples there is a discussion of their differing features. The next chapter then begins development of a rigorous setup that incorporates these examples.

\section{Example 1: Heuristic BRST-QEM}\label{sec:exheuEM}
The following is a version of heuristic BRST-Quantum Electromagnetism (BRST-QEM) as given in Scharf \cite{Schf2001} Chapter 1. This is an example of Lagrangian BRST, where the basic idea is, given a gauge field theory, we `replace the gauge parameter by ghost parameter' K\&O \cite{KuOj79} p16. This is a fairly vague concept given that rigorous interacting quantum gauge field theory is still an open problem, but we describe the process for electromagnetism below. 

Note that \cite{Schf2001} does not make reference to BRST but has all features described in the preceeding sections. Moreover this treatment is essentially the same as that in K\&O \cite{KuOj79} Chapter 2, for abelian gauge group ($C^c_{ab}=0$), the differences essentially being in labelling. We use Scharf's treatment for which the labelling is straightforward, and it makes use of the Hilbert space involution $*$.

\subsection{Basic fields}\label{sbsc:alg}
For our basic ingredients, we follow \cite{Schf2001} p9-p11, and assume that we have a bosonic field on $\mathbb{R}^3$ (momentum space) with the following CCR relations,
\begin{gather}\label{eq:hccr}
[a^{\mu}(\bq)\,,\,a^{\nu}(\bp)] = 0, \qquad [a^{\mu}(\bq)^*\,,\,a^{\nu}(\bp)^{*}] = 0, \qquad \bp,\bq \in \mathbb{R}^3\\ 
[a^{\mu}(\bq)\,,\,a^{\nu}(\bp)^{*}] = \delta_{\mu}^{\nu} \delta^3(\bp-\bq),
\end{gather}
and a fermionic field with the following CAR relations,
\begin{gather}\label{h:car}
\{ c_j(\bp),c_k(\bq) \} = 0,\qquad \{ c_j(\bp)^*,c_k(\bq)^* \} = 0,\qquad \bp,\bq \in \mathbb{R}^3\\
\{ c_j(\bp),c_k(\bq)^* \} = \delta_{jk} \delta^3(\bp-\bq), \qquad j,k=1,2.
\end{gather}
These act on the respective Fock spaces $\cH_0, \, \cH_g$ in the usual way, with the $*$ denoting Hilbert involution. Note that the above are actually operator valued distributions and so to make sense of them mathematically we will smear them over a one-particle test function space. As this chapter is aimed at describing formal structures, we continue in the above manner leaving the mathematical precise treatment for Chapter \ref{ch:BRSTQEM}. 

We assume further that we can take tensor products of the fields above acting on  $\cH=\cH_0 \otimes \cH_g$. A Krein space structure and involution $\dag$ will be defined on $\cH$ below. Following the usual convention we will not explicitly use  $\otimes$ to denote the operator products. 

In this example, we will not define $\cA_0$ and $\cA_g$ as they would correspond to pointwise products of operator valued distributions.

\subsection{The QEM field}\label{sbs:heuQEMfield}
We follow \cite{Schf2001} p10 {but} with different metric convention $g:=diag(-1,1,1,1,)$, 
\begin{align}\label{eq:Adef}
A^{k}(x)&=\FC\int \frac{d^3 \bp}{\sqrt{ p_0}} ( a^{k}(\bp)e^{-ip x} + a^{k}(\bp)^{*} e^{ip x} ), \quad p_0=|\mathbf{\bp}|,\, k =1,2,3, \\
\intertext{and,}
A^{0}(x)&=\FC \int \frac{d^3 \bp}{\sqrt{ p_0}} ( a^{0}(\bp)e^{-ip x} - a^{0}(\bp)^{*} e^{ip x} ), \quad p_0=|\mathbf{\bp}|.
\notag
\end{align}
Note that $A^{0}(x)^*=-A^{0}(x)$ which reflects the fact that the $*$-involution is the Hilbert involution rather than the Krein involution.

So far we have been considering $\cH_p$ as the Hilbert-Fock space with the  fields $a^{\mu}(\bp),a^{\mu}(\bq)^{*}$ acting as annihilation and creation operators. However the Gupta-Bleuler (GB) version of EM is usually defined in terms of creation and annihilation operators acting on a Krein-Fock space (see Schweber \cite{Sch64} p246) . To see the connection we define a fundamental symmetry $J \in B(\cH_0)$ so that  $\cH_0$ becomes the GB Krein-Fock space and $A^{0}(x)^{\dag}=JA^{0}(x)^*J=A^{0}(x)$. 

As we are in Fock space we  assume that we have an $N_0$ which is a number operator for the $a^{0}(\bp)$ and that $J= (-1)^{N_0}$ (so $J$ is bounded). It is easy to see that,
\[
J^*=J,\quad J^2=\one.
\]
Let $\ip{\cdot}{\cdot}$ be the Hilbert inner product on $\cH_0$, for which $*$ is the adjoint, then we define a new inner product as,
\[
\iip{\cdot}{\cdot}:=\ip{\cdot}{J\cdot},
\]
and denote its adjoint by $\dag$, and we have that $F^{\dag}=JF^*J$, $F\in \op(\cH_0)$ (see Appendix \eqref{ap:IIP}).

Now it is straightforward to see that
\[
a^{0}(\bp)^*J= -Ja^{0}( \bp)^*, \qquad a^{k}(\bp)^*J= Ja^{k}(\bp)^* ,\quad k = 1,2,3, 
\]
and so in terms of $\dag$-involutions we have
\begin{gather}\label{eq:daginv}
 a^{0}(\bp)^{\dag}=Ja^{0}(\bp)^*J= -a^{0}(\bp)^*  \\
a^{k}(\bp)^{\dag}=Ja^{k}(\bp)^*J= a^{k}(\bp)^*,\quad k = 1,2,3. \notag
\end{gather}
From the definitions \eqref{eq:Adef} we get that
\begin{align}\label{eq:Adefdag}
A^{\mu}(x)&=\FC \int \frac{d^3\bp}{\sqrt{p_0}} ( a^{\mu}(\bp)e^{-ipx} + a^{\mu}(\bp)^{\dag} e^{ipx} ), \quad p_0=|\mathbf{\bp}|,\,, \mu =0,1,2,3, 
\end{align}
All the $A^{\mu}(x)$'s are Krein ($\dag$) hermitian, but we now have that
\begin{equation}\label{eq:dagheuccr}
[a_{\mu}(\bp)\,,\,a_{\nu}(\bq)^{\dag}] =  g_{\mu \nu}\delta(\bp-\bq)
\end{equation}
and so we do not have the usual CCR's for the $a^{\mu}(p)$'s with involution $\dag$.

 If,
\begin{equation}\label{eq:JPdist}
D_0(x):=(2\pi)^{-3} \int_{C_+}\frac{d^3\bp}{p_0} e^{-i\bp.\mathbf{x}}\sin(p_0 x_0),
\end{equation}
then we can formally verify that: 
\begin{gather}
[ A^{\mu}(x)\,,\, A^{\nu}(y)] = -i g^{\mu \nu} D_0(x-y), \label{eq:guv} \\
\square A^{\mu}(x)=0, \notag \\
[Q_0\,,\, A^\nu(x)] = -i\partial^\nu \Lambda (x), \notag
\end{gather}
where $\Lambda \in \mathcal{S}(\mathbb{R}^4), \, \square \Lambda=0$ and, 
\begin{align}\label{eq:hemch}
Q_0&= \int_{x_0=const} d^3\mathbf{x}\,[ (\partial_\nu A^\nu ) \partial_0 \Lambda - ( \partial_0 \partial _\nu A^\nu) \Lambda],\\
&:=\int_{x_0=const} d^3\mathbf{x} \,(\partial_\nu A^\nu )\overleftrightarrow{\partial_0} \Lambda .\notag 
\end{align}

\begin{proof}See Appendix \eqref{ap:herpfs}
\end{proof}
That is, we have the correct commutation relations for $A^{\mu}(x)$ and $Q_0$ generates the gauge transformations. Note that the definition of $A^{0}(x)$ is made to ensure we get a $g^{\mu \nu}$ factor in \eqref{eq:guv} rather than $\delta_{\nu}^{\mu}$ (as pointed out in \cite{Schf2001}, p12).

To get \index{QEM} we have to further impose Maxwell's equations. Let $F^{\mu \nu}(x)=A^{\mu,\, \nu}(x)-A^{\nu,\, \mu}(x)$ then the Maxwell's equations are ${F^{\mu \nu}}_{,\nu}=0$. In Gupta-Bleuler QEM we also impose the Lorentz condition by by first selecting $\cH'=\ker p^{\mu}a_{\mu}(\bp)$ for all $\bp=0$, $p_0=\np$ and then defining the physical state space as $\cH^{GB}_{phys}:=\cH'/\cH''$ where $\cH''$ is the neutral part of $\cH'$ with respect to the indefinite inner product on $\cH_0$. See Schweber \cite{Sch64} p240-246 for details.

\subsection{BRST Extension of QEM}\label{sbs:hrext}
The idea in Lagrangian BRST is to replace the `gauge parameter' $\Lambda(x)$ in \eqref{eq:hemch} by an anti-commuting `ghost parameter' $\gh(x)$, such that,
\begin{equation}\label{eq:gcon}
\{\gh(x)\,, \, \gh(y) \}=0, \qquad \square \gh(x) =0,
\end{equation}
The BRST charge $Q$ is then the generator of the gauge transformations with `ghost parameter', and we will see below that $Q^2=0$. The condition $\square \gh(x) =0$ is motivated by $\square \Lambda=0$.

To this end, using the fields $c(\bp)$ in \eqref{h:car}, we define:
\begin{equation}\label{eq:udef}
\gh(x)=\FC \int \frac{d^3\bp}{\sqrt{p_0}} ( c_2(\bp)e^{-ipx} + c_1(\bp)^{*} e^{ipx} ), \quad p_0= |\bp |. 
\end{equation}
We can check that $\gh(x)$ satisfies \eqref{eq:gcon}.
\begin{rem}\label{rm:hgfcomrel}
\begin{itemize}
\item[(i)]  Note that we do not use $c_1(p)$ and $c_1(p)^*$ in the definition of $\gh(x)$ as this will not give $\{\gh(x) \,,\, \gh(y) \}=0$. This is in contrast to the CCR field $A_{\mu}(x)$.
\item [(ii)] We want a Krein structure such that $\gh(x)$ is Krein Hermitian, to get features such as $\bch^{\dag}=\bch$. To this end we choose a fundamental $J_g$ symmetry on $\cH_g$ such that $c_2(\bp)^{\dag}=c_1(\bp)^{*}$. We will postpone the construction of $J_g$ to section \eqref{sec:ghost}.
\item [(iii)]Scharf defines a second ghost field as,
\begin{equation}\label{eq:halcg}
\tilde{\gh}(x)=\FC \int \frac{d^3\bp}{\sqrt{p_0}} (- c_1(\bp)e^{-ipx} + c_2(\bp)^{*} e^{ipx} ), \quad p_0= |\bp |,
\end{equation}
(see \cite{Schf2001}, p9) and notes that this is \textit{not} the $*$-adjoint of $\gh(x)$. Using the assignment $c_2(\bp)^{\dag}=c_1(\bp)^{*}$ we see that $\tilde{\gh}(x)^{\dag}=-\tilde{\gh}(x)$. Moreover, we have the anticommutation relations
\begin{equation}\label{eq:gfheq}
\{\gh(x),\tilde{\gh}(y)\}=-iD_0(x-y)
\end{equation}
\begin{proof} See Appendix \eqref{ap:herpfs}\end{proof}
Also note that $\tilde{\gh}(x)$ is not the conjugate ghost field, $\cgh(x)$, in our treatment, but it is related as will be seen in Remark \eqref{rm:concghsegh}.
\end{itemize}
\end{rem}
So now, replacing $\Lambda(x)$ by $\gh(x)$ in \eqref{eq:hemch}, we define, 
\begin{align*}
Q&:=\int_{x_0=const} d^3\mathbf{x} (\partial_\nu A^\nu )\overleftrightarrow{\partial_0} \gh, 
\end{align*}
and note that $Q^{\dag}=Q$, and that now, 
\begin{equation}\label{eq:HQ}
Q= \int d^3\bp\; p_0 [(a_{\Vert}(\bp)^* + a_0(\bp)^*)c_2(\bp) + (a_{\Vert}(\bp) - a_0(\bp))c_1(\bp)^*],\quad p_0= |\bp |,
\end{equation}
where, 
\[
a_{\Vert}(\bp)= (p_k/p_0) a^k(\bp), \qquad k=1,2,3
\]
where we are using the convention $p_k x^k=\mathbf{p}.\mathbf{x}$.
\begin{proof} See Appendix \eqref{ap:herpfs}.\end{proof} 
From this it follows that,
\begin{align}
[Q\,,\, A^\nu(x)] &= -i\partial^\nu \gh(x), \label{eq:hdA}\\
\{Q\,,\, \gh(x)\} &= 0, \label{eq:hdgh}\\
\{Q\,,\, \tilde{\gh}(x)\} &= -i\partial_\nu A^\nu(x), \label{eq:hdcgh}\\
Q^2&=0 \label{eq:h2np}
\end{align}
\begin{proof}See Appendix \eqref{ap:herpfs}\end{proof}

Now we have formally defined BRST-QEM and replaced the gauge parameter by a ghost parameter. To complete the analysis we need to find the physical state space. Usually we would also want to calculate the physical observables, but as we have not defined the BRST superderivation $\drb$ in this formal example we will leave this till we have developed the appropriate machinery (subsection \eqref{sbs:phalopem}).

\subsection{$\Delta$ and State space}
We now turn to the task of calculating the physical state space. Since the $dsp$-decomposition (equation \eqref{eq:hdsp} ) gives $\ker Q= \ker\Delta\oplus \overline{\ran Q}$, we now calculate $\Delta$. 
To simplify notation we define,
\begin{align*}
b_1(\bp) &= \frac{1}{\sqrt{2}}(a_{\Vert}(\bp) - a_0(\bp))\\
b_2(\bp) &= \frac{1}{\sqrt{2}}(a_{\Vert}(\bp) + a_0(\bp)),
\end{align*}
which satisfy the usual commutation relations.
\begin{equation}\label{eq:bcomm2}
[b_i(\mathbf{q})\,,\,b_j(\bp)^{*}] = \delta_{ij} \delta^3(\bp-\mathbf{q})
\end{equation}
and so we get that,
\begin{align}\label{eq:Qpdef}
Q&= \sqrt{2} \int d^3\bp\; p_0 [b_{2}(\bp)^* c_2(\bp) + c_1(\bp)^* b_1(\bp)], \\
Q^*&= \sqrt{2} \int d^3\bp\; p_0 [b_{2}(\bp) c_2(\bp)^* + c_1(\bp) b_1(\bp)^*],\notag
\end{align}
This is essentially the form of $\bch$ given in K\&O \cite{KuOj79} p19,  H\&T\cite{HenTei92} p464, \cite{Sl1989}p584,  \cite{HoVo92,HoVo292}, \cite{RazRyb1990}p218, and we use it now to calculate,
\begin{align}\label{eq:HDel}
\Delta= & \;\{Q,Q^*\},\notag\\
= & \;2\int d^3\bp \;p_{0}^{2} (b_1(\bp)^*b_1(\bp) +b_2(\bp)^*b_2(\bp) +c_1^*(\bp)c_1(\bp) +c_2^*(\bp)c_2(\bp)) 
\end{align}
\begin{proof} See Appendix \eqref{ap:herpfs}. \end{proof}

We see that $\Delta$ looks like a number operator with $\ker \Delta= \cap_{i=1,2}\{\ker b_i(\bp) \cap \ker c_i(\bp)\,|\, \bp \in \mathbb{R}^3 \}$, that is kernel consisting of states with no $b(\bp)$'s or $c(\bp)$'s. 

More explicitly, let 
\[
\cF:=\{ f:C(\mathbb{R}^3) \to C(\mathbb{R}^4)   \,|\, f_0=0,\,\bp. \mathbf{f}(\bp)=0,  \, \mathbf{f}(\bp). \mathbf{f}(\bp)=1 \}.
\]
Let $a_f(\bp)= f(\bp)_\mu a^{\mu}(\bp)=  \mathbf{f}(\bp). \mathbf{a}(\bp)$. Then we have  that for $f\in \cF$,
\begin{align*}
[a_f(\bp),a_{f}^{*}(\bq)]&=\mathbf{f}(\bp). \mathbf{f}(\bp)\delta(\bp-\bq) =\delta(\bp-\bq),\\
[a_f(\bp),b_i^*(\bq)]&=(1/(p_0\sqrt{2}))\bp.\mathbf{f}(\bp) \delta(\bp-\bq)=0, \qquad i=1,2,
\end{align*}
and all other commutators are zero. This is not suprising as $a_f(\bp), a^*_f(\bp)$ are the `transversal photons'. Now we also have that formally $\cH$ is the Fock space  generated by on the vacuum state $\Omega_0 \otimes \Omega_g$ with $b_1^*(\bp), b_2^*(\bq), a_f^*(\bp), c_1^*(\bp), c_2^*(\bp)$'s. This statement is made precise in Subsection~\eqref{sbs:testfunc} using Proposition \eqref{lm:stdec2} where the bosonic fields are smeared over $\fH=\fH_t \oplus \fHL \oplus \fHJ$ with $\fH_t,\fHL,\fHJ $ corresponding to the $a_{f}(\bp)$'s, $b_1(\bp)$'s, $b_2(\bp)$'s above and the ghost fields smeared over $\fL=\fL_1 \oplus \fL_2$ where $\fL_1, \fL_2$ correspond to $c_1(\bp), c_2(\bp)$ respectively.

So we see from \eqref{eq:HDel} that $\ker \Delta$ is the Fock space generated by the $a_f^*(\bp)$'s, $f\in \cF$, or the space of transversal photons which we would expect from QEM.

To see that the above construction corresponds to QEM we need show that the Maxwell equations are satisfied. We have formally,
\begin{align}
\drb(F^{\mu \nu}(x))= & \;0, \label{eq:MxOb}\\
{F^{\mu \nu}}_{,\nu}(x)= & \;i \drb( \partial^{\mu} \tilde{\gh}(x)). \label{eq:Mxcond}
\end{align}
with $\tilde{gh}(x)$ as in equation \eqref{eq:halcg}.
\begin{proof} See Appendix\eqref{ap:herpfs}. \end{proof}
Therefore we have that formally $F^{\mu \nu}(x)$ factors to $\ker Q / \ran Q$ and that ${F^{\mu \nu}}_{,\nu}(x)$ factors trivially to $\ker Q / \ran Q$ and thus the Maxwell conditions hold on the BRST physical state space. 

\section{Example 2: Hamiltonian BRST} \label{sec:hamBRST}
\subsection{Setup}
Now we give the example of Hamiltonian BRST with closed constraint algebra, which  follows the treatment in H\&T \cite{HenTei92} chapter 14, as well as many others, e.g. \cite{vHol2005,LanLin1992,McPat1989I,McPat1989II}. 

We begin with a field algebra $\cA_0$  acting on the Hilbert space $\cH_0$ with inner product $\ip{\cdot}{\cdot}_{0}$. We are told that we should `not pay attention to the the functional subtleties arising from the infinite dimensionalitiy of the space of states' and `proceed as if the space of states were finite-dimensional', \cite{HenTei92} p302. Hence we will assume that $\cA_0$ is a $C^*$-algebra and so acts on all of $\cH_0$. Let $\cC=\{G_i\,|\, i=1,\dots,n\}\subset \cA_0$ be a finite set of self-adjoint, linearly independent elements in $\cA_0$ that form a Lie algebra, ie,
\begin{equation}\label{eq:conbrak}
[G_a,G_b]=iC^c_{ab}G_c, \qquad C^c_{ab} \in \mathbb{R}.
\end{equation}
Note $C^c_{ab}=-C^c_{ab}$, and self-adjointness of the constraints imply that we must have $C^c_{ab} \in \mathbb{R}$. 

The class of Hamiltonian BRST systems includes the case where the constraints do not `close', \ie the `structure constants' $C^c_{ab}$ are actually operators \cite{HenTei92} p317, but this will not be discussed further as the example given is sufficient for our purposes.

We assume now that $\cC$ is our constraint set, which selects the physical state space by $\cH^0_p:=\cap_{i=1}^{n}\ker G_i$. We then extend $\cA_0$ by tensoring on a ghost algebra (cf. \eqref{eq:hghstrel}), with a ghost for every constraint, ie
\[
\cA=\cA_0 \otimes \cA_g, \qquad \cA_g=\alg{ \gh_a, \cgh_b \,|\, a=1, \dots, n},
\] 
which act on a Krein space $\cH_g$. Let the Krein-inner product on $\cH_g$ be given by $\iip{\cdot}{\cdot}_g$ which has fundamental symmetry $J_g$. Therefore $\ip{\cdot}{\cdot}_g:=\iip{\cdot}{J\cdot}_g$ is a Hilbert-inner product, with Hilbert involution denoted by $*$. We let $\cA= \cA_0 \otimes \cA_g$ act on the Hilbert space $\cH:=\cH_0 \otimes \cH_g$. We assume that $\cH$ has the usual tensor inner product and define the Krein inner product by using $J:=\one \otimes J_g$ as a fundamental symmetry.

Following \cite{vHol2005} equation (2.116), we assume the Hilbert-adjointness properties for $\cA_g$ as
\begin{equation}\label{eq:hgadj}
\gh_a^*=\cgh_a.
\end{equation} 
That this is consistent will be shown in Chapter \eqref{sec:ghost}, and it is impicitly assumed in the usual Berezin representation which will be discussed in Section \eqref{sbs:brz}.

We now construct the BRST charge as,
\begin{equation}\label{eq:hnonabQ}
Q=G_a \otimes \gh_a - (i/2) C^c_{ab}\one \otimes \gh_a \gh_b \cgh_c,
\end{equation}
where we use convention of summing over repeated indices. This formula is given in numerous places, e.g. \cite{HenTei92} p197, \cite{vHol2005}, p48 \cite{LanLin1992}, p423. 

It is obvious that $gh(Q)=1$, it is straightforward to check that $Q^{\dag}=Q$, and with some work we check that 
\begin{equation}\label{eq:hHam2nilQ}
Q^2=0
\end{equation}
\begin{proof}see Appendix \eqref{ap:herpfs}. 
\end{proof}
Therefore $Q$ fulfils the criteria for a BRST charge.  For motivation for this definition of $Q$ via classical theory see H\&T \cite{HenTei92} Chapters 1-12. 

Another motivation for $Q$ is that we begin with constraints $G_a$ and we want to construct a 2-nilpotent constraint to develop BRST-cohomology. We then take  $G_a \otimes \gh_a$, but as  $[G_a,G_b]=iC^c_{ab}G_c$ we get 2-nilpotence only when the structure constants vanish, hence we need to add a term to cancel the terms coming from the non-commuting constraints $G_a$ and this is $- (i/2) C^c_{ab}\one \otimes \gh_a \gh_b \cgh_c$.

\subsection{Physical subspace and $\ker \Delta$}\label{sbs:hhamexphsp}
Recall that $\Delta:= \{Q, Q^* \}$, and that the $dsp$-decomposition gives a natural isomorphism  $\ker Q/ \overline{\ran Q}\cong \ker \Delta$, by $\ker Q= \ker \Delta \oplus \overline{\ran Q}$. Hence we would like to calculate $\ker \Delta$. 

\medskip
\noindent{\textbf{Case 1: Abelian Constraints, $C^c_{ab}=0$.}}
\smallskip

We consider first the simplest case where all the constraints commute, \ie $C^c_{ab}=0$. In this case we have that,
\[
Q=G_a\otimes \gh_a, \qquad Q^*= G_a \otimes \cgh_a.
\]
\begin{rem}\label{rm:indhpssp}
We show that $\ker Q$ is \textit{not} a semidefinite subspace of $\cH$: Note that $\cH^0_p \otimes \cH_g \subset \ker Q$. So take a $\psi_g \in \cH_g$ such that $\iip{\psi_g}{\psi_g}>0$. As we have shown in Remark \eqref{rm:heugh} (ii) we have that $\iip{A\psi_g}{A\psi_g}_g<0$ where $A=i(\gh_1 \cgh_1-\cgh_1 \gh_1)$, and hence for $\psi=\psi_0\otimes A\psi_g$, $\iip{\psi}{\psi}<0$. This is problematic as the final physical space  $\ker Q/ \overline{\ran Q}$ now has an indefinite inner product. This problem is noted in \cite{HenTei92} p311, where it is stated that an extra constraint condition maybe needed in specific examples. One possibility given is to restrict to states of ghost number $0$. However, if we take $\psi_g$ to have $gh(\psi_g)=0$ then we see that $gh(A \psi_g)=0$  and hence the same problem arises.
\end{rem}

Now $\cH^0_p \otimes \cH_g \subset \ker Q$, and we will now show that in fact $\cH^0_p \otimes \cH_g =\ker \Delta$. Note that:
\begin{align}
\Delta=\{Q,Q^*\}&=\{G_a \otimes \gh_a,G_b \otimes \cgh_b\}, \label{eq:calhdelham}\\
		&=[G_a \otimes \gh_a,G_b\otimes \one](\one \otimes \cgh_b)+(G_b\otimes \one)\{G_a\otimes \gh_a,\one \otimes \cgh_b\},\notag \\
		&=0+G_bG_b\otimes \one \notag ,
\end{align}
where the second equality follows from $\{A,BC\}=[A,B]C+B\{A,C\}$, and the third from $\{\gh_a, \cgh_b\}= \delta^a_b \one$. So we see that,
\begin{equation}\label{eq:hcomhamdel}
\ker \Delta= \ker (\sum_{a=1}^{n} G_a^2\otimes \one)=\bigcap_{a=1}^{n}\ker (G_a\otimes \one)= \cH^0_p \otimes \cH_g.
\end{equation}

This is a problematic result since if we identify $\cH^0_p$ with $\cH^0_p\otimes \psi$ for any nonzero $\psi \in \cH_g$ then we have a copy of $\cH^0_p$ in $\ker \Delta$ for every nonzero $\psi \in \cH_g$, and so $\ker \Delta$ selects a larger constraint space than $\cC$. We will refer to this problem as the \emph{Multiple Copies of the Physical Space Problem} (MCPS problem). The MCPS problem is noted in \cite{vHol2005} p50 (who proves the same for the non-abelian case) and suggests that we need to fix the ghost number of physical states to the correct identification of the physical state space. However, as already noted at the end of Subsection \ref{sec:heursetup}, $G^{\dag}=-G$ implies that subspaces with nonzero ghost number, are neutral with respect to $\iip{\cdot}{\cdot}$. Hence if we restrict to vectors of a particular ghost number then this ghost number must be zero. As discussed in the remark above, this is not a good restriction. We will return to the MCPS problem and other choices for the physical space in Subsection \eqref{sbs:brz} Remark \eqref{rm:spfstuff}.

Also in \cite{HenTei92} p311 it is mentioned that extra conditions beyond $Q\psi=0$ and $\one \otimes G\psi=0$ may be needed, a simple example of which is discussed in \cite{HenTei92} p330 Excercise 14.16. We will return to these issues in section \eqref{sbs:brz}  remark \eqref{rm:spfstuff}. 

\medskip
\noindent{\textbf{Case 2: Non-abelian Constraints, $C^c_{ab}\neq0$.}}
\smallskip

We consider as above in \eqref{eq:conbrak} the case that the constraints form a Lie algebra. For simplicity, further assume that $C^c_{ab}\neq 0$ is antisymmetric in all indices as this makes the following proofs simpler. The more general case of semi-simple Lie algebra is in \cite{vHol2005,vHo1990} and we are essentially following the treatment there. Following \cite{vHol2005} p50, we define,
\begin{equation}\label{eq:sigantbrs}
\Sigma_a:=-\one \otimes iC^c_{ab} \gh_b \cgh_c, \qquad \Sigma_a^{\dag}=\Sigma_a
\end{equation}
We also have,
\[
\Sigma_a^*=\one \otimes iC^c_{ab} \gh_c \cgh_b=-\one \otimes iC^b_{ac}\gh_c \cgh_b=\Sigma_a,
\]
where we used $\gh_a^{*}=\cgh_a$. Then,
\begin{align}\label{eq:nonabdel}
\Delta= & \;(1/2)G_aG_a \otimes \one + (1/2)( G_a \otimes \one +\Sigma_a)(G_a \otimes \one +\Sigma_a),\notag\\
\end{align}
\begin{proof}
See Appendix \eqref{ap:herpfs}.
\end{proof}
As $G_a$ and $(G_a \otimes \one +\Sigma_a)$ are $*$-self adjoint operators, we see that $\Delta$ is the sum of $2n$ positive operators and hence,
\begin{align*}
\ker \Delta= & \; \cap_{a} (\ker (G_a\otimes \one) \cap \ker( G_a \otimes \one +\Sigma_a)),\\
= & \; \cap_{a} (\ker (G_a\otimes \one) \cap \ker (\Sigma_a)),\\
= & \; (\cH^0_p \otimes \cH_g)\cap (\cap_{a}\ker (\Sigma_a)).
\end{align*}
So we see that in the non-abelian case we have the extra constraints $\Sigma_a$ on $\ker \Delta$. As $\ker \Delta \cong \cH^{BRST}_{phys}$ through the \emph{dsp}-decomposition this is an extra constraint on the physical space. We will see in section \eqref{sbs:fndmgh} that $\cap_{a}\ker (\Sigma_a)|_{\cH_g}$ is not one dimensional so again we have many copies of the $\cH^0_p \otimes \cH_g)$ in $\ker \Delta$ and also that the restriction of $\iip{\cdot}{\cdot}$ to $\ker \Delta$  is indefinite (see Remark \eqref{rm:spfstuff}).

\section{Discussion of the Examples}
We now list the major positives and negatives of the above examples and then compare them to each other. The Lagrangian-BRST-QEM example has positive features: 
\begin{itemize}
\item It formally corresponds to other versions of QEM. It gives a state space isomorphic to the transversal photons, which will will be made precise in Section \ref{sc:EM} Theorem \eqref{th:GBalg}, where we see that this space is isomorphic to that given by the Gupta-Bleuler method in \cite{Hendrik2000} p36. We have that $\drb(F^{\mu \nu}(x))=0$ and ${F^{\mu \nu}}_{,\nu}(x)=i \drb( \partial^{\mu} \tilde{\gh}(x))$, so the field operators $F^{\mu \nu}(x)$ factor to the physical subspace, and the Maxwell equations hold as ${F^{\mu \nu}}_{,\nu}(x)$ factors trivially. 
\item It models the example of free QEM, an example containing many physical structures.
\end{itemize}
and difficulties:
\begin{itemize}
\item Many of the basic objects are defined only heuristically: e.g. $Q$ is defined as an integral of products of operator valued distributions; the fields $A(x)$ are unbounded distributions and domains of definition need to be specified; we need a sensible definition of our ghost algebra to give meaning to odd and even elements; the algebras involved involve products of operator valued distributions so defining $\drb$ and its domain needs care.  
\item Once defined on a specific representation, to move to a $C^*$-algebraic setting, we need to encode the unbounded operator algebras in some bounded fashion. We also need to encode $\drb$ in a bounded fashion which is not a straightforward problem as it maps bounded elements to unbounded elements, see Definition \eqref{df:brstsd} below. The natural involution is the Krein-$\dag$ involution with associated objects such as $\dag$-algebras not Hilbert-$*$-algebras. Hence general $C^*$-theory will not always apply directly but has to be arrived at by through the relation of $\dag$-involution to $*$-involution via $Ad(J)$. 
\item The BRST extension of replacing the `gauge parameter' $\Lambda$ by `ghost parameter' $\gh$ relies on gauge theory structure and is not well defined in an algorithmic sense. Hence generalization is needed to be able to apply BRST to other theories with constraints.
\end{itemize}
The Hamiltonian BRST example has the following positive features:
\begin{itemize}
\item It only depends on the number and bracket structure of the constraints for the ghost algebra extension and  construction of $Q$ and so lends itself to a general constraint algorithm.
\item The operators involved are in $B(\cH)$ and so the structures can be encoded straightforwardly in a $C^*$-algebraic setting.
\item  We have an algorithm for dealing with the case where the constraints do not commute.
\end{itemize}
It has the following negative features:
\begin{itemize}
\item The constraint condition does not remove the ghosts from the physical subspace. Therefore extended and unextended systems will not be equivalent.
\item  The physical constraint space does not have  positive definite or semidefinite inner product, and also suffers from the MCPS problem.
\item This leads to the BRST-physical algebra being not equivalent to the unextend constrained algebra for basic examples, as we shall see in Subsection \ref{sbs:egAHC}. 
\item Adding the extra condition that physical states and algebras must be of ghost number zero does not solve the above problems.  
\end{itemize}

Comparing the two examples we see that they both fit the general description outlined in Section \ref{sec:heurcha}, Section \ref{sec:heursetup}, Section \ref{sbs:hedsperv}, but differ in crucial aspects: 

A major difference between these two examples is that in the BRST-QEM example $\Delta$ removes the ghost space, whereas in the Hamiltonian BRST example $\Delta$ does not. This stems from the fact that while the $Q$'s in the different frameworks are cosmetically similar, they are derived from different reasoning. We will see below that this leads to BRST-QEM removing the ghosts from the observable algebra, but not for Hamiltonian BRST in general.

The Hamiltonian example is more clearly defined than BRST-QEM example from a mathematical standpoint. In Hamiltonian BRST example both $Q$, $\drb$ are constructed as finite sums of bounded operators, the original field algebra is a $C^*$-algebra and as we shall see in the  Section \ref{sec:ghost} the ghost algebra is easily seen to be a finite dimensional CAR algebra.  In contrast, the basic objects in the BRST-QEM example, such as $Q$, involve integrals of products of operator valued distributions and the field algebras are not clearly defined, however the BRST-QEM example models many physical structures and hence needs a mathematically rigorous treatment.

\section{General Discussion and Directions}
In light of the previous discussion we can discuss the issues related to quantum BRST in greater depth as well as the general directions of this thesis.
 
We need to give an account of the general BRST structures descriped at the beginning of this chapter that is general enough to incorporate both examples given above. In particular, we need to give a description of the Ghost Algebra that will  describe the smeared ghosts of the BRST-QEM example as well as the connect this to the finite dimensional superfunction representations used in Hamiltonian BRST. We need to give the \emph{dsp}-decomposition for the case of unbounded BRST charge $Q$, and also show its realation to the operator cohomology for the case of unbounded Krein-$\dag$ algebras. There are results in these directions as discussed in the introduction, but we need extend these and unify them in order to incorporate both examples.

With the BRST structures well defined, in particular the ghost structures, we can discuss problems for Hamiltonian BRST related to the ghost space representations: ghost number zero restriction, neutral subspace and the MCPS problem. Although discussions related to these problems can be found in the literature in several places e.g. (\cite{LanLin1992} p425-429, \cite{McPat1989II} p489,  \cite{Hen1988} Section 8, H\&T \cite{HenTei92} Sections 14.2.5, 14.2.6, 14.5.3, 14.5.4), I feel that these need to be collected, highlighted, and expanded as it is not always clear from the literature that these are \emph{always} problems faced by Hamiltonian BRST with heuristic set-up as in Section \ref{sec:exheuEM}. 
The example of BRST applied to a constraint set consisting of a finite number of canonical momentum variables is often assumed to fit into heuristic Hamiltonian BRST e.g. H\&T \cite{HenTei92} Section 14.4, and we do not have neutrality or MCPS problems with this example at the heuristic level. However we cannot ignore boundedness issues in this example and it turns out that this example fits better into a Lagrangian framework rather than Hamiltonian. This is for reasons such as that the momentum constraints are Krein-symmetric rather than Hilbert-symmetric and quantum momentum variables with canonical conjugates have trivial kernels. See Subsection \ref{sbs:findimbos} for more details.

We need a mathematically precise formulation of the BRST-QEM example so we can compare the results to other rigorous versions of quantum QEM. We do this first at the level of Krein representations in Fock-Krein space, where we actually give an account that covers more general free abelian bosonic theories with setup similar to \KO.

With BRST done for Hamiltonian and BRST-Bose-Fock theories, it is natural to see if we can combine the two approaches in a way that captures the positive features of both approaches, but avoids the negative.  We can do this by combining abelian Hamiltonian BRST with a finite number constraints with the abelian BRST-Bose-Fock set-up to get an general BRST algorithm such that selects $\Delta$ the correct state space.

We want to give a $C^*$-algebraic formulation of BRST so as to connect with the algebraic approach to QFT. We have many difficulties as already mentioned in the introduction, but we can do this in terms of an unbounded superderivation acting on a subalgebra $D_2(\drb)$ of a $C^*$-algebra $\cA$. A key to this construction is the ghost algebra as developed in Section \ref{sec:ghost} , and the Resolvent Algebra as developed in \cite{HendrikBuch2007}. We identify the interesting states on the algebra, and how the general structures from Chapter \ref{ch:GenStruct}   are incorporated in the $C^*$-framework. We then apply these structures to the examples developed.

There are also major equivalence issues: how do Lagrangian and Hamiltonian quantum BRST relate to each other and when do they give equivalent results; and when will quantum BRST give equivalent results to the Dirac algorithm. Using the general BRST structures above, we can answer these questions on a case by case basis for given examples where the ghost extensions have been made and $\drb$, $Q$ explicitly constructed. However these are difficult questions to anwer in general and we discuss this further at the end of this thesis when we have developed and explored the examples and structures to a greater extent.

\section{Appendix: Proofs}\label{ap:herpfs}

\begin{proof}(Proof of equations \eqref{eq:guv}).

Recall $g=diag(-1,1,1,1)$, $[a_{\mu}(\bp)\,,\,a_{\nu}(\bq)^{\dag}] =  g_{\mu \nu}\delta(\bp-\bq)$, and
\[
A^{\mu}(x)=\FC\int \frac{d^3 \bp}{\sqrt{ p_0}} ( a^{\mu}(\bp)e^{-ip x} +  a^{\mu}(\bp)^{\dag} e^{ip x} ), \quad p_0=|\mathbf{\bp}|,\,, \mu =0,1,2,3. 
\]
In the following calculations we assume that $p_0=|\mathbf{\bp}|,q_0=|\mathbf{\bq}|$.
\begin{align*}
(i) \,[ A^{\mu}(x)\,,& \, A^{\nu}(y)]  = \\
&=\FCC \int \frac{d^3 \bp}{\sqrt{ p_0}} \int \frac{d^3 \bq}{\sqrt{ q_0}}([a^{\mu}(\bp),a^{\nu}(\bq)^{\dag}]e^{-ipx+iqy} - [a^{\nu}(\bq), a^{\mu}(\bp)^{\dag}] e^{-iqy+ip x} ) \\
&=\FCC g^{\mu \nu}\int \frac{d^3 \bp}{{ p_0}} (e^{-ip(x-y)} -  e^{ip(x- y)} ), \\
&=\FCC g^{\mu \nu}\left(   \int \frac{d^3 \bp}{{ p_0}}e^{ip_0(x_0-y_0)}e^{-i\bp.(\mathbf{x- y})}- \int \frac{d^3 \bp}{{ p_0}} e^{-ip_0(x_0-y_0)}e^{i\bp.(\mathbf{x-y})}\right), \\
&=-i(2\pi)^{-3}g^{\mu \nu}\int \frac{d^3 \bp}{{ p_0}}e^{-i\bp.(\mathbf{x-y})}\sin(p_0(x_0-y_0)),\\
&=-ig^{\mu \nu}D_0(x-y)
\end{align*}
 
\noindent(ii) This follows from $\square e^{-ipx}= -(p_0^2 -  |\mathbf{\bp}|^2) e^{-ipx}$
 
\noindent(iii) We use the notation $\partial_{y\mu}:=(\partial/\partial y^{\mu})$.
\begin{align*}
[Q_0,A^{\nu}(x)]= & \;\int_{y_0=const} d^3\mathbf{y}\,\partial_{y\mu}[  A^\mu(y),A^{\nu}(x)] \partial_{y0} \Lambda(y) -  \partial_{y0} \partial _{y\mu} [A^\mu(y), A^{\nu}(x)]\Lambda(y),\\
= & \;i\int_{y_0=const} d^3\mathbf{y}\,\big(\partial_{y\mu} g^{\mu \nu} D_0(x-y)\big) \partial_{y0} \Lambda(y) -  \big(\partial_{y0} \partial _{y\mu}  g^{\mu \nu} D_0(x-y)\big)\Lambda(y),\\
= & \;i\int_{y_0=const} d^3\mathbf{y}\,\big(\partial^{y\nu} D_0(x-y)\big) \partial_{y0} \Lambda(y) -  \big(\partial_{y0} \partial ^{y\nu}  D_0(x-y)\big)\Lambda(y),\\
= & \;-i\int_{y_0=const} d^3\mathbf{y}\, D_0(x-y)\big( \partial_{y0}\partial^{y\nu} \Lambda(y)\big) -  \big(\partial_{y0}  D_0(x-y)\big)\partial ^{y\nu} \Lambda(y).
\end{align*}
This is the solution to the Cauchy problem of the 3-dimensional wave equation with initial data $f_1(\mathbf{x})=\partial ^{\nu} \Lambda(0,\mathbf{x})$ and $f_2(\mathbf{x})=\partial_{0} f_1(\mathbf{x})=\partial_{0}\partial ^{\nu} \Lambda(0,\mathbf{x})$ and so we get that,
\[
[Q_0,A^{\nu}(x)]=-i\partial^{\nu}\Lambda(x),
\]
(see Pinsky \cite{Pin1998} p325 with $c=1$).
\end{proof}
\begin{proof}(Proof of equation \eqref{eq:gfheq})Recall the defintions of $\gh(x)$ and $\tilde{\gh}(x)$ in equation \eqref{eq:udef} and equation \eqref{eq:halcg} respectively, and the CAR's 
\[
\{ c_j(\bp),c_k(\bq)^* \} = \delta_{jk} \delta^3(\bp-\bq), \qquad j,k=1,2
\]
in equation \eqref{h:car}. In the following calculations we assume that $p_0=|\mathbf{\bp}|,q_0=|\mathbf{\bq}|$:
\begin{align*} 
\{\gh(x),\tilde{\gh}(y)\}&=\FCC \int \frac{d^3 \bp}{\sqrt{ p_0}} \int \frac{d^3 \bq}{\sqrt{ q_0}}(-\{c_1(p)^*,c_1(q)\}e^{ipx-iqy} + \{c_2(p),c_2(q)^*\} e^{-ipx+iq y} ) \\
&=\FCC \int \frac{d^3 \bp}{{ p_0}} (e^{-ip(x-y)} -  e^{ip(x- y)} ), \\
&=-iD_0(x-y)
\end{align*}
where the last equality follows as the calculation of equation \eqref{eq:guv} (i) above.
\end{proof}
\begin{proof}(Proof of equation \eqref{eq:HQ}) Note that as $qx:=-x_0q_0+(\mathbf{x}.\bq)$ we have $\partial_{\mu}a^{\mu}(\bq)e^{-iqx}=-iq_{\mu}a^{\mu}(\bq)e^{-iqx}$. Therefore, using  \eqref{eq:udef} and \eqref{eq:Adefdag} we get,
\begin{align*}
Q= & \;\int_{x_0=const} d^3\mathbf{x} \,(\partial_\nu A^\nu ){\partial_0} \gh-{\partial_0} (\partial_\nu A^\nu )\gh, \\
&=\FCC \int_{x_0=const}\frac{d^3\bp}{\sqrt{p_0}} \frac{d^3\bq}{\sqrt{q_0}}d^3 \mathbf{x} \, \,I(x,q,p), \qquad p_0= |\bp |, q_0= |\bq |
\end{align*}
where,
\begin{align*}
I(x,q,p)= & \;\partial_{\mu} (a^{\mu}(\bq)e^{-iqx} + a^{\mu}(\bq)^{\dag} e^{iqx}) \partial_0( c_2(\bp)e^{-ipx} + c_1(\bp)^{*} e^{ipx} )\\
&-\partial_0( \partial_{\mu} (a^{\mu}(\bq)e^{-iqx} + a^{\mu}(\bq)^{\dag} e^{iqx})) ( c_2(\bp)e^{-ipx} + c_1(\bp)^{*} e^{ipx} ),\\
= & \; (-iq_{\mu}a^{\mu}(\bq)e^{-iqx} + iq_{\mu}a^{\mu}(\bq)^{\dag} e^{iqx}) ( -ip_0c_2(\bp)e^{-ipx} + ip_0c_1(\bp)^{*} e^{ipx} )\\
&-(-q_0q_{\mu}a^{\mu}(\bq)e^{-iqx} - q_0q_{\mu}a^{\mu}(\bq)^{\dag} e^{iqx}) ( c_2(\bp)e^{-ipx} + c_1(\bp)^{*} e^{ipx} ),\\
= & \; (p_0+q_0)q_{\mu}a^{\mu}(\bq)c_1(\bp)^{*}e^{-i(q-p)x} + (p_0+q_0)q_{\mu}a^{\mu}(\bq)^{\dag}c_2(\bp) e^{i(q-p)x}) \\
&+(q_0-p_0)q_{\mu}a^{\mu}(\bq)c_2(\bp)e^{-i(q+p)x} + (q_0-p_0)q_{\mu}a^{\mu}(\bq)^{\dag}c_1(\bp)^{*} e^{i(q+p)x}).
\end{align*}
Now using $ (2\pi)^{-3} \int_{x_0=const}d^3 \mathbf{x} \, e^{i(p\pm q)x} = e^{i(p_0\pm q_0)x_0}\delta(\bp \pm \bq)$, and the condition $p_0=|\bp|$, $q_0=|\bq|$, we get,
\begin{align*}
Q= & \;\FCC 
\int \frac{d^3\bp}{\sqrt{p_0}} \frac{d^3\bq}{\sqrt{q_0}}\int_{x_0=const}d^3 \mathbf{x} \, \,I(x,q,p), \qquad p_0= |\bp |, q_0= |\bq |,\\
= & \;\FCC \int \frac{d^3\bp}{{p_0}}2(2\pi)^{3} (p_0p_{\mu}a^{\mu}(\bq)c_1(\bp)^{*}+p_0p_{\mu}a^{\mu}(\bp)^{\dag}c_2(\bp)),\qquad p_0= |\bp |\\
\intertext{and by the relations \eqref{eq:daginv},}
= & \;\int d^3\bp\, p_0[  (a_{\Vert}(\bp) - a_0(\bp))c_1(\bp)^*+(a_{\Vert}(\bp)^* + a_0(\bp)^*)c_2(\bp)], \qquad p_0= |\bp |.
\end{align*}
\end{proof}

\begin{proof}(Proof of equations \eqref{eq:hdA}, \eqref{eq:hdcgh}, \eqref{eq:hdgh}, \eqref{eq:h2np}.)
If we use equation \eqref{eq:HQ} to write,
\begin{equation*}
Q= \int d^3\bp\; p_0 [(a_{\Vert}(\bp)^* + a_0(\bp)^*)c_2(\bp) + (a_{\Vert}(\bp) - a_0(\bp))c_1(\bp)^*],
\end{equation*}
where $a_{\Vert}(\bp)= (p_j/p_0) a^j(\bp)$. Then it follows from the defining (anti-)commutation relations \eqref{h:car}, \eqref{eq:hccr}, that
\begin{gather}\label{eq:hbrsderosc}
[Q\,,\,a^{\mu}(\bp)] = -p^{\mu}c_2(\bp), \qquad [Q\,,\,a^{\mu}(\bp)^{\dag}] =p^{\mu}c_1(\bp)^*,\\ 
\{Q,c_1(\bp)\}= p_0(a_{\Vert}(\bp) - a_0(\bp)),\qquad \{Q,c_1^*(\bp)\}=0,\\
\{Q,c_2^*(\bp)\}= p_0(a_{\Vert}^*(\bp) + a_0^*(\bp)),\qquad \{Q,c_2(\bp)\}=0.
\end{gather}

To prove $Q^2=0$ \eqref{eq:h2np} we write $Q$ as,
\begin{align*}
Q=\int \frac{d^3\bp}{{p_0}} (p_0p_{\mu}a^{\mu}(\bq)c_1(\bp)^{*}+p_0p_{\mu}a^{\mu}(\bp)^{\dag}c_2(\bp)),\qquad p_0= |\bp |
\end{align*}
and so,
\begin{align*}
2Q^2=\{Q,Q\}=\int {d^3\bp} (\{Q,p_{\mu}a^{\mu}(\bq)c_1(\bp)^{*}\}+\{Q,p_{\mu}a^{\mu}(\bp)^{\dag}c_2(\bp)\}),\qquad p_0= |\bp |.
\end{align*}
Now using $\{Q,AB\}=[Q,A]B+A\{Q,B\}$  and the above relations we get,
\begin{align*}
\{Q,p_{\mu}a^{\mu}(\bq)c_1(\bp)^{*}\}= & \;[Q,p_{\mu}a^{\mu}(\bq)]c_1(\bp)^{*}+p_{\mu}a^{\mu}(\bq)\{Q,c_1(\bp)^{*}\},\\
= & \;-p_{\mu}p^{\mu}c_2(\bp)c_1(\bp)^{*}
\end{align*}
and similarly $\{Q,p_{\mu}a^{\mu}(\bp)^{\dag}c_2(\bp)\}=p_{\mu}p^{\mu}c_1^{*}(\bp)c_2(\bp)$. Substituting back into $2Q^2$ we get,
\begin{align*}
2Q^2= & \;\int {d^3\bp} (-p_{\mu}p^{\mu}c_2(\bp)c_1(\bp)^{*}+ p_{\mu}p^{\mu}c_1^{*}(\bp)c_2(\bp)),\qquad p_0= |\bp |,\\
= & \;0
\end{align*}
where the last line follows as $p_{\mu}p^{\mu}=0$ for $p_0= |\bp |$.

Equations \eqref{eq:hdA}, \eqref{eq:hdcgh}, \eqref{eq:hdgh}, follow in a similar manner. Note that \cite{Schf2001} p20 has a different proof of $Q^2=0$.
\end{proof}
\begin{proof}(Proof of \eqref{eq:HDel})
Using equations\eqref{eq:Qpdef} we have,
\[
\Delta = \{Q\,,\,Q^* \} =2 \int \int d^3\bp \, p_0 \, d^3\bq  \,q_0\, (\{b_2(\mathbf{q})^*c_2(\mathbf{q}) + c_1(\mathbf{q})^*b_1(\mathbf{q})\,,\,b_2(\bp)c_2(\bp)^* + c_1(\bp)b_1(\bp)^*\}).
\]
Now,
\[
\{b_2(\mathbf{q})^*c_2(\mathbf{q})\,,\,c_1(\bp)b_1(\bp)^*\} =b_2(\mathbf{q})^*b_1(\bp)^*\{c_2(\mathbf{q})\,,\,c_1(\bp)\}=0
\]
so,
\[
\Delta =2\int \int d^3\mathbf{q}\, d^3\bp \, q_0 \, p_0 (\{b_2^*(\mathbf{q})c_2(\mathbf{q}) \,,\,b_2(\bp)c_2(\bp)^*\} + \{c_1(\bp)b_1(\bp)^*\,,\, c_1^*(\mathbf{q})b_1(\mathbf{q})\}),
\]
and,
\begin{align}\label{eq:horrnd}
\{b_2(\mathbf{q})^*c_2(\mathbf{q}) &\,,\,b_2(\bp)c_2(\bp)^*\}= \\ = & \;b_2(\mathbf{q})^*b_2(\bp)c_2(\mathbf{q})c_2(\bp)^* +b_2(\bp)b_2(\mathbf{q})^*c_2(\bp)^*c_2(\mathbf{q}) \notag \\
= & \;b_2(\mathbf{q})^*b_2(\bp)c_2(\mathbf{q})c_2(\bp)^* +b_2(\mathbf{q})^*b_2(\bp)c_2(\bp)^*c_2(\mathbf{q})
-b_2(\mathbf{q})^*b_2(\bp)c_2(\bp)^*c_2(\mathbf{q}) \notag \\ &+b_2(\bp)b_2(\mathbf{q})^*c_2(\bp)^*c_2(\mathbf{q})
\notag \\
= & \;b_2(\mathbf{q})^*b_2(\bp) \{c_2(\mathbf{q})\,,\,c_2(\bp)^* \} +[b_2(\bp)\,,\,b_2(\mathbf{q})^*]c_2(\bp)^*c_2(\mathbf{q}) \notag \\
= & \;[b_2(\mathbf{q})^*b_2(\bp)+ c_2(\bp)^*c_2(\mathbf{q})]\delta( \mathbf{q} -\bp) \notag
\end{align}
Note that the above calculation uses the CCR's \eqref{eq:bcomm2}. Substituting back into $\Delta$ gives,
\begin{equation}
\Delta=2\int d^3\bp\, p_{0}^{2}\, (b_1(\bp)^*b_1(\bp) +b_2(\bp)^*b_2(\bp) +c_1^*(\bp)c_1(\bp) +c_2^*(\bp)c_2(\bp)) 
\end{equation}
\end{proof}

\begin{proof}(Proof of equations \eqref{eq:MxOb} and \eqref{eq:Mxcond})
First,
\begin{align*}
F^{\mu \nu}(x)=A^{\mu,\, \nu}(x)-A^{\nu,\, \mu}&(x),\\
=\FC  &\left( \int \frac{d^3\bp}{\sqrt{p_0}} (-ip^{\nu} a^{\mu}(\bp)e^{-ipx} + ip^{\nu}a^{\mu}(\bp)^{\dag} e^{ipx} ) \right. ,\\
- &  \left.  \int \frac{d^3\bp}{\sqrt{p_0}} ( -ip^{\mu}a^{\nu}(\bp)e^{-ipx} + ip^{\mu}a^{\nu}(\bp)^{\dag} e^{ipx} )\right) ,  \qquad p_0= |\bp |.
\end{align*}
Using the action of $\drb$ given in equations \eqref{eq:hbrsderosc} it follows straightfowardly that $\drb(F^{\mu \nu}(x))=0$ which is equation \eqref{eq:MxOb}.

To prove equation \eqref{eq:Mxcond} we differentiate the above equation to get 
\begin{align*}
{F^{\mu \nu}}_{,\nu}(x)= & \;-\FC \left( \int \frac{d^3\bp}{\sqrt{p_0}} ( p^{\mu}p_{\nu}a^{\nu}(\bp)e^{-ipx} + p^{\mu}p_{\nu}a^{\nu}(\bp)^{\dag} e^{ipx} )\right.\\
&\left.- \int \frac{d^3\bp}{\sqrt{p_0}} (p^{\nu}p_{\nu} a^{\mu}(\bp)e^{-ipx} + p^{\nu}p_{\nu}a^{\mu}(\bp)^{\dag} e^{ipx} )\right),  \qquad p_0= |\bp |,\\
= & \;-\FC \left(\int \frac{d^3\bp}{\sqrt{p_0}} ( p^{\mu}p_{\nu}a^{\nu}(\bp)e^{-ipx} + p^{\mu}p_{\nu}a^{\nu}(\bp)^{\dag} e^{ipx} )\right),  \qquad p_0= |\bp |,\\
\intertext{using $p^{\nu}p_{\nu}=0$ for $p_0=|\bp|$. Furthermore using equations \eqref{eq:hbrsderosc} and $p_{\nu}a^{\nu}=p_0(a_{\Vert}-a_0)$,}
= & \;i\drb \left(\FC \int \frac{d^3\bp}{\sqrt{p_0}} ( ip^{\mu}c_1(\bp)e^{-ipx} + ip^{\mu}c_2^*(\bp) e^{ipx} )\right) ,  \qquad p_0= |\bp |,\\
= & \;i\drb( \partial^{\mu} \tilde{\gh}(x)).
\end{align*} 
where $\tilde{\gh}$ is as in equation \eqref{eq:halcg}.
\end{proof}

\begin{proof}(Proof $Q^2=0$ \eqref{eq:hHam2nilQ})

This is the proof as in masters thesis \cite{Zw1998}. We will repeatedly use the properties of superbrackets $\sbr{\cdot}{\cdot}$ as in Appendix\eqref{ap:SS} and for ease of notation we will drop the $\otimes$ for this proof. Recall that $Q$ is defined by equation\eqref{eq:hnonabQ},
\[
Q=G_a \gh_a - (i/2) C^c_{ab}\gh_a \gh_b \cgh_c,
\]
then,
\begin{align*}
2Q^2=\sbr{Q}{Q}= & \;\sbr{G_a \gh_a}{ G_j \gh_j}-i\sbr{G_a \gh_a}{C^j_{kl}\gh_k \gh_l \cgh_j} -(1/4)\sbr{C^c_{ab}\gh_a \gh_b \cgh_c}{C^j_{kl}\gh_k \gh_l \cgh_j},\\
 = & \; A+B+C,
\end{align*}
where we have that,
\begin{align*}
A=\sbr{G_a \gh_a}{ G_j \gh_j}=\{G_a \gh_a, G_j \gh_j\}
= & \; G_aG_j \gh_a \gh_j +  G_jG_a\gh_j \gh_a,\\
= & \; [G_a,G_j]\gh_a \gh_j+ G_jG_a \{\gh_j, \gh_a\}
=iC^{c}_{ab}G_c\gh_a \gh_j,
\end{align*}
and,
\begin{align*}
B&=-iC^j_{kl} \sbr{ G_a\gh_a}{ \gh_k \gh_l \cgh_j},\\
&=-iC^j_{kl}G_a \sbr{ \gh_a}{ \gh_k \gh_l \cgh_j},\\
&=-iC^j_{kl}G_a (\sbr{ \gh_a}{ \gh_k \gh_l} \cgh_j+\gh_k \gh_l\sbr{ \gh_a}{\cgh_j}),\\
&=-iC^j_{kl}G_a(\gh_k\gh_l \delta_{aj}),\\
&=-iC^a_{kl}G_a\gh_k\gh_l,\\
&=-A.
\end{align*}
Therefore we have that $2Q^2=C$,
\begin{align*}
-4C= & \;C^c_{ab}C^i_{jk}\sbr{\gh_a \gh_b \cgh_c}{\gh_j \gh_k \cgh_i},\\
= & \;C^c_{ab}C^i_{jk} \left( \sbr{\gh_a \gh_b \cgh_c}{ \gh_j \gh_k}\cgh_i+ \gh_j \gh_k \sbr{\gh_a \gh_b \cgh_c}{\cgh_i} \right),\\
= & \;C^c_{ab}C^i_{jk} \left( \gh_a \gh_b\sbr{ \cgh_c}{ \gh_j \gh_k}\cgh_i+ \sbr{\gh_a \gh_b }{ \gh_j \gh_k}\cgh_c\cgh_i  +  \gh_j \gh_k \gh_a \gh_b \sbr{ \cgh_c}{\cgh_i}- \gh_j \gh_k \sbr{\gh_a \gh_b}{ \cgh_i}\cgh_c \right),\\
= & \;C^c_{ab}C^i_{jk} \left( \gh_a \gh_b\sbr{ \cgh_c}{ \gh_j \gh_k}\cgh_i- \gh_j \gh_k \sbr{\gh_a \gh_b}{ \cgh_i}\cgh_c \right),\\
= & \;C^c_{ab}C^i_{jk} \left( \gh_a \gh_C(\delta^c_j \gh_k - \delta^c_k \gh_j)\cgh_i- \gh_j \gh_k (\delta^a_i \gh_b - \delta^b_i \gh_a)\cgh_c \right),\\
= & \;\left(C^c_{ab}C^i_{ck}\gh_a \gh_b\gh_k \cgh_i- (1/4)C^i_{jk} C^c_{ib} \gh_j \gh_k  \gh_b \cgh_c\right)\\
& +\left( -C^c_{ab}C^i_{jc}  \gh_a \gh_b \gh_j\cgh_i +C^i_{jk}C^c_{ai}\gh_j \gh_k  \gh_a\cgh_c \right),\\
= & \;0.
\end{align*}
The last equality comes from equating the coefficients of each term in the brackets.
\end{proof}

\begin{proof}(Proof of \eqref{eq:nonabdel}) We will drop the $\otimes$ notation for this proof. First we let,
\begin{equation}\label{eq:nonaD}
\Delta=\{Q,Q^*\}=A+B+C+D,
\end{equation}
where,
\begin{gather*}
A=\{G_a  \gh_a,G_b  \cgh_b\}, \qquad B=\{G_j  \cgh_j,(-i/2) C^c_{ab}\gh_a \gh_b \cgh_c \},\\
C=\{G_j  \gh_j,(i/2) C^c_{ab}\gh_c\cgh_b \cgh_a \}, \qquad D=(1/4)\{ C^c_{ab}\gh_a \gh_b \cgh_c,  C^i_{jk}\gh_i\cgh_k \cgh_j \}.
\end{gather*} 
Now recall that,
\[
\Sigma_a := -iC^c_{ab} \gh_b \gh_c. 
\] 
Using the ghost commutation relations \eqref{eq:hghstrel}, the constraint brackets \eqref{eq:conbrak} and antisymmetry of $C^c_{ab}$ in every index, we calculate $A,B,C$ and $D$.
\begin{align*}
A&=[G_a ,G_b]\gh_a\cgh_b+G_bG_a\{ \gh_a,\cgh_b\},\\
&=iC^c_{ab}G_c \gh_a \cgh_b +G_bG_b,\\
&=G_a\Sigma_a+G_bG_b,,
\end{align*}
and,
\begin{align*}
B&=-(i/2) C^c_{ab}G_j\sbr{ \cgh_j}{\gh_a \gh_b \cgh_c },\\
&=-(i/2) C^c_{ab}G_j (\sbr{\cgh_j}{\gh_a \gh_b }\cgh_c+\gh_a \gh_b\sbr{\cgh_j}{\cgh_c }),\\
&=-(i/2) C^c_{ab}G_j (\delta^a_j \gh_b \cgh_c- \delta^b_j \gh_a\cgh_c),\\
&=-(i/2) C^c_{ab}G_a\gh_b \cgh_c+(i/2) C^c_{ab}G_b \gh_a  \cgh_c,\\
&=G_a\Sigma_a.
\end{align*}
By a similar computation we calculate,
\[
C= -G_a \Sigma_a.
\]
Next we calculate,
\begin{align}\label{eq:DelHen}
4D= & \;C^c_{ab}C^i_{jk}\sbr{\gh_a \gh_b \cgh_c}{  \gh_i\cgh_k \cgh_j },\notag\\
= & \;C^c_{ab}C^i_{jk} \left( \sbr{\gh_a \gh_b \cgh_c}{\gh_i}\cgh_k \cgh_j - \gh_i \sbr{\gh_a \gh_b \cgh_c}{\cgh_k \cgh_j} \right),\notag\\
= & \;C^c_{ab}C^i_{jk} \left( \gh_a \gh_b \sbr{\cgh_c}{\gh_i}\cgh_k\cgh_j - \sbr{\gh_a \gh_b}{\gh_i} \cgh_c \cgh_k \cgh_j - \gh_i \gh_a \gh_b \sbr{\cgh_c}{ \cgh_k \cgh_j} - \gh_i \sbr{\gh_a \gh_b}{\cgh_k \cgh_j} \cgh_c \right),\notag\\
= & \;C^c_{ab}C^i_{jk} \left( \delta_c^i\gh_a \gh_b \cgh_k\cgh_j  - \gh_i \sbr{\gh_a \gh_b}{\cgh_k \cgh_j} \cgh_c \right),\notag\\
= & \;C^c_{ab}C^i_{jk} \left( \delta_c^i\gh_a \gh_b \cgh_k\cgh_j  - \gh_i \sbr{\gh_a \gh_b}{\cgh_k} \cgh_j \cgh_c - \gh_i \cgh_k \sbr{\gh_a \gh_b}{\cgh_j} \cgh_c \right),\notag\\
= & \;C^c_{ab}C^i_{jk} \left( \delta_c^i\gh_a \gh_b \cgh_j\cgh_k  -\delta_b^k \gh_i \gh_a \cgh_j \cgh_c + \delta_k^a \gh_i \gh_b \cgh_j \cgh_c - \delta_b^j \gh_i \cgh_k \gh_a \cgh_c + \delta_a^j \gh_i \cgh_k \gh_b \cgh_c \right),\notag\\
= & \;[C^c_{ab}C^c_{jk} \gh_a \gh_b \cgh_k\cgh_j  - C^c_{ab}C^i_{jb} \gh_i \gh_a \cgh_j \cgh_c + C^c_{ab}C^i_{ja} \gh_i \gh_b \cgh_j \cgh_c \notag\\
&- C^c_{ab}C^i_{bk} \gh_i \cgh_k \gh_a \cgh_c + C^c_{ab}C^i_{ak} \gh_i \cgh_k \gh_b \cgh_c].
\end{align}
Now,
\begin{align}
- C^c_{ab}C^i_{bk} \gh_i \cgh_k \gh_a \cgh_c + C^c_{ab}C^i_{ak} \gh_i \cgh_k \gh_b \cgh_c&=-2C^c_{ab}C^k_{ai} \gh_i \cgh_k \gh_b \cgh_c \notag,\\
& =2 \Sigma_a \Sigma_a \label{eq:delhen2}
\end{align}
Also by the antisymmetry of the structure constants and the Jacobi identity,
\begin{align}
C^c_{ab}C^c_{jk} \gh_a \gh_b \cgh_k\cgh_j&=-C^b_{ac}C^c_{jk} \gh_a \gh_b \cgh_k\cgh_j,\notag\\
&=(C^b_{kc}C^c_{aj}+C^b_{jc}C^c_{ka})\gh_a \gh_b \cgh_k\cgh_j,\notag \\
&=C^b_{jc}C^c_{ak}\gh_a \gh_b \cgh_j\cgh_k+C^b_{jc}C^c_{ka}\gh_a \gh_b \cgh_k\cgh_j, \notag \\
&=-C^b_{jc}C^c_{ak}\gh_a \gh_b \cgh_k\cgh_j+C^b_{jc}C^c_{ka}\gh_a \gh_b \cgh_k\cgh_j,\notag \\
&=-2C^b_{jc}C^c_{ak}\gh_a \gh_b \cgh_k \cgh_j \label{eq:delhen3},
\end{align}
and,
\begin{align}
- C^c_{ab}C^i_{jb} \gh_i \gh_a \cgh_j \cgh_c + C^c_{ab}C^i_{ja} \gh_i \gh_b \cgh_j \cgh_c&= C^a_{cb}C^b_{ij} \gh_i \gh_a \cgh_j \cgh_c + C^a_{ij}C^b_{ca} \gh_i \gh_b \cgh_j \cgh_c,\notag \\
= & \;C^b_{jc}C^c_{ak} \gh_a \gh_b \cgh_k \cgh_j + C^c_{ak}C^b_{jc} \gh_a \gh_b \cgh_k \cgh_j,\notag \\
= & \;2C^b_{jc}C^c_{ak} \gh_a \gh_b \cgh_k \cgh_j. \label{eq:delhen4}
\end{align}
Substituting equations \eqref{eq:delhen2}, \eqref{eq:delhen3}, \eqref{eq:delhen4} in equation\eqref{eq:DelHen} gives that,
\[
D=(1/2) \Sigma_a \Sigma_a
\]

Substituting back in \eqref{eq:nonaD} we calculate,
\begin{align*}
\Delta= & \;(1/2)G_aG_a \otimes \one + (1/2)( G_a \otimes \one +\Sigma_a)(G_a \otimes \one +\Sigma_a),\\
\end{align*}
\end{proof}

\chapter{General BRST Structures}\label{ch:GenStruct} 
The aim of this chapter is to obtain a mathematically consistent definition of the heuristic BRST structures outlined above including the ghosts, the charge $Q$, the \emph{dsp}-decomposition, and the operator cohomology.
\section{Ghost Algebra}
\label{sec:ghost}
We first need to properly define the ghost algebra (cf. Section \ref{sec:heursetup}). Recall the key heuristic features:
\begin{enumerate}
\item There is a representation of an algebra $\cA_g:= \alg{\gh_j, \, \cgh_j \,| j \in I }$ acting on a Krein space $\cH_g$ with indefinite inner product $\iip{\cdot}{\cdot}_g$ such that the ghosts $\gh_j$ and conjugate ghosts $\cgh_j$ are Krein-hermitian.
\item There is a ghost number operator $G$ such that $G^{\dag}=-G$, $\cH_g=\oplus_{n=-\infty}^{\infty}\cH_{g,n}$ where $\cH_{g,n}=\{ \psi \in \cH_g| G\psi= n \psi \}$, and $[G,A]=nA$, $A \in \cG_n$.
\end{enumerate}

Conditions for the existence of a ghost number operator have been studied in \cite{AzKh89}, where the structure of Krein spaces $\cH_g$ that admit such a ghost number operators is determined. We will connect with these results after we have constructed the ghost algebra in a well-defined mathematical framework.

We start with a construction for even and infinite dimensional $\cA_g$. This construction will correspond to the smeared version of the ghosts in the BRST-QEM example of the previous section (example \eqref{sec:exheuEM}). We will then treat the odd degree case by restriction and finally, connect to superfunction spaces as given in Berezin \cite{Be1966} sections 3.1, 3.2, \HT \cite{HenTei92} Chapter 20.2,  \cite{McPat1989I} p482.

We define the Ghost Algebra as a CAR algebra over a test function with the following structures. Let $\cH$ be a Krein space with indefinite inner product $\iip{\cdot}{\cdot}$, a fundamental symmetry $J$ (see Appendix \eqref{ap:IIP}), \ie $\cH$ is a Hilbert space with norm coming from the inner product $\ip{\cdot}{\cdot}$,  $\cH=\cH_{+}\oplus \cH_{-}$ where $\cH_{+},\cH_{-}$ are closed subspaces of $\cH$ and $\oplus$ denotes Hilbert orthogonality, and $\iip{\cdot}{\cdot}=\ip{\cdot}{J\cdot}$ where $J=P_{+}-P_{-}$ where $P_{+},P_{-}$ are the projections on $\cH_{+}, \cH_{-}$ respectively. 
Note that $J$ is a unitary on $\cH$ such that $J^*=J$, $J^2=\one$.  We will denote $[\perp]$, $[\oplus]$ for Krein orthogonality, etc. and $\perp$, $\oplus$, for Hilbert orthogonality, etc.

\begin{lemma}\label{lm:ghtsfext}
Let $\cH_{+}$ and $\cH_{-}$ have the same cardinality. Then there exist closed subspaces $\cH_1,\cH_2 \subset \cH$ such that $\cH=\cH_1\oplus \cH_2$, $J\cH_1=\cH_2$.

Both $\cH_1$ and $\cH_2$ are neutral with respect to the Krein inner product $\iip{\cdot}{\cdot}$ and $P_1=JP_2J$ where $P_1,P_2$ are the projections on $\cH_1,\cH_2$ respectively.  
\end{lemma}
\begin{proof}
Let $\cH_{+}$ and $\cH_{-}$ have orthonormal basis $(e^1_{\lambda})_{\lambda\in \Lambda}$, $(e^2_{\lambda})_{\lambda\in \Lambda}$ respectively for some index set $\Lambda$. Now define,
\[
f^1_{\lambda}=\frac{e^1_{\lambda}+e^2_{\lambda}}{\sqrt{2}},\qquad f^2_{\lambda}=\frac{e^1_{\lambda}-e^2_{\lambda}}{\sqrt{2}},\qquad \lambda \in \Lambda
\]
We have that $\{f^1_{\lambda},f^2_{\lambda}\}_{\lambda \in \Lambda}$ is an orthonormal basis of $\cH$ and $Jf^1_{\lambda}=f^2_{\lambda}$. 

Let $\cH_1:=[f^1_{\lambda}\:|\, {\lambda \in \Lambda}]$ and $\cH_2:=[f^2_{\lambda}\:|\, \lambda \in \Lambda]$. Then $J\cH_1=\cH_2$. As $\iip{\cdot}{\cdot}=\ip{\cdot}{J\cdot}$ by definition, it follows that $\cH_1$ is neutral in the Krein inner product by $\cH_1\perp \cH_2$. Also for $\psi \in \cH_1$ we have $P_2J\psi=J\psi$ from which it follows that $JP_2J=P_1$.
\end{proof}

The Fermi-Fock space for $\cH$ has a natural Krein structure associated to $J$:
\begin{lemma}\label{lm:JKsp2}
Let $\cH$ and $J$ be as above, that $J \neq \pm \one$, and let $\fF^{-}(\cH)$ be the Fermi-Fock space with Hilbert inner product $\ip{\cdot}{\cdot}_g$. Let $J_g:=\Gamma(J)$ be the second quantisation of $J$  and define the indefinite inner product:
\[
\iip{\cdot}{\cdot}_g:=\ip{\cdot}{J_g\cdot}_g,
\]
on $\fF^{-}(\cH)$. Then $\fF^{-}(\cH)$ is a Krein space with indefinite inner product $\iip{\cdot}{\cdot}_g$, and fundamental symmetry $J_g$.
\end{lemma}
\begin{proof} We have that $J$ is $*$-unitary implies that $J_g$ is $*$-unitary and furthermore $J=J^*$ and $J \neq \pm \one$ implies that $J_g=J_g^*$ and $J_g \neq \pm \one$. Hence by lemma \eqref{lm:JKsp} we have that $\iip{\cdot}{\cdot}_g$ is indefinite and the $\fF^{-}(\cH)$ is a Krein space with indefinite inner product $\iip{\cdot}{\cdot}_g$ and fundamental symmetry $J_g$.
\end{proof}

Let $c(f),c^*(f)$ be the usual annihilators and creators on $\fF^{-}(\cH)$ and define:
\[
\cA(\cH):= \osalg{c(f)\,,c^*(g)\,|\, f,g \in \cH }\subset B(\fF^{-}(\cH)).
\]
We have:
\begin{equation}\label{eq:CARac}
\{c(f),c(g)\}=\{c^*(f),c^*(g)\}=0,\qquad \{c(f),c^*(g)\}=\ip{f}{g}\one,
\end{equation}
(cf. \cite{BraRob21981} section 5.2). As $\fF^{-}(\cH)$ is a Krein space, we can calculate Krein adjoints by $A^{\dag}=J_gA^*J_g$ (cf. lemma \eqref{lm:khad}). So:
\begin{equation}\label{eq:ccomm}
c^{\dag}(f)=\Gamma(J) c^*(f) \Gamma(J)=c^*(Jf),
\end{equation}
and hence:
\[
\{c(f),c^{\dag}(g)\}=\{c(f),c^{*}(Jg)\}=\ip{f}{Jg}\one=\iip{f}{g}\one.
\]
Define the Krein-hermitian fermionic quantum field operators:
\begin{equation}\label{eq:Cliffel}
C(f)=\fst(c(f)+c^{\dag}(f))=\fst(c(f)+c^{*}(Jf)),
\end{equation}
on $\fF^{-}(\cH)$. Note that $C(f)^*=C(Jf)$ and so $C(f)$ is Hilbert hermitian only when $Jf=f$. Then,
\begin{equation}\label{eq:Cliffcom}
\{C(f), C(g) \}= \re \iip{f}{g}\one,
\end{equation}
which corresponds to a Clifford algebra with indefinite inner product. Define the ghost fields as:
\begin{equation}\label{eq:gfield}
\gh(f):=C(f) \qquad  \forall f \in \cH_2,
\end{equation}
then $\forall f,g \in \cH_2$,
\begin{align}
\gh^*(f)&=C(Jf) \notag \\
\gh^{\dag}(f)&= \gh(f), \quad (\gh^*)^{\dag}(f)= \gh^*(f),  \notag \\
\{ \gh(f), \gh(g) \}& = Re \iip{f}{g}\one= 0 \quad \text{as $\cH_2$ is neutral} \notag \\
\{ \gh(f), \gh^*(g) \}& = \{ C(f), C(Jg) \}= Re \iip{f}{Jg}\one= Re \ip{f}{g}\one \label{eq:Cliffhilcom}.
\end{align}

Now take an orthonormal basis $(f_j)_{j\in \Lambda}$ of $\cH_2$, where $\Lambda$ is some index set. The above relations justify the identification of $\gh(f_j)$ as our ghosts and 
\begin{equation}\label{eq:congfield}
\cgh(f_j):=\gh^*(f_j)=C(Jf_j)\qquad  \forall f \in \cH_2,
\end{equation}
as its conjugate ghost. That is, they are both $\dag$-hermitian and, 
\begin{align}
\{\gh(f_j),\gh(f_k) \} & =\{\cgh(f_j), \cgh(f_k) \} =Re \iip{f}{f}\one = 0 \quad \text{(as $\cH_2$ is neutral)} \label{eq:ghcomrelrig}\\
\{\gh(f_j), \cgh(f_k) \} &= \{\gh(f_j), \gh^*(f_k) \}= Re \ip{f_j}{f_k}\one=\delta_{kj}\one \notag
\end{align}
which are the relations \eqref{eq:hghstrel}. Moreover $J_g\gh(f)J_g=\cgh(f)$. Note that we have associated the Clifford algebra elements $C(f)$ with $\cH_2$ arguments as ghosts and with $\cH_1$ arguments as conjugate ghosts. Using $c(if)=-ic(f)$, $c^*(if)=ic^*(f)$ and $[J,i]=0$ for all $f\in \cH$, we can recover the creators and annihilators from the ghosts by:
\begin{align}
c(f)&=\fst(\gh(f)+i\gh(if)), \qquad f \in \cH_2, \notag \\ 
c^*(f)&=\fst(\cgh(f)-i\cgh(if)), \qquad  f \in \cH_2, \label{eq:ghs}\\ 
\intertext{and,}
c^*(f)&=\fst(\gh(Jf)-i\gh(iJf)), \qquad f \in \cH_1,\label{eq:cghs} \\ 
c(f)&=\fst(\cgh(Jf)+i\cgh(iJf)), \qquad  f \in \cH_1. \notag 
\end{align}

We define our ghost algebra and Fock ghost space as:
\begin{definition}\label{df:GA}
The Fock ghost space is,
\[
\cH_g:= \fF^{-}(\cH)=\overline{\cA_g(\cH_2) \Omega}.
\]
where the ghost algebra is,
\begin{align*}
\cA_g(\cH_2)&:=\osalg{ \gh(f), \gh^*(g)\, |\, f,g \in \cH_2 } \\
&=\osalg{ C(f)|\, f \in \cH=\cH_1 \oplus \cH_2 } \\
&=\cA(\cH)\subset B(\fF^{-}(\cH)).
\end{align*}
Denote $\cA_{g0}(\cH_2):=\alg{ \gh(f), \gh^*(g)\, |\, f,g \in \cH_2 }$ and note that $\cA_{g0}(\cH_2)$ is dense in $\cA_g(\cH_2)$ in the norm topology. Also $\cA_g(\cH_2)$ and $\cA_{g0}(\cH_2)$ will be denoted by $\cA_g$ and $\cA_{g0}$ respectively when no confusion will arise.

\end{definition}

\begin{rem}
\begin{itemize}
\item [(i)] There is an automorphism $\alpha \in \Aut(\cA_g(\cH_2))$ given by 
\begin{equation}\label{eq:autghalg}
\alpha(C(f)):=C(Jf).
\end{equation}
That this is an automorphism follows from $\{C(Jf),C(Jg)\}=\ip{Jf}{Jg}\one=\ip{f}{g}\one$ and \cite{BraRob21981} Theorem 5.2.5. We will denote this automorphism by $\alpha$ in the sequel unless specified otherwise. Note that we have $A^{\dag}=J_gA^*J_g= \alpha(A^*)$ for $A \in \cA_g$ by which we define the Krein involution on the algebra independently of the representation $\cH_g$.
\item[(ii)] Given an orthonormal basis $(f_j)_{j\in \Lambda}$ of $\cH_2$, the linearly independent set:
\begin{align}\label{eq:setS}
S= & \;\{ \gh(f_j), \gh^*(f_j), \gh(if_j) ,\gh^*(if_j) \,| \, j\in \Lambda \} ,\\
= & \;\{ \gh(f_j), \cgh(f_j), \gh(if_j) ,\cgh(if_j) \,| \, j \in \Lambda \} ,
\end{align}
generates, via equations \eqref{eq:ghs} and \eqref{eq:cghs}, a dense subalgebra of $\cA_g$, whereas the set,
\[
\{ \gh(f_j), \gh^*(f_j) \,| \, j\in \Lambda \}=\{ \gh(f_j), \cgh(f_j)\,| \, j\in \Lambda \}
\]
does not. This is easy to see by dimensional comparison in the finite dimensional case, which extends to the infinite dimensional case as for each $j \in \Lambda$ we have that the set $\alg{ \gh(f_j), \cgh(f_j) }$ generates a proper norm closed subspace of \\
$\cA_g=\osalg{ \gh(f_j), \gh^*(f_j), \gh(if_j) ,\gh^*(if_j) \,| \, j\in \Lambda }$.
\end{itemize}
\end{rem}

\subsection{Gradings}\label{sbs:ghgrad}
Next, we define useful gradings for $\cA_g$. Recall that $P_i$ is the projection onto $\cH_i$ ($i=1,2$), and define,
\[
N_1:=d\Gamma(P_1), \qquad N_2:=d\Gamma(P_2), \qquad \text{where} \qquad D(N_1):=D(N_2)=\fF_0^{-}(\cH).
\]
where $D(N_1)$ and $D(N_2)$ denote the domains of $N_1$ and $N_2$ respectively and $\fF_0^{-}(\cH)$ is the finite particle subspace of the Fock space $\cH_g=\fF^{-}(\cH)$. We have:
\[
J_gN_1\psi=\Gamma(J)d\Gamma(P_1)\psi=d\Gamma(JP_1)\psi=d\Gamma(P_2J)\psi=N_2J_g\psi
\]
for all $\psi \in \fF_0^{-}(\cH)$, making use of $JP_1=P_2J$. Thus by $J_g\fF_0^{-}(\cH)=\fF_0^{-}(\cH)$:
\begin{equation}\label{eq:N1N2}
J_gN_1J_g = N_2.
\end{equation}
Now $N_1$ and $N_2$ on $\fF^{-}_0(\cH_1)$ and $\fF^{-}_0(\cH_2)$ respectively are just the number operators which we combine:
\[
N:=N_1+N_2=d\Gamma(P_1+P_2)=d\Gamma(\one),\qquad D(N)=\fF_0^{-}(\cH),
\]
\ie $N$ is the number operator on $\fF_0^{-}(\cH)$. Note that $N$ is $*$-symmetric on $\fF_0^{-}(\cH)$, so by $N^\dag=J_gN^*J_g$ (see Appendix \eqref{ap:IIP}) and equation \eqref{eq:N1N2}, we have that $N$ is $\dag$-symmetric on $\fF_0^{-}(\cH)$. Now $N$ induces the usual grading on $\cH_g$:
\begin{proposition}\label{pr:CARgr}
\[
\cA_g=\sqb{\cup_{n=-2m}^{2m}\cA_n}, \qquad \cH_g={\oplus_{n=0}^{2m}\fF_n}
\]
where $m=\dim(\cH_1)$ (possibly infinite), and,
\begin{align*}
\cA_n:&=\{ A\in \cA_g \,|\, AD(N)\subset D(N), \; [N, A]\psi=n A\psi, \,\psi \in D(N) \},\\
\fF_n:&=\{ \psi \in \cH_g \,|\, N\psi= n \psi \}.
\end{align*}
Moreover $\fF_n= \overline{\cA_n \Omega}$ for all $n\in \mathbb{N}$.
\end{proposition}
\begin{proof}
The standard definition of the number operator is to define $N\psi_n:=n\psi_n$ for $\psi_n$ in the $n$-particle subspace of $\fF^{+}(\cH)$, and extend linearly to the domain \\
$\{ \psi \in \fF^{+}(\cH)\,|\, \sum_{i=1}^{\infty} |n \psi_n|^2 < \infty\}$ \cite{BraRob21981} p7. It is then proved that $N=d\Gamma(\one)$ on $\fF_0^{+}(\cH)$, and so the above decomposition of $\cH_g=\oplus_{n=0}^{2m}\fF_n$ follows.

Let $(f_{\lambda})_{\lambda\in \Lambda}$ be an orthonormal basis for $\cH$, and let,
\[
M^n_{\bla,\bmu}=c^*(f_{\lambda_1})\ldots c^*(f_{\lambda_k})c(f_{\mu_1})\ldots c(f_{\mu_l})\neq 0
\]
be a monomial of creators and annihilators where $\bla=(\lambda_1,\dots,\lambda_k)\in \Lambda^k$, $\bmu =(\mu_1,\dots,\mu_l)\in \Lambda^l$ and $(k-l)=n\in \mathbb{Z}$. As $M^n_{\bla,\bmu} \neq 0$ we have that $\lambda_1\ldots \lambda_k$ are distinct and $\mu_1\ldots \mu_l$ are distinct. Also as $c(f)^2=c^*(f)^2=0$ for all $f \in \cH$ we have that, 
\begin{equation}\label{eq:overkill}
M^n_{\bla,\bmu}\neq 0 \Rightarrow n \leq 2m.
\end{equation} 
Let $j \in \mathbb{N}$ and take $\psi \in \fF_j$. Now $\fF_n$ is the subspace of $\fF_0^{-}(\cH)$ spanned by the $n$-particle vectors, hence by the definition of the creators and annihilators:
\[
c^*(f)\fF_n\subset \fF_{n+1}, \qquad c(f)\fF_n\subset \fF_{n-1},  \qquad c(f)\fF_0=\{0\}, \qquad n\geq 1, \, \forall f \in \cH
\]
From this we have that $M^n_{\bla,\bmu} \psi \in \fF_{j+n}$ for $(j+n)\geq 0$ and $M^n_{\bla,\bmu} \psi=0$ for $(j+n)<0$. Therefore, 
\begin{equation}\label{eq:comblah}
[N,M^n_{\bla,\bmu}]\psi=(j+n)M^n_{\bla,\bmu}\psi-jM^n_{\bla,\bmu}\psi=nM^n_{\bla,\bmu}\psi,
\end{equation}
where equation \eqref{eq:comblah} is trivially true when $M^n_{\bla,\bmu}\psi=0$. As $j$ was arbitrary and $D(N)=\fF_0^{-}(\cH)=\oplus_{n=0}^{2m}\fF_n$ we get that \eqref{eq:comblah} is true for all $\psi \in D(N)$. Hence $M^n_{\bla,\bmu}\in \cA_{n}$. Now by \cite{BraRob21981} Theorem 5.2.5 p15,
\[
\cA_g =[\{\one, M^n_{\bla,\bmu} \,|\, n\in \mathbb{Z},\, \bla\in \Lambda^k, \,\bmu \in \Lambda^l\, \text{and}\, (k-l)=n\} ],
\]
and we obviously have $\one \in \cA_{0}$, hence,
\begin{align*}
&\left[ \left\{\one \cup \left( \bigcup \{ M^n_{\bla,\bmu} \,|\,   n\in \mathbb{Z},\, \bla\in \Lambda^k, \,\bmu \in \Lambda^l\, \text{and}\, (k-l)=n \} \right)\right\}\right]\subset {\sqb{\cup_{n=-2m}^{2m}\cA_n}}  \subset \cA_g,\\
&= { \left[\left\{ \one \cup \left(\bigcup \{M^n_{\bla,\bmu} \,|\,  n\in \mathbb{Z},\, \bla\in \Lambda^k, \,\bmu \in \Lambda^l\, \text{and} \,(k-l)=n\}\right)\right\}\right]},
\end{align*}
where we used equation \eqref{eq:overkill} for the first inclusion. Therefore ${\sqb{\cup_{n=-2m}^{2m}\cA_n}} = \cA_g$. 

By the definition of $N$, $\fH_n$ is the $n$-th component of a vector in Fock space and so for $n\geq 0$,
\begin{align*}
\fF_n=\sqb{M^n_{\bla,\bmu}\Omega_g \,|\,  n\in \mathbb{Z}\backslash \{0\},\, \bla\in \Lambda^k, \,\bmu \in \Lambda^l\, \text{and} \,(k-l)=n} \subset \overline{\cA_n \Omega}\subset \fF_n
\end{align*}
where the last inclusion follows as $N\Omega=0$ implies $\cA_n\Omega \subset \fF_n$ and that $\fF_n$ is closed. That is we have $\fF_n=\overline{\cA_n \Omega}$ and we are done.
\end{proof}
We will refer to these as the \textit{{CAR gradings}}. 

To the Ghost Gradings below we introduce the notation:

We define the ghost number operator as,
\[
G:=N_1-N_2=d\Gamma(P_1-P_2),\qquad D(G)=\fF_0^{-}(\cH),
\]
and note that $G$ is $*$-symmetric on $\fF_0^{-}(\cH)$. By $J_gN_1J_g \psi = N_2 \psi$ for $\psi \in \fF_0^{-}(\cH)$, we have that $G$ is $\dag$-skew symmetric on $ \fF_0^{-}(\cH)$. This gives a grading of $\cA_g$ and $\cH_g$ as,
\begin{proposition}\label{pr:ghgrad}
\[
\cA_g:={[\cup_{n=-2m}^{2m}\cG_n]}, \qquad \cH_g:={\oplus_{n=-m}^{m}\fH_n}
\]
where $m=\dim(\cH_1)$ (possibly infinite), and
\begin{align*}
\cG_n:&=\{ A\in \cA_g \,|\,AD(G)\subset D(G),\, [G, A]\psi=n A\psi, \,\psi \in D(G) \},\\
\fH_n:&=\{ \psi \in \cH_g \,|\, G\psi= n \psi \}=\overline{\cG_n \Omega},
\end{align*}
\end{proposition}

\begin{proof}
Suppose that $\Omega_1$ and $\Omega_2$ are the vacuum vectors of  $\fF_0^{-}(\cH_1)$ and $\fF_0^{-}(\cH_2)$ respectively, then  $\fF^{-}_0(\cH_1\oplus \cH_2)$ is unitarily equivalent to $\fF_0^{-}(\cH_1)\otimes\fF_0^{-}(\cH_2)$ by defining the map,
\begin{align*}
U:\fF_0^{-}(\cH_1)\otimes\fF_0^{-}(\cH_2) &\to \fF^{-}_0(\cH_1\oplus \cH_2),\\
U(\Omega_1\otimes \Omega_2):=\Omega, \hspace{.9cm}&\\
\begin{split}
U(a^*(f_{\lambda_1})\ldots a^*(f_{\lambda_n})\Omega_1&\otimes a^*(f_{\mu_1})\ldots a^*(f_{\mu_l})\Omega_2):=\\
&=\left(\sqrt{\frac{n!l!}{(n+l)!}}\right)a^*(f_{\lambda_1})\ldots a^*(f_{\lambda_n})a^*(f_{\mu_1})\ldots a^*(f_{\mu_l})\Omega,
\end{split}\\
\end{align*}
where $f_{\lambda_1},\ldots,f_{\lambda_n}\in \cH_1$, $f_{\mu_1},\ldots,f_{\mu_l}\in \cH_2$, $n,l \in \mathbb{Z}$. We then extend linearly to all of $\fF_0^{-}(\cH_1)\otimes\fF_0^{-}(\cH_2)$. It is easily checked that $U$ is well defined and unitary, the only thing to note is that the $({n!l!})/({(n+l)!})$ comes from the antisymmetrization of components of elements in $\fF_0^{-}(\fH)$. Now we identify $\fF_0^{-}(\cH_1)$ with the subspace $U(\fF_0^{-}(\cH_1)\otimes \Omega_2)$ of $\fF^{-}_0(\cH_1\oplus \cH_2)$ and  $\fF_0^{-}(\cH_1)$ with the subpsace $U(\fF_0^{-}(\cH_1)\otimes \Omega_2)$ of $\fF^{-}_0(\cH_1\oplus \cH_2)$.

The result now follows similar to Proposition \ref{pr:CARgr} if we recall that $N_1$ and $N_2$ on $\fF^{-}_0(\cH_1\oplus \cH_2)=U(\fF_0^{-}(\cH_1)\otimes\fF_0^{-}(\cH_2))$ are just the number operators when restricted to $\fF_0^{-}(\cH_1)$ and $\fF_0^{-}(\cH_2)$ and that $G=N_1-N_2$. 
\end{proof}
We refer to these as the {\textit{Ghost gradings}}. Note that $G\Omega =0$. 

The CAR and Ghost gradings of the spaces and Algebras give the following structure
\begin{lemma}\label{lm:gradprops}
Let $\cA_g(\cH_2)$ be the ghost algebra with $\dim{\cH_2}=m$ (possibly infinite). Let $k,l \in \{-2m,\ldots,2m\}$, then:
\begin{itemize}
\item[(i)] $\cA_k \cA_l \subset \cA_{k+l}, \qquad \cG_k \cG_l \subset \cG_{k+l}$.
\item[(ii)] $\cA_k^*=\cA_k^{\dag}=\cA_{-k}$ and $\cG_k^*=\cG_{-k}$ and $\cG_{k}^{\dag}=\cG_{k}$. 
\item[(iii)] $\gh(f) \in \cG_1$ and $\cgh(f) \in \cG_{-1}$ for all $f \in \cH$.
\item[(iv)] $\fF_{k} \perp \fF_{l}$  for $k,l \geq 0$. $\fH_{k}\perp \fH_{l}$ for $k \neq l$ where $\perp$ denotes Hilbert orthogonality and $k,l \in \mathbb{Z}$. Moreover $\fH_{k} [\perp] \fH_{-l}$ for $k\neq l$  where $[\perp]$ denotes Krein orthogonality with respect to $\iip{\cdot}{\cdot}_g$,
\end{itemize}
where we take the convention that $\cA_{k}=\cG_{k}=\{0\}$ for $k < -2m$, or $k>2m$.
\end{lemma}
\begin{proof}

\noindent(i): Take $A \in \cA_k$, $B \in \cA_l$ and $\psi \in \fF_0^{-}(\cH)$. Then using $\cA_n D(G) \subset D(G)$ for all $n \in \mathbb{Z}$, we get $[N,AB]\psi=[N,A]B\psi+A[N,B]\psi=kAB\psi+lAB\psi=(k+l)AB\psi$. Therefore $AB \in \cA_{k+l}$ hence $\cA_j \cA_k \subset \cA_{k+l}$. Similarly $\cG_j \cG_k \subset \cG_{j+k}$.

\smallskip \noindent (ii): Let $A \in \cA_{k}$ and $\psi, \xi \in \fF_0^{-}(\cH)$. Then as $N$ is $*$-symmetric, \[
\ip{kA^*\psi}{\xi}_g=\ip{\psi}{kA\xi}_g=\ip{\psi}{[N,A]\xi}_g=\ip{-[N,A^*]\psi}{\xi}_g.
\]
Thus $[N,A^*]\psi=-kA^*\psi$ for all $\psi \in \fF_0^{-}(\cH)$ and hence $A^*\in \cA_{-k}$. So $\cA_k^*=\cA_{-k}$ and similar arguments produce the rest of (ii) using the fact that $N$ is $\dag$-symmetric, $G$ is $*$-symmetric, and $G$ is $\dag$-antisymmetric.

\smallskip \noindent (iii): We have,
\begin{align*}
\gh(f)&=\fst(c(f)+c^{*}(Jf)), \qquad f \in \cH_2\\
\cgh(f)=\gh^*(f)&=\fst(c(Jf)+c^{*}(f)), \qquad f \in \cH_2
\end{align*}
from which we see that, for all $f \in \cH_2$:
\begin{align}
[G,\gh(f)]\psi= & \;[(N_1-N_2),\fst(c(f)+c^{*}(Jf))]\psi=\gh(f)\psi, \label{eq:ghcghgradcom}\\ 
[G,\cgh(f)]\psi= & \;[(N_1-N_2),\fst(c(Jf)+c^{*}(f))]\psi=-\cgh(f)\psi,\notag
\end{align}
for all $\psi \in \fF_0^{-}(\cH)$, by $[N_1,c(f)]\psi=0$, $[N_2,c(f)]\psi=-c(f)\psi$, etc.

\smallskip \noindent (iv): Let $\psi \in \fF_{k}$, $\xi \in \fF_{l}$. Then as $N$ is $*$-symmetric,  $l\ip{\psi}{\xi}_{g}=\ip{\psi}{N\xi}_{g}=\ip{N\psi}{\xi}_g=k\ip{\psi}{\xi}_{g}$ and so $\psi \perp \xi$ if $k \neq l$. As $G$ is $*$-symmetric, the analogous argument shows that $\fH_{k}\perp \fH_{l}$ for $k \neq l$, and as $G$ is $\dag$-antisymmetric we also get $\fF_{k} [\perp] \fF_{-l}$ for $k\neq l$.
\end{proof}

\begin{rem}\label{rm:ghnad}
\begin{itemize}
\item[(i)]  Let $(f_j)_{j\in \Lambda}$ be an orthonormal basis of $\cH_2$, and let $A$ be a monomial of the $\gh(f_j)$'s and $\cgh^*(f_j)$'s. Then  $A \in \cG_{k}$ where $k$ is the difference between the number of $\cgh(f)$'s and the number of $\gh(f)$'s in $A$. By lemma \eqref{lm:gradprops} (i) and (iii) we get that $A \in \cG_{k}$. Thus the definition of the ghost grading above agrees with the heuristic definition given in Section \ref{sec:heursetup}.
\item[(ii)] The identity $\gh(f)=\fst(c(f)+c^{*}(Jf))$ for $f\in \cH_2$ shows that the CAR grading and ghost grading are different as $\gh(f)$ has ghost number one but no definite CAR grading number.
\end{itemize}
\end{rem}
Also we have,
\begin{definition}\label{df:Z2grad}
The Ghost grading induces the $\Z_2$-grading: 
\[
\cA_g:=\cA_g^{+}\oplus \cA_g^{-},
\]
where $m=\dim(\cH_1)$ (possibly infinite) and
\[
\cA_g^{+}:=[\cup_{n=-2m}^{2m}\cG_{2n}], \qquad \cA_g^{-}:=[\cup_{n=-2m}^{2m}\cG_{2n+1}].
\]
We refer to this as the \emph{$\Z_2$-grading} in $\cA_g$. We define
\[
\gamma(A)=A, \qquad \forall A \in \cA_g^{+}, \qquad \gamma(A)=-A, \qquad \forall A \in \cA_g^{-}.
\]
It is easy to see that $\gamma$ extends to a $*$-automorphism on $\cA$ such that $\gamma^2=\iota$. We refer to $\gamma$ as the \emph{$\Z_2$-grading automorphism} on $\cA_g$.
\end{definition}
\begin{rem}\label{rm:Zwgr} It is easy to see that $c(f), c^*(f) \in \cG_{1}\subset \cA_g^{-}$ for $f \in \cHL$ and $c(g), c^*(g) \in \cG_{-1}\subset \cA_g^{-}$ for $g \in \cHJ$. Therefore $c(f), c^*(f), C(f) \in \cA_g^{-}$ for all $f \in \cH$.
\end{rem}
These algebras and gradings are defined in the Fock representation. We want to find conditions on states so that the ghost grading structures can also be defined in the state space of their GNS representations. We follow \cite{Hendrik1991} lemma 6.1, p26,
\begin{proposition}\label{pr:ghsp} Let $\cA_g(\cH_2)$ be the ghost algebra with $dim(\cH_2)=m$ (possibly infinite), and suppose that $\omega\circ \alpha = \omega \in \fS(\cA_g)$ where $\alpha$ is defined in equation \eqref{eq:autghalg}. We define:
\[
\cH_n:= \overline{\pi_{\omega}(\cG_n) \Omega_{\omega}}\subset \cH_{\omega}, \qquad n \in \{-m,\ldots,m\}
\]
where $(\pi_{\omega},\Omega_{\omega},\cH_{\omega})$ is the GNS-data for $\omega$, with Hilbert inner product $\ip{\cdot}{\cdot}_{\omega}$.
Define $\iip{\cdot}{\cdot}_{\omega}=\ip{\cdot}{J_{\omega}\cdot}_{\omega}$ where $J_{\omega}$ is the implementer for $\alpha$ in $\cH_{\omega}$. Then we have
\[
\omega(\cG_n)=0 \; \forall n \neq 0 \quad \text{iff} \quad  \cH_n [\perp] \cH_k  \; \text{when} \; n \neq -k.
\]           
In this case we also have that $\cH_n$ is neutral with respect to $\iip{\cdot}{\cdot}_{\omega}$ (\ie $\iip{\psi}{\psi}_{\omega}=0$ for all $\psi \in \cH_n$) , $J_{\omega} \cH_n = \cH_{-n}$, $\cH_{\omega}= \oplus_{n=-m}^{m} {\cH_n}$, $\cH_n \perp \cH_k$ for $k\neq n$, and we have the following decomposition,
\[
\cH_{\omega}=\overline{ \cH_0 [ \oplus]^{m}_{n=1}(\cH_n \oplus \cH_{-n})}.
\]
Furthermore there exists a $*$-symmetric, and $\dag$-symmetric operator $G$ with domain\\
$D(G)=\cH_0 [ \oplus]^{m}_{n=1}(\cH_n \oplus \cH_{-n})$ such that $G\psi= n\psi$ for $\psi \in \cH_n$.
\end{proposition}
\begin{proof}
First for $\psi= \pi_{\omega}(A) \Omega_{\omega}$, $\xi=\pi_{\omega}(B) \Omega_{\omega}$ we have that, 
\[
\iip{\psi}{\xi}_{\omega}=\ip{\pi_{\omega}(A) \Omega_{\omega}}{J_{\omega}\pi_{\omega}(B) \Omega_{\omega}}_{\omega}=\omega(A^*\alpha(B))=\omega(\alpha(A^*)B)
\]
Let $\omega(\cG_n)=0$ for $n\neq 0$, and let $A \in G_{k}$, $B\in G_{n}$ and $k\neq -n$. By lemma \eqref{lm:gradprops} (ii)  we have that $\alpha(A^*)=A^{\dag} \in \cG_{k}$ and so by lemma \eqref{lm:gradprops} (i) $ \alpha(A^*)B \in G_{k+n} \neq G_{0}$. Therefore $\iip{\psi}{\xi}_{\omega}=\omega(\alpha(A^*)B)=0$. As such $\psi, \xi$ are dense in $\cH_k$ and $\cH_{n}$ respectively, we get that $\cH_{k} [\perp] \cH_{n}$.

Conversely suppose that $\cH_{n}[\perp] \cH_{k}$ for $n \neq -k$. Now $\Omega_{\omega} \in \cH_{0}$ and so if $n \neq 0$, then  $\cH_{n}[\perp] \Omega_{\omega}$. So if $A \in \cG_n$ then $\psi:=\pi_{\omega}(A) \Omega_{\omega} \in \cH_{n}$ and we get that $\omega(A)=\ip{\Omega_{\omega}}{\pi_{\omega}(A)\Omega_{\omega}}_{\omega}=\iip{\Omega_{\omega}}{\pi_{\omega}(A)\Omega_{\omega}}_{\omega}=\iip{\Omega_{\omega}}{\psi}_{\omega}=0$, where we used $J_{\omega}^*=J_{\omega}$ and $J_{\omega}\Omega_{\omega}=\Omega_{\omega}$.
 
Observe that $\cH_{n}[\perp] \cH_{k}$ for $k \neq -n$ implies that $\cH_{n}[\perp]\cH_{n}$ for all $n \neq 0$ and so $\cH_n$ is a neutral subspace for all $n\neq 0$. Next since $A^{\dag}=\alpha(A^*)$ we have:
\[
J_{\omega}\cH_n=J_{\omega}\overline{\pi_{\omega}(\cG_n)\Omega_{\omega}}=\overline{\pi_{\omega}(\alpha(\cG_n^*)^*)\Omega_{\omega}}=\overline{\pi_{\omega}((\cG_n^{\dag})^*)\Omega_{\omega}}=\overline{\pi_{\omega}(\cG_{-n})\Omega_{\omega}}=\cH_{-n}
\]
where we used $J_{\omega}\Omega_{\omega}=\Omega_{\omega}$ in the first equality and lemma \eqref{lm:gradprops} (ii) in the second last equality. Next,
\[
\cH_{\omega}=\overline{\pi_{\omega}(\cA_g)\Omega_{\omega}}=\overline{\pi_{\omega}(\sqb{\cup_{n=-m}^{m}\cG_n})\Omega_{\omega}}=\bigoplus_{n=-m}^{m} {\cH_n},
\]
where we used Proposition \ref{pr:ghgrad} in the third equality.

Now we have that $\cH_n\perp \cH_k$ for $k \neq n$ by $\cH_n [\perp]\cH_k$ for $k \neq -n$, $\ip{\cdot}{\cdot}_{\omega}=\iip{\cdot}{J_{\omega}\cdot}_{\omega}$ and $J_{\omega}\cH_n=\cH_{-n}$.

The decomposition, 
\[
\cH_{\omega}=\overline{ \cH_0 [ \oplus]_{n=1}^{m}(\cH_n \oplus \cH_{-n})}.
\]
now follows from the above relations.

Finally, $G\psi= n\psi$ for $\psi \in \cH_n$ defines a linear operator $G$ on $\cH_n$. As $\cH_n \perp \cH_k \Rightarrow \cH_n \cap \cH_k =\{0\}$ for $k \neq n$ this extends to a linear operator on $D(G)=\cH_0 [ \oplus]_{n=1}^{m}(\cH_n \oplus \cH_{-n})$. That $G$ is $*$-symmetric and $\dag$-antisymmetric follows by direct computation in the inner products $\ip{\cdot}{\cdot}_{\omega}$ and $\iip{\cdot}{\cdot}_{\omega}$ (This proof of the existence of $G$ is based on that in Proposition 1 \cite{AzKh89} p673-674).
\end{proof}

Hence we define,
\begin{definition}\label{df:ghst}
Let $\fS_{g}\subset \fS(\cA_g)$ be the set of states, 
\[
\fS_{g}=\{ \omega \in \fS(\cA_g)\,|\, \omega\circ \alpha = \omega, \, \omega(\cG_n)=0\; \text{for $n\neq 0$}\}.
\]
\end{definition}
That is, $\fS_{g}$ is the set of states with the correct ghost grading structure, and with GNS cyclic vector with ghost number $0$. 

We connect the above abstract structures to the Fock structures at the beginning of this chapter as follows. Let $\pi_F:\cA_g \to B(\cH_g)$ be the Fock representation of $\cA_g$, and the Fock state $\omega_F\in \fS(\cA_g)$ be defined by  $\omega_F(A):=\ip{\Omega}{\pi_F(A)\Omega}$. Then as $J_g\Omega=\Omega$ and $J_g\pi_F(A)J_g=\pi_F(\alpha(A))$ we have that $\omega_F \circ \alpha \in \fS_g$, hence $\omega_F \in \fS_g$. Moreover the GNS-representation $\pi_{\omega_F}$ is unitarily equivalent to $\pi_F$ by a unitary $U:\cH_g \to \cH_{\omega_F}$ with $U\Omega=\Omega_{\omega_F}$. Hence we recover all the concrete Fock structures from $\cA_g$ and $\omega_F$. The explicit spatial connection of the Fock structures to Proposition \eqref{pr:ghsp} is, 
\begin{gather*}
\cH_{\omega_F}=U\cH_g, \quad
\cH_n=\overline{\pi_{\omega_F}(\cG_n)\Omega_{\omega_F}}=U\fH_n=U\{ \psi \in \cH_g \,|\, G\psi= n \psi \}.
\end{gather*}
The ghost number operators are connected by $U^*GU=G$, where we have used the symbol $G$ for the ghost number operator on $\cH_g$ and for the ghost number operator defined in Proposition \eqref{pr:ghsp}.

\subsection{Finite Dimensional Ghost algebras}\label{sbs:fndmgh}
Here we consider the case of finite dimensional ghost algebras. For finite dimensional $\cH_2$ we have that all representations $\cA_g(\cH_2)$ are multiples of the Fock representation by \cite{BraRob21981} Theorem 5.2.14. Hence we will only consider the Fock-ghost representation of $\cA_g(\cH_2)$ on the Fock space $\cH_g$ in this section, and not denote the representation explicitly.

Suppose that $\cH_2$ has a complex orthonormal basis $(f_i)_{i=1}^{m}$ where $m < \infty$. So $\dim \cH=2m$ and $\cH$ has orthonormal basis $\{ Jf_i,  f_i|,\, i=1\dots \MM \}$. Then $\fF^{-}(\cH)=\fF_0^{-}(\cH)$ is finite dimensional,
\begin{align}\label{eq:NGdef}
N&=N_1+N_2=\sum_{j=1}^{\MM} (c^*(f_j)c(f_j) + c^*(Jf_j) c(Jf_j)) 
\intertext{and}
G&= N_1-N_2 \notag\\
&=\sum_{j=1}^{\MM} [  c^*(Jf_j) c(Jf_j)-c^*(f_j) c(f_j) ]   \notag\\
&=\frac{1}{2}\sum_{j=1}^{\MM} [\gh(f_j)\cgh(f_j)-\cgh(f_j)\gh(f_j)  +\gh(if_j)\cgh(if_j)-\cgh(if_j)\gh(if_j) ] 
\notag \\
\end{align}
where the last equality for $G$ follows from equations \eqref{eq:ghs} and \eqref{eq:cghs} and the brackets relations in equation \eqref{eq:ghcomrelrig}.

\begin{rem}
\begin{itemize}
\item[(i)] The formula for $G$ in terms of the ghosts corresponds to the heuristic ghost number operator in \HT \cite{HenTei92} p313.
\item[(ii)] The reason why we did not choose 
\[
G_1=\sum_{j=1}^{\MM} [\gh(f_j)\cgh(f_j) +\gh(if_j)\cgh(if_j) ]
\]
as our ghost number operator is that while it gives the correct commutation relation with elements in $\cA_g$ it  does not annihilate the vacuum, and $G_1^{\dag}\neq -G_1$. That is it will not serve as a good number operator on the state space. In fact by the commutation relations we get $G_1+1=G$.
\end{itemize}
\end{rem}

As $\cH=\cH_1\oplus J\cH_1$ is even dimensional, we have an even number of linearly independent ghosts, \ie $\{\gh(f_j), \gh(if_j)\,|\, j=1,\dots,\MM\}$. 
To deal with with an odd number of ghosts define,
\[
\rga(\cH_2):= \osalg{ \gh(f_j), \cgh(f_j)\,| \, (f_j)_{j=1}^m\,\text{a complex orthonormal basis of $\cH_2$}},
\]
which we denote $\rga$ when no confusion will arise. Note that $\rga$ makes sense for $m$ being either even or odd but in the case that $m$ is odd we get $\rga$ has an odd number of ghosts $\gh(f_j)$. Hence when we want an odd number of ghosts we use $\rga$.
\begin{rem}\label{rm:oddghrem}
\begin{itemize}
\item[(i)] Note that $\rga$ is a proper subalgebra of $\cA_g=\osalg{C(f) \,|\, f \in \cH=\cH_1 \oplus \cH_2}$ since $f \to \gh(f)$ and $f \to \cgh(f)$ are only real linear and $(f_j)_{j=1}^m$ is a complex basis. A problem with using $\rga$ is that we cannot recover the  $c(f)$'s and $c(f)^*$'s from it. 
\item[(ii)] A natural representation to consider is the action of $\rga$ on 
\begin{equation*}
\rgs_{\psi}:=\rga\psi,\qquad \psi \in \cH_g \backslash \{0\},
\end{equation*}
for $\psi$ chosen such that $J_g\rgs_{\psi}=\rgs_{\psi}$, \ie such that $\rgs_{\psi}$ is a Krein space (note that $\rgs$ is complete as it is finite dimensional). We can choose $\psi$ such that $\rga$ acts irreducibly on $\rgs_{\psi}$ and the corresponding representation is isomorphic to the Berezin representation as will be discussed in the next section. However, we will see below that the choice $\psi=\Omega$ does not give an irreducible representation. 
\end{itemize}
\end{rem}

\subsection{Connection to Berezin Superfunctions}\label{sbs:brz}
A common representation for the ghost algebra is of operators acting on Berezin superfunctions, and is described formally for the case of finite ghosts in Berezin \cite{Be1966} sections 3.1, 3.2, \HT \cite{HenTei92} Chapter 20.2,  \cite{McPat1989I} p482. In this section, we give a well-defined interpretation of these formal definitions by constructing a representation of $\rga$ that can be identified with the formal Berezin superfunction representation's vector space structure and indefinite inner product. This Berezin representation is useful as it gives an irreducible representation of the ghost algebra for an odd number of ghosts (cf. \eqref{lm:rghlm}). We will also discuss the formal formulas in the literature for the Berezin integral, products of superfunctions and involution on superfunctions, but we will \emph{not} give rigorous interpretations of these as they are used to construct the formal indefinite inner product and will not be needed explicitly in this thesis. 

The motivation and purpose for this section is so that we can connect to the literature and discuss related problems, as done Remark \eqref{rm:spfstuff} at the end of the section. After this section we will only use the formal Berezin terminology when referring to the heuristic literature. We summarize in Definition \eqref{df:ghsp} below, the definitive ghost stuctures used in the rest of this thesis.

Another rigorous treatment of the heuristic Berezin calculus is in \cite{Rob1999} but we do not discuss this here since this will take us too far afield, as we would need to connect it to both our approach and the heuristic approach. The Ghost-Fock algebras used here are sufficient for our purposes.

In this subsection we consider both the cases of an odd and an even number of ghosts. Note however that in the case of an even number of ghosts the Berezin superfunctions are redundant as we can use the full ghost algebra as discussed in the previous section. 

Let $(f_j)_{j=1}^{m}$ be an orthonormal basis of $\cH_2$ where $m\in \mathbb{N}$. Define $\gh_j:=\gh(f_j)$ and $\cgh_j:=\cgh(f_j)=\gh_j^*$. We will construct an irreducible representation for $\rga$  with a cyclic vector $\Omega_{sf}$ such that $\cgh_j$ annihilates $\Omega_{sf}$ for all $1 \leq j \leq \MM$. 
\begin{lemma}\label{lm:sfvac}
Let $\psi:=\cgh_{1} \dots \cgh_{\MM}( \gh_{1} \dots \gh_{\MM} + (i)^{\MM(\MM-1)/2}\one) \Omega$. Then the unit vector\\
$\Omega_{sf}:=(\sqrt{2}\norm{\gh_{1} \dots \gh_{\MM}\cgh_{1} \dots \cgh_{\MM}\Omega}_g)^{-1}\psi$ satisfies: 
\begin{itemize}
\item[(i)] $\cgh_{j} \Omega_{sf}=0$ for $1 \leq j \leq \MM$, \qquad \text{and}\label{eq:cgann}
\item[(ii)]  \label{eq:spsgint} $\iip{\Omega_{sf}}{\Omega_{sf}}_g=0, \qquad \iip{ \Omega_{sf}}{ \gh_{1} \ldots \gh_{j}\Omega_{sf}}_g= 
\begin{cases}
(i)^{\MM(\MM-1)/2}, \qquad \text{for}\quad j= \MM\\
0, \qquad \text{for}\quad j<\MM
\end{cases},
$ 
\end{itemize}
\end{lemma}
\begin{proof}
(i): follows from $\{\cgh_j,\cgh_k\}=0$ for all $1 \leq j,k \leq m$, in particular $\cgh_j^2=0$ for all $1\leq j \leq m$. Using $\cgh_j^{\dag}=\cgh_j$ and $\cgh_j^2=0$ we get $\iip{\Omega_{sf}}{\Omega_{sf}}_g=0$.

Now using $\fH_{n}[\perp]\fH_m$ for $n\neq -m$ we calculate,
\begin{align*}
\iip{\psi}{\gh_{1} \dots \gh_{\MM}\psi}_g= & \; \iip{\cgh_{1} \dots \cgh_{\MM}\gh_{1} \dots \gh_{\MM}\Omega}{(i)^{\MM(\MM-1)/2}\gh_{1} \dots \gh_{\MM}\cgh_{1} \dots \cgh_{\MM}\Omega}_g\\
&+\iip{(i)^{\MM(\MM-1)/2}\cgh_{1} \dots \cgh_{\MM}\Omega}{\gh_{1} \dots \gh_{\MM}\cgh_{1} \dots \cgh_{\MM}\gh_{1} \dots \gh_{\MM}\Omega}_g,\\
= & \;\ip{J_g\cgh_{1} \dots \cgh_{\MM}\gh_{1} \dots \gh_{\MM}\Omega}{(i)^{\MM(\MM-1)/2}\gh_{1} \dots \gh_{\MM}\cgh_{1} \dots \cgh_{\MM}\Omega}_g\\
&+\iip{(i)^{\MM(\MM-1)/2}\cgh_{1} \dots \cgh_{\MM}\Omega}{\gh_{1} \dots \gh_{\MM}\cgh_{1} \dots \cgh_{\MM}\gh_{1} \dots \gh_{\MM}\Omega}_g,\\
= & \;(i)^{\MM(\MM-1)/2}\norm{\gh_{1} \dots \gh_{\MM}\cgh_{1} \dots \cgh_{\MM}\Omega}_g^2\\
&+\overline{(i)^{\MM(\MM-1)/2}}\iip{\gh_{\MM} \dots \gh_{1}\cgh_{1} \dots \cgh_{\MM}\Omega}{\cgh_{1} \dots \cgh_{\MM}\gh_{1} \dots \gh_{\MM}\Omega}_g,\\
= & \;(i)^{\MM(\MM-1)/2}\norm{\gh_{1} \dots \gh_{\MM}\cgh_{1} \dots \cgh_{\MM}\Omega}_g^2+{(i)^{\MM(\MM-1)/2}}\norm{\gh_{1} \dots \gh_{\MM}\cgh_{1} \dots \cgh_{\MM}\Omega}_g^2,\\
= & \;2(i)^{\MM(\MM-1)/2}\norm{\gh_{1} \dots \gh_{\MM}\cgh_{1} \dots \cgh_{\MM}\Omega}_g^2,
\end{align*}
where we used $J_g\Omega=\Omega$ and $J_g\gh_jJ_g=\cgh_j$ for $1\leq j \leq \MM$ for the third equality, and $\gh_{\MM} \dots \gh_{1}=(-1)^{\MM(\MM-1)/2}\gh_{1} \dots \gh_{\MM}$ for the fourth equality, and applying the method of the previous step.

Now define, 
\begin{equation}\label{eq:BzGNop}
G_{sf}:=\left(\sum_{j=1}^{\MM} \gh_j\cgh_j\right) - (\MM/2) \one .
\end{equation}
As $\cgh_j \Omega_{sf}=0$ $\forall j$, we have that $G_{sf}\Omega_{sf}= - (\MM/2)\Omega_{sf}$. Moreover $G_{sf}^{\dag}=-G_{sf}$, $[G_{sf}, \gh_j]= \gh_j$ and $[G_{sf}, \cgh_j]= -\cgh_j$ $\forall j$. That is $G_{sf}$ acts as a ghost number operator with $\Omega_{sf}$ at ghost number $-(\MM/2)$.

We calculate, for $1\leq j < \MM$: 
\begin{align*}
 -(\MM/2) \iip{ \Omega_{sf}}{\gh_{1} \dots \gh_{j} \Omega_{sf}}_g= & \;\iip{G \Omega_{sf}}{\gh_{1} \dots \gh_{j} \Omega_{sf}}_g,\\
= & \;-\iip{ \Omega_{sf}}{G_{sf}\gh_{1} \dots \gh_{j} \Omega_{sf}}_g,\\
= & \;((\MM/2)-j)\iip{ \Omega_{sf}}{\gh_{1} \dots \gh_{j} \Omega_{sf}}_g,\\
\end{align*}
hence $\iip{ \Omega_{sf}}{\gh_{1} \dots \gh_{j} \Omega_{sf}}=0$.
\end{proof}

We define the Berezin ghost representation as,
\begin{definition}
Let $\cA_g(\cH_2)$ be a ghost algebra with $dim(\cH_2)=\MM<\infty$, and $\Omega_{sf} \in \cH_g$ be as above. Then $\cH_{bz}:=\rgs_{\Omega_{sf}}=\rga\Omega_{sf}$ and the Berezin representation $\pi_{bz}$ of $\rga$ is this action of $\rga$ on $\cH_{bz}$, with inner product $\ip{\cdot}{\cdot}_g$. 
\end{definition}
\begin{lemma}\label{lm:rghlm}
Let $\cA_g(\cH_2)$ be a finite dimensional ghost algebra with $\dim(\cH_2)=m$. Then: 
\begin{itemize}
\item[(i)] $S:=\{\Omega_{sf}, \gh_1\Omega_{sf},\ldots,\gh_1\ldots \gh_\MM \Omega_{sf}\}$ is a $\mathbb{C}$-linear basis for $\cH_{bz}$, hence $\dim(\cH_{bz})=2^m$.
\item[(ii)] $J_g \Omega_{sf}= (-i)^{m(m-1)/2}\pi_{bz}(\gh_1\ldots \gh_m) \Omega_{sf}$, and  $J_g\cH_{bz}=\cH_{bz}$, where $J_g$ is the fundamental symmetry on $\cH_g$. Hence $\cH_{bz}$,is a Krein space. 
\item[(iii)]$\pi_{bz}:\rga \to B(\cH_{bz})$ is  irreducible. 
\end{itemize}
\end{lemma}
\begin{proof}(i): For this proof we assume that we are in the representation $\pi_{bz}$ and not use the notation explicitly. For any $\{k_1,\ldots,k_l\} \subset \{ 1, \ldots ,\MM\}$ with $k_1 < k_2 < \ldots < k_l$ define
\[
\xi_{k_1,\ldots,k_l}:=\gh_{k_1}\ldots \gh_{k_l}\Omega_{sf}, \qquad M_{k_1,\ldots,k_l}:=\cgh_{k_1}\ldots\cgh_{k_l}\gh_{1} \ldots \widehat{\gh_{k_1}} \ldots \widehat{\gh_{k_l}} \ldots \gh_m,
\]
where $\widehat{\gh_j}$ denotes omission. Now using $\gh_j^2=0$, $\cgh_j\Omega_{sf}=0$ and the commutation relations \eqref{eq:ghcomrelrig} we calculate, $M_{j_1,\ldots,j_n}\xi_{k_1,\ldots,k_l}$. We have three cases.

\smallskip
\noindent Case 1: $\exists k_i \in \{k_1, \dots, k_l\}$ such that $k_i \notin \{j_1,\ldots,j_n\}$. In this case we have that $M_{j_1,\ldots,j_n}$ has a ghost term $\gh_{k_i}$, and so using the ghost anticommutation relations:
\begin{align*}
M_{j_1,\ldots,j_n}\xi_{k_1,\ldots,k_l}= & \;\pm A \gh_{k_i}\gh_{k_i}\gh_{k_1}\ldots \widehat{\gh_{k_i}}\ldots \gh_{k_l} \Omega_{sf}=0,
\end{align*}
where $A$ is a monomial of conjugate ghosts with indices $j_1,\dots j_n$ and ghosts $k_1, \dots, \hat{k_i}, \ldots, {k_l}$.

\smallskip
\noindent Case 2: $\exists j_i \in \{j_1,\ldots,j_n\}$ such that $j_i \notin \{k_1,\ldots,k_l\}$. Then $M_{j_1, \ldots, j_n}$ contains a conjugate ghost $\cgh_{j_i}$ and so,
\begin{align*}
M_{j_1,\ldots,j_n}\xi_{k_1,\ldots,k_l}= & \;\pm B \cgh_{j_i} \gh_{k_1}\ldots \gh_{k_l} \Omega_{sf},\\
= & \;\pm B  \gh_{k_1}\ldots \gh_{k_l} \cgh_{j_i} \Omega_{sf},\\
= & \;0
\end{align*}
where $B$ is a monomial of conjugate ghosts with indices $j_1,\dots,\hat{j_i},\dots j_n$ and ghosts $k_1, \dots, {k_l}$, and we used $\cgh_{j_i}\Omega_{sf}=0$.

\smallskip
\noindent Case 3: $\{j_1,\dots,j_n\}=\{k_1,\ldots,k_l\}$.
\begin{align*}
M_{j_1,\ldots,j_n}\xi_{k_1,\ldots,k_l}= & \;\cgh_{j_1}\ldots\cgh_{j_n}\gh_{1} \ldots \widehat{\gh_{j_1}} \ldots \widehat{\gh_{j_n}} \ldots \gh_m\gh_{j_1}\ldots \gh_{j_m}\Omega_{sf},\\
= & \;\pm \cgh_{j_1}\ldots \cgh_{j_n} \gh_{1}\ldots \gh_{m} \Omega_{sf},\\
= & \;\pm \gh_1\ldots \widehat{\gh_{j_1}} \ldots \widehat{\gh_{j_n}} \ldots \gh_m \Omega_{sf},\\
\neq&\;0
\end{align*}
where we used the fact that the ghosts and conjugate ghosts with different indices anticommute in the first line, that $\{\gh_i,\cgh_i\}=\one$ and that the conjugate ghosts annihilate the vacuum in the second, and lemma \eqref{lm:sfvac} (ii)  in the last line. Summarizing these three cases:
\begin{equation}\label{eq:linbzmon}
M_{j_1,\ldots,j_n}\xi_{k_1,\ldots,k_l}=
\begin{cases}
\pm \gh_1\ldots \widehat{\gh_{j_1}} \ldots \widehat{\gh_{j_n}} \ldots \gh_m \Omega_{sf} \qquad \text{for $\{j_1,\ldots,j_m\}=\{k_1,\ldots,k_l\}$ }\neq \emptyset,\\
0  \qquad \text{otherwise.}
\end{cases}
\end{equation}
To see that $S$ is linearly independent, let:
\begin{align*}
0= & \;\alpha_0\Omega_{sf} + \sum_{j=1}^{m} \alpha_{j} \gh_j \Omega_{sf}+ \sum_{1=j_1<j_2}^{m}\alpha_{j_1, j_2} \gh_{j_1}\gh_{j_2}\Omega_{sf}+ \ldots + \alpha_{1,\ldots, m} \gh_{1}\dots \gh_{m}\Omega_{sf},\\
= & \;\alpha_0\Omega_{sf} + \sum_{j=1}^{m} \alpha_{j} \xi_j+ \sum_{1=j_1<j_2}^{m}\alpha_{j_1 ,j_2}\xi_{j_1,j_2}+ \ldots + \alpha_{1\ldots, m}\xi_{1,\dots,m},
\end{align*}
where $\alpha_{0},\ldots,\alpha_{1,\ldots, m} \in \mathbb{C}$. Now act with $M_{j_1,\ldots,j_n}$ on this linear combination where $j_1 < \ldots < j_n \in \{1,\ldots,m\}$, then using equation \eqref{eq:linbzmon} gives
\[
0=\pm \alpha_{j_1 \ldots j_n} \gh_1\ldots \widehat{\gh_{j_1}} \ldots \widehat{\gh_{j_n}} \ldots \gh_m \Omega_{sf} \implies \alpha_{j_1 \ldots j_n}=0, \qquad \forall j_i, n \in \{1,\ldots,m\}
\]
Hence $0=\alpha_0 \Omega_{sf}$ and hence $\alpha_0=0$. Therefore the set $S$ is linearly independent. Also $S$ spans $\cH_{bz}$ as $\cgh_j \Omega_j=0$ and $\cH_{bz}=\rga\Omega_{sf}$, and so we have $\mathrm{dim}(\cH_{bz})$ $=$ the number of vectors in $S$ $=$ $2^\MM$. 

\smallskip
\noindent (ii): We have $\gh_j \cgh_j \gh_j=\gh_j(\one-\gh_j \cgh_j)=\gh_j-\gh_j^2\cgh_j=\gh_j$ and the brackets $\{\gh_j,\gh_k\}=\{\gh_j,\cgh_k\}=0$ for $k \neq j$, and so
\begin{align}
\gh_{1} \dots \gh_{\MM} \cgh_{1} \dots \cgh_{\MM}\gh_{1} \dots \gh_{\MM}= & \;\gh_{2} \dots \gh_{\MM}\gh_1 \cgh_{1}\gh_1 \cgh_2 \dots \cgh_{\MM}\gh_2 \dots \gh_{\MM},\notag\\
= & \;\gh_{2} \dots \gh_{\MM}\gh_1 \cgh_{2} \dots \cgh_{\MM}\gh_2 \dots \gh_{\MM},\notag\\
= & \;(-1)^{(\MM-1)}\gh_{1} \dots \gh_{\MM} \cgh_2  \dots \cgh_{\MM}\gh_2 \dots \gh_{\MM},\notag\\
\vdots&\notag\\
= & \;(-1)^{\sum_{k=1}^{\MM}(\MM-k)}\gh_{1} \dots \gh_{\MM}, \notag \\
= & \;(-1)^{\MM(\MM-1)/2}\gh_{1} \dots \gh_{\MM} \label{eq:ghred}
\end{align}
Using $J_g \gh_k J_g=\cgh_k$, $J_g\Omega=\Omega$, and the brackets \eqref{eq:ghcomrelrig} we see that,
\begin{align*}
J_g\Omega_{sf}= & \;C\gh_{1} \dots \gh_{\MM}( \cgh_{1} \dots \cgh_{\MM} + (i)^{\MM(\MM-1)/2}\one) \Omega,\\
= & \;C[\gh_{1} \dots \gh_{\MM} \cgh_{1} \dots \cgh_{\MM}\Omega + (i)^{\MM(\MM-1)/2}\gh_{1} \dots \gh_{\MM}\one)] \Omega,\\
= & \;C[\gh_{1} \dots \gh_{\MM} \cgh_{1} \dots \cgh_{\MM}\Omega + (-i)^{\MM(\MM-1)/2}\gh_{1} \dots \gh_{\MM} \cgh_{1} \dots \cgh_{\MM}\gh_{1} \dots \gh_{\MM}\one)] \Omega,\\
= & \;(-i)^{\MM(\MM-1)/2}\gh_{1} \dots \gh_{\MM} [C\cgh_{1} \dots \cgh_{\MM}(\gh_{1} \dots \gh_{\MM}+ (i)^{\MM(\MM-1)/2}\one) ]\Omega,\\
= & \;(-i)^{\MM(\MM-1)/2}\gh_{1} \dots \gh_{\MM} \Omega_{sf},
\end{align*}
where $C \in \mathbb{R}_{+}$ is the normalization constant, and we used \eqref{eq:ghred} in the third equality. Now if we take any vector $\psi \in [S]=\cH_{bz}$, then it follows from the above equality, $J_g \gh_k J_g=\cgh_k$, $\cgh_k \Omega_{sf}=0$, and the ghost commutation relations that $J_g \psi \in [S]=\cH_{bz}$.

\smallskip \noindent (iii): Now let $0\neq \psi\in \cH_{bz}$. Using equation \eqref{eq:linbzmon} and the decomposition of $\psi$ in terms of the basis $S$ as above, there exists $M_{j_1,\dots,j_n}$ such that, 
\[
M_{j_1,\dots,j_n}\psi=\lambda \gh_1\ldots \widehat{\gh_{j_1}} \ldots \widehat{\gh_{j_n}} \ldots \gh_m \Omega_{sf}\neq 0, \qquad \lambda \in \mathbb{C}.
\]
Thus,
\[
\cgh_1\ldots \widehat{\cgh_{j_1}} \ldots \widehat{\cgh_{j_n}} \ldots \cgh_m M_{j_1,\dots,j_n}\psi= \lambda \Omega_{sf}.
\]
where we used the ghost anticommutation relations and $\cgh_j\Omega_{sf}=0$ for the last equality. Therefore we have that for all $\psi \in \cH_{bz}$ there exists $A \in \rga$ such that $A\psi=\Omega_{sf}$ and as $\Omega_{sf}$ is cyclic we have that every vector in $\cH_{bz}$ is cyclic. Hence $\pi_{bz}:\rga \to B(\cH_{bz})$ is  irreducible.  
\end{proof}
In equation \eqref{eq:BzGNop} we defined $G_{sf}$ and we will use the same notation for its restriction to $\cH_{bz}$. Denote $G_{sf}$ restricted to $\cH_{bz}$ the \emph{{Berezin Ghost Number Operator}} and it gives a grading on $\cH_{bz}$: 
\begin{proposition}\label{pr:bsghgrad}
Let $G_{sf}$ be the Berezin Ghost Number Operator. Then
\begin{itemize}
\item[(i)]$G_{sf}\Omega_{sf}= - (\MM/2)\Omega_{sf}$, $G_{sf}^{\dag}=-G_{sf}$, $[G_{sf}, \gh_j]= \gh_j$, and $[G_{sf}, \cgh_j]= -\cgh_j$ $\forall j$.
\item[(ii)] We have the decomposition:
\[
\cH_{bz}=\bigoplus_{k=0}^{\MM} \cH_{bz}^{(k-(\MM/2))}
\]
where $\cH_{bz}^{j}:=\{\psi \in \cH_{bz}\,|\, G_{sf}\psi=j\psi\}$ and $j
\in \{-\MM/2, \MM/2+1,\ldots, \MM/2\}$.
\item[(iii)] 
Furthermore we have that $\cH_{bz}^{j}[\perp]\cH_{bz}^{k}$ for $j \neq -k$ with respect to $\iip{\cdot}{\cdot}_g$, in particular each $\cH_{bz}^{j}$ is a neutral space with respect to $\iip{\cdot}{\cdot}_g$ for $j\neq 0$.
\item[(iv)] Let $\xi:=\cgh_{1} \dots \cgh_{\MM}( \gh_{1} \dots \gh_{\MM} - (i)^{\MM(\MM-1)/2}\one) \Omega\in \rgs_{\Omega}=\rga\Omega$. Then $\xi \perp \cH_{bz}$.
\end{itemize}
\end{proposition} 
\begin{proof}
(i): First $G_{sf}\Omega_{sf}= - (\MM/2)\Omega_{sf}$, $G_{sf}^{\dag}=-G_{sf}$, $[G_{sf}, \gh_j]= -\gh_j$, $[G_{sf}, \cgh_j]= -\cgh_j$ $\forall j$ follow immediately from lemma \eqref{lm:sfvac} equation \eqref{eq:BzGNop} and the ghost anticommutation relations \eqref{eq:ghcomrelrig}. From these relations we see that $G_{sf}\xi_{j_1,\ldots,j_n}=-(\MM/2-n)\xi_{j_1,\ldots,j_n}$ where we are using the notation as in the proof of lemma \eqref{lm:rghlm}. 

\smallskip
\noindent(ii): Now $G_{sf}^*=G_{sf}$ and so eigenspaces of $G_{sf}$ with different eigenvalues will be orthogonal with respect to $\ip{\cdot}{\cdot}_g$. By lemma \eqref{lm:rghlm} (i), the set $S:=\{ \Omega_{sf}, \xi_{j_1,\dots,j_n}\,|\, 1\leq j_1 < \ldots < j_n \leq n \}$ spans $\cH_{bz}$ and so the grading decomposition,
\[
\cH_{bz}=\oplus_{j=0}^{\MM} \cH_{bz}^{(j-(\MM/2))},
\]
follows.

\smallskip
\noindent(iii): This follows from $G_{sf}^{\dag}=-G_{sf}$. 

\pfit (iv): It follows from $\cgh_j^2=0$ for all $1\leq j \leq m$ and the definition of $G_{sf}$  (equation \eqref{eq:BzGNop}) that $G_{sf}\xi=-(\MM/2)\xi$. Hence as eigenvectors with distinct eigenvalues of $G_{sf}$ are orthogonal, it follows that $\xi \perp \bigoplus_{k=1}^{\MM} \cH_{bz}^{(k-(\MM/2))}$. 

To complete the proof we need to show that $\xi \perp \Omega_{sf}$. Recall from lemma \eqref{lm:sfvac} that $\Omega_{sf}=C\psi$ where $\psi=\cgh_{1} \dots \cgh_{\MM}( \gh_{1} \dots \gh_{\MM} + (i)^{\MM(\MM-1)/2}\one) \Omega$ and $C$ is the normalisation constant. Using $\fH_{n}[\perp]\fH_j$ for $n\neq -j$ (Proposition \eqref{pr:ghgrad}) we calculate,
\begin{align*}
\iip{\xi}{\psi}_g= & \; \iip{\cgh_{1} \dots \cgh_{\MM}\gh_{1} \dots \gh_{\MM}\Omega}{\cgh_{1} \dots \cgh_{\MM}\gh_{1} \dots \gh_{\MM}\Omega}_g\\
&+\iip{-(i)^{\MM(\MM-1)/2}\cgh_{1} \dots \cgh_{\MM}\Omega}{(i)^{\MM(\MM-1)/2}\cgh_{1} \dots \cgh_{\MM}\Omega}_g,\\
= & \; \iip{\gh_{1} \dots \gh_{\MM}\Omega}{\gh_{\MM} \dots \gh_{1}\cgh_{1} \dots \cgh_{\MM}\gh_{1} \dots \gh_{\MM}\Omega}_g\\
&-\iip{\cgh_{1} \dots \cgh_{\MM}\Omega}{\cgh_{1} \dots \cgh_{\MM}\Omega}_g,\\
\intertext{using $\gh_{\MM} \dots \gh_{1}=(-1)^{\MM(\MM-1)/2}\gh_{1} \dots \gh_{\MM}$ and equation \eqref{eq:ghred} gives,}
= & \; \iip{\gh_{1} \dots \gh_{\MM}\Omega}{\gh_{1} \dots \gh_{\MM}\Omega}_g -\iip{\cgh_{1} \dots \cgh_{\MM}\Omega}{\cgh_{1} \dots \cgh_{\MM}\Omega}_g,\\
= &\;0
\end{align*}
where we used $J_g\Omega=\Omega$, $J_g\gh_jJ_g=\cgh_j$ for $1\leq j \leq \MM$ and that $J_g$ is unitary for the last equality.
\end{proof}
We call the above decomposition $\cH_{bz}$ the \emph{{Berezin Ghost Grading}}. As $G$ and $G_{sf}$ have the same commutation relations with $\gh_j$ and $\cgh_j$, we have that $[G,A]=nA$ iff $[G_{sf},A]=nA$ for $A \in \rga$, $n\in \mathbb{Z}$, we do not need to define a Berezin ghost grading of $\rga$ as it would coincide with the ghost grading already defined. Note however that the spatial and algebra gradings do not match in the sense that,
\[
\cH_{bz}^{j-(m/2)}=\sqb{(\cG_j \cap \rga)\Omega_{sf}}, \qquad j=0,\dots m,
\]
which is due to the fact that $\Omega_{sf}$ is the state with `no ghosts' but has Berezin ghost number $-(m/2)$.

For the remainder of this section we will assume that we are in the Berezin representation of $\rga$ and not denote it explicitly. Using the Berezin representation of $\rga$, we now define the space of Berezin superfunctions which corresponds to the usual heuristic definitions (cf. references at the beginning of this subsection). Assume that $\cH_0$ is a fixed Hilbert space with inner product $\ip{\cdot}{\cdot}_0$ and define,
\[
\sfu:= \cH_0 \otimes \cH_{bz}.
\]
which has the tensor Hilbert inner product, $\ip{\cdot}{\cdot}_{sf}$, and Krein inner product $\iip{\cdot}{\cdot}_{sf}:=\ip{\cdot}{(\one\otimes  J_g) \cdot}_{sf}$. Now by lemma \eqref{lm:rghlm}(i), we have that a vector in $\cH_{sf}$ has the form,
\begin{align}
\psi= & \;\psi_0\otimes \Omega_{sf} + \sum_{a=1}^{\MM} \psi_a\otimes \gh_a \Omega_{sf}+ \sum_{1=a_1<a_2}^{\MM}\psi_{a_1 a_2}\otimes \gh_{a_1 a_2}\Omega_{sf}+ \ldots + \psi_{1 \ldots \MM} \otimes \gh_{1}\dots \gh_{\MM}\Omega_{sf},\notag \\
=:&\,\psi_0 + \psi_a \gh_a +  \ldots + \psi_{1 \ldots \MM} \gh_{1}\dots \gh_{\MM} \label{eq:dfsfu},
\end{align}
where $\psi_0, \ldots, \psi_{1 \ldots \MM} \in \cH_0$, and the last line is the heuristic expression for a superfunction, where the operators $\gh_k$ are now reinterpreted as Grassman variables in a formal polynomial with coefficients in $\cH_0$.
\begin{rem}
In the above definitions we have written the ghost terms in order of increasing index. Another common convention is writing superfunctions $\psi$ in sums of all permutations of indices, ie,
\[
\psi=\psi_0 + \sum_{a=1}^{\MM} \psi_a \gh_a +  \ldots + \sum_{1=a_1,\dots, a_\MM}^{\MM}\psi_{a_1 \ldots a_\MM} \gh_{a_1}\dots \gh_{a_\MM}.
\]
(see \cite{Be1966} p50,  \HT \cite{HenTei92} p319). These are equivalent by using the ghost commutation relations, but will have coefficient functions $\psi_{a_1 \ldots a_j}$ differing by a constant. Also using increasing indices will give slightly different formulas in the definition of the integrals and inner products below.
\end{rem}
We can check that lemma \eqref{lm:sfvac} \eqref{eq:spsgint} gives for $\psi, \xi \in \cH$, via the decomposition in equation \eqref{eq:dfsfu} that:
\begin{align}
\iip{\psi}{\xi}_{sf}= & \;(i)^{\MM(\MM-1)/2} \Big( \ip{\psi_0}{ \xi_{1 \dots \MM}}_0 + (-1)^{\MM(\MM-1)/2} \ip{\psi_{1 \dots \MM}}{ \xi_0}_0 \notag  \\
& +  \sum_{j=1}^{\MM-1} \sum_{a_1=1 < \ldots <a_j}^{\MM} (-1)^{((\sum_{k=1}^{j} a_k)-j)} \ip{\psi_{a_1 \dots a_j}} {\xi_{1 \ldots \hat{a_1} \ldots \hat{a_j} \ldots \MM}}_0 \Big) \label{eq:dfsfuiip}.
\end{align}
We call the representation $\pi_{sf}:\rga \to B(\sfu)$  defined by $\pi_{sf}:=\one \otimes \pi_{bz}$ the \textit{Berezin superfunction representation} of $\rga$, and for the remainder of this section we will not denote $\pi_{sf}$ explicitly.

We connect the above construction with the usual description of the Berezin superfunctions as follows. Formally, the Berezin superfunctions are defined by the last line of equation \eqref{eq:dfsfu} where $\cH_0$ is assumed to be an $L^2(\mu)$ space, and act as the coefficients of the Grassman variables $\gh_j$ (cf. references at the beginning of this subsection). Then a formal integration rule, formal multiplication, and formal conjugation is assumed on $\sfu$ and are used to give an indefinite inner product on the $\sfu$. We give these formal rules below but do \emph{not} try to make sense of them rigorously, we just show that they produce the same inner product as the formula \eqref{eq:dfsfuiip}. 

The Berezin integral is formally defined on $\sfu$ as,
\[
\int  d \gh_1\ldots d \gh_\MM\,\psi:= \int  d \mu \,\psi_{1 \ldots \MM}.
\]
and is referred to as integrating the `top function'. This is a formal definition and we are not assuming that $d \gh_1\ldots d \gh_\MM$ defines a measure on the ghosts. Let $\psi, \xi \in \sfu$. Formal multiplication $\psi\xi$ is defined to be multiplication of the terms in $\psi$ and $\xi$ with respect to the decomposition given by equation \eqref{eq:dfsfu}, with the coefficients $\psi_{a_1 \ldots a_n}$, $\xi_{a_1 \ldots a_l}$ being multiplied pointwise as ordinary $L^2(\mu)$ functions producing an $L^1(\mu)$ function, and the $\gh_{a_1}\dots \gh_{a_\MM}$ terms being multiplied and rearranged into increasing index order using the ghost anticommutation relations. Note that repeated indices in multiplied terms (denoted $a_k$) will give $\gh_{a_k}^2=0$ in the ghost part of that term, hence that term will equal $0$. See Example \eqref{ex:bzm2} below for a sample calculation. Formal conjugation on $\sfu$ is defined by:
\begin{align*}
\psi^{\dag}= & \;\overline{\psi_0} + \overline{\psi_a} \gh_a^{\dag} + \ldots + \overline{\psi_{1 \ldots \MM}} (\gh_{1}\ldots \gh_{\MM})^{\dag},\\
= & \;\overline{\psi_0} + \overline{\psi_a} \gh_a + \ldots + \overline{\psi_{1 \ldots \MM}} \gh_{\MM}\ldots \gh_{1},\\
= & \;\overline{\psi_0} + \overline{\psi_a} \gh_a + \ldots + (-1)^{\MM(\MM-1)/2}\overline{\psi_{1 \ldots \MM}} \gh_{1}\ldots \gh_{\MM},
\end{align*}
where $\overline{\psi}$ is complex conjugation in $L^2(\mu)$ and the formal conjugation $\dag$ acts as $(\gh_{j_1}\ldots \gh_{j_n})^{\dag}=(\gh_{j_n}\ldots \gh_{j_1})$. 

Using the formal Berezin integral, multiplication and conjugation, the formal inner product on $\sfu$ is defined by:
\begin{align*}
\iip{\psi}{\xi}_{sf}:=  (i)^{\MM(\MM-1)/2}\int  d \mu\, d\gh_1\ldots d \gh_\MM \, \psi^{\dag} \xi,
\end{align*}
A straightforward calculation shows that the above inner product on $\sfu$ defined using the formal Berezin integral agrees with the rigorous equation \eqref{eq:dfsfuiip}. We give a sketch of the calculation below for the case of $m=2$ (Example \eqref{ex:bzm2}) with the general case following in a similar way. 
\begin{eje}\label{ex:bzm2}
Take $\MM=2$, $\cH_0=L^2( \mu)$, and $\psi, \xi \in \sfu$, i.e,
\[
\psi= \psi_0 + \psi_1 \gh_1 + \psi_2 \gh_2 + \psi_{12} \gh_1 \gh_2, \qquad \xi= \xi_0 + \xi_1 \gh_1 + \xi_2 \gh_2 + \xi_{12} \gh_1 \gh_2.
\]
Then:
\begin{align*}
\psi^{\dag}\xi=\overline{\psi_0}\xi_0+(\overline{\psi_1}\xi_0+&\overline{\psi_0}\xi_1)\gh_1\\
+&(\overline{\psi_2}\xi_0+\overline{\psi_0}\xi_2)\gh_2+(\overline{\psi_0} \xi_{1 2} - \overline{\psi_{1 2}}\xi_0 + \overline{\psi_1}{\xi_2}- \overline{\psi_2}{\xi_1})\gh_1\gh_2.
\end{align*}
Therefore, the top function of $\psi^{\dag}\xi$ is:
\[
(\overline{\psi_0} \xi_{1 2} - \overline{\psi_{1 2}}\xi_0 + \overline{\psi_1}{\xi_2}- \overline{\psi_2}{\xi_1})\gh_1\gh_2.
\]
Hence the formal inner product is:
\begin{align*}
\iip{\psi}{\xi}_{sf}= & \;i\int \rd \mu \,(\overline{\psi_0} \xi_{1 2} - \overline{\psi_{1 2}}\xi_0 + \overline{\psi_1}{\xi_2}- \overline{\psi_2}{\xi_1}),\\
= & \;(i)^{\MM(\MM-1)/2} \Big( \ip{\psi_0}{ \xi_{1 2}}_0- \ip{\psi_{1 2}}{ \xi_0}_0  +   \ip{\psi_1}{ \xi_{ 2}}_0- \ip{\psi_{ 2}}{ \xi_1}_0 \Big),
\end{align*}
which we see is the same as equation \eqref{eq:dfsfuiip} with $\MM=2$. The factor of $i$ in front of the integral is needed to ensure that $\overline{\iip{\psi}{\xi}}_{sf}=\iip{\xi}{\psi}_{sf}$.
\end{eje}
In the calculation of $\iip{\psi}{\xi}_{sf}$ for general $m$, an important point is that the only term that contributes $\iip{\psi}{\xi}_{sf}$ is the integral of the `top function' of $\psi^{\dag} \xi$, and so this the only term we have to calculate in the product. The sign factor $(-1)^{((\sum_{k=1}^{j} a_k)-j)}$ in front of the $\ip{\psi_{a_1 \dots a_j}} {\xi_{1 \ldots \hat{a_1} \ldots \hat{a_j} \ldots \MM}}_0= \int \rd \mu \overline{\psi_{a_1 \dots a_j}}{\xi_{1 \ldots \hat{a_1} \ldots \hat{a_j} \ldots \MM}}$ term in equation \eqref{eq:dfsfuiip} formally comes from rearranging the indices of the ghost terms in the top function into increasing order. This allows us to identify the heuristic superfunctions with the Hilbert space $\sfu=\cH_0 \otimes \cH_{bz}$.

Formally define left differentiation by Grassman variables by: 
\[
\frac{\partial}{\partial \gh_a}\psi_{j_1,\ldots,j_n}\gh_{j_1}\ldots\gh_{j_n} =
\begin{cases}
(-1)^{i-1}\psi_{j_1,\ldots,j_n}\gh_{j_1}\ldots\widehat{\gh_{j_i}}\ldots\gh_{j_n},\qquad a=j_i,\\
0,\qquad \text{otherwise}
\end{cases}
\]
and extending linearly to $\sfu$ (cf. \cite{Be1966} p51,57). Then $\rga$ is heuristically assumed to act on $\sfu$ via,
\begin{equation}\label{eq:ghsfopr}
\hat{\gh_a}\psi=\gh_a\psi, \qquad \hat{\cgh_a}= \frac{\partial}{\partial \gh_a}.
\end{equation}
From equation \eqref{eq:dfsfu} we see that this corresponds to the representation $\pi_{sf}:\rga \to B(\sfu)$, hence we know that it extends to a well-defined representation of $\rga$. This provides a rigorous interpretation of the Berezin superfunction formalism. What is more, from lemma \eqref{lm:rghlm} we know that $\sfu$ is actually a Krein space, we have the explicit form of its fundamental symmetry and we know that $\pi_{sf}:\rga \to B(\sfu)$ is an irreducible representation of $\rga$.

\begin{rem}\label{rm:spfstuff}
\begin{itemize}

\item[(i)] Another natural candidate for a space on which the restricted ghost algebra can act is $\rgs_{\Omega_{g}}=\rga \Omega_{g}$, where $\Omega_{g}$ is the Ghost-Fock vacuum vector. However this will \emph{not} be an irreducible representation. This can be seen as $\Omega_{sf} \in \rgs_{\Omega_{g}}$ implies that $\cH_{bz}=\rga \Omega_{sf} \subset \rgs_{\Omega_{g}}$. But by Proposition \eqref{pr:bsghgrad} (iv) we have that the inclusion is proper. As $\rga$ preserves $\cH_{bz}$ we get that action $\rga$ on $\rgs_{\Omega}$ is reducible, hence we see that the restriction of a pure state of $\cA_g(\cH_2)$ to $\rga$ is no longer pure. 
\item[(ii)] One of the main features of the Berezin representation is that
\begin{equation}
\cgh_a \Omega_{sf} =0, \qquad 1\leq i \leq \MM,
	\label{eq:gstvah}
\end{equation} 
which is an attractive feature since we can think of the ghost and conjugate ghosts as acting as creators and annihilators with vacuum $\Omega_{sf}$. Hence the superfunction representation is often assumed in Hamiltonian BRST to identify a natural physical subspace amongst the `multiple copies of the physical space' (see section \eqref{sbs:hhamexphsp}). If we take a Hamiltonian BRST model with commuting constraints (cf. Section \ref{sec:hamBRST}) with Berezin ghost space then we saw that $\ker \Delta= \cH^0_p\otimes \cH_{bz}$ (equation \eqref{eq:hcomhamdel}), \ie the span of vectors of the form
\[
 \psi_0 + \psi_a \gh_a +  \ldots + \psi_{1 \ldots \MM}\gh_{1}\dots \gh_{\MM}
\] 
where $\psi_0,\psi_a,\psi_{1 \ldots \MM} \in \cH^0_p\otimes \{\mathbb{C}\Omega_{sf}\}$ are Dirac physical vectors. 

The choice of representative for the physical space is then the space spanned by the physical vectors with no ghosts, \ie the $\psi_0$ vectors. Note that the subspace of $\cH_{bz}$ with Berezin ghost number $-(\MM/2)$ is spanned  by $\Omega_{sf}$ and so we can also see this choice of representative space, \ie $\cH^0_p\otimes \mathbb{C}\Omega_{sf}$, as restricting to the subspace with Berezin Ghost number $-(\MM/2)$, by which we get BRST physical space naturally isomorphic to the Dirac physical space. A significant problem with this choice is that if we have a representation of the ghost algebra with  $\dag$-adjoint ghosts and conjugate ghosts, and a vector $\psi$ such that $\cgh_j \psi=0$ for some $j$, then:
\begin{gather*}
\iip{\psi}{\psi}=\iip{ \{\gh_j, \cgh_j\}\psi}{\psi}= \iip{ (\gh_j+ \cgh_j)^2\psi}{\psi}=\\
=\iip{ (\gh_j+ \cgh_j)\psi}{(\gh_j+ \cgh_j)\psi}=\iip{ \gh_j\psi}{\gh_j\psi}=\iip{ \gh_j^2\psi}{\psi}=0.
\end{gather*}
So the span of such vectors $\psi$ will not be a good physical subspace as it is neutral with respect to the indefinite inner product. This is exactly the case for $\cH^0_p\otimes \{\mathbb{C}\Omega_{sf}\}$ in the superfunction representation. This problem is noted in \cite{vHo2006} p7, and a correction is suggested. This problem also occurs due to the Berezin grading structure on $\cH_{bz}$ which we discuss in the next item.

\item[(iii)] We have a grading structure on $\cH_{bz}$ by using the Berezin ghost number operator $G_{sf}$ (cf. equation \eqref{eq:BzGNop}  and Proposition \ref{pr:bsghgrad}(ii)). As $G_{sf}^{\dag}=-G_{sf}$ we get that any eigenspace of $G_{sf}$ with nonzero eigenvalue will be neutral. Hence the only subspaces with definite Berezin ghost number that are positive with respect to $\iip{\cdot}{\cdot}_{sf}$ (hence can serve as a physical subspace) must have Berezin ghost number zero. This shows again that the usual choice of physical space $\cH^0_{p}\otimes \{\mathbb{C} \Omega_{sf}\}$ is problematic, as it has Berezin ghost number $-(m/2)$.

Restricting to the zero ghost subspace is problematic when there are an odd number of ghosts $\dim(\cH)=\MM=2k+1$. Proposition \eqref{pr:bsghgrad}(ii) shows that the eigenvalues of $G_{sf}$ are fractional and hence the non-trivial ghost numbered subspaces have fractional number, \ie the ghost number $0$ subspace is $\{0\}$. Hence restricting to ghost number zero subspaces gives a trivial model (as noted in \cite{LanLin1992} p425).

\item[(iv)] A second method to deal with the neutrality of the natural physical subspace $\cH^0_p \otimes \mathbb{C}\Omega_{sf}$ has been suggested by \cite{LanLin1992} p426 and \cite{vHo2006} p7. One takes a subspace of $\ker \Delta$ spanned by the ghost number $-\MM/2$ and ghost number $\MM/2$ space, ie,
\[
(\ker \Delta)_{+}:=[ \psi_0\otimes(\one +\gh_{1}\dots \gh_{\MM})\Omega_{sf}\,|\,\psi_0 \in \cH^0_p]=\cH^0_p\otimes \{\mathbb{C}(\one +\gh_{1}\dots \gh_{\MM})\Omega_{sf}\}
\]
(where $(\ker \Delta)_{+}=V_{+}$ in the notation of \cite{LanLin1992} p426). $(\ker \Delta)_{+}$ a subspace of indeterminate ghost number, which can be seen as a problem discussed in \cite{vHo2006} p8. 

To the author, indeterminate ghost space number does not seem to be a serious objection for these reasons: First we have already graded our algebra so as to apply the BRST superderivation and so do not need the ghost space to do this. Second, the choice of zero ghost number for the physical subspace was arbitrary in the first place. Third and most importantly, $(\ker \Delta)_{+}$ is a subspace which is positive with respect to $\ip{\cdot}{\cdot}$ and corresponds to $\cH^p_0$ without the MCPS problem (cf. Subsection \ref{sbs:hhamexphsp}).

\item[(v)] Consider the case of Hamiltonian BRST with a finite number of non-abelian constraints as in Subsection \ref{sbs:hhamexphsp}(2). We have that $(C^c_{ab}\one \otimes \gh_a \gh_b \cgh_c)\Omega_{sf} = 0$. In fact we also have $\Omega_{sf}, \gh_1\ldots \gh_\MM \Omega_{sf} \in  \cap_{a}\ker (\Sigma_a)$ as defined in \eqref{eq:sigantbrs}, and so $\cH^p_0 \otimes \mathbb{C}\Omega_{sf}\subset \ker \Delta$ and $\cH^p_0 \otimes \mathbb{C}\Omega_{sf} \subset \ker \Delta$.  Both $\Omega_{sf}$ and $\gh_1\ldots \gh_\MM \Omega_{sf}$ are neutral vectors and so we have that $\cap_{a}\ker (\Sigma_a)$ is not one-dimensional and that $\ker \Delta$ has indefinite inner product. So we have the MCPS problem for non-abelian Hamiltonian BRST as claimed in the end of heuristic Subsection \ref{sbs:hhamexphsp}. Note that  $\ker \Delta_{+}$ as defined in the preceding item is in $\ker \Delta$ and so is a good candidate for the choice to solve the MCPS problem, as long as we are not concerned with definite ghost number for the physical space.
\item[(vi)] The above construction for the Berezin space does not work in the infinite dimensional space as the construction of $\Omega_{sf}$ requires a finite number of $\cgh_j$'s. Also the heuristic definition of integral does not work as we have no top function. This means that a choice such as $(\ker \Delta)_{+}$ to solve the MCPS for any formulation of Hamiltonian BRST with infinite constraints will become problematic.
Infinite dimensional Berezin representations are dealt with in \cite{Rob1999} p136 onwards. We do not discuss these further as the infinite dimensional ghost algebras are sufficient for us and come from the motivating physical example.
\end{itemize}
\end{rem}
Above, we found the two conventions for the ghost algebra:
\begin{itemize}
	\item $\cA_g$ on $\cH_g$ with ghost operator $G=-G^{\dag}$ and vacuum $\Omega$.
	\item $\rga \subset \cA_g$ on $\cH_{bz} \subset \cH_g$ with ghost operator $G_{sf}=-G_{sf}^{\dag}$ and vacuum $\Omega_{sf}$.
\end{itemize}
Both are used in the heuristic literature, but for the analysis below we would like a unified terminology. So we summarize their essential similarities and differences:
\begin{itemize}
\item As remarked in Proposition \eqref{pr:bsghgrad}, the ghost grading of $\cA_g$ with respect to $G$ coincides with the ghost grading of $\rga \subset \cA_g$ with respect to $G_{sf}$.
\item The spatial ghost gradings do \emph{not} coincide;- whereas
\[
\cH_g=\oplus_{n=-m}^{m}\fH_n, \qquad \fH_n=\overline{\cG_n\Omega}, \qquad G\Omega=0,
\]
we have:
\[
\cH_{bz}=\oplus_{k=0}^{m}\cH_{bz}^{(k-m/2)}, \qquad \cH_{bz}^{(k-m/2)}=\overline{(\cG_k\cap \rga)\Omega_{sf}}, \qquad G_{sf}\Omega_{sf}=-(m/2)\Omega_{sf},
\]
\end{itemize}
Hence the following unified terminology is appropriate with respect to the general structure required to define the BRST superderivation.
\begin{definition}\label{df:ghsp}
A \textbf{Ghost Space} $\cH_g$ is either:
\begin{itemize}
\item[(i)] A Hilbert space $\cH_{\omega}$ where $\omega \in \fS_g \subset \fS(\cA_g)$, and it is equipped with the grading structure, $D(G)$, and $G$ as given by Proposition\eqref{pr:ghsp}.
\item[(ii)] The Berezin ghost space $\cH_g=\cH_{bz}$ and it is equipped with the Berezin Ghost Number Operator $G_{sf}$ and Berezin grading structure as given by Proposition \eqref{pr:bsghgrad}. 
\end{itemize}
Once it has been established that we are using the Berezin ghost space, we will also refer to the restricted ghost algebra $\rga$ as the ghost algebra and drop the `$sf$' subscript when no confusion will arise.
\end{definition}
\begin{rem}
In constructions it is preferable to use the full ghost algebra rather than the  Berezin representation of $\rga$ as we can recover the CAR creators and annihilators, as discussed in Remark \eqref{rm:oddghrem}(i). However the Berezin representation of $\rga$ is used in the odd number of ghost case as this is the finite irreducible representation of $\rga$. For the even number of ghosts case we can use the full ghost algebra, but may in examples still use the Berezin representation of $\rga$ so as to stay as close as possible to the heuristic literature.
\end{rem}

\section{Charge and \emph{dsp}-decomposition} \label{sec:dsprig}

In physical examples, such as BRST-QEM, the BRST charge is given as a formal integral of unbounded operators (e.g. equation \eqref{eq:HQ}) which suggests that any well-defined formulation will need to include the case of an unbounded $Q$. The main aim of this section is to prove rigorously a version of the \emph{dsp}-decomposition which is general enough to be applied to the subsequent QEM example. This has been done in \cite{HoVo89}, and here we give slightly different proofs. This includes the case of bounded $Q$. 

Initially we do not assume any Krein or ghost structure. Let $\cH$ be a Hilbert space with inner product $\ip{\cdot}{\cdot}$.

\begin{theorem}[dsp-decomposition]\label{th:Hdsp1}
Let $Q$ be a closed operator, $\ran Q \subset D(Q) \subset \cH$, and $Q^2=0$. Then we have the following decomposition,
\[
\cH=\cH_d \oplus \cH_s \oplus \cH_p = \ker Q \oplus \cH_p, 
\]
where $\cH_s:=\ker Q \cap \ker Q^*$, $\cH_d:=\overline{\ran Q}$, $\cH_p:=\overline{\ran Q^*}$, $\ker Q=\cH_s\oplus \cH_d$
\end{theorem}
\begin{proof}
First, as $Q$ is closed we have that $\ker Q=(\ran Q^*)^{\perp}$ (\cite{Con1985} Proposition X 1.13 p310) and hence we get the first decomposition. As $Q^2=0$ we have that $\ran Q \subset \ker Q$ and so if $\cL$ is the orthogonal complement of $\ran Q$ in $\ker Q$ we have,
\[
\cH= \overline{\ran Q}\oplus \cL \oplus \overline{\ran Q^*}.
\]
Now $Q^2=0$ implies $(Q^*)^2=0$ and so the above decomposition holds for $Q^*$, therefore we have that $\cL \subset \ker Q \cap \ker Q^*=\cH_s$. It is straightforward to check that $\cH_s= \ker Q \cap \ker Q^* \subset  (\overline{\ran Q} \oplus \overline{\ran Q^*})^{\perp}=\cL$ and so $\cH_s=\ker Q \cap \ker Q^*$, which gives the decomposition.
\end{proof}
For the remainder of this subsection $\cH_s,\cH_p,\cH_d$ denote the same spaces as in Theorem \eqref{th:Hdsp1} and $P_s,P_s,P_d$ are their corresponding orthogonal projections. 

The heuristics in Chapter \ref{ch:heu} use $Q$ to select the physical subspace. Correspondingly we define,
\begin{definition}\label{df:opbrsphsp}
The BRST-physical space of $Q$ is,
\begin{equation}\label{df:BRSTsp}
\cH^{BRST}_{phys}:=\ker Q/ \overline{\ran Q}=\ker Q/ \cH_d
\end{equation}
Let $\vp: \ker Q \to \cH^{BRST}_{phys}$ be the factor map, and denote $\hat{\psi}:=\vp(\psi)$ for all $\psi \in \ker Q$. 
\end{definition}
Further motivated by the heuristics, we assume that $\cH$ has a Krein structure with fundamental symmetry $J$ and denote the indefinite inner product by $\iip{\cdot}{\cdot}=\ip{\cdot}{J\cdot}$. Then we have,
\begin{lemma}\label{lm:ranQnll}
Let $Q\in \op(\cH)$ be a closed operator which is Krein symmetric on its domain $D(Q)$, \ie $Q \subset Q^{\dag}$. Moreover, let $\ran Q \subset D(Q)$ and $Q^2=0$. Then 
\[
\ker Q = \cH_s \oplus \cH_d =\cH_s [\oplus] \cH_d
\]
where $\cH_s:=\ker Q \cap \ker Q^*$, $\cH_d:=\overline{\ran Q}$ and $[\oplus]$ denotes orthogonality with respect to the indefinite inner product $\iip{\cdot}{\cdot}$. 
\end{lemma}
\begin{proof}
The first decomposition is given by Theorem \eqref{th:Hdsp1}. Next, assume that $\psi \in \ker Q$, $\xi=Q\phi \in \cH_d$. Then as $\psi \in \ker Q \subset D(Q) \subset D(Q^{\dag})$, we have $\iip{\psi}{Q\phi}=\iip{Q^{\dag}\psi}{\phi}=\iip{Q\psi}{\phi}=0$. That is, $\ran Q [ \oplus] \ker Q$. As $\iip{\cdot}{\cdot}=\ip{\cdot}{J\cdot}$ we have that it is continuous in both arguments and hence $\cH_{d}=\overline{\ran Q} [\oplus] \ker Q$.
\end{proof}

The Krein inner product on $\cH$ now induces an indefinite inner product on $\cH^{BRST}_{phys}$, which is well defined by lemma \eqref{lm:ranQnll}.
\begin{definition}\label{df:BRSTphysindef}
Define an indefinite inner product on  $\cH^{BRST}_{phys}$ as,
\begin{equation*}
\iip{\hat{\psi}}{\hat{\xi}}_p:=\iip{\psi}{\xi}=\iip{P_s\psi}{P_s\xi}, \qquad \forall \psi,\xi \in \ker Q 
\end{equation*}
Note that $\hat{\psi}=\widehat{P_s\psi}$ for all $\psi \in \ker Q$. Moreover, as $\ker Q$ and $\cH_d$ are complete with respect to the Hilbert space topology, the factor space $\cH^{BRST}_{phys}$ is complete with respect to the factor topology.
\end{definition}
An interesting case of the above constructions is when $\cH^{BRST}_{phys}$ becomes a Krein space with the inner product $\iip{\cdot}{\cdot}_p$. The following example shows that the assumptions of lemma \eqref{lm:ranQnll} are not sufficient to imply that $\cH^{BRST}_{phys}$ is a Krein space.
\begin{eje}\label{ex:spdcount}
Let $\cL$ be a Hilbert space and $A \in \op(\cL)$ be a closed, symmetric operator with dense domain $D(A)$. Furthermore, let $A$ be such that $\ker A= \{0\}$ and $\ker A^* \neq \{0\}$ (\cite{ReSi1972v1}). Define $\cH=\cL \oplus \cL$, $D(Q)=D(A^*) \oplus D(A)$ and

\begin{align*}
Q:=\left(
\begin{matrix}
0 & A \\
0 & 0 \\
\end{matrix}
\right),\qquad 
J:=\left(
\begin{matrix}
0 & I \\
I & 0 \\
\end{matrix}
\right),
\end{align*}
From these definitions it follows that $Q^2=0$ and,
\begin{align*}
Q^*=\left(
\begin{matrix}
0 & 0 \\
A^* & 0 \\
\end{matrix}
\right), \qquad
D(Q^*)=\left(
\begin{matrix}
D(A^*) \\
\cL \\
\end{matrix}
\right) 
\end{align*}
and that
\begin{align*}
\ker Q=\left(
\begin{matrix}
D(A^*) \\
0 \\
\end{matrix}
\right), \qquad
\ker Q^*=\left(
\begin{matrix}
\ker A^* \\
\cL \\
\end{matrix}
\right),
 \end{align*}
Hence 
\begin{align*}
\ker Q \cap \ker Q^*=\left(
\begin{matrix}
\ker A^* \\
0 \\
\end{matrix}
\right), 
\qquad
J\left( \ker Q \cap \ker Q^* \right)=\left(
\begin{matrix}
0 \\
\ker A^* \\
\end{matrix}
\right). 
\end{align*}
Hence $J\cH_s =J(\ker Q\cap \ker Q^*)\perp \cH_s$. This shows that $\cH_s$ is a neutral space and that the inner product $\ip{\cdot}{\cdot}_p$ is neutral on $\cH^{BRST}_{phys}$.

\end{eje}
\begin{rem}
In general, $J$ does not preserve $\ker Q$ (cf. Theorem \eqref{th:Hdsp} below) and so $J$ does not factor to $\cH^{BRST}_{phys}$. However, if $J$ preserves a linear section of the factoring map $\pi:\ker Q \to \cH^{BRST}_{phys}$ then the restriction of $J$ to that section can be used to define a fundamental symmetry on $\cH^{BRST}_{phys}$. A fundamental symmetry defined this way depends on the choice of section.
\end{rem}  
The following lemma is an example of the above remark, and will hold when $Q$ is Krein-selfadjoint (cf. Theorem \eqref{th:Hdsp}).  
\begin{lemma}\label{lm:phspksp}
Let $J\cH_s=\cH_s$ and define $J_p \in B(\cH^{BRST}_{phys})$ by $J_p\hat{\psi}:=\widehat{JP_s\psi} $ for all $\psi \in \ker Q$. Then: 
	\begin{itemize}
		\item[(i)]$\cH^{BRST}_{phys}$ is a Krein space with fundamental symmetry $J_p$, and Hilbert inner product:
			\begin{equation}\label{eq:BRSTpHilip}
				\ip{\hat{\psi}}{\hat{\xi}}_p:=\iip{\hat{\psi}}{J_p\hat{\xi}}_p=\ip{P_s\psi}{P_s \xi}, \qquad \forall \psi,\xi \in \ker Q
			\end{equation}
		\item[(ii)] The map $\vp|_{\cH_s}:\cH_s \to \cH^{BRST}_{phys}$ is a Krein and Hilbert isometric isomorphism. 
		\item[(iii)]The space $\cH^{BRST}_{phys}$ is a Hilbert space with respect to the inner product 			$\iip{\hat{\psi}}{\hat{\xi}}_p$ if and only if $JP_s=P_s$.
	\end{itemize}
\end{lemma}
\begin{proof}
(i) As $J$ is a fundamental symmetry $J^*=J$. Combining this with the assumption $J\cH_s=\cH_s$ gives $J\cH_s^{\perp}=\cH_s^{\perp}$, hence $[P_s,J]=0$ and so $J_p$ is well defined. Therefore $(JP_s)^2=P_s$ and $(JP_s)^*=JP_s$ and so $J_p^*=J_p$ and $J_p^2=\one$. From:
\begin{equation*}
\ip{\hat{\psi}}{\hat{\xi}}_p=\iip{\hat{\psi}}{J_p\hat{\xi}}_p=\iip{P_s\psi}{JP_s\xi }=\ip{P_s\psi}{P_s\xi},
\end{equation*} 
we see that $\ip{\hat{\psi}}{\hat{\xi}}_p$ is a positive definite and hence $\cH^{BRST}_{phys}$ is a Krein space, as $\cH^{BRST}_{phys}$ is complete with respect to the factor topology.

\smallskip
\noindent (ii): Obvious from (i)

\smallskip
\noindent (iii):
As $\cH^{BRST}_{phys}$ is a Krein space then $J_p=J=P^p_{+}-P^p_{-}$ where $P^p_{+}$ and $P^p_{-}$ are the projections onto the positive and negative subspaces of $\cH^{BRST}_{phys}$ with respect to $\iip{\hat{\psi}}{\hat{\xi}}_p$.  $\cH^{BRST}_{phys}$ is a Hilbert space with respect to $\iip{\hat{\psi}}{\hat{\xi}}_p$ iff $P^p_{-}=0$ iff $J_p=\one$ iff $JP_s=P_s$ by the definition of $J_p$.
\end{proof}
The above sufficient condition (lemma \eqref{lm:phspksp}) is satisfied in many examples, in particular it will be satisfied whenever $Q$ is Krein selfadjoint. Note however, that it is not satisfied in Example \eqref{ex:spdcount}.
\begin{theorem}\label{th:Hdsp}
Let $Q\in \op(\cH)$ be Krein selfadjoint, \ie $Q = Q^{\dag}$. Moreover, let $\ran Q \subset D(Q)$, and $Q^2=0$. Then:
\begin{itemize}
\item[(i)] We have the following decompositions, 
\[
\cH=\cH_d \oplus \cH_s \oplus \cH_p = \cH_s [\oplus] (\cH_p \oplus \cH_d)=\ker Q \oplus \cH_p, 
\]
where $\cH_s:=\ker Q \cap \ker Q^*$, $\cH_d:=\overline{\ran Q}$, $\cH_p:=\overline{\ran Q^*}$, $\ker Q=\cH_s\oplus \cH_d$ and $[\oplus]$ denotes Krein orthogonality.  Furthermore $\cH_s=J\cH_s$ and $\cH_p=J\cH_d$, and hence $\cH_s$ and $\cH_p\oplus \cH_d$ are Krein spaces, with respect to $J$.
\item[(ii)] The inner product $\iip{\cdot}{\cdot}_p$ makes $\cH^{BRST}_{phys}$ a Krein space. 
\end{itemize}

\end{theorem}
\begin{proof}
(i) Theorem \eqref{th:Hdsp1} gives the first Hilbert decomposition. 

Let $\psi \in \ker Q$ and $\xi \in D(Q)$, then by $Q^{\dag}=Q$ we have $\iip{\psi}{Q \xi}=\iip{Q\psi}{\xi}=0$, \ie $\ran Q [\perp] \ker Q$. Similarly $\ran Q^* [\perp] \ker Q^*$ and so $\cH_s=(\ker Q \cap \ker Q^*) [\perp] (\ran Q\oplus \ran Q^*)$. Now the second decomposition follows from the first.

From Appendix \eqref{ap:IIP} we have $Q^*=JQ^{\dag}J=JQJ$, hence $\psi \in \ker Q \cap \ker Q^*$ iff $ QJ^2\psi = Q^*J^2 \psi=0$ iff $ (JQJ)J\psi = (JQ^*J)J \psi=0$ iff $J\psi \in \ker Q \cap \ker Q^*$ hence $J\cH_s=\cH_s$. Also from $Q^*=JQ^{\dag}J=JQJ$ we have $\ran Q^*= \ran (JQJ)=J (\ran Q)J$ and it follows that $\overline{\ran Q^*}=J\overline{ \ran Q}$, \ie $\cH_d=J \cH_p$.

\pfit (ii): We have $J\cH_s=\cH_s$ from (i) and so (ii) follows from lemma \eqref{lm:phspksp} (i).
\end{proof}
\begin{rem}
\begin{itemize}
\item[(i)]Any bounded $Q$ satisfying the hypothesis of lemma \eqref{lm:ranQnll} can be extended to a Krein selfadjoint operator on $\cH$. Therefore Theorem \eqref{th:Hdsp} holds for bounded BRST charges and so $\cH^{BRST}_{phys}$ is a Krein space in the case of bounded charges.
\item[(ii)] When $Q$ is Krein-selfadjoint, Theorem \eqref{th:Hdsp} gives $J\cH_d=\cH_p$ and $J \cH_s =\cH_s$, hence $Q$ and $J$ have the following decompositions with respect to the $dsp$-decomposition:
\begin{equation*}
Q=\left(
\begin{matrix}
0 & 0 & M \\
0 & 0 & 0 \\
0 & 0 & 0
\end{matrix}
\right) 
\; \text{and} \;
J=\left(
\begin{matrix}
0 & 0 & L^{-1} \\
0 & W & 0 \\
L & 0 & 0
\end{matrix}
\right) 
\; \text{on} \;
\cH= 
\left(
\begin{matrix}
\cH_d\\
\cH_s\\
\cH_p
\end{matrix}
\right).
\end{equation*}
where $M\in \op(\cH_p,\cH_d)$ is a (possibly unbounded) closed operator, $W\in B(\cH_s)$ is a unitary, $L^{-1}=L^{*}$ and $L:\cH_d \to \cH_p$.

Note that the above decomposition for $J$ is not guaranteed in the case where $Q$ is not Krein-selfadjoint, as can be seen by example \eqref{ex:spdcount}.
\end{itemize}
\end{rem}

Now $\cH^{BRST}_{phys}$ is the physical state space and the inner product $\iip{\cdot}{\cdot}_p$ is usually taken to be the physical inner product. A physicality requirement is then that $\iip{\cdot}{\cdot}_p$ is positive definite, hence implying $\cH^{BRST}_{phys}$ is a Hilbert space. By lemma \eqref{eq:BRSTpHilip} (ii) this leads to:
\begin{definition}\label{df:spphcond}
The condition of physicality on $\cH^{BRST}_{phys}$ is:
\begin{equation}\label{eq:spphcond}
JP_s=P_s
\end{equation}
\end{definition}
In the heuristic picture a key object used to obtain $\cH_s$ is $\Delta =(QQ^*+Q^*Q)$. However it follows from (Ota \cite{Ota1984} Theorem 3.3 p232) that in the case of unbounded $Q$ we do not have $\ran Q \subset D(Q^*)$, so we need to take some care in defining $\Delta$. 

\begin{lemma}\label{lm:kDelHs}
Let $Q$ be as in Theorem \eqref{th:Hdsp} and let there be a space $D(\Delta)$ dense in $\cH$ and such that:
\[
D(\Delta) \subset (D(Q)\cap D(Q^*)), \qquad QD(\Delta) \subset D(\Delta) \supset Q^*D(\Delta),
\]
and,
\[
\cH_s=\ker Q \cap \ker Q^* \subset D(\Delta).
\]
Define,
\[
\Delta\psi:=(QQ^*+Q^*Q)\psi , \qquad \psi \in D(\Delta).
\]
Then $\Delta$ is $*$-symmetric on $D(\Delta)$ and $\ker \Delta=\cH_s$.
\end{lemma}
\begin{proof}
That $\Delta$ is $*$-symmetric is obvious. 

By definition of $\Delta$ we have $\ker Q \cap \ker Q^* =\cH_s \subset \ker \Delta$. Conversely, if $\psi \in \ker \Delta$ then $0=\ip{\psi}{\Delta \psi}=\norm{Q\psi}^2+\norm{Q^*\psi}^2$ and so $\ker \Delta \subset \cH_s$, thus $\ker \Delta = \cH_s$. 
\end{proof}
The assumptions above are of course realised in the case that $Q$ is bounded and will be realized in the example of \KOB QEM below (cf. lemma \eqref{lm:hdecDq}). 

An interesting question is, given an operator $Q$ as in Theorem \eqref{th:Hdsp1}, can we find a fundamental symmetry on $\cH$ such that $Q$ is Krein self adjoint? 
\begin{lemma} Let $Q$ be a closed operator on $\cH$ such that, $\ran Q \subset D(Q)$, and $Q^2=0$. Then there exists a fundamental symmetry $J\in B(\cH)$ (\ie  $J^2=\one$ and $J^*=J$) such that $Q$ is Krein-selfadjoint with respect to the indefinite inner product $\iip{\cdot}{\cdot}:=\ip{\cdot}{J\cdot}$ where $\ip{\cdot}{\cdot}$ is the original inner product. 
\end{lemma}
\begin{proof}
Let $\cH=\cH_d\oplus \cH_s \oplus \cH_p$ be the decomposition given by Theorem \eqref{th:Hdsp1}. Let $Q=W|Q|$ be the polar decomposition (\cite{Con1985} Theorem VIII 3.11 p242) for $Q$, \ie  $W$ is a partial isometry such that $W: ((\ker Q)^{\perp}=\cH_p) \to (\overline{\ran Q}=\cH_d)$ is isometric, and $W \ker Q=0$. Furthermore let $P_s$ be the orthogonal projection onto $\cH_s$, then,
\[
W^*W=P_p, \qquad WW^*=P_d, \qquad 0=W^2=(W^*)^2=P_sW^*=WP_s=W^*P_s.
\]
Now define $J:=P_{s}+W+W^*$, then using the above identities we get that $J^2=P_s^2+WW^*+W^*W=P_s+P_d+P_p=\one$. Obviously $J^*=J$.

Now $Q^*=|Q|W^*$ hence $|Q|=Q^*W$ hence $WQ^*W=W|Q|=Q$. Also, $Q^2=0$ implies $(Q^*)^2=0$ hence Theorem \eqref{th:Hdsp1} applied to $Q^*$ gives that $\ker Q^* = \cH_s \oplus \cH_p$. Therefore $Q^*J=Q^*(P_{s}+W+W^*)=Q^*W$. Furthermore $Q^*(D(Q^*))=\ran Q^* \subset \cH_p$ and so $JQ^*=(P_{s}+W+W^*)Q^*=WQ^*$. 

Therefore $JQ^*J=WQ^*W=Q$. Now by Appendix \eqref{ap:IIP} lemma \eqref{lm:khad} we have that $Q^{\dag}=JQ^*J=Q$ and so $Q$ is Krein-self adjoint. 
\end{proof}

By this lemma we see that the assumptions in Theorem \eqref{th:Hdsp1} are enough to ensure that there is a Krein structure that makes $Q$ Krein self-adjoint and gives all the structure in Theorem \eqref{th:Hdsp}.

\section{Operator Cohomology and $dsp$-decomposition}\label{sbs:AbcuQ}

The formal operator cohomology is defined as $\cP^{BRST}=\ker \drb/ \ran \drb$, but we need to consider domain issues when making sense of this rigorously. When we have done this we will show the connection between the $dsp$-decomposition and the operator cohomology. We will decompose $D(\drb)$ into blocks with respect to $\cH_d\oplus \cH_s \oplus \cH_p$, then use the block decomposition of $\ker \drb$ to  construct elements in $\ran \drb$. This is basically the same technique as in \cite{Hen1989} p285, \cite{Schf2001} p127, and K\&O \cite {KuOj79} Proposition 5.9 p68.

We describe here the basic mathematical structures used for the BRST-constraint method. In practice, the actual spaces and operators are constructed from the given field theory and constraints, where the number of ghosts equals the number of constraints.

To investigate the operator cohomology, we first have to define $\drb$ and the algebra it acts on explicitly. 
\begin{itemize}
\item[(i)] Let  $\cL$ be a Krein space with indefinite inner product $\iip{\cdot}{\cdot}_0$ with Krein involution $\dag$, and definite inner product $\ip{\cdot}{\cdot}_0=\iip{\cdot}{J_0\cdot}_0$ with Hilbert involution $*$, and let $\cD_0\subset \cL$ be a dense subspace of $\cL$.
\item[(ii)] Let $\cH_g$ be a Ghost space with ghost gradings and ghost number operator $\tilde{G}$ with domain $D(\tilde{G})$ as in Definition \eqref{df:ghsp}, \ie $\tilde{G}$ can be either the full ghost space of Berezin ghost space. In the case that $\cH_g$ is the full ghost space, we use the full ghost algebra $\cA_g$ with dense subalgebra $\cA_{g0}$. In the case of the Berezin ghost space we use the restricted ghost algebra $\rga$ and by an abuse of notation use $\cA_{g0}=\rga$ for this subsection only as the results it contains depend only upon the $\bZ_2$-grading of the Ghost algebra.
\item[(iii)] Let $\cH=\cL \otimes \cH_g$ with natural indefinite inner product $\iip{\cdot}{\cdot}$, positive definite inner product $\ip{\cdot}{\cdot}$ and  fundamental symmetry $J$ induced by  the indefinite inner products, positive definite inner products and  fundamental symmetries on $\cL$ and $\cH_g$.
\item[(iv)] Let $\cA_{0}$ be a Krein and Hilbert involutive unital subalgebra of $\op(\cL)$ defined on a common dense invariant domain $\cD_0\subset \cL$, \ie  for $A \in \cA_0$ we have $D(A)=\cD_0$, $A\cD_0 \subset \cD_0$. In addition we will assume that all $A\in \cA_0$ are closable on $\cD_0$. Let $\cA_{0} \otimes \cA_{g0}$ be  the algebraic tensor product where we are assuming no topology.
\item [(v)] The ghost grading of $\cA_{g0}$ naturally extends to $\cA_{0} \otimes \cA_{g0}$, however we will need to enlarge our algebra by operators not in the tensor product (such as the charge $Q$) and so we make the following definitions. Let $D(Q):=\cD_0 \otimes D(\tilde{G})$,  $G:=\one \otimes \tilde{G}$ with $D(G)=D(Q)$ and define:
\begin{align*}
(\cA_{0} \otimes \cA_{g0})^{+}:= & \;\{ A\in (\cA_{0} \otimes \cA_{g0}) \,|\, [ G,A]\psi=2kA\psi, \; \psi \in D(Q), k\in \mathbb{Z}\},\\
(\cA_{0} \otimes \cA_{g0})^{-}:= & \;\{ A\in (\cA_{0} \otimes \cA_{g0}) \,|\, [G,A]\psi=(2k+1)A\psi, \; \psi \in D(Q), k\in \mathbb{Z}\}
\end{align*}
ie, $(\cA_{0} \otimes \cA_{g0})^{+}$ are the elements with difference between ghosts and conjugate ghosts is even, and similarly $(\cA_{0} \otimes \cA_{g0})^{-}$ corresponds to elements with an odd difference. As $\gh(f)^*=\cgh(f)$, $\gh(f)^{\dag}=\gh(f)$, $\cgh(g)^{\dag}=\cgh(g)$ for $f \in \cH_2$, $g\in \cH_1$ it follows that $(\cA_{0} \otimes \cA_{g0})^{+}$ and $(\cA_{0} \otimes \cA_{g0})^{-}$ are both $*$-closed and $\dag$-closed subspaces of $\cA_{0} \otimes \cA_{g0}$, moreover $(\cA_{0} \otimes \cA_{g0})^{+}$ is a $*$-closed and $\dag$-closed subalgebra of $\cA_{0} \otimes \cA_{g0}$. Define a $\mathbb{Z}_2$-grading of $\cA_0\otimes \cA_{g0}$ by,
\[
\cA_0\otimes \cA_{g0}:=(\cA_{0} \otimes \cA_{g0})^{+}\oplus (\cA_{0} \otimes \cA_{g0})^{-}
\]
and define the grading automorphism on $\cA_{0} \otimes \cA_{g0}$ as:
\[
\gamma(A_{+}+A_{-})=A_{+}-A_{-},
\]
where $A_{+}\in (\cA_{0} \otimes \cA_{g0})^{+}, A_{-}\in (\cA_{0} \otimes \cA_{g0})^{-}$.
\item[(vii)] Now let $\cA$ be a $\dag$-involutive subalgebra of $\cA_{0} \otimes \cA_{g0}$ such that $\gamma(\cA)=\cA$. As $\gamma(\cA)=\cA$ we get that,
\[
\cA=\cA^{+}\oplus \cA^{-},
\]
where $\cA^{+}:=\cA\cap (\cA_{0} \otimes \cA_{g0})^{+}$ and $\cA^{-}:=\cA\cap (\cA_{0} \otimes \cA_{g0})^{-}$. We need to include the case of proper inclusion $\cA \subset \cA_{0} \otimes \cA_{g0}$ to handle the $C^*$-algebraic BRST-QEM below where, for technical reasons, the BRST superderivation cannot be defined on all of $\cA_{0} \otimes \cA_{g0}$ (cf. Definition \eqref{df:sddom1} and Definition \eqref{df:sddom2} below).
\end{itemize}
\begin{lemma}\label{lm:QessaCsa}
Let $Q\in \op(Q)$ with domain $D(Q)$ be such that $Q:D(Q)\to D(Q)$ be a Krein symmetric and 2-nilpotent operator, \ie $Q^2\psi=0$ for all $\psi \in D(Q)$. Then:
\begin{itemize}
\item[(i)] $Q$ is closable.
\item[(ii)] $\ran \cQ \subset D(\cQ)$ and $\cQ^2\psi =0$ for all $\psi \in D(\cQ)$. 
\item[(iii)] $(Q^*)^2\psi=0$ for all $\psi \in D(Q^*)$  
\end{itemize}
\end{lemma}
\begin{proof}
(i): This follows from Proposition \eqref{pr:Krclos}.
\smallskip 
\noindent(ii): Let $\xi =\cQ \psi$ for $\psi \in D(\cQ)$. As $\cQ$ is the closure of $Q$, we have a sequence $(\xi_n) \subset \ran Q \subset D(Q)$ such that $\xi_n \to \xi$. By the assumptions on $Q$ we have that $Q\xi_n=0$ and so if we denote an element in the graph of $\cQ$ as $(x,y)$, then
\[
(\xi_n,Q\xi_n)=(\xi_n,0)\to (\xi,0).
\]
As the graph of $\cQ$ is closed we get that $\xi \in D(\cQ)$, \ie $\ran \cQ \subset D(\cQ)$ and $0=\cQ\xi=\cQ^2\psi$, \ie $\cQ^2\psi=0$ for all $\psi \in D(\cQ)$.
 
\smallskip
\noindent(iii): 
First recall that $Q^*=\cQ^*$. Let $\psi \in D(Q^*)=D(\cQ^*)$ and $\xi \in D(\cQ)$. Then by (ii) $0=\ip{\psi}{\cQ^2\xi}=\ip{Q^*\psi}{\cQ\xi}$. Hence it follows that $Q^* \psi \in D(Q^*)$ and hence that $(Q^*)^2=0$.
\end{proof}
By this lemma, $\cQ$ satisfies the conditions of Theorem \eqref{th:Hdsp1} and so we have the $dsp$-decomposition for $\cQ$:
\[
\cH=\cH_d \oplus \cH_s \oplus \cH_p= \overline{\ran \cQ} \oplus \cH_s  \oplus \overline {\ran \cQ^*}.
\]
Let $P_i$ be the projection onto $\cH_i$ for $i=d,s,p$. As $\bch \subset \cQ$ there is no guarantee that these projections preserve $D(Q)$. To decompose the $\cA$ into a convenient form we assume:
\begin{ass}\label{ass:DQdec} Let $Q$ be as in lemma \eqref{lm:QessaCsa} and let:
\begin{itemize} 
\item[(i)] $P_iD(Q) \subset D(Q)$ for $i=p,s,d$, and $Q^*D(Q)\subset D(Q)$. Then the following are well defined:
\begin{itemize}
\item[(a)] $M:=Q|_{P_pD(Q)}$.  
\item[(b)] As $\ker M =\{0\}$ there exists an inverse $M^{-1}:\ran Q \to (\cH_p\cap D(Q))$. Let $K:=M^{-1}P_d$, with domain $D(K)=\{\psi\in D(Q)\,|\, P_d\psi \in \ran Q\}$.
\end{itemize}
\item[(ii)] $P_d \ran Q \subset \ran Q$. 
\item[(iii)] Let:
\begin{gather}\label{eq:progncom}
 [G,Q]\psi=Q\psi,\quad[G,Q^*]\psi=-Q^*\psi, \quad [G,K]\psi=-K\psi \\ 
 [G,P_d]\psi=[G,P_p]\psi=[G,P_s]\psi=0, \notag
\end{gather}
for $\psi \in D(Q)$.
\end{itemize}
\end{ass}
Note that $D(KQ)=\{\psi \in D(Q)\,|\, Q\psi \in D(K)\}=\{\psi \in D(Q)\,|\, P_dQ\psi \in \ran Q\}$, and so by Assumption \eqref{ass:DQdec} (ii) we get $D(KQ)=D(Q)$ and hence $KQ=P_p|_{D(Q)}$.

Assumption (iii) will allow us to enlarge $\cA$ by $Q,Q^*,K,P_s,P_p,P_d$, e.g. below in equation \eqref{eq:Aextdef}. In examples, the assumptions need verification, and we now give conditions on $D(Q)$ for these to be true. These conditions will hold in the QEM example below and the next lemma will show that as well as in the case of bounded $Q$ and $G$ (\ie finitely many bounded constraints).

\begin{lemma}\label{lm:projdspinc}
Let $Q \in \op(\cH)$ preserve its domain $D(Q)$, be Krein symmetric and let $Q^2\psi=0$ for all $\psi \in D(Q)$. Assume that: 
\[
D(Q)= \ran Q \oplus \cH_s\oplus \ran Q^*,\quad \text{and} \quad[G,Q]\psi=Q\psi. \;\forall \psi \in D(Q).
\]
Then, 
\begin{itemize}
\item[(i)]$P_iD(Q) \subset D(Q)$ for $i=p,s,d$, and $Q^*D(Q)\subset D(Q)$,
\item[(ii)] $P_d \ran Q \subset \ran Q$. 
\item[(iii)]We have that:
\begin{gather*}
 [G,Q^*]\psi=-Q^*\psi,\quad [G,K]\psi=-K\psi, \\
 [G,P_d]\psi=[G,P_p]\psi=[G,P_s]\psi=0,\notag
\end{gather*}
for $\psi \in D(Q)$. 
\end{itemize}
\end{lemma}
\begin{proof}
\noindent(i): $Q$ is closable by lemma \eqref{lm:QessaCsa}. Hence $Q^*=\cQ^*$ by \cite{ReSi1972v1} Theorem VIII.1 (c) p253, and so $\ran Q^* \subset \overline{\ran \cQ^*}= \cH_d$, moreover $P_p|_{\ran Q^*}=\one$. As $\ran Q \subset \ran \cQ$, $P_iD(Q) \subset D(Q)$ for $i=p,s,d$ now follows. Furthermore $Q^*\cH_s=\{0\} \subset D(Q)$, $Q^*\ran Q \subset \ran Q^* \subset D(Q)$ and $Q^*\ran Q^*=\{0\}\subset D(Q)$ by lemma \eqref{lm:QessaCsa} (iii), hence $Q^*D(Q) \subset D(Q)$ by the assumed decomposition of $D(Q)$.

\smallskip \noindent(ii): $\ran Q \subset \ran \cQ \subset \cH_d$.

\smallskip \noindent(iii): Let $\psi\in D(Q)$. It follows from $[G,Q]\psi=Q\psi$ and $*$-adjoints that $[G,Q^*]\psi=-Q^*\psi$. Hence we have $[G,Q^*Q]\psi=0$.

It is sufficient to prove the remaining identities for $\psi \in D(Q)$ such that $G\psi=g\psi$, $g \in \mathbb{R}$ if we recall that $D(Q)=\cD_0\otimes D(\tilde{G})$ is spanned by linear combinations of such vectors where $g$ runs over all the possible ghost numbers we are considering (cf. Proposition \eqref{pr:ghgrad} and Proposition \eqref{pr:bsghgrad} (ii)). So for the remainder of the proof we let the chosen $\psi \in D(Q)$ be such that $G\psi=g\psi$, $g \in \mathbb{R}$ and so by the assumption on $D(Q)$ we have $\psi=P_s\psi +Q\vp +Q^* \xi$ where $\vp \in \ran Q^* \perp \ker Q$ and $\xi \in \ran Q \perp \ker Q^*$.

We want to show that $[G,P_d]\psi=0$ for all $\psi \in D(Q)$. 
By (i) we have $P_dG\psi=gP_d\psi=gQ\vp$ for $\vp \in \ran Q \perp \ker Q^*$ as above. So we need now only show $GP_d\psi=GQ\vp=gQ\vp$. To do this we calculate,
\[
Q^*QG\vp=GQ^*Q\vp=GQ^*\psi=Q^*(G-1)\psi=(g-1)Q^*\psi=Q^*Q(g-1)\vp.
\]
where the first inequality comes from $[G,Q^*Q]=0$ on $D(Q)$ using $QD(Q)\subset D(Q)$ and $[G,Q^*]=-Q^*$ and its adjoint on $D(Q)$. Therefore:
\[
Q^*Q[(G-(g-1))\vp]=0.
\]
Now $\ran Q \perp \ker Q^*$ and $\ran Q^*\perp \ker Q$ imply that $\ker (Q^*Q)|_{\ran Q^*}=\{0\}$. Also $[G,Q^*]=-Q^*$ on $D(Q)$ implies that $G\ran Q^* \subset \ran Q^*$. Hence as $\vp, G\vp \in \ran Q^*$ we have that $Q^*Q[(G-(g-1))\vp]=0$ implies that $G\vp=(g-1)\vp$. By $[G,Q]=Q$ on $D(Q)$ we have $GQ\vp=gQ\vp$, and so $[G,P_d]\psi=0$. 

Now by Theorem \eqref{th:Hdsp} we have $\cH_p=J\cH_d$ which implies $P_p=JP_dJ$ which implies $P_p^{\dag}=JP_d^{*}J=P_p$. By Proposition \eqref{pr:ghsp} (or Proposition \eqref{pr:bsghgrad} (i) for the Berezin ghost number operator) we have  $G^{\dag}=-G$ on $D(Q)$, hence $[G,P_p]\psi=[G,P_d]^{\dag}\psi=0$. 

We have shown that $[G,P_d]\psi=[G,P_p]\psi=0$ and since $P_s=\one-P_p-P_d$ we have that $[G,P_s]\psi=0$.

Lastly, $GK\psi= GM^{-1}Q\vp=G\vp=(g-1)\vp$ where the last equality was proven above. Therefore $[G,K]\psi=-\vp=-K\psi$.
\end{proof}

\begin{rem}\begin{itemize}
\item[(i)] For unbounded $Q$ the above condition that $D(Q)= \ran Q \oplus \cH_s\oplus \ran Q^*$ needs verification, and will be done for the forthcoming BRST-QEM example by lemma \eqref{lm:hdecDq}. We can always restrict the domain of $Q$ so that the above condition holds, however we then would have to check that $\cA$ preserves the restricted domain if we use $Q$ to generate $\drb$. 
\item[(ii)] For the case of bounded $Q$ and finitely many bounded constraints (\ie bounded $G$) we can restrict the domain of $Q$ to $D(Q)=\ran Q \oplus \cH_s \oplus \ran Q^*$, then extend \eqref{eq:progncom} to all of $\cH$ by continuity.
\end{itemize}
\end{rem}
With the above assumptions we enlarge $\cA$  to 
\begin{equation}\label{eq:Aextdef}
\cAe:= \salg{\cA \cup\{Q,K,P_p,P_d, P_s\}},
\end{equation}
and using the commutation relations \eqref{eq:progncom} we can extend the $\mathbb{Z}_2$ grading to $\cAe$ which is a $*$-algebra but not a $\dag$-algebra. Note that by  assumptions, $D(Q)$ is an invariant domain for all the elements in $\cAe$. By a slight abuse of notation we use the same notation $\cA^{-}$ and $\cA^{+}$ for the extended even and odd algebras. Note that $Q,Q^*,K \in \cA^{-}$ and $P_s ,P_p,P_d\in \cA^{+}$. We define the BRST superderivation as:

\begin{definition} With definitions as above, the BRST superderivation is
\[
\drb(A)\psi:=\sbr{Q}{A}\psi =(QA-\gamma(A)Q)\psi, \qquad A\in \cAe, \psi \in D(Q).
\]
\end{definition}

An interesting result is that we have a unique $Q$ which generates $\drb$.
\begin{lemma}\label{lm:unch}
Let $Q_1,Q_2 \in \cA^{-}$ have common dense invariant domain $D(Q)$ (as they are in $\cAe$), and be such that $Q_1^*,Q_2^*\in \cA^{-}$ and $\drb(A)\psi:=\sbr{Q_1}{A}\psi=\sbr{Q_2}{A}\psi$ for all $A\in \cAe$ and $\psi \in D(Q)$. Then $Q_1=Q_2$. 
\end{lemma}
\begin{proof}
Let $\psi \in D(Q)$, and calculate
\begin{align*}
[(Q_1-Q_2)^*(Q_1&-Q_2)+(Q_1-Q_2)(Q_1-Q_2)^*]\psi= \\
&=[\{Q_1,Q_1^*\}+\{Q_2,Q_2^*\}-\{Q_1,Q_2^*\}-\{Q_2,Q_1^*\}]\psi,\\
&=[\drb(Q_1^*)+\drb(Q_2^*)-\drb(Q_1^*)-\drb(Q_2^*)]\psi,\\
&=0
\end{align*}
From this we have that $\norm{(Q_1-Q_2)\psi}=\norm{(Q_1-Q_2)^*\psi}=0$, \ie $Q_1=Q_2$.
\end{proof}
This result is interesting as it relies on $\drb$ being a \emph{super}derviation, in contrast to the case of derivations where the generators are non-unique up to a central term. 

Since the original algebra $\cA$ contains the physical information, we will define the BRST-observable algebra using the restriction of $\drb$ to $\cA$.
\begin{proposition}\label{pr:BRSTopal}
Define the BRST-physical observable algebra as:
\begin{equation*}
\cP^{BRST}:=(\ker \drb \cap \cA) /(\ran \drb \cap \cA).
\end{equation*}
Let $\vp:\ker \cQ \to \cH^{BRST}_{phys}=\ker \cQ /\cH_d$ and $\tau: \ker \drb  \to \ker \drb /\ran \drb$ be the factor maps and note that $\tau: (\ker \drb \cap \cA) \to \cP^{BRST}=(\ker \drb \cap \cA) /(\ran \drb \cap \cA)$. Define the dense subspace $\cD^{BRST}_{phys}:=\vp(\ker Q)$ and for ease of notation,
\begin{equation}\label{eq:hfacmap}
\hat{\psi}:=\vp(\psi),\, \forall \psi \in \ker Q, \qquad \hat{A}:=\tau(A), \, \forall A \in \ker \drb.
\end{equation}
Then:
\begin{itemize}
\item[(i)] Both $\ker \drb / \ran \drb$ and $\cP^{BRST}$ have a natural actions on $\cD^{BRST}_{phys}$ by,
\begin{equation}\label{eq:brphacsbsp}
\hat{A}\hat{\psi}=\widehat{A\psi}, 
\end{equation}
for all $A \in \ker \drb$ respectively $A \in (\ker \drb \cap \cA)$ and all $\psi \in \ker Q$.
\item[(ii)] The $\dag$-involution  on $\ker \drb$ and $\ker \drb \cap \cA$ factors to the $\dag$-involution on $\ker \drb/ \ran \drb$ respectively $\cP^{BRST}$ with respect the action defined in equation \eqref{eq:brphacsbsp}, \ie
\[
\widehat{A^{\dag}}\hat{\psi}=\hat{A}^\dag\hat{\psi}
\]
for all $A \in \ker \drb$ respectively $\ker \drb \cap \cA$ where $\hat{A}^{\dag}$ is the the adjoint of $\hat{A}$ with respect to the inner product $\iip{\cdot}{\cdot}_p$.
\end{itemize}
\end{proposition}
\begin{proof}
(i): We check that the action given by equation \eqref{eq:brphacsbsp} is well defined. Now $\ker \drb$ preserves $\ker Q$ and $\ran Q$, hence factors to $\cD^{BRST}_{phys}$. Furthermore $(\ran \drb) \ker Q \subset \ran Q$ and so factors trivially to $\cD^{BRST}_{phys}$ hence the action is well defined.

\smallskip \noindent (ii): First $Q^{\dag}=Q$ implies that $\ker \drb$ and $\ran \drb$ are $\dag$-algebras, hence $\dag$ factors as an involution on $\cP^{BRST}$. 
Moreover it factors to the involution with respect to the indefinite inner product on $\cH^{BRST}_{phys}$ by equation \eqref{eq:brphacsbsp} which can be seen by calculation:
\[
\iip{\widehat{A^{\dag}}\hat{\psi}}{\hat{\xi}}_p=\iip{\widehat{A^{\dag}\psi}}{\hat{\xi}}_p=\iip{A^{\dag}\psi}{\xi}=\iip{\psi}{A\xi}=\iip{\hat{\psi}}{\hat{A}\hat{\xi}}_p,
\]
for all $\psi,\xi \in \cD^{BRST}_{phys}$ and all $A \in \ker \drb$, where we have used that by definition $\iip{\hat{\psi}}{\hat{\xi}}_p=\iip{\psi}{\xi}$ for all $\psi,\xi \in \ker Q$ (cf. Definition \eqref{df:opbrsphsp}). 
\end{proof}
The fact that $\dag$ factors to an involution on $\cP^{BRST}$ in Proposition \eqref{pr:BRSTopal} (ii) relies on $Q^{\dag}=Q$. But $Q\neq Q^*$ implies that $\ker \drb$ and $\ran \drb$ are \emph{not} $*$-algebras in general, hence the $*$-involution on $\ker \drb \cap \cA$ does not factor to $\cP^{BRST}$. This is problematic from the mathematical standpoint as $*$-algebras have more tools available e.g. $C^*$-algebra theory in their bounded representations. However as $\cH^{BRST}_{phys}$ is a Krein space it has a Hilbert inner product which we can use to define a $*$-involution on $\ker \drb / \ran \drb$. 

\begin{lemma}\label{lm:facsal} Given notation above:
\begin{itemize}
\item[(i)]  For $A \in \ker \drb$ we have:
\begin{equation*}
\hat{A}^*\hat{\psi} =\widehat{P_sA^*P_s\psi}=\widehat{P_sA^*P_s}\hat{\psi}, \qquad \forall \hat{\psi} \in \cD^{BRST}_{phys}.
\end{equation*}
and 
\begin{equation}\label{eq:keranhilad}
 \hat{A}^*|_{\cD^{BRST}_{phys}} \in (\ker \drb /\ran \drb).
\end{equation}
Hence $\hat{A}\to \hat{A}^*:= \widehat{P_sA^*P_s}$ is an involution on $\ker \drb /\ran \drb$ which coincides with the $*$-involution on $\ker \drb /\ran \drb$ with respect to the indefinite inner product $\iip{\cdot}{\cdot}_p$. With respect to this involution, $\ker \drb /\ran \drb$ is a $*$-algebra.
\item[(ii)] If the physicality condition $JP_s=P_s$ holds (cf. equation \eqref{eq:spphcond}) then the $\dag$-adjoint on $\ker \drb $ factors to the $*$-adjoint on $\ker \drb /\ran \drb$ as defined in (i). 
\end{itemize}
\end{lemma}
\begin{proof}
(i):  Let $A \in \ker \drb$ and ${\psi},{\xi} \in \ker Q$. Using equation \eqref{eq:BRSTpHilip} we have,
\[
\ip{\hat{\xi}}{\hat{A}\hat{\psi}}_p=\ip{P_s\xi}{P_sA\psi}=\ip{P_s\xi}{P_sAP_s\psi}=\ip{P_sA^*P_s\xi}{\psi},
\]
where we used the fact that $A \in \ker \drb$ implies that $A \ran Q \subset \ran Q$ in the second last inequality. It is easy to check that $P_s \cAe P_s \subset \ker \drb$ hence $\widehat{P_sA^*P_s}\in (\ker \drb/\ran \drb)$ and we have that,
\begin{equation*}
\hat{A}^*\hat{\psi} =\widehat{P_sA^*P_s\psi}=\widehat{P_sA^*P_s}\hat{\psi}, \qquad \forall \hat{\psi} \in \cD^{BRST}_{phys}.
\end{equation*}
As $P_s$ is $*$-self adjoint and $\cAe$ is a $*$-algebra we get that $\ker \drb /\ran \drb$ is a $*$-algebra.

\smallskip \noindent (ii): Let $JP_s=P_s$. Then lemma \eqref{lm:phspksp} (iii)  gives that $\iip{\hat{\psi}}{\hat{\xi}}_p=\ip{\hat{\psi}}{\hat{\xi}}_p$ and it follows that $\widehat{A^{\dag}}\hat{\psi}=\hat{A}^{\dag}\hat{\psi}=\hat{A}^*\hat{\psi}$ for all $\hat{\psi} \in \cD^{BRST}_{phys}$.
\end{proof}

\begin{rem} \begin{itemize} \item[(i)]
If we use $\cA_{ext}$ as the domain of $\drb$ then we have by lemma \eqref{lm:facsal} (i) that $\ker \drb /\ran \drb$ is a $*$-algebra. So although $\ker \drb$ was not necessarily a $*$-algebra, factoring out $\ran \drb$ from $\ker \drb$ produced a $*$-algebra. Notice that we had to extend to $\cAe$ to get this result in general, \ie we do not have that $(\ker \drb \cap \cA)/(\ran \drb \cap \cA)$ is a $*$-algebra in general and we would need to put more restrictions on the domain of $\drb$. 

We will return to this point in the $C^*$-theory (Proposition \eqref{pr:repban}) as in the abstract theory we want the  BRST physical algebra to be a $C^*$-algebra with respect the $*$-involution and norm defined in Proposition \eqref{pr:bddcshinv}. We will find that a sufficient condition is for the physicality condition $JP_s=P_s$ to hold (cf. Definition \eqref{df:spphcond} and Example \eqref{ex:pscs} (i)). In this case we also have that the $\dag$-involution factors to the $*$-involution on the final physical algebra.
\item[(ii)] When $Q$ is Krein-selfadjoint, Theorem \eqref{th:Hdsp} gives $J\cH_d=\cH_p$ and so the physicality condition is $JP_s=P_s$, hence in terms of the \emph{dsp}-decomposition, the physicality condition states that:
\begin{equation*}
J=\left(
\begin{matrix}
0 & 0 & L^{-1} \\
0 & \one & 0 \\
L & 0 & 0
\end{matrix}
\right) 
\; \text{on} \;
\cH= 
\left(
\begin{matrix}
\cH_d\\
\cH_s\\
\cH_p
\end{matrix}
\right).
\end{equation*}
where $L^{-1}=L^{*}$, $L:\cH_d \to \cH_p$.
\end{itemize}
\end{rem}

We want to decompose $\ker \drb$ into a convenient form. To do this we introduce the following algebra homomorphism.
\begin{lemma}\label{lm:alkerhom}
Define,
\begin{align*}
\Phi^{ext}_s:\cA_{ext} \to \cA_{ext}, \qquad \text{by} \qquad \Phi^{ext}_s(A):=P_sAP_s, 
\end{align*}
Then $\Phi^{ext}_s$ is an algebra homomorphism on $\ker \drb$.
\end{lemma}
\begin{proof}
That $\Phi^{ext}_s$ For $A \in \ker \drb$ we have $A\cH_s \subset \cH_s \oplus \cH_d$ and $A \cH_d \subset \cH_d$ and so for $A, B \in \ker \drb$ we have
\[
\Phi^{ext}_s(AB)=P_sA(P_s+P_d)BP_s=P_sAP_s^2BP_s=\Phi^{ext}_s(A)\Phi^{ext}_s(B),
\]
and so it follows that $\Phi^{ext}_s$ is an algebra homomorphism on $\ker \drb$.
\end{proof}
 
The next theorem shows that the kernel of $\Phi^{ext}_s$ is $\ran \drb$. 
\begin{theorem}\label{pr:krdel}
Let $Q$, $D(Q)$ satisfy Assumptions \eqref{ass:DQdec} and let $\cAe$,  $\Phi^{ext}_s$ be as in equation \eqref{eq:Aextdef} and lemma \eqref{lm:alkerhom}. In particular $Q\in \cA^{-}$ is a 2-nilpotent Krein-symmetric (hence closable by Proposition \eqref{pr:Krclos}) operator with domain $D(Q)$, and  $\cQ$ satisfies the hypothesis of Theorem \eqref{th:Hdsp1}. Let $\drb(A)=[Q,A]_s$ for $A \in \cAe$. If $B \in \ker \delta$ and $P_d \ran B \subset \ran Q$, then:
\begin{itemize}
\item[(i)]  $B \in \ran \delta$ iff $\Phi^{ext}_s(B)=0$, \ie $B \in \ran \delta$ iff $(P_sBP_s)D(Q)=0$. This implies that:
\[
B=\Phi^{ext}_s(B) + C,
\]
where $C\in \ran \drb$.
\item[(ii)] Suppose that $P_d \ran A \subset \ran Q$ for all $A \in \cA$. Then:
\[
\cP^{BRST}=(\ker \drb \cap \cA)/ (\ran \drb \cap \cA)\cong  \Phi^{ext}_s(\ker \drb \cap \cA).
\]
where $\cong$ denotes an algebra isomorphism.
\end{itemize}
\end{theorem}
\begin{proof}
\noindent(i): Given that $P_i D(Q)\subset D(Q)$ for $i=s,p,d$, where $P_i$ is the projection onto $\cH_i$, then with respect to the $dsp$-decomposition of $\cQ$, $Q$ has the representation:
\begin{equation*}
\left(
\begin{matrix}
0 & 0 & M \\
0 & 0 & 0 \\
0 & 0 & 0
\end{matrix}
\right) 
\; \text{on} \;
\cH= 
\left(
\begin{matrix}
\cH_d\\
\cH_s\\
\cH_p
\end{matrix}
\right).
\end{equation*}
Let $B\in \cAe$. Then it has the representation: 
\begin{equation*}
B=\left(
\begin{matrix}
B_{11} & B_{12} & B_{13} \\
B_{21} & B_{22} & B_{23}\\
B_{31} & B_{32} & B_{33}
\end{matrix}
\right)
\; \text{on} \;
\cH= 
\left(
\begin{matrix}
\cH_d\\
\cH_s\\
\cH_p
\end{matrix}
\right)
\end{equation*}
Now as $P_s,P_p,P_d \in \cAe$ we find:
\begin{equation}\label{eq:ubderd}
\delta(B)=
\left(
\begin{matrix}
MB_{31} & MB_{32} & MB_{33}-\gamma(B_{11})M \\
0 & 0 & -\gamma(B_{21})M\\
0 & 0 & -\gamma(B_{31})M
\end{matrix}
\right).
\end{equation}
Therefore if $B \in \ker \delta$, then $MB_{31} = MB_{32} =0$ which implies $B_{31}=B_{32}=0$ as $\ker M=0$. As $\ran M= \ran Q$, $\gamma(B_{21})M=0$ implies $\gamma(B_{21})\ran Q=0$. Also $\ran Q$ is dense in $P_dD(Q)= (D(Q) \cap \overline{\ran \cQ})$. Let $\xi \in P_dD(Q)$ and a sequence $(Q\psi_n)_{n\in \mathbb{Z}}$ such that $ Q\psi_n \to \xi$. As $\cAe$ is $*$-involutive and $\gamma(B_{21})M=0$, we get for any $\psi \in D(Q)$ that: 
\[\ip{\psi}{ \gamma(B_{21})\xi}=\lim_n\ip{\gamma(B)^*P_s\psi}{Q\psi_n}=\lim_n\ip{\psi}{\gamma(B_{21})Q\psi_n}=0,
\] 
where we used the fact that  $\gamma(B_{21})=\gamma(P_s)\gamma(B)\gamma(P_d)=P_s\gamma(B)P_d$ and that $\gamma$ is a $*$-automorphism in the first equality. As $D(Q)$ is dense in $\cH$, we have $\gamma(B_{21})D(Q)=\gamma(B_{21})P_dD(Q)=0$. As $\gamma$ is an automorphism,  $B_{21}=0$.

Now $MB_{33}-\gamma(B_{11})M=0$ implies $B_{11} \ran Q \subset \ran Q$ and $B_{33}= M^{-1}\gamma (B_{11})M $. Therefore,
\begin{equation} \label{eq:ubkerd}
B\in \ker \delta \quad \text{iff} \quad B=\left(
\begin{matrix}
B_{11} & B_{12} & B_{13} \\
0 & B_{22} & B_{23}\\
0 & 0 & M^{-1}\gamma(B_{11})M
\end{matrix}
\right),
\end{equation}

Now suppose $B\in \ran \delta \subset \ker \delta$. If we compare equations \eqref{eq:ubderd} and \eqref{eq:ubkerd} we see that $B_{22}=0$. Conversely suppose that
that $B_{22}=0$. Then we would like to choose a $C \in \cAe$ to produce $\drb(B)$. Consider:
\begin{equation}\label{eq:Cran}
C=\left(
\begin{matrix}
0 & 0 & 0 \\
-\gamma( B_{23}M^{-1}) & 0 & 0\\
M^{-1}B_{11} & M^{-1}B_{12} & \quad M^{-1}B_{13}
\end{matrix}
\right)=KB-P_s\gamma(B)K \in \cAe
\end{equation}
As $B\in \cAe$ preserves $D(Q)$, and by the assumption $P_d\ran B \subset \ran Q$ we see that $\ran(P_dB)=P_d\ran(B)\subset \ran Q= D(M^{-1})$ and it follows that the bottom row of $C$ is well defined. Also $\ran M^{-1} \subset D(Q)$ and so $-\gamma( B_{23})M^{-1}$ is well defined. Hence all the entries of $C$ are well defined. 

Assumption \eqref{ass:DQdec} (iii) gives that $\gamma(K)=-K$ and $\gamma(P_d)=P_d$, so as $B_{11}=P_dB_{11}$ we have 
\[
\gamma(M^{-1}B_{11})=\gamma(M^{-1}P_dB_{11})=\gamma(KB_{11})=-K\gamma(B_{11})=-M^{-1}\gamma(P_d)\gamma(B_{11})=-M^{-1}\gamma(B_{11}).
\]
Substituting $C$ into equation \eqref{eq:ubderd} gives,
\begin{equation}
\delta(C)=
\left(
\begin{matrix}
B_{11} & B_{12} & B_{13} \\
0 & 0 & \gamma(\gamma( B_{23}M^{-1}))M\\
0 & 0 & -\gamma(M^{-1}B_{11})M
\end{matrix}
\right)=
\left(
\begin{matrix}
B_{11} & B_{12} & B_{13} \\
0 & 0 &  B_{23} \\
0 & 0 &  M^{-1}\gamma(B_{11})M
\end{matrix}
\right),
\end{equation}
where we used $\gamma^2=\iota$ and $\gamma(M^{-1}B_{11})=-M^{-1}\gamma(B_{11})$ in the second equality. Hence we can choose $C$ as above to get $\delta(C)=B$. 

\medskip \noindent (ii): Let $B \in \ker \drb$. As $P_s^2=P_s$ we have $(\Phi^{ext}_s)^2(B)=\Phi^{ext}_s(B)$, hence $C=(B-\Phi^{ext}_s(B))\in (\ker \Phi^{ext}_s(B)\cap \ker \drb)$. Therefore by part (i) we have $C\in \ran \drb$ and so part (ii) follows.

\medskip \noindent (iii): Let $A \in (\ker \drb \cap \cA)$ and  $\Phi^{ext}_s(A)=0$. Then by part (i) $A \in (\ran \drb \cap \cA)$. Conversely, let $A=(QC -\gamma(C)Q) \in (\ran \drb \cap \cA)$ for some $C \in \cAe$. As $P_sQ=QP_s=0$ we have that $\Phi^{ext}_s(A)=0$ hence $(\ran \drb \cap \cA)= (\ker \Phi^{ext}_s \cap \ker \drb \cap \cA)$. By lemma \eqref{lm:alkerhom}, $\Phi^{ext}_s$ is an algebra homomorphism on $\ker \drb$, and so part (iii) follows.
\end{proof}
We have the important corollary:
\begin{corollary}\label{cr:oneinker}
Given the hypothesis of Theorem \eqref{pr:krdel} then:
\[
\one \in \ran \delta\qquad \text{iff} \qquad \cH_s =\{0\}
\]
\end{corollary}
\begin{rem}\label{rm:ext}
\begin{itemize}\item[(i)] Note that while the definition $\cP^{BRST}$ involves only the unextended algebra $\cA$, to prove the isomorphism  $\cP^{BRST}\cong \Phi^{ext}_{s}(\ker \drb \cap \cA)$ needed the extension to $\cAe$. This is unsatisfactory as it requires us add extra elements to our algebra. In the $C^*$-algebraic approach it desirable to do this as little as possible as it effects the representation theory of the algebras involved. If we want to take the minimal extension of $\cA$ for the same construction to work in Theorem \eqref{pr:krdel}, we see from equation \eqref{eq:Cran}, that we still need extend to $\alg{Q,K, P_s, \cA}$.

Note however that $\ker \delta/ \ran \delta$  can be obtained algebraically in specific examples (cf. \cite{HoVo92}) without extending $\cA$, but the price paid is much more intensive calculations which the \emph{dsp}-decomposition construction avoids (cf. \cite{HoVo92} p1322-1326).
\item[(ii)] Motivated by the last section, a natural direction to investigate is to see if there is an `operator \emph{dsp}-decomposition'.

This is the approach taken in \cite{HoVo92} where the authors examine the structures as defined in the BRST-QEM example. They define a new superderivation $\drb^*(A):=[Q^*,A]_s$ and motivated by $\Delta$ define $d:=\{\drb, \drb^*\}$. It is straightforward to calculate from superbracket properties (Appendix\eqref{ap:SS}) that $d$ is a \emph{derivation} and that $d(A)=[\Delta,A]$. Unfortunately it turns out that $\ker d + \ran \drb \subset \ker \drb$, where the containment is proper (\cite{HoVo92} p1321). So $d$ is used  to calculate a large part of $\ran \drb$, the remaining part being calculated directly using the properties of $Q,Q^*,\Delta,G$ specific to the BRST-QEM setup. While this solves the problem, the last calculation involved is not straightforward, and it is not clear how to extend to the Hamiltonian BRST with general constraints. The calculations in \cite{HoVo92} are algebraic and it is conjectured at the end of \cite{HoVo92} that using the $dsp$-decomposition of $\cH$, the calculations could be done in a more economic way. We answer this in the affirmative using Theorem \eqref{pr:krdel} below (cf. Proposition \eqref{pr:fdkob} (iii)).
\end{itemize}
\end{rem}

\subsection{Example - Abelian Hamiltonian Constraints}\label{sbs:egAHC}
We have seen for abelian Hamiltionian BRST with a finite number of constraints that $\cH^{BRST}_{phys}$ is larger than the Dirac state space, and that adding the extra condition of restriction to ghost number zero states does not fix the problem (cf. MCPS problem and discussion in Subsection \ref{sbs:hhamexphsp}). We may however conjecture that the BRST physical algebra is isomorphic to the Dirac physical algebra. An application of Theorem \eqref{pr:krdel} shows that this is not the case, \ie $\cP^{BRST}$ is larger than physical algebra associated to the Dirac algorithm and adding the extra condition of restricting to elements with zero ghost number does not fix the problem.

Let $\cA_0=B(\cH_0)$ for a Hilbert space $\cH_0$. Take $\cC=\{P\}$ where $P\in B(\cH_0)$ is a projection and so the physical subspace is $\cH^0_{phys}=\ker P$. If we now apply the $T$-procedure to $(\cA_0,\{P\})$ (cf. Appendix \eqref{app:Tp}) then for the Dirac physical algebra we have the following isomorphism $\cP_{phys}\cong (1-P)P'(1-P)$ where $P'$ is the commutant in $\cA_0$ of $P'$ (see example \eqref{ex:DC} (i)).

Following the Hamiltonian BRST method in Section \ref{sec:hamBRST} we extend $\cA_0$ by tensoring on a ghost algebra $\rga=\osalg{\gh, \cgh}$ with one ghost for the constraint $P$.   The extended algebra is $\cA=\cA_0\otimes \rga$, and we identify $\cA_0$ with $\cA_0 \otimes \one$ (the norm on $\cA$ is unique as $\rga$ finite dimensional, hence nuclear). As all operators involved are bounded, we take $D(Q)=\cH_0\otimes\cH_g=\cH$ where $\cH_g$ is a ghost space for $\rga$. Note that by Definition \eqref{df:ghsp} and lemma \eqref{lm:rghlm} $\dim(\cH_g)=2$ and $\rga$ acts as an irreducible representation $\pi_{bz}:\rga \to B(\cH_g)$, therefore $\pi_{bz}(\rga)=B(\cH_g)$. Now \cite{Tak2002} p185 equation (10) gives that $B(\cH_1)\overline{\otimes} B(\cH_2)=B(\cH_1 \otimes \cH_2)$ for Hilbert space $\cH_1, \cH_2$ where $B(\cH_1)\overline{\otimes} B(\cH_2)$ is defined in \cite{Tak2002} Definition 1.3 p185. In fact as $B(\cH_2)$ is finite dimensional we have $B(\cH_1)\overline{\otimes} B(\cH_2)=B(\cH_1){\otimes} B(\cH_2)$. So as the norm on $\cA$ is unique we have $\cA=B(\cH_0){\otimes}B(\cH_g)=B(\cH)$, hence we do not need to extend $\cA$ to $\cAe$.

The BRST charge is $Q=P\otimes \gh\in B(\cH)$, and $\drb(\cdot)=\sbr{Q}{\cdot}$ with $D(\drb)=\cA$.  As $Q$ is bounded we do not have domain problems with $\Delta$ and we calculate as in equation \eqref{eq:hcomhamdel} that $\ker \Delta=\cH^0_p\otimes \cH_g$. Therefore we have that $P_s=(\one-P)\otimes \one$. From this we see that $\Phi^{ext} (\cA_0 \otimes \one)=P_s(\cA_0 \otimes \one)P_s= (\one -P)\cA_0 (\one-P)\otimes\one  =(\one -P)P' (\one-P)\otimes\one\cong \cP_{phys}$ with the obvious isomorphism from above. Now $P_s (\cA_0 \otimes \one) P_s \subset \ker \delta$ so by Theorem \eqref{pr:krdel},
\begin{equation}\label{eq:DirBRSTalgcomp}
\cP_{phys}\cong \Phi^{ext}_s (\cA_0\otimes\one ) \subset \Phi^{ext}_s (\ker \delta) \cong \cP^{BRST},
\end{equation}
and so the physical algebra produced by the $T$-procedure is contained in the physical BRST algebra. 

We see this containment is proper as follows. Take $P_2 \in \cB(\cH_0)$ such that $PP_2=0$ (e.g. take $P_2$ a projection orthogonal to $P$), and let $A=P_2\otimes \gh \cgh$. Then:
\begin{align*}
\drb(A)&=\drb(P_2\otimes \one)(\one \otimes \gh \cgh) +(P_2\otimes \one) \drb(\one \otimes \gh) (\one \otimes \cgh) -(P_2\otimes  \gh) \drb(\one \otimes \cgh),\\
&=[P,P_2]\otimes \gh^2 \cgh +P_2P\otimes \{\gh,\gh\} \cgh-P_2P\otimes \gh \{\gh, \cgh\},\\
&=0,
\end{align*}
therefore $A \in \ker \drb$. Now take $B\in \Phi^{ext}_s (\cA_0 \otimes \one)$. Then we have that $(A-B)\in \ker \drb$ and as $\gh,\, \cgh$ are linearly independent, $(A-B) \neq 0$. Now we have $P_s (A-B) P_s=((\one-P)P_2(\one-P)\otimes \gh \cgh)- B=(A-B)\neq 0$, and so by Theorem \eqref{pr:krdel} we have $(A-B) \notin \ran \drb$. Therefore $A \notin \Phi^{ext}(\cA_0 \otimes \one)P_s + \ran \drb$ and so the BRST physical algebra properly contains $\cP_{phys}$ if we take the natural isomorphisms as in \eqref{eq:DirBRSTalgcomp}.

An important point is that $A \in \cG_0$ and so restricting the BRST physical algebra to ghost number zero elements will still imply that the physical Dirac algebra is properly contained in the physical BRST algebra.

\subsection{Alternative Physical Algebra} \label{sbs:alalg}
So far we used the heuristic prescription of taking $\ker \drb/ \ran \drb$ as the physical algebra. We would like to motivate this from a more structural point of view. Recall the previous definitions:
\begin{itemize}
\item[(i)] $\cH^{BRST}_{phys}:=\ker \cQ /\cH_d$ and has indefinite inner product $\iip{\hat{\psi}}{\hat{\xi}}_p=\iip{P_s\psi}{P_s \xi}$ (cf. Definition \eqref{df:opbrsphsp} and Definition \eqref{df:BRSTphysindef}),
\item[(ii)] $\vp:\ker \cQ \to \cH^{BRST}_{phys}$ and $\tau: \ker \drb  \to \ker \drb /\ran \drb$ are the factor maps,
\item[(iii)] $\tau: (\ker \drb \cap \cA) \to \cP^{BRST}$,
\item[(iv)] $\cD^{BRST}_{phys}:=\vp(\ker Q)$ is dense in $\cH^{BRST}_{phys}$,
\item[(v)] $\hat{\psi}:=\vp(\psi),\, \forall \psi \in \ker Q$ and $\hat{A}:=\tau(A), \, \forall A \in \ker \drb,$ 
\end{itemize}
 
Now $\ker \drb$ preserves $\ker Q$ and $\ran Q$, hence factors to $\cD^{BRST}_{phys}$. Furthermore $(\ran \drb) \ker Q \subset \ran Q$ and so factors trivially to $\cD^{BRST}_{phys}$. Therefore $\ker \drb$ has a natural action on $\cD^{BRST}_{phys}$ by,
\begin{equation*}
\pi_{phys}\hat{\psi}:=\widehat{A\psi} , \qquad \forall A \in (\ker \drb \cap \cA), \, \psi \in \ker Q.
\end{equation*}
Define,
\begin{align*}
\Phi_s(A):\cA \to \cAe,\qquad \text{by} \qquad \Phi_s(A):=P_sA P_s, 
\end{align*}
that is $\Phi_s=\Phi^{ext}_s|_{\cA}$. By lemma \eqref{lm:alkerhom} we have that $\Phi_s$ is a homomorphism on $\ker \drb \cap \cA$, and:
\begin{proposition}\label{pr:phrepaldf}
We have that:
\[
\pi_{phys}(\ker \drb \cap \cA) \cong  \Phi_s( \ker \drb \cap \cA),
\] 
where $\cong$ denotes an algebra isomorphism. 
\end{proposition}
\begin{proof}
As $\cH^{BRST}_{phys}$ is a Krein space and $\cD^{BRST}_{phys}$ is dense we have that $\pi_{phys}(A)=0$ iff $\iip{\hat{\xi}}{ \pi_{phys}(A)\hat{\psi}}_p=0$ for all $\xi, \psi \in \ker Q$ iff $\iip{\hat{\xi}}{ \hat{A}\hat{\psi}}_p=\iip{P_s\xi}{P_sA\psi}=\iip{P_s\xi}{P_sAP_s\psi}=0$ for all $\xi, \psi \in \ker Q$ iff $P_s A P_s=0$, where $P_s$ is the projection onto $\cH_s$. From this it follows that $(\ker \Phi_s \cap \ker \drb \cap \cA)= (\ker \pi_{phys}\cap \ker \drb \cap \cA)$, hence $\Phi_s$ and $\pi_{phys}$ are homomorphisms on $\ker \drb \cap \cA$ with the same kernel, and hence have isomorphic images. This proves the first isomorphism, the second being obvious.
\end{proof}
By Theorem \eqref{pr:krdel} (iii) and the Proposition \eqref{pr:phrepaldf} we get that $\cP^{BRST}\cong \pi_{phys}(\ker \drb \cap \cA)$ and so we see that the cohomological definition of the BRST physical algebra is natural with repect to the constraint structure corresponding to the physical representation $\pi_{phys}:\ker \drb \cap \cA \to \op(\cH^{BRST}_{phys})$.

To prove Theorem \eqref{pr:krdel} (iii) however, we had to extend $\cA$ to $\cAe$. In particular we had to include $Q$ in the extended algebra. In examples, such as QEM below, $Q$ is an unbounded operator and so constructing a proof similar to Theorem \eqref{pr:krdel} (iii) in an abstract $C^*$-algebraic setting is not straightforward. So we take the point of view from now on that the abstract definition of the BRST-physical algebra is:
\begin{definition}\label{df:alphal}
{The alternative definition of the BRST physical algebra} is: 
\[
\tilde{\cP}^{BRST}:=(\ker \drb \cap \cA) / (\ker \drb \cap \ker \Phi_s \cap \cA).
\]
\end{definition}
The heuristic cohomological definition of $\cP^{BRST}=(\ker \drb \cap \cA)/(\ran \drb \cap \cA)$ is of physical interest because of Theorem \eqref{pr:krdel} (iii). Now we have that $\tilde{\cP}^{BRST}=\cP^{BRST}$ by Theorem \eqref{pr:krdel} (iii), and so the notation $\tilde{\cP}^{BRST}$ is redundant at this point. However when we move to the $C^*$-setting, $\cA$ will be a $C^*$-algebra with superderivation $\drb$ acting on domain $D_2(\drb)\subset \cA$ and $\Phi_s$ can be constructed using the universal representation of $\cA$, hence $\tilde{\cP}^{BRST}$ will be the natural definition the BRST physical algebra, and it will only be in representations were we can construct the extra unbounded elements in $\cAe$ (e.g. $Q$) in which the heuristic cohomological definition will become equivalent.

\begin{rem}\label{rm:alus}
Note the similarity of the definition of $\tilde{\cP}^{BRST}$ to equation \eqref{eq:Phkphys} in Appendix \eqref{app:Tp}. This suggests we can take a more Dirac-like approach to physical algebra selection for BRST. This natural question is to ask is how a Dirac-like constraint procedure ($T$-procedure Appendix \eqref{app:Tp}) can be applied to using $Q$ or $P_s$ as a constraint, and how this relates to the BRST-observable algebra $\cP^{BRST}$. 

If we were to follow a Dirac-like constraint procedure using $P_s$ as a constraint, we would then take observables as self-adjoint operators in the field algebra preserving $\ker P_s$, and take the observable algebra as the $*$-algebra generated by these. $\Phi_s$ would then be a homomorphism on this $*$-algebra and we would factor out $\ker \Phi_s$ to get the constrained algebra. The $T$-procedure is the abstract $C^*$-version of this process.

There is an extra complication with BRST in that we take $\iip{\cdot}{\cdot}$ as our physical inner product. Hence following a Dirac-BRST like procedure, observables will be $\dag$-algebras and the trivial observables will be those which map $\ker Q$ to $\overline{\ran Q}$. To incorporate this at the abstract level means that we would have to modify the $T$-procedure as well. This can be done, but we will not pursue this here.

Note however that the Dirac-BRST approach to selection of observables as described in the preceding paragraph is taken in K\&O \cite{KuOj79} p58. The result is then that \emph{local observables}, \cite{KuOj79} p59,  are in $\ker \drb$, which is proven by the authors using the Reeh-Schlieder theorem \cite{KuOj79} Proposition 5.4 p59. This is why we are motivated to define the constrained observable algebra as $\Phi_s(\ker \drb)$.
\end{rem}

\chapter{\KOB and Electromagnetism} \label{ch:BRSTQEM}

In this chapter we give a mathematically consistent interpretation of the BRST-QEM example in Section \eqref{sec:exheuEM}. We will smear the fields $A(x)$ and ghost field $\gh(x)$ over appropriate test function spaces and give a well-defined meaning to `replacing the gauge parameter by a ghost parameter'. This is open to interpretation and we take the view that the object of primary importance is the BRST superderivation $\drb$ as this selects the physical algebra. We will define a superderivation $\drb$ which gives the smeared versions of  \eqref{eq:hdA}, \eqref{eq:hdgh}, \eqref{eq:hdcgh}: 
\begin{align*} 
\drb(A^\nu(x))&= -i\partial^\nu \gh(x), \\
\drb(\gh(x))&= 0, \\
\drb(\tilde{\gh}(x))& = -i\partial_\nu A^\nu(x), 
\end{align*}
will find a BRST charge $Q$ which generates the the superderivation $\drb(\cdot)=\sbr{Q}{\cdot}$, show that $Q^2=0$, and that we have we have a smeared analogue of equation \eqref{eq:HDel}:-
\begin{align*}
\Delta= & \;\{Q,Q^*\},\\
= & \;2\int d^3\bp \;p_{0}^{2} (b_1(\bp)^*b_1(\bp) +b_2(\bp)^*b_2(\bp) +c_1^*(\bp)c_1(\bp) +c_2^*(\bp)c_2(\bp)).
\end{align*} 
Using the results of the previous chapter, such as the \emph{dsp}-decomposition, we will construct the physical space and algebra and compare them to other approaches.

One could try to make sense of the formal definition,
\begin{align*}
Q&:=\int_{x_0=const} d^3\mathbf{x} \,(\partial_\nu A^\nu )\overleftrightarrow{\partial_0} \gh, 
\end{align*}
but as the integrand involves  a pointwise multiplication of operator valued distributions, this is not straightforward to define. Such a definition of $Q$ need not generate the correct superderivation, as was found in the example \cite{HorVo1992} Theorem 1 p2827. Given the uniqueness of the BRST charge associated to $\drb$ (lemma \eqref{lm:unch}) and the central importance of $\drb$, it seems that the definition of $Q$ is best used as a guide for the form of $\drb$.

Section \ref{sbs:testfunc} gives an account of the necessary structures for the BRST method in the one particle test function space for QEM and the associated Ghost algebra. In Section \ref{sec:FKBRST} we assume a general abstract test function space with the relevant structures of the QEM and give a description of the smeared bosonic fields in a Fock-Krein representation as in the CCRs \cite{Min1980}. We do not fix the specific test function space, so we can extend the method to other bosonic field theories with a similar one particle test function space to QEM. We then tensor on an appropriate ghost algebra, construct the BRST superderivation and show that this gives the correct heuristic BRST-QEM superderviation for the QEM test function space. We derive the physical subspace and algebra for the general abstract test function space. We refer to this abstract model as \KOB as it based on the BRST constructions developed in the foundational work \KO \cite{KuOj79} in the abelian case. Using the QEM test function space we compare BRST-QEM  with Gupta-Bleuler QEM as in \cite{Hendrik2000}.  Using the construction based on the abstract test function space we do the examples of \KOB with a finite number of bosonic constraints and BRST for massive abelian gauge theory and compare these examples with the literature. Finally, we combine \KOB with a finite number of bosonic constraints with finite abelian Hamiltonian BRST to avoid the MCPS problem of Hamiltonian BRST (cf. Subsection \ref{sbs:hhamexphsp}).

\section{Test functions}\label{sbs:testfunc} 

To make the formal BRST-QEM example in Section \ref{sec:exheuEM} well-defined, we first smear the fields over test function spaces to get operator valued distributions.
 
\subsection{QEM Test Functions}\label{sbs:QEMtf}
 Recall the formal definition of the QEM gauge potential in equation \eqref{eq:Adefdag}:
\begin{equation*}
A^{\mu}(x)=\FC \int_{C_{+}} \frac{d^3\bp}{\sqrt{p_0}} ( a^{\mu}(\bp)e^{-ipx} + a^{\mu}(\bp)^{\dag} e^{ipx} ), \, \mu =0,1,2,3, 
\end{equation*}
where $\dag$ is the Krein involution with respect to the heuristic Krein-Fock space $\cH_0$, and has CCR's (equation \eqref{eq:guv}),
\[
[ A^{\mu}(x)\,,\, A^{\nu}(y)] = -i g^{\mu \nu} D_0(x-y), 
\]
where, 
\begin{equation*}
D_0(x):=(2\pi)^{-3} \int_{C_+}\frac{d^3\bp}{p_0} e^{-i\bp.\mathbf{x}}\sin(p_0 x_0),
\end{equation*}
is the Pauli-Jordan distribution.

For the following mathematical developments, the previous pointwise objects need to be understood in the sense of distributions. Hence we next define the appropriate 1-particle test function space over which to smear $A(x)$. We use the notation for inner product spaces and symplectic spaces as in Appendices  \eqref{ap:IIP} and \eqref{ap:symp}.

Let $f\in \cS(\mathbb{R}^4,\mathbb{C}^4)$ and denote its Fourier transform by $\hat{f}_{\mu}(p):=(2\pi)^{-2}\int d^4x f_{\mu}(x)e^{-ipx}$. Note the unconventional inner product $px=p_0x_0-\bp.\mathbf{x}$ in the exponential term. We want to smear our gauge potentials of over real vector valued test functions, so we define:
\begin{gather}
\widehat{\cS}:= \{\hat{f} \,| \,f\in \cS(\mathbb{R}^4, \mathbb{R}^4) \})=\{ f\in \cS(\mathbb{R}^4, \mathbb{C }^4)\, |\, \overline{f(p)}=f(-p)\}) \notag \\
\text{with IIP:} \quad  \iip{f}{g}:=  2\pi\int_{C_{+}} \frac{d^3p}{p_0} \overline{f_{\mu}}(p)g^{\mu}(p), \qquad \forall f, g \in \widehat{\cS}. \label{eq:iipqem}
\end{gather}
Smear $A(x)$ as:
\begin{align}
A(\hat{f}):= & \;\int d^4x A_{\mu}(x)f^{\mu}(x), \label{eq:Asmear}\\
= & \; \sqrt{\pi}\int  \frac{d^3\bp}{\sqrt{p_0}} ( a^{\mu}(\bp)\hat{f}_{\mu}(p) + a^{\mu}(\bp)^{\dag} \overline{\hat{f}_{\mu}(p)} )), \qquad p_0=|\mathbf{\bp}|, \notag \\
= & \;(a(\hat{f})+a(\hat{f})^{\dag})/\sqrt{2}, \notag \\
\text{with} \qquad a(f):= & \;\sqrt{2\pi}\int_{C_{+}}  \frac{d^3\bp}{\sqrt{p_0}} ( a^{\mu}(\bp){f}_{\mu}(p)), \notag
\end{align}
where $\hat{f} \in \widehat{\cS}$, \ie $f\in\cS(\mathbb{R}^4,\mathbb{R}^4)$ . The commutation relations \eqref{eq:dagheuccr} mean in our present context,
\begin{align}
[A(f),A(g)]= & \;\pi\int_{C_{+}} \frac{d^3\bp}{{p_0}}\left(f_{\mu}(p)\overline{g^{\mu}(p)}-\overline{f_{\mu}(p)}g^{\mu}(p)\right), \label{eq:Acommsym}\\
= & \; -i\, \im \iip{f}{g} \notag
\end{align}
The symplectic form for the smeared CCR's for QEM is then
\begin{equation}\label{eq:covsym}
\smp{f}{g}{1}:= -\im \iip{f}{g} \qquad \forall f,g \in \widehat{\cS},
\end{equation}
and we call this the \emph{covariant symplectic form}. Note that $\sigma_1$ is the Fourier transform of the Pauli-Jordan distribution, ie,
\begin{equation}\label{eq:JPsym}
\smp{\hat{f}}{\hat{g}}{1}= \int \int d^4x\,d^4y\,f_{\mu}(x)g^{\mu}(y)D_0(x-y), \qquad \forall f,g \in \widehat{\cS}.
\end{equation}
We define ${J}$ on $\cS(\mathbb{R}^4, \mathbb{C}^4)$ by:
\begin{equation}\label{eq:Jdef}
({J}f)_0:=-f_0, \quad \text{and} \quad ({J}f)_l=f_{l}, \; l=1,2,3.
\end{equation}
Using this we define a second inner product and symplectic form on $\widehat{\cS}$ by:
\begin{align}
\ip{f}{g}:= & \;\iip{f}{Jg}=2\pi \int_{C_{+}} \frac{d^3p}{p_0} \overline{f_{\mu}}(p)g_{\mu}(p), \label{eq:auxsymp}\\
\smp{f}{g}{2}:= & \;\smp{f}{Jg}{1}=- \im \ip{f}{g}, \notag
\end{align}
and we see that $\ip{\cdot}{\cdot}$ is positive semidefinite. We call $\sigma_2$ the \emph{\ncsf}. We want our inner products and symplectic forms to be non-degenerate so we define, 
\[
\fXp:= \widehat{ \cS}/\ker \sigma_1,
\]
and let $\rho_{\fXp}$ be the factor map. By slight abuse of notation denote $\sigma_1$ factored to $\fXp$ also by $\sigma_1$.

Note that $\ker \sigma_1$ consists of all the functions in $\widehat{\cS}$ which vanish on the forward light cone $C_{+}$, so $\rho_{\fXp}$ can be thought of as restricting $\widehat{\cS}$ to $C_{+}$. $\ker \sigma_1$ is also the degenerate part of $\widehat{\cS}$ with respect to both $\iip{\cdot}{\cdot}$, $\ip{\cdot}{\cdot}$ and $\sigma_2$, and so the latter forms factor to $\fXp$ and we still denote the factored forms by $\iip{\cdot}{\cdot}$, $\ip{\cdot}{\cdot}$ and $\sigma_2$ respectively. The inner product $\ip{\cdot}{\cdot}$ is now positive definite on $\fXp$.

It is easy to check that ${J}$ preserves $\ker \sigma_1$, hence it factors to an operator on $\fXp$, which we denote still by $J$. Then $J$ preserves both inner products, is symplectic in both symplectic forms and:
\[
\ip{f}{g}=\iip{f}{Jg}, \qquad \smp{f}{g}{1}=\smp{f}{Jg}{2}, \qquad \forall f,g \in \fXp.
\] 
As it stands, $\fXp$ is a real symplectic space with respect to $\sigma_1$, however in subsequent constructions we will need a complexified version of it. We have that $\cS(\R^4,\C^4)=\widehat{\cS}+i\widehat{\cS}$ where $i$ is the usual multiplication by $i=\sqrt{-1}$. We extend $\iip{\cdot}{\cdot}$ from $\widehat{\cS}$ to $\cS(\mathbb{R}^4,\mathbb{C}^4)$ by
\[
\iip{f}{g}:=  2\pi\int_{C_{+}} \frac{d^3p}{p_0} \overline{f_{\mu}}(p)g^{\mu}(p), \qquad \forall f, g \in \cS(\R^4,\C^4),
\]
and,
\[
\smp{f}{g}{1}=:- \im \iip{f}{g}, \qquad \forall f,g \in \cS(\R^4,\mathbb{C}^4).
\]
Now $\ker \sigma_1$ is the subspace of functions in $ \cS(\R^4,\mathbb{C}^4)$ vanishing on $C_{+}$. Let $\rho$ be the factor map $\rho:\cS(\R^4,\C^4) \to \cS(\R^4,\C^4) / \ker \sigma_1$ and define,
\[
\fDp:=\rho(\cS(\R^4,\C^4)).
\]
Multiplication by $i$ in $ \cS(\R^4,\C^4)$ preserves $\ker \sigma_1$ and so factors to an operator also denoted by $i$ on $\fDp$. Note that $i$ is a complex structure on $\fDp$. Now $\rho_{\fXp}:\widehat{\cS}\to \fXp$ and $\rho|_{\widehat{\cS}}:\widehat{\cS} \to  (\cS(\R^4,\C^4)/\ker \sigma_1)$ are both homomorphisms on $\widehat{\cS}$ with the same kernel and hence $\fXp\cong \rho(\widehat{\cS})$. We will always assume this identification and not denote it explicitly. Therefore we have that,
\[
\fDp=\rho(\cS(\R^4,\C^4))=\rho(\widehat{\cS}+i\widehat{\cS})=\fXp+i\fXp.
\]
We also factor the indefinite inner product $\iip{\cdot}{\cdot}$ and the symplectic form $\sigma_1$ to $\fDp$, and denote the factored forms with the same symbol as their corresponding unfactored counterpart.

Now $J$ extends to $\fDp$ in the obvious way and $[J,i]=0$. On $\fDp$ we can define the positive definite inner product,
\begin{align*}
\ip{f}{g}:&=\iip{f}{Jg}, \qquad \forall f, g \in \fDp,\\
\intertext{that is,}
\ip{\rho(f)}{\rho(g)}&=2\pi\int_{C_{+}} \frac{d^3p}{p_0} (\overline{f_{\mu}}(p)g_{\mu}(p)), \qquad \forall f, g \in \cS(\R^4,\C^4),
\intertext{and define the symplectic form}
\smp{f}{g}{2}:&= -\im \ip{f}{g}= \smp{f}{Jg}{1}, \qquad \forall f, g \in \fDp.
\end{align*}
Let $\fH=\overline{\fDp}$ where closure is with respect to the norm given by $\ip{f}{f}^{1/2}$, hence $\fH$ is a Hilbert space. Now $J$ is $\ip{\cdot}{\cdot}$-isometric on $\fD$ and $J^2=\one$. Therefore $J$ extends to a $*$-unitary on $\fH$ such that $J^*=J$ and $J^2=\one$. Therefore by lemma \eqref{lm:JKsp} $(\fH,\,\iip{\cdot}{\cdot})$ is a Krein space with fundamental symmetry $J$. 
\begin{rem}
Note that $\sigma_1$ is causal on $\fXp$ but not its complexification $\fDp$. That is $\smp{\hat{f}}{\hat{g}}{1}=0$ for all $\hat{f}=\rho(\hat{f}_1)\in\fXp$ and $\hat{g}=\rho(\hat{g}_1)\in \fXp$ such that $f_1$ and $g_1$ have spacelike separated supports (this follows from the fact that the Jordan-Pauli distribution has support on $\overline{V}_{+}\cup (- \overline{V}_{+})$ where $\overline{V}_{+}$ is the positive light cone, cf. \cite{ReSi1975v2} Theorem IX.34 and p71-72 where the Pauli-Jordan distribution is denoted $C_2(x)$). This is not the case for all $f,g \in \fDp$. Although we will make constructions using field algebras with test functions in $\fDp$, when we restrict finally to the physical objects, the resulting physical fields will have test functions in $\fXp$ and so causality will not be violated for physical objects.
\end{rem}

Let $L^2(C_{+}, \mathbb{C}^4,\lambda)$ be the Hilbert space of square integrable functions on the positive light cone with respect to unique $\cL^{\uparrow}_{+}$-invariant measure, $\rd \lambda= 2\pi(\rd^3 \bp)/\np$ on $C_{+}$ 
as defined in \cite{ReSi1975v2} p70, (\ie $L^2(H_0, \mathbb{C}^4,\Omega_0)$ in \cite{ReSi1975v2} notation). We identify $\fDp$ with a dense subspace of $L^2(C_{+}, \mathbb{C}^4,\lambda)$ as follows. Recall that $\rho:\widehat{\cS} \to \fDp$ is the factor map and identify $g=\rho(f)\in \fDp$ with  $f|_{C_{+}}\in L^2(C_{+}, \mathbb{C}^4,\lambda)$. This identification is independent of the representative $f$ for a given $g$ and it is linear, injective, and is isometric with respect to the $\fH$ norm. Hence we can identify $\fDp$ with a dense subspace of $L^2(C_{+}, \mathbb{C}^4,\lambda)$, and so by taking closures we have that $\fH \cong L^2(C_{+}, \mathbb{C}^4,\lambda)$. In effect, we are thinking of $\fDp$ as the `Schwartz functions on $C_{+}$'. We will assume this identification in the sequel, will not denote the isomorphism explicitly, and  will identify $\fXp$ with a subspace of $L^2(C_{+}, \mathbb{C}^4,\lambda)$. 

Also $J$ acts on $L^2(C_{+}, \mathbb{C}^4,\lambda)$ as in equation \eqref{eq:Jdef} and by lemma \eqref{lm:JKsp} acts as a fundamental symmetry that makes $\fH$ into a Krein space. The Hilbert inner product on $L^2(C_{+}, \mathbb{C}^4,\lambda)$ is the extension of equation \eqref{eq:auxsymp} on $\fD$ which we also denote by $\ip{\cdot}{\cdot}$.  The indefinite inner product $\iip{\cdot}{\cdot}=\ip{\cdot}{J\cdot}$ is the extension of equation \eqref{eq:iipqem} from $\fD$ to $L^2(C_{+}, \mathbb{C}^4,\lambda)$, and $\iip{f}{g}=\ip{f}{Jg}$ for $f,g \in L^2(C_{+}, \mathbb{C}^4,\lambda)$. The symplectic forms $\sigma_1$ and $\sigma_2$ are the extended from $\fD$ to  $L^2(C_{+}, \mathbb{C}^4,\lambda)$ via equation \eqref{eq:covsym} and equation \eqref{eq:auxsymp}.

We need scalar valued test functions for the ghosts, and obtain these as follows. Let:
\begin{align*}
P_0 &:  \widehat{\cS} \to \cS(\mathbb{R}^4,\mathbb{C})\quad \text{be:}\quad P_0  (f_0 , f_1 , f_2 , f_3) \to f_0,
\end{align*}
\ie $P_0$ is the projection onto the first component of $f \in \cS$. So $P_0$  factors to $\fXp$ and we define:
\begin{gather}\label{eq:D0}
\widehat{\cS}_0:=P_0 \widehat{\cS} \subset \cS(\mathbb{R}^4,\mathbb{C}), \qquad \fXp_0:=P_0\fXp,\\ 
\fDp_0:=\fXp_0+i\fXp_0. \notag
\end{gather}
The inner product on $\widehat{\cS}$ induces one on $\widehat{\cS}_0$ by:
\[
\ip{f}{g}_0:=\ip{P_0f}{P_0g}, \qquad \forall f,g \in \widehat{\cS},
\]
Similarly we have symplectic forms $\sigma_0$ and inner products $\ip{\cdot}{\cdot}_0$ induced by $\sigma_2$ and $\ip{\cdot}{\cdot}$ respectively on $\widehat{\cS}_0$, $\cS(\R^4,\C)$, $\fXp_0$ and $\fDp_0$. 
Then $\ip{\cdot}{\cdot}$ is positive definite on $\fDp_0$ hence generates a norm $\norm{\cdot}_0$, and we define $\fH_0 =\overline{\widehat{\cS}_0}$ where closure is with respect to $\norm{\cdot}_0$. The map $P_0:\fDp \to \fDp_0$ is continuous with respect to $\norm{\cdot}_0$ and so $P_0$ extends to an orthogonal projection $P_0:\fH \to \fH_0$. As in the vector valued case, we identify $\fH_0$ with $L^2(C_{+}, \mathbb{C},\lambda)$ and $\fDp_0$ with a dense subspace of $L^2(C_{+}, \mathbb{C},\lambda)$.
 
Summarizing the discussion above:
\begin{itemize}
\item[(i)] $\fXp$ is a real vector space corresponding to the Fourier transforms of $\cS(\mathbb{R}^4,\mathbb{R}^4)$ restricted to the light cone $C_{+}$. $\fXp_0$ is a real vector space corresponding to the Fourier transforms of $\cS(\mathbb{R}^4,\mathbb{R})$ restricted to the light cone $C_{+}$. 
\item[(ii)] $\fDp=\fXp + i \fXp$ and $\fDp_0=\fXp_0 +i \fXp_0$ are complex vector spaces.
\item[(iii)] $\fXp$ and  $\fDp$ have a non-degenerate indefinite inner product $\iip{\cdot}{\cdot}$, a positive definite inner product $\ip{\cdot}{\cdot}$ and non-degenerate symplectic forms $\sigma_1$ and $\sigma_2$. Furthermore, $\fXp_0$ and  $\fDp_0$ have a non-degenerate positive definite inner product $\ip{\cdot}{\cdot}_0$ and a non-degenerate symplectic form $\sigma_0$.
\item[(iv)] $\sigma_1$ and $\sigma_0$ correspond to the Fourier transformed smeared  Jordan-Pauli distribution, 
\[
D_0(x)=(2\pi)^{-3} \int_{C_+}\frac{d^3\bp}{p_0} e^{-i\bp.\mathbf{x}}\sin(p_0 x_0),
\]
and are causal but $\sigma_2$ is not causal.
\item[(v)] $\fH$ and $\fH_0$ correspond to the closures of $\fDp$ and $\fDp_0$ respectively in $\ip{\cdot}{\cdot}, \ip{\cdot}{\cdot}_0$ respectively. There exists a unitary $J\in B(\fH)$ such that $J^2=\one$, $J$ preserves $\fXp$, $\fDp$, and $\iip{\cdot}{\cdot}=\ip{\cdot}{J\cdot}$, $\smp{\cdot}{\cdot}{1}=\smp{\cdot}{J\cdot}{2}$.
\item[(vi)] $\fH$ and $\fH_0$ can be identified with $L^2(C_{+}, \mathbb{C}^4,\lambda)$ and $L^2(C_{+}, \mathbb{C},\lambda)$ respectively, and so $\fDp$ and $\fDp_0$ are identified with dense subspaces of $L^2(C_{+}, \mathbb{C}^4,\lambda)$ and $L^2(C_{+}, \mathbb{C},\lambda)$ respectively. 
\end{itemize}
Pointwise multiplication of the components of vectors in $\widehat{\cS}$ and $\cS(\mathbb{R}^4, \mathbb{C}^4)$ by $ip_{\mu}$ preserves $\ker \sigma_1$, and so factors to $\fXp$ and $\fDp$. Hence we have that multiplication by $ip_{\mu}$ is an unbounded operator on $\fH=L^2(C_{+}, \mathbb{C}^4,\lambda)$ which preserves the dense domain $\fDp$. Note that we need the $i$ factor if we want multiplication by $ip_{\mu}$ to preserve $\fXp$. Similarly we have that pointwise multiplication by $ip_{\mu}$ is an unbounded operator on $\fH_0$ that preserves $\fD_0$.

To analyze the Lorentz condition, we need the following notation: 
\begin{definition}\label{df:mulop} Let $(X,\kappa)$ be a measure space with positive measure $\kappa$,  $1 \leq q \leq \infty$ and $n \in \mathbb{N}$.
\begin{itemize}
\item[(i)] Let $f\in L^{q}(X, \mathbb{C}^n,\kappa)$ with components $f_j$, $j=1,\ldots,n$. When $f\in L^{q}(C_{+}, \mathbb{C}^4,\lambda)$ we use the usual convention that $f_{\mu} \in L^{q}(C_{+}, \mathbb{C},\lambda)$ is the $\mu$-th component of $f$ for $\mu=0,1,2,3$, and $f^{\mu}$ is the $\mu$-th component of $Jf$.
\item[(ii)] Let $g:X \to \mathbb{C}^n$ be measurable with respect to $\kappa$. Define:
\[
D(M_g)=\{ f\in  L^{q}(X,\mathbb{C}^n,\kappa)\,|\, \int_{X}\rd \kappa |g_j(p)f_j(p)|^q < \infty \},
\]
and the multiplication operator $M_g \in \op(L^{q}(X,\mathbb{C}^n,\kappa))$ by:
\[
M_gf(p):=g_j(p)f_j(p) \quad a.e,
\]
for $f \in D(M_g)$. Note that $M_g\in B(L^{q}(X,\mathbb{C}^n,\kappa))$ when $g\in L^{\infty}(X, \mathbb{C}^n,\kappa)$.
\item[(iii)] Define the bounded functions $b_j:\mathbb{R}^n \to \mathbb{C}$ by $b_j(p):= p_j/\np$ for $p\neq 0$ and $b(0)=0$, where $\nps=p_1^2+\ldots +p_n^2$. Define $b \in L^{\infty}(C_{+}, \mathbb{C}^4,\lambda)$ by:
\[
b(p):=(1/\np)(p_0,p_1,p_2,p_3) \quad a.e.
\]
Then $b_{\mu}(p)=p_{\mu}/\np$ and $b^{\mu}(p)=p^{\mu}/\np$ $a.e$. These are bounded since\\
$p \in C_{+} \Rightarrow \nps =p_0^2=\bp.\bp$.
\end{itemize}
\end{definition}
We now decompose $\fH$:
\begin{proposition}\label{lm:stdec2}
Define a map $P_{\JJ}: \fH \to \fH=L^{2}(C_{+},\mathbb{C}^4,\lambda)$ by:
\begin{align}\label{eq:P2def}
(P_{\JJ} f)_{\mu}:= & \; (1/2)M_{b^{\mu}}M_{b_{\nu}}f^{\nu}, 
\end{align}
{where we are summing over $\nu$. That is,}
\begin{align*}
(P_{\JJ} f(p))_{\mu}:= & \; \left(\frac{p^{\mu}p_{\nu}}{2\norm{\mathbf{p}}^2}\right)f^{\nu}(p), \qquad a.e. \qquad \text{for $p \in C_{+}$}.
\end{align*}
Then $P_{\JJ}$ is a projection on $\fH$, $P_{\LL}:=JP_{\JJ}J$ is a projection and 
\[
(P_{\LL} f(p))_{\mu}:= \left(\frac{p_{\mu}p_{\nu}}{2\norm{\mathbf{p}}^2}\right)f_{\nu}(p), \qquad a.e \qquad \text{for $p \in C_{+}$},
\]
for all $f \in \fH$, $p \in C_{+}$. Furthermore $P_{\LL} P_{\JJ}=0$, $JP_{\LL}=P_{\JJ}J$, and we have the following decomposition of $\fH$:
\[
\fH=\fH_t\oplus \fHL \oplus \fHJ,
\] 
where $\fHL:=\ran P_{\LL}$, $\fHJ:=\ran P_{\JJ}$, $\fH_t:=(\fHL\oplus \fHJ)^{\perp}$. Furthermore we have,
\begin{align}
\fH_t &=\{ f \in \fH \,| \, \mathbf{p.f(p)}=0, f_0(p)=0 \},\label{eq:tsHidec}\\
\fHL&= \{ f \in \fH \,| \, f_{\mu}(p)=(ip_{\mu}/\norm{\mathbf{p}})h(p)\; \text{\textrm{for}}\; h\in \fH_0 \},\notag\\
\fHJ&=J\fHL= \{ f \in \fH \,| \, f_{\mu}(p)=(ip^{\mu}/\norm{\mathbf{p}})h(p)\; \text{\textrm{for}}\; h \in \fH_0\},\notag\\
\fH_3&=\fHJ^{\perp}=\fH_{t}\oplus \fHL =\{ f \in \fH\,| \, p_\mu f^{\mu}(p)=0 \}\notag,
\end{align}
where the above multiplications are understood to be a.e.
\end{proposition}
\begin{proof}
For  $\mu=0,1,2,3$,  we have that $b_{\mu} \in L^{\infty}(C_{+}, \mathbb{C},\lambda) $ hence $M_{b_{\nu}} \in B(\fH_0)$ so the operator $Tf: \fH \to \fH_0$ defined by,
\[
Tf:=(P_{\JJ}f)_{\mu}=(1/2)M_{b^{\mu}}M_{b_{\nu}}f^{\nu},
\]
is bounded and so $P_{\JJ} \in B(\fH)$. Now for $p \in C_{+}$ we get that $p_{\mu}p_{\mu}=p_0^2+\bp.\bp=2\norm{\mathbf{p}}^2$ from which it follows that $P_{\JJ}^2=P_{\JJ}$. As $b_{\mu}$ is real valued we have that $P_{\JJ}^*=P_{\JJ}$. 

Now as $J^2=1$ we have that $P_{\LL}$ is a projection and $JP_{\LL}=J^2P_{\JJ}J=P_{\JJ}J$. The explicit form of $P_{\LL}$ follows directly by calculating $(JP_{\JJ}Jf)_{\mu}(p)$ for $f\in \fH$.

Since $p_{\mu}p^{\mu}=0$ for $p\in C_{+}$, we have for $f,h \in \fH$ that:
\begin{align*}
\ip{P_{\LL}f}{P_{\JJ}h}_{\fH}= & \; 2\pi  \int_{C_{+}} \frac{d^3p}{p_0}(p_{\mu}p^{\mu})\overline{(p_{\nu}f_{\nu}(p))}(p_{\nu}h^{\nu}(p))/(4\norm{\mathbf{p}}^4),\\
= & \;0,
\end{align*}
hence $P_{\LL} P_{\JJ}=0$. To prove the relations \eqref{eq:tsHidec} we proceed as follows. By definition,
\[
\fHJ=\ran P_{\JJ} = \{ f \in \fH\,|\, f_{\mu}=(1/2)M_{b^{\mu}}M_{b_{\nu}}g^{\nu}, \, g \in \fH \} \subset \{ f \in \fH\,|\, f_{\mu}=iM_{b^{\mu}}h, \, h \in \fH_0 \}.
\] 
For the reverse inclusion let $f_{\nu}=iM_{b^{\nu}}h$ for some $h \in \fH_0$ then we calculate:
\[
(P_{\JJ}f)_{\mu}=i(1/2)M_{b^{\mu}}M_{b_{\nu}}M_{b_{\nu}}h=iM_{b^{\mu}}h=f_{\mu},
\]
using again that $p_{\mu}p_{\mu}=2\norm{\mathbf{p}}^2$ for $p \in C_{+}$, hence $f \in \ran P_{\JJ} = \fHJ$. 

The equality for $\fHL$ follows from $\fHL=J \fHJ$ and the action of $J$ on $\fHJ$.

For the $\fH_3$ equality we take $f\in \fH$ and let $h=f-P_{\JJ}f$ then, 
\[
p_{\mu}h^{\mu}(p)= p_{\mu}f^{\mu}(p)- p_{\mu}p_{\mu}/(2\norm{\mathbf{p}}^2)(p_{\nu}f^{\nu}(p))=p_{\mu}f^{\mu}(p)-p_{\nu}f^{\nu}(p)=0,
\]
and hence $\fH_3=\fHJ^{\perp} \subset  \{ f \in \fH\,: \, p_\mu f^{\mu}(p)=0 \}$. For the reverse inclusion let $f \in \fH$ and $p_\mu f^{\mu}(p)=0$. Then for all $h \in \fH$ we have,
\begin{align*}
\ip{f}{P_{\JJ}h}= & \; 2 \pi \int_{C_{+}} \frac{d^3p}{p_0}(\overline{f^{\mu}}(p)(p_{\mu}p_{\nu}/\nps)h^{\nu}(p)),\\
= & \; 2\pi \int_{C_{+}} \frac{d^3p}{p_0}\overline{((p_{\mu}/\np)f^{\mu}(p))}((p_{\nu}/\np)h^{\nu}(p)),\\
=\;0
\end{align*}
hence $f\perp \fHJ$ and hence $\{ f \in \fH\,: \, p_\mu f^{\mu}(p)=0 \} \subset \fHJ^{\perp}=\fH_3$.

For the $\fH_t$ equality we adapt the argument in \cite{Hendrik2000} (p1193-1194). Let $f \in \fH_3 $ then $p_{\mu}f^{\mu}(p)=0$ implies $p_0f_0(p)=\bp.\mathbf{f}(p)$, hence
\begin{align*}
(P_{\LL} f(p))_{\mu}= p_{\mu}(p_{\nu}f_{\nu}(p))/(2\norm{\mathbf{p}}^2)
= & \;p_{\mu}(p_{0}f_{0}(p)+\mathbf{p.f}(p))/(2\norm{\mathbf{p}}^2),\\
= & \;p_{\mu}(\mathbf{p.f}(p))/\norm{\mathbf{p}}^2, 
\end{align*}
Therefore if we let $h=(f-P_{\LL}f)\in \fH_t$ we get that 
\[
\mathbf{p.h}(p)=\bp .(\mathbf{f}(p)- \bp(\bp .\mathbf{f}(p)/\nps))=0.
\]
Also as $h \in \fH_3$ this implies that $h_0(p)=0$, and so $\fH_t\subset \{ f \in \fH \,: \, \mathbf{p.f(p)}=0, f_0(p)=0 \}$. The reverse inclusion follows by direct verification of the orthogonality relations, using the equalities in \eqref{eq:tsHidec} already obtained.
\end{proof}

We would now like to obtain an analogous decomposition of the original test function space $\fDp$. However we have the complication that for $f=\rho(g) \in \fDp$ with $f(0)=\rho(g(0))\neq 0$, the $(1/2\norm{\mathbf{p}}^2)$ factor in $(P_{\JJ}f)_{\mu}$ poses the problem of whether there exists  $h \in \cS(\R^4,\C^4)$ such that $P_{\JJ}f=\rho(h)$ due to differentiability problems at $0$. So it is not evident that $P_{\JJ}$ preserves $\fDp$. To deal with this, we extend $\fDp$ to a more mathematically convenient test function space dense in $L^2(C_{+},\C^4,\lambda)$. Criteria that need to be taken into account when extending are:
\begin{itemize}
\item[(i)] Poincar{\'e} invariance: The Poincar{\'e} transformations on $\cS(\mathbb{R}^4,\mathbb{C}^4)$ are given by:
\begin{equation}\label{eq:Ptran}
(V_gf)(p):=e^{ipa}\Lambda f(\Lambda^{-1}p) \quad \forall f \in \cS(\mathbb{R}^4,\mathbb{C}^4), \quad g=(\Lambda,a) \in \cP^{\uparrow}_{+},
\end{equation}
which obvioulsy preserve $\cS(\mathbb{R}^4,\mathbb{C}^4)$. As $\Lambda C_{+} = C_{+}$ for all $\Lambda \in \cL^{\uparrow}_{+}$ we have that   $f|_{C_{+}}=0$ iff $V_gf|_{C_{+}}=0$ for all $g=(\Lambda,a) \in \cP^{\uparrow}_{+}$ and hence $V_g$ factors to and preserves $\fDp$. We denote the factored transformations with the same symbol, hence we have that the Poincar{\'e} transformations preserve $\fDp$. Now we can use formula \eqref{eq:Ptran} to define the Poincar{\'e} transformations on any function in $L^2(C_{+},\mathbb{C}^4,\lambda)$  (and adding the appropriate `a.e.'). This is because $V_g$ maps $L^2(C_{+},\mathbb{C}^4,\lambda)$ to $L^2(C_{+},\mathbb{C}^4,\lambda)$ functions for all $g=(\Lambda,a) \in \cP^{\uparrow}_{+}$ as $\lambda$ is $\Lambda$ invariant for all $\Lambda \in \cL^{\uparrow}_{+}$. We require that any extension of $\fDp$ will be preserved by all the Poincar{\'e} transformations so defined.
\item[(ii)] Causality: Causality of $\sigma_1$ on $\fXp$ means that for any $f,g \in \cS(\mathbb{R}^4,\mathbb{R}^4)$ with spacelike separated support, we have that $\smp{\hat{f}}{\hat{g}}{1}=0$. That is, spacelike separated support  refers to the support of $f,g\in \cS$ in $x$-space rather that support of $\hat{f}|_{C_{+}}$ in $p$-space. So to be able to give meaning to causality of $\sigma_1$ on an extension of $\fXp$ we will define the extension by first extending $\cS(\mathbb{R}^4,\mathbb{C}^4)$ to an appropriate space for which the inverse Fourier transform exists, and then restrict this space to the light cone.   \end{itemize}
The following extension incorporates both these criteria and produces  a physically reasonable model for QEM after applying the BRST procedure in Section \ref{sc:EM}. Let
\[
sgn(p):\mathbb{R}\to \mathbb{R}, \quad \text{be,} \quad sgn(p):=
\begin{cases}
1, \qquad p>0,\\
0, \qquad p=0,\\
-1, \qquad p<0
\end{cases}
\]
and define the bounded functions:
\[
c_k(p):\mathbb{R}^4\to \mathbb{R},\quad \text{by} \quad c_k(p):=\begin{cases}
\frac{sgn(p_0)}{\np}p_k, \qquad p\neq 0,\, k=1,2,3,\\
0, \qquad p=0,
\end{cases}, 
\]
and,  
\[
c_{k,g}(p):\mathbb{R}^4\to \mathbb{R},\quad \text{by} \quad c_{k,g}(p):=(V_gc_k)(p)=e^{ipa}c_k(\Lambda^{-1} p) \quad \text{for} \quad g=(\Lambda,a)\in \cP^{\uparrow}_{+},\qquad k=1,2,3.
\]
Note that:
\begin{itemize}
\item[(i)] $c_k(-p)=c_k(p)$.
\item[(ii)] $c_k=c_{k,\one}$ where $\one$ is the identity in $\cP^{\uparrow}_{+}$.
\item[(iii)] $M_{c_{k,g}}\in B(L^2(\mathbb{R}^4,\mathbb{C}))$ for all $g\in \cP^{\uparrow}_{+}$ as $c_{k,g}$ is bounded.
\end{itemize}
\begin{rem}\label{rm:schal}
The space of Schwartz functions $\cS(\R^4,\C)$ is an algebra with product defined by pointwise multiplication, as follows from \cite{ReSi1975v2} Theorem IX.3 (a) and  (b) p6 and the fact that $\cS(\R^4,\C)$ is preserved by the Fourier transform (\cite{ReSi1975v2} Theorem X.1 p3)
\end{rem}

\begin{proposition}\label{pr:extdec}
Let $M^m_{\mathbf{c_{j,g}}}:=M_{c_{j_1,g_1}}\ldots M_{c_{j_m,g_m}}$ for \,$m \in \mathbb{N}$,  $\mathbf{j}=(j_1,\ldots,j_m)\in \{1,2,3\}^m$ and $\mathbf{g}=(g_1,\ldots,g_m) \in (\cP^{\uparrow}_{+})^m$. Let:
\begin{align*}
\cSet_0&:=\alg{ f, M^m_{\mathbf{c_{j,g}}}h\,|\, f,h \in \cS(\mathbb{R}^4,\mathbb{C}), \; m \in \mathbb{N},\; \mathbf{j}\in \{1,2,3\}^m,\; \mathbf{g} \in (\cP^{\uparrow}_{+})^{m}, },\\
\cSet&:=\{ f \in L^2(\mathbb{R}^4,\mathbb{C}^4)\,|\, f_{\mu} \in \cSet_0,\,  \mu=0,1,2,3\}
\end{align*}
which is well defined as an algebra by Remark \eqref{rm:schal}. Furthermore, let
\begin{align*}
\fD&:=\{f|_{C_{+}}\,|\, f \in \cSet \}\supset \fDp, \text{hence $\fD$ is dense in $\fH=L^2(C_{+},\C^4, \lambda)$}\\
\fX&:=\{f|_{C_{+}}\,|\, f \in \cSet,\, \overline{f(p)}=f(-p)\quad a.e.\,\text{with respect to the Lebesgue measure on $\mathbb{R}^4$} \}\supset \fXp
\end{align*}
Let $P_{\LL}^{\cSet}:=JP_{\JJ}^{\cSet}J$ where $P_{\JJ}^{\cSet}\in B(L^2(\mathbb{R}^4,\mathbb{C}^4))$ is given by:
\[
(P_{\JJ}^{\cSet}f)_0:=(1/2)(f_0+M_{c_j}f_j), \qquad (P_{\JJ}^{\cSet}f)_k:=(1/2)M_{c_k}(f_0+M_{c_j}f_j)
\]
where $k=1,2,3$ and there is a summation over the repeated index $j=1,2,3$. Note that $P_{\JJ}^{\cSet}$ is bounded as $|c_j(p)|\leq 1$ for all $p\in \mathbb{R}^4$. 

Then:
\begin{itemize}
\item[(i)] $\cSet \subset (L^1(\mathbb{R}^4,\mathbb{C}^4)\cap L^2(\mathbb{R}^4,\mathbb{C}^4))$, hence $\check{f}\in (C_0(\mathbb{R}^4,\mathbb{C}^4)\cap L^2(\mathbb{R}^4,\mathbb{C}^4))$  for all $f \in \cSet$ where  $\check{f}_{\mu}:=(2\pi)^{-2}\int d^4x f_{\mu}(x)e^{ipx}$ is the inverse Fourier transform of $f\in \cSet$.
\item[(ii)] $f|_{C_{+}} \in L^2(C_{+},\mathbb{C}^4,\lambda)$ for $f \in \cSet$.
\item[(iii)] $V_g \fD \subset \fD$ and $V_g \fX \subset \fX$ for all $g\in \cP^{\uparrow}_{+}$.
\item[(iv)] $P_{\LL}^{\cSet}$ and $P_{\JJ}^{\cSet}$ are self-adjoint projections such that $P_{\LL}^{\cSet}\cSet \subset \cSet$ and $P_{\JJ}^{\cSet}\cSet \subset \cSet$.
\item[(v)] $(P_{i}^{\cSet}f)|_{C_{+}}=P_{i}(f|_{C_{+}})$ for $i=1,2$ and all $f \in \cSet$ where $P_{\LL}$ and $P_{\JJ}$ are the projection from Proposition \eqref{lm:stdec2}. Hence $P_{\LL}$ and $P_{\JJ}$ preserve $\fD$ and $\fX$.
\end{itemize}
\end{proposition}
\begin{proof}
(i) and (ii): As $|c^g_j(p)|\leq 1$ we get that $|M^m_{\mathbf{c_{j,g}}}h(p)|\leq |h(p)|$ for all $h \in \cS(\mathbb{R}^4,\mathbb{C})$, $m \in \mathbb{N}$,  $\mathbf{j}\in \{1,2,3\}^m$, $p\in \mathbb{R}^4$, $j=1,2,3,$ and all $g \in \cP^{\uparrow}_{+}$. As $\cS(\mathbb{R}^4,\mathbb{C}) \subset (L^1(\mathbb{R}^4,\mathbb{C}))\cap L^2(\mathbb{R}^4,\mathbb{C})$ we get that the components of elements in $\cSet$ are in $L^1(\mathbb{R}^4,\mathbb{C}))\cap L^2(\mathbb{R}^4,\mathbb{C})$ hence (i) follows. 

Similarly $f|_{C_{+}}\in L^2(C_{+},\mathbb{C}^4,\lambda)$ for all $f \in \cS(\mathbb{R}^4,\mathbb{C}^4)$ and  $|c_{j,g}(p)|\leq 1$ for all $p\in C_{+}$, $j=1,2,3,$ and all $g \in \cP^{\uparrow}_{+}$ which implies (ii).

By \cite{Rud1987} Theorem 9.6 p182 and that the Fourier transform is a unitary in $B(L^2(\mathbb{R}^4,\mathbb{C}^4))$, we have $\check{f}\in(C_0(\mathbb{R}^4,\mathbb{C}^4)\cap L^2(\mathbb{R}^4,\mathbb{C}^4))$.

\smallskip \noindent (iii): Suppose $g_1,g_2 \in \cP^{\uparrow}_{+}$. Then $V_{g_1}c_{k,g_2}=c_{k,g_1g_2}$ and as $\cS(\mathbb{R}^4,\mathbb{C}^4)$ is preserved by the Poincar{\'e} transformations, we get that $\cSet$ is preserved by the Poincar{\'e} transformations hence $V_g \fD \subset \fD$ for all $g\in\cP^{\uparrow}_{+}$. 

Let $ f \in \cSet$ and $\overline{f(p)}=f(-p)$ $a.e.$ with respect to the Lebesgue measure on $\mathbb{R}^4$. Then $\overline{(V_gf)(p)}=e^{-ipa}\Lambda \overline{f(\Lambda^{-1}p)}=(V_gf)(-p)$ $a.e.$ with respect to the Lebesgue measure on $\mathbb{R}^4$ for all $g=(\Lambda,a) \in \cP^{\uparrow}_{+}$. Hence $V_g\fX \subset \fX$ for all $g=(\Lambda,a) \in \cP^{\uparrow}_{+}$.

\smallskip \noindent (iv):
That $P_{\JJ}^{\cSet}\cSet \subset \cSet$ is obvious from the definitions. 
As $c_{j,g}$ is a real valued function, 
\begin{align*}
\ip{P_{\JJ}^{\cSet}f}{g}= &\;(1/2) \int 
\rd p^4\, \overline{(P_{\JJ}^{\cSet}f)_{\mu}(p)}g_{\mu}(p),\\
=&\; (1/2) \int 
\rd p^4\, [\overline{f_0(p)}g_0(p)+\overline{f_j(p)}c_j(p)g_0(p)+\overline{f_0(p)}c_k(p)g_k(p)+c_j(p)\overline{f_j(p)}c_k(p)g_k(p)],\\
 =&\;\ip{f}{P_{\JJ}^{\cSet}g}
\end{align*}
for all $f,g \in  L^2(\mathbb{R}^4,\mathbb{C}^4)$, and so $P_{\JJ}^{\cSet}$ is self-adjoint.  We get that $P_{\JJ}^2=P_{\JJ}$ by using $c_j(p)c_j(p)=\nps/\nps=1 \Rightarrow M_{c_j}M_{c_j}=1$ and direct substitution in the defining formulas. We give the calculation for the $0$-th component:
\begin{align*}
((P_{\JJ}^{\cSet})^2f)_0:= & \;(1/2)((P_{\JJ}^{\cSet}f)_0+M_{c_j}(P_{\JJ}^{\cSet}f)_j),\\
= & \;(1/4)(f_0+M_{c_j}f_j+M_{c_j}M_{c_j}(f_0+M_{c_k}f_k)),\\
= & \;(2/4)((P_{\JJ}^{\cSet}f)_0+M_{c_j}(P_{\JJ}^{\cSet}f)_j),\\
= & \;(P_{\JJ}^{\cSet}f)_0.
\end{align*} 
The calculations for the other components are similar.
 
It follows from $P_{\LL}^{\cSet}=JP_{\JJ}^{\cSet}J$, $J^{*}=J$ and $J^2=\one$ that $P_{\LL}^{\cSet}$ is a self-adjoint projection. As $J\cSet\subset \cSet$ we get $P_{\LL}^{\cSet}\cSet \subset \cSet$.

\smallskip \noindent (v): Now $b_0(p):=p_0/\np=1$ for $p \in C_{+}$ and so $M_{b_0}=\one$. Substituting this into the defining formula for $P_{\JJ}$ (equation \eqref{eq:P2def}) and substitution $c_j(p)=p_j/\np=b_j(p)$ for $p\in C_{+}$ in the formula for $P_{2}^{\cSet}$ we get that $(P_{2}^{\cSet}f)|_{C_{+}}=P_{2}(f|_{C_{+}})$, e.g. for the $k$-th component,
\[
(P_{2}^{\cSet}f)_k|_{C_{+}}=(1/2)M_{b_k}(f_0+M_{b_j}f_j)|_{C_{+}}=(1/2)M_{b_k}M_{b_{\nu}}f^{\nu}|_{C_{+}}=P_2(f|_{C_{+}})
\]
Similarly $(P_{1}^{\cSet}f)|_{C_{+}}=P_{1}(f|_{C_{+}})$. By (iv) $P_{\LL}^{\cSet}\cSet \subset \cSet \supset P_{\JJ}^{\cSet}\cSet $, hence $P_{\LL}\fD \subset \fD \supset P_{\JJ}\fD$.

Let $ f \in \cSet$ and $\overline{f(p)}=f(-p)$ $a.e.$ with respect to the Lebesgue measure on $\mathbb{R}^4$. Then, 
\begin{align*}
\overline{(P_{\JJ}^{\cSet}f)_k(p)}:= & \;(1/2)c_k(p)(\overline{f_0(p)}+c_j(p)\overline{f_j(p)}),\\
= & \;(1/2)c_k(-p)({f_0(-p)}+c_j(-p){f_j(-p)}), \qquad a.e.\\
= & \;\overline{(P_{\JJ}^{\cSet}f)_k(-p)}, \qquad a.e.
\end{align*}
where $a.e.$ above is with respect to Lesbegue measure, where we used that $c_j(p)=c_j(-p)$ $a.e.$ with respect to the Lesbegue measure on $\mathbb{R}^4$. Therefore $P_{\JJ} \fX \subset \fX$. As $J\fX \subset \fX$ we also get $P_{\LL}\fX \subset \fX$.
\end{proof}

Given Proposition \eqref{pr:extdec} (v) we get a decomposition of $\fX$ and $\fD$ analogous to $\fH$ in Proposition \eqref{lm:stdec2}.
\begin{theorem}\label{th:fdtstdec} Let
\begin{gather*}
\fX_j:=P_j\fX, \qquad \fD_j:=P_j\fD, \qquad j=t,\LL,\JJ,\\
\fX_3:=\fX \cap \fH_3 \qquad \fD_3:=\fD\cap \fH_3,\\
\fX_0:=\{f_0 \in \fH_0\,|\, f \in \fX\}, \qquad \fD_0:=\{f_0 \in \fH_0\,|\, f \in \fD\}
\end{gather*}
where $P_{t}\in B(\fH)$ is the projection on $\fH_t$.
Then:
\begin{itemize}
\item[(i)]We have the decompositions:
\begin{gather*}
\fX= \fX_t\oplus \fXL \oplus \fXJ \supset \fXp, \qquad \fD= \fD_t\oplus \fDL \oplus \fDJ \supset \fDp,\\
\fX_3=\fX_t \oplus \fXL, \qquad \fD_3 = \fD_t \oplus \fDL,
\end{gather*}
where $\oplus$ denotes orthogonality with respect to the $\ip{\cdot}{\cdot}$ inner product on\\
 $\fH=L^2(C_{+},\mathbb{C}^4,\lambda)\supset \fD$. 
\item[(ii)] We have that $\fDL=J\fDJ$ and $\fXL,\fXJ,\fDL,\fDJ$ are all $\sigma_1$-null subspaces.
\item[(iii)] $\fXp_0 \subset \fX_0$ and $\fDp_0 \subset \fD_0$. Also $M_{c_k}\fX_0 \subset \fX_0$, and $M_{c_k}\fD_0 \subset \fD_0$ where $k=1,2,3$. 
\end{itemize}
\end{theorem}
 \begin{proof}
(i): As $P_t,P_{\LL},P_{\JJ}$ are projections on $\fH$, we get for $f \in \fD$ that $f=P_tf+P_{\LL}f+P_{\JJ}f$, so by Proposition \eqref{pr:extdec} (v) we get the decompositions of $\fX$ and $\fD$. As $\fH_3=\fHL\oplus \fHJ$ we similarly get the decomposition of $\fX_3$ and $\fD_3$.

\pfit(ii): We have that $\fDL=J\fDJ$ follows from $P_{\JJ}=JP_{\LL}J$. From (i) we get $0=-\im \ip{\fXL}{\fXJ}=\smp{\fXL}{\fXJ}{2}=\smp{\fXL}{J\fXJ}{1}=\smp{\fXL}{\fXL}{1}$ hence $\fXL$ is a $\sigma_1$-null subspace. Similarly $\fXL,\fXJ,\fDL,\fDJ$ are all $\sigma_1$-degenerate subspaces.

\pfit (iii): $\fXp_0 \subset \fX_0$ and $\fDp_0 \subset \fD_0$ follow as $\fXp\subset \fX$ and $\fDp \subset \fD$.

Let $f \in \fD_0$. Then $f=g|_{C_{+}}$ for some $g \in \cSet_0$. By the definition of $\cSet_0$ in Proposition \eqref{pr:extdec} that $M_{c_k}\cSet_0\subset \cSet_0$ and hence $M_{c_k}f=(M_{c_k}g)|_{C_{+}}\in \fD_0$ for $k=1,2,3$. Hence $M_{c_k}\fD_0 \subset \fD_0$ and $M_{c_k}\fX_0 \subset \fX_0$ follows similarly.
\end{proof}
Lastly, we verify that extending $\fDp$ by $\fD$ preserves causality.
\begin{lemma} Let $\cSetin:=\{ \check{f}\,|\, f\in \cSet\}$ where $\check{f}_{\mu}:=(2\pi)^{-2}\int d^4x f_{\mu}(x)e^{ipx}$ is the inverse Fourier transform of $f\in \cSet$, and define $\tilde{f}:=\hat{f}|_{C_{+}}$ for all $f \in \cSetin$. Let $f, g \in \cSetin$ be such that $\tilde{f},\tilde{g} \in \fX$ and $f,g$ have spacelike seperated supports, then $\smp{\tilde{f}}{\tilde{g}}{1}=0$. That is $\sigma_1$ is causal on $\fX$.
\end{lemma}
\begin{proof}
Proposition \eqref{pr:extdec} (i) gives that $\cSetin$ is well defined. Let $f, g \in \cSetin$ be such that $\tilde{f},\tilde{g} \in \fX$ and $f,g$ have spacelike seperated supports. Now $\tilde{f}\in \fX$ implies that $\hat{f} \in \cSet$ and $\overline{\hat{f}(p)}=\hat{f}(-p)\quad a.e.$ with respect to the Lesbegue measure on $\mathbb{R}^4$, which in turn implies that $f=\check{\hat{f}}$ is real valued. That is $f \in  (C_0(\mathbb{R}^4,\mathbb{R}^4)\cap L^2(\mathbb{R}^4,\mathbb{R}^4))$. 

Using notation $x^2=x_{\mu}x^{\mu}$, we have by \cite{ReSi1975v2} Theorem IX.48 p107 that
\begin{align*}
\smp{\tilde{f}}{\tilde{g}}{1}=&\;\int \int dx^4\,dy^4\, f(x)g(y)D_0(x-y),\\
=&\;\int \int d^4x\,d^4y\, f_{\mu}(x)g^{\mu}(y)D_0(x-y),\\
=&\;\int \int d^4x\,d^4y\, f_{\mu}(x)g^{\mu}(D_0^{+}(x-y)-D_0^{+}(y-x)),\\
=&\;\int \int d^4x\,d^4y\, f_{\mu}(x)g^{\mu}[h( (x-y)^2)- h((y-x)^2)],\\
=0
\end{align*}
where $D_0(x),D_0^{+}(x)$ are symbolic notations for the Jordan-Pauli distribution and advance Jordan-Pauli distributions ($\Delta(x,0)$ and $\Delta_{+}(x,0)$ in \cite{ReSi1975v2} notation), and the fourth equality follows as $f(x)$ and $g(y)$ have spacelike seperated supports hence $D_0^{+}(x-y)$ is given as integration against the function $h:(0,\infty) \to \R$ given in \cite{ReSi1975v2} Theorem IX.48 (c).
\end{proof}

\subsection{Lorentz condition}\label{sbs:lorentz}
The Lorentz condition of classical electromagnetism is $\partial_\mu  A^\mu(x)=0$, which is well known to be problematic in the quantum context (\cite{Str1967}). Smearing against $ {h}\in \mathcal{S}(\mathbb{R},\mathbb{R})$ gives:
\begin{align*}
A(\hat{f})&=\int dx^4  \partial_\mu A^\mu(x)   {h}(x)\\
&=-\int dx^4   A^\mu(x)  \partial_\mu {h}(x),
\end{align*}
where $\hat{f}_{\mu}(p)=(\widehat{ \partial_{\mu} {h}})(p)=ip_{\mu} \hat{h}(p)$. Let
\[
\fXp_L= \{ f \in \fXp \,| \, f_{\mu}(p)=ip_{\mu}h(p)\; \text{for}\; h\in \fXp_0 \},
\]
and so $\{ A(f) \,|\, f\in \fXp_L\}$ corresponds to the smeared $\partial_\mu  A^\mu(x)$. 

Now by Theorem \eqref{th:fdtstdec} and the expression for $\fHL$ in equation \eqref{eq:tsHidec} we have  
\[
\fXL=\fX\cap \fHL= \{ f \in \fX \,| \, f_{\mu}(p)=ip_{\mu}h(p)\; \text{for}\; h\in \fH_0 \}.
\]
As we extend the test function space to $\fX \supset \fXp$ in the following QEM examples and so from the above expressions we associate $\{A(f)\,|\, f \in \fXL \}$ with the smeared $\partial_\mu  A^\mu(x)$. Similarly we will associate $\{A(f)\,|\, f \in \fDL\}$ with extended Lorentz condition when we are considering $\fD=\fX+i\fX$, the complexified version of $\fX$. 
\begin{rem}
Note that,
\[
\{ A(f) \;| \; f \in \fXp, f_{\mu}(p)=p_{\mu}{p^{\nu}}g_{\nu}(p), \quad p\in C_{+}, \quad g \in \fDp \}
\]
is properly contained in $\fXp_{L}$, as shown in \cite{Hendrik2000} p1193 remark 5.4 (ii). As the above set corresponds to the smeared Maxwell equations, we see that the Maxwell conditions are a proper subset of the Lorentz conditions.
\end{rem}

\subsection{Ghost Test Functions}\label{sbs:QEMGHtf}
To construct the BRST superderivation we first need to choose the test function space to use in the ghost algebra in Section \ref{sbs:ghgrad}. Recall the formal superderivation $\drb$ defined by equations \eqref{eq:hdA}, \eqref{eq:hdgh}, \eqref{eq:hdcgh}: 
\begin{align*} 
\drb(A^\nu(x))&= -i\partial^\nu \gh(x), \\
\drb(\gh(x))&= 0, \\
\drb(\tilde{\gh}(x))& = -i\partial_\nu A^\nu(x), 
\end{align*}
Smearing equation \eqref{eq:hdA} against $f\in \cS(\R^4,\R^4)$:
\begin{equation}\label{eq:smAdrb}
\drb(A(\hat{f}))=\int d^4x \,\drb(A^\nu(x))f_{\nu}(x)= i\int d^4x \, \gh(x)\partial^\nu f_{\nu}(x)=\gh(\hat{g}),
\end{equation}
where $\hat{g}(p):=ip^{\nu}\hat{f}_{\nu}(p)\in \cS(\R^4,\C)$. 

Smearing equation \eqref{eq:hdgh} is trivial. 

Smearing equation \eqref{eq:hdcgh} with $h\in \cS(\R^4,\R)$ gives: 
\begin{equation}\label{eq:smghdrb}
\drb(\tilde{\gh}(\hat{h}))=\int d^4x \, \drb(\tilde{\gh}(x))h(x) = i\int d^4x \, A^\nu(x)\partial_\nu h(x)=iA(\hat{k}), 
\end{equation}
where $k \in \cS(\R^4,\C^4)$ is defined by $\hat{k}_{\mu}(p):=ip_{\mu}\hat{h}(p)$.

These relations suggest using a scalar valued test function space contained in $\cS(\R^4,\C)$ to smear the heuristic fields $\gh$ and $\tilde{\gh}$ against. While this is a plausible direction to take, if we consider the expression for $P_{\JJ}$ in lemma \eqref{lm:stdec2}, then for $f\in \fD$:
\[
(P_{\JJ} f(p))_{\mu}:= \left(\frac{-ip_{\mu}}{2\norm{\mathbf{p}}^2}\right)(ip^{\nu}f_{\nu}(p))=\left(\frac{-ip_{\mu}}{2\norm{\mathbf{p}}^2}\right)g(p), \qquad a.e \qquad \text{for $p \in C_{+}$,}
\]
where $g(p):=ip^{\nu}{f}_{\nu}(p)$ is the same as in the RHS of equation \eqref{eq:smAdrb}. So we see that there is a connection between $\fDJ$ and the test functions on the RHS of $\drb(A(\hat{f}))=\gh(\hat{g})$. We will make this connection explicit below in lemma \eqref{lm:Tunit} where we show that there exists unitary $T:\overline{\fD_0} \to \overline{\fDJ}^{\mH}$ where the closure in $\overline{\fDJ}^{\mH}$ is with respect to the norm defined below in Proposition~\eqref{lm:mH}. As $T$ is unitary, the Ghost Algebra $\cA_g(\overline{\fDJ}^{\mH})=CAR(\overline{\fDL}^{\mH} \oplus \overline{\fDJ}^{\mH})$ will have the same CAR's as if we constructed a Ghost Algebra using the scalar test function space $\overline{\fD_0}$ for the ghosts. We will use the Ghost Algebra $\cA_g(\overline{\fDJ}^{\mH})$ in the BRST extension below and show, after defining the BRST superderivation as a map, that it gives the correct smeared relations above, cf. the discussion after Definition \eqref{df:brstsd}.

Two reasons for prefering to smear the ghosts over the vector valued function space $\overline{\fDJ}^{\mH}$ rather than a scalar test function space as the heuristics suggest are:
\begin{itemize}
\item[(i)] As $\fDL$ corresponds to the smeared $\partial^{\mu}A_{\mu}(x)$ (Subsection \eqref{sbs:lorentz}), we get that using $\cA_g(\overline{\fDJ}^{\mH})=CAR(\overline{\fDL}^{\mH} \oplus \overline{\fDJ}^{\mH})$ as the ghost algebra follows the philosophy of associating a ghost-conjugate ghost pair to each constraint.
\item[(ii)] Generalisation will be easier and we will use the constructions below for other examples where there is no corresponding scalar test funcion space $\fD_0$, e.g. finite \KOB in Subsection \ref{sbs:findimbos}.
\end{itemize}
We now analyze the relation between $\fDL\oplus\fDJ$ and the scalar valued function space $\fD_0$ as defined in Theorem \eqref{th:fdtstdec} and give explicit formulas that connect them.
\begin{definition}\label{df:T} 
\begin{itemize}
\item[(i)] Let $y(p):=\norm{\bp}^{-1}$ and define \begin{equation}\label{eq:dfscghis}
\fB:=\{ h\in \fD_0 \, | \, M_{y}h \in \fD_0  \}  \subset \fD_0
\end{equation}
\item[(ii)]                               
Let $T:\fB \rightarrow \fDJ$ be defined by:
\begin{align*}
(T h)(p):= \frac{h(p)}{{2}\nps} (-p_0 , p_1 , p_2 , p_3 ).
\end{align*}
\end{itemize}
\end{definition}
\begin{rem}
\begin{itemize}
\item[(i)] 
Note that the containment $\fB \subset \fH_0$ is proper since for any $f \in \fD_0$ with $f(0)\neq 0$ we have $f \notin \fB$. Moreover, $\fB$ is dense in $\fH_0$ in the $\norm{\cdot}_{\fH_0}$-topology by a simple approximation argument smooth bump functions with support outside outside an open set containing the origin. 
\item[(ii)] We have that $T:\fB\to \fD$ is well defined by the definition of $\fB$ and Theorem \eqref{th:fdtstdec} (iii). That $T(\fB) \subset \fDJ$ follows as $P_{\JJ}Th=h$ for all $h\in \fB$ which verified using the defining formula for $T$ above, $p_0^2=\nps$ for $p\in C_{+}$ and the formula for the components of $P_{\JJ}$ in Proposition \eqref{lm:stdec2}.
\end{itemize}
\end{rem}
Let $f \in \cS(\R^4,\R^4)$ and $g\in \cS(\R^4,\C)$ be defined by ${g}(p):=ip^{\nu}{f}_{\nu}(p)$. Then using equation \eqref{eq:P2def} it is immediate that $P_{\JJ}f =Tg$, that is $T$ maps the smearing function for $\drb(A^\nu(x))$ (cf. equation \eqref{eq:smAdrb}) to a function in $\fDJ$. To use $\fDJ$ as the test function space with which to smear the ghosts we would like $T$ to have dense range and to be isometric so that ghost algebras generated using test functions from $\fB$ or $\fDJ$ have equivalent CAR's. However it is straightforward to see that $T$ is \emph{not} isometric with respect to the $\fH=L^2(C_{+},\C^4,\lambda)$ norm on $\fD$ and $\fH_0=L^2(C_{+},\C,\lambda)$ norm on $\fD_0$. Hence we will define a new inner product on $\fDJ$ with respect to which $T$ is isometric and has dense range in the closure of $\fDJ$ with respect to this inner product.

\begin{proposition}\label{lm:mH} Let $z(p):=2\nps$ and define, 
\[   
\mH:=M_z, \qquad  D(\mH):=D(M_z)=\{ f\in \fH_0\,|\, \int_{C_{+}} \frac{d \bp^3}{p_0} |\nps f(p)|^2 < \infty \}.
\]
Then:
\begin{itemize}
\item[(i)]$\mH$ is a positive ($*$) and a Krein($\dag$) self-adjoint operator with  $\ker \mH= \{0\}$. Furthermore we have $\fD \subset D(\mH)$ and $[\mH, J]=[\mH, P_{t}]=[\mH, P_{1}]=[\mH, P_{2}]=0$ on $\fD$. Hence $\mH \fD_j \subset \fD$ for $j=t,\LL,\JJ$.  
\item[(ii)] Define a new inner product on $\fD$ by:
\[
\ip{\cdot}{\cdot}_{\mH}:=\ip{\cdot}{\mH \cdot}.
\]
We have that $\ip{\cdot}{\cdot}_{\mH}$ is positive definite and
\[
\fD=\fD_t \oplus_{\mH} \fDL \oplus_{\mH} \fDJ= \fD_3 \oplus_{\mH} \fDJ,
\]
where $\oplus_{\mH}$ denotes $\ip{\cdot}{\cdot}_{\mH}$ orthogonality. Let $\fL$ be the closure of $\fD$ with respect to this inner product, then:
\[
\fL=\fL_t \oplus_{\mH} \fLL \oplus_{\mH} \fLJ= \fL_3\oplus_{\mH} \fLJ,
\]
where $\overline{\fD}_j:= \fL_j$, $j=t,\LL,\JJ$ and the closure is with respect to the $\ip{\cdot}{\cdot}_{\mH}$ topology.
\item[(iii)] We have that $J$ is $\ip{\cdot}{\cdot}_{\mH}$-isometric on $\fD$ and $J$ extends to an $\ip{\cdot}{\cdot}_{\mH}$-unitary operator on $\fL$ such that $J^2=\one$. Define an indefinite inner product on $\fL$ by:
\[
\iip{\cdot}{\cdot}_{\mH}:=\ip{\cdot}{J \cdot}_{\mH}
\]
which makes $(\fL, \iip{\cdot}{\cdot}_{\mH})$ into a Krein space. Furthermore $\fLL$, $\fLJ$ are neutral subspaces. We denote $P^{\fL}_i$ for $i=t,\LL,\JJ,3$ for orthogonal projections on $\fL$ and note that $f \in \fD$ we have $P_if=P^{\fH}_if=P^{\fL}_if$, $i=t,\LL,\JJ,3$.
\end{itemize}
\end{proposition}
\begin{proof}
(i):That $[\mH, J]=0$ on $\fD$ is obvious. We have that $*$-self-adjointness follows from Proposition 1 \cite{ReSi1972v1} p259. Hence $\mH^{\dag}=J\mH^* J=\mH J^2=\mH$. $\ker \mH= \{0\}$ as $\nps=0$ iff $\bp=0$ for $\bp \in \R^3$, \ie the multiplier function $z(p)=0$ only for a set of measure zero. 

As $\mH$ is a multiplication operator of a polynomial in $p_1,p_2,p_3$ it follows that $\fD \subset D(\mH)$ and $[\mH, P_{1}]=[\mH, P_{2}]=[\mH, P_{t}]=0$ on $\fD$. 

\smallskip \noindent(ii): As $\mH$ is positive,  $\ker \mH =\{ 0 \}$ and $\fD \subset D(\mH)$, we get that $\ip{\cdot}{\cdot}_{\mH}$ is a positive definite inner product on $\fD$.

As $[\mH, P_{t}]=[\mH, P_{1}]=[\mH, P_{2}]=0$, we get a $\ip{\cdot}{\cdot}_{\mH}$-decomposition from $\fD=\fD_t \oplus \fDL \oplus \fDJ$.

\smallskip \noindent(iii): As $J$ is $\ip{\cdot}{\cdot}$-isometric on $\fD$ and $[J,\mH]=0$, it follows that $J$ is $\ip{\cdot}{\cdot}_{\mH}$-isometric on $\fD$.
Therefore $J$ extends isometrically to $\fL$ and it follows that this extension is $\ip{\cdot}{\cdot}_{\mH}$-unitary
as $J^2=\one$. Hence $(\fL, \iip{\cdot}{\cdot}_{\mH})$ is a Krein space by lemma \eqref{lm:JKsp}. Using $[\mH, J]=[\mH, P_{t}]=[\mH, P_{1}]=[\mH, P_{2}]=0$ on $\fD$, that $\fDL=J\fDJ$ and that $\fDL \perp \fDJ$ we get that $\fDL \perp_{\mH} \fDJ$ and hence $\fLL$ and $\fLJ$ are $\iip{\cdot}{\cdot}_{\mH}$-neutral spaces.
\end{proof}
The reason why we give $M_z$ a special symbol $Y$ is that we want to use the same notation for more general test function spaces below. We show that $\fL$ is the appropriate ghost test function space:
\begin{lemma}\label{lm:Tunit}
Let $T$ be defined as in Definition \eqref{df:T}. Then $T$ is an $\fL$-isometric isomorphism  (but is not $\fH$-isometric) with $\ran T$ dense in $\fLJ$. Hence it extends to a unitary $T:\fH_0 \to \fLJ$. Furthermore,
\begin{equation}\label{eq:Tinv}
(T^{-1}f)(p)={p_{\mu}f^{\mu}(p)},
\end{equation}
for all $f \in \fDJ$ and, 
\begin{align}
i(T^{-1}P_{\JJ}f)(p)=\,&ip_{\mu}f^{\mu}(p), \qquad \forall f\in \fD, \label{eq:sp2stuff}\\
i(\mH P_{\LL}JTh)_{\mu}(p)=\,& ip_{\mu}h(p)\qquad \forall h \in \fD_0, \qquad \mu=0,1,2,3. \notag
\end{align}
\end{lemma}
\begin{proof}
Let $h_1,h_2 \in \fB$, then using $p_{\mu}p_{\mu}=2\nps$ on $C_{+}$:
\begin{align*}
\ip{Th_1}{Th_2}_{\mH}&=\int_{C_{+}} \frac{d \bp^3}{p_0}\, \overline{Th_1(p)}\, { (\mH Th_2)(p)},\\
 &=\int_{C_{+}} \frac{d \bp^3}{p_0} z(p) \, \frac{2 \nps }{4 \np^{4}} \overline{h_1(p)}{h_2(p)},\\
 &=\int_{C_{+}} \frac{d \bp^3}{p_0} \overline{h_1(p)}\,{h_2(p)},\\
 &=\ip{h_1}{h_2}_0,
\end{align*}
and so $T$ is a $\fL$-isometric (but not $\fH$ isometric). 

We want to show that $\ran T$ is dense in $\fDJ$. By the definition of $T$ we know that $T(\fB)\subset \fDJ$ and we now show that $T(\fB)= \fDJ$. Let $f \in \fDJ$ and $g(p):=p_{\mu}f^{\mu}(p)$. Then as multplication by polynomials preserves Scwhartz space and commutes with the operators $M_{c_k}\in B(\fH_0)$ it follows that $g\in \fD_0$. Now as $y(p)=\np^{-1}$, $c_k(p)=p_k/\np$ and $p_0=\np$ for $p\in C_{+}$, we get that $M_y g= f_0+M_{c_k}f_k \in \fD_0$ where we have summed over $k=1,2,3$ and we use Theorem \eqref{th:fdtstdec} (iii). Therefore by the definition of $\fB$ we get that $g \in \fB$. Using equation \eqref{eq:P2def} it is immediate that $Tg=P_{\JJ}f=f$ and so $\fDJ \subset \ran T$, hence $T(\fB)=\fDJ$. 

As $T$ is $\fL$-isometric and $T(\fB)=\fDJ$ we get that it extends to a unitary $T:\fH_0 \to \fLJ$ and so $\overline{\ran T}=\fLJ$ in the $\norm{\cdot}_{\fL}$-topology.

The equation for the $T^{-1}$ can be checked by direct calculation of $T^{-1}Th=h$ for $h\in \fB$ using the definitions of $T$ and $T^{-1}$ and that $p_0^2=\nps$ for $p\in C_{+}$. 

Equations \eqref{eq:sp2stuff} follow as $p_{\mu}p_{\mu}=2\nps$ for $p \in C_{+}$ and for all $f\in \fD$ Proposition \eqref{lm:stdec2} gives
\begin{equation*}(T^{-1}P_{\JJ}f)(p)= p_{\mu}(P_{\JJ}f)^{\mu}(p)=\frac{2\nps}{{2}\nps}(p_{\nu}f^{\nu}(p))=p_{\nu}f^{\nu}(p)
\end{equation*}
and for all $h \in \fD_0$, $\mu=0,1,2,3$,
\[
i(\mH P_{\LL}JTh)_{\mu}(p)= i\frac{2\nps}{{2}\nps} p_{\mu} h(p)= ip_{\mu}h(p).
\]
where we used $JTh \in (J\fDJ=\fDL)$.
\end{proof}

We will now show how $T$ connects the heuristic scalar ghosts in Subsection \ref{sbs:hrext} with the Ghost Algebra $\cA_g(\fLJ)$. Recall that the heuristic ghost field was given by equation \eqref{eq:udef}:
\begin{equation*}
\gh(x)=\FC \int_{C_{+}} \frac{d^3p}{\sqrt{p_0}} ( c_2(\bp)e^{-ipx} + c_1(\bp)^{*} e^{ipx} ), , 
\end{equation*}
where $c_2(\bp),c_2(\bp)^*$ and $c_1(\bp),c_1(\bp)^*$ correspond to the creators and annihilators of distinct scalar fermionic fields  in $p$-space (cf. equations \eqref{h:car}). Accordingly, we would usually smear each of these scalar fields over a separate copy of $\fH_0=L^2(C_{+},\C,\lambda)$. However, by lemma \eqref{lm:Tunit} and $J\fLJ = \fLL$, we have the following unitary equivalences 
\[
\fLL \cong \fH_0 \cong \fLJ
\]
using the unitaries $JT:\fH_0 \to \fLL$ and $T:\fH_0 \to \fLJ$. So instead of using separate copies of $\fH_0$ to smear the fields over, we use $\fLL$ and $\fLJ$ and identify the heuristic fermionic fields $c_1(p),c_1(p)^*$ with the CAR algebra $CAR(\fLL)$ and $c_2(p),c_2(p)^*$ with $CAR(\fLJ)$. We combine these into the single fermionic field $CAR(\fLL \oplus \fLJ)=\cA_g(\fLJ)$. 
Explicitly, using equation \eqref{eq:udef} the smeared the heuristic ghost field $\gh(x)$ is:
\begin{align}
\int dx^4\, \gh(x) h(x)= &\;\FC \int_{C_{+}} \frac{d^3p}{\sqrt{p_0}} ( c_2(\bp)\hat{h}(p) + c_1(\bp)^{*} \overline{\hat{h}(p)}),\notag\\
=&\;\fst(c(f)+c^{*}(Jf)),\notag \\
=&\;C(f), \label{eq:hgsm1}
\end{align}
where  $h \in \cS(\R^4,\R)$ and $f=T (\hat{h}|_{C_{+}}) \in \fDJ$ and we use equation \eqref{eq:Cliffel} in the last equality. Hence we define the smeared scalar ghost field by:
\begin{align}\label{eq:cgsgf}
\tu:\fH_0 \to \cA_g, \qquad \text{by} \qquad \tu(h):=  \gh(Th)=C(Th),
\end{align}
(note that $Th \in \fLJ$)

Also recall in formal BRST-QEM the second ghost field $\tilde{\gh}(x)$ (equation \eqref{eq:halcg}) which is \emph{not} the conjugate ghost and does not anticommute with $\gh(x)$, given by the formula:
\begin{align*}
\tilde{\gh}(x)= & \;\FC \int_{C_{+}}\frac{d^3p}{\sqrt{p_0}} (- c_1(\bp)e^{-ipx} + c_2(\bp)^{*} e^{ipx} ), 
\end{align*} 
Smearing we obtain:
\begin{align}
\int dx^4\, \tilde{\gh}(x) h(x)= &\;\FC \int_{C_{+}} \frac{d^3p}{\sqrt{p_0}} (- c_1(\bp)\hat{h}(p) + c_2(\bp)^{*} \overline{\hat{h}(p)}),\notag\\
=&\;\left(\FC \int_{C_{+}} \frac{d^3p}{\sqrt{p_0}} ( c_2(\bp)\widehat{h}(p)  -c_1(\bp)^{*} \overline{\widehat{h}(p)})\right)^*,\notag\\
=&\;\left(\fst(c(f)-c^{*}(Jf))\right)^*,\notag \\
=&\;-i\left(\fst(c(if)+c^{*}(iJf))\right)^*,\notag \\
=&\;-iC(if)^* \label{eq:hgsm2}, 
\end{align}
where  $h \in \cS(\R^4,\R)$, $f=T (\hat{h}|_{C_{+}}) \in \fDJ$ and we used that $f\to c(f)$ is antilinear and $f\to c(f)^*$ is linear in the fourth equality.

Hence we define the smeared second ghost field as:
\begin{equation}\label{eq:cgsgf2}
\tilde{\tu}:\fH_0 \to \cA_g, \qquad \text{by} \qquad \tilde{\tu}(h):=-iC(iTh)^*=-i\tu(ih)^*,
\end{equation}
\begin{rem}\label{rm:concghsegh}
Note that $JTh \in \fLL$ for $h \in \fH_0$, hence $\tilde{\tu}(h)=-iC(iTh)^*=-iC(iJTh)=-i\cgh(iJTh)$ where $\cgh(iJh)$ is the conjugate ghost as defined in equation \eqref{eq:congfield}. Hence we see that we can recover the second ghost field from the conjugate ghost field as stated in Remark \eqref{rm:hgfcomrel} (iii).
\end{rem}
With the above definitions we get
\begin{proposition}\label{pr:consvgh}
The definitions of the ghost fields in equation \eqref{eq:cgsgf} and equation \eqref{eq:cgsgf2} give that for all $g,h \in \fH_0$:
\begin{itemize}
\item[(i)] $\tu(h)^{\dag}=\tu(h)$ and $\tilde{\tu}(h)^{\dag}=-\tilde{\tu}(h)$ .
\item[(ii)] $\osalg{\tu(h),\, \tilde{\tu}(h)\,|\, h \in \fH_0}=\cA_g(\fLJ)$.
\item[(iii)] The CAR's:\begin{align*}
\{\tu(h),\tu(g) \} = \{\tilde{\tu}(h), \tilde{\tu}(g) \}  = 0,\qquad \{\tu(h), \tilde{\tu}^*(g) \} = i \smp{h}{g}{0}\one,
\end{align*} 
where $\smp{\cdot}{\cdot}{0}=-\im \ip{\cdot}{\cdot}_0$, where we recall $\ip{\cdot}{\cdot}$ is the inner product on $\fH_0$..  
\end{itemize}
\end{proposition}
\begin{proof}
(i): These follow from $C(f)^{\dag}=C(f)$ (cf. equation \eqref{eq:Cliffel}) and the defining equation \eqref{eq:cgsgf} and equation \eqref{eq:cgsgf2}.

\pfit (ii): This follows by:
\begin{align*}
\osalg{\tu(h),\, \tilde{\tu}(h)\,|\, h \in \fH_0}&\;=\osalg{C(f),\, C(f)^*\,|\, f \in \fLJ},\\
&\;=\osalg{C(f),\, C(Jf)\,|\, f \in \fLJ},\\
&\;=\cA_g(\fLJ)
\end{align*}
where the first equality used that $T:\fH_0\to \fLJ$ is unitary, the second that $C(f)^*=C(Jf)$ for all $f \in \fL$ and the last that $J\fLJ=\fLL$ and $\fL=\fLL\oplus_{\mH}\fLJ$.

\pfit(iii): Using the CAR's given by equation \eqref{eq:Cliffcom}, equation \eqref{eq:cgsgf} and the $T$ is isometric gives for $g,h \in \fH_0$:
\[
\{\tu(h),\tu(g) \} =\{C(Th), C(Tg)\}= \re \iip{Th}{Tg}_{\mH} \one =  \re\ip{Th}{JTh}_{\mH} \one=0,
\]
where we used $Th \in \fLJ$, $JTh\in \fLL$ and $\fLL \perp_{\mH} \fLJ$ in the last equality. Similarly using $J^2=\one$ we get
\[
\{\tilde{\tu}(h), \tilde{\tu}(g) \} = - \re \ip{JTh}{Tg}_{\mH} \one=0
\]
Lastly we have
\[
\{\tu(h), \tilde{\tu}^*(g) \} = -i\{C(Th), C(iJTg)\}=-i \re \ip{Th}{iTg}_{\mH} \one= -i \im \ip{h}{g}_0\one= i \smp{h}{g}{0}\one,
\]
where we used that $T$ is isometric in the second equality.
\end{proof}
Note that $\smp{\cdot}{\cdot}{0}$ is the Fourier transform of the Pauli-Jordan distribution. That is, let $g,h \in \cS(\R^4,\R)$ and by a slight abuse of notation denote $\hat{g}|_{C_{+}},\hat{h}|_{C_{+}}$ by $\hat{g},\hat{h}$. Then 
\begin{align*}
\smp{\hat{h}}{\hat{g}}{0}= &\; \pi\int_{C_{+}} \frac{d^3\bp}{{p_0}}\left({\hat{h}(p)}\overline{\hat{g}(p)}-\overline{\hat{h}(p)}{\hat{g}(p)}\right) \\
=&\;\int \int dx^4 \,dy^4\, h(x)g(y)D_0(x-y).
\end{align*} 
where $D_0(x-y)$ is given by equation \eqref{eq:JPdist}. 

Hence, comparing Proposition \eqref{pr:consvgh} with the formal equation \eqref{eq:gcon} and Remark \eqref{rm:hgfcomrel} (iii) shows that the Ghost Algebra $\cA(\fLJ)$ has the correct anticommutation relations and ghost hermicity assignments when using the identifications given by equation \eqref{eq:cgsgf} and equation \eqref{eq:cgsgf2}. Therefore we are justified using $\cA(\fLJ)$ as the Ghost Algebra for QEM.

The formal smeared formula for $\drb$ in equation \eqref{eq:smAdrb} and equation \eqref{eq:smghdrb} will be connected below to the well-defined construction of QEM-BRST in Subsection \ref{sbs:dspv}. 

\section{Fock-Krein BRST}\label{sec:FKBRST}

\subsection{Fock-Krein CCR's}\label{sec:FKCCR}
We now give an account of the smeared Fock-Krein CCR's, which follows the treatment as given in \cite{Min1980, Hendrik2000}. 

Assume $\fH$ is a Krein space with Hilbert inner product $\ip{\cdot}{\cdot}$ indefinite inner product $\iip{\cdot}{\cdot}$ and fundamental symmetry $J$. The Bose-Fock space $(\fF^{+}(\fH), \ip{\cdot}{\cdot})$ of the Hilbert space $(\fH, \ip{\cdot}{\cdot})$ is a Krein space with respect to the indefinite inner product $\iip{\cdot}{\cdot}:=\ip{\cdot}{\Gamma_{+}(J)\cdot}$, cf. lemma \eqref{lm:JKsp}.

We define the creators and annihilators on the finite particle space $\fF_{0}^{+}(\fH) \subset \fF^{+}(\fH)$ with respect to the Krein inner product, as in \cite{Min1980}. That is we define the creation and annihilation operators as usual, except for the replacement of the Hilbert inner product $\ip{\cdot}{\cdot}$ on the 1-particle space $\fH$ with the Krein inner product $\iip{\cdot}{\cdot}$. On the symmetric $n$-particle space $\fH^n$ they are, for $f,h_1, \ldots, h_n \in \fH$:
\begin{align}
a^{\dag}(f)&S_n^{+} h_1 \otimes \cdots \otimes h_n = \sqrt{n+1} S_{n+1}^{+} f\otimes h_1 \otimes \cdots \otimes h_n, \label{eq:npcan}\\
a(f)&S_n^{+} h_1 \otimes \cdots \otimes h_n =\frac{1}{\sqrt{n}}\sum_{i=1}^{n} \iip{ f}{h_i}  S_{n-1}^{+}  h_1 \otimes \cdots \widehat{h_i} \cdots \otimes h_n, \notag
\end{align}
where $S^{+}_n$ is the symmetrization projection on $\fH^n$.
Note that on $\fF_0^{+}(\fH)$, $a^{\dag}(f)$ is the Krein adjoint of $a(f)$ (analogous to the usual $*$-adjoint case). 
\begin{lemma}
Let $a(f), a^{\dag}(f)\in \op(\fF_0^{+}(\fH))$ for $f,g \in \fH$ as above. Then for all $\psi \in  \fF_0^{+}(\fH)$:
\begin{itemize}
\item[(i)] We have, 
\begin{equation}\label{eq:adhadcon}
a^{\dag}(f)\psi=a(Jf)^{*}\psi, 
\end{equation}
where $a(f)^*$ is the adjoint of $a(f)$ with respect to the Hilbert inner product $\ip{\cdot}{\cdot}$ on $\fF_0^{+}(\fH)$. 
\item[(ii)] We have,
\begin{gather}\label{eq:acomm}
[a(f),a^{\dag}(g)]\psi=[a(f),a(Jg)^{*}]\psi=\ip{ f}{Jg}\psi= \iip{f}{ g}\psi,\\
[ a(f),a(g)]\psi=0. \notag
\end{gather}
\end{itemize}
\end{lemma}
\begin{proof}
\noindent (i): Let $b(f)=a(Jf)$ on $\fF_{0}^{+}(\fH)$. Then by equations \eqref{eq:npcan} and $\ip{\cdot}{\cdot}=\iip{\cdot}{J\cdot}$, $b(f)$ is the usual annihilator with respect to the Hilbert inner product on $\fH$. Hence $a^{\dag}(f)=b(f)^*=a(Jf)^*$ on $\fF_{0}^{+}(\fH)$.

\smallskip
\noindent (ii): These follow from (i) and the standard commutation relations for smeared Fock space creators and annihilators as in \cite{BraRob21981} p8-p10.
\end{proof}

We now define the Krein $\dag$-symmetric field operators as
\begin{definition}\label{df:Af}Let
\begin{align*}
A({f}):&= \frac{1}{\sqrt{2}} (a(f) + a^{\dag}(f) ), \qquad D(A(f))=\fF^{+}(\fH)_0\\
	&=\frac{1}{\sqrt{2}} (a(f) + a^{*}(Jf) ), \qquad D(A(f))=\fF^{+}(\fH)_0.
\end{align*}
The field algebra is then the non-normed $\dag$-algebra $\cA_{0,\fH}:= \alg{ A(f) \,: \, f \in \fH }$ which preserves the common dense domain $\fF^{+}(\fH)_0$.
\end{definition}

With these definition we get all the properties in \cite{Min1980} section 4, Theorem 1.
\begin{theorem}\label{th:ccrsa}
We have for $f,g \in \fH$, $\psi \in \fF^{+}(\fH)_0$:
\begin{enumerate}
\item The operator $A(f)$ is closable ,
\item $D(A(f))$ is a set of analytic vectors for $A(f)$,
\item \label{eq:stkocont} If $(f_j)\subset \fH$ and $\lim_{j\to \infty}f_j=f$ then,
\[
\lim_{j\to \infty}A(f_j)\psi=A(f)\psi, 
\]
\item The vacuum vector $\Omega \in \fF^{+}(\fH)_0$ is cyclic with respect to $\cA_{0,\fH}$.
\item $A(f)^*\psi=A(Jf)\psi$ \label{eq:Asad}
\item We have,
\[
[A(f),A(g)]\psi=i \smp{f}{g}{1}\psi,
\]
where $\smp{f}{g}{1}:=-\im \iip{f}{g}$ for all $f,g \in \fH$.
\end{enumerate}
\end{theorem}
\begin{proof}
The first four statements are proved in \cite{Min1980} Theorem 1. The fifth is immediate from the definition. The commutation relations \eqref{eq:acomm} give, 
\[
[A(f),A(g)]\psi=iIm\iip{f}{g}\psi=i \smp{f}{g}{1}\psi.
\]
\end{proof}
\begin{rem}\begin{itemize}
\item[(i)] From Theorem \eqref{th:ccrsa}, \eqref{eq:Asad} $A(f)$ is $*$-symmetric  if and only if $Jf=f$. 
\item[(ii)] If we let $\fH=L^2(C_{+},\mathbb{C}^4,\lambda)$ with inner products and symplectic form $\sigma_1$ as in Section \ref{sbs:testfunc}, then $\{A(f)\,|\, f \in \fXp\}$ corresponds to the smeared gauge potentials for QEM.
\end{itemize}
\end{rem}

\subsection{\KOB test function space}\label{sbs:abstf}
The BRST model we develop below will be able to be used for other abelian bosonic theories e.g.: the case of a finite number of Krein-symmetric bosonic constraints as in Subsection \ref{sbs:findimbos} and H\&T\cite{HenTei92} p313-316, and also for the case of Massive Abelian Gauge Theory as in Subsection \ref{sbs:Mabga} and \cite{Schf2001} p28.  We will prefer to use a general test function space in the $C^*$-algebraic theory (Chapter \ref{ch:CsBRST}). Hence we will not assume the specific QEM test function space as above, but instead use a general test function space with useful structures (already present on the QEM test function space). We refer to this as the \KO test funtion space.
\begin{enumerate}
\item Let $(\fX, \sigma_1)$ be a non-degenerate real symplectic space with symplectic form $\sigma_1$, and suppose there exists two idempotents, $P_{+}, P_{-}$ on $\fX$ such that
\[
P_{+}+P_{-}=\one, \qquad [P_{+}, P_{-}]=0, \qquad \smp{P_{+}\fX}{P_{-} \fX}{1}=0.
\]
We let $J=P_{+}-P_{-}$. Note that $J^2=\one$, and that $\smp{P_{+}\fX}{P_{-} \fX}{1}=0 \Rightarrow J \in \Sp(\fX, \sigma_1)$.
\item Let $\fX$ have the decomposition:
\[
\fX=\fX_t\oplus \fXL \oplus \fXJ= \fX_3 \oplus \fXJ 
\]
where $J\fX_t=\fX_t$, $\fXJ=J \fXL$ and $\fXL$ is the degenerate part of $\fX_3=\fX_t \oplus \fXL$. As $J$ is symplectic we get that, $\fXJ$ is the degenerate part of $J\fX_3= \fX_t\oplus \fXJ$.  Denote the algebraic projections onto $\fX_i$ as $P_i$ where $i=t,\LL,\JJ,3$. Note that $\oplus$ denotes the algebraic direct sum, and that $P_{\LL}=JP_{\JJ}J$. We will call $(\cX,\sigma_1)$ with the above structure the \emph{covariant symplectic space} and $\sigma_1$ the \emph{covariant symplectic form}. This terminology is motivated by the fact that the Poincar{\'e} transformations in subsection \eqref{sbs:hcov} are symplectic with respect to the corresponding QEM test function space.
\item Define a second symplectic form:
\begin{equation}\label{eq:gensmpdef}
\smp{\cdot}{\cdot}{2}:=\smp{\cdot}{J\cdot}{1}
\end{equation}
on $\fX$,  and note that $J$ is also $\sigma_2$-symplectic. We call a symplectic space with structure such as $(\fX, \sigma_2)$, where $\sigma_2$ is related to the \emph{covariant} form via $J$, the \emph{\ncss} and the form the \emph{\ncsf}.
\item Let $K$ be a complex structure of type $J$ (cf. \cite{Jak1985} p327 but where the roles of $J$ and $K$ are reversed), \ie $K$ is a symplectic operator on $\fX$ such that:
\begin{enumerate}
\item $K^2=-\one$ 
\item $[J,K]=0$ and for any $f\in \cX$, $\smp{f}{KJf}{1} \geq 0$, where equality holds if and only if $f=0$ \label{df:K2}
\end{enumerate}
Suppose further that $K \fX_i= \fX_i$ for $i=t,\LL,\JJ$. We define scalar multiplication and a complex IIP on $\fX$ by,
\begin{align*}
(\lambda_1 + i\lambda_2)f:= \lambda_1 f + \lambda_2Kf, \qquad \lambda_i \in \mathbb{R}, \; f \in \fX \\
\iip{f}{g}:= \smp{f}{Kg}{1} + i \smp{f}{g}{1}, \qquad f,g \in \fX
\end{align*}
We denote the complexified $\fX$ as $\fD$. As $[J,K]=0$ and $J\in \Sp(\fX,\sigma_1)$ and $J^2=\one$ on $\fX$ we have that $J$ is isometric and $J^2=\one$ on $\fD$. Due to property \eqref{df:K2} we see that,
\begin{equation}\label{eq:ipcst}
\ip{\cdot}{\cdot}:=\iip{\cdot}{J\cdot},
\end{equation}
is a positive definite inner product on $\fD$, and so $\fD$ is a pre-Hilbert space. Denote the closure of $\fD$ with respect to $\ip{\cdot}{\cdot}$ as $\fH$. As $J$ is isometric and $J^2=\one$ on $\fD$, it follows that $J$ extends to a unitary on $\fH$ such that $J^2=\one$. Hence by lemma \eqref{lm:JKsp}, $\iip{\cdot}{\cdot}$ extends to $\fH$ and  $(\fH, \iip{\cdot}{\cdot})$ is a Krein space. Note that the covariant and auxiliary symplectic forms extend to $\fH$ via $\smp{\cdot}{\cdot}{2}=\im \ip{\cdot}{\cdot}$ on $\fH$ and equation \eqref{eq:gensmpdef}.
\begin{rem}\label{rm:abstinp}
Note that for the QEM test function space we have a different sign convention in the correspondence between the symplectic form and inner product, \ie we have $\smp{\cdot}{\cdot}{1}=- \im \iip{\cdot}{\cdot}$. We have done this to satisfy the correspondence with the Jordan-Pauli distribution and to mantain the convention of the second argument in the inner product $\iip{\cdot}{\cdot}$ being linear. In the abstract case we have used the usual convention to avoid confusion in the many calculations in the following Chapters below. To get the explicit correspondence of the QEM test function space with the \KOB test function space, let $\sigma_{JP,1}=-\sigma_{KO,1}$ where the $JP$ subscript correspondends to the Jordan-Pauli symplectic form and $KO$ to the abstract $\KOB$ symplectic form.  
\end{rem}
Let $\fD_i:=(\fX_i+ K\fX_i)$ for $i=t,\LL,\JJ$.  
\begin{lemma} We have:
\begin{itemize}
\item[(i)] The decomposition:
\begin{equation}\label{eq:orthdsum}
\fD=\fD_t\oplus \fDL \oplus \fDJ
\end{equation} 
where $\oplus$ denotes $\ip{\cdot}{\cdot}$-orthogonal sum, and $J\fD_t=\fD_t$, $J\fDL=\fDJ$.
\item[(ii)] Let $\fH:=\overline{\fD}$ in the norm coming from $\ip{\cdot}{\cdot}$. Then 
\[
\fH=\fH_t\oplus \fHL \oplus \fHJ,
\]
where $\fH_{i}=\overline{\fD_i}$ where $i=t,\LL,\JJ$. Furthermore $P_{\LL},P_{\JJ}$ extend to $\ip{\cdot}{\cdot}$-orthogonal projections and $P_{\LL}=JP_{\JJ}J$.
\end{itemize}
\end{lemma}
\begin{proof}
(i): For $i=t,\LL,\JJ$, $K\fX_i\subset \fX_i$ implies $K\fD_i\subset \fD_i$, and  we get the decomposition in equation \eqref{eq:orthdsum} as an algebraic direct sum, and that $\fDL$ is the $\sigma_1$-isotropic part of $\fD_t\oplus\fDL$. To prove orthogonality, note $[K,J]=0$ implies $J\fD_t=\fD_t$ and $J\fDL=\fDJ$, hence,
\[
\smp{\fDL}{\fDJ}{2}=\smp{\fDL}{\fDL}{1}=0=\smp{\fDL}{\fD_t}{2}=\smp{\fDL}{\fD_t}{1},
\]
and $\smp{\fDJ}{\fD_t}{2}=\smp{\fDJ}{\fD_t}{1}=0$. As $K$ preserves $\fD_t,\fDL, \fDJ$ we have by the definitions of the inner products that $\fD_t\perp \fDL\perp\fDJ$ with respect to $\ip{\cdot}{\cdot}$. Hence $\oplus$ in equation \eqref{eq:orthdsum} is an $\ip{\cdot}{\cdot}$-orthogonal sum.

\smallskip \noindent (ii): The decomposition of $\fH$ is obvious. As $\fDL \perp \fDJ$ we get that $\fHL \perp \fHJ$ hence $P_{\LL},P_{\JJ}$ extend to the orthogonal projections on $\fHL$ and $\fHJ$. As $J$ is bounded we get $J\fHL=\fHJ$ by (i) and so $JP_{\LL}J=P_{\JJ}$ on $\fH$.
\end{proof}
Hence we we see that $\fX, \fD, \fH$ have the same decompositions as in QEM test function case (cf. Theorem \eqref{th:fdtstdec} ) with respect to the definite and indefinite inner products, and symplectic forms.
\item We define the ghost test function space as follows. 
\begin{definition}\label{df:mH}Let $\mH\in \op(\fH)$ be a possibly unbounded positive (w.r.t. $\ip{\cdot}{\cdot}$) and Krein (w.r.t. $\iip{\cdot}{\cdot}$ ) self-adjoint operator with domain $D(\mH)\subset\fH$, $\fD \subset D(\mH)$, $\ker \mH= \{0\}$, $\overline{\mH \fD}=\fH$. Suppose further that $[\mH, J]=[\mH, P_{t}]=[\mH, P_{1}]=[\mH, P_{2}]=0$ on $\fD$. Define a new inner product on $\fD$:
\[
\ip{\cdot}{\cdot}_{\mH}:=\ip{\cdot}{\mH \cdot}.
\]
\end{definition}
As $\mH$ is strictly positive,  $\ker \mH =\{ 0 \}$, and $\fD \subset D(\mH)$, we get that $\ip{\cdot}{\cdot}_{\mH}$ is a positive definite inner product on $\fD$. As  $[\mH, P_t]=[\mH, P_{\LL}]=[\mH, P_{\JJ}]=0$, we get that 
\[
\fD=\fD_t \oplus_{\mH} \fDL \oplus_{\mH} \fDJ= \fD_3 \oplus_{\mH} \fDJ
\]
where $\oplus_{\mH}$ signifies $\ip{\cdot}{\cdot}_{\mH}$-orthogonality. 
\begin{lemma} Let $\fL$ be the closure of $\fD$ with respect to $\ip{\cdot}{\cdot}_{\mH}$. Then:
\begin{itemize}
\item[(i)] $\fL$ is a Hilbert space and we have:
\[
\fL=\fL_t \oplus_{\mH} \fLL \oplus_{\mH} \fLJ= \fL_3\oplus_{\mH} \fLJ,
\]
where $\overline{\fD_i}:= \fL_i$ for $i=t,\LL,\JJ$ and the closure is with respect to the $\ip{\cdot}{\cdot}_{\mH}$-topology. 
\item[(ii)] $J$ extends to a unitary on $\fL$ such that $(\fL, \iip{\cdot}{\cdot}_{\mH})$ is a Krein space where
\[
\iip{\cdot}{\cdot}_{\mH}:=\ip{\cdot}{J \cdot}_{\mH}.
\]
\item[(iii)] We have that $\fLL,\fLJ$ are neutral subspaces of $\fL$ with respect to $\iip{\cdot}{\cdot}_{\mH}$.
\end{itemize}
\end{lemma} 
\begin{proof}
(i):Obvious as $\ip{\cdot}{\cdot}_{\mH}$ is positive definite on $\fD$.

\smallskip \noindent (ii): Now $J$ is $\ip{\cdot}{\cdot}$-isometric and $J^2=\one$ and $[J,\mH]=0$ on $\fD$, hence $J$ is $\ip{\cdot}{ \cdot}_{\mH}$-isometric and invertible on $\fD$, hence extends to a unitary on $\fL$. That $(\fL, \iip{\cdot}{\cdot}_{\mH})$ is a Krein space, follows from lemma \eqref{lm:JKsp}.

\end{proof}

We will use the same notation $A^*$ for the adjoint of an operator $A \in \op(\fD)$ with respect to the inner products $\ip{\cdot}{\cdot}$ and $\ip{\cdot}{\cdot}_{\mH}$, explicitly making note of the inner product where confusion may arise. Similarly for $A^\dag$ where $A \in \op(\fD)$.

We denote $P^{\fL}_i$ for $i=t,\LL,\JJ,3$ for orthogonal projections on $\fL$. Note that for $f \in \fD$ we have $P_if=P^{\fH}_if=P^{\fL}_if$, $i=t,\LL,\JJ,3$. Also as $[\mH, P_{t}]=[\mH, P_{1}]=[\mH, P_{2}]=0$ on $\fD$ we have that the $*$-adjoints of $P_i$ with respect to $\ip{\cdot}{\cdot}$ and $\ip{\cdot}{\cdot}_{\mH}$  coincide on $\fD$ for $i=t,\LL,\JJ,3$. Similarly for the $\dag$-adjoints with respect to $\ip{\cdot}{\cdot}$ and $\ip{\cdot}{\cdot}_{\mH}$. 
\end{enumerate}

A useful lemma for the future,
\begin{lemma}\label{lm:plmnJ}
Given $T \in \op(\fD)$ we have,
\[
T+JTJ=2(P_{+}TP_{+}+P_{-}TP_{-}), \qquad T-JTJ=2(P_{+}TP_{-}+P_{-}TP_{+}).
\]
\end{lemma}
\begin{proof}
Expand $\one$ and $J$ in terms of $P_{+}$ and $P_{-}$ and collect the terms.
\end{proof}
\subsection{BRST Extension}\label{sec:brstext}
We next define the BRST-extension of this system (as in Chapter \ref{ch:GenStruct}). Assume that we have a one-particle test function space $\fD$ with the structures in Subsection \ref{sbs:abstf}. Recall that $\fH=\overline{\fD}$ and take the Fock-Krein CCR algebra $\cA_{0,\fH}$ as in Subsection \ref{sec:FKCCR} with one particle Krein space $(\fH, \iip{\cdot}{\cdot})$. Take the Fock representation of the ghost algebra with one-particle space $(\fLL \oplus_{\mH} \fLJ, \iip{\cdot}{\cdot}_{\mH})$ (cf. Definiton \eqref{df:ghsp}). We are assuming the full ghost algebra here and so will have integer valued ghost spaces.

Let
\begin{gather*} 
\cA_0:=\alg{A(f)\,|\, f \in \fD}, \qquad \cD_0:= \fF^{+}_0(\fD),  
\end{gather*}
where $\fD \subset \fH$ and
\begin{gather*}
\cA_{g0}:=\alg{C(f)\,|\, f \in (\fDL\oplus_{\mH} \fDJ)}, \qquad D(G):= \fF^{-}_0(\fDL \oplus_{\mH} \fDJ), 
\end{gather*}
where $\fDL \oplus_{\mH} \fDJ \subset \fHL \oplus_{\mH} \fHJ$ and we have used notation of Section \ref{sbs:AbcuQ}. Then:
\[
D(Q)=\cD_0\otimes D(G)=\fF^{+}_0(\fD)\otimes \fF^{-}_0(\fDL \oplus_{\mH} \fDJ),
\]
and
\[
\cH=\overline{D(Q)}=\fF^{+}(\fH)\otimes \fF^{-}(\fLL \oplus_{\mH} \fLJ),
\]
where $\cD_0\otimes D(G)$ is the algebraic tensor product of $\cD_0$ and $D(G)$, \ie the linear span in of the elementary tensors in $\cH$. 

Let $\Omega:=\Omega_0\otimes \Omega_g$ where $\Omega_0$ and $\Omega_g$ are the cyclic vacuum vectors for $\fF^{+}(\fH)$ and $\fF^{-}(\fLL \oplus_{\mH} \fLJ)$ respectively. Denote the Hilbert inner product on $\cH$ by $\ip{\cdot}{\cdot}_{T}$ and Krein inner product by $\iip{\cdot}{\cdot}_{T}=\ip{\cdot}{J_{T}\cdot}_{T}$, where $J_T=\Gamma_{+}(J) \otimes \Gamma_{-}(J)$, where we have used the $\pm$ subscripts on the $\Gamma$'s to denote if they are the second quantization of $\fF^{\pm}(\fH)$. We use the subscript $T$ to avoid confusion with the one particle inner products. 

Then the BRST field algebra is the algebraic tensor product,
\[
\cA=\cA_0\otimes\cA_{g0}= \alg{A(f)\otimes \one,\, \one \otimes C(g)\, |\, f \in \fD,\, g\in (\fDL \oplus \fDJ)}\subset \cA_0 \otimes \cA_g,
\]
where which acts on $D(Q)$.  Then $\cA$ has the gradings described in Section \ref{sbs:AbcuQ}.
\begin{rem} Important points to keep track of for the remainder of Chapter \ref{ch:BRSTQEM} are:
\begin{itemize}
\item[(i)] The 1-particle indefinite inner product for the bosonic fields is $\iip{\cdot}{\cdot}$ while the 1-particle indefinite inner product on the ghost fields is $\iip{\cdot}{\cdot}_{\mH}=\iip{\cdot}{\mH \cdot}$.
\item[(ii)] For subalgebras $\cB \subset \cA_0$ and $ \cC \subset \cA_{g0}$, $\cB \otimes \cC$ denotes the algebraic tensor product and does not assume any topology. $C^*$-algebraic tensor products will be used in Chapter \ref{ch:CsBRST}. 
\end{itemize}
\end{rem}

\subsection{Superderivation}\label{sbs:dspv}
We now want to make a well-defined version of the BRST superderivation that models the formal example BRST-QEM as in Section \ref{sec:exheuEM}. To define the superderivation we proceed with a method similar to \cite{HendrikBuch2006} Section 4. First define a map $\drb$ on the generating elements of $\cA$ by:
\begin{definition}\label{df:brstsd}
We define a map on the generating elementary tensors of $\cA$ by:
\begin{align*}
\delta(A(g)\otimes \one)&= -i \one \otimes C( P_{\JJ}  ig), \qquad g \in \fD, \\
\delta(\one \otimes C(g))&= A(\mH P_{\LL} g)\otimes \one \qquad g \in (\fDL\oplus \fDJ),
\end{align*}
\end{definition}
This map extends to a superderivation on $\cA$ as proved below in Theorem \eqref{th:hbrstsd}. Before we proceed we have to show that $\drb$ defined above is a smeared version of the formal superderivation in the case of QEM. That is, where we identify the data $\fX,\sigma_i,K,J,\fD,\mH$, $i=1,2$ with the objects with the same labels as in Subsection \ref{sbs:QEMtf}. 

Recall that by Proposition \eqref{pr:consvgh}, we identified the Ghost Algebra $\cA_g(\fHJ)$ with the formal scalar ghosts via:
\begin{align*}
\tu(h)= \gh(Th)=C(Th),
\qquad \tilde{\tu}(h)=-i\tu(ih)^*=-iC(iJTh), 
\end{align*}
for all $h \in \fH_0=L^2(C_{+},\C,\lambda)$, where $T:=\fH_0 \to \fLJ$ is the unitary given in lemma \eqref{lm:Tunit}. Also recall for $f\in \fD$, lemma~\eqref{lm:Tunit} gives
\begin{equation*}
(iT^{-1}P_{\JJ}f)(p)=ip_{\nu}f^{\nu}(p) \quad a.e, \qquad i(\mH P_{\LL}JTh)_{\mu}(p)= ip_{\mu}h(p)\quad a.e, 
\end{equation*}
for all $f \in \fD$, $h \in \fD_0$ and $\mu=0,1,2,3$. 
Therefore, using the definition of $\drb$ (Definition \eqref{df:brstsd}) we calculate,
\begin{align*}
\drb (A(f) \otimes \one )&= -i\one \otimes C( P_{\JJ}  i(f))
=-i\tu(i T^{-1}P_{\JJ}(f))
=-i\tu(i p_{\mu}f^{\mu}(p))
\intertext{and for the ghost algebra}
\drb (\one \otimes \tu(h)) &= \drb  (\one \otimes C(P_{\JJ}Th))= A(\mH P_{\LL}  (P_{\JJ}Th)) \otimes \one = 0 \notag \\
\drb(\tilde{\tu}(h))&=-i\drb (\one \otimes \tu^*(ih))=-i \drb(\one \otimes C(iJTh)),\\ 
&= -iA(i\mH P_{\LL}(J P_{\JJ}Th)) \otimes \one = -iA(i\mH P_{\LL}^2 (JTh)) \otimes \one= -iA(i\mH J(Th)) \otimes \one \notag\\
&=-iA\left(ip_0{h}(p),ip_1{h}(p),ip_2{h}(p),ip_3{h}(p))\right)\otimes \one \notag
\end{align*}
Restricting to $\fX$ and $\fX_0\subset \fH_0$ and by a slight abuse of notation, denoting $\rho(\hat{f})$ and $\rho(\hat{h})$ by $\hat{f}$ and $\hat{h}$, for $f \in \cS(\mathbb{R}^4, \mathbb{R}^4)$,  $h \in \cS(\mathbb{R}^4, \mathbb{R})$:  
\begin{align*}
\drb (A(\hat{f}) \otimes \one )&= -i \one \otimes \tu(\widehat{ \partial_{\nu}f^{\nu}}) \\
 \drb (\one \otimes \tu(\hat{h})) &=0\\
 \drb (\one \otimes \tilde{\tu}(\hat{h})) &= -i A(\widehat{\partial_{0} h},\widehat{\partial_{1} h},\widehat{\partial_{2} h},\widehat{\partial_{3} h}) \otimes \one 
\end{align*}
If we compare the above with the formal superderivation $\drb$ defined by equations \eqref{eq:hdA}, \eqref{eq:hdgh}, \eqref{eq:hdcgh}: 
\begin{align*} 
\drb(A^\nu(x))&= -i\partial^\nu \gh(x), \\
\drb(\gh(x))&= 0, \\
\drb(\tilde{\gh}(x))& = -i\partial_\nu A^\nu(x), 
\end{align*}
and consider the ghost smearing formulas given in equation \eqref{eq:hgsm2} and equation \eqref{eq:hgsm1}, then we see that these are the correct smeared relations for $\drb$. Hence we see that Definition \eqref{df:brstsd} gives the correct mathematically well-defined relations. 
\begin{rem} 
Note that we could have used a scalar test function space in  Definition \eqref{df:brstsd}, but this will make calculations below much more cumbersome and obscure the fact the \KOB is really a constraint theory that uses test function spaces with structures such as $\fD$ with neutral subspace $\fDL$ which corresponds to the smearing functions for the constraint set. 
\end{rem}
Now that we have checked that the map $\drb:\cA \to \cA$ agrees with the formal superderivation for QEM on the generating tensors, we want to show that $\drb$ extends to a superderivation on all of $\cA$. So we return to the general context of the data $\fX,\sigma_i,K,J,\fD,\mH$, $i=1,2$ as in Subsection \ref{sbs:abstf}.
\begin{lemma}\label{lm:2nQs} Let $\vL=(f_j)_{j=1}^{n}$ be a finite $\fH$-orthonormal basis of a subspace $\cXs\subset \fDL$. Define
\[
\Qs:= \sum_{j=1}^{n} \left( A( f_j)\otimes C( Jf_j) + A(i f_j) \otimes C(i Jf_j) \right), \quad D(\Qs)=D(Q),
\]
Then:
\begin{itemize}
\item[(i)] $\Qs$ is Krein symmetric, and 2-nilpotent.
\item[(ii)] We have
\[
\Qs=\sum_{j=1}^{n}[a^*( Jf_j)\otimes c( Jf_j) + a( f_j)\otimes c^*( f_j)].
\] 
Furthermore $\Qs$ is independent of the choice of basis $(f_j)$ of $\cXs$.
\end{itemize}
\end{lemma}
\begin{proof}
(i): Krein symmetry is obvious. Now $\fDL \subset \cXs$ is a neutral space in both the respective indefinite inner product's on $\fH$ and $\fL$, so it follows that, 
\[
[A(f_j),A(f_k)]=[A(f_j),A(if_k)]=\{ C(Jf_j), C(Jf_k) \}= \{C(Jf_j), C(iJf_k) \}=0
\]
for $j,k=1,\ldots,n$. Note that this also implies that $C(Jf_j)^2=C(iJf_j)^2=0$. Using these we get,
\begin{align*}
\Qs^2&=\sum_{j,k=1}^{n} \left[ A( f_j)A(f_k)\otimes C( Jf_j)C(Jf_k) + A(i f_j)A(i f_k) \otimes C(i Jf_j)C(i Jf_k) \right.\\
&+ \left. A( f_j)A(if_k)\otimes C( Jf_j)C(iJf_k) + A(i f_j)A( f_k) \otimes C(i Jf_j)C( Jf_k) \right],\\
&=\sum_{1\leq k<j}^{n} \left( A( f_j)A(f_k)\otimes \{C( Jf_j),C(Jf_k)\} + A(i f_j)A(i f_k) \otimes \{C(i Jf_j),C(i Jf_k)\} \right) \\
& +\sum_{j=1}^{n}\left( A( f_j)A(f_j)\otimes C( Jf_j)^2 + A(i f_j)A(i f_k) \otimes C(i Jf_j)^2 \right) \\
& + \sum_{j,k=1}^{n} \left( A( f_j)A(if_k)\otimes \{C( Jf_j),C(iJf_k)\} \right),\\
&=0
\end{align*}
\smallskip \noindent (ii): Using $a(if)=-ia(f)$, $a^*(if)=ia^*(f)$, $c(if)=-ic(f)$, $c^*(if)=ic^*(f)$ and $[J,i]=0$ for all $f\in (\fDL\oplus \fDJ)$, then in terms of creators and annihilators we have,
\begin{align*}
\Qs&=\sum^{n}_{j=1}[A( f_j)\otimes C( Jf_j)  + A(i f_j)\otimes C(i Jf_j) ],\\
&=\frac{1}{2}\sum^{n}_{j=1}[(a( f_j)+a^*( Jf_j))\otimes (c( Jf_j)+c^*( f_j))   +(a(i f_j)+a^*(i Jf_j))\otimes (c(i Jf_j)+c^*(i f_j)) ],\\
&=\sum^{n}_{j=1}[a^*( Jf_j)\otimes c( Jf_j) + a( f_j)\otimes c^*( f_j)],
\end{align*}
As orthonormal bases of $\cXs$ are related by unitary matrices, it follows directly from the above formula that $\Qs$ is independent of the choice of basis $(f_j)_{j=1}^{n}$.
\end{proof}
Now we have the following superderviations on $\cA$:
\begin{lemma}\label{lm:Qndv}
Let $\vL=(f_j)_{j=1}^{n}$ be a finite $\fH$-orthonormal basis of $\cXs\subset \fDL$. Define the superderivaton $\drb_s$ by:
\[
\drb_s(A):= \sbr{\Qs}{A }, \qquad A\in \cA
\] 
Then we have:
\begin{itemize}
\item[(i)] For $g \in \fD$, $h \in (\fDL \oplus_{\mH} \fDJ)$ and $\psi \in D(Q)$ we get,
\begin{align*}
\delta_s(A(g)\otimes \one)\psi&=[\Qs,A(g)\otimes \one ]\psi= -i \one \otimes C( \sum_{j=1}^{n} i  \ip{Jf_j}{ g} Jf_j)\psi,  \\
\delta_s(\one \otimes C(h))\psi&=\{ \Qs,\one \otimes C(h) \}\psi= A(  \sum_{j=1}^{n}  \ip{f_j}{\mH h}f_j)\otimes \one\psi 
\end{align*}
\item[(ii)] For $g \in (\fD_t\oplus \cXs \oplus J\cXs)$, $h\in (\cXs \oplus_{\mH} J\cXs)$ we get,
\[
\delta_s(A(g)\otimes \one)= -i \one \otimes C( P_{\JJ}  ig), \qquad \delta_s(\one \otimes C(h))= A(\mH P_{\LL} h)\otimes \one,
\]
\ie in this case $\drb_s$ coincides with the map in Definition \eqref{df:brstsd}.
\end{itemize}
\end{lemma}
\begin{proof}
(i): These follow easily by writing everything in terms of creators and annihilators via lemma \eqref{lm:2nQs} and direct calculation. We give the calculation for $\delta_s(A(g)\otimes \one)\psi$. Let $\psi \in D(Q)$ and $g \in \fD$, then
\begin{align*}
[\Qs,A(g)]\psi= & \;\sum_{j=1}^{n} \left[ a^*(Jf_j)\otimes c(Jf_j)+a(f_j)c^*(f_j), \fst(a(g)+a^*(Jg))\right],\\
= & \;\fst \sum_{j=1}^{n} \left( [a^*(Jf_j),a(g)]\otimes c(Jf_j)+[a(f_j),a^*(Jg)]\otimes c^*(f_j)\right),\\
= & \;\fst \sum_{j=1}^{n}\left( -\ip{g}{Jf_j}\one \otimes c(Jf_j)+\ip{f_j}{Jg}\one \otimes c^*(f_j)\right),\\
&=-i\one \otimes \left(c(\sum_{j=1}^{n}i \ip{ Jf_j}{g}  Jf_j) + c^*( \sum_{j=1}^{n}i  \ip{ Jf_j}{g} f_j) \right)\psi \\
&=-i\one \otimes C(\sum_{j=1}^{n}i \ip{ Jf_j}{g}  Jf_j)  \psi \\
\end{align*}
where the fourth used $c(\lambda f)=\overline{\lambda}c(f)$ and $c^*(\lambda f)={\lambda}c^*(f)$, and the last used $C(f)=\fst(c(f)+c^{\dag}(f))=\fst(c(f)+c^*(Jf))$ for all $f\in (\fDL \oplus \fDJ)$.

The identity for $\delta_s(\one \otimes C(h))\psi$ follows similarly but with the $\mH$ factor appearing as $\{c(f),c^*(g)\}=\ip{f}{g}_{\mH}=\ip{f}{\mH g}$ for $f,g \in \fD$.

\smallskip \noindent(ii):
Let $g \in (\fD_t\oplus \cXs \oplus J\cXs)$. Then as $(f_j)_{j=1}^n$ is an orthonormal basis for $\cXs$ we get
\[
\sum_{j=1}^{n}i \ip{ Jf_j}{g}  Jf_j=J(\sum_{j=1}^{n}i \ip{ f_j}{Jg}f_j)=iJP_{\LL}Jg=iP_{\JJ}g.
\]
and for $h \in (\cXs \oplus J\cXs)$ we get $\sum_{j=1}^{n}  \ip{f_j}{\mH h}f_j=P_{\LL}\mH h= \mH P_{\LL} h$ as $[\mH,P_{\LL}]=0$. The result follows from (i).
\end{proof}
Using this we have that 
\begin{theorem}\label{th:hbrstsd}
There is a superderivation $\drb:\cA \to \cA$ with respect to the $\Z_2$-grading on $\cA$ as described in Section \ref{sbs:AbcuQ}, which coincides with $\drb$ in Definition \eqref{df:brstsd} on the generating $G:=\{A(f)\otimes \one,\one \otimes C(g)\,|\, f\in \fD, g\in \fDL\oplus_{\mH} \fDJ\}$. Furthermore $\drb^2=0$.
\end{theorem}
\begin{proof}
Let $\cXs \subset \fDL$ be any finite dimensional subspace and consider the subalgebra,
\[
\cA(\cXs):=\salg{A(g),\, C(h)\, |\, g \in \fD_t\oplus \cXs \oplus J\cXs,\, h\in \cXs \oplus J\cXs}.
\] 
Every elementary tensor $A$ is in $\cA(\cXs)$ for some $\cXs$. Hence lemma \eqref{lm:Qndv} (ii) gives that $\drb_s(A)=\drb(A)$ for all $A\in G$ and this $\cXs$. As elements $\cA$ are finite polynomials in the elementary tensors, we see that $\drb$ given by definition \eqref{df:brstsd} extends to a superderivation on $\cA$.

Also $\Qs^2=0$ (cf. lemma \eqref{lm:Qndv} (i)) implies $\drb^2=0$.
\end{proof}
\begin{rem}  
 Note that the appearance of the $\mH$ appears in $\delta_s(\one \otimes C(g))\psi$, is due to the fact that on the ghost space we are using the inner product $\ip{\cdot}{\cdot}_{\mH}$.

\end{rem}

\subsection{BRST charge}\label{sbsc:brch}
We want to construct $Q$ such that $Q$ generates $\drb$. Once this is done we can connect to results in Section \ref{sec:dsprig} such as the $dsp$-decomposition, and Section \ref{sbs:AbcuQ} such as the extended algebra $\cAe$ and Theorem \eqref{pr:krdel}. First a useful lemma,
\begin{lemma}\label{lm:thbas}
Let $\cXs$ be a finite dimensional subspace of $\fDL$. Then there exists an $\fH$-orthonormal basis $(f_j)_{j=1}^{\MM}$ of $\cXs$ where $m=\dim(\cXs)$ such that $ \{h_j:=  f_j /\norm{ f_j}_{\fL}\,|\, j=1,\ldots,\MM\}$ is a $\fL$-orthonormal basis of $ \cXs$,
\end{lemma}
\begin{proof}
As $\mH$ is positive, we have that $\ip{\cdot}{\cdot}_{\mH}$ defines a positive quadratic form (cf. \cite{BraRob21981} p27) on $\cXs$. Therefore we have that there exists a positive self adjoint operator $T \in B(\cXs)$ such that 
\begin{equation}\label{eq:thbas}
\ip{ g}{f}_{\mH}=\ip{ g}{\mH f}= \ip{f}{Tg},
\end{equation}
for $f,g \in \cXs$. As $\ker \mH= \{0\}$, we get that $\ker T=0$. Now $T$ is a self adjoint operator on a finite dimensional space $\cXs$ and so has a complete set of orthonormal eigenvectors $(f_j)$ that span $\cXs$. By \eqref{eq:thbas} we see that $(f_j)$ is the basis we are looking for.
\end{proof}

\begin{rem} Note that as $\mH$ does not in general preserve the finite dimensional subspace $\cXs$, $T$ in the above proof need not be equal to $\mH$. 
\end{rem}
Now define:
\begin{lemma}\label{lm:dfpscs}
Let $\cXs$ be a finite dimensional subspace of $\fDL$, let $\vL=(f_j)_{j=1}^{n}$ be an $\fH$-orthonormal basis of $\cXs$, and let,
\[
\cP^{\cXs}_{\vL}=\alg{ a^*(g)\otimes\one ,a^*(f_j)\otimes \one,a^*(Jf_j)\otimes \one, \one \otimes c^*(f_j), \one \otimes c^*(Jf_j) | g \in \fD_t},
\] 
by the algebra generated by the above creation operators. Furthermore let,
\[
\cS^{\cXs}_{\vL}:=\cP^{\cXs}_{\vL}\Omega
\] 
where $\Omega$ is the cyclic vacuum vector for $D(Q)$. Then
\begin{itemize}
\item[(i)]  $\cP^{\cXs}_{\vL}$ and $\cS^{\cXs}_{\vL}$ do not depend on the choice of $\vL$, and so we will only use the $\vL$ subscript when we want to emphasize the basis we are using in a proof.
\item[(ii)]  For all $\psi \in D(Q)$ there exists a finite dimensional subspace $\cXsp\subset\fDL$ such that $\psi \in \cS^{\cXsp}$. Moreover
\[
D(Q)=\sn{\cup \cS^{\cXsp} \,|\, \text{finite dimensional subspaces $\cXsp \subset \fDL$} }
\]
\item[(iii)] $J_T \cS^{\cXs}=\cS^{\cXs}$ where $J_T$ is the fundamental symmetry on $\cH=\overline{D(Q)}$. 
\end{itemize}
\end{lemma}
\begin{proof}
\pfit (i): Follows as $\cXs$ is finite dimensional and $f\to a^*(f)$ and $g\to c^*(g)$ are linear.

\pfit (ii): Let $\psi \in D(Q)=\fF_0^{+}(\fD)\otimes\fF_0^{-}(\fDL\oplus_{\mH}\fDJ)$. Then $\psi=A\Omega$ where $A$ is a polynomial of a finite number of creators $a^*(f_j)\otimes \one$ and $\one \otimes c^*(g_k)$ where $f_j \in \fD$, $g_k \in (\fDL\oplus_{\mH} \fDJ)$, $j,k\in \Z^{+}$. As $\fD=\fD_{t}\oplus \fDL \oplus \fDJ$ we can choose  a  finite dimensional subspace $\cXsp\subset\fDL$ such that $\fD_t \oplus \cXsp \oplus J\cXsp$ contains all the arguments of the creators and annihilator in the polynomial $A$. Hence $A\in \cP^{\cXs}$ and $\psi=A\Omega \in \cS^{\cXsp}$.
 
\pfit(iii): We have $\Gamma_{-}(J) c^*(f) \Gamma_{-}(J)=c(Jf)$ for $f\in (\fDL\oplus_{\mH} \fDJ)$, $\Gamma_{+}(J) a^*(f) \Gamma_{+}(J)\psi=a(Jf)\psi$ for $f \in \fD$ and $\psi \in \fF^{+}_0(\fD)$ and $J_T=\Gamma_{+}(J)\otimes \Gamma_{-}(J)$. As $J(\fD_t \oplus \cXs \oplus J\cXs)=\fD_t \oplus \cXs \oplus J\cXs$ it follows that $J_T\cP^{\cXs}J_T=\cP^{\cXs}$, hence as $J_T\Omega=\Omega$ we get that $J_T \cS^{\cXs}=\cS^{\cXs}$.
\end{proof}
We want to define $Q$ so that it is independent of a basis of $\cXs$. 
\begin{lemma}\label{lm:Qsind}
Let $\cXs \subset \cXsp$ be finite dimensional subspaces of $\fDL$ then:
\begin{itemize}
\item [(i)] For all  $\psi \in \cS^{\cXs}$
\[
\Qs \psi = \Qsp  \psi
\]
\item[(ii)] $\Qs \cS^{\cXs} \subset \cS^{\cXs}$ and $\Qs^* \cS^{\cXs} \subset \cS^{\cXs}$.
\end{itemize}
 
\end{lemma}
\begin{proof}
(i): Take $\psi \in \cS^{\cXs}$, and suppose that $\dim(\cXs)=m< \dim(\cXsp)=n$. Now by lemma\eqref{lm:thbas}, take an $\fH$-orthonormal basis $\vL=(f_j)_{j=1}^{m}$ of $\cX_{s}$ such that it is also a $\fL$-orthogonal basis for $\cX_{s}$, and take an  $\fH$-orthonormal basis $\vL''=(f_j)_{j=m+1}^{n}$ of $\cXsp\ominus \cXs$ such that it is also a $\fL$-orthogonal basis for $\cXsp\ominus \cXs$.  Therefore $\vL'=(f_j)_{j=1}^{n}$ is an $\fH$-orthonormal basis and $\fL$-orthogonal basis for  $\cXsp$. Now,
\begin{align}\label{eq:Qbind}
(\Qs-\Qsp )&=\sum_{j=m+1}^{n}[a^*( Jf_j)\otimes c( Jf_j) + a( f_j)\otimes c^*( f_j)].
\end{align}
By the way we chose $\vL'$ we have that $f_i \perp_{\fH} f_j$ and $f_i \perp_{\fL} f_j$ for $i\leq n$, $j>n+1$. So as $\psi=A \Omega$ where $A \in  \cP^{\cXs}$, we see that we can (anti-)commute all the terms in the RHS of \eqref{eq:Qbind} through the terms in $A$ to $\Omega$, which they annihilate. Therefore
\[
(\Qs-\Qsp )\psi=0
\]

\pfit (ii): Let $\psi  \in \cS^{\cXs}$. Then $\psi=A\Omega$ where $A \in \cP^{\cXs}$, and so using the (anti)commutation relations for $a(f)$ and $c(g)$ (cf. \eqref{eq:acomm} and \eqref{eq:CARac}) and that the annihilators annihilate $\Omega$, we have $\Qs\psi \in \cS^{\cXs}$. Also we have
\[
\Qs^*\cS^{\cXs}=J_T\Qs^{\dag}J_T\cS^{\cXs}\subset J_T\cS^{\cXs}=\cS^{\cXs}
\]
where we have used lemma \eqref{lm:dfpscs} (iii) and $\Qs^{\dag}=\Qs$ on $\cS^{\cXs}$ (cf. lemma\eqref{lm:2nQs} (i)).
\end{proof}
Using this lemma we construct the BRST charge $Q$.
\begin{theorem} \label{th:HQ1}
Let $\psi \in D(Q)$, let $\cXs$ be a finite dimensional subspace of $\fDL$ such that $\psi \in \cS^{\cXs}$, and define
\[
Q\psi:= \Qs \psi. 
\]
Then $Q$ extends to well defined operator with domain $D(Q)$ and:
\begin{itemize}
\item[(i)] $Q$ preserves $D(Q)$, is Krein symmetric (hence closable by Proposition \eqref{pr:Krclos}) and $Q^2\psi=0$ for all $\psi \in D(Q)$. 
\item[(ii)] $\cQ$ is Krein symmetric (hence closable) and preserves $D(\cQ)$ and $\cQ^2\psi=0$ for all $\psi \in D(\cQ)$. Hence $\cQ$ has \emph{dsp}-decomposition (cf. Theorem \eqref{th:Hdsp1})
\[
\cH=\cH_d\oplus \cH_s \oplus \cH_p,
\]
where $\cH_d= \overline {\ran \cQ}$, $\cH_s = \ker \cQ \cap \ker Q^* $, $\cH_p= \overline{\ran Q^*}$.
\item[(iii)] $(Q^*)^2\psi=0$ for all $\psi \in D(Q^*)$.
\end{itemize}
\end{theorem}
\begin{proof}
Let $\psi \in D(Q)$, then by lemma \eqref{lm:dfpscs} (ii) we have there exists a finite dimensional subspace $\cXs \subset \fDL$ such that $\psi \in \cS^{\cXs}$, and so $Q\psi= \Qs \psi$. We check that $Q\psi=\Qs \psi$ is independent of the choice of this $\cX_s$. Suppose that $\cXsp\subset \fDL$ is another finite dimensional subspace such that $\psi \in \cS^{\cXsp}$ and suppose that $\cXs,\cXsp\subset\cXspp $ where $\cXspp$ is a finite dimensional subspace of $\fDL$. Then by lemma \eqref{lm:Qsind}
\[
\Qs\psi=\Qspp \psi=\Qsp \psi,
\]
and so $Q\psi=\Qs \psi$ is independent of the choice of this $\cX_s$. That $Q$ extends to a well defined linear operator on $D(Q)$ is now obvious. 
 
\pfit(i): Let $\psi \in D(Q)$. Then $\psi \in \cS^{\cXs}$ for some finite dimensional subspace of $\cXs \subset \fDL$.  So we have that $Q\psi=\Qs\psi \in \cS^{\cXs}$, hence $Q D(Q)\subset D(Q)$. Furthermore by lemma \eqref{lm:Qsind} (ii) we have that $Q^2\psi=\Qs^2\psi=0$ hence $Q^2\psi=0$, and $Q$ is $\dag$-symmetric as $\Qsp $ is for all finite dimensional subspace $\cXsp\subset \fDL$. Also, $\cQ$ is closable by Proposition \eqref{pr:Krclos}.

\pfit (ii) and (iii) now follow from lemma \eqref{lm:QessaCsa} and Theorem \eqref{th:Hdsp1}. 
\end{proof}
\subsection{State space}\label{sbs:statespace}
By Theorem \eqref{th:HQ1} (ii) $\cQ$ and $\cH$ satisfy the conditions of the \emph{dsp}-decomposition (Theorem \eqref{th:Hdsp1}). We want to be able to calculate $\cH_s$ explictly since $\cH_s \cong \cH^{BRST}_{phys}$ via the natural isomporphism $\vp$ in lemma \eqref{lm:phspksp} (ii). Also to use Theorem \eqref{pr:krdel}, we need to extend $\cA$ to $\cAe$ to include the above projections, etc, as in lemma \eqref{lm:projdspinc}. For this, we need the following decomposition of $D(Q)$:
\begin{lemma}\label{lm:hdecDq} 
\begin{itemize} \item[(i)] $D(Q)$ has the decomposition
\[
D(Q)= (\fF_0(\fD_t)\otimes \{\mathbb{C}\Omega\} )\oplus \ran Q \oplus \ran Q^*.
\]
\item[(ii)] We have
\begin{gather*}
P_sD(Q)=(\fF_0(\fD_t)\otimes \{\mathbb{C}\Omega\})\subset D(Q), \quad P_dD(Q)= \ran Q\subset D(Q),\\
P_pD(Q)= \ran Q^* \subset D(Q),
\end{gather*}
where $P_i$ is the projection onto $\cH_i$ for $i=s,p,d$.
\end{itemize}
\end{lemma}
\begin{proof}
(i): Let $\vp \in D(Q)$, then by lemma \eqref{lm:dfpscs} (ii) there exists a finite dimensional subspace $\cXs\subset \fDL$ such that $\vp=A\Omega$ where $A \in \cP^{\cXs}$. Hence it suffices to show that for all monomials of generating tensors $A_2 \in \cP^{\cXs}$ we have
\[
\psi=A_2\Omega\in (\fF_0(\fD_t)\otimes \{\mathbb{C}\Omega\} )\oplus \ran Q \oplus \ran Q^*
\]
We drop the tensor product $\otimes$ notation for the remainder of the proof.

\pfit \textbf{Case 1:} 

Now by lemma \eqref{lm:thbas}, we can choose an $\fH$-orthonormal basis $\Lambda=(f_j)_{j=1}^m$ of $\cXs$ such that $(f_j)_{j=1}^m$ is also a $\fL$-orthogonal basis of $ \cXs$. Suppose that
\begin{align*}
A_2&=\\
 & a^*(g_{i_1})
\ldots a^*(g_{i_j})a^*(f_{k_1})\ldots  a^*(f_{k_l})a^*(Jf_{p_1})\ldots a^*(Jf_{p_q}) c^*(f_{r_1})\ldots c^*(f_{r_s})c^*(Jf_{t_1})\ldots c^*(Jf_{t_u}),\\ &\neq 0
\end{align*}
where
\begin{itemize}
\item $ g_{i_1},\ldots,g_{i_j}\in \fD_t$,  
\item $f_{k_1},\ldots, f_{k_l}, f_{p_1},\ldots,f_{p_q},f_{r_1},\ldots,f_{r_s},f_{t_1},\ldots,f_{t_u} \in  (f_j)_{j=1}^m$. Note that the arguments of the $a^*(\cdot)$'s may be repeated, but that the arguments of the $c^*(\cdot)$'s may not as the CAR's imply $(c^*(h))^2=0$ for all $h \in (\fDL\oplus_{\mH} \fDJ)$.
\item $(l+q+s+u)=n \in \N$ and $n \neq 0$,
\item We use the convention that $l,q,s,u=0$ means that there are no corresponding creator terms in the monomial $A_2$.
\end{itemize}
That is $A_2$ is a monomial of `$n$-unphysical creators'.

Now motivated by the formal calculation of $\Delta$ in equation \eqref{eq:horrnd} we calculate as follows, where we rely heavily on the (anti)commutation relations for $a(f)$ and $c(g)$ (cf. \eqref{eq:acomm} and \eqref{eq:CARac}) and that the $a(f)$'s and $c(g)$'s commute:
\begin{align}
\Qs(\Qs)^*+ (\Qs)^*\Qs = 
 \sum_{j,i=1}^{m} \{a^*(Jf_i)c(Jf_i)+a(f_i)c^*(f_i), a(Jf_j)c^*(Jf_j)+a^*(f_j)c(f_j)\}
\end{align}
Now $\{a^*(Jf_i)c(Jf_i), a^*(f_j)c(f_j)\}=a^*(Jf_i)a^*(f_j)\{c(Jf_i),c(f_j)\}=0$ and hence,
\begin{align}
\Qs(\Qs)^*+ &(\Qs)^*\Qs = \;\sum_{j,i=1}^{m} \large(\{a^*(Jf_i)c(Jf_i), a(Jf_j)c^*(Jf_j)\} +\{a(f_i)c^*(f_i), a^*(f_j)c(f_j)\}\large) \notag\\
=&\;\sum_{j,i=1}^{m} ( a^*( J f_i)a( J f_j)\{c( J f_i), c^*(  J f_j) \} + [a( J f_j), a^*(  J f_i)] c^*( J f_j) c( J f_i)) \notag \\
&+( a^*(  f_j)a(  f_i)\{c(  f_j), c^*(   f_i) \} + [a(  f_i), a^*(   f_j)] c^*(  f_i) c( J f_j)) \notag \\
= & \;\sum_{j=1}^{m}( \norm{f_j}^2_{\fL}[a^*( J f_j)a( J f_j)+ a^*( f_j)a( f_j)] + [c^*(J f_j) c(J  f_j) + c^*( f_j) c( f_j)] )\label{eq:kerD}
\end{align}
where the second equality can be checked by expanding the respective brackets and the third follows from the (anti)commutation relations for $a(f)$ and $c(f)$ (cf. \eqref{eq:acomm} and \eqref{eq:CARac}), the $\fH$-orthonormality and $\fL$-orthogonality of $(f_j)_{j=1}^m$, and the unitarity of $J$. 

Now by the definition of $Q$ (Theorem \eqref{th:HQ1}) and lemma \eqref{lm:Qsind} we  have that $\psi, Q\psi=\Qs\psi \in \cS^{ \cXs}$. By lemma \eqref{lm:dfpscs} (ii) and (iii) it follows that $J_TD(Q)=D(Q)$. Hence, by Theorem \eqref{th:HQ1} (ii) and lemma \eqref{lm:khad} it follows that $Q^*=J_TQ^{\dag}J_T=J_TQJ_T$ on $D(Q)\supset \cS^{\cXs}$. Furthermore, $J_T \cS^{ \cXs}\subset \cS^{ \cXs}$ (cf. lemma \eqref{lm:dfpscs} (iii)), hence we have that $Q^*\psi= J_TQJ_T\psi \in \cS^{ \cXs}$. We calculate,
\begin{align}
&(Q(Q)^*+ (Q)^*Q)\psi=(\Qs(\Qs)^*+ (\Qs)^*\Qs)\psi \notag \\ 
&= \left(\sum_{j=1}^{m}  \norm{f_j}^2_{\fL}[a^*( J f_j)a( J f_j)+ a^*( f_j)a( f_j)]\right)\psi + \left(\sum_{j=1}^{m} [ c^*(J f_j) c(J  f_j) + c^*( f_j) c( f_j)]\right)\psi \notag \\
&= \left(\sum_{j=1}^{m}  \norm{f_j}^2_{\fL}[a^*( J f_j)a( J f_j)+ a^*( f_j)a( f_j)]\right)\psi + \left(\sum_{j=1}^{m} \norm{f_j}^2_{\fL}[ c^*(J h_j) c(J  h_j) + c^*( h_j) c( h_j)]\right)\psi, \label{eq:Delrigop}
\end{align}
where $(h_j)_{j=1}^m=(f_j/\norm{f_j}_{\fL})_{j=1}^m$ is a $\fL$-orthonormal basis of $\cX$. 

Note that the above formula looks similar to a number operator but with the $\norm{f_j}^2_{\fL}$ factors. This motivates the following calculation. Let $f_{v} \in (f_j)_{j=1}^m$ and suppose that $A_2$ has $w$ creators $a^*(f_{v})$ where $w \in \mathbb{Z}^{+}$. Then
\begin{equation}\label{eq:asded1}
\norm{f_v}^2_{\fL}a^*( f_v)a( f_v)\psi=\norm{f_v}^2_{\fL}a^*( f_v)a( f_v)A_2\Omega=w\norm{f_v}^{2}_{\fL}A_2\Omega=w\norm{f_k}^{2}_{\fL}\psi,
\end{equation}
where we calculated similar as for number operators in the last equality using the CCR's (cf. equation \eqref{eq:acomm}). We also have
\begin{equation}\label{eq:asded2}
\norm{f_k}^2_{\fL}c^*( h_k)c( h_k)\psi=\norm{f_k}^2_{\fL}c^*( h_k)c( h_k)A_2\Omega=\norm{f_k}^{2}_{\fL}A_2\Omega=\norm{f_k}^{2}_{\fL}\psi,
\end{equation}
where we used the CAR's (cf. equation \eqref{eq:CARac}) $c^*(f_k)=\norm{f_k}^2_{\fL}c^*(h_k)$, that $(h_j)=(f_j/\norm{f_j}_{\fL})$ is an orthonormal basis of $\fL$ and that $A_2$ can have at most one $c^*(f_k)$ term with the identity being trivial when there is no $c^*(f_k)$.

Using equation \eqref {eq:asded1} and equation \eqref{eq:asded2}, the similar equations for the orthonormal basis $(Jf_j)_{j=1}^{m}$ of $J\cXs$, and that $n>0$ and $\norm{f_j}^2_{\fL}>0$ we can calculate using equation \eqref{eq:Delrigop},
\[
(Q(Q)^*+ (Q)^*Q)\psi= C_{\psi} \psi.
\]
where $C_\psi>0$ is some positive constant. Therefore we have,
\begin{equation}\label{eq:kerran}
\psi=(1/C_{\psi} )(Q(Q)^*+ (Q)^*Q)\psi \subset \ran Q \oplus \ran Q^*
\end{equation}

\pfit \textbf{Case 2:} 

Suppose $A_2\in \cP^{\cXs}$ is a monomial of $a^*(g)$'s with arguments in $g\in \fD_t$, then $\psi=A_2\Omega \in (\fF_0(\fD_t)\otimes \{\mathbb{C}\Omega_g\})$. By \eqref{eq:kerD} we also see that $(\fF_0(\fD_t)\otimes \{\mathbb{C}\one\})\subset \ker (QQ^*+ Q^*Q) \perp (\ran Q \oplus \ran Q^*)$.

Combining Case 1 and Case 2 gives that for any monomial $A_2\in \cP^{\cXs}$ we have that 
\[
A_2\Omega\in (\fF_0(\fD_t)\otimes \{\mathbb{C}\Omega\} )\oplus \ran Q \oplus \ran Q^*,
\]
hence, as dicsussed at the beginning of the proof, $D(Q)$ has the decomposition we require.

\smallskip \noindent (ii): Immediate from (i).
\end{proof}
\begin{rem}
\begin{itemize}\item [(i)] The calculation in equation \eqref{eq:Delrigop} is the rigorous version of the formal equation \eqref{eq:horrnd}. This is the equation that shows that the ghosts are removed in the final physical space $\cH^{BRST}_{phys}$ in BRST-QEM case, and is the crucial calculation in showing that BRST gives the correct physical results for QEM below, e.g. Theorem \eqref{th:kerdrand}. 
\item[(ii)] Also $D(Q)= \ran Q \oplus \cH_s\oplus \ran Q^*$  is assumed to follow in a slightly different context in \cite{HoVo93} Proposition 6 p203 in  a very similar manner to lemma \eqref{lm:hdecDq}. The difference with the approach taken in \cite{HoVo93} is that only representations of the algebra $l(1,1):=\alg{Q,Q^*,G,\Delta}$ is considered rather than the whole algebra $\cA$ which will correspond to the smeared gauge potentials in the BRST-QEM example.
\end{itemize}
\end{rem}

Using this we can calculate $\cH_s$.
\begin{corollary}\label{cr:kerran} 

\begin{itemize}
\item[(i)] We have,
\begin{align*}
\cH_s &=\fF(\fH_t) \otimes \mathbb{C}\Omega_g, \qquad (\cH_d \oplus \cH_p)= \cH \ominus(\fF( \fH_t) \otimes \mathbb{C}\Omega_g),
\end{align*}
\item[(ii)] If $\vp:(\ker \cQ=(\cH_s \oplus \cH_d)) \to (\cH^{BRST}_{phys}=\ker \cQ/ \cH_d)$ is the factor map then,
\begin{align*}
\cH^{BRST}_{phys}\cong \cH_s=\fF(\fH_t) \otimes \mathbb{C}\Omega_g,
\end{align*}
where the above isometric isomorphism is given by $\vp|_{\cH_s}:\cH_s \to \cH^{BRST}_{phys}$. Moreover,
\[
\cD^{BRST}_{phys}:=\vp(\ker Q) \cong \fF_0(\fD_t) \otimes \mathbb{C}\Omega_g,
\]
where the above isometric isomorphism is also given by $\vp|_{\cH_s}$.
\item[(iii)] If $J|_{\fD_t}=\one$ then $J_TP_s=P_s$ where $J_T=\Gamma_{+}(J)\otimes \Gamma_{-}(J)$ is the fundamental symmetry on $D(Q)$. That is the physicality condition given by equation \eqref{eq:spphcond} holds.
\end{itemize}
\end{corollary}
\begin{proof}
(i): From lemma \eqref{lm:hdecDq} we have,
\[
\cH_s \cap D(Q) = \fF_0(\fD_t) \otimes \mathbb{C}\Omega_g,
\]
and,
\[
(\cH_d \oplus \cH_p) \cap D(Q) = D(Q) \ominus (\fF_0( \fD_t) \otimes \mathbb{C}\Omega_g).
\]  
As $\cH_s$ and $(\cH_d \oplus \cH_p)$ are closed we have
\begin{align*}
\cH_s &\supset \overline{\fF_0(\fD_t) \otimes \mathbb{C}\Omega_g}=\fF(\fH_t) \otimes \mathbb{C}\Omega_g\\
(\cH_d \oplus \cH_p) &\supset \overline{D(Q) \ominus (\fF_0( \fD_t) \otimes \mathbb{C}\Omega_p)} = \cH \ominus(\fF( \fH_p) \otimes \mathbb{C}\Omega_g),
\end{align*}
and so, as $\cH_s \perp (\cH_d \oplus \cH_p)$,
\begin{align*}
\cH_s &=\fF(\fH_t) \otimes \mathbb{C}\Omega_g, \qquad (\cH_d \oplus \cH_t)= \cH \ominus(\fF( \fH_t) \otimes \mathbb{C}\Omega_g).
\end{align*}

\pfit (ii): Obvious.

\pfit (iii): If $J|_{\fH_t}=\one$ we have that $J_T=(\Gamma_{+}(J)\otimes \Gamma_{-}(J))=\one$ on $\cH_s=(\fF(\fH_t) \otimes \Omega_g)$. Therefore $J_TP_s=P_s$.
\end{proof}

\subsection{Physical algebra}\label{sbs:phalopem}
By lemma \eqref{lm:hdecDq} and lemma \eqref{lm:projdspinc} we can extend the $\mathbb{Z}_2$-grading from $\cA$ to $\cAe:= \salg{\cA \cup\{Q,K,P_p,P_d, P_s\}}$ and extend the domain $\drb$ to  $D(\drb)=\cAe$ as done in Section \ref{sbs:AbcuQ}. To investigate the BRST operator cohomology, recall that from Theorem \eqref{pr:krdel} we have an algebra isomorphism: 
\[
\cP^{BRST}:=(\ker \drb \cap \cA) /(\ran \drb \cap \cA)\cong \Phi^{ext}_s(\ker \drb \cap \cA).
\]
To calculate this we need make some definitions. Let,
\[
\cA(\fM)=\alg{  \one, a(f),a^*(g): f,g \in \fM  }, 
\]
and $\fM$ is a subspace of $\fH$. Note that $\cA(\{0\})= \{ \mathbb{C}\one \}$. 
Now recall that $\cA_0=\cA(\fD)$, $\cA_{g0}= \salg{ C(f)\,|, f \in \fDL\oplus_{\mH}  \fDJ}$ and $\cA=\cA_0\otimes \cA_{g0}$ (cf. Subsection \ref{sec:brstext}). Define
\begin{align*}
\cA_b:= \cA(\fDL\oplus\fDJ), \qquad
\cA_{ph}:= \cA(\fD_t), \qquad \cA_{u}:= \cA_b \otimes \cA_{g0},
\end{align*}
and note that $\cA_b  \cA_{ph}= \cA_0$, and $\cA= (\cA_{ph}\otimes \one) \cA_{u}$, since $\fD=\fD_t\oplus \fDL \oplus \fDJ=\fD_t\oplus_{\mH} \fDL \oplus_{\mH} \fDJ$.

Now let,
\[
\cB:=\alg{ a(f)\otimes \one , \one \otimes c(g) : f,g \in \fDL \oplus \fDJ  }, 
\]
Note that $\cB$ is not a $*$-algebra and $\one \notin \cB$. Furthermore, it is not hard to see that
\[
\cA_u=\sn{\one+\cB+\cB^* + \cB^*\cB } 
\]
\begin{rem} We need to include $\one$ in RHS above, otherwise we can't (anti)commute elements and we don't get a $*$-algebra. This is not done in \cite{HoVo92} p1321 for the 1-cell case, and is why the authors have the result that $\ker \drb = \ran \drb$ which is Theorem 1 of that paper. This is not to say that the result in \cite{HoVo92} is wrong, just that the domain of $\drb$ is in that paper is not the whole field algebra $\cA$. 
\end{rem}
\begin{lemma}\label{lm:indFh}
We have:
\begin{itemize}
\item[(i)] Let $A_1 \in \cA_{ph}\otimes \one$ and  $A_2 \in \cA_{u}$, and  $\Omega=\Omega_0 \otimes \Omega_g \in \fF^{+}(\fH) \otimes \fF_g$. Then
\[
\ip{\Omega}{A_1A_2\Omega}=\ip{\Omega}{A_1\Omega}\ip{\Omega}{A_2\Omega}
\]
\item[(ii)] For all $T \in \cA$ there exists $S \in \cA_{ph}\otimes \one$ such that,
\[
\ip{\Omega}{R_1(T-S) R_2\Omega} )=0, \qquad \forall R_1, R_2 \in \cA_{ph}\otimes \one
\]
\end{itemize}
\end{lemma}
\begin{proof}
(i): First note that by normal ordering, we can write $A_2= c \one + C$ where $c \in \mathbb{C}$, $C \in (\cB+\cB^* + \cB^*\cB)$. Now $\cB \Omega= 0$, so we get that $A_2 \Omega= c \Omega + C_{2} \Omega$ where $C_{2} \in \cB^*$, and 
\[
\ip{\Omega}{A_2\Omega}=c+\ip{\Omega}{C_{2}\Omega}=c+\ip{C_{2}^*\Omega}{\Omega}=c.
\]
Now we calculate,
\begin{align*}
\ip{\Omega}{A_1A_2\Omega}= & \;\ip{\Omega}{A_1(c + C_{2})\Omega},\\
= & \; \ip{\Omega}{A_1\Omega}c +\ip{\Omega}{A_1C_{2}\Omega},\\
= & \;\ip{\Omega}{A_1\Omega}\ip{\Omega}{A_2\Omega}+\ip{C_{2}^*\Omega}{A_1\Omega},\\
= & \;\ip{\Omega}{A_1\Omega}\ip{\Omega}{A_2\Omega}
\end{align*}
where we used that $C_{2} \in \cA_u$ and $[\cA_u,\cA_{ph}\otimes \one]=0$ on $D(Q)$.

\pfit (ii): Let $T=\sum_{i=1}^{n} a_i A_i B_i \in \cA=\cA_{ph}\cA_u$, where $a_i \in \C$. $A_i \in \cA_{ph}\otimes \one$ and $B_i \in \cA_u$.  Then we have that for all $ R_1, R_2 \in \cA_{ph}\otimes \one$,
\begin{align*}
\ip{\Omega}{R_1TR_2\Omega}= & \;\sum_{i=1}^{n}\ip{\Omega}{R_1a_iA_iR_2B_i)\Omega},\\
= & \;\sum_{i=1}^{n}\ip{\Omega}{R_1a_iA_iR_2\Omega}\ip{\Omega}{B_i\Omega}.
\end{align*}
where we have used that $[\cA_u,\cA_{ph}\otimes \one]=0$ on $D(Q)$ for the first line and (i) for the second. Now if we let $b_i= \ip{\Omega}{B_i\Omega}$ and $S=\sum_{i=1}^{n} a_i b_i A_i \in \cA_{ph}\otimes \one$, then it is clear from the above calculation that
\[
\ip{\Omega}{R_1(T-S)R_2\Omega}=\sum_{i=1}^{n}\ip{\Omega}{R_1a_iA_iR_2\Omega}\ip{\Omega}{(B_i-b_i)\Omega}=0.
\]
As all elements in $\cA$ are of the form of $T$ we are done.
\end{proof}
From the above lemma
\begin{theorem}\label{th:kerdrand} We have,
\[
\cP^{BRST}\cong \Phi^{ext}_s(\ker \drb \cap \cA)=\cA_{ph}\otimes \one.
\]
where the above algebra isomorphism is as given in Theorem \eqref{pr:krdel}.
\end{theorem}
\begin{proof}
Take $T=\sum_{i=1}^{n} a_i A_i B_i \in (\ker \drb\cap \cA)$, $A_i \in \cA_{ph}$ $B_i \in \cA_u$. Then by lemma \eqref{lm:indFh} (ii) we have that there exists $S \in \cA_{ph}$ such that, 
\[
\ip{\Omega}{R_1(T-S) R_2\Omega} =0, \qquad \forall R_1, R_2 \in \cA_{ph}.
\]
By Corollary \eqref{cr:kerran} $\fH_s= \overline{\cA_{ph} \Omega}$, hence by the above identity $P_sTP_s=P_sSP_s$. As $S$ is a polynomial in creators and annihilators over $\fD_t$, we get that $\delta(S)=0$ by the definition of $\drb$ (cf. Definition \eqref{df:brstsd}) and $P_{\LL}\fD_t=P_{\JJ}\fD_t=\{0\}$.

Hence
\[
\Phi^{ext}_s(\ker \drb \cap \cA)=P_s (\ker \drb \cap \cA)P_s= P_s(\cA_{ph}\otimes \one)P_s=\cA_{ph}\otimes \one,
\]
where the last equality follows from $P_sD(Q) = \fF_0(\fD_t)\otimes\{ \mathbb{C} \one \} $ and $\cA_{ph}:= \cA(\fD_t)$. By Theorem \eqref{pr:krdel} the result follows.
\end{proof}

\subsection{Example: Finite \KOB }\label{sbs:findimbos}
The following example is \KOB for a finite number of bosonic constraints. The corresponding test function spaces for the fields and the ghosts for this situation are as follows: 
\begin{itemize}
\item Let 
\[
\fD=\fD_t\oplus \fDL \oplus \fDJ
\]
be an abstract test function space with all the structures as in Subsection \eqref{sbs:abstf}, and with $\dim(\fDL)=\MM$ where $\MM < \infty$. It is easy to see that such test function spaces exist if we recall the construction in lemma \eqref{lm:ghtsfext}. In particular, $\fD$ has positive definite inner product $\ip{\cdot}{\cdot}$ giving $\fH_{j}=\overline{\fD_{j}}$ for $j=t,\LL,\JJ$ and
\[
\fH=\fH_t\oplus \fHL \oplus \fHJ
\]
Also there exists a fundamental symmetry $J\in B(\cH)$ such that $(\cH,\iip{\cdot}{\cdot})$ is a Krein space, where $\iip{\cdot}{\cdot}:=\ip{\cdot}{J\cdot}$. Furthermore $J\fDL=\fDJ$ and $J \fD_t=\fD_t$.
\item For the ghost test function space let $\mH=\one$ and so as $\dim{\fDL}=\dim{\fDJ}=\MM < \infty$ we have $\fL=\overline{\fDL\oplus \fDJ}=\fDL\oplus \fDJ$ where $\fL$ is defined as in Subsection \eqref{sbs:abstf}.
\end{itemize} 
Using the constructions and results of previous sections we construct the BRST extension, superderivation and charge  and calculate  $\cH^{BRST}_{phys}$ and easily $\cP^{BRST}$ for this $\fD$.
\begin{proposition}\label{pr:fdkob}
Let $\fD$ and $\fL=\fDL\oplus \fDJ$ be as above. Then:
\begin{itemize}
\item [(i)]Let $(f_j)_{j=1}^{\MM}$ be an orthonormal basis for $\fDL$. Then
\begin{align*}
Q:= & \; \sum_{j=1}^{\MM} \left( A( f_j)\otimes C( Jf_j) + A(i f_j) \otimes C(i Jf_j) \right), \quad D(\Qs)=D(Q),\\
=& \; \sum_{j=1}^{\MM}[a^*( Jf_j)\otimes c( Jf_j) + a( f_j)\otimes c^*( f_j)]
\end{align*}
\item[(ii)] We have $\cH^{BRST}_{phys}\cong (\fF(\fH_t)\otimes \C \Omega_g)$ where the isomorphism is given as in Corollary \eqref{cr:kerran}, and  $\cP^{BRST}\cong (\cA_{ph}\otimes \one)$ where the isomorphism is as in Theorem \eqref{th:kerdrand}.
\item[(iii)] If $\fD_t=\{0\}$ then $\cH^{BRST}_{phys}\cong  \C \Omega$ and $\cP^{BRST}\cong \C\one$, \ie we have a trivial theory.
\end{itemize}
\end{proposition}
\begin{proof}
(i): As $\fDL$ is finite dimensional we let $\cXs=\fDL$ and it follows from Theorem \eqref{th:HQ1} that $Q=\Qs$. The first formula for $Q$ is then the defining formula for $\Qs$, the second comes from lemma \eqref{lm:2nQs} (ii).

\pfit (ii): Follows by Corrollary \eqref{cr:kerran} and Theorem \eqref{th:kerdrand}. 

\pfit (iii): Let $\fD_t=\{0\}$. Then $\fH_t=\{0\}$ and so by (ii) $\cH^{BRST}_{phys}\cong \fF(\fH_t)\otimes \C \Omega_g = \C \Omega$. Futhermore, $\cA_{ph}=\cA(\fD_t )=\cA(\{ 0 \} )= \mathbb{C}\one $ and so (ii) gives  $\cP^{BRST}\cong \cA_{ph}\otimes \one = \C \one$.
\end{proof}
\begin{rem}\label{rm:fdkob}
\begin{itemize}\item[(i)] Proposition \eqref{pr:fdkob} (iii) (\ie $\fD_t=\{0\}$) corresponds to H\&T\cite{HenTei92} p313-316 (though the operator algebra is not calculated here), and \cite{HoVo92} though with the different result that $\one \notin \ran \delta$. This also answers in the affirmative the conjecture at the end of \cite{HoVo92}, \ie the \emph{dsp}-decomposition can be used to efficiently calculate the physical algebra $\cP^{BRST}$ (cf. Remark \eqref{rm:ext} (ii)).
\item[(ii)] Also note that while the formula for $Q$ looks cosmetically very similar to that given for abelian Hamiltionian BRST there are key differences:
\begin{enumerate}\item We have that $\ker A(f_j)=\{0\}$, and $\ker A(if_j)=\{0\}$ for $j=1 \ldots \MM$, hence we get that $\cH^p_0=\{0\}$ if we were to assume that these were constraints. This is why the Gupta-Bleuler approach is used in the Dirac constraint version of QEM. 
\item The $A(f)$'s are Krein symmetric, but \emph{not} Hilbert symmetric. This effects the calculation of $\ker \Delta$ in equation \eqref{eq:calhdelham} and this is why we do not have the result that $\ker \Delta= \cH^p_{0}\otimes \cH_g= \{0\}$ as may be expected (wrongly) from the abelian Hamiltionian BRST example (c.f equation \eqref{eq:hcomhamdel}).
\item We did not use the Berezin representation of the restricted ghost algebra. We used the Ghost algebra that comes from smearing the CAR creators and annihilators over the appropriate test function space.
\end{enumerate}
\end{itemize}
\end{rem}
We see that finite \KOB differs from Hamiltonian BRST in several significant ways.

\section{Electromagnetism}\label{sc:EM}
Now that we have constructed the appropriate well-defined BRST structures, we substitute the QEM test function space $\fD$ and ghost test function space $\fL$ (cf. Subsection \eqref{sbs:QEMtf} and Subsection \eqref{sbs:QEMGHtf}) to calculate the BRST physical space and BRST physical algebra for BRST-QEM and compare this to Gupta-Bleuler QEM as done in  $\cite{Hendrik2000}$ Section 5.6.

Recall for QEM: 
\[
\fH_t =\{ f \in \fH \,| \, \mathbf{p.f(p)}=0, f_0(p)=0 \},
\]
and $\fD_t=\fD\cap \fH_t$ (cf. \eqref{eq:tsHidec}). 
\begin{proposition}\label{pr:QEMres} We have that:
\begin{itemize}
\item[(i)] The physical state space for BRST-QEM is:
\begin{align*}
(\cH^{BRST}_{phys}=\ker Q/ \cH_d)\cong (\cH_s=\fF(\fH_t) \otimes \mathbb{C}\Omega_g),\qquad \cD^{BRST}_{phys}\cong \fF_0(\fH_t) \otimes \mathbb{C}\Omega_g,
\end{align*}
where the above isomorphisms are as Corrollary \eqref{cr:kerran}.
\item[(ii)] 
The physical algebra for BRST-QEM is:
\[
\cP^{BRST}\cong P_s (\ker \drb \cap \cA) P_s= \cA_{ph} \otimes \one=\cA(\fD_t)\otimes \one. 
\]
where the isomorphism as in Theorem \eqref{th:kerdrand}.
\item[(iii)] The physicality condition $J_TP_s=P_s$ holds where $J_T$ is the fundamental symmetry for the Krein space $(\cH, \iip{\cdot}{\cdot}_T)$ (cf. Definition \eqref{df:spphcond} and Subsection \ref{sec:brstext}). Hence the Hilbert and Krein inner products coincide on $\cH^{BRST}_{phys}$, \ie $\iip{\cdot}{\cdot}_p=\ip{\cdot}{\cdot}_p$. Moreover, the Krein $\dag$-adjoint on $\cP^{BRST}$ corresponds to the Hilbert $*$-adjoint.
\end{itemize}
By the expression for $\fH_t$ above, the BRST physical space and algebra correspond to the transversal photons and have Hilbert inner product and adjoint respectively.
\end{proposition}
\begin{proof}
(i):  By Corollary \eqref{cr:kerran} (ii).

\pfit (ii): By Theorem \eqref{th:kerdrand}.

\pfit (iii): By equation \eqref{eq:tsHidec} we have that $J\fH_t=\fH_t$ hence by Corollary \eqref{cr:kerran} (iii) the physicality condition $J_TP_s=P_s$ holds. By lemma \eqref{lm:phspksp} (iii) the inner product is positive definite and by lemma \eqref{lm:facsal} (ii) the Krein $\dag$-adjoint on $\cP^{BRST}$ corresponds to the Hilbert $*$-adjoint.
\end{proof}

We compare this to QEM using Gupta-Bleuler constraints as done in $\cite{Hendrik2000}$ Section 5.6. A technicality in making this comparison is that $\cite{Hendrik2000}$ defines the unconstrained fields on Fock-Krein space with 1-particle test function space $\fDp$ rather than $\fD \supset \fDp $ that we have been using (see Subsection \ref{sbs:testfunc}). We assume that this technicality is only a minor difference due to  $\overline{\fD}=\overline{\fDp}=\fH$, and the fact that by Theorem \eqref{th:ccrsa} \eqref{eq:stkocont} gives $f \to A(f)\psi$ is continuous for all $\psi \in \fF_0^{+}(\fH)$. The results and proofs in $\cite{Hendrik2000}$ Section 5.6 are virtually identical using $\fD$ as 1-particle space and we now give a summary using our notation.

Let $\cC=\{ \chi(h):=a(f)\,|\, h \in \cS(\mathbb{R}^4,\mathbb{R}), \, f_{\mu}(p)=ip_{\mu}\widehat{h}(p)\}$, where we are identifying $f$ with its image in $\fH=L^2(C_{+},\mathbb{C}^4,\mu)$ as done in Section \ref{sbs:testfunc}.  That is, $\cC$ is the set of smeared Gupta-Bleuler constraints as in \cite{Sch64} p246. Let
\[
\cH'=\{ \psi \in \fF_0(\fD)\,|\, \psi \in D(\chi(h))\quad \text{and} \quad \chi(h)\psi=0 \,\forall h \in \cS(\mathbb{R}^4,\mathbb{R})\}
\]
Let $\cH''$ be the neutral part of $\cH'$ with respect to the indefinite inner product $\iip{\cdot}{\cdot}_0$ on $\cD_0=\fF_0(\fD)$. If we recall  $\fD_3=\fD_t \oplus \fDL$ (cf. Theorem \eqref{th:fdtstdec}) then,
\begin{proposition}\label{pr:phguspa} We have,
\begin{itemize}
\item [(i)] $\cH'= \fF_0(\fD_3)$
\item [(ii)] $ \cH''=\{\psi \in \fF_0(\fD_3)\,|\, \psi^{(n)}\in S_n(\fDL \otimes \fD_3 \otimes \ldots \otimes \fD_3)\}$
\item [(iii)] Let $P_t$ be the projection on $\fD_t$ and let $\hat{\psi}=\vp{\psi}$ for $\psi \in \cH'$ where $\vp:\cH'\to \cH'/\cH''$ is the factor map. Then 
\begin{align*}
U:\;\cH'/\cH'' \to \fF_0(\fD_t), \qquad  U\hat{\psi}=\Gamma(P_t)\psi
\end{align*}
for all $\psi \in \cH'$ is a well defined isomorphism. 
\end{itemize}
\end{proposition}
\begin{proof}
By Proposition \eqref{lm:stdec2} and Theorem \eqref{th:fdtstdec} (i),
\begin{align}
\fD_3&=\fD_t\oplus \fDL=\fD\cap \fH_3 =\{ f \in \fD\,: \, p_\mu f^{\mu}(p)=0 \}\notag,\\
\fDL&= \{ f \in \fD \,: \, f_{\mu}(p)=(ip_{\mu})h(p)\; \text{for}\; h\in \fH_0 \},\notag
\end{align}
(i) and (ii): The Proposition follows from Proposition 5.19 and Proposition 5.21 in $\cite{Hendrik2000}$.

\pfit (iii): As $P_t|_{\fD_t}=\one$ it follows that $\Gamma_{+}(P_t)|_{\fF_0(\fD_t)}=\one$, and as $P_t \fDL=0$ it follows that
\[
\cH'' \subset \ker \Gamma_{+}(P_t)\qquad\text{and}\qquad \Gamma_{+}(P_t)\fF_0(\fD_t)=\fF_0(\fD_t)\subset \ran \Gamma_{+}(P_t)
\]
Now $P_t^*=P_t$ and $\Gamma_{+}(P_t)^*=\Gamma_{+}(P_t)$ and so by the above we see $\cH'' \perp \fF_0(\fD_t)$. Now from $\fD_3=\fD_t \oplus \fDL$ it follows that
\[
\cH'=\fF_0(\fD_3)=\fF_0(\fD_t)\oplus \cH''
\]
From the above we see that $\Gamma_{+}(P_t)|_{\cH'}$ is the projection on $\fF_0(\fD_t)\subset \cH'$ and so the map $U$ is a well defined isomorphism.
\end{proof}
The physical state space in the Gupta-Bleuler approach is $\cH'/\cH''\cong \fF_0(\fD_t)\cong \cD^{BRST}_{phys}$ and so we see that the BRST phsycial state space is naturally isomorphic.

We now define and calculate the physical algebra for the Gupta-Bleuler approach.
\begin{theorem}\label{th:GBalg}
The observables of $\cA_0=\alg{A(f) \,|\, f \in \fD}$ in the Gupta-Bleuler approach is the subalgebra of $\cO^{GB}\subset \cA_0$ that factors to $\cH'/\cH''$ \ie 
\[
\cO^{GB}:=\{ A \in \cA_0 \,|\, A\cH'\subset \cH',\, A\cH''\subset \cH'' \}.
\]
The trivial observables are
\[
\cD^{GB}:=\{ A \in \cO^{GB} \,|\, A\cH'\subset \cH''\},
\]
The natural physical representation of $\cO^{GB}$ is defined by
\[
\pi_{GB}:\cO^{GB}\to \op(\cH'/\cH''), \qquad \pi_{GB}(A)\hat{\psi}=\widehat{A\psi}
\]
for all $\psi \in \cH'$. 
With these definitions we have:
\begin{itemize}
\item[(i)] $\ker \pi_{GB} =\cD^{GB}$. Hence $\cD^{GB}$ is a two-sided ideal.
\item[(ii)] The following algebra isomorphisms
\[
\cO^{GB}/\cD^{GB} \cong \pi_{GB}(\cO^{GB}) \cong \cA(\fD_t)
\]
where  we are assuming no topology.
\end{itemize}
\end{theorem}
\begin{proof}
(i): Obvious.

\pfit(ii): First, $\cO^{GB}/\cD^{GB} \cong \pi_{GB}(\cO^{GB})$ is obvious from (i).
 
Next recall from the end of the proof of Proposition \eqref{pr:phguspa} (iii) that $P_{Dt}:=\Gamma_{+}(P_t)|_{\cH'}$ is the projection on $\fF_0(\fD_t)\subset \cH'$. Define
\[
\Phi_{GB}:\cO^{GB} \to \op( \fF_0(\fD_t)) \qquad \text{by} \qquad \Phi_{GB}(A)=P_{Dt} A P_{Dt}
\]
If we replace $\ker \drb$ by $\cO^{GB}$, $\cH_s$ by $\fF_0(\fD_t)$, $\cH_p$ by $\cH''$, and $\Phi^{ext}_s$ by $\Phi_{GB}$ in the proof of lemma \eqref{lm:alkerhom} we see that $\Phi_{GB}$ is an algebra homomorphism. 

From the proof of  Proposition \eqref{pr:phguspa} (iii) we know
\begin{equation}\label{eq:Hddecmp}
\cH'=\fF_0(\fD_3)=\fF_0(\fD_t)\oplus \cH''.
\end{equation}
Moreover $P_{Dt}\cD^{GB}\cH'\subset P_{Dt} \cH'' =\{0\}$ and so it follows $\cD^{GB}\subset \ker\Phi_{GB}$. Conversely, let $A \in \ker \Phi_{GB}$. We want to show $A \cH' \subset \cH''$. By the definition of $\cD^{GB}$ we know $A \cH'' \subset \cH''$. So given equation \eqref{eq:Hddecmp} above we need only show $A\fF_0(\fD_t) \subset \cH''$. But as $A \in \ker \Phi_{GB}$ we have $A\fF_0(\fD_t)=AP_{Dt}\cH' \perp P_{Dt}\cH'=\fF_0(\fD_t)$ hence $A\fF_0(\fD_t) \subset \cH''$ hence $A \cH' \subset \cH''$ hence $A \in \cD^{GB}$. Therefore we have that
\[
\cD^{GB} =\ker\Phi_{GB}.
\]
But $P_{Dt}|_{\fF_0(\fD_t)}=\one$ we see that $\Phi_{GB}(\cA(\fD_t))=\cA(\fD_t)$, and so we have that 
\[
\cO^{GB}/\cD^{GB} \cong \Phi_{GB}(\cA(\fD_t))=\cA(\fD_t)
\]
where the above is an algebra isomorphism. Hence the statement of (ii) follows.
\end{proof}

If we compare Theorem \eqref{th:GBalg} (ii) and Proposition \eqref{pr:QEMres} (ii) we see Gupta-Bleuler algebra, $\cO^{GB}/\cD^{GB}$, is isomorphic to the BRST-physical algebra above. So the BRST approach and Gupta-Bleuler approach to QEM give equivalent results at the level of the CCR unbounded field operators acting on Fock-Krein space.

\begin{rem} 
Different choices of $\mH$ can give equivalent results for the final physical objects but do necessarily give the same smeared $\drb$. To see this take the QEM test function space $\fD$ but let $\mH=\one$, hence we have ghost test function space $\fL=\fHL\oplus \fHJ$ as $\iip{\cdot}{\cdot}_{\mH}=\iip{\cdot}{\mH\cdot}=\iip{\cdot}{\cdot}$ ( cf. Subsection \ref{sbs:QEMGHtf}). We still get from Corollary \eqref{cr:kerran} that:
\[
\cH^{BRST}_{phys}=\ker Q/ \cH_d\cong \cH_s=(\fF(\fH_t) \otimes \mathbb{C}\Omega_g).
\]
and from \eqref{th:kerdrand},
\[
\cP^{BRST}\cong P_s (\ker \drb \cap \cA) P_s= \cA_{ph} \otimes \one 
\]
the difference being that the ghost spaces, ghost algebras, and $Q$ are now different. That is using $\mH=\one$ gives equivalent physical results to the usual BRST-QEM with $\mH=M_z$ defined in Subsection \ref{sbs:QEMtf}.

However using $\mH=\one$ does \emph{not} give the usual smeared heuristic formulas for $\drb$. We can calculate this directly, but first we show what using $\mH=\one$ corresponds to heuristically. If we look at how $\Qs$ is constructed, we see that it is by summing smeared creators and annihilators over finite orthonormal basis $(f_j)_{j \in \Lambda}$ of $\fDL$. However these basis vectors correspond to
\[
(f_j(p))_{\mu}=\frac{-ip_{\mu}{h_j(p)}}{\sqrt{2}\norm{\mathbf{p}}},
\]
where $(h_j)_{j \in \Lambda}$ is an orthonormal basis of $\fH_0$ rather than the smeared constraints in $\fDL$, ie
\[
(f(p))_{\mu}=-ip_{\mu}{h(p)},
\]
where $h(p)\in \fH_0$ is has $\fH_0$ norm equal one. That is we have constructed a rigorous version of
\begin{align*}
 Q_h'&=\sqrt{2}\int d^3p (p_0)\left[\frac{b_{2}(\bp)^*}{\sqrt{2}p_0} \frac{c_2(\bp)}{\sqrt{2}p_0} + \frac{c_1(\bp)^*}{\sqrt{2}p_0} \frac{b_1(\bp)}{\sqrt{2}p_0}\right],\\
 &=(1/2)\int d^3p\; (1/p_0) [b_{2}(\bp)^* c_2(\bp) + c_1(\bp)^* b_1(\bp)]
\end{align*}
rather than \eqref{eq:Qpdef}
\begin{equation*}
Q_h= \sqrt{2} \int d^3p\; (p_0) [b_{2}(\bp)^* c_2(\bp) + c_1(\bp)^* b_1(\bp)].
\end{equation*}

To see the problem rigorously let $\mH=\one$. Then $\fL=\fH$ and $T:\fB \to \fDJ$ though still well defined (cf. lemma \eqref{lm:Tunit}), is no longer isometric and does not extend to a unitary from $\fH_0$ to $\fL$ (now $=\fH$). Hence we are no longer justified in using the identifications of $\tu(h)=C(Th)$ and $\tilde{\tu(h)}=-i\tu(ih)$ for $h\in \fH_0$ with the formal scalar ghosts as the anticommutation relations will not be preserved (cf. Proposition \eqref{pr:consvgh} (iii) requires that $T$ be isometric).
 
To get an $\fH$-isometric operator $T_2$ that goes from the scalar $\fD_0$ to $\fD$ we define,
\begin{definition}\label{df:T2}
Let $T_2:\fD_0 \rightarrow \fDJ$ be
\begin{align*}
T_2:& h(p) \to f(p)=\frac{ h(p)}{\sqrt{2}\np} (-p_0, p_1 , p_2 , p_3 ).
\end{align*}
\end{definition}
then we define $\tu_2(h):=\one \otimes \gh(T_2h)$. From  $(P_{\JJ}f)_{\mu}(p)=p^{\mu}(p_{\nu}f^{\nu}(p))/(2\nps)$, $f\in \fD$,
 we can calculate:
\begin{equation*}
(T_2^{-1}P_{\JJ}f)(p)= \frac{p_{\nu}}{\sqrt{2}\np}f^{\nu}(p)
\end{equation*}
Therefore Definition \eqref{df:brstsd}  and the above equation give:
\begin{align*}
\delta (A(f) \otimes \one )&= -i \one \otimes C(P_{\JJ}  {if})
=-i \one \otimes \tu_2( T_2^{-1}P_{\JJ}  {if})\notag \\
&=-i\one \otimes \tu_2 \left(-i\frac{p_{\mu}}{\sqrt{2}\np}{f^{\mu}(p)}\right)\\
& \neq-i\one \otimes \tu_2 \left(-i{p_{\mu}{f^{\mu}(p)}}\right)
\end{align*}
where the last line is what we would need for the previous formal correspondence. 
\end{rem}

\subsection{Covariance}\label{sbs:hcov}

Consider the Poincar{\'e} transformations generated by the transformations on $\fD$:
\[
(V_gf)(p):=e^{ipa}\Lambda f(\Lambda^{-1}p) \quad \forall f \in \cS(\mathbb{R}^4,\mathbb{C}^4), \quad g=(\Lambda,a) \in \cP^{\uparrow}_{+}.
\]
As $\Lambda$ is a constant matrix, $V_g$ preserves $\fD$. By the covariance of the formula for $\ip{\cdot}{\cdot}$ (cf. Subsection \ref{sbs:QEMtf})  $V_g$ is $\iip{\cdot}{\cdot}$-Krein-unitary on $\fD$. It is not however $\iip{\cdot}{\cdot}_{\mH}$-Krein unitary on $\fL$. As $\fL$ is the ghost test function space, and the ghosts are unphysical, we are free to construct any unitary representation of $\cP^{\uparrow}_{+}$ on $\fL$. We do this by using the Poincar{\'e} transformations on the scalar space $\fD_0$, \ie let
\[
(U_gf)(p):=:=e^{ipa} f(\Lambda^{-1}p) \quad \forall f \in \cS(\mathbb{R},\mathbb{C}), \quad g=(\Lambda,a) \in \cP^{\uparrow}_{+},
\]
then $g\to U_g$ defines a unitary representation $\cP^{\uparrow}_{+} \to B(\fH_0)$, (recall $\fH_0=\overline{\fD_0}=L^2(C_{+},\C, \lambda_0)$). Using this and the fact that $T\in B(\fH_0, \fLJ)$ is also unitary, we define $S_g\in B(\fLL\oplus \fLJ)$ by
\begin{equation}\label{eq:Sg}
S_g:= TU_gT^{-1}P^{\fL}_{\JJ} + JTU_gT^{-1}P^{\fL}_{\JJ}J,\qquad  \forall g=(\Lambda,a) \in \cP^{\uparrow}_{+}
\end{equation}
Now 
\begin{lemma}\label{lm:SgKHsa} We have:
\begin{itemize}
\item[(i)] $S_g$ is both $\iip{\cdot}{\cdot}_{\mH}$-unitary and 
$\ip{\cdot}{\cdot}_{\mH}$-unitary, for all $g=(\Lambda,a) \in \cP^{\uparrow}_{+}$.
\item[(ii)] $g \to S_g$ defines a representation $\cP^{\uparrow}_{+} \to B(\fLL\oplus \fLJ)$, which is both $\iip{\cdot}{\cdot}_{\mH}$-unitary and 
$\ip{\cdot}{\cdot}_{\mH}$-unitary.
\end{itemize}
\end{lemma}
\begin{proof}
(i): Let $f \in \fLL\oplus \fLJ$. By the definition of $T$ (Definition \eqref{df:T}) we have that $TU_gT^{-1}P^{\fL}_{\JJ}f \in \fLJ$ and $JTU_gT^{-1}P^{\fL}_{\JJ}Jf\in \fLL$. As $T$, $U_g, J$ are isometric and $\fLL \perp \fLJ$ we get that,
\begin{align*}
\norm{S_g f}_{\fL}^2= & \;\norm{TU_gT^{-1}P^{\fL}_{\JJ}f}_{\fL}^2+\norm{JTU_gT^{-1}P^{\fL}_{\JJ}Jf}_{\fL}^2,\\
= & \;\norm{P^{\fL}_{\JJ}f}_{\fL}^2+\norm{JP^{\fL}_{\JJ}Jf}_{\fL}^2,\\
= & \;\norm{P^{\fL}_{\JJ}f}_{\fL}^2+\norm{P^{\fL}_1f}_{\fL}^2,\\
= & \;\norm{ f}_{\fL}^2,
\end{align*}
and so $S_g$ is $\ip{\cdot}{\cdot}_{\mH}$-unitary.

From $J^2=\one$ and the definition of $S_g$ it follows that $JS_gJ=S_g$, hence $JS_g S_g^{\dag}=JS_g JS_g^{*}J=S_g S_g^{*}J=J$ therefore $S_g S_g^{\dag}=\one$ and similarly $S_g^{\dag} S_g=\one$, hence $S_g$ is  $\iip{\cdot}{\cdot}_{\mH}$-unitary .

\pfit (ii): Let $g_1,g_2 \in \cP^{\uparrow}_{+}$. As $g \to U_g$ is a representation and as $T:\fH_0 \to \fHJ$ implies $T^{-1}P^{\fL}_{\JJ}T=\one_{\fH_0}$, we get:
\begin{equation}\label{eq:sgrep1}
TU_{g_1g_2}gT^{-1}P^{\fL}_{\JJ}=TU_{g_1}U_{g_2}gT^{-1}P^{\fL}_{\JJ}=(TU_{g_1}T^{-1}P^{\fL}_{\JJ})(TU_{g_2}gT^{-1}P^{\fL}_{\JJ}).
\end{equation}
Similarly using $T^{-1}P^{\fL}_{\JJ} J^2 T=T^{-1}P^{\fL}_{\JJ}T=\one_{\fH_0}$ we calculate:
\begin{align}
JTU_{g_1g_2}T^{-1}P^{\fL}_{\JJ}J=JTU_{g_1}U_{g_2}T^{-1}P^{\fL}_{\JJ}J
=&\;JTU_{g_1}(T^{-1}P^{\fL}_{\JJ} J^2 T)U_{g_2}T^{-1}P^{\fL}_{\JJ}J, \notag \\
=&\;(JTU_{g_1}T^{-1}P^{\fL}_{\JJ} J)(J TU_{g_2}T^{-1}P^{\fL}_{\JJ} J) \label{eq:sgrep2}
\end{align}
Combining equation \eqref{eq:sgrep1} and equation \eqref{eq:sgrep2} with the definition of $S_g$ gives that
\[
S_{g_1g_2}=S_{g_1}S_{g_2}
\]
and hence $g \to S_g$ is a representation $\cP^{\uparrow}_{+} \to B(\fLL\oplus \fLJ)$, which is both $\iip{\cdot}{\cdot}_{\mH}$-unitary and 
$\ip{\cdot}{\cdot}_{\mH}$-unitary by (i).
\end{proof}
 
So given that $V_g$ and $S_g$ are Krein unitary in their respective IIP's, we see via \cite{Min1980} section 4 that $R_g=\Gamma_{+}(V_g)\otimes \Gamma_{-}(U_g)$ is well defined on $\fF_0^{B}$, preserves $\fF_0^{B}$ and is $\iip{\cdot}{\cdot}_T$-unitary (\ie Krein-unitary on $\fF_0^{B}$ cf. Section \ref{sec:brstext}). Moreover as $g\to V_g$ and $g \to U_g$ are representations of $\cP^{\uparrow}_{+}$ we get that and $g \to R_g$ is a $\iip{\cdot}{\cdot}_T$-unitary representation   $\cP^{\uparrow}_{+}\to \op(\fF_0^{B})$. We define:
\[
\alpha_g(A(f)\otimes C(f))\psi:=R_g^{\dag}(A(f)\otimes C(f)))R_g\psi=  (A(V_gf)\otimes C(S_gf))\psi, \qquad \psi \in \fF_0^{B}.
\]
Note that the for $g \in \cP^{\uparrow}_{+}$ above formula extends to a $\dag$-automorphism on all of $\cA=\cA_0\otimes \cA_g$ (cf. Section \ref{sec:brstext}), and that as $g \to R_g$ is a represention of $\cP^{\uparrow}_{+}$ we have
\[
g \to \alpha_g
\]
defines a representation of $\cP^{\uparrow}_{+} \to \aut(\cA)$.

We see that $\delta$ commutes with the relativistic transformations. 
\begin{lemma}\label{lm:comdrel} We have for all $g=(\Lambda,a) \in \cP^{\uparrow}_{+}$, and all $A\in \cA$.
\[
(\drb\circ \alpha_g)(A)\psi = (\alpha_g \circ \drb)(A) \psi, \qquad \psi \in \fF_0^{B}.  
\]
\end{lemma}
\begin{proof}
Clearly it suffices to check the statement on the generating elementary tensors in $\cA$. Suppose that $g=(\Lambda,a) \in \cP^{\uparrow}_{+}$, $f \in \fD$, and $\psi \in \fF_0^{B}$. Recall that $P^{\fL}_2f=P^{\fH}_2f=P_{\JJ}f$ for $f \in \fD$. Now from Definition \eqref{df:brstsd} we have that for $f \in \fD$,
\[
(\drb\circ \alpha_g)(A(f)\otimes \one)\psi=\drb(A(V_g f )\otimes \one)\psi=- i (\one \otimes C( P_{\JJ} V_g if))\psi, 
\]
and,
\[
(\alpha_g \circ \drb)(A(f)\otimes \one)\psi= -i (\one \otimes C( S_g P_{\JJ} if))\psi, 
\]
Also,
\[
(\drb\circ \alpha_g)(\one \otimes C(f))\psi=\drb(\one \otimes C(S_gf))\psi=  (A(\mH P_{\LL} S_g f)\otimes \one)\psi, \qquad f \in \fDL\oplus_{\mH} \fDJ
\] 
and,
\[
(\alpha_g \circ \drb)(\one \otimes C(f))\psi=  (A(V_g\mH P_{\LL} f)\otimes \one)\psi, 
\]
To complete the proof we must show that $P_{\JJ} V_g f=S_g P_{\JJ} f$ for $f \in \fD$ and $\mH P_{\LL} S_g f=V_g\mH P_{\LL} f$ for $f \in \fDL\oplus_{\mH} \fDJ$, which we do by direct calculation.

\medskip
\noindent{$\mathbf{P_{\JJ} V_g f=S_g P_{\JJ} f}$:}
\smallskip

Using Proposition \eqref{lm:stdec2} we calulate,
\begin{align*}
(P_{\JJ} V_gf(p))_{\mu}= e^{ipa}p^{\mu}(p_{\nu}(\Lambda f(\Lambda^{-1}p))^{\nu})/(2\norm{\mathbf{p}}^2).
\end{align*}
Also,
\begin{align*}
(S_gP_{\JJ}f(p))_{\mu}= & \;(TU_gT^{-1}P_{\JJ} f(p))_{\mu},\\
= & \;(TU_g) (p_{\nu}f^{\nu}(p))_{\mu}, &&\text{(by equation \eqref{eq:Tinv})}\\
= & \;e^{ipa}T((\Lambda^{-1}p)_{\nu}f^{\nu}(\Lambda^{-1}p))_{\mu},\\
= & \;e^{ipa} T(p_{\nu}(\Lambda f(\Lambda^{-1}p))^{\nu})_{\mu},&& \text{(by $\Lambda^{-1}=\Lambda^{\dag}$)}\\
=& e^{ipa}p^{\mu}(p_{\nu}(\Lambda f(\Lambda^{-1}p))^{\nu})/(2\norm{\mathbf{p}}^2), && \text{(by Defintion \eqref{df:T} (ii))}
\end{align*}
hence $P_{\JJ} V_g=S_gP_{\JJ}$.

\medskip
\noindent{$\mathbf{\mH P_{\LL} S_g f=V_g\mH P_{\LL} f}$ for $f \in \fDL\oplus \fDL$:}
\smallskip

Now $T:\fH_0 \to \fHJ$ implies $P^{\fL}_{\LL}Th=0$ for $h\in \fH_0$. Hence by the definition of $S_g$ we have $\mH P_{\LL} S_g f=\mH P_{\LL} JTU_gT^{-1}P_{\JJ}J f$, from which it follows similar to the above calculation that,
\[ 
(\mH P_{\LL} S_g f(p))_{\mu}=e^{ipa}p_{\mu}((\Lambda^{-1}p)_{\nu}( f(\Lambda^{-1}p))_{\nu})
\]
and also,
\begin{align*}
( V_g\mH P_{\LL}f(p))_{\mu}= & \;e^{ipa}\Lambda_{\mu \sigma}(\Lambda^{-1}p)_{\sigma}((\Lambda^{-1}p_{\nu})( f(\Lambda^{-1}p))_{\nu}),\\
= & \;e^{ipa}p_{\mu}((\Lambda^{-1}p)_{\nu}( f(\Lambda^{-1}p))_{\nu}),
\end{align*}
and so $\mH P_{\LL} S_g =V_g\mH P_{\LL}$.
\end{proof}
Using this lemma we show that $Q$ is relativistically invariant and hence that  BRST-QEM is `manifestly covariant'. 
\begin{proposition}\label{pr:invQ}
We have for all $g=(\Lambda,a) \in \cP^{\uparrow}_{+}$ :
\begin{itemize}
\item[(i)] For all $ \psi \in \fF_0^{B}$,
\begin{equation}\label{eq:defcoop}
\alpha_g(Q)\psi:=R_g^{\dag}QR_g\psi=Q\psi. 
\end{equation}
hence $R_g \ker Q = \ker Q$ and $R_g \ran Q= \ran Q$.
\item[(ii)] We have $R_g$ factors to a $*$-unitary operator $\hat{R_g}$ on the dense set $\cD^{BRST}_{phys}\subset \cH^{BRST}=\ker Q/ \overline{\ran Q}$, which extends to a $*$-unitary on $\cH^{BRST}_{phys}$. Hence $g \to \hat{R_g}$ is a Hilbert unitary representation of $\cP^{\uparrow}_{+}$.
\item[(iii)] We have:
\[
\alpha_g(\ker \drb \cap \cA) =(\ker \drb \cap \cA) , \qquad \alpha_g(\ran \drb \cap \cA)=(\ran \drb \cap \cA).
\]
Hence $\alpha_g$ factors to an automorphism, $\hat{\alpha_g}$ on $\cP^{BRST}:=(\ker \drb \cap \cA) /(\ran \drb \cap \cA)$, that is implemented by $\hat{R_g}$ when $\cP^{BRST}$ acts on $\cD^{BRST}_{phys}$ as described in Proposition \eqref{pr:BRSTopal} (i). Furthermore $g \to \hat{\alpha_g}$ is a representation $\cP^{BRST} \to \aut(\cP^{BRST})$.
\end{itemize}
\end{proposition}
\begin{proof}
(i): Let $\psi=P\Omega\in \fF_0^{B}$, where $T$ is a polynomial in $A(f)$'s and $C(g)$'s. Then using the fact that the vacuum is invariant, \ie $R_g \Omega=R_g^{\dag}\Omega=\Omega$ for all $g=(\Lambda,a) \in \cP^{\uparrow}_{+}$ and $Q \Omega=0$, we get that
\[
R_g^{\dag}QR_g\psi=R_g^{\dag}QR_gT\Omega=R_g^{\dag}\drb(R_gTR_g^{\dag})\Omega=\drb(T)R_g^{\dag}\Omega=QT\Omega=Q\psi
\]
where we use lemma \eqref{lm:comdrel} in the third last equality.

\pfit (ii): Part (i) implies that $R_g$ factors to a $\dag$-unitary operator $\hat{R_g}$ on the dense set $\cD^{BRST}_{phys}\subset \cH^{BRST}_{phys}=\ker Q/ \overline{\ran Q}$. Proposition \eqref{pr:QEMres} (iii) gives that $\iip{\cdot}{\cdot}_p=\ip{\cdot}{\cdot}_p$ on $\cH^{BRST}_{phys}$, hence the $\dag$-adjoint and $*$-are the same on $\cH^{BRST}_{phys}$ and so $\hat{R_g}$ is a $*$-unitary operator on $\cD^{BRST}_{phys}$ hence extends to a $*$-unitary operator on $\cH^{BRST}_{phys}$. Furthermore, $g \to  \hat{R_g}$ is a representation of $\cP^{\uparrow}_{+}$ as $g \to R_g$ is.

\pfit (iii): By lemma \eqref{lm:comdrel} we have for all $A\in \cA$.
\[
(\drb\circ \alpha_g)(A)\psi = (\alpha_g \circ \drb)(A) \psi, \qquad \psi \in \fF_0^{B}.  
\]
and so 
\[
\alpha_g(\ker \drb \cap \cA) =(\ker \drb \cap \cA) , \qquad \alpha_g(\ran \drb \cap \cA)=(\ran \drb \cap \cA).
\]
The rest of the statements are obvious.
\end{proof}
Part (iii) of the above Proposition says that the physical BRST representation is relativistically covariant.

Another useful result is,
\begin{lemma}\label{lm:covX1}
We have that $V_g\fXL \subset \fXL$, for all $g=(\Lambda,a) \in \cP^{\uparrow}_{+}$.
\end{lemma}
\begin{proof}
By Theorem \eqref{th:fdtstdec} and Proposition \eqref{lm:stdec2} we know that,
\[
\fXL= \fX\cap \fHL =\{ f \in \fX \,: \, f_{\mu}(p)=i(p_{\mu}/\np) h(p)\; \text{for}\; h\in \fH_0 \}
\]
and $V_g \fX \subset \fX$. Therefore for $f \in \fXL$
\[
(V_gf(p))_{\mu}=e^{ipa}\Lambda_{\mu \nu} (i(\Lambda^{-1}p)_{\nu}h(\Lambda^{-1}p)/\norm{\mathbf{\Lambda^{-1}p}})=e^{ipa} (ip_{\mu}/\norm{\mathbf{\Lambda^{-1}p}})h(\Lambda^{-1}p)
\]
Now using the formula for $P_{\LL}$ given in lemma \eqref{lm:stdec2},
\begin{align*}
(P_{\LL} V_gf(p))_{\mu}=&\; \left(\frac{p_{\mu}p_{\nu}}{2\norm{\mathbf{p}}^2}\right)(V_gf)_{\nu}(p), \\
=&\; e^{ipa}\left(\frac{p_{\mu}p_{\nu}}{2\norm{\mathbf{p}}^2}\right)(ip_{\nu}/\norm{\mathbf{\Lambda^{-1}p}})h(\Lambda^{-1}p),\\
=&\; e^{ipa}\left(\frac{p_{\mu}p_{\nu}}{2\norm{\mathbf{p}}^2}\right)(ip_{\nu}/\norm{\mathbf{\Lambda^{-1}p}})h(\Lambda^{-1}p),\\
=& (V_gf(p))_{\mu}
\end{align*}
where we used $p_{\nu}p_{\nu}=2\nps$ for $p \in \C_{+}$. As $V_g$ and $P_{\LL}$ preserve $\fX$ we get thatand so $V_gf\in (\fX \cap \ran P_{\LL})=\fXL$.
\end{proof}
The above result is unsurprising as $\fXL$ is the test function space corresponding to the Lorentz condition (subsection \eqref{sbs:lorentz}) and so it shows that the Lorentz condition is covariant in our setup.

\subsection{Example: Constraints for Massive Abelian Gauge Theory}\label{sbs:Mabga}
We show that massive abelian gauge theory has test function space with structure as above, and that applying the preceding \KOB gives the smeared version of the superderivation as in \cite{Schf2001} p28. We do not go through all the steps in detail as they are similar to those done for BRST-QEM.

Suppose again that $\widehat{\cS}_0$ and $\widehat{\cS}$ are as given in subsection \eqref{sbs:testfunc}, \ie 
\[
\widehat{\cS}= \{\hat{f} \,| \,f\in \cS(\mathbb{R}^4, \mathbb{R}^4) \}), \qquad \widehat{\cS}_0= \{\hat{f} \,| \,f\in \cS(\mathbb{R}^4, \mathbb{R}) \})
\]
Let the mass hyperboloid be for $m>0$:
\[
C_{+}^{m}:=\{ x \in \mathbb{R}^4\,|\, x_{\mu}x^{\mu}=-m^2, \,x_0>0\}.
\]
\begin{rem}
Note that the convention of $g_{\mu \nu}=diag(-1,1,1,1)$ gives the $-m^2$ term in the above definition of the mass hyperboloid. To define the mass hyperboloid using a positive $m^2$ term we can use the metric $-g_{\mu \nu}=diag(1,-1,-1,-1)$ as is commonly done, however this conflicts with the definition of $J$ (cf. equation \eqref{eq:Jdef}), hence our convention.
\end{rem}
Let $\lambda_m$ be the unique $\cL^{\uparrow}_{+}$ invariant measure on $C_{+}^{m}$ (see \cite{ReSi1975v2} p74). Define the inner product on $\widehat{\cS}_0$ by:
\begin{gather*}
 \ip{f}{h}^m:=  2\pi \int_{C_{+}^{m}} d\lambda_m \, \overline{f(p)}h(p), \qquad \forall f, h \in \widehat{\cS}_0 
 \end{gather*}
and define the symplectic form
\begin{align*}
 \smp{f}{g}{0}^m:=-\im ( \ip{f}{g}^m)= & \;i\pi \int_{C_{+}^{m}} d\lambda_m \left( \overline{f(p)}h(p)-f(p) \overline{h(p)} \right),
\end{align*}
Let $\fX_0^m:=\widehat{\cS}_0/\ker \sigma_0^m$. These are now the `Schwartz functions on $C_{+}^{m}$' and $\fX_0^m$ can be completed to $\fH_0^m:=L^2(C_{+}^{m},\C^4, \lambda_m)$. To define $\fD$ for the massive theory, we need to add in an extra field as will be seen by the following. Define the inner product on $\widehat{\cS}$ by:
\begin{gather*}
\ip{f}{h}^m:= 2\pi \int_{C^{+_m}} d\lambda_m \, \overline{f_{\mu}}(p)h_{\mu}(p), \qquad \forall f, h \in \widehat{\cS}
\end{gather*}
and define the symplectic form
\[
\smp{f}{h}{2}^m:= -\im(\ip{f}{h}_2)=i\pi\int_{C_{+}^{m}} d\lambda_m\left( \overline{f_{\mu}}(p)h_{\mu}(p)-{f_{\mu}}(p)\overline{h_{\mu}(p)}\right).
\]
We define ${\fX}_{m}:= \widehat{\cS}/\ker \sigma_2^m$. Let $Q_{\mu}(p)=p_{\mu}$ for $p \in C_{+}^m$, $\mu=0,1,2,3$, then $J$ and $M_{Q_\mu}$ are defined on $\widehat{\cS}$ as in Subsection \ref{sbs:testfunc} and factor to $\fX_m$. Using $\fX_m$ we get a massive gauge theory with 4-component vector field $A_{\mu}(f)$. We would like to impose the Lorentz condition $\partial_\mu  A^\mu(x)=0$ which corresponds to the subspace (cf. $\fXp_L$ in Subsection \ref{sbs:lorentz}) ,
\[
\fXL^m= \{ f \in \fX \,: \, f_{\mu}(p)=ip_{\mu}h(p)\; \text{for}\; h\in \fX_0^m \},
\]
However this presents a problem as for $p \in C_{+}^{m}$ we have $p_{\mu}p^{\mu}=-m^2 \neq 0$ and so $\fXL^m$ is not a $\sigma_1^m$-degenerate subspace. Hence there exist $f,g \in \fXL^m$ such that,
\[
[A(f),A(g)]\psi=c\psi, \qquad \psi \in \fF_0^{+}(\fX_m), \, 0\neq c \in \mathbb{C}
\]
if we are in the Fock-Krein representation in subsection \eqref{sec:FKCCR}. This tells us that the Lorentz condition is heuristically a \emph{second class} constraint set in Dirac sense (\cite{HenHur1988}) and so problematic. It is also problematic for rigorous theory as it means that we cannot apply \cite{Hendrik2000} Theorem 4.5 in the case of the Weyl Algebra or Proposition \eqref{pr:conRA} in the case of the Resolvent Algebra below.  

The way around this problem, both formally and rigorously, is to introduce a new scalar field with mass $m$, $B(x)$, called the Stueckelberg field \cite{RuRuiAl2004} p3272 and replace the Lorentz condition with $\partial_\mu  A^\mu(x)+mB(x)=0$. We now have that
\[
[\partial_\mu  A^\mu(x)+mB(x),\partial_\mu  A^\mu(y)+mB(y)]=0
\]
and so this is formally first class. Actually in the formal non-BRST version ( \cite{RuRuiAl2004} p3272) we use the constraints $(\partial_\mu  A^\mu(x)+mB(x))^{(-)}$ similar to the Gupta-Bleuler version of QEM. In BRST we impose the whole condition.

Rigorously, adding the new field involves adding another component to $\widehat{\cS}$, \ie we define:

\begin{gather*}
\widehat{\cSe}:= \{\hat{f} | f\in \cS(\mathbb{R}^4, \mathbb{R}^5) \}=\{ f\in \cS(\mathbb{R}^4, \mathbb{C}^5) | \overline{f(p)}=f(-p)\}\\
 \text{with IP:} \quad  \ip{f}{h}^{m}_{ext}:= 2 \pi \int_{C_{+}^{m}} \rd \lambda_m \, \big( \overline{f_{\mu}}(p)h_{\mu}(p)+\overline{f_4}(p)h_{4}(p) \big), \qquad \forall f, h \in \widehat{\cSe}
 \end{gather*}
and the symplectic form,
\[
\smp{f}{h}{2,ext}:= -\im \ip{f}{h}^{m}_{ext}
\]
and ${\fX}^m_{ext}:= \widehat{\cSe}/\ker\sigma_{2,ext}$. Now define $J_{ext}(f_0,f_1,f_2,f_3,f_4)=(-f_0,f_1,f_2,f_3,f_4)$ on $\cSe$. It follows that $J_{ext}$ factors to ${\fX}^m_{ext}$ and we denote it by the same symbol, and we define
\[
\smp{f}{g}{1,ext}:=\smp{f}{Jg}{2,ext}, \qquad f,g \in {\fX}^m_{ext}.
\]
Now the smeared version of $\partial_\mu  A^\mu(x)+mB(x)$ corresponds to,
\[
\fX_{ext}^{m,\LL}:= \{ f \in \fX \,: \, f_{\mu}(p)=ip_{\mu}h(p), f_{4}(p)=mh(p) \,|\, \mu=0,1,2,3 \quad h\in \fX_0^m \},
\]
It is easy to check using $p_{\mu}p^{\mu}=-m^2$ that $\fX^{m,\LL}_{ext}$  is a $\sigma_{1,ext}$-neutral subspace and so we can now use this in the $T$-procedure \cite{Hendrik2000} Theorem 4.5 or Resolvent Algebra Proposition \eqref{pr:conRA}, or in \KOB as above. For \KOB we define $\fD^m_{ext}:=\fX^m_{ext}+i\fX^m_{ext}$ and $\fH_{ext}^m:=\overline{\fD^m_{ext}}$. Let $s(p):=1/(2\norm{\mathbf{p}}^2+m^2)$ and define
\begin{align*}
(P_{\JJ}^m f)_{\mu}:&= M_sM_{P^{\mu}}(M_{P_{\nu}}f^{\nu}+m f^4), \qquad \mu=0,1,2,3 \quad f\in \fD^m_{ext}, \\
(P_{\JJ}^m f(p))_{4}:&= mM_s(M_{P_{\nu}}f^{\nu}+m f^4), \qquad \mu=0,1,2,3 \quad f\in \fD^m_{ext}, \\
\end{align*}
that is,
\begin{align*}
(P_{\JJ}^m f(p))_{\mu}:&= (p^{\mu}/(2\norm{\mathbf{p}}^2+m^2)(p_{\nu}f^{\nu}(p)+mf^4(p)), \qquad \mu=0,1,2,3 \quad f\in \fX^m_{ext}, \\
(P_{\JJ}^m f(p))_{4}:&= m/(2\norm{\mathbf{p}}^2+m^2)(p_{\nu}f^{\nu}(p)+mf^4(p)), \qquad \mu=0,1,2,3 \quad f\in \fX^m_{ext}, \\
\end{align*}
Note above that we are deviding by the factor $2\norm{\mathbf{p}}^2+m^2>0$ and so we do not have differentiability problems at $p=0$ in the massive case, hence we do not have to enlarge $\fX^m_{ext}$ as in Proposition \eqref{pr:extdec} to decompose the test function space analogously to Theorem \eqref{th:fdtstdec}. Using $p_{\mu}p_{\mu}=2\nps +m^2$ for $p \in C_{+}^m$, it is easy to see that $P_{\JJ}^m$ is a an algebraic projection on $\fX^m_{ext}$ and extends to a Hilbert space projection on $\cH$. Now defining $P_{\LL}^m:=J_{ext} P_{\JJ}^m J_{ext}$ and the projection $P_t^m:=\one-(P_{\LL}^m+P_{\JJ}^m) $ gives the following decomposition ,
\[
\fX^m_{ext}:=\fX^{m,t}_{ext} \oplus\fX^{m,\LL}_{ext} \oplus \fX^{m,\JJ}_{ext}
\]
where $\fX^{m,j}_{ext}:=P_j^m \fX^m_{ext}$ for $j=t,\LL,\JJ$. 
Similarly we have the decompositions:
\[
\fD^m_{ext}:=\fD^{m,t}_{ext} \oplus\fD^{m,\LL}_{ext} \oplus \fD^{m,\JJ}_{ext}, \qquad \fH^m_{ext}:=\fH^{m,t}_{ext} \oplus\fH^{m,\LL}_{ext} \oplus \fH^{m,\JJ}_{ext}
\]
where $\fD^{m,j}_{ext}:=P_j^m\fD^m_{ext}$, $\fH^{m,j}_{ext}:=P_j^m\fH^m_{ext}$ for $j=t,\LL,\JJ$. We use $\fD_{ext}^m$ as the test function space for the smeared massive gauge fields, and note that it has all the structures of the abstract test function space discussed in Subsection \ref{sbs:abstf}. For the ghosts we define 
\begin{gather*}
\mH:\cS(\R^4, \C^5) \to \cS(\R^4, \C^5)\\
(\mH f(p))_{\mu}:=2\norm{\bp}^2(f_{\mu}(p)), \qquad (\mH f(p))_4:=mf_4(p), \qquad \mu=0,1,2,3.
\end{gather*}
This factors to $\fD^m_{ext}$ which we also denote by $\mH$, and using this define and define $\fL{}_{ext}^m= \overline {\fD^{m,\LL}_{ext} \oplus \fD^{m,\JJ}_{ext}}$ in the $\ip{\cdot}{\cdot}_{\mH}$-topology as in Subsection \ref{sbs:abstf}.

Using these test function spaces we construct the BRST extension, superderviation, charge as in Section \ref{sec:FKBRST}  and calculate the BRST-physical subspace and BRST-physical algebra using Corollary \eqref{cr:kerran} (ii) and Theorem \eqref{th:kerdrand}:
\[
\cH^{BRST}_{phys}\cong \cH_s= (\fF^{+}(\fH_{ext}^{mt})\otimes \Omega_g), \qquad \cP^{BRST}\cong P_s (\ker \drb \cap \cA) P_s= (\cA(\fD^{mt}_{ext}) \otimes \one).
\]
To connect with the formal picture in \cite{Schf2001} p28 we want a scalar ghost (similar to QEM) so we define $T_m:\fX_0^m \rightarrow \fX_{ext}^{m2}$ by,
\begin{align*}
T_m: h(p) \to f(p)=\frac{h(p)}{({2}\nps+m^2)} (-p_0 , p_1 , p_2 , p_3 , m).
\end{align*}
Note that in the massive case we do not have to restrict $\fX_0^m$ to define $T_m$ as $({2}\nps+m^2)>0$ and we have:
\begin{itemize}
\item[(i)] $T_m^{-1}f(p)= (p_{\mu}f^{\mu} +mf^4(p))$ for $f \in \fX_{ext}^{m2}$, as can be verified by substitution.
\item[(ii)]  $T_m$ extends to $\fD^m_{ext}$ using the same defining formula and  $T_m$ can easily checked to be $\fL{}_{ext}^m$-isometric using $\ip{\cdot}{\cdot}_{\mH}$. As $T_m$ is also invertible we get that it extends to a unitary $T_m:\fH_0 \to \fL{}_{ext}^m$. 
\end{itemize}
Let
\[
\tu(h):=C(T_m h),\qquad  \tilde{\tu}(h):=i\tu(ih), \qquad h \in \fH_0,
\]
As $T_m$ is an isometry we have that using $\osalg{\tu(h),\, \tilde{\tu}(h)\,|\, h \in \fH_0}$ for the ghost algebra is equivalent to using $\cA_g(\fH^{m,\JJ}_{ext})$ as $T_m$ preserves the anticommutation relations (cf. Proposition \eqref{pr:consvgh}). Let $f_{\mu}\in{\fX}^m_{ext}$  be that vector with all $0$ entries except $f_{\mu} \in \fX_0^m$ in the $\mu$-th entry. We now find using the definition of $P_{\LL}^m,P_{\JJ}^m,T_m$ and Definition \eqref{df:brstsd}, the action of $\drb$ is:
\begin{align*}
\drb (A(\hat{f_{\mu}}) \otimes \one )&= -i \one \otimes \tu (\widehat{ \partial_{\mu}f^{\mu}}), \qquad \mu=0,1,2,3 \\
 \drb (A(\hat{f_{4}}) \otimes \one )&= -i \one \otimes \tu(m\widehat{f_4}) \\
 \drb (\one \otimes \tu(\hat{h})) &=0\\
 \drb (\one \otimes \tilde{\tu}(\hat{h})) &= -i A(\widehat{\partial_{0} h},\widehat{\partial_{1} h},\widehat{\partial_{2} h},\widehat{\partial_{3} h},\widehat{m h}) \otimes \one 
\end{align*}
If we associate $A(0,0,0,0,\hat{f_4})$ with the smeared Stueckelberg field $B(x)$ then these are the smeared version of the formal equations in \cite{Schf2001} p28, \ie 
\begin{align*} 
\drb(A^\nu(x))&=[Q\,,\, A^\nu(x)] = -i\partial^\nu \gh(x), \\
\drb(B(x))&=[Q\,,\, B(x)] = -im\,\gh(x), \\
\drb(\gh(x))&=\{Q\,,\, \gh(x)\} = 0, \\
\drb(\tilde{\gh}(x))&=\{Q\,,\, \tilde{\gh}(x)\} = -i(\partial_\mu  A^\mu(x)+mB(x)), 
\end{align*} 
Covariance can be treated similarly as in the QEM case, cf. Subsection \eqref{sbs:hcov}.
 
\section{Finite \KOB and Finite Abelian Hamiltonian BRST}\label{sec:Bsfinabhm}

We have seen in the previous sections that finite abelian Hamiltonian BRST (Subsection \ref{sbs:hhamexphsp}) has attractive features as it assumes little structure, only  beginning with a finite set of commuting self adjoint constraints. However the BRST-physical subspace it selects is larger than that selected by the original constraints and has indefinite inner product. Restricting to ghost number zero states does not fix the problem. On the other hand \KOB selects the correct physical state space, however much more structure is assumed and there are difficulties associated with unboundedness of the operators involved. The aim of this section is to combine finite abelian BRST with finite \KOB (Example \eqref{sbs:findimbos} ) to get a general BRST constraint algorithm that selects the correct physical state space. This is done by enlarging the orignal constraint system by tensoring on a ghost algebra and a bosonic field algebra. $Q$ is then constructed in a natural way with all the usual properties and we see that $\cH_s$ from the \emph{dsp}-decomposition is naturally isomorphic with the original Dirac physical state space.

As before begin with a unital $C^*$-algebra $\cA^{H}_0$  acting on the Hilbert space $\cH_0$ with inner product $\ip{\cdot}{\cdot}_{0}$, and let the constraint set $\cC=\{G_i\,|\, i=1,\dots,m\}\subset \cA^{H}_0$ be a finite linearly independent set of self-adjoint commuting operators. The physical state space is $\cH^0_p:=\cap_{i=1}^{m}\ker G_i$.  We extend by tensoring on a ghost algebra and a CCR field algebra. We assume the structures of finite \KOB as in Subsection \eqref{sbs:findimbos}.

Suppose that $\MM< \infty$ is a positive integer, and let $\fD=\fH$ be a Hilbert space with $\dim(\fD)=2\MM$ and Hilbert inner product $\ip{\cdot}{\cdot}$. If we recall the construction in lemma \eqref{lm:ghtsfext} then:
\begin{itemize}
\item[(i)] There exists a decomposition $\fD=\fDL \oplus \fDJ$
\item[(ii)] There exists a unitary $J\in B(\fD)$ with $J^2=\one$ such that $J\fDL=\fDJ$ and $\fD$ has Krein inner product $\iip{\cdot}{\cdot}:=\ip{\cdot}{J\cdot}$, hence $\fDL [\oplus] \fDJ$. 
\item[(iii)] Let $\mH=\one$ so $\fL=\fD$ where $\fL$ is defined as in Subsection \eqref{sbs:abstf}.
\end{itemize}
It follows now that $\fD$ has all the structure of $\fD$ in Subsection \eqref{sbs:abstf} with $\fD_t=\{0\}$ and with ghost test function space $\fL=\fD$.
Also let $(f_j)_{j=1}^{\MM}$ be an orthonormal basis of $\fDL$.

Now let,
\[
D(Q)=\cH_0\otimes \fF^{+}(\fH)_0\otimes \fF^{-}(\fH).
\]
Note $\fF^{-}(\fH)=\fF^{-}(\fH)_0$ as $\fH$ is finite dimensional, let $\cH:=\overline{D(Q)}=\cH_0\otimes\fF^{+}(\fH)\otimes \fF^{-}(\fH) $, and let $\ip{\cdot}{\cdot}$ be the its usual Hilbert space inner product. Let $J_g=\one \otimes \Gamma_{+}(J)\otimes \Gamma_{-}(J)$. Then $J_g$ is a fundamental symmetry that makes $\cH$ a Krein space. Let the indefinite inner product be $\iip{\cdot}{\cdot}=\ip{\cdot}{J_g\cdot}$ which has associated involution $\dag$. Let
\[
\cA=\cA^{H}_0 \otimes \cA^{B}_0\otimes \cA_{g}
\]
where $\cA^{B}_0:= \alg{A(f)\otimes \one\, |\, f \in \fH}$ and $\cA_g:=\alg{\one \otimes C(g)\,|\,g \in \fH}$ and we are assuming the algebraic tensor product. Note that $\cA_g=\cA_{g0}$ as $\fH$ finite dimensional (\cite{BraRob21981} Theorem 5.2.5). To connect to the terminology in Chapter \ref{ch:GenStruct}  we define $\cA_0:=\cA_0^H\otimes \cA_0^B$ and grade $\cA$ as in Section \ref{sbs:AbcuQ}. Define,  
\begin{align*}
Q^B:=\sqrt{2} \sum_{j=1}^{\MM}\left[ \one \otimes A(f_j)\otimes C(Jf_j) + \one \otimes  A(if_j) \otimes C(iJf_j) \right], \qquad D(Q^B)=D(Q)
\end{align*}
where $\Lambda:=(f_j)_{j=1}^{\MM}$ is an orthonormal basis of $\fDL$ and,
\begin{align*}
Q^H:= & \;\sum_{j=1}^{\MM}\left[G_j\otimes \one \otimes C(Jf_j) +G_j\otimes \one \otimes C(iJf_j) \right], \qquad D(Q^H)=D(Q)\\
= & \;\sum_{j=1}^{\MM}\left[G_j\otimes \one \otimes C((1+i)Jf_j)\right], \qquad D(Q^H)=D(Q)
\end{align*}
We have that $Q^B$ and $Q^H$ are the the \KOB and Hamiltonian BRST charges respectively. Following the proofs as in lemma \eqref{lm:2nQs} and Section \ref{sec:hamBRST} we see that these are both Krein symmetric and 2-nilpotent with dense domains. Now we define the BRST charge as the sum of these two charges, ie
\begin{align}
Q:=\sum_{j=1}^{\MM}\left[(G_j\otimes \one + \one \otimes \sqrt{2} A(f_j))\otimes C(Jf_j) +(G_j\otimes \one + \one \otimes \sqrt{2} A(if_j))\otimes C(iJf_j) \right)\label{eq:QHamBos}
\end{align}
with domain $D(Q)$.

For the remainder of this section we drop the tensor product $\otimes$ for ease of notation.
\begin{lemma}
$Q$ is an odd 2-nilpotent and Krein symmetric.
\end{lemma}
\begin{proof}
That $Q$ is odd and Krein symmetric is obvious. We have that for all $f_j,f_k \in \Lambda \subset \fDL$ 
\[
0=[G_j+A(f_j),G_k+A(f_k)]\psi=[G_j+A(f_j),G_k+A(if_k)]\psi=[G_j+A(f_j),G_k+A(if_k)]\psi
\]
for all $\psi \in D(Q)$. It follows that $Q^2\psi=0$ for all $\psi \in D(Q)$ similarly as in lemma \eqref{lm:2nQs}.
\end{proof}
Therefore  we have that $Q$ is Krein symmetric, hence closable by Proposition \eqref{pr:Krclos}, and that $\cH$ has an \emph{dsp}-decomposition $\cH=\cH_d \oplus \cH_s \oplus \cH_p$ with respect to $\cQ$ (cf. lemma \eqref{lm:QessaCsa} and Theorem \eqref{th:Hdsp1}). We want to calculate $\cH_s$. As $Q$ is densely defined, $Q^*$ exists. We can easily check that
\[
Q^*\psi=\sum_{j=1}^{\MM}\left[(G_j + \sqrt{2} A(Jf_j)) C(f_j) +(G_j +  \sqrt{2} A(iJf_j)) C(if_j)\right]\psi
\]
for $\psi \in D(Q)$ and hence $Q^*$ preserves $D(Q)$. To calculate $\cH_s$ we will calculate $\ker \Delta$ on $D(Q)$ as in lemma \eqref{lm:kDelHs}.  To do this we first need:
\begin{lemma}\label{lm:tired1}
Let $1\leq k \leq \MM$ and $g \in \fH$. Then
\[
\ker (G_k+a(g))\cap D(Q)= (\ker G_k\cap D(Q)) \cap (\ker a(g) \cap D(Q))
\]
where $a(g)$ is the annihilator in $\cA_0^B$.
\end{lemma}
\begin{proof}
Let $N_b$ be the number operator on $\fF^{+}(\fH)_0$ (\cite{BraRob21981} Chapter 5.2.1). Let $N=\one \otimes N_b \otimes \one$ be the extension to $D(Q)$. Then $D(Q)= \sn{ \bigcup_{j} \cH_{j}\,|\, j \in \Z^{+}}$ where $\cH_{j}$ are the eigenspaces of $N$ with integer eigenvalue $j$. It clear that $G_k \cH_j \subset \cH_j$ and $a(g)\cH_{j} \subset \cH_{j-1}$ for $j>0$, $1\leq k \leq \MM$ and $a(g) \cH_{0}=0$. 

Let $\psi \in \ker (G_k+a(g))\cap D(Q)$. As $\psi \in D(Q)$ we have for some $n\geq 0$, that $\psi=\sum_{j=0}^{n} \psi_{j}$ where $\psi_{j} \in \cH_{j}$. Therefore,
\begin{align*}
0=(G_k+a(g))\psi= & \;G_k\psi_n + \sum_{j=1}^{n}(G_k\psi_{j-1} + a(g)\psi_{j}) +a(g) \psi_0,\\
&=G_k\psi_n + \sum_{j=1}^{n}(G_k\psi_{j-1} + a(g)\psi_{j}).
\end{align*}
As $\cH_{j}\perp \cH_{i}$ for $i \neq j$ we get
\begin{align}
G_k\psi_n= & \;0, \label{eq:han1}\\
G_k\psi_{j-1}+a(g)\psi_{j}= & \;0. \label{eq:han2}
\end{align}
Combining $[G_k,a(g)]\psi=0$ for $\psi \in D(Q)$ and  equation \eqref{eq:han2} gives 
\[
G_k(G_k\psi_{j-1})=-a(g)G_k\psi_{j}.
\]
Substituting $j=n$ gives $G_k^2 \psi_{n-1}=-a(g)G_k\psi_{n}=0$ by \eqref{eq:han1}. Now $G_k^*=G_k$ hence $\ker G_k \perp \ran G_k$, and so $G_k G_k \psi_{n-1}=0$ implies $G_k \psi_{n-1}=0$. Using equation \eqref{eq:han2} again gives $a(g)\psi_{n}=0$. 

Summarising we have $\psi_n \in (\ker a(g) \cap \ker G_k)$ and $\psi_{n-1} \in \ker G_k$. We can iterate the above arguments using equations \eqref{eq:han1} and \eqref{eq:han2} to get that $\psi_j \in (\ker a(g) \cap \ker G_k)$ for $0\leq j \leq n$ and hence $\psi \in (\ker a(g) \cap \ker G_k)$.

This shows that $\ker(G_k +a(g))\cap D(Q) \subset ( \ker G_k \cap D(Q)) \cap ( \ker a(g) \cap D(Q))$. The reverse inclusion is obvious. 

\end{proof}
Using this lemma we can prove that
\begin{proposition}\label{pr:HamBosstsp}
Let $\cH_s$ be defined as in the \emph{dsp}-decomposition. Then 
\[
\cH_s= \cH^0_p\otimes \mathbb{C}\Omega_b \otimes \mathbb{C} \Omega_g,
\] 
where $\Omega_b$ is the vacuum vector in $\fF^{+}(\fH)_0$ and $\Omega_g$ is the vacuum vector in $\fF^{-}(\fH)_0$, and $\cH^0_p= \cap_{k=1}^{\MM}\ker G_k$.
\end{proposition}
\begin{proof}
Let $\psi \in D(Q)$. Then
\begin{align*}
\{Q,Q^*\}\psi= & \; \{Q^{B*}, Q^B\}\psi +\sum_{jk} G_jG_k\{C(Jf_j)+C(iJf_j),C(f_k)+C(if_k)\} \\
& + \sum_j \{ Q^B, G_j(C(f_j)+C(if_j)) \} \\
&+\sum_j \{ Q^{B*}, G_j(C(Jf_j)+C(iJf_j)) \},\\
= & \; A +B +C +D
\end{align*}
where the obvious correspondence between the terms in the two RHS identities are made. Since $Q^B$ has the same structure as $\Qs$ in lemma \eqref{lm:2nQs}, the proof of equation \eqref{eq:kerD} adapts immediately to give:
\[
A\psi=2\sum_{j} (a^*(Jf_j)a(Jf_j)+ a^*(f_j)a(f_j)+ c^*(Jf_j)c(Jf_j)+ c^*(f_j)c(f_j) )\psi
\]
Using the CAR relations $\{C(f),C(g)\}=\re \iip{f}{g}\one$ for all $f,g \in \fD$  (cf. equation \eqref{eq:Cliffcom}) we get that
\[
B\psi=(\sum_{jk=1}^{\MM} G_jG_k\{C(Jf_j)+C(iJf_j),C(f_k)+C(if_k)\})=2(\sum_{j=1}^{\MM} (G_j)^2)\psi.
\]
Now from the action of the BRST superderivation for Bose-Fock theories (cf. Definition \eqref{df:brstsd}) we have that $\{ Q^B, G_j(C(f_j)+C(if_j)) \}\psi=G_j\{ Q^B, (C(f_j)+C(if_j)) \}\psi=G_j(A(f_j)+A(if_j))\psi$ and by taking $*$-adjoints $\{ (Q^{B})^*, G_j(C(Jf_j)+C(iJf_j)) \}\psi=G_j(A(Jf_j)+A(iJf_j))\psi$. Therefore using $A(g)=\fst(a(g)+a^*(Jg))$ for $g\in \fD$ gives:
\[
(C+D)\psi=\sum_{j}(G_j(a^*(f_j+if_j)+a^*(Jf_j+iJf_j)+a(f_j+if_j)+a(Jf_j+iJf_j)))\psi
\]
Also, using the linearity of $f\to a^*(f)$ and antilinearity $f \to a(f)$ we get that
\begin{equation}\label{eq:tired1}
a^*(f_j+if_j)a(f_j+if_j)\psi=2a^*(f_j)a(f_j)\psi.
\end{equation}
Therefore we get that
\begin{align*}
\{Q^*,Q\}\psi= & \;\sum_{j} \Big(2[c^*(Jf_j)c(Jf_j)+ c^*(f_j)c(f_j) ] +[(G_j+a^*(f_j+if_j))(G_j+a(f_j+if_j))] \\
&+[(G_j+a^*(J(f_j+if_j))(G_j+a(J(f_j+if_j)))]\Big)\psi
\end{align*}
which can be checked by expanding, using equation \eqref{eq:tired1} and the expressions for $A,B,C,D$.

Now this is a sum of positive operators and so:
\begin{align*}
\ker \{Q^*,Q\} = \bigcap_{j} &( \ker [c^*(Jf_j)c(Jf_j)+ c^*(f_j)c(f_j) ]\\
& \cap \ker [(G_j+a^*(f_j+if_j))(G_j+a(f_j+if_j))] \\
&\cap \ker[(G_j+a^*(J(f_j+if_j))(G_j+a(J(f_j+if_j)))] ), \\
=\bigcap_{j} &\Big[ \ker c(Jf_j) \cap \ker c(f_j)  \cap \ker (G_j+a(f_j+if_j))\cap \ker(G_j+a(J(f_j+if_j))) \Big] 
\end{align*}
where the last equality follows from $\psi \in \ker T^*T \Leftrightarrow \psi \in \ker T$. By lemma \eqref{lm:tired1} to get that $\ker (G_j+a(f_j+if_j))= \ker G_j \cap \ker a(f_j+if_j)$. But  $\ker a(f_j+if_j)= \ker (1-i)a(f_j)=\ker a(f_j)$. Therefore we have that, 
\begin{align*}
\ker \{Q^*,Q\} = & \; \bigcap_{j} (\ker G_j \cap \ker c(Jf_j) \cap \ker c(f_j) \cap \ker a(Jf_j) \cap \ker a(f_j) ),\\
= & \; \cH^0_p\otimes \mathbb{C}\Omega_b \otimes \mathbb{C} \Omega_g,
\end{align*}
and so by lemma \eqref{lm:kDelHs}  we have that $\cH_s=\cH^0_p\otimes \mathbb{C}\Omega_b \otimes \mathbb{C} \Omega_g$.
\end{proof}
The above Proposition tells us that the BRST charge constructed above will select the correct physical space for an arbitrary finite commuting set of constraints $ \cC=\{G_i\,|\, i=1,\dots,m\}\subset \cA^{H}_0$ without adding ghost number restrictions and does not suffer the MCPS problem and neutrality problems of usual Hamiltonian BRST  (cf. Subsection \ref{sbs:hhamexphsp}, Remark \eqref{rm:spfstuff} (ii) ).

\begin{eje}
We revisit the example in subsection \eqref{sbs:egAHC}, and we want to show that adding this extended BRST now selects the correct physical algebra. Let $\cH_0$ be a Hilbert space and that $\cA^H_0=B(\cH_0)$. Let  $\cC=\{C_j\,|\, j=1,\ldots,\MM\}$ be a finite self adjoint set of commuting constraints and so $\cH^0_p=\cap_{j=1}^{\MM}\ker C_j $. Let $(\one-P)$ be the projection on $\cH^0_p$, and apply the $T$-procedure (section \eqref{app:Tp}) to $(\cA_0,\cC)$. As $(\one - P) \in \cA^H_0$ it is straightforward to see that for the Dirac physical algebra we have the following isomorphism $\cP_{phys}\cong (1-P)\cA^H_0 (1-P)$ (see example \eqref{ex:DC} (i)). 

We have $\MM$ constraints and so $\dim(\fD)=2\MM$ and we construct the BRST model above.
\begin{proposition}\label{pr:BoHaalex}
Let $\cA=\cA^H_0\otimes \cA^B_0\otimes \cA_g$ be as above, and let $Q$ have domain $D(Q)$ and be defined by equation \eqref{eq:QHamBos}. Then:
\begin{itemize}
\item[(i)] We have $\cH_s=\cH^0_p\otimes \mathbb{C}\Omega_b \otimes \mathbb{C} \Omega_g=\ker P\otimes \mathbb{C}\Omega_b \otimes \mathbb{C} \Omega_g$. Therefore 
\[
P_s=(\one -P)\otimes P_{\Omega_b} \otimes P_{\Omega_g}\in \cA,
\]
where $P_{\Omega_b},P_{\Omega_g}$ are the projections onto $\mathbb{C}\Omega_b$, $\mathbb{C} \Omega_g$.
\item[(ii)] Define $\drb(\cdot):=\sbr{Q}{\cdot}$ for $A \in \cA$ as in Subsection \eqref{sbs:AbcuQ}, then 
\[
\tilde{\cP}^{BRST}\cong \Phi_s(\ker \drb) = (\one - P) \cA_0^H (\one -P)\otimes \one \otimes \one \cong \cP_{phys},
\]
where we are using the alternative definition of the BRST observables in Subsection \ref{sbs:alalg}, and $\Phi_s(\ker \drb)$ for all $A \in \cA$.
\end{itemize}
\end{proposition}
\begin{proof}
(i): That $\cH_s=\cH^0_p\otimes \mathbb{C}\Omega_b \otimes \mathbb{C} \Omega_g$ follow by Proposition \eqref{pr:HamBosstsp}. And as $(\one-P)$ is the projection   on $\cH^0_p$ this implies $\cH_s=\ker P\otimes \mathbb{C}\Omega_b \otimes \mathbb{C} \Omega_g$ hence $P_s=(\one -P)\otimes P_{\Omega_b} \otimes P_{\Omega_g}\in \cA$ where $P_{\Omega_b},P_{\Omega_g}$ are the projections onto $\mathbb{C}\Omega_b$, $\mathbb{C} \Omega_g$.

\pfit (ii): Let $T \in \ker \drb$. Using the commutation relations for $a(f), c(f)$ we can write,
\[
T=A_0\otimes \one \otimes \one + \sum_{j=1}^{M} A_j\otimes M_j
\]
where $A_j \in \cA^H_0$, $M_j$ is normally ordered monomial of creators and annihilators corresponding the ghosts and bosonic field. As the $M_j$ are normally ordered we have that $M_j(P_{\Omega_b} \otimes P_{\Omega_g})\psi=0$ for $\psi \in D(Q)$ and so,
\[
P_sTP_s\psi=P_s(A_0\otimes \one \otimes \one)P_s=(\one -P)A_0(\one-P)\otimes \one \otimes \one.
\]
Also, $QP_s\psi=P_sQ\psi=0$ for $\psi \in \ker Q$ and so as $P_s \in \cA$, we have that $P_s\cA P_s \subset \ker \drb$. Combining this with the above calculation gives that,
\[
P_s\cA P_s= \Phi_s(P_s\cA P_s) \subset \Phi_s(\ker \drb) = (\one - P) \cA_0^H (\one -P)\otimes \one \otimes \one =P_s \cA P_s,
\]
\end{proof}
Proposition \eqref{pr:BoHaalex}(ii) shows that the combined \KOB and Hamiltonian BRST charge $Q$ selects the Dirac physical space and algebra in this example. It is important to note that we chose our original algebra to be $\cA^H_0=B(\cH_0)$ so that $(\one-P) \in \cA^H_0$ which in turn gave $P_s \in \cA$. This was crucial for the above isomorphism.
\end{eje}

\begin{rem}
\begin{itemize}
\item[(i)] It would be nice to extend the above procedure in two directions: The Hamiltonian case where the constraints do not commute, and the case of an infinite set of constraints $\cC$. This is work still in progress.
\end{itemize}
\end{rem}

\chapter{$C^*$-BRST}\label{ch:CsBRST}

In this chapter we will cast the above structures into a $C^*$-algebraic context.

So far we have given an account of the general structures of quantum BRST in a given representation in Chapter \ref{ch:GenStruct} and constructed and examined explicitly rigorous examples for basic abelian Hamiltionian BRST (cf. Subsection \ref{sbs:egAHC}), and abelian Bose-Fock theories, such as QEM, in Chapter \ref{ch:BRSTQEM}. In the algebraic approach we would like to be able to move beyond a given representation and construct BRST structures at the $C^*$-algebraic level. For the case of Hamiltionian BRST given in Section \ref{sec:hamBRST}, this is straightforward as all have already assumed that our operators, etc are bounded. However for BRST-QEM we have that the basic objects, such as the fields, the BRST charge $Q$, etc. are unbounded and so we need more elaborate constructions. 

We first have to interpret the structures of Chapter \ref{ch:GenStruct} in an abstract setting. As already discussed in the introduction to Chapter \ref{ch:BRSTQEM}, we take that the BRST superderivation $\drb$ to be of primary interest and so aim to make a formulation of BRST as a superderivation acting on a $C^*$-algebra $\cA$. Once this is done we need to identify the correct states on $\cA$. To motivate our definitions we investigate what a BRST theory with bounded $Q$ looks like and what the physical states on this theory will be.

With bounded BRST as a guide we aim to construct BRST for QEM. For this, we need to find a $C^*$-form for the superderivation in definition \eqref{df:brstsd}:
\begin{align*}
\delta(A(g)\otimes \one)&= i \one \otimes C( P_{\JJ}  ig), \qquad g \in \fD, \\
\delta(\one \otimes C(g))&= A(\mH P_{\LL} g)\otimes \one \qquad g \in \fD.
\end{align*}
Due to unboundedness, there are several technical hurdles to overcome, such as domain issues related to the unboundedness of the fields $A(g)$. The main tool used for dealing with these issues is the Resolvent Algebra \cite{HendrikBuch2007,HendrikBuch2006}, which encodes the CCR relations in bounded form similar to the Weyl algebra. Using the QEM test function space we perform the $T$-procedure (cf. Appendix \ref{app:Tp}) on the Resolvent Algebra, calculate the physical algebra and see that this gives the results we expect from other examples \cite{Hendrik2000}.

We then encode the \KOB structures in bounded form, following an approach similar to \cite{HendrikBuch2006}. A surprise is that \KOB using the symplectic test function space with covariant symplectic form $(\fX,\sigma_1)$ (cf. Subsection \eqref{sbs:abstf}), produces more BRST-observables than what we get from the $T$-procedure. This is because the BRST-procedure does not remove the ghosts.

To resolve this issue we study a second the Resolvent Algebra using the auxiliary symplectic form, \ie $(\fD,\sigma_2)$ as in Subsection \ref{sbs:abstf}. In this case we get that BRST-procedure and the Dirac procedure produce the same physical algebra, but we pay the price that we have to work harder to encode the Poincar{\'e} transformations.

Finally, we give a general formulation of a BRST-theory motivated by these results, and show how the examples we have seen so far fit into this framework.

We use the following notation with respect to $C^*$-algebras:
\begin{definition}\label{df:csksrepterm}
Let $\cA$ be a $C^*$-algebra: 
\begin{itemize}
\item $\fS(A)$ is the set of states of $\cA$.
\item Let $\cB$ be a $C^*$-subalgebra of $\cA$ and $\omega \in \fS(\cA)$. Then $\omega_{\cB}:=\omega|_{\cB} \in \fS(\cB)$.
\item $\pi_{\omega}:\cA\to \cH_{\omega}$ is the GNS-representation associated $\omega \in \fS(\cA)$, $\Omega_{\omega}\in \cH_{\omega}$ is its cyclic generating vector, and we denote its inner product by $\ip{\psi}{\xi}_{\omega}$ for all $\psi, \xi \in \cH_{\omega}$.
\item Let $\alpha \in \aut(\cA)$ be $*$-automorphism and $\omega \circ \alpha=\omega$. Then there exists a unitary $U_{\omega} \in \cH_{\omega}$ such that $\pi_{\omega}(\alpha(A))=U_{\omega}^*\pi_{\omega}(A)U_{\omega}$ for all $A \in \cA$ and $U_{\omega}\Omega_{\omega}=\Omega_{\omega}$. We call $U_{\omega}$ the \emph{implementer} of $\alpha$ in $\pi_{\omega}:\cA\to B(\cH_{\omega})$.
\item Let $\alpha \in \aut(\cA)$ be $*$-automorphism, $\omega \circ \alpha=\omega$, and $\alpha^2=\iota$. Let $J_{\omega}$ be the implementer of $\alpha$. It follows that $J_{\omega}^2=\one$ and $J_{\omega}^*=J_{\omega}$. We define the indefinite inner product $\iip{\cdot}{\cdot}_{\omega}$,
\[ 
\iip{\psi}{\xi}_{\omega}:=\ip{\psi}{J_{\omega}\xi}_{\omega}, \qquad \psi, \xi \in \cH_{\omega}. 
\]
By lemma \eqref{lm:JKsp} we have that $\cH_{\omega}$ with $\iip{\cdot}{\cdot}_{\omega}$ is a Krein space with fundamental symmetry $J_{\omega}$. This will always be the Krein structure we associate to a $*$-automorphism such as $\alpha$. We define an involution on $\cA$ by 
\begin{equation}\label{eq:daginvalg}
A^{\dag}:=\alpha(A^*), 
\end{equation}
and note that by lemma \eqref{lm:khad}, 
\[
\pi_{\omega}(A)^{\dag}=J_{\omega}\pi_{\omega}(A)^*J_{\omega}=J_{\omega}\pi_{\omega}(A^*)J_{\omega}=\pi_{\omega}(\alpha(A^*))=\pi_{\omega}(A^{\dag}),
\]
where $\pi_{\omega}(A)^{\dag}$ is the Krein adjoint of $\pi_{\omega}(A)$ in $(\pi_{\omega},\cH_{\omega})$. Therefore equation \eqref{eq:daginvalg} is the natural way to define the Krein involution with respect to $\alpha$ on $\cA$ in a representation independent way.
\end{itemize}
\end{definition}

\section{Bounded $\drb$}\label{sbs:bddcbrst}
We first cast the BRST structures of Chapter \ref{ch:GenStruct} in a $C^*$-algebra setting, for the case where $\cA_0$ is a $C^*$-algebra. We assume that $\cA_0$ is a unital $C^*$-algebra and that $\beta \in Aut(\cA_0)$ is a $*$-automorphism such that $\beta^2=\iota$. $\beta$ encodes any Krein structure present in $\cA_0$, and we define an involution on $\cA_0$ by $A^{\dag}:=\beta(A^*)$ for all $A\in \cA_0$. When we have no relevant Krein structures we will set$\beta=\iota$.

The \emph{Unextended Field Algebra} is $\cA_0$ and we assume that it has some kind of degeneracy, such as constraints. Depending on the example we are modelling we will take the $C^*$-tensor product of $\cA_0$ with either the ghost algebra $\cA_g(\cH_2)$ or the restricted ghost algebra $\rga(\cH_2)$  (cf. Subsection \ref{sbs:fndmgh}), where $\cH_2$ corresponds to the degrees of degeneracy, \ie $\dim(\cH_2)=$ number of linear independent constraints in the Hamiltonian case (cf. Definition \ref{df:GA}). As $\cA_g(\cH_2)$ is a CAR algebra and $\rga(\cH_2)$ is finite dimensional we have that these are both nuclear $C^*$-algebras, the norm on the tensor product is unique.  Let $\alpha'\in \Aut(\cA_g)$ be the automorphism that corresponds to the Krein-ghost stucture (cf. equation \eqref{eq:autghalg}), hence $(\alpha')^2=\iota$. 
\begin{definition}\label{df:fagbbrst}
The \emph{BRST-Field Algebra} is either $\cA:=\cA_0\otimes \cA_g$ or $\cA:=\cA_0\otimes \rga$ and we let $\cA$ have the $\mathbb{Z}_2$-grading with grading automorphism $\iota\otimes\gamma$ (cf. Definition \eqref{df:Z2grad}). Let $\alpha:=\beta\otimes \alpha' \in \Aut(\cA)$ and note that $\alpha^2=\iota$. Define the involution on $\cA$ by:
\[
A^{\dag}:=\alpha(A^*), \qquad \forall A\in \cA.
\]
We assume that there exists a BRST-charge $Q \in \cA$ such that $Q$ is 2-nilpotent, $Q^{\dag}=\alpha(Q^*)=Q$ and $Q \in \cA^{-}$, \ie $\gamma(Q)=-Q$. 
\end{definition}
\begin{rem}
The tensor $C^*$-norm on $\cA$ is unique as $\cA_g$ is a CAR algebra hence nuclear (cf. \cite{Bla2006} Example II.8.2.2 (iii), Example II.9.4.2 and II.9.4.5). Also $\rga$ is only used in the case of a finite number of ghosts, hence $\rga$ is finite dimensional and hence nuclear.
\end{rem}

This is the situation as for Hamiltonian BRST in Section \ref{sec:hamBRST}. Now as $Q$ is bounded we have:
\begin{lemma}\label{lm:propbsd}
Let $Q$ be as in Definition \eqref{df:fagbbrst}. Then $Q$ generates a bounded superderivation $\drb:\cA \to \cA$:
\begin{align*}
\drb(A):=\sbr{Q}{A}=QA-\gamma(A)Q,
\end{align*}
such that:
\begin{itemize}
\item[(i)] $\drb^2=0$.
\item[(ii)] $\gamma \circ \drb \circ \gamma= -\drb$. 
\item[(iii)] $\drb(A)^*=- \alpha \circ \drb \circ \alpha \circ \gamma(A^*)$. This identity encodes that $\drb(A^{\dag})=\drb(A)^{\dag}$.
\end{itemize}
\end{lemma}
\begin{proof}
First, as $\cA$ is a $C^*$-algebra $\norm{\drb(A)}\leq \norm{Q}\norm{A}+\norm{\gamma(A)}\norm{Q}=2\norm{Q}\norm{A}$, hence $\drb$ is bounded as a linear map $\drb:\cA \to \cA$.

\pfit(i) Let $A \in \cA$. Then
\begin{align*}
\drb^2(A)=\drb(QA-\gamma(A)Q)=&\;Q^2A-\gamma(QA)Q-Q\gamma(A)Q+\gamma(\gamma(A)Q)Q,\\
=&\; Q\gamma(A)Q-Q\gamma(A)Q-AQ^2,\\
=&\;0
\end{align*}
where used  $\gamma^2=\iota$, $\gamma(Q)=-Q$ and $Q^2=0$.

\pfit(ii):  Follows as $Q \in \cA^{-} \Rightarrow \gamma(Q)=-Q$. 

\pfit(iii): We calculate using $\alpha(Q^*)=Q$ and $\alpha^2=\iota$ that:
\[
\drb(A)^*=A^*Q^*-Q^*\gamma(A^*)=-\alpha( Q\alpha(\gamma(A^*))-\alpha(A^*)Q)=- \alpha \circ \drb \circ \alpha \circ \gamma(A^*).
\] 
\end{proof}
We now want to select states from which we can construct structures as in Chapter \ref{ch:GenStruct}. In the usual heuristic setup, the extended state space is the original state space tensored with a ghost state space. In terms of states on $\cA$ this corresponds to the set:
\[
\fS_T:=\{ \omega \in \fS(\cA) \,|\, \omega=\omega_1 \otimes \omega_2, \; \omega_1 \in \fS(\cF), \omega_2\in \fS_g \},
\]
where $\fS_g$ is definied in Definition \ref{df:ghst}, and so $\omega_2 \in \fS_g$ implies that $\omega_2 \circ \alpha' =\omega$ hence $\cH_{\omega_2}$ is a Krein space with fundamental symmetry, $J_{\omega_2}$, implementing $\alpha$ and $J_{\omega_2}\Omega_{\omega_2}=\Omega_{\omega_2}$. Thus $\Omega_{\omega_2}$ is positive in the Krein inner product on $\cH_{\omega_2}$ (cf. Definition \eqref{df:csksrepterm}). Assuming $\omega_1\circ \beta =\omega_1$ gives analogous Krein structures for $\cH_{\omega_1}$ (cf. Definition \eqref{df:csksrepterm}). We choose states that have GNS-cyclic vectors in $\ker \pi_{\omega}(Q)$. 
\begin{definition}\label{df:bsgsbdd}
The BRST-states $\fS_{\drb}\subset \fS_T$ are states of the form $\omega_1 \otimes \omega_2$, where $\omega_1 \circ \beta= \omega_1$, $\omega_2 \in \fS_{g}$ and 
\[
\omega(AQ)=0 \qquad \forall A \in \cA
\]
\end{definition}
Since, for $\omega \in \fS_{\drb}$ then $\omega \circ \alpha = \omega$ and so $\omega(AQ^*)= \omega(\alpha(A)Q)=0$ for all $A \in \cA$, hence $\omega \in \fS_{\drb}$ corresponds to a state such that 
\[
\Omega_{\omega} \in \ker \pi_{\omega}(Q) \cap \ker \pi_{\omega}(Q^*)= \ker \pi_{\omega}(\{Q, Q^* \}).
\]
We can rephrase the $\ker \pi_{\omega}(Q)$ condition in terms of $\delta$ alone \cite{Hendrik1991}:-
\begin{lemma}\label{lm:bQkdel}
Let $\omega=\omega_1 \otimes \omega_2$, where $\omega_1 \circ \beta= \omega_1$, $\omega_2 \in \fS_{g}$. Then the following are equivalent:
\begin{itemize}
\item[(i)] $\omega(AQ)=0 \quad \forall A \in \cA$.
\item[(ii)] $\omega(\drb(A))=0$ $\forall A \in \cA$.
\item[(iii)] $ \omega(B)=0$ $\forall B \in (\ran \drb \cap \ran \drb^*)$.
\end{itemize}
\end{lemma}
\begin{proof}
Let $\omega(AQ)=0$ $\forall A \in \cA$ then,
\[
\omega(\drb(A))= \omega(QA)= - \overline{\omega(A^*Q^*)}= - \overline{\omega(\alpha(A^*)Q)}=0,
\]
which proves one direction for both equivalence statements.

Conversely, let  $\omega(\drb(A))=0\, \forall A \in (\ran \drb \cap \ran \drb^*)$. We have that $\Delta= QQ^*+Q^*Q= \drb(Q^*) \in (\ran \drb \cap \ran \drb^*)$, hence $0=\omega(\drb(Q^*))=\omega( QQ^*+Q^*Q)=0$, and hence $\omega(QQ^*)=\omega(Q^*Q)=0$. Therefore by Cauchy-Schwartz we get that $|\omega(AQ)|\leq \omega(AA^*)\omega(Q^*Q)=0$, and so $\omega(AQ)=0$ $\forall A \in \cA$. Therefore (ii) $\Rightarrow$ (iii) $\Rightarrow$ (i) and we are done.
\end{proof}
With the appropriate states chosen, we would like to describe the structures of Chapter \ref{ch:GenStruct} algebraically. 
\begin{definition}\label{df:bdspdec}
\begin{itemize}
\item[(i)] Define the representation $\pi_{\drb}:\cA \to B(\cH_{\drb})$ by:
\[
\cH_{\drb}:=\bigoplus\{ \cH_{\omega}\,|\, \omega \in \fS_{\drb}\}, \qquad \pi_{\drb}:\bigoplus \{ \pi_{\omega}\,|\, \omega \in \fS_{\drb}\},  
\]
and denote the Hilbert inner product on $\cH_{\drb}$ by $\ip{\cdot}{\cdot}_{\cH_{\drb}}$. 
\item[(ii)] We have   $\omega \circ \alpha =\omega$ for all $\omega\in \fS_{\drb}$, hence $\alpha$ is unitarily implemented in each $\pi_{\omega}$ so $\alpha$ is unitarily implemented in $\cH_{\drb}$. Denote the implementer for $\alpha$ in $\pi_{\drb}$ by $J^{\drb}$. As $J^{\drb}|_{\cH_{\omega}}=J_{\omega}$ it follows from $\alpha^2=\iota$ that $(J^{\drb})^2=\one$ and $(J^{\drb})^*=J^{\drb}$. By lemma \eqref{lm:JKsp} $\cH_{\drb}$ is a Krein space with fundamental symmetry $J^{\drb}$ and indefinite inner product $\iip{\cdot}{\cdot}_{\cH_{\drb}}:=\ip{\cdot}{J^{\drb}\cdot}_{\cH_{\drb}}$. 
\item[(iii)] Let $\cH_{\drb}=\cH^d_{\drb}\oplus \cH^s_{\drb}\oplus \cH^p_{\drb}$ and $\cH_{\omega}=\cH^d_{\omega}\oplus \cH^s_{\omega}\oplus \cH^p_{\omega}$ be the \emph{dsp}-decompositions with respect to $\pi_{\drb}(Q)$ and $\pi_{\omega}(Q)$ where $\omega \in \fS_{\drb}$ (cf. Theorem \eqref{th:Hdsp1}). Let $P^{k}_{j}$, $k=\omega, \drb$, $j=d,s,p$ be the corresponding projections on $\cH_{\drb},\cH_{\omega}$. 
\item[(iv)] Let $\omega \in \fS_{\drb}$. For $j=\drb, \omega$ define,
\[
\cH^{BRST}_{phys,j}:=\ker \pi_{j}(Q)/\cH^d_{j}.
\]
and let $\vp_{j}:\ker \pi_{j}(Q) \to  \cH^{BRST}_{phys,j}$ be the factor map. Denote $\hat{\psi}:=\vp_{\drb}(\psi)$ for $\psi \in \cH_{\drb}$ and $\hat{\psi}{}_{\omega}=\vp_{\omega}(\psi)$ for $\psi \in \cH_{\omega}$. 
\end{itemize}
\end{definition}
\begin{rem}\label{rm:obvious}
As $\cH_{\drb}=\oplus_{\omega \in \fS_{\drb}}\cH_{\omega}$ it is obvious that $\cH^{j}_{\drb}=\oplus_{\omega \in \fS_{\drb}}\cH^{j}_{\omega}$ for $j=d,s,p$, and hence $P^{\drb}_{j}=\oplus_{\omega \in \fS_{\drb}}P^{\omega}_{j}$, $j=d,s,p$.
\end{rem}

To get the spatial structures of Chapter \ref{ch:GenStruct}:
\begin{proposition}\label{pr:cbiksp}
We have that $\cH^{BRST}_{phys,\drb}$ is a Krein space with indefinite inner product:
\begin{equation*}
\iip{\hat{\psi} }{\hat{\xi}}_p:=\iip{\psi}{\xi}_{\cH_{\drb}}=\iip{P^{\drb}_s\psi}{P^{\drb}_s\xi }_{\cH_{\drb}},
\end{equation*} 
fundamental symmetry $J^{\drb}_p\hat{\psi}:=\vp_{\drb}(J^{\drb}P^{\drb}_s\psi)$ and Hilbert inner product,
\begin{equation}\label{eq:BRSTpHilip2}
\ip{\hat{\psi}}{\hat{\xi}}_p:=\iip{\hat{\psi} }{J^{\drb}_p\hat{\xi}}_p= \ip{P^{\drb}_s\psi}{P^{\drb}_s\xi }_{\cH_{\drb}},
\end{equation} 
and norm $\norm{\hat{\psi}}_p:=\ip{\hat{\psi}}{\hat{\psi}}_p^{1/2}$. The space $\cH^{BRST}_{phys,\drb}$ is a Hilbert space with respect to the inner product $\iip{\hat{\psi}}{\hat{\xi}}_p$ if and only if $J_p^{\drb}=\one$ if and only if \begin{equation}\label{eq:ppossub}
J^{\drb}P^{\drb}_s=P^{\drb}_s,
\end{equation}
Hence the physicality condition for the abstract BRST system is equation \eqref{eq:ppossub}.
\end{proposition}
\begin{proof}
As $Q=\alpha(Q^*)=Q^{\dag}$ we have that $\pi_{\drb}(Q)$ is Krein self adjoint with respect to $\ip{\cdot}{\cdot}_{\cH_{\drb}}$. Hence, we have that $\pi_{\drb}(Q)$ and $\cH_{\drb}$ satisfy the hypothesis of Theorem \eqref{th:Hdsp}, and so the Proposition follows from Theorem \eqref{th:Hdsp}, Definition \eqref{df:opbrsphsp} and  lemma \eqref{lm:phspksp}.
\end{proof}

To get the algebraic structures as in Chapter \ref{ch:GenStruct} we define:
\begin{definition}\label{df:bddbrshom}
Define the linear map:
\begin{align*}
\Phi^{\drb}_s:\cA \to B(\cH_{\drb}),\qquad \text{by} \qquad \Phi^{\drb}_s(A):=P^{\drb}_s\pi_{\drb}(A)P^{\drb}_s,
\end{align*}
By lemma \eqref{lm:alkerhom} we have that $\Phi^{\drb}_s$ is an algebra homomorphism on $\ker \drb\subset \cA$. We define the \emph{BRST-physical algebra} as
\[
\cP_0^{BRST}:=\ker \drb / (\ker \drb \cap \ker \Phi^{\drb}_s)\cong \Phi^{\drb}_s(\ker \drb).
\]
Let the factor map be $\tau: \ker \drb \to \cP_0^{BRST}$, and denote $\hat{A}:=\tau(A)\in \cP_0^{BRST}$ for $A \in \ker \drb$.  
\end{definition}
\begin{rem} The above definition corresponds to the alternative definition of the phsyical algebra in Subsection \ref{sbs:alalg}. To connect with the cohomological definition of the physical algebra we extend $\pi_{\drb}(\cA)$ to $\cAe=\salg{\cA, P^{\drb}_s, P^{\drb}_p,P^{\drb}_d, \pi_{\drb}(Q), K}$ as in Section \eqref{sbs:AbcuQ}. Then $(\ker  \drb \cap \cA)/ (\ran \drb \cap \cA) \cong {\cP}_{\omega}^{BRST}$ by Theorem \eqref{pr:krdel} (see also Remark \eqref{rm:alus}).
\end{rem}
As $\ker \drb$ is not in general a $*$-subalgebra of $\cA$ (\eg $Q\in \ker \drb$ but $\drb(Q^*)=Q^*Q+QQ^* \neq 0$ for $Q \neq 0$) we do not get that $\cP^{BRST}$ is a $C^*$-algebra using the usual factor norm. To get a natural norm for $\cP_0^{BRST}$ we use the norm on $\cH^{BRST}_{phys,\drb}$. 
\begin{proposition}\label{pr:repban}
Define the representations,
\begin{align*}
\pi_{j, p}:\cP_0^{BRST}&\to B(\cH^{BRST}_{phys,j}),\\
\pi_{j, p}( \hat{A})\hat{\psi}:&= \widehat{\pi_{j}(A)\psi}, \qquad \psi \in \ker \pi_j(Q)
\end{align*}
where $j=\omega, \delta$ and $\omega \in \fS_{\drb}$. Define the seminorm on $\cP_0^{BRST}$ by: 
\[
\norm{\hat{A}}_p:=\norm{\pi_{\drb, p} ( \hat{A})}_{B(\cH^{BRST}_{phys,\drb})}.
\]
Then $\norm{\hat{A}}_p$ is a norm on $\cP_0^{BRST}$,
\begin{equation}\label{eq:nmfact}
\norm{\hat{A}}_p= \norm{\Phi^{\drb}_s(A)}_{\cH_{\drb}},
\end{equation}
and $\cP^{BRST}:=\overline{\cP_0^{BRST}}$ is a Banach algebra where closure is with respect to $\norm{\hat{A}}_p$. Moreover $\cP^{BRST}\cong \overline{\Phi^{\drb}_s(\ker \drb)}^{\cB(\cH_{\drb})}$ where the isomorphism is isometric. Thus $\pi_{\drb,p}$ is a faithful representation of $\cP^{BRST}$ and all calculations can be done in this representation.
\end{proposition}
\begin{proof}
First $A \in \ker \drb \Rightarrow A \ran \pi_{\drb}(Q) \subset \ran \pi_{\drb}(Q)$, hence
\begin{equation}\label{eq:Akerrz}
A \in \ker \drb \Rightarrow \pi_{\drb}(A)\cH^d_{\drb}\subset \cH^{d}_{\drb} \Rightarrow P^{\drb}_s \pi_{\drb}(A) \ker \pi_{\drb}(Q)=P^{\drb}_s \pi_{\drb}(A)P^{\drb}_s \ker \pi_{\drb}(Q)
\end{equation}
Using this we see that $\norm{\hat{A}}_p$ is a norm by the calculation:
\begin{align*}
\norm{\hat{A}}_p= \sup_{\norm{\hat{\psi}}_p=1}\norm{\pi_{\drb, p}( \hat{A})\hat{\psi}}_p=&\; \sup\{ \norm{ P^{\drb}_s \pi_{\drb}(A)\psi}_{\cH_{\drb}}\;|\; \psi \in \ker  \pi_{\drb}(Q),\,\norm{{P^{\drb}_s} \psi}_{\cH_{\drb}}=1\} \\
=&\;\sup_{\norm{{P^{\drb}_s} \psi}_{\cH_{\drb}}=1}\norm{ P^{\drb}_s \pi_{\drb}(A)P^{\drb}_s\psi}_{\cH_{\drb}},\\
=&\;\norm{\Phi^{\drb}_s(A)}_{\cH_{\drb}},
\end{align*}
for $A \in \ker \drb$, where we have used that $\norm{\hat{\psi}}_p=\norm{P^{\drb}_s\psi}_{\cH_{\drb}}$ for $\psi \in \ker \pi_{\drb}(Q)$ by equation \eqref{eq:BRSTpHilip2}  for the second equality, and  equation \eqref{eq:Akerrz} for the third. Therefore $\norm{\hat{A}}_p=0$ iff $\norm{\Phi^{\drb}_s(A)}_{\cH_{\drb}}=0$ iff $\Phi^{\drb}_s(A)=0$  iff $\hat{A}=0$, and so $\norm{\cdot}_p$ is a norm. Let $A,B \in \ker \drb$. As $\norm{P^{\drb}_s}_{B(\cH_{\drb})}=1$ we have that:
\[
\norm{\hat{A}\hat{B}}_p=\norm{\Phi^{\drb}_s(AB)}_{\cH_{\drb}}\leq\norm{\Phi^{\drb}_s(A)}_{\cH_{\drb}} \norm{\Phi^{\drb}_s(B)}_{\cH_{\drb}}=\norm{\hat{A}}_p\norm{\hat{B}}_p.
\]
Thus $\cP^{BRST}$ is a Banach algebra where closure is with respect to $\norm{\hat{A}}_p$. The last isomorphism comes from equation \eqref{eq:nmfact}.
\end{proof}
Now for $A \in \ker \drb$ we have that $A^*$ need not be in $\ker \drb$ (\eg $Q$). Hence we have to be careful of how we define a $*$-involution on  $\cP^{BRST}$ (cf. lemma \eqref{lm:facsal}).
\begin{proposition}\label{pr:bddcshinv}
We have:
\begin{itemize}
\item[(i)] $(\ker \drb)^{\dag}=\ker \drb$ and $(\ker \drb \cap \ker \Phi^{\drb}_s)^{\dag}=(\ker \drb \cap \ker \Phi^{\drb}_s)$. Hence $\dag$ on $\ker \drb$ factors to the $\dag$-involution on $\cP^{BRST}$ which coincides with the $\dag$-involution with respect to the representation $\pi_{\drb,p}$, \ie
\[
\pi_{\drb,p}(\widehat{A^{\dag}})\hat{\psi}=\pi_{\drb,p}(\hat{A})^\dag\hat{\psi}
\]
for all $A \in \ker \drb$ where $\pi_{\drb,p}(\hat{A})^{\dag}$ is the the adjoint of $\pi_{\drb,p}(\hat{A})$ with respect to the inner product $\iip{\cdot}{\cdot}_p$. Furthermore,
\[
\cP^{BRST} \cong \overline{\Phi_s^{\drb}(\ker \drb)}
\]
where the above is an isometric $\dag$-isomorphism.
\item[(ii)] Let $\cM \in \ker \drb$ be a subalgebra such that $\Phi^{\drb}_s(\cM)=\Phi^{\drb}_s(\cM^*)$. Given $A \in \cM$, define
\begin{equation}\label{eq:dfhadpa}
\hat{A}^*:=\hat{B},
\end{equation} 
where $\Phi^{\drb}_s(A^*)=\Phi^{\WW}_s(B)$ for some $B \in \cM$. This defines an involution on $\cM/ (\cM \cap \ker \Phi^{\drb}_s)$ such that $\overline{\cM/ (\cM \cap \ker \Phi^{\drb}_s)}$ is a $C^*$-algebra where closure is with respect to $\norm{\cdot}_p$.
\item [(iii)] Let the physicality condition $J^{\drb}P^{\drb}_s=P^{\drb}_s$ hold (equation \eqref{eq:ppossub}). Then $\ker \drb$ satisfies the conditions on $\cM$ in (ii), hence $\cP^{BRST}$ is a $C^*$-algebra with respect to the norm $\norm{\cdot}_p$. Moreover, the $\dag$-involution from (i) and $*$-involution from (ii) coincide.
\end{itemize}
\end{proposition}
\begin{proof}
(i): $(\ker \drb)^{\dag}=\ker \drb$ follows from $Q^{\dag}=Q$. By Theorem \eqref{th:Hdsp}, we have $[P^{\drb}_s, J^{\drb}]=0$, hence $(P^{\drb}_s)^{\dag}=J^{\drb}P^{\drb*}_sJ^{\drb}=P^{\drb}_s$ and so it follows that $(\ker \drb \cap \ker \Phi^{\drb}_s)$ is a $\dag$-subalgebra of  $\ker \drb$. Hence the involution $\dag$ factors to an involution on $\cP^{BRST}$. Moreover, it coincides with the $\dag$-involution with respect to the representation $\pi_{\drb,p}$ which can be seen by the calculation:
\[
\iip{\pi_{\drb,p}(\widehat{A^{\dag}})\hat{\psi}}{\hat{\xi}}_p=\iip{\widehat{\pi_{\drb}(A^{\dag})\psi}}{\hat{\xi}}_p=\iip{\pi_{\drb}(A^{\dag})\psi}{\xi}=\iip{\psi}{\pi_{\drb}(A)\xi}=\iip{\hat{\psi}}{\pi_{\drb,p}(\hat{A})\hat{\xi}}_p,
\]
for all $\psi,\xi \in \cH^{BRST}_{phys, \drb}$ and all $A \in \ker \drb$, where we have used that by definition $\iip{\hat{\psi}}{\hat{\xi}}_p=\iip{\psi}{\xi}$ for all $\psi,\xi \in \ker \pi_{\drb}(Q)$ (cf. Proposition \eqref{pr:cbiksp}). 

Let $A \in \ker \drb$. Then using $[P^{\drb}_s, J^{\drb}]=0$,
\[
\Phi_s^{\drb}(A^{\dag})=\Phi_s^{\drb}(\alpha(A)^*)=P_s^{\drb}J_{\drb}\pi_{\drb}(A)^*J_{\drb}P_s^{\drb}=(J_{\drb}P_s^{\drb}\pi_{\drb}(A)P_s^{\drb}J_{\drb})^*=\Phi_s^{\drb}(A)^{\dag}
\]
where we used $J_{\drb}^*=J_{\drb}$ in the last line. Combining this with Proposition \eqref{pr:urepban} gives that there is a isometric $\dag$-isomorphism such that $\cP^{BRST} \cong \overline{\Phi_s^{\drb}(\ker \drb)}$.  
\smallskip
\noindent(ii): Let $A\in \cM$ then by assumption there exists $B\in \cM$ such that $\Phi_s(A^*)=\Phi_s(B)$. Now by equation \eqref{eq:Akerrz} we have
\[
A \in \ker \drb \Rightarrow  P^{\drb}_s \pi_{\drb}(A) \ker \pi_{\drb}(Q)=P^{\drb}_s \pi_{\drb}(A)P^{\drb}_s \ker \pi_{\drb}(Q),
\]
hence for all $\xi,\psi \in \ker \pi_{\drb}(Q)$ we have by equation \eqref{eq:BRSTpHilip2},
\begin{align}
\ip{\hat{\xi}}{\pi_{\drb, p}(\hat{A}) \hat{\psi}}_p= & \;\ip{P^{\drb}_s\xi}{P^{\drb}_s \pi_{\drb}(A)\psi}_{\cH_{\drb}},\notag \\
= & \;\ip{P^{\drb}_s\xi}{P^{\drb}_s \pi_{\drb}(A)P^{\drb}_s\psi}_{\cH_{\drb}},\notag \\
= & \;\ip{\Phi^{\drb}_s(A^*)\xi}{\psi}_{\cH_{\drb}},\notag \\
= & \;\ip{\Phi^{\drb}_s(B)\xi}{\psi}_{\cH_{\drb}}, \notag \\
= & \;\ip{\pi_{\drb, p}(\hat{B})\hat{\xi}}{\hat{\psi}}_p,\label{eq:ipad}
\end{align}
and hence 
\[
\pi_{\drb, p}(\hat{A})^*=\pi_{\drb, p}(\hat{B})=\pi_{\drb, p}((\hat{A})^*)\in \cM/ (\cM \cap \ker \Phi^{\drb}_s)
\]
where we used equation \eqref{eq:dfhadpa} and that $B \in \cM$. This shows that $\cM/ (\cM \cap \ker \Phi^{\drb}_s)$ is a $*$-algebra. Furthermore, Proposition \eqref{pr:repban} gives that $\pi_{\drb, p}:\cM/ (\cM \cap \ker \Phi^{\drb}_s)\to B(\cH^{BRST}_{phys, \drb})$ is an isometric isomorphism, and so we have proved that $\pi_{\drb, p}$ is a $*$-isometric isomorphism. Hence as $B(\cH^{BRST}_{phys, \drb})$ is a $C^*$-algebra so is  $\overline{\cM/ (\cM \cap \ker \Phi^{\drb}_s)}$.

\pfit (iii): Now $A \in \ker \drb$ implies $A^{\dag} \in \ker \drb$ and by equation \eqref{eq:ppossub} $\Phi^{\drb}_s(A^{\dag})= P^{\drb}_s J^{\drb}\pi_{\drb}(A^*)J^{\drb} P^{\drb}_s=P^{\drb}_s \pi_{\drb}(A^*) P^{\drb}_s=\Phi^{\drb}_s(A^{*})$. Therefore we can apply (ii) with $\cM=\ker \drb$ and $\cM_2=(\ker \drb)^{\dag}=\ker \drb=\cM$ to get that ${\cP^{BRST}}=\overline{\ker \drb/ (\ker \drb \cap \ker \Phi^{\drb}_s)}$ with norm $\norm{\cdot}_p$ is a $C^*$-algebra.  Moreover by Proposition \eqref{pr:bddcshinv} (i) and the defining equation \eqref{eq:dfhadpa}, $\hat{A}^{\dag}=\widehat{A^{\dag}}=\hat{A}^*$.  
\end{proof}
\begin{rem} \begin{itemize}
\item[(i)] Note that the above Proposition does \emph{not} assume that $\cM$ is a $C^*$-algebra. If is factoring out by $(\cM \cap \ker \Phi^{\drb}_s)$ and using the the Hilbert $*$-involution coming from the BRST physical space $\cH^{BRST}_{phys,\drb}$ that gives a $C^*$-algebra.
\item[(ii)] Proposition \eqref{pr:bddcshinv} (iii) shows that equation \eqref{eq:ppossub} is a good physicality condition as it ensures that the $\dag$-involution factors to a $C^*$-involution on the physical algebra $\cP^{BRST}$.
\end{itemize}
\end{rem}
\begin{eje}\label{ex:pscs} Extend $\cA$ by defining $\tilde{\cA}=\osalg{P^{\drb}_s, \pi_{\drb}(\cA)}$. We drop the notation $\pi_{\drb}$ here, and extend all the structures such as $\drb$ and $\alpha$ to $\tilde{\cA}$ by using $Q$ and $J^{\drb}$, etc.

As $Q P^{\drb}_s=P^{\drb}_sQ=0$ we have that $\Phi^{\drb}_s(\tilde{\cA}) \subset \ker \drb$ hence we have 
\[
\Phi^{\drb}_s(\tilde{\cA}) \subset\Phi^{\drb}_s( \ker \drb) \subset \Phi^{\drb}_s(\tilde{\cA}),
\]
hence $\Phi^{\drb}_s( \ker \drb) = \Phi^{\drb}_s(\tilde{\cA})$. By taking adjoints we see 
\[
\Phi^{\drb}_s( (\ker \drb)^*)= \Phi^{\drb}_s(\tilde{\cA}^*)=\Phi^{\drb}_s(\tilde{\cA})=\Phi^{\drb}_s( \ker \drb)
\]
and so by extending $\cA$ by $P^{\drb}_s$, we get from Proposition \eqref{pr:bddcshinv} (ii) that $\cP^{BRST}$ is a $C^*$-algebra with norm $\norm{\cdot}_p$. Notice again that we did not assume that $\ker \drb$ was $*$-algebra to begin with.
\end{eje}
Summarising the above discussion we get:
\begin{theorem}\label{th:bdbrststruc}
Let $\cA,\alpha,Q,\drb,\fS_{\drb}, \pi_{\drb}:\cA \to B(\cH_{\drb})$ be as in Definitions \eqref{df:fagbbrst},  \eqref{df:bsgsbdd} and  \eqref{df:bdspdec}, and let,
\[
\Phi^{\drb}_s(A):=P^{\drb}_s\pi_{\drb}(A) P^{\drb}_s, \qquad A \in \cA.
\] 
as in Definition \eqref{df:bddbrshom}. Then:
\begin{itemize}
\item[(i)]$\Phi^{\drb}_s$ is a homomorphism on $\ker \drb$ and on $(\ker \drb)^*$. 
\item[(ii)]Let  $\cP_0^{BRST}=\ker \drb / (\ker \drb \cap \ker \Phi^{\drb}_s)$ and $\pi_{\drb, p}:\cP_0^{BRST} \to B(\cH^{BRST}_{phys,\drb})$ be as in Proposition \eqref{pr:repban}. Then $\cP_0^{BRST}$ has norm 
\[
\norm{\hat{A}}_p:=\norm{\pi_{\drb, p} (\hat{A})}_{B(\cH^{BRST}_{phys,\drb})}=\norm{\Phi^{\drb}_s(A)}_{\cH_{\drb}}, \qquad \forall A \in \ker \drb
\]
with respect to which $\cP^{BRST}:=\overline{\cP_0^{BRST}}$ it is a $\dag$-Banach algebra and we have a $\dag$-isometric isomorphism such that ${\cP^{BRST}}\cong \overline{\Phi^{\drb}_s(\ker \drb)}^{\cB(\cH_{\drb})}$ 
\item[(iii)] If $\Phi^{\drb}_s(\ker \drb)= \Phi^{\drb}_s((\ker \drb)^*)$, then $\cP^{BRST}$ is a $C^*$-algebra with norm $\norm{\cdot}_p$ and involution denoted by $*$ as defined in equation \eqref{eq:dfhadpa}. 
\begin{enumerate}
\item When $J^{\drb}P_s^{\drb}=P_s^{\drb}$ we have that $\Phi^{\drb}_s(\ker \drb)= \Phi^{\drb}_s((\ker \drb)^*)$ and so $\cP^{BRST}$ is a $C^*$-algebra.  
\item When we extend $\cA$ to $\tilde{\cA}=\osalg{P^{\drb}_s, \cA}$ as we have that $\Phi^{\drb}_s(\ker \drb)= \Phi^{\drb}_s((\ker \drb)^*)$ and so $\cP^{BRST}$ is a $C^*$-algebra.  
\end{enumerate}
\end{itemize}
\end{theorem}
\begin{proof}
(i): Follows from lemma \eqref{lm:alkerhom} applied to  for $\pi_{\drb}(Q)$ and $\pi_{\drb}(Q^*)$. (ii) is Proposition \eqref{pr:repban} and Proposition \eqref{pr:bddcshinv}. (iii) follows from Proposition \eqref{pr:bddcshinv} and Example \eqref{ex:pscs}.
\end{proof}
We may want a more intrinsic description of $\cP^{BRST}$, \ie one that does not involve $\Phi^{\drb}_s$.Define,
\[
\fS^{\drb}_{S}:=\{ \omega_{\psi}\,|\, \psi \in \cH^{\drb}_s=\ker \pi_{\drb}(\Delta) \},\qquad  \cT_{S}:=\cap_{\omega \in \fS^{\drb}_S}\ker \omega,
\]
where $\omega_{\psi}$ denotes the vector state of $\psi \in \cH_{\drb} \backslash \{0\}$, \ie $\omega_{\psi}(\cdot):=\ip{\psi}{\cdot \,\psi}/ \norm{\psi}^2$. 
The above definition implies that if $A \in \cT_S$ then $\pi_{\omega}(A)|_{\cH^{\omega}_s}=0$ for $\omega\in \fS_{\drb}$ (\cite{Con1985} II Proposition 2.15 p34). Hence if $\tau:\ker \drb \to \cP_0^{BRST}$ is the factor map $(\pi_{\drb, p}\circ \tau)(\ker \drb \cap \cT_{S})=\{0\}$, \ie  $\ker \drb \cap \cT_{S}$ are trivial BRST observables. We will see that these are all the trivial observables. 
\begin{proposition}\label{pr:bddinttriv}
We hav, 
\[
\ker \drb \cap \cT_{S}=\ker \drb \cap \ker \Phi^{\drb}_s,
\]
hence $\ker \drb \cap \cT_{S}$ is a two-sided ideal in $\ker \drb$ and,
\[
\cP_0^{BRST} =\ker \drb/ (\ker \drb \cap \cT_S). 
\]
\end{proposition}
\begin{proof}
Let $A \in (\ker \drb\cap \cT_S)$, $\omega \in \fS_{\drb}$, $\psi \in \cH_{\omega}$. Since $P_s^{\omega}\psi \in \cH^{\omega}_s$ we have $\omega_{P_s^{\omega}\psi}\in \fS^{\drb}_S$. As $A \in \cT_S$ we have $\ip{\psi}{P_s^{\omega}\pi_{\omega}(A)P_s^{\omega} \psi}=\omega_{P_s^{\omega}\psi}(A)=0$ and as $\psi\in \cH_{\omega}$ was arbitrary $P_s^{\omega}\pi_{\omega}(A)P_s^{\omega}=0$. As $P^{\drb}_{s}=\oplus_{\omega \in \fS_{\drb}}P^{\omega}_{s}$ (cf. Remark \eqref{rm:obvious}) it follows that $A \in \ker \Phi^{\drb}_s$.

Conversely, let $A \in \ker \Phi_s^{\drb}$ and $\omega_{\psi} \in \fS^{\drb}_{S}$. $A \in \ker \Phi_s^{\drb} \Rightarrow P^s_{\omega}\pi_{\omega}(A) P^s_{\omega}=0$ and so $\omega_{\psi}(A)=\omega_{P_s^{\omega}\psi}(A)= \ip{\psi}{P_s^{\omega}\pi_{\omega}(A)P_s^{\omega}\psi}_{\omega}=0$. Hence $\ker \drb \cap \ker \Phi_s^{\drb}= \ker \drb \cap \cT^{\drb}_S$. 
\end{proof}
\begin{eje}\label{ex:csfham}
We are now in a position to formulate Hamiltonian BRST with finite constraints (cf.Section \ref{sec:hamBRST}) in a $C^*$-algebra format. Let $(\cA_0,\cC)$ be a quantum system with constraints (cf. Appendix \ref{app:Tp}) where $\cC=\{G_{j}\}_{j=1}^{\MM}$ is a set of finite set of self-adjoint constraints that form a Lie algebra,
\begin{equation*}
[G_a,G_b]=iC^c_{ab}G_c, \qquad C^c_{ab} \in \mathbb{R}.
\end{equation*}
where $C^c_{ab}$ is antisymmetric in all indices. Let $\fS_D(\cA_0)$ be the set of Dirac states and let $P \in \cA_0''$ be the open projection from Theorem \eqref{th:vnT}. Therefore $P$ is such that $\omega(P)=0$ iff $\omega \in \fS_D(\cA_0)$. 

Let $\cH$ be a complex Hilbert space with $\dim(\cH)=2\MM$. If we recall the construction in lemma \eqref{lm:ghtsfext} then given any decomposition 
\[
\cH=\cH_1 \oplus \cH_2
\]
such that $\dim(\cH_2)=\MM$, there exists a unitary $J\in B(\cH)$ with $J^2=\one$ such that $J\cH_1=\cH_2$ and $\cH$ has Krein inner product $\iip{\cdot}{\cdot}:=\ip{\cdot}{J\cdot}$, hence $\cH_1 [\oplus] \cH_2$.

Given such a choice of decomposition $\cH=\cH_1\oplus \cH_2$ and $J$, then $\cH$ has all the structure of the ghost test function space $\cH$ in Section \ref{sec:ghost}.

Let $\rga(\cH_2)$ be the restricted ghost algebra as in Definition \eqref{df:GA} and let $\gh_j:=\gh(f_j)$, $\cgh_j:=\cgh(f_j)$. Let the BRST-Field Algebra be $\cA=\cA_0\otimes \rga$ with $\alpha=\iota \otimes \alpha'$ where $\alpha' \in \aut(\cA_g(\cHJ))$ is the automorphism associated to the ghost-Krein structure (cf. Defintion \eqref{df:fagbbrst}). Note that the tensor norm on $\cA$ is unique as $\rga$ is finite dimensional. Define the BRST charge as in equation \eqref{eq:hnonabQ}:
\begin{equation}
Q:=G_a \otimes \gh_a - (i/2) C^c_{ab}\one \otimes \gh_a \gh_b \cgh_c,
\end{equation}
and we have as in \eqref{eq:hHam2nilQ} that $Q^{\dag}=\alpha(Q^*)$, $Q^2=0$, $Q\in \cA^{-}$. Then we have all the structures as in Theorem \eqref{th:bdbrststruc} as well as the representations $\pi_{\omega, p}:\cP^{BRST} \to B(\cH^{BRST}_{\omega, phys})$ as in Proposition \eqref{pr:repban} for $\omega \in \fS_{\drb}$

\begin{lemma}\label{lm:abcsst}
In the case of abelian constraints, \ie $C^c_{ab}=0$ $\forall a,b,c,$ we have:
\begin{itemize}
\item[(i)] $\fS_{\drb}=\fS_D\otimes \fS_g$,
\item[(ii)] $P^{\drb}_s=(\one -P) \otimes \one$,
\end{itemize}
where $P$ is the restriction to $\cH_{\drb}$ of the projection from Theorem \eqref{th:vnT}.
\end{lemma}
\begin{proof}
\noindent (i) Let $\omega=\omega_1 \otimes \omega_2 \in \fS(\cA)$ and $\omega_2 \in \fS_g$. Recall that $\Delta = \{Q^*, Q\}$. From lemma \ref{lm:bQkdel} we know that if $\omega \in \fS_{\drb}$ then $\omega(\drb(Q^*))=\omega(\Delta)=0$. Conversely, if $\omega(Q^*Q+QQ^*)=0$ then by positivity we have $\omega(Q^*Q)=0$ hence $\omega \in \fS_{\drb}$. Now, as in Subsection \ref{sbs:hhamexphsp}, we have $\Delta=\{Q^*,Q\}=G_jG_j\otimes \one$ and therefore for any $\omega=\omega_1\otimes \omega_2 \in \fS_{\drb}$ iff $\omega(\Delta)=\omega_1(G_bG_b)\otimes \omega_2(\one)=0$, hence $\omega_1 \in \fS_D$.

\medskip
\noindent{(ii)} Let $\omega=\omega_1 \otimes \omega_2 \in \fS_{\drb}$. Then $P^{\omega}_s$ is the projection onto
\[
\ker \pi_{\omega}(Q) \cap \ker \pi_{\omega}(Q^*)=\ker \pi_{\omega}(\Delta)= \bigcap_{j}\big(\ker \pi_{\omega_1}(G_j) \otimes \cH^g_{\omega_2}\big)= \ker \pi_{\omega}(P)\otimes \cH_g,
\] 
hence $P^{\drb}_s=(\one -P) \otimes \one$.
\end{proof}
From this it follows that is we restrict ourselves to $\ker \drb \cap (\cA_0 \otimes \one)$ then we get that the restricted BRST observables are equivalent to the traditional Dirac observables.
\begin{proposition}\label{pr:rabbdbral}
For abelian constraints we have:
\begin{itemize}
\item[(i)] $\ker \drb \cap (\cA_0\otimes \one)= \cC' \otimes \one$,
\item[(ii)] $ \Phi^{\drb}_s(\ker \drb \cap (\cA_0 \otimes \one))\cong (\cC'/(\cC'\cap \cD)))\otimes \one$
\end{itemize} 
where $\cC'$ is the commutant of $\cC'$ in $\cA_0$ and $\cD$ is defined as in Appendix \ref{app:Tp}, and $\cong$ is a $*$-isomorphism.
\end{proposition}
\begin{proof}
\noindent(i) Take $A \otimes \one \in \cA_0\otimes \one$. Then $\drb(A)=\sum_{j}[G_j,A]\otimes \gh_j$. As $\gh_j$ are linearly independent we get that $A\otimes \one \in \ker \drb$ iff $A \in \cC'$.

\medskip
\noindent{(ii)} Using lemma \eqref{lm:abcsst} we have that $\Phi^{\drb}_s(\ker \drb \cap (\cA_0 \otimes \one))=[(1-P)\cC'(1-P)] \otimes \one \cong (\cC'/(\cC'\cap \cD))\otimes \one$ by Theorem \eqref{th:vnT} (v). This is obviously a $*$-isomorphism.
\end{proof}

\end{eje}
\begin{rem}\label{rm:wkcom}
\begin{itemize}
\item[(i)] From lemma \eqref{lm:abcsst} we see that $\cH_{\omega}=\bigcap_{j}\big(\ker \pi_{\omega_1}(G_j)\big) \otimes \cH^g_{\omega_2}$ for $\omega=\omega_1\otimes \omega_2\in \fS_{\drb}$. As $\cA_g$ is simple we get that $\pi_{\omega_2}$ is faithful and we still get the MCPS problem for abelian BRST as in Subsection \ref{sbs:hhamexphsp} and neutrality problems as in Remark \eqref{rm:spfstuff}). It is not hard to see that this also follows for non-abelian constraints from equation \eqref{eq:nonabdel} in the same way as described in Subsection \ref{sbs:hhamexphsp}.
\item[(ii)] Proposition \eqref{pr:rabbdbral} tells us that if we restrict $\drb$ to our original algebra that the BRST method selects a physical algebra equivalent to the one obtained using the traditional Dirac method, \ie the commutant of the constraints. The problem with this approach is that it does not handle equivalent sets of constraints well because different sets of constraints which select the same set of physical states may have different commutants. This is the reason why the $T$-procedure defines the observable algebra as the abstract version of the \emph{weak commutant} of the constraints \cite{Hendrik2006} p100, and a significant advantage which this generalized Dirac method of constraints has over Hamiltonian BRST. 
\end{itemize}
\end{rem}

\section{Resolvent Algebra}\label{sec:RA}
To cast the structures of Chapter \ref{ch:BRSTQEM} into a $C^*$-algebraic setting we have several hurdles to overcome. The main difficulty is that the objects so far defined are unbounded operators, hence analytically difficult. Here we want to find a $C^*$-algebra that which encodes the algebraic information of the operators involved, \ie it can reproduce the above structures in the appropriate representations. We will use notations for inner product spaces and symplectic spaces as in Appendices \ref{ap:IIP},  \ref{ap:symp}.

The main relation we would like to model is the superderivation action, given in Definition \eqref{df:brstsd},
\begin{align*}
\delta(A(g)\otimes \one)&= i \one \otimes C( P_{\JJ}  ig), \qquad g \in \fD, \\
\delta(\one \otimes C(g))&= A(\mH P_{\LL} g)\otimes \one \qquad g \in \fD.
\end{align*}
These are awkward relations because:
\begin{itemize}
\item [(i)] They involve Krein-symmetric unbounded operators such as $A(g) \otimes \one$. Unboundedness gives domain technicalities and spectral theory for Krein-symmetric operators is more complicated and less well understood than for the Hilbert space case.
\item [(ii)] The equation $\delta(\one \otimes C(g))= A(\mH P_{\LL} g)\otimes \one$ maps bounded operators to unbounded operators, and so $\drb$ will be hard to  interpret even after the CCR relations have been cast into a bounded form. This problem has been previously analyzed for the supersymmetry superderivation (cf.  \cite{HendrikBuch2007}, instead leading to the Resolvent algebra. 
\end{itemize}

We start with a brief description of the structure of Resolvent Algebras of CCR's. Proofs can be found in \cite{HendrikBuch2007}. One observes that resolvents of self-adjoint operators can be used as `mollifiers' in the following sense. Let $A$ be a self-adjoint operator acting on a Hilbert space $\cH$ with domain $D(A)$. Since $A$ has real spectrum, $(i\lambda\one -A)^{-1} \in B(\cH)$ for $\lambda \in \mathbb{R}\backslash \{0\}$, and by functional calculus we get that $A(i\lambda\one -A)^{-1} =\overline{(i\lambda\one -A)^{-1}A}=i\lambda(i\lambda\one -A)^{-1}-\one \in B(\cH)$. Thus $M:=(i\lambda\one -A)^{-1}$ carries the information of $A$ in bounded form and `mollifies' $A$, \ie $AM \in B(\cH) \ni \overline{MA}$. Recall from \cite{DamGeo2004} that $A$ is \emph{affiliated} with a $C^*$-algebra $\cA\subset B(\cH)$ if $M=(i \lambda \one-A)^{-1} \in \cA$. It makes sense to look for a $C^*$-algebra $\cA \subset B(\cH)$ which contains the resolvents $(i\lambda \one -A)^{-1} $ of all self-adjoint operators $A$ which we need. 
The mollifying property will prove key ingredient for BRST in making sense of the identities $\drb(\one \otimes C(g))= A(\mH P_{\LL} g)\otimes \one$ as will be seen in Subsection \ref{sbs:sdI}.

A bosonic field is often described by operators satisfying the CCR as in \cite{BraRob21981}. That is given $(\fX,\sigma)$ a real symplectic space with symplectic $\sigma$ a field $\phi$ is a linear map from $\fX$ to a linear space of self-adjoint operators on some common dense invariant core $\cD$ in a Hilbert space $\cH$, satisfying the relations
\[
[\phi(f),\phi(g)]=i\smp{f}{g}{1}\one,\qquad \text{on} \quad \cD.
\]
A common way to encode the CCR-fields in a bounded way is to study the \emph{Weyl Algebra} \cite{BraRob21981} which is the algebra we get from exponentiating the $\phi(f)$, ie
\[
\osalg{\exp(i\phi(f))\,|\, f\in \fX}\subset B(\cH)
\]
This is can be defined abstractly and is commonly denoted by $\overline{\Delta(\\fX,\sigma)}$. It is not optimal for modelling the $\drb$ as it does not contain mollifiers \cite{HendrikBuch2007} Proposition 2.1. More useful is the algebra we get from taking resolvents of $\phi(f)$, ie
\[
\osalg{(i\lambda\one-\phi(f))^{-1}\,|\, f\in \fX,\,\lambda \in \mathbb{R}\backslash \{0\}}\subset B(\cH).
\]
Note that as $\phi(f)^*=\phi(f)$,  the spectrum of $\phi(f)$ is real and hence the resolvents $(i\lambda\one-\phi(f))^{-1}$  are well defined. 

The Resolvent Algebra can be defined abstractly as follws:
\begin{definition}\label{df:RA}
Given a symplectic space $(\fX,\sigma)$, define $\cR_0(\fX,\sigma)$ to be the universal $*$-algebra generated by the set $\{R(\lambda,f)\,|\,  \lambda \in \mathbb{R} \backslash \{0\},\,f\in \fX\}$ subject to the relations:
\begin{enumerate}
\item $R(\lambda,0)=-(i/\lambda)\one$,
\item $R(\lambda,f)^*=R(-\lambda,f)$,
\item $R(\lambda,f)=(1/\lambda)R(1,f/\lambda)$, \label{eq:R3}
\item $R(\lambda,f)-R(\mu,f)=i(\mu-\lambda)R(\lambda,f)R(\mu,f)$, \label{eq:R4}
\item $[R(\lambda,f),\,R(\mu,g)]=i\smp{f}{g}{}R(\lambda,f)R(\mu,g)^2R(\lambda,f)$,\label{eq:R5}
\item $R(\lambda,f)R(\mu,g)=R(\lambda+\mu,f+g)[R(\lambda,f)+R(\mu,g)+i\smp{f}{g}{} R(\lambda,f)R(\mu,g)^2R(\lambda,f)]$
\end{enumerate}
where $\lambda, \mu \in\mathbb{R} \backslash \{0\}$ and $f,g \in \fX$, and for (6) we require $\lambda +\mu \neq 0$. That is we start with the free unital algebra generated by $\{R(\lambda,f)\,|\,  \lambda \in \mathbb{R} \backslash \{0\},\,f\in \fX\}$  and factor out by the ideal generated by the relations (1) to (6). We also denote $\cR_0(\fX,\sigma)$ by $\cR_0$ when no confusion will arise.
\end{definition}

\begin{rem}
\begin{itemize}
	\item $\cR_0$ has non-trivial representations as can be seen by taking resolvents of the fields $\phi(f)$ in the Fock representation.
	\item let $\mu=-\lambda$ in equation \eqref{eq:R4}. Then we have that
	\begin{equation}\label{eq:R42}
	R(\lambda,f)-R(\lambda,f)^*=-2i\lambda R(\lambda,f)R(\lambda,f)^*
	\end{equation} 
\end{itemize}
\end{rem}

We want a norm on $\cR_0$ such that we can complete it to a $C^*$-algebra. Let $\fS$ be the set of functionals $\omega:\cR_0\to \mathbb{C}$ such that $\omega(A^*A)\geq 0$ and $\omega(\one)=1$. Then for each $\omega \in \fS$ its GNS-representations is bounded (\cite{HendrikBuch2007} Proposition 3.3) and hence we can define:
\begin{definition}
The universal representation $\pi_u:\cR_0 \to B(\cH_0)$ is given by
\[
\pi_u(A):= \oplus_{\omega \in \fS} \{\pi_{\omega}(A)\} \quad \text{and} \quad \norm{A}_u:= \norm{\pi_u(A)}=\sup_{\omega \in \fS}\norm{\pi_{\omega}(A)}
\]
denotes the enveloping $C^*$ seminorm of $\cR_0$. Define the resolvent algebra $\cR(\fX, \sigma)$ as the abstract $C^*$-algebra generated by $\pi_u(\cR_0)$, \ie we factor $\cR_0$ by $\ker \pi_u$ and complete with respect to the eveloping $C^*$ seminorm $\norm{\cdot}_{u}$.
\end{definition}

Useful properties of the Resolvent Algebra are:
\begin{theorem}\label{th:symautRA}
Let $(\fS, \sigma)$ be a given nondegenerate symplectic space, and define $\cR(\fX, \sigma)$ as above. Then for all $\lambda, \mu \in \mathbb{R} \backslash \{0\}$ we have:
\begin{itemize}
\item[(i)]$[R(\lambda,f),R(\mu,f)]=0$. Substitute $\mu=-\la$ to see that $R(\la,f)$ is normal.
\item[(ii)] $R(\lambda,f)R(\mu,g)^2R(\lambda,f)=R(\lambda,g)R(\mu,f)^2R(\lambda,g)$.
\item[(iii)] $\norm{R(\la,f)}_u=|\la|^{-1}$.
\item[(iv)] $R(\la,f)$ is analytic in $\la$. Explicitly, the series expansion (von Neumman series)
\[
R(\la,f)=\sum_{n=0}^{\infty} (\la_0-\la)^n i^n R(\la_0,f)^{n+1}, \qquad \la, \la_0 \neq 0
\]
converges absolutely in norm whenever $|\lambda_0 -\la|< |\la_0|$.
\item[(v)] Let $T\in \Sp(\fX,\sigma)$ be a symplectic transformation. Then $\alpha(R(\la,f)):=R(\la,Tf)$ extends to an automorphism $\alpha \in \Aut(\cR(\fX,\sigma))$.
\end{itemize}
\end{theorem}
In constrast to the Weyl algebra, the Resolvent algebra is not simple. We have the following ideal structure.

\begin{theorem}
\label{Ideals0}
Let $(\fX,\,\sigma)$ be a given nondegenerate symplectic space. 
Then for each $\lambda\in\R\backslash 0$
and $f\in \fX\backslash 0$ we have that the
 closed two--sided ideal generated by $\rlf$ in 
$\al R.(\fX,\,\sigma)$ is
\[
\big[\rlf\rsl\big]=\big[\rsl\rlf\big]=\big[\rsl\rlf\rsl\big]
\]
where ${[\, \cdot \,]}$ indicates the closed linear span 
of its argument. This ideal is proper.
Moreover the intersection of the ideals 
$\big[R(\lambda_i,f_i)\rsl\big]$,
 $i=1,\ldots,\,n$ for distinct $f_i\in \fX\backslash 0$ is the ideal
${\big[R(\lambda_1,f_1)\cdots R(\lambda_n,f_n)\rsl\big]}$.
\end{theorem}
{}From these ideals we can build other ideals, \eg for a set $S\subseteq \fX$
we can define the ideals ${\mathop{\bigcap}\limits_{f\in S}\big[\rsl R(\lambda,f)\big]}$
as well as ${\Big[\mathop{\bigcup}\limits_{f\in S}\big[\rsl R(\lambda,f)\big]\Big]}$. 
Ideals of a different structure will occur in the following sections.
Thus $\rsl$ has a very rich ideal structure.

\subsection{States, representations and regularity}
\label{StatesRep}

Any operator family $R_{\lambda}$, $\lambda \in \R \backslash 0$ 
on a Hilbert space which satisfies the resolvent equation
\eqref{eq:R4} is called by Hille a pseudo-resolvent
and for such a family we know (cf.\ Theorem~1
in~\cite[p 216]{Yos1980}) that:
\begin{itemize}
\item{} All $R_\lambda$ have a common range and a common null space.
\item{} A pseudo-resolvent is the resolvent for an operator
$B$ iff $\ker R_{\lambda} =\{0\}$ for some (hence for
all) $\lambda \in \R \backslash 0$, and in this case
$\dom B=\ran R_\lambda$ for all $\lambda \in \R \backslash 0$.
\end{itemize}
This leads us to an examination of $\ker\pi\big(\rlf\big)$ for representations
$\pi$.
\begin{theorem}
\label{Ideals1}
Let $(\fX,\,\sigma)$ be a given nondegenerate symplectic space, and define
$\al R.(\fX,\,\sigma)$ as above. Then for $\lambda\in\R\backslash 0$
and $f\in \fX\backslash 0$ we have:
\begin{itemize}
\item[(i)] If for a representation $\pi$ of $\rsl$ we have
$\ker\pi\big(\rlf\big)\not=\{0\} $,  
then $\ker\pi\big(\rlf\big)$ reduces $\pi(\rsl)$. Hence
there is a unique orthogonal decomposition
$\pi=\pi_1\oplus\pi_2$ such that $\pi_1(\rlf)=0$ and
$\pi_2(\rlf)$ is invertible.
\item[(ii)] Let $\pi$ be any nondegenerate representation of $\rsl$,
then  $$P_f:=\slim_{\lambda\to\infty}i\lambda\, \pi\big(\rlf\big) $$
exists, defines a central projection of $\pi(\rsl)''\!$, and it is the
range projection of $\pi\big(\rlf\big)$ as well as the projection
of the ideal ${\pi\left(\big[\rsl R(\lambda,f)\big]\right)}$.
\item[(iii)] If $\pi$ is a factorial representation of $\rsl$,
then $P_f=0$ or $\un$ and such $\pi$ are classified by the sets
${\big\{f\in \fX\backslash 0\,\mid\, P_f=\un\big\}}$.
\item[(iv)] There is a state $\omega\in\ot S.\big(\rsl\big)$
such that $\rlf\in\ker\omega$. Moreover, given a state $\omega$
with  $\rlf\in\ker\omega$, then $\rlf\in\ker\pi_\omega$. 
\end{itemize}
\end{theorem}
Given a $\pi\in\rep\big(\al R.(\fX,\,\sigma),\al H._\pi\big)$ with
$\ker\pi\big(R(1,f)\big)=\{0\},$ we define a field
operator by
\[
\j_\pi(f):=i\un-\pi\big(R(1,f)\big)^{-1}\;
\]
with domain $\dom \j_\pi(f)=\ran\pi\big(R(1,f)\big)$, and it has the following
properties:
\begin{theorem} 
\label{RegThm}
Let $\al R.(\fX,\,\sigma)$ be as above, and let $\pi\in\rep\big(\al R.(\fX,\,\sigma),\al H._\pi\big)$
satisfy $\ker\pi\big(R(1,f)\big)=\{0\}= \ker\pi\big(R(1,h)\big)$ for given
$f,\; h\in \fX.$ Then
\begin{itemize}
\item[(i)]  $\j_\pi(f)$ is selfadjoint, and $\pi(\rlf)\dom \j_\pi(h)
\subseteq\dom \j_\pi(h)$.
\item[(ii)]
$\lim\limits_{\lambda\to\infty}i\lambda\pi(\rlf)\psi=\psi$ for all
$\psi\in\al H._\pi$.
\item[(iii)]
$\lim\limits_{\mu \to 0}i\pi(R(1, \mu f))\psi=\psi$ for all
$\psi\in\al H._\pi$.
\item[(iv)]
The space $\al D.:={\pi\big(R(1,f)R(1,h)\big)\al H._\pi}$
is a joint dense domain 
for $\j_\pi(f)$ and $\j_\pi(h)$ 
and we have:
$[\j_\pi(f),\,\j_\pi(h)]=i\sigma(f,h) \un $ on $\al D..$
\item[(v)] $\ker\pi\big(R(1,\nu f+h)\big)=\{0\}$
for $\nu \in\R$. 
Then $\j_\pi(\nu f+h)$ is defined,  
$\al D.$ 
is a core for  $\j_\pi(\nu f+h)$ and
$\j_\pi(\nu f+h)=\nu \j_\pi(f)+\j_\pi(h)$
 on $\al D..$ Moreover
${\pi\big(R(1,\nu f+h)\big)}\in
{\big\{\pi\big(R(1,f)\big),\,\pi\big(R(1,h)\big)\big\}''.}$
\item[(vi)]
$\j_\pi(f)\pi(\rlf)=\pi(\rlf)\j_\pi(f)=i\lambda\pi(\rlf)-\un$ on
$\dom \j_\pi(f)$.
\item[(vii)]
$\big[\j_\pi(f),\pi(R(\lambda,h))\big]=i\sigma(f,h)\pi(R(\lambda,h)^2)$
on $\dom \j_\pi(f)$.
\end{itemize}
\end{theorem}
Thus we define:
\begin{definition}\label{df:regrep}
A  representation $\pi\in\rep\big(\al R.(\fX,\,\sigma),\al H._\pi\big)$
is {regular on} $S\subset \fX$ if
\[
\ker\pi\big(R(1,f)\big)=\{0\}\qquad\mbox{for all} \ f\in
S\;.
\]
A state $\omega$ of $\rsl$ is {regular on} $S\subset \fX$ if its GNS--representation
$\pi_\omega$ is regular on $S\subset \fX$.
A {regular representation (resp.\ state)} is a
representation (resp.\ state) which is regular on $\fX$. 
Given a Hilbert space $\al H.,$ we denote 
the set of (nondegenerate) regular representations $\pi:\rsl\to\al B.(\al H.)$
by ${{\rm Reg}\big(\rsl,\al H.\big)}.$
The set of regular states of $\rsl$ is denoted by
$\ot S._r\big(\rsl\big).$
\end{definition}
Obviously many regular representations are known, \eg the
Fock representation.
The class of all regular representations of $\rsl$ is not a set,
hence the necessity to fix $\al H..$
Thus for $\pi\in{{\rm Reg}\big(\rsl,\al H.\big)},$ all the field operators
$\j_\pi(f),$ $f\in
\fX$ are defined, and we have the resolvents 
$\pi(\rlf)=(i\lambda\un-\j_\pi(f))^{-1}$.

{}From Theorem~\ref{RegThm}, we can now establish a bijection between the
regular representations of $\rsl$ and the regular representations
of the Weyl algebra $\CCRX:$
\begin{corollary}
\label{RegBij} Let $\rsl$ be as above.
Given a regular representation $\pi\in{\rm Reg}\big(\rsl,\al H.\big),$ define a
regular representation $\wt\pi\in{\rm Reg}\big(\CCRX,\al H.\big)$
by $\wt\pi(\delta_f):=\exp(i\j_\pi(f))=W(f)$
(using Theorem~\ref{RegThm}(viii)). 
This correspondence establishes a bijection 
between ${\rm Reg}\big(\rsl,\al H.\big)$ and ${\rm Reg}\big(\CCRX,\al H.\big)$
which respects irreducibility and direct sums. Its inverse is given by
the Laplace transform,
\begin{equation}
\label{Laplace1}
\pi(R(\lambda,f)):= - i\int_0^{\, \infty}  e^{-\lambda t}\pi
( \delta\s -tf.)\,dt\;, \qquad \sigma := \mbox{sign} \, \lambda \,.
\end{equation}
By an application of this to the GNS--representations of regular states, 
we also obtain an affine bijection between $\ot S._r\big(\rsl\big)$
and the regular states $\ot S._r\big(\CCRX\big)$ of $\CCRX,$
and it restricts to a bijection between the pure regular states of
$\rsl$ and the pure regular states of $\CCRX.$
\end{corollary}
Note that whilst we have a bijection between the regular states of
$\rsl$ and those of $\CCRX$, there is no such map between the
nonregular states of the two algebras. In fact, fix a nonzero $f\in \fX$ and
consider the two commutative subalgebras 
$C^*\{\rlf\, , \un \, \mid\,\lambda\in\R\backslash 0\} \subset\rsl$ and
$C^*\{\delta_{tf}\,\mid\,t\in \R\}\subset\CCRX$, then these are isomorphic
respectively to the continuous functions on
the one point compactification of $\R$, and the continuous functions on
the Bohr compactification of $\R$.
Note that the point measures on
the compactifications without $\R$
 produce nonregular states (after extending to the full C*--algebras by
 Hahn--Banach) and there are
many more of these for the Bohr compactification than for the one 
point compactification of $\R$,
(cf. Theorem~5 in~\cite[p 949]{DS2}).
So the Weyl algebra has many more nonregular states than
the resolvent algebra.

Some further properties of regular representations and states are:
\begin{proposition} 
\label{RegAlg}
Let $\al R.(\fX,\,\sigma)$ be as above. 
\begin{itemize}
\item[(i)] If a representation $\pi$ of $\rsl$ is faithful and
factorial, it must be regular.
\item[(ii)] If a representation $\pi:\rsl\to\al B.(\al H.)$ is regular
then ${\|\pi(\rlf)\|}=\|\rlf\|=|\lambda|^{-1}$ for all $\lambda\in\R\backslash 0,$
$f\in \fX$.
 \item[(iii)] A state $\omega$ of $\rsl$ is regular iff
 $\omega(A)=\lim\limits_{\lambda\to\infty}i\lambda\,\omega\big(\rlf A\big)$
 for all $A\in\rsl$ and $f\in \fX$.
\end{itemize}
\end{proposition} 
Thus regular states restrict to regular states on subalgebras generated by the Resolvents.

We also find the the following decomposition of the test function space with regards to a given representation useful,
\begin{proposition}
\label{Xdecomp}
Let $\pi:\rsl\to\al B.(\al H.)$ be a nondegenerate representation. Then
\begin{itemize}
\item[(i)]
the set $\fX_R:=\set f\in \fX,\ker\pi\big(R(1,f)\big)=\{0\}.$ is a linear space.
Hence if $f\in \fX_S:=\fX\backslash \fX_R,$ then $f+g\in \fX_S$ for all $g\in \fX_R$.
\item[(ii)]
The set $\fX_T:=\set f\in \fX,\ker\pi\big(R(1,f)\big)=\{0\}\;\hbox{and}\;\pi\big(R(1,f)\big)^{-1}\in
{\al B.(\al H.)}.\subset \fX_R$ is a linear space. Moreover if $f\in \fX_T$ then
$\pi\big(R(1,g)\big)=0$ for all $g\in \fX$ with $\sigma(f,g)\not=0$.
Thus $\sigma(\fX_T,\fX_R)=0$.
\item[(iii)]
If $\pi$ is factorial, then $\pi\big(R(1,f)\big)=0$ for all $f\in \fX_S,$ and
$\pi\big(R(1,f)\big)\in\C\un\backslash 0$ for all $f\in \fX_T$. 
Moreover $\fX_T=\fX_R\cap \fX_R^\perp$.
\item[(iv)] Let $\fX$ be finite dimensional 
and let ${\big\{q_1,\ldots,\,q_n\big\}}$ be a basis for  $\fX_T$.
If $\pi$ is factorial, we can augment this basis of  $\fX_T$ by
${\big\{p_1,\ldots,\,p_n\big\}}\subset \fX_S$ into a symplectic
basis of $Q:={\rm Span}{\big\{q_1,p_1;\ldots;\,q_n,p_n\big\}},$
\ie $\sigma(p_i,q_j)=\delta_{ij}$, $0=\sigma(q_i,q_j)=\sigma(p_i,p_j)$.
Then we have the decomposition
\begin{equation}
\label{Qdecomp}
\fX=Q\oplus(Q^\perp\cap \fX_R)\oplus(Q^\perp\cap \fX_R^\perp)
\end{equation}
into nondegenerate spaces  such that
$Q^\perp\cap \fX_R\subset\{0\}\cup(\fX_R\backslash \fX_T)$ and
$Q^\perp\cap \fX_R^\perp\subset\{0\}\cup \fX_S$.
\end{itemize}
\end{proposition}
Clearly $\fX_R$ is the part of $\fX$ on which $\pi$ is regular,
$\fX_T$ is the part on which it is ``trivially regular'', 
$\fX_S$ is the part on which it is singular, and these have a particularly nice form when
$\pi$ is factorial. This proposition can be used to prove the following theorem.
\begin{theorem}
\label{UniqueR}
Let $\fX$ be a nondegenerate symplectic space of arbitrary dimension. Then
\begin{itemize}
\item[(i)]
The norms of $\al R.(\fX,\,\sigma)$ and $\al R.(S,\,\sigma)$ coincide on 
 ${\hbox{*-alg}\{\rlf\,\mid\,f\in S,\,\lambda\in\R\backslash 0\}}$ for each
finite dimensional nondegenerate subspace $S\subset \fX$.
Thus we obtain a containment ${\al R.(S,\,\sigma)}\subset\al R.(\fX,\,\sigma).$
\item[(ii)]$\al R.(\fX,\,\sigma)$ is the inductive limit of the net of all
${\al R.(S,\,\sigma)}$ where $S\subset \fX$ ranges over all finite dimensional
nondegenerate subspaces of $\fX.$
\end{itemize}
\end{theorem}
It follows therefore from Fell's theorem (cf. Theorem~1.2 in \cite{Fe1960}) 
that {\it any} state of $\rsl$ is in the w*-closure of the convex hull of
the vector states of $\pi_r$, hence of
the regular states. The following result 
is relevant for physics.
\begin{theorem}
\label{RegFaith}
Let $(\fX,\,\sigma)$ be any nondegenerate symplectic space, and 
$\rsl$ as above. Then every regular representation of $\rsl$ is faithful.
\end{theorem}
 The importance of this result lies in the fact that the regular representations
 are taken to be the physically relevant ones, and the images of $\rsl$
 in all regular representations are isomorphic. Thus, since we can obtain
 the quantum fields from $\rsl$ in these representations, we are justified 
 in taking $\rsl$ to be the observable algebra for bosonic fields.
 Usually one argues that for a C*--algebra $\al A.$ to be an  
 observable algebra of a physical system,
 it must be simple (cf.~\cite[p 852]{HaKa}). 
 The argument is that by Fell equivalence of
 the physical representations, the image of $\al A.$ in all
 physical representations must be isomorphic. However, if one restricts the class
 of physical representations (as we do here to the regular representations
 of $\rsl),$
 then the latter isomorphism does not imply that $\al A.$ must be simple.

This theorem also has structural consequences, \eg it implies that
$\rsl$ has faithful irreducible representations, hence that its centre must
be trivial. 
{}For many applications one needs regular representations where there is
a dense invariant joint domain for all the fields  $\j_\pi(f),$ and this 
leads us to a subclass of the regular representations as follows.
We will say that a state $\omega$ on the Weyl algebra $\CCRX$
is {strongly regular} if the functions
\[
\R^n\ni(\lambda_1,\ldots,\lambda_n)\mapsto\omega\big(\delta_{\lambda_1f_1}
\cdots\delta_{\lambda_nf_n}\big)
\]
are smooth for all $f_1,\ldots,\,f_n\in \fX$ and all 
$n\in\N$. Of special importance is that the GNS-representation 
of a strongly regular state has a common dense invariant domain
for all the generators $\j_{\pi_\omega}(f)$ of the one parameter groups
$\lambda\to\pi_\omega(\delta_{\lambda f})$
(this domain is obtained by applying the polynomial algebra of the Weyl operators
$\set\pi_\omega(\delta_f),f\in \fX.$ to the cyclic GNS-vector).
By the bijection of Corollary~\ref{RegBij},
we then obtain the set of strongly regular states on $\rsl,$ 
and we denote this by ${\ot S._{sr}\big(\rsl\big)}.$

\subsection{Further structure.}
\label{FurStruc}

Here we want to explore the algebraic structure of $\al R.(\fX,\,\sigma)$.
\begin{theorem}
\label{TensorAlg}
Let $(\fX,\,\sigma)$ be a given nondegenerate symplectic space, and let 
$\fX=S\oplus S^\perp$ for $S\subset \fX$ a nondegenerate subspace. Then
$$\al R.(\fX,\,\sigma)
 \supset{\rm C}^*\big(\al R.(S,\,\sigma)\cup\al R.(S^\perp,\,\sigma)\big)
\cong\al R.(S,\,\sigma)\otimes\al R.(S^\perp,\,\sigma)$$
where the tensor product uses the minimal (spatial) tensor norm.
The containment is proper in general.
\end{theorem}
Thus we cannot generate $\rsl$ from a basis alone, \ie if $\{q_1,p_1;\,q_2,p_2;\cdots\}$
is a symplectic basis of $\fX$,  then 
${\rm C}^*\big\{R(\lambda_i,q_i),\;R(\mu_i,p_i)\,\mid\,\lambda_i,\,\mu_i\in\R\backslash 0,\;
i=1,\,2,\ldots\big\}$ is in general a proper subalgebra of $\rsl,$ though in any
regular representation $\pi$ it is strong operator dense in
$\pi\big(\rsl\big)$ by Theorem~\ref{RegThm}(v).

Note that since ${C^*(\{\rlf\,\mid\,\lambda\in\R\backslash 0\})}\cong C_0(\R)$
(easily seen in any regular representation), and we have that 
$C_0(\R^{n+m})=C_0(\R^n)\otimes C_0(\R^m),$ it follows from Theorem~\ref{TensorAlg}
that any $C_0\hbox{--function}$ of a finite commuting set of variables 
is in $\rsl$. More concretely, we have the following result which will
be used later.
\begin{proposition}
\label{Czero}
Let $\{q_1,\ldots,\,q_k\}\subset \fX$ satisfy $\sigma(q_i,q_j)=0$
for all $i,\,j$. Then for each
$F \in C_0(\R^k)$ there is a (unique) $R_F \in\rsl$ such that
in any regular representation $\pi$
we have $\pi(R_F)={F \big(\j\s\pi.(q_1),\ldots,\j\s\pi.(q_k)\big)}$.
\end{proposition}
Thus the resolvent algebra contains in abstract form 
all $C_0\hbox{--functions}$ of commuting fields. 
Note that such a result neither holds for the 
Weyl algebra nor for the corresponding 
twisted group algebra (in the case of finite
dimensional $\fX$).

\begin{theorem}
\label{Nonsep}
Let $(\fX,\,\sigma)$ be a given nondegenerate symplectic space 
and let $f,\,h\in \fX\backslash 0$ such that  
$f\not\in\R h$. Then
\begin{itemize}
\item[(i)] $R(1,f)\not\in\big[\rsl R(1,h)\big]$, \ie the ideals separate the rays 
of $\fX$,
\item[(ii)] $\big\|R(1,f)-R(1,h)\big\|\geq 1$, and if $\sigma(f,h)=0$ we have equality.
\item[(iii)] $\rsl$ is nonseparable.
\end{itemize}
\end{theorem}

\subsection{Constraint theory.}\label{sbs:racons}
In this subsection we assume the structures associated to the $T$-procedure (quantum Dirac constraint procedure) as described in Appendix \eqref{app:Tp}.

{}For linear bosonic constraints, we start with a nondegenerate symplectic space
${(\fX,\sigma)}$ and specify a nonzero {\it constraint subspace} $C\subset \fX$.
Our task is to implement the heuristic constraint conditions 
\[
\j(f)\,\psi=0\quad f\in C
\]
to select the subspace spanned by the physical vectors $\psi$.
There are many examples where these occur, \eg in quantum electromagnetism,
cf. \cite{HendrikHu1985,Hendrik1988,Hendrik2000}.
Now in a representation $\pi$ of $\rsl$ for which $\ker\pi\big(\rlf\big)=\{0\}$
we have by Theorem~\ref{RegThm}(vi) that
$\pi(\rlf)\j_\pi(f)=i\lambda\pi(\rlf)-\un$ on $\dom \j_\pi(f)$.
Hence the appropriate form in which to impose the heuristic constraint condition
in the resolvent algebra is to select the set of physical (``Dirac'') states by
\begin{equation}
\label{Rconstraint}
\ot S._D:=\left\{\omega\in\ot S.(\rsl)\,\mid\,
\pi_\omega\big(i\lambda\rlf-\un\big)\,\Omega_\omega=0\quad 
 f\in C, \, \lambda \in \R \backslash 0  \right\} \, ,
\end{equation}
where $\pi_\omega$ and $\Omega_\omega$ denote the GNS--representation and 
GNS--cyclic vector of $\omega.$ Thus $\omega\in\ot S._D$ iff
$\al C.\subset
{\cal N}_\omega:={\big\{A\in\rsl\,\mid\,\omega(A^*A)=0\big\}}$,
where $\al C.:={\big\{i\lambda\rlf-\un\,\mid\, f\in
  C, \, \lambda \in \R \backslash 0 \big\}}$. Note that 
${\cal C}^* = {\cal C}$. 

\begin{proposition}
\label{RDirac}
Given the data above, we have:
\begin{itemize}
\item[(i)] $\ot S._D=\left\{\omega\in\ot S.(\rsl)\,\mid\,
\omega\big(R(1,f)\big)=-i, \ f\in C\right\}$.
\item[(ii)] If $\omega\in\ot S._D$, then it is not regular. 
In particular, if $\sigma(g,C)\not=0$ for some $g\in \fX,$ 
then $\pi_\omega(R(\lambda,g))=0$ for all $\lambda\in\R\backslash 0$.
\item[(iii)]  $\ot S._D\not=\emptyset$ iff $\sigma(C,C)=0$.
\end{itemize}
\end{proposition}
Henceforth we will assume that $\sigma(C,\,C)=0$ and hence $\ot S._D\not=
\emptyset$.

\begin{proposition}\label{pr:conRA} With notation 
\[
\cD:=\cN\cap \cN^*,\qquad \cO:=\{A \in \cR(\fX,\sigma)\,|\, [A,\cD]\subset \cD\}, \qquad \cN:=[\cR(\fX,\sigma)\cC]
\]
from the $T$-procedure we have:
\chop
$\WD=\rsl$ with the proper ideal  $\al D.=[\rsl\al C.]=[\al C.\rsl]$,
and $\Ob =\al C.'\big/(\al C.'\bigcap\al D.)$ where $\cP:=\cO/\cD$.
\end{proposition}
So Dirac constraining of linear bosonic constraints is considerably simpler
in the resolvent algebra $\rsl$ than in the CCR--algebra $\CCRX$~cf.\cite{HendrikHu1985}.

\section{Constraints I: Symplectic Form $\sigma_1$} \label{sbs:cnst1}
We want to use the Resolvent Algebra to model a Fock-Krein bosonic field algebras of the type discussed in Section \ref{sec:FKCCR}. Hence for the remainder of this chapter we will assume that $\fD$ and $\fL$ are test function spaces with all the structures as in Subsection \ref{sbs:abstf}. Now we have a technical difficulty that $A(g)$ as defined in Section \ref{sec:FKCCR} are Krein symmetric but not Hilbert selfadjoint in general, and hence their spectrum need not be a subset of the real line, hence $(i \lambda -A(g))^{-1}$ is need not be defined for all $\lambda \in \mathbb{R}\backslash\{0\}$.
This means that although the Resolvent Algebra exists using the QEM test function space, there is no reason it maps to `resolvents' $(i \lambda -A(g))^{-1}$ and this poses the problem of what Resolvent algebra to use to model the $A(g)$'s. 

The Resolvent Algebra exists for any non-degenerate symplectic space. So given the CCR's in Section \ref{sec:FKCCR} are: 
\[
[A(f),A(g)]=i\smp{f}{g}{1}, \qquad f,g \in \fX
\]
we will first try using the Resolvent Algebra $\cR(\fX,\sigma_1)$ to model the Fock-Krein fields given in Section \ref{sec:FKCCR}. 

Before we use $\cR(\fX,\sigma_1)$ to develop a $C^*$-algebraic model for the structures in in Chapter \ref{ch:BRSTQEM}, we apply the $T$-procedure to $\cR(\fX,\sigma_1)$ as in Subsection \ref{sbs:racons}. This will establish a reference point for comparison with subsequent results.

\subsection{$T$-procedure I}\label{sbs:DcI}
In this subsection we assume the structures associated to the $T$-procedure (quantum Dirac constraint procedure) as described in Appendix \eqref{app:Tp}.

Motivated by the Lorentz condition in Subsection \ref{sbs:lorentz}, we want implement the heuristic constraint $\{A(f)\psi=0\,|\, f \in \fXL\}$. As described in Subsection \ref{sbs:racons}, this corresponds in the Resolvent Algebra to the T-procedure  with constraint set,
\[
\cC:=\{ i \lambda R(\lambda,f)-\one \,|\, f \in \fXL,\, \lambda \in \R \backslash 0 \}.
\]
As $\smp{\fXL}{\fXL}{1}=0$ we have that Dirac states exist by Proposition \eqref{RDirac} and that the physical (constrained algebra) is:
\[
\cP= \cC'/(\cC'\cap \cD)
\]
by Proposition \eqref{pr:conRA}. We next construct $\cP$ more explicitly for comparison with other results. By equation \eqref{Rconstraint} we have that for $\omega \in \fS_D$,
\[
\pi_{\omega}(i \lambda R(\lambda,f))\Omega_{\omega}=\Omega_{\omega}.
\]
However we can get the stronger statement:
\begin{corollary}\label{cr:trivcons}
Let $\omega \in \fS_D$, then for all $f \in \fXL$, $\lambda \in \R \backslash 0$:
\[
\pi_{\omega}(i \lambda R(\lambda,f))=\one.
\]
\end{corollary}
\begin{proof}
Since $\cD$ is an ideal of $\cR(\fX,\sigma_1)$ (Proposition \eqref{pr:conRA}) and $\cD \subset \ker \omega$ for all $\omega \in \fS_D$, it follows by \cite{Dix77} 2.4.10 that $\cD \subset \ker \pi_{\omega}$ hence that $i \lambda R(\lambda,f)-\one \in \cC \subset \cD \subset \ker \pi_{\omega}$ for all $f \in \fXL$.
\end{proof}

Now  $\fX_t$ and $\fXL\oplus \fXJ$ are non-degenerate and $\fX_t \perp (\fXL\oplus \fXJ)$ $\sigma_1$-symplectically. Hence we have by Theorem \eqref{TensorAlg} that
\[
\cR_{T}:=C^*(\{ R(\lambda, f),\,| f \in \fX_t \cup \fXL \oplus \fXJ, \, \lambda \in \R \backslash 0 \} )\cong \cR(\fX_t, \sigma_1)\otimes \cR(\fXL \oplus \fXJ, \sigma_1).
 \]
Let $\vp$ be this $*$-isomorphism, and define $\cR_{ph}:= \vp^{-1}(\cR(\fX_t, \sigma_1) \otimes \one)$ and $\cR_{u}:= \vp^{-1}(\one \otimes \cR(\fXL \oplus \fXJ, \sigma_1))$.

To characterize the physical observable algebra $\cP$, we prove the following lemmma. Recall the notation $\omega_{\cB}=\omega|_{\cB}$ where $\omega \in \fS_{\cA}$ and $\cB$ is a $C^*$-subalgebra of $\cA$ (cf. Definition \eqref{df:csksrepterm}).

\begin{lemma}\label{lm:resD2} We have:
\begin{itemize} 
\item[(i)]For every ${\omega_1} \in \fS(\cR(\fX_t, \sigma_1))$ there exists an $\omega \in \fS_D(\cR(\fX, \sigma_1))$ such that $\omega_{\cR_{ph}}= {\omega_1 \circ \vp}$. Furthermore we have that,
\[
\cR_{ph} \cap \cD= \{0 \}.
\]
\item[(ii)]Define the representation of the observables as $\pi_D:\cR(\fX,\sigma_1) \to B(\cH_D)$ by: 
\[
\cH_{D}:=\bigoplus\{ \cH_{\omega}\,|\, \omega \in \fS_{D}\}, \qquad \pi_{D}:\bigoplus \{ \pi_{\omega}\,|\, \omega \in \fS_{D}\}.  
\]
Then $\pi_D(\cR_{ph})=\overline{\pi_D(\cR_{0}(\fX_t, \sigma_1))}^{B(\cH_D)} \cong \cR(\fX_t,\sigma_1)$, the closure in the first equality with respect to the uniform norm of $B(\cH_D)$.
\item[(iii)] Let $\omega \in \fS^{P}_D(\cR(\fX, \sigma_1))$, \ie $\omega$ is a pure Dirac state. Then for all $f\in \fX_t$, $g\in \fXL$ and $\lambda \in \R\backslash \{0\}$ we have:
\[
\pi_{\omega}(R(\lambda, f+g))= \pi_{\omega}(R(\lambda, f))
\]
\end{itemize}
\end{lemma}
\begin{proof}
(i): We have that $\cC:=\{ i \lambda R(\lambda,f)-\one \,|\, f \in \fXL,\, \lambda \in \R \backslash 0 \}\subset \cR (\fXL \oplus \fXJ, \sigma_1)$ and so we can use $\cC$ as the constraint set for the $T$-procedure in $\cR (\fXL \oplus \fXJ, \sigma_1)$. As $\smp{\fXL}{\fXL}{1}=0$ Proposition \eqref{RDirac} implies that Dirac states $\fS_D(\cR (\fXL \oplus \fXJ, \sigma_1))$ exist. Take $\omega_1 \in \fS(\cR(\fX_t, \sigma_1))$ and $\omega_2 \in \fS_D(\cR (\fXL \oplus \fXJ, \sigma_1))$ and let $\tilde{\omega}:= (\omega_1 \otimes \omega_2)\circ \vp \in \fS(\cR_T)$. We can extend $\tilde{\omega}$ to $\omega \in \fS(\cR(\fX, \sigma_1))$ by the Hahn-Banach theorem, and it is easy to check that $\omega|_{\cR_{ph}}=\omega_1 \circ \vp$ and $\cC \subset \cN_{\omega}$.

For the last statement in (i) let $A\in \cR_{ph}$ and ${\omega_1} \in \fS(\cR(\fX_t, \sigma_1))$ such that ${(\omega_1\circ \vp)(A^*A)}\neq 0$. By the above there exists $\omega \in \fS_D$ such that $\omega(A^*A)= {(\omega_1\circ \vp)(A^*A)} \neq 0$, hence $\cR_{ph} \cap \cN_\omega = \{ 0 \}$. As $\cR_{ph}$ is a $*$-algebra, $\cR_{ph} \cap \cD = \{0 \}$.

\smallskip
\noindent (ii): By the definition of $\pi_D$ we see that $\ker \pi_D= \cap \{ \ker \pi_{\omega}\,|\, \omega \in \fS_D(\cR(\fX, \sigma_1))\}=\cD$ where the last  equality follows by Theorem \eqref{th:DsD}. Hence by (i) we see that $(\ker \pi_{D}\cap \cR_{ph})=(\cD \cap \cR_{ph})=\{0\}$, and hence $\pi_{D}|_{\cR_{ph}} $ is isometric.
 As $\cR_{ph} \cong \cR(\fX_t, \sigma_1)$ the result follows.

\pfit(iii): Let $\omega \in \fS^{P}_D(\cR(\fX,\sigma_1))$ and let $f \in \fX_t$. By Corollary \eqref{cr:trivcons} we have $\pi_{\omega}(R(\lambda, g))=-(i/\lambda)\one$ for all $g \in \fXL$ and $\lambda \in \mathbb{R} \backslash \{0\}$. As $\,\omega\,$ is pure, $\pi_{\omega}(R(\lambda, f))=0$ or $\pi_{\omega}(R(\lambda, f))$ is invertible by Theorem \eqref{Ideals1} (i). 

Consider the case that $\pi_{\omega}(R(\lambda, f))$ is invertible. Note that $\pi_{\omega}(R(\lambda, g))=-(i/ \lambda)\one$ is invertible, hence  $\phi_{\pi_{\omega}}(g)$ exists and $\phi_{\pi_{\omega}}(g)=0$. Therefore,
\begin{align*}
\pi_{\omega}(R(\lambda, f))= & \;(i \lambda \one +\phi_{\pi_{\omega}}(f))^{-1},\\
= & \;(i \lambda \one +\phi_{\pi_{\omega}}(f)+\phi_{\pi_{\omega}}(g))^{-1},\\
= & \;(i \lambda \one +\phi_{\pi_{\omega}}(f+g))^{-1},\\
= & \;\pi_{\omega}(R(\lambda, f+g))
\end{align*}
where we used Theorem \eqref{RegThm} (v) in the third equality. 

Next consider the case that $\pi_{\omega}(R(\lambda, f))=0$. Note that $\pi_{\omega}(R(1, f))$ is invertible if and only if $\pi_{\omega}(R(\lambda, f))$ is, since in each case $\phi_{\pi_{\omega}}(f)$ exists, hence since $\omega$ is pure $\pi_{\omega}(R(1, f))=0$ if and only if $\pi_{\omega}(R(\lambda, f))=0$. So by assumption we have $\pi_{\omega}(R(1, f))=0$. Proposition \eqref{Xdecomp} (i) gives $(f+g) \in \fX_S$ where $\fX_S$ defined as in Proposition \eqref{Xdecomp} and so $\pi_{\omega}(R(1, f+g))=0$. Therefore
\[
\pi_{\omega}(R(\lambda, f))=0 \Rightarrow \pi_{\omega}(R(\lambda, f+g))=0.
\]
A similar argument shows that $\pi_{\omega}(R(\lambda, f+g))=0$ implies $\pi_{\omega}(R(\lambda, f))=0$ for all $\lambda \in \mathbb{R} \backslash \{0\}$, and we are done.
\end{proof}
From these lemmas we can characterize the Dirac physical observable algebra:
\begin{proposition}\label{pr:RA1phalg}
We have
\[
\cP:=\cO/\cD \cong \cR(\fX_t, \sigma_1) \cong \cR((\fXL \oplus \fX_t)/\fXL, \sigma_1)
\]
where `$\cong$' denotes (isometric) $*$-isomorphism.
\end{proposition}
\begin{proof}
First, by a slight abuse of notation we will write $\cR(\fX_t, \sigma_1)$  for $\cR_{ph}$ as by Theorem \eqref{UniqueR} (i), we have that the norms of $\cR_{ph}$ and $\cR(\fX_t, \sigma_1)$ coincide. Let $\omega \in \fS_D(\cR(\fX, \sigma_1))$ and let $f \in \fX_t$, $g \in \fXL$, $h \in \fXJ$. Consider
\[
R(\lambda, f+g+h) \in \cR(\fX, \sigma_1)
\]
where $\lambda \in \R \backslash 0$. By assumption we have that $\fX$ is $\sigma_1$-nondegenerate and that $\smp{\fXL}{\fX_t\oplus\fXL}{1}=0$, so if $h\neq0$ then there exists a $k \in \fXL$ such that $\sigma_1(k,h) \neq 0$. Hence by Proposition \eqref{RDirac} we have for all $\lambda \in \R \backslash 0$ that $\pi_{\omega}(R(\lambda, f +g +h ))=0$ and hence $R(\lambda, f +g +h ) \in \cD$. Now assume $h=0$, \ie we consider
\[
R(\lambda, f + g)  \in \cR(\fX, \sigma_1).
\]
By lemma \eqref{lm:resD2} (iii) and that $\pi_{\omega}$ is a representation we have that $R(\lambda,f)-R(\lambda, f + g)\in \ker \pi_{\omega}$ for all $\omega \in \fS_D^{P}(\cR(\fX,\sigma_1))$, and hence that $R(\lambda,f)-R(\lambda, f + g)\in  \cD$ by Corollary 3.13.8 \cite{Ped79}, where  we recall that $\cD \triangleleft \cO=\cR(\fX, \sigma_1)$ where the last equality is by Proposition \eqref{pr:conRA}. 
Let $\tau:\cR(\fX, \sigma_1) \to \cP=\cO/\cD$ be the factor map. By the arguments in the preceding paragraph, and as $\tau$ is a continuous $*$-homomorphism, we get that:
\[
\cP\cong \tau(\cR(\fX,\sigma_1))\subset \overline{\tau(\cR_0(\fX,\sigma_1))} \subset \overline{\tau(\cR_0(\fX_t,\sigma_1))}  
\]
By lemma \eqref{lm:resD2} (i) we have that $\ker \tau \cap \cR_{ph}= \cD \cap \cR_{ph}=\{0\}$ where $\cR_{ph}=\cR(\fX_t,\sigma_1)$, hence \\
$\tau:\cR_0(\fX_t,\sigma_1) \to \cP$ is injective hence isometric (cf. \cite[Theorem 3.1.5, p80]{Mur90}). Therefore 
\[
\overline{\tau(\cR_0(\fX_t,\sigma_1))}\cong \overline{\cR_0(\fX_t,\sigma_1)}=\cR(\fX_t,\sigma_1)
\]
where `$\cong$' denotes $*$-isomorphism (automatically isometric by \cite[Theorem 3.1.5, p80]{Mur90}). Putting these together we get that
\[
\cP \cong  \cR(\fX_t,\sigma_2)
\]
The last isomorphism statement follows as $\fXL$ is the $\sigma_1$ degenerate part of $\fX_t \oplus \fXL$.
\end{proof}

\subsection{QEM and Covariance I}\label{sbs:QEMcovI}
We take $\fX$ as in Subsection \ref{sbs:testfunc} and hence we get that
\[
\cP=\cO/\cD \cong \cR(\fX_t, \sigma_1) \cong \cR((\fXL \oplus \fX_t)/\fXL, \sigma_1)
\]
Recall that for $\fX$ as in Subsection \ref{sbs:testfunc}, the Poincar{\'e} transformations defined by 
\begin{equation*}
(V_gf)(p):=e^{ipa}\Lambda f(\Lambda^{-1}p) \quad \forall  f\in L^2(C_{+},\mathbb{C}^4), \quad g=(\Lambda,a) \in \cP^{\uparrow}_{+},
\end{equation*}
(cf. equation \eqref{eq:Ptran}). As the measure $\lambda$ is Lorentz invariant, it is straightforward to see that $V_g$ is $\iip{\cdot}{\cdot}$-unitary on $L^2(C_{+},\mathbb{C}^4,\lambda)$ and so $\sigma_1$-symplectic on $\fX$. Hence we have by Theorem \eqref{th:symautRA} (v) that,
\[
\alpha_{g}(R(\lambda,f))=R(\lambda,V_gf)
\]
extends to an automorphism on $\cR(\fX,\sigma_1)$ for all $g=(\Lambda,a) \in \cP^{\uparrow}_{+}$. Also as $V_g$ preserves $\fXL$ (lemma \eqref{lm:covX1}), we have that $\alpha_{g}(\cC)=\cC$ and hence $\alpha_{g}$ factors to an automorphism on $\cP$. Hence the Poincar{\'e} transformations are defined naturally on the constrained algebra and we use representations that satisfy the spectral condition for physical representations.

\subsection{$C^*$-BRST I}\label{sbsc:cbrstv1}
We want to construct a $C^*$-algebraic version of Fock-Krein BRST as in Section \ref{sec:FKBRST}. The definition of $\drb$ in \eqref{df:brstsd}  uses the complexified test function space $\fD=\fX+i\fX$ (cf. Subsection \ref{sbs:dspv}), hence we use will use the Resolvent algebra $\cR(\fD, \sigma_1)$ to define the $C^*$-algebraic $\drb$ and restrict to $\cR(\fX, \sigma_1)$ once done.

A complication with this approach is, as already discussed at the beginning of Subsection \ref{sbs:cnst1}, that $A(f)$ as defined in Subsection \ref{sec:FKCCR} is $\dag$-symmetric but not necessarily $*$-symmetric. Hence the operators $(i\lambda\one-A(f))^{-1}$ need not necessarily exist nor does  $\alg{ (i\lambda\one-A(f))^{-1}\,|\, f \in \fX}$ into which $\cR(\fD, \sigma_1)$ naturally maps. We can however think of $R(\lambda, f)$ `$=$' $(i\lambda\one-A(f))^{-1}$ as a heuristic formula and use it as guideline to construct a rigorous BRST superderivation and investigate the results, comparing the final constrained BRST-physical system with that in Subsection \ref{sbs:cnst1}.

To begin the construction to define the abstract BRST-Field Algebra: 
\begin{itemize}
\item Let $\fD$ and $\fL$ have all the structures as in Subsection \ref{sbs:abstf}. Let $ \cR(\fD, \sigma_1)$  and let $\cA_g(\fLJ)$  be the ghost algebra (Section \ref{sec:ghost}) with all the definitions as there, in particular
\[
C(f):=\fst( c( f)+ c^*(J f)), \qquad f \in \fLL\oplus \fLJ.
\]
and $\{C(f),C(g)\}=\re\,\iip{f}{g}\one$ for all $f,g \in \fLL \oplus \fLJ$
\item We define the BRST-Field Algebra as
\[
\cA:=\cR(\fD, \sigma_1)\otimes \cA_g(\fLJ).
\]
The tensor norm on $\cA$ is unique as the CAR algebra is nuclear. We define a grading on $\cA$ by extending the ghost grading on $\cA_g$ (cf. Definition \eqref{df:Z2grad}), \ie we define a grading automorphism $\gamma$ on $\cA$ as $\gamma$ equal to the identity on $\cR(\fD, \sigma_1)$, and equal to the $\Z_2$-grading automorphism on $\cA_g(\fLJ)$ (cf. Definition \eqref{df:Z2grad}). 
\item As $J$ is symplectic, we have that $\beta'(R(\lambda, f)):= R( \lambda, Jf)$ defines a unique automorphism on $\cR(\fD, \sigma_1)$ (cf. Theorem \eqref{th:symautRA} (v)), and we let  $\beta= \beta' \otimes \one \in \Aut(\cA)$. Also, as $J$ defines a unitary on $\fLL\oplus \fLJ$, we get that $\alpha'(C(f)):=C(Jf)$ defines a unique automorphism on $\cA_g$, and we let $\alpha:=\alpha' \otimes \one \in \Aut(\cA)$. 
\end{itemize}

\subsection{Superderivation I}\label{sbs:sdI}
We next define the BRST superderivation. As it will be unbounded we first specify its domain. 
\begin{definition}\label{df:sddom1}
Let
\[
D_1(\drb):=\alg{\one, \, R(\lambda, f), \, \zeta_1(h), \, \zeta_1(h)^* \,| \, f \in \fD, \,h \in \fDL \, \lambda\in \mathbb{R} \backslash 0},
\]
where $\fD$ and $\fDL$ are as in Subsection \ref{sbs:abstf}, and 
\[
D_2(\drb):=\alg{D_1(\drb),\, \zeta_2(f) \,| \, f \in \fD,\, \lambda\in \mathbb{R} \backslash 0},
\]
where $\zeta_1(h):=R(1, \mH h)\otimes C( h)$ and $\zeta_2(f):=R(1,  f)\otimes C(P_{\JJ} K f)$.
\end{definition}
\begin{rem}
\begin{itemize}
\item[(i)] Note that $D_1(\drb)\subset D_2(\drb)$ but that $D_1(\drb)$ is a $*$-algebra whereas $D_2(\drb)$ is not, and that neither is norm dense in $\cA$.
\item[(ii)] Note that $R(\lambda,  f)\otimes C(P_{\JJ} K f)=R(1,  f/\lambda)\otimes C(P_{\JJ} K f/\lambda)=\zeta_2(f/\lambda)$ for $\lambda\in \mathbb{R} \backslash 0$. This implies that $\alpha(\zeta_i^*(f))=\zeta_i(-f)$ for $i=1,2$.
\item[(iii)] Recall that $\gamma$ is the grading automorphism on $\cA$, hence $\gamma(R(\lambda,f)\otimes \one)=R(\lambda,f)\otimes \one$ and $\gamma(\one \otimes C(g))=-\one \otimes C(g)$ for all $f \in \fD$ and $g \in \fDL \oplus_{\mH} \fDJ$ (cf. Remark \eqref{rm:Zwgr}). Hence $\gamma(D_1(\drb))=D_1(\drb)$ and $\gamma(D_2(\drb))=D_2(\drb)$.
\end{itemize}
\end{rem}
Now recall the Definitions \eqref{df:brstsd},
\begin{align*}
\delta(A(g)\otimes \one)&=- i \one \otimes C( KP_{\JJ}  g), \qquad g \in \fD, \\
\delta(\one \otimes C(g))&= A(\mH P_{\LL} g)\otimes \one \qquad g \in \fD,
\end{align*}
The second identity is problematic and we want to encode it in bounded form. To do so we construct a mollified version of $\drb$. To motivate this we make the heuristic identification of $R(\lambda, f)$`$=$'$(i\lambda\one-A(f))^{-1}$, and use the mollifying property of resolvents in the following heuristic calculation: For $h \in \fDL$
\begin{align*}
\drb(\zeta_1(h))=&\; \drb(R(1, \mH h)\otimes C( h)),\\
=&\; \drhb(i\lambda\one-A(f))^{-1})C(h)+R(1,  \mH h))\drb(C(h)),\\
=&\; -R(1,  \mH h)^2C(iP_{\JJ}P_{\LL} \mH h)C(P_{\LL} f)+R(1,  \mH h))A(\mH h),\\
=&\; iR(1, \mH h )- \one, 
\end{align*}
where we used the mollifying properties of the resolvent in the last equality. The LHS and RHS are both well defined elements in $D_2(\drb)$. Similarly we can use heuristic calculations on all the generators of $D_2(\drb)$ to rigorously define $\drb$ as map. We can then recover the algebraic structure of superdrivation as given in Definitions \eqref{df:brstsd} in regular representations. Note however this will \emph{not} correspond exactly as resolvents in the Resolvent Algebra correspond to Hilbert essentially selfadjoint fields, while the $A(f)$ in Chapter \ref{ch:BRSTQEM} are Krein symmetric. 
\begin{theorem}\label{th:calgdrb}
Define a map on the elementary tensors in $D_2(\drb)$ by:
\begin{align*}
\drb(R(\lambda, f))&= -i R(\lambda, f)^2 C( KP_{\JJ} f)\\
&=-i R(\lambda, f)\zeta_2(f/\lambda) \quad\in D_2(\drb),  \\
\drb(\zeta_1(h))&= iR(1, P_{\LL} \mH h)- \one \quad \in D_1(\drb), \\
\drb(\zeta_1(h)^*)&=0, \\
\drb(\zeta_2(f))&=0,
\end{align*}
for $f\in \fD$, $h\in \fDL$, $\lambda\in \mathbb{R} \backslash 0$. This extends to a superderivation on $\drb:D_2(\drb)\to D_2(\drb)$ such that:
\begin{itemize}
\item[(i)] $\gamma \circ \drb \circ \gamma= -\drb$ on $D_2(\drb)$.
\item[(ii)] $\drb(A)^*=- \alpha \circ \drb \circ \alpha \circ \gamma(A^*)$ on $D_2(\drb)$.
\item[(iii)] $\drb^2=0$ on $D_2(\drb)$.
\end{itemize}
\end{theorem}
\begin{proof}
First we verify that $\drb$ is a superderivation on $D_2(\drb)$. We follow \cite{HendrikBuch2006} p708. Let $\pi_S$ be a strongly regular (hence faithful) representation of $\cR(\sigma_1, \cX)$. As $\pi_S$ is regular we have that $\phi_{\pi_S}(f)\in \op(\cH_s)$ exists for all $f\in \fD$ and have the properties given by Theorem \eqref{RegThm}. Furthermore, as $\pi_S$ is strongly regular, there exists a dense invariant domain $\cD_{\infty}$ for all $\phi_{\pi_S}(f)$, $f\in \fD$. Thus by Theorem \eqref{RegThm} (i), the resolvents $\pi_{S}(R(\lambda,f))$ map $\cD_{\infty}$ back into $D(\phi_S(f))$ for all $f \in \fD$. Thus we can define a second dense invariant domaint $\cD_S$ by applying all polynomials in $\phi_{\pi_S}(f)$ and $\phi_{\pi_S}(R(\lambda,f))$ to $\cD_{\infty}$.

We define the non-normed $*$-algebra
\[
\cE_0:= \salg{\phi_{\pi_S}(f), \, \pi_S(R(\lambda, f)) \, | \, f\in \fD, \, \lambda\in \mathbb{R} \backslash 0}
\]
which acts on the common dense invariant domain $\cD_S $, and we have the CCR's
\begin{equation}\label{eq:CCRRA1}
[\phi_{\pi_S}(f),\phi_{\pi_S}(g)]\psi=i\smp{f}{g}{1}\psi, \qquad \psi \in \cD_S
\end{equation}
Now let $\pi_0$ be any faithful representation of $\cA_g$, and so $\pi_S \otimes \pi_0$ is a faithful representation of $\cA$. Furthermore $\cD:= \cD_S \otimes \cH_g$ is a common dense invariant domain for $\pi_S(R(\lambda, f)) \otimes \one$, $\phi_{\pi_S}(f)\otimes \one$, and $\one \otimes \pi_0(C(f))$. For convenience of notation we will drop the $\otimes$ and $\pi_S$, $\pi_0$ for the remainder of this proof and define
\[
\cE:=\salg{\phi(f),\, R(\lambda,f), C(g) \,|\,f\in \fD, \, g \in \fDL \oplus \fDJ \, \lambda\in \mathbb{R} \backslash 0 }.
\]
We now define the map, $\drhb$ from the generators of $\cE$ to $\cE$ by:
\begin{align*}
\drhb(\phi(f))&= -i  C( KP_{\JJ}  f),  \qquad f \in \fD\\
\drhb(R(\lambda, f))&= -i R(\lambda, f)^2 C( KP_{\JJ}  f), \qquad f \in \fD\\
\drhb(C(f))&= \phi(\mH P_{\LL} f), \qquad f \in \fLL \oplus \fLJ
\end{align*}
and show that this extends to a well defined superderivation on $\cE$. To do this  we show that $\drhb$ is linear and satisfies the graded Leibniz rule on any finite polynomial in the operators $\phi(f)\, R(\lambda,f), C(f)$ where $f\in \fD, \lambda\in \mathbb{R} \backslash 0$.

Let $\cXs$ be a finite-dimensional subspace of $\fDL$ and let $\fD_s= \cXs \oplus J \cXs$. Let 
\[
\cE(\fD_s):=\alg{\phi(f)\, R(\lambda,f), C(g) \,|\, f \in \fD_t \oplus \fD_2,\,g \in \fD_s, \, \lambda\in \mathbb{R} \backslash 0} \subset \cE,
\] 
and let $(f_j)_{j\in \Lambda}$ be a finite orthonormal basis for $\cXs$ and define,
\[
\Qs:= \sum_{j\in \Lambda} \left( \phi( f_j)C( Jf_j) + \phi(K f_j)  C(K Jf_j) \right).
\]
Let 
\[
\drhb_s(\cdot)\psi:=\sb{\Qs}{ \cdot}\psi, \qquad \psi \in \cD
\] 
This is the same formula as the one for $\Qs$ in lemma \eqref{lm:2nQs}, but we have substituted $\phi(f)$ for $A(f)$. Since $A(f)$ and $\phi(f)$ satisfy the same CCR's, we obtain by the same calculations in the proof in lemma \eqref{lm:2nQs} that $\Qs^2=0$.
 
Given $g\in \fD_s, \lambda\in \mathbb{R} \backslash 0$ then $\sum_{j\in \Lambda} \ip{f_j}{g}f_j=P_{\LL}g$ and $\sum_{j\in \Lambda} \ip{Jf_j}{g}Jf_j=P_{\JJ}g$, and we calculate:
\begin{align*}
\drhb_s(\phi(g))=[\Qs,\phi(g)]&= \sum_{j\in \Lambda} (i\smp{f_j}{g}{1}C(Jf_j)+i\smp{Kf_j}{g}{1}C(KJf_j)) \\
&=-i\sum_{j\in \Lambda} C(K(K\smp{f_j}{g}{1}Jf_j+\smp{f_j}{Kg}{1}Jf_j)), \\
&=-i\sum_{j\in \Lambda} C(K\ip{Jf_j}{g}Jf_j), \\
&=-iC(KP_{\JJ}g),
\intertext{where we have used the CCR's (equation \eqref{eq:CCRRA1}) for $\phi(f)$, the definition of $\fD$ (cf Subsection \ref{sbs:abstf}) and $\ip{\cdot}{\cdot}$ (equation \eqref{eq:ipcst}) for the first identity. To get} 
\drhb_s(R(\lambda, g))&= -i R(\lambda, g)^2 C( KP_{\JJ}  g), \qquad g \in \fDJ
\intertext{we use $\big[\j(f),R(\lambda,h)\big]=i\sigma_1(f,h)R(\lambda,h)^2$
on $\cD_S$ (Theorem \eqref{RegThm} (vii)). Noting that the CAR's for the $C(f)$ use the inner product $\ip{\cdot}{\cdot}_{\mH}=\ip{\cdot}{\mH \cdot}$ on $\fLL\oplus \fLJ$, we get (similar to calculation for $\drhb_s(\phi(g))$)}
\drhb_s(C(g))&= \{\Qs, C(g) \}= \phi(\mH P_{\LL} g), \qquad g \in \fDL
\end{align*}
So we see that $\drhb$ agrees with $\drhb_s$ on the generating elements of $\cE(\fD_s)$ and as $\cXs$ was arbitrary we see that $\drhb_s$ extends to all of $\cE$ as a graded derivation and coincides with $\drhb_s$ on each $\cE(\cD_s)$.

For $h \in \fDL$ we calculate the identity:
\begin{align*}
\drb(\zeta_1(h))=&\; \drhb(R(1, \mH h))C(h)+R(1,  \mH h))\drhb(C(h)),\\
=&\; -iR(1,  \mH h)^2C(KP_{\JJ} \mH h)C(P_{\LL} h)+R(1,  \mH h))\phi(\mH h),\\
=&\; R(1,  \mH h))\phi(\mH h),\\
=&\; iR(1, \mH h )- \one, 
\end{align*}
where we used use the mollifying properties of the Resolvent Algebra in the last line.

The identities $drb(\zeta_1(h)^*)=\drb(\zeta_2(f))=0$ for $h \in \fDL$, $f \in \fD$ are obvious. Now $\drhb$ preserves $\pi_S \otimes \pi_0(D_2(\drb))\subset \cE$, hence defines a superderivation $\drb:D_2(\drb) \to D_2(\drb)$ by $\drhb \circ (\pi_s \otimes \pi_0)=(\pi_S \otimes \pi_0)\circ \drb$ which agrees on the generating elements $R(\lambda, f), \zeta_1(f), \zeta_1(h)^*, \zeta_2(f)$ with the given equations.

\pfit(i) and (ii): These follow from $\Qs^{\dag}:=J_{\omega}\Qs^*J_{\omega}=\Qs$ and calculations as in lemma \eqref{lm:propbsd} (ii) and (iii). Alternately, these can be easily verified directly for the generating tensors on of $D_2(\drb)$ and so extend to all $D_2(\drb)$

\pfit (iii): For all finite dimensional subspaces $\cX_s \subset \fDL$ we have $\Qs^2=0$ on $\cD_S$, hence $\drhb_s^2=0$,
hence it follows that $\delta^2=0$ on $D(\delta_2)$
\end{proof}
\begin{rem} Some important points to note about the above superderivation:
\begin{itemize} 
\item[(i)] From the definition above we see that $\drb$ preserves $D_2(\drb)$ but \textit{not} $D_1(\drb)$. So $\drb^2$ makes sense on $D_2(\drb)$ but not on $D_1(\drb)$. The reason we define $D_1(\drb)$ is that it is a $*$-algebra whereas $D_2(\drb)$ is not.
\item[(ii)]  Property (iii) above states 
\begin{equation}\label{eq:delka}
\drb(A)^*= - \alpha\circ\delta\circ\alpha \circ \gamma(A^*)
\end{equation}
for $A \in D_2(\drb)$. It should be pointed out that although $(\zeta_2(f))^* \notin D_2(\drb)$, $\alpha((\zeta_2(f))^*) =\zeta_2(-f) \in D_2(\drb)$ and so the $*$ on the RHS of equation \eqref{eq:delka} does not give us domain problems.
\end{itemize}
\end{rem}
To define the Poincar{\'e} transformations on $\cA$ recall the representations of $\cP^{\uparrow}_{+}$
\[
g \to V_g \in \sp(\fD,\sigma_1)  \qquad{and}\qquad g \to S_g \in U(\fDL\oplus_{\mH}\fDJ),
\]
where $U(\fHL\oplus_{\mH}\fLJ)\subset B(\fHL\oplus_{\mH}\fLJ) $ is the group of unitaries in $B(\fHL\oplus_{\mH}\fLJ)$ (cf. Subsection \ref{sbs:hcov} and equation \eqref{eq:Sg}). Using these we define the Poincar{\'e} transformations by the following action:
\begin{proposition}\label{pr:RA1keragcom}
There exists a homomorphism $\alpha_{(\cdot)}:\cP^{\uparrow}_{+} \to \aut(\cA)$ such that for all $f \in \fX$, $h \in \fDL\oplus_{\mH}\fDJ$ and $g \in \cP^{\uparrow}_{+}$:
\begin{equation}\label{eq:autr1}
\alpha_g(R(\lambda,f)\otimes C(h))=R(\lambda,V_gf)\otimes C(S_gh)
\end{equation}
Moreover:
\begin{itemize}
\item[(i)] $\alpha_g(D_1(\drb))= D_1(\drb)$ and   $\alpha_g(D_2(\drb))= D_2(\drb)$.
\item[(ii)]$\alpha_{g^{-1}} \circ \drb \circ \alpha_g = \drb$ on $D_2(\drb)$ and $D_2(\drb)$, hence $\alpha_g(\ker \drb) =\ker \drb$.
\item[(iii)] $\alpha_{g^{-1}} \circ \alpha \circ \alpha_g = \alpha$, where we recall that $\alpha=\iota \otimes \alpha'\in \Aut(\cA)$ for $\alpha'\in \Aut(\cA_g)$ is such that $\alpha'(C(h))=C(Jh)$ for all $h \in \fL$.
\end{itemize}
\end{proposition}
\begin{proof}
Let $g \in \cP^{\uparrow}_{+}$. Before we begin the proof we recall several properties of $V_g$ and $S_g$:
\begin{gather}
V_g \fXL = \fXL, \qquad S_g \fDL = \fDL, \qquad S_g \fDJ= \fDJ, \label{eq:SVcomP}\\
V_{g}^{\dag}=V_{g^{-1}}, \qquad S_{g}^*=S_{g}^{\dag}=S_{g^{-1}},\qquad [S_g,J]=0,\notag \\
[\mH, P_{\LL}]=[\mH, P_{\JJ}]=0 , \qquad \mH P_{\LL} S_g = V_g \mH P_{\LL}, \qquad P_{\JJ}V_g = S_gP_{\JJ}, \notag
\end{gather}
but also recall $[\mH,V_g]$ and $[\mH,S_g]$ need not be zero. The first line above follows from lemma \eqref{lm:covX1} and the definition of $S_g$ (equation \eqref{eq:Sg}). The
second follows from the definition of $V_g$, lemma \eqref{lm:SgKHsa} (ii) which proves that $S_g$ is both $\ip{\cdot}{\cdot}_{\mH}$-unitary and $\iip{\cdot}{\cdot}_{\mH}$-unitary, and $[S_g,J]=0$ follows directly from the definition. The statements in the last line were proved as cases in lemma \eqref{lm:comdrel}.
  
Since $V_g \in \sp(\fX, \sigma_1)$ we have by Theorem \eqref{th:symautRA} (v) that there exists $\alpha_g' \in \Aut(\cR(\fX,\sigma_1))$ such that $\alpha_g'(R(\lambda,f)):=R(\lambda,V_gf)$. As $S_g$ is $\iip{\cdot}{\cdot}_{\mH}$-unitary we get that the exists $\alpha_g''\in \Aut \cA_g$ such that $\alpha_g''(C(h)):=C(S_g h)$ (as it preserves the CAR's). Let $\alpha_g:=\alpha_g'\otimes \alpha_g'' \in \Aut(\cA)$ which satisfies equation \eqref{eq:autr1}. Furthermore $g \to \alpha_g$ is a homomorphism of $\cP^{\uparrow}_{+}$ as $g \to V_g$ and $g \to S_g$ are representations.

\pfit (i): For $f \in \fX$, $h \in \fDL$ and $g \in \cP^{\uparrow}_{+}$ we calculate $\alpha_g$ on the generators of $D_2(\drb)$:
\begin{align*}
\alpha_g(R(\lambda,f)\otimes &\one ) = R(\lambda,V_g f)\otimes \one  \in D_1(\drb),\\
\alpha_g(\zeta_1(h))=&\;R(1, V_g\mH h)\otimes C( S_gh)= R(1,\mH S_g h)\otimes C(S_g h)=\zeta_1(S_g h) \in D_1(\drb),\\
\alpha_g(\zeta_1(h)^*)=&\;R(-1, V_g\mH h)\otimes C( S_gJh)= R(-1,\mH S_g h)\otimes C(JS_g h)=\zeta_1(S_g h)^* \in D_1(\drb),\\
\alpha_g(\zeta_2(f))=&\; R(1, V_g f)\otimes C(S_g P_{\JJ} K f)=R(1, V_g f)\otimes C( P_{\JJ}K V_g  f) = \zeta_2(V_g f)\in D_2(\drb)
\end{align*}
where we have used the equations \eqref{eq:SVcomP} above, in particular $V_g\mH h =V_g\mH P_{\LL} h =\mH P_{\LL} S_g h =\mH S_g h$ in the second and third identities, and $P_{\JJ}V_g = S_gP_{\JJ}$ in the last identity.

As $\alpha_g \in \Aut(\cA)$ and $\alpha_g$ preserves the generators of $D_1(\drb)$ and $D_2(\drb)$, we have that (i) follows.

\pfit (ii): For $f \in \fX$, $h \in \fDL$ and $g \in \cP^{\uparrow}_{+}$ we calculate on the generators of $D_2(\drb)$:
\begin{align*}
(\alpha_{g^{-1}} \circ \drb \circ \alpha_g)(R(\lambda, f)\otimes \one )=&\; (\alpha_{g^{-1}} \circ \drb )(R(\lambda, V_g f)\otimes \one),\\
=&\;-i\alpha_{g^{-1}}( R(\lambda, V_g f)^2 C( KP_{\JJ}V_g f))\\
=&\;-i\alpha_{g^{-1}} (R(\lambda, V_g f)^2 C( S_g KP_{\JJ} f))\\
=&\;-i R(\lambda,  f)^2 C(  KP_{\JJ} f)\\
=&\; \drb(R(\lambda,  f)\otimes \one)
\end{align*}
where we used $KP_{\JJ}V_g=S_g KP_{\JJ}$ in the second line.
\begin{align*}
(\alpha_{g^{-1}} \circ \drb \circ \alpha_g)(\zeta_1(h))=&\; (\alpha_{g^{-1}} \circ \drb)(\zeta_1(S_g h)),\\
=&\;\alpha_{g^{-1}}(iR(1, P_{\LL} \mH S_g h)- \one),\\
=&\; \alpha_{g^{-1}}(iR(1, V_g P_{\LL} \mH  h)- \one),\\
=&\; iR(1, P_{\LL} \mH  h)- \one,\\
=&\; \drb(\zeta_1(h)),
\end{align*}
where we have used $\alpha_g(\zeta_1(h))=\zeta_1(S_g h)$ as calculated above. Furthermore using the calculation of $\alpha_g$ on the generators in (i) we get:
\begin{align*}
(\alpha_{g^{-1}} \circ \drb \circ \alpha_g)(\zeta_1(h)^*)=&\; (\alpha_{g^{-1}} \circ \drb )(\zeta_1(S_g h)^*)=0=\drb(\zeta_1(h)^*)), \\
(\alpha_{g^{-1}} \circ \drb \circ \alpha_g)(\zeta_2(f))=&\; (\alpha_{g^{-1}} \circ \drb )(\zeta_2(V_g f))=0=\drb(\zeta_2(f)).
\end{align*}
Hence we have proved that $\alpha_{g^{-1}} \circ \drb \circ \alpha_g=\drb$ on the generators of $D_1(\drb)$ and $D_2(\drb)$. As $\alpha_{g}$ is an automorphism that preserves $D_1(\drb)$ and $D_2(\drb)$, and as $\drb$ is a superderivation we get that this extends to $\alpha_{g^{-1}} \circ \drb \circ \alpha_g=\drb$.

As this is true for all $g \in \cP^{\uparrow}_{+}$ we have that $\alpha_g(\ker \drb)=\ker \drb$.

\pfit (iii):  As $[J,S_g]=0$ (iii) follows.
\end{proof}

\subsection{States}{\label{sbs:brstIstates}
We want to define the set of states that will give us the correct ghost gradings on the ghost space and also serve as vacuum vectors, hence be Krein positive. These states should be in the kernel of the BRST charge if this exists in a given representation. Now in the above construction, we do not have a BRST charge $Q$ to select physical states, but motivated by the bounded case (lemma \eqref{lm:bQkdel}) we define 
\begin{definition}
Let $\fS_\drb$ be the set of states on $\cA$ of the form $\omega_1\otimes \omega_2$ such that  $\omega_2 \in \fS_g(\cA_g)$, and $(\omega_1\otimes \omega_2)(\delta(A))=0 \quad \forall A \in D(\drb)$.
\end{definition}
Note that we have defined the Krein involution on $\cA$ via $(\iota \otimes \alpha')(A^*)=A^{\dag}$ and so we do not require that $\omega_1 \circ \beta= \omega$. The reason that we do not associate any Krein structure in the physical representations with $\beta \in \aut (\cR(\fX,\sigma_1))$ is that the Hilbert involution of the resolvents in $\cR(\fX,\sigma_1)$ corresponds to the Krein involution on the fields $A(f)$ that we have heuristically associated them to. 

We immediately get the following
\begin{lemma}\label{lm:BRSTalgbas}
Let $f\in \fDL$ and $\omega=\omega_1\otimes \omega_2 \in \fS_\drb$. Then:
\begin{enumerate}
\item [(i)] $\omega(R(1, \mH f)\otimes \one)=-i$
\item[(ii)]  We have $\omega_{1}|_{\cR(\fX,\sigma_1)}\in \fS_D(\cR(\fX,\sigma_1))$ where we have assumed Dirac constraints $\cC$ as in Subsection \ref{sbs:DcI}. Furthermore, $\omega_1$ is nonregular.
\item [(iii)] $\pi_{\omega}(R(1,\mH f)\otimes \one))= -i\one $.
\item [(iv)]$\pi_{\omega}(R(1, \mH Jf)\otimes \one)=0,$ 
\item[(v)] $\omega_1 \neq \omega_1 \circ \beta $.

\end{enumerate}
\end{lemma}
\begin{proof}
(i): This follows from $0=\omega(\drb(\zeta_1(f)))=\omega(iR(1,\mH f)-1)$.

\smallskip \noindent
(ii): It follows from (i) that $\omega_1( R(1, \mH f))=-i$ for all $f \in \fD$. Hence if we take $\fXL$ as our test function constraint space as in Subsection \ref{sbs:DcI}, then by Proposition \eqref{RDirac} (i) $\omega_{1}|_{\cR(\fX,\sigma_1)}\in \fS_D(\cR(\fX,\sigma_1))$. Moreover, by Proposition \eqref{RDirac} (ii) $\omega_1$ is nonregular.
  
\pfit (iii): Follows by (ii) and Corollary \eqref{cr:trivcons}. 

\smallskip \noindent (iv): As $\fXJ=J\fXL$ we have that $g:=\mH Jf= J\mH f \in \fXJ$. By the definition of the inner products on $\fX$ (cf. Subsection \ref{sbs:abstf} (4)) we have for this $g$ and $f$:
\begin{equation}\label{eq:d2ncomd1}
\smp{g}{\mH Kf}{1}=\smp{J\mH f}{K\mH f}{1}=\re \iip{\mH f}{J\mH f}=\ip{\mH f}{\mH f}>0.
\end{equation}
By (ii) we have $\omega_1\in \fS_D(\cR(\fX,\sigma_1))$, and so by Proposition \eqref{RDirac} (ii) and $\smp{g}{\mH Kf}{1}\neq0$ we get $\pi_{\omega_1}(R(1, \mH J  f))=0$. 

\smallskip \noindent
(v): Assume $\omega = \omega \circ \beta $. By part (i) and part (iii),  $-i=\omega(R(1, \mH f))=\omega(\beta(R(1, \mH f)))=\omega(R(1, \mH Jf))=0$ for all $f\in \fDL$.
\end{proof}
Lemma \eqref{lm:BRSTalgbas} (v) says that state in $\fS_{\drb}$ do not have GNS-cyclic vectors that are positive with respect to the Krein-structure associated to $\beta$. As discussed at the beginning of this subsection, $\beta$ does not encode the involution on $\cR(\fX,\sigma_1)$ associated to the heuristic correspondence with $A(f)$, hence lemma \eqref{lm:BRSTalgbas} (v) is not unreasonable.

\subsection{{A problem with the Physical Algebra}}\label{sbs:problem}
The above structure gives
\begin{proposition}\label{pr:rabrstd1}
We have:
\begin{itemize}
\item[(i)] $\fS_{\drb}=\fS_D\otimes \fS_g$, where $\fS_D$ are the Dirac states of $\cR(\fD,\sigma_1)$ using $\fDL$ as a constraint test function space as in Subsection \ref{sbs:racons}.
\item[(ii)] Let $\omega \in \fS_{\drb}$. Then $\pi_{\omega}\circ \drb=0$ on $D_2(\drb)$.
\end{itemize}
\end{proposition}
\begin{proof}
\noindent (i): First $\fS_{\drb}\subset \fS_D\otimes \fS_g$ follows from lemma \ref{lm:BRSTalgbas} (ii). Conversely, let $\omega_1 \otimes \omega_2 \in \fS_D\otimes\fS_g$. Suppose $f \in \fD$ and $P_{\JJ} f \neq 0$. Then for $g:=JKP_{\JJ}f \in \fDL$ we have $\smp{f}{g}{1}=\ip{P_{\JJ}f}{P_{\JJ}f}\neq 0$ (cf.  equation \eqref{eq:d2ncomd1}). So by Proposition \ref{RDirac} (ii) (taking $C=\fDL$), we get that $\pi_{\omega}(R(\lambda, f))=0$ for all $\lambda \in \mathbb{R} \backslash 0$. Using Theorem \eqref{th:calgdrb}, we get that, 
\[
\pi_{\omega}(\drb(R(\lambda, f)))=i\pi_{\omega}(R(\lambda, f))^2\pi_{\omega}( C( P_{\JJ}  Kf))=0.
\]
But if we suppose that $P_{\JJ}f=0$ then $C(P_{\JJ} Kf)=0$ and we can see from the first equality above that again $\pi_{\omega}(\drb(R(\lambda, f)))=0$. Therefore we have that, 
\begin{align*}
\pi_{\omega}\circ \drb(R(\lambda, f)))&= 0,\\
\pi_{\omega}\circ \drb(\zeta_1(f)))&= i\pi_{\omega}(R(1, P_{\LL} \mH f))- \one=0, \\
\pi_{\omega}\circ \drb(\zeta_1(f)^*)&=0, \\
\pi_{\omega}\circ \drb(\zeta_2(f))&=0.
\end{align*}
where the second line follows from lemma \ref{lm:BRSTalgbas} (ii). Now $\pi_{\omega}\circ\delta=0$ on all the generators of $D_2(\drb)$, so $\pi_{\omega}\circ \drb=0$ and hence $\omega \circ \drb =0$ on $D_2(\drb)$ which proves (i). (ii) is immediate from the argument above.
\end{proof}
Proposition \eqref{pr:rabrstd1} (ii) is a big problem as it says that if $\drb_{\omega}(\pi_{\omega}(A)):=\pi_{\omega} (\drb(A))$ for $\omega \in \fS_{\drb}$, then $\drb_{\omega}=0$ and so  $\ker \drb_{\omega}= \pi_{\omega}(D_2(\drb))$ which contains ghost terms, etc, and $\ran \drb_{\omega}=0$ so when passing to the physical algebra, we won't factor out anything! To be more precise and also consistent with previous definitions we make the following definition.
\begin{definition}
Suppose $\omega=\omega_1\otimes \omega_2 \in \fS_\drb$. Let $Q_{\omega}$ be the generator of $\drb_{\omega}$ given above and let
\[
\cH_{\omega}=\cH_d^{\omega}\oplus \cH_s^{\omega}\oplus \cH_p^{\omega}
\]
be the \emph{dsp}-decomposition with respect to $\cQ_{\omega}$. Then we define,
\[
\tilde{\cP}_{\omega}^{BRST}:=\ker \drb_{\omega}/(  \ker \drb_{\omega} \cap \ker \Phi^{\omega}_s) \cong P^{\omega}_s\ker \drb_{\omega} P^{\omega}_s,
\]
where the above isomorphism is assumed to be algebraic, and we have assumed no topology on $\tilde{\cP}_{\omega}^{BRST}$.
\end{definition}

At this stage we do not assume a topology or involution on $\tilde{\cP}_{\omega}^{BRST}$ and so the above isomorphism are purely algebraic. We address these issues in the next section (cf. Theorem \eqref{th:ubdbrststruc}). As $\drb_{\omega}=0$ we have that $Q_{\omega}=0$ hence $\cH_s^{\omega}=\overline{\ker Q_{\omega} \cap \ker Q_{\omega}^*}=\cH_{\omega}$ and $ P^{\omega}_s=\one$. Therefore $P^{\omega}_s\ker \drb_{\omega} P^{\omega}_s= \ker \drb_{\omega}= \pi_{\omega}(D_2(\drb))$. As $\omega(R(1,\mH f))=-i$ we have that $\pi_{\omega}(R_0(\fX,\sigma_1)\otimes \one))=\pi_{\omega}(R_0(\fX_t,\sigma_1)\otimes \one))$ by the same arguments as in Proposition \eqref{pr:RA1phalg}, hence 
\begin{equation}\label{eq:omg1}
\pi_{\omega}(D_2(\drb))=\pi_{\omega}(\cR_0(\fX_t, \sigma_1)\otimes \cA_g)
\end{equation}
However as $\cA_g$ is simple it only has faithful representations and there do not exist representations such that $\pi_{\omega}(\one \otimes C(P_{\LL}f))=0$. Hence  if $\omega \in \fS_{\drb}$ then $\pi_{\omega}(\zeta_1(f))=\pi_{\omega_1}(R(1, P_{\LL} \mH f))\otimes \pi_{\omega_2}(C(P_{\LL}f))=-i\one \otimes\pi_{\omega_2}(C(P_{\LL}f))\neq 0$, and hence there are non-zero ghost number elements in $\pi_{\omega}(\cA)$.

We see that $\omega \in \fS_{\drb}$ selects the same resolvent part of the algebra as the Dirac method as in Proposition \eqref{pr:RA1phalg}, but does not remove the ghost terms. Hence BRST for QEM  defined as in Subsection \ref{sbs:sdI} and Subsection \ref{sbs:brstIstates} does \emph{not} give the same results as the Dirac method as it does not remove the ghosts. This is somewhat surprising as we found that BRST gave the correct results for QEM when viewed from the operator point of view in Chapter \ref{ch:BRSTQEM}. 

The explanation for the difference is that we used the heuristic identification of $R(\lambda, f)$`$=$'$(i\lambda\one-A(f))^{-1}$ when defining $\drb$ and as already discussed, there is no reason to think that  $(i\lambda\one-A(f))^{-1}$ is well defined in a representation having BRST structures.

We will leave the discussion of relativistic covariance for BRST Section \ref{sec:ubbcsbrst} (cf. Proposition \eqref{pr:RA1keragcom2}) .
\begin{rem}
It is of course possible to do BRST purely at an algebraic level, \ie to take $\ker \drb/ \ran \drb$ as the physical algebra, and considering its representations. However, this does not correspond to any BRST-constraining at an operator level.
\end{rem}
\section{$C^*$-BRST II}\label{sec:cstbrstii}

\subsection{$T$-prodedure II}\label{sbs:qemdii}
As we have seen in the previous section, using the algebra $\cR(\fD, \sigma_1)$ in the construction of BRST-QEM does not lead to the correct results for QEM as the ghosts were not removed when passing to the final BRST-Physical Algebra $
\tilde{\cP}_{\omega}^{BRST}$. One reason for this is that we have assumed that $R(\lambda,f) \in\cR(\fD, \sigma_1)$ should in some sense represent the resolvent of the Krein-symmetric field $A(f)$ and so have defined our BRST superderivation accordingly. It turned out to give results different from operator model in Chapter \ref{ch:BRSTQEM} and the $T$-procedure in Subsection \ref{sbs:DcI}, in that it did not remove the ghosts.

To remedy this problem we take an approach that more directly connects to the Krein-symmetric field operator construction in Chapter \ref{ch:BRSTQEM}. We use for the algebra of bosonic fields, $\cR(\fD, \sigma_2)$ where $\smp{f}{g}{2}:=\smp{f}{Jg}{1}$ is the \ncsf  (cf. Subsection \ref{sbs:abstf}). Heuristically the fields producing the symplectic form $\sigma_2$ from its commutators are of the form
\[
\phi(f)= \fst(a(f)+a^*(f)),
\]
where $[a(f), a^*(g)]= \ip{f}{g}$, rather than
\[
A(f)= \fst(a(f)+a^{\dag}(f))=\fst(a(f)+a^*(Jf)),
\]
which produces $\sigma_1$ from its commutators. Since $\phi(f)=\phi(f)^*$, its resolvents $(i \lambda \one - \phi(f))^{-1}$ make sense, so we can use $\cR(\fD,\sigma_2)$ to model the field $\phi(f)$. Solving $A(f)$ and $\phi(f)$ for each other gives:
\begin{align}
\phi(f)&=\frac{1}{2}(A(f+Jf)+iA(K(f-Jf))=A(P_{+}f)+iA(K(P_{-}f)), \label{eq:Aphc}\\
A(f)&=\phi(P_{+}f)+i\phi(K(P_{-}f)). \notag
\end{align}
using $J=P_{+}-P_{-}$. In the Krein space formulation BRST-QEM of  $\delta$ is defined in terms of its action on the $A(f)$'s and $C(f)$'s (Definition \ref{df:brstsd}) and so we will use the above equations to make a rigorous definition of $\drb$ on $\cR(\fD, \sigma_2)\otimes \cA_g$. The idea of using auxiliary fields to generate a $C^*$-algebra for a quantum system was used in a different context in \cite{HendrikHur1993} p505 where the $C^*$-algebra under consideration in that case was the Weyl Algebra.

First, we need to apply the $T$-procedure to constraints in $\cR(\fD, \sigma_2)$. For this we need to specify the constraints corresponding to the smeared Lorentz condition as follows.

Observe that for $f\in \fD$ we have,
\begin{align*}
A(f)^*A(f)=A(Jf)A(f)= & \;(\phi(P_{+}f)-i\phi(KP_{-}f))(\phi(P_{+}f)+i\phi(KP_{-}f),)\\
= & \;\phi(P_{+}f)^2+\phi(KP_{-}f)^2
\end{align*}
via $\smp{P_{+}\fD}{P_{-}\fD}{2}=0$.  
Now $\phi(f)^*=\phi(f)$ \,so\,  $\ker (A(f)^*A(f))= \ker (\phi(P_{+}f) \cap \phi(KP_{-}f))$ hence we can restate the heuristic condition $A(f)\psi=0$ as $\phi(P_{+}f)\psi=\phi(KP_{-}f)\psi=0$. Define the subspaces of $\fX$ as,
\begin{align}\label{eq:ctf2}
 \fT_1:= & \; P_{+}P_{\LL} \fD ,\\
 \fT_2:= & \;  P_{-}P_{\LL} \fD
\end{align}
then we associate $\{ \phi(f), \phi(g)\,|\, f \in \fT_1, g\in \fT_2\}$ with the smeared complexified Lorentz condition for the symplectic form $\sigma_2$. 
We would like to perform the $T$-procedure using Proposition \eqref{pr:conRA} using $\fT_1$ and $\fT_2$ as the test function spaces corresponding to the constraints. Unfortunately $\fT_1 \cup \fT_2$ is not $\sigma_2$-symplectically neutral since $\smp{f}{Kf}{2}=\norm{f}^2_{\fH}$ and $K\fT_1=\fT_11$, $K\fT_2=\fT_2$. Hence Proposition \eqref{RDirac} (iii) shows that are no Dirac states for these constraints.

We would like to find $\sigma_2$-neutral subspaces of $\fT_1$ and $\fT_2$ to use for constraints in the $T$-procedure, cf. Proposition \eqref{pr:conRA}. In QEM there exists a conjugation on $\fX$ by which we can do this.

Let $\widehat{\cS},\widehat{\cS_0},\fX, \fX_0,\fD,\fD_0$ be the  test function spaces for QEM as in Subsection \ref{sbs:testfunc}, with all the structures assumed there. It follows from Proposition \eqref{lm:stdec2} and Theorem \eqref{th:fdtstdec} that we have
\begin{align}
\fT_{\LL}:= & \; P_{+}P_{\LL} \fD =\{ f \in \fD \,|\, f(p)=(0, p_1h(p), p_2h(p),p_3h(p)),\,h \in \fH_0 \}\label{eq:c2ra}\\
 \fT_{\JJ}:= & \; P_{-}P_{\LL} \fD=\{ f \in \fD \,|\, f(p)=(ip_0h(p), 0, 0,0),\,h \in \fH_0 \}\notag .
\end{align}
Note that $\fDL \oplus \fDJ= \fT_1 \oplus \fT_2$ where $\oplus$ is Hilbert orthogonality with respect to $\ip{f}{g}=\smp{f}{Kg}{2}+i\smp{f}{g}{2}$. Now define for the scalar functions $f\in \widehat{\cS_0}$:
\[
(C_0 f )(p):=\overline{f(p_0, -p_1, -p_2, -p_3)},
\]
so that $C_0$ is a well defined real linear (but not complex linear) operator such that $C_0^{2}=\one$. Furthermore, as $p=(p_0, p_1, p_2, p_3) \in C_{+} \iff p=(p_0, -p_1, -p_2, -p_3) \in C_{+}$ we have that if $f|_{C_{+}}=0$ then $C_0f|_{C_{+}}=0$. Hence we have that $C_0$ factors to a (real) linear operator on $\fX_0$ which we will still denote by $C_0$. 
\begin{proposition}\label{pr:QM2tfdec}
Define the real linear operator $C:\fX \to \fX$ by:
\[
Cf:=(C_0f_0,C_0f_1, C_0f_2, C_0f_3)
\]
for $f\in \fX$. Then
\begin{itemize}
\item[(i)] $C^2=\one$.
\item[(ii)] $CK=-KC$ and so $C$ extends to $\fD$ as an \emph{antilinear} operator.
\item[(iii)] $\smp{Cf}{Cg}{2}=\smp{g}{f}{2}=-\smp{f}{g}{2}$ for all $f,g \in \fD$.
\item[(iv)] Let $\fM \subset \fD$ be a complex linear subspace such that $C\fM=\fM$ and define  
\[
\fM_C:=\{f \in \fM\,|\, Cf=f\}.
\]
Then $\fM=\fM_C \oplus K\fM_C$  and  $K\fM_C=\{f \in \fM\,|\, Cf=-f\}$ where $\oplus$ means algebraic direct sum of real linear spaces. Furthermore $\fM_C$ and $K\fM_C$ are $\sigma_2$-neutral real linear spaces. 
\item[(v)] $C\fT_{\LL}=\fT_{\LL}$ and $C\fT_{\JJ}=\fT_{\JJ}$ hence we have the decomposition:
\[
\fDL\oplus \fDJ=\fT_{\LL}\oplus\fT_{\JJ}= (\fT_{\LL C} \oplus K\fT_{\JJ C})\oplus (K\fT_{\LL C} \oplus \fT_{\JJ C}),
\]
where $\oplus$ means algebraic direct sum of real linear spaces. Furthermore if we define $\fT:=\fT_{\LL C} \oplus K\fT_{\JJ C}$ then:
\[
\fD=\fD_t\oplus \fT \oplus K\fT 
\]
and $\smp{\fT}{\fT}{2}=\smp{K\fT}{K\fT}{2}=\smp{\fD_t}{\fT}{2}=\smp{\fD_t}{K\fT}{2}=0$ 
\end{itemize}
\end{proposition}
\begin{proof}
(i) and (ii): Obvious.

\pfit (iii): By $CK=-KC$ and the polarization identity (R\&S \cite{ReSi1972v1} p63) we have that,
\[
\ip{Cf}{Cg}=\ip{g}{f},
\] 
where we recall $\ip{\cdot}{\cdot}$ is the inner product given by $\ip{\cdot}{\cdot}=\smp{\cdot}{K\cdot}{2}+i\smp{\cdot}{\cdot}{2}$, and so 
(iii) follows.

\pfit (iv): Let $f\in \fM$. Now $f= (f+Cf)/2+ K(-K(f-Cf)/2)$ and as $C\fM \subset \fM \supset K\fM$ we get $\fM=\fM_C\oplus  K\fM_C$ (uniqueness of the decomposition is a trivial exercise). Furthermore, as $CK=-KC$ and $C^2=\one$ we get $K\fM_C=\{f \in \fM\,|\, Cf=-f\}$. From (iii) we have for $f,g \in \fM_C$ or $f,g\in K\fM_C$ that,
\[
\smp{f}{g}{2}=\smp{Cf}{Cg}{2}=-\smp{f}{g}{2},
\]
and so $\fM_C$ and $K\fM_C$ are $\sigma_2$-neutral spaces. 

\pfit (v): Using equation \eqref{eq:c2ra} above we see that $C\fT_{\LL}=\fT_{\LL}$, $C\fT_{\JJ}=\fT_{\JJ}$ hence 
the decomposition:
\[
\fDL\oplus \fDJ=\fT_{\LL}\oplus\fT_{\JJ}= (\fT_{\LL C} \oplus K\fT_{\JJ C})\oplus (K\fT_{\LL C} \oplus \fT_{\JJ C})=\fT\oplus K\fT,
\]
follows by (iv), hence $\fD=\fD_t\oplus \fT \oplus K\fT$. 

Finally, it follows from $\smp{\fD_t}{\fT_{\LL}}{2}=\smp{\fD_t}{\fT_{\JJ}}{2}=\smp{\fT_{\LL}}{\fT_{\JJ}}{2}=0$, that $K\fD_t=\fD_t$, $K\fT_{\LL}=\fT_{\LL}$, $K\fT_{\JJ}=\fT_{\JJ}$, that $K$ is $\sigma_2$-symplectic and the above decompositions of $\fD$.

\end{proof}
\begin{rem}
Note that the use of $C$ to decompose $\fM$ into neutral subspaces as above requires that $C$ preserve $\fM$. By Proposition \eqref{lm:stdec2} and Theorem \eqref{th:fdtstdec} we have,
\[
\fXL=P_{\LL} \fX =\{ f \in \fX \,|\, f(p)=i((p_0, p_1, p_2,p_3)/\np) h(p),\,h \in \fH_0 \}
\]
and so $C$ does \emph{not} preserve $\fXL$ or $\fDL$ and hence we do \emph{not} get an analogous lemma using $\fXL$ in place of $\fT_{\LL}$.
\end{rem}
By the argument leading to  equation \eqref{eq:ctf2} we have that the use of $\fT$ as test function constraint space corresponds heuristically to the constraint set $\{ \ker A(f)\,|\, f \in \fT^1_C\oplus \fT^2_C\}$. So rigorously we take the quantum system with constraints to be $(\cR(\fD, \sigma_2), \cC_2)$ where $\cC_2$ is the set, 
\[
\cC_2:=\{ i \lambda R(\lambda,f)-\one \,|\, f \in \fT \}.
\]
\begin{proposition}\label{pr:physalgQEMii}
Let $(\cR(\fD, \sigma_2),\cC_2)$ be a quantum system with constraints as in Appendix \ref{app:Tp}. Then the Dirac observable algebra is:
\[
\cP \cong \cR(\fD_t, \sigma_2),
\]
where $\cP=\cO/\cD$ is the maximal $C^*$-algebra of physical observables.
\end{proposition}
\begin{proof}
By Proposition \eqref{pr:QM2tfdec} (v)  we have $\smp{\fT}{\fT}{2}=\smp{\fD_t}{\fT}{2}=\smp{\fD_t}{K\fT}{2}=0$, and $f \in \fT\backslash 0$ does not symplectically commute with $K\fT$. These are the analogous inputs in the arguments that lead to Proposition \eqref{pr:RA1phalg} with $\fDL$ replaced by $\fT$. Hence we get
\[
\cP \cong \cR(\fD_t, \sigma_2).
\]
\end{proof}
Therefore we have that the $T$-procedure applied to $(\cR(\fD, \sigma_1), \cC_1)$ or $(\cR(\fD, \sigma_2), \cC_2)$ give the same result.

\subsection{Superderivation II}

To construct BRST-QEM using the \ncss we follow a similar procedure to Subsection \ref{sbsc:cbrstv1} except now we use the Resolvent Algebra with the \ncsf $\sigma_2$. For notational efficiency, we use much of the same notation as in Subsection \ref{sbsc:cbrstv1}, \eg $R(\lambda, f)$ will denote the generating elements of $\cR(\fD, \sigma_2)$:
\begin{itemize}
\item Let $\fD$ and $\fL$ have all the structures as in Subsection \eqref{sbs:abstf}. Let $\cA_g(\fLJ)$  be the ghost algebra  with all the structures defined in Section \ref{sec:ghost}, e.g.
\[
C(f):=\fst( c( f)+ c^*(J f)), \qquad f \in \fL.
\]
\item We define the BRST-Field Algebra as
\[
\cA:=\cR(\fD, \sigma_2)\otimes \cA_g(\fLJ).
\]
The tensor norm on $\cA$ is unique as the CAR algebra is nuclear. We define a grading on $\cA$ by extending the ghost grading on $\cA_g$, \ie we define a grading automorphism $\gamma$ on $\cA$ as $\gamma$ equal to the identity on $\cR(\fD, \sigma_2)$, and equal to the $\Z_2$-grading on $\cA_g$ (cf. Definition \eqref{df:Z2grad}). 
\item As $J$ is symplectic, we have that $\beta'(R(\lambda, f)):= R( \lambda, Jf)$ defines a unique automorphism on $\cR(\fD, \sigma_2)$ (cf. Theorem \eqref{th:symautRA} (v)). As $J$ also defines a unitary on $\fL$, we get that $\alpha'(C(f)):=C(Jf)$ defines a unique automorphism on $\cA_g$ and we let $\alpha=\alpha' \otimes \beta' \in \Aut(\cA)$ which encoded the Krein structure on $\cA$. We define the involution:
\[
A\to A^{\dag}:=\alpha(A^*), \qquad A \in \cA.
\] 
\end{itemize}
 
To define the BRST superderivation rigorously we first use equations \eqref{eq:Aphc} and Definition \eqref{df:brstsd} to get the heuristic equations:
\begin{align*}
\drb(\phi(f) \otimes \one)=&\; \drb(A(P_{+}f)\otimes \one)+i\drb(A(K(P_{-}f))\otimes \one)= -\one \otimes (i C(KP_{\JJ}P_{+}f)+C(P_{\JJ}P_{-}f)),\\
\drb(\one \otimes C(g))=&\; A(\mH P_{\LL} g) \otimes \one =( \phi(\mH P_{+}P_{\LL}g)+i\phi(K\mH P_{-}P_{\LL}g))\otimes \one,
\end{align*}
where $f \in \fD$ and $g \in \fDL\oplus_{\mH}\fDJ$. We can interpret this rigorously using the Resolvent Algebra, but first we must choose an appropriate domain in $\cA$ for $\drb$. Again there is a problem of defining bounded operators corresponding to $\drb(\one \otimes C(f))$ since the $\phi$'s on the RHS of this expression are unbounded operators. We will use mollifiers to encode this expression in $\cA$. 
\begin{definition}\label{df:sddom2}
Let
\[
D_1(\drb)=\salg{\one, \, R(\lambda, f)\otimes \one , \, \zeta_1(g),\zeta_2(g)\,| \, f \in \fD, g \in \fDL \, \lambda\in \mathbb{R} \backslash 0},
\]
and, 
\[
D_2(\drb)=\alg{ D_1(\drb)\,,\, \one \otimes C(g) | \, f \in \fD, g \in \fDJ},
\]
where $\zeta_1(g):=R(1, \mH P_{+}g)R(1,K\mH P_{-}g)\otimes C(g)$ and $\zeta_2(g)= R(1, -\mH P_{+}g)R(1,K\mH P_{-}g)\otimes C( g)$ for $g \in \fDL$.
\end{definition}
\begin{rem} \label{rm:DBRSII}
We have:
\begin{itemize}
\item[(i)]  $D_1(\drb)\subset D_2(\drb)$ but $D_1(\drb)$ is a $*$-algebra  whereas $D_2(\drb)$ is not. Neither is norm dense in $\cA$.
\item[(ii)] Note that by $R(\lambda,g)^*=R(-\lambda,g)=-R(\lambda,-g)$ $Jf \in \fDJ$ for $f \in \fDL$, and $\smp{P_{+}}{K\mH P_{-}f}{2}=0$ we have 
\begin{align*}
\zeta_1(f)^*=&\;R(-1, \mH P_{+}f)R(-1,K\mH P_{-}f)\otimes C(J f),\\
=&\;R(1, -\mH P_{+}f)R(1,-K\mH P_{-}f)\otimes C(J f),\\
 \in &\; D_2(\drb),
\end{align*}
for all $f \in \fDL$.
\item[(iii)] Using $JP_{+}=P_{+}$ and $JP_{-}=-P_{-}$ we have for all $f \in \fDL$ that $\zeta_2(f)=\alpha(\zeta_1(f)^*)=\zeta_1(f)^{\dag}$. This implies that 
\[
\zeta_2(f)^*=\alpha(\zeta_1(f))= R(1, -\mH J P_{+}f)R(1,-K\mH JP_{-}f)\otimes C( Jf) \in D_2(\drb)
\]
as $Jf \in \fD_2$. Furthermore $\alpha(\zeta_1(f))=\alpha(\zeta_1(f)^*)^*=\zeta_2(f)^*\in D_2(\drb)$.
\item[(iv)] From (ii) and (iii) and $\alpha(C(f)^*)=C(f)$, we have that $\alpha(A^*) \in D_2(\drb)$ for all $A$ that are generators of $D_2(\drb)$ and hence all $A \in D_2(\drb)$.
\item[(iii)] Recall that $\gamma$ is the $\Z_2$-grading automorphism on $\cA$, hence $\gamma(R(\lambda,f)\otimes \one)=R(\lambda,f)\otimes \one$ and $\gamma(\one \otimes C(g))=-\one \otimes C(g)$ for all $f \in \fD$ and $g \in \fDL \oplus_{\mH} \fDJ$ (cf. Remark \eqref{rm:Zwgr}). Hence $\gamma(D_1(\drb))=\gamma(D_1(\drb)$ and $\gamma(D_2(\drb))=D_2(\drb)$. This says that $D_2(\drb)$ is a $\dag$-algebra.
\end{itemize}
\end{rem}
To reduce notation we define for all $f \in \fDL$:
\[
R_{+}(f):= R(1, \mH P_{+}f), \qquad R_{-}(f)=R(1,\mH P_{-}f),
\]
and so 
\begin{gather*}
\zeta_1(f)=R_{+}(f)R_{-}(Kf)\otimes C(f), \quad \text{and,} \quad \zeta_1(f)^*=R_{+}(-f)R_{-}(-Kf)\otimes C(Jf),\\
\zeta_2(f)=R_{+}(-f)R_{-}(Kf)\otimes C(f), \quad \text{and,} \quad \zeta_2(f)^*=R_{+}(f)R_{-}(-Kf)\otimes C(Jf).
\end{gather*}
Also note that as $P_{+}P_{-}=0$ we get using Definition \eqref{df:RA} \eqref{eq:R5} that
\[
[R_{+}(f),R_{-}(g)]=0, \qquad \forall f,g \in \fDL,
\]
which we will use frequently in calculations. We will also drop the tensor product $\otimes$ in Theorem \eqref{th:calgdrb2} below. We define the mollified version of the BRST superderivation for this context, the proof of which follows the method in Theorem \eqref{th:calgdrb}.
\begin{theorem}\label{th:calgdrb2}
Define a map on the generating elements of $D_2(\drb)$ as follows, for $f\in \fD$, $g \in \fDL$, $h\in \fDJ$,  $\lambda\in \mathbb{R} \backslash 0$:
 \begin{align*}
\drb(R(\lambda, f)\otimes \one)= & \; - R(\lambda, f)^2 (iC( K P_{\JJ}P_{+} f)+ C(P_{\JJ}P_{-}f))  \\
= & \;- (1/\sqrt{2}) R(\lambda, f)^2(c(P_{\JJ}f)-c^*(P_{\JJ}Jf)) ,  \\
\drb(\one \otimes C(h)) = & \;0,\\
\drb(\zeta_1(g))= & \; -R_{+}(g)R_{-}(Kg)[iR_{+}(g)C(KP_{\JJ}\mH P_{+}g)-R_{-}(Kg)C(KP_{\JJ}\mH P_{-}g)]C(g)  \\ 
 &+ (i-1)R_{+}(g)R_{-}(K g))-iR_{+}(g)-R_{-}(Kg)\\
= & \; -[iR_{+}(g)C(KP_{\JJ}\mH P_{+}g)-R_{-}(Kg)C(KP_{\JJ}\mH P_{-}g)]\zeta_1(g)  \\ 
 &+ (i-1)R_{+}(g)R_{-}(K g))-iR_{+}(g)-R_{-}(Kg) \\
\drb(\zeta_1(g)^*)= & \;-R_{+}(g)R_{-}(Kg)[iR_{+}(g)C(KP_{\JJ}\mH P_{+}g)-R_{-}(Kg)C(KP_{\JJ}\mH P_{-}g)]  C( J g) ), \\
\drb(\zeta_2(g))= & \; R_{+}(g)R_{-}(-Kg)[iR_{+}(g)C(KP_{\JJ}\mH P_{+}g)+R_{-}(-Kg)C(KP_{\JJ}\mH P_{-}g)]C(g)  \\
 & +(i+1)R_{+}(g)R_{-}(-K g))-iR_{+}(g)-R_{-}(Kg)\\
=& \; [iR_{+}(g)C(KP_{\JJ}\mH P_{+}g)+R_{-}(-Kg)C(KP_{\JJ}\mH P_{-}g)]\zeta_2(g)  \\
 & +(i+1)R_{+}(g)R_{-}(-K g))-iR_{+}(g)-R_{-}(Kg) \\
\drb(\zeta_2(g)^*)= & \;-R_{+}(g)R_{-}(Kg)[iR_{+}(g)C(KP_{\JJ}\mH P_{+}g)-R_{-}(Kg)C(KP_{\JJ}\mH P_{-}g)]  C( J g) ) , 
\end{align*}
The image of this map is a subset of $D_2(\drb)$, and furthermore this extends to a superderivation on $\drb:D_2(\drb)\to D_2(\drb)$ such that:
\begin{itemize}
\item[(i)] $\drb^2=0$ on $D_2(\drb)$.
\item[(ii)] $\gamma \circ \drb \circ \gamma= -\drb$ on $D_2(\drb)$.
\item[(iii)] $\drb(A)^*=- \alpha \circ \drb \circ \alpha \circ \gamma(A^*)$ on $D_2(\drb)$.
\end{itemize}
\end{theorem}
\begin{proof}
First we verify that $\drb$ is a superderivation on $D_2(\drb)$. We follow \cite{HendrikBuch2006} p708 (cf. proof of Theorem \eqref{th:calgdrb} also). Let $\pi_S$ be a strongly regular (hence faithful) representation of $\cR(\sigma_1, \cX)$. As $\pi_S$ is regular we have that $\phi_{\pi_S}(f)\in \op(\cH_s)$ exists for all $f\in \fD$ and have the properties given by Theorem \eqref{RegThm}. Furthermore, as $\pi_S$ is strongly regular, there exists a dense invariant domain $\cD_{\infty}$ for all $\phi_{\pi_S}(f)$, $f\in \fD$. Thus by Theorem \eqref{RegThm} (i), the resolvents $\pi_{S}(R(\lambda,f))$ map $\cD_{\infty}$ back into $D(\phi_S(f))$ for all $f \in \fD$. Thus we can define a second dense invariant domaint $\cD_S$ by applying all polynomials in $\phi_{\pi_S}(f)$ and $\phi_{\pi_S}(R(\lambda,f))$ to $\cD_{\infty}$. In particular, we can form the non-normed $*$-algebra,
\[
\cE_0:= \salg{\phi_{\pi_S}(f), \, \pi_S(R(\lambda, f)) \, | \, f\in \fD, \, \lambda\in \mathbb{R} \backslash 0}
\]
which acts on the common dense invariant domain $\cD_S $.

Let $\pi_0$ be any representation of $\cA_g$ (which is faithful as $\cA_g$ is simple), and so $\pi_S \otimes \pi_0$ is a faithful representation of $\cA$. Furthermore $\cD:= \cD_S \otimes \cH_g$ is a common dense invariant domain for $\pi_S(R(\lambda, f)) \otimes \one$, $\phi_{\pi_S}(f)\otimes \one$, and $\one \otimes \pi_0(C(f))$. For convenience of notation we will drop the $\otimes$ and $\pi_S$, $\pi_0$ for the remainder of this subsection. Let,
\[
\cE:=\salg{\phi(f),\, R(\lambda,f), C(g) \,|\,f\in \fD, \,g\in \fDL\oplus \fDJ, \, \lambda\in \mathbb{R} \backslash 0 }
\]
and define the map $\drhb$ from the generators of $\cE$ to the generators of $\cE$ as follows: For $f \in \fD$, $g \in \fDL \oplus \fDJ$
\begin{align*}
\drhb(\phi(f))&=  (iC( K P_{\JJ}P_{+} f)+ C(P_{\JJ}P_{-}f)),  \\
\drhb(C(h))&= \phi(\mH P_{\LL} P_{+} g)+i \phi(\mH KP_{\LL} P_{-}g),\\
\end{align*}
and show that this extends to a well defined superderivation on $\cE$. To do this all we have to is show that $\drhb$ is linear and satisfies the graded Leibniz rule on any finite polynomial in the operators $\phi(f),\, R(\lambda,f), C(h), \zeta_2(g), \zeta_2^*(g)$ where $f\in \fD, g\in \fDL, h\in \fDL \oplus \fDJ, \lambda\in \mathbb{R} \backslash 0$.

Let $\cXs$ be a finite-dimensional subspace of $\fDL$ and let $\fD_s= \cXs \oplus J \cXs$. Let 
\[
\cE(\fD_s):=\alg{\phi(f),\, R(\lambda,f), C(g) \,|\, f \in \fD_t \oplus \fD_2,\,g \in \fD_s, \, \lambda\in \mathbb{R} \backslash 0} \subset \cE,
\] 
and let $(f_j)_{j\in \Lambda}$ be a finite orthonormal basis for $\cXs$  and define,
\[
\bch_s:= \sum_{j} \left( (\phi(P_{+} f_j)+i \phi(KP_{-}f_j))C( Jf_j) + (\phi(KP_{+} f_j)-i\phi(P_{-}f_j))  C(K Jf_j) \right).
\]
Let 
\[
\drhb_s(A)\psi:=[\Qs, A \}\psi, \qquad \psi \in \cD, \qquad A \in \cE(\fD_s)
\] 
This is the same formula which defined $\Qs$ in lemma \eqref{lm:2nQs}, but we have substituted $(\phi(P_{+}f) + i \phi (KP_{-}f))$ for $A(f)$. Thus by the same calculations as in the proof of lemma \eqref{lm:2nQs} we get $\Qs^2=0$

Theorem \eqref{RegThm} (iv) gives that $[\phi(f), \phi(g)]=i \smp{f}{g}{2}\one$ on $\cD$, so given $g\in (\fD_t \oplus \fD_s), \lambda\in \mathbb{R} \backslash 0$ then $\sum_j \ip{f_j}{g}f_j=P_{\LL}g$ and $\sum_j \ip{Jf_j}{g}Jf_j=P_{\JJ}g$, and we calculate 
\begin{align*}
\drhb_s(\phi(g))= & \; [\Qs, \phi(g)],\notag\\
& \;\sum_j [i\smp{P_{+}f_j}{g}{2}-\smp{KP_{-}f_j}{g}{2}]C(Jf_j)\\
&+[i\smp{KP_{+}f_j}{g}{2}+\smp{P_{-}f_j}{g}{2}]C(KJf_j), \notag\\
= & \;-i\sum_j C(K[K\smp{P_{+}f_j}{g}{2}+\smp{P_{+}f_j}{Kg}{2}]Jf_j), \notag \\
&+\sum_j C([\smp{P_{-}f_j}{Kg}{2}+K\smp{P_{-}f_j}{g}{2}]Jf_j), \notag \\
= & \;\sum_j (-i C(K\ip{P_{+}f_j}{g}Jf_j) +  C(\ip{P_{-}f_j}{g}Jf_j), \notag\\
= & \;-iC(KJ P_{\LL} P_{+}g)+ C(JP_{\LL}P_{-}g), \notag\\
= & \;-(iC(K P_{\JJ} P_{+}g) +C(P_{\JJ}P_{-}g)), 
\end{align*}
where we have used that $P_{\pm}^*=P_{\pm}$ in the second last line, and $JP_{\LL}=P_{\JJ}J$, $JP_{+}=P_{+}$ and $JP_{-}=-P_{-}$ in the last. Noting that the CAR's for the $C(f)$ use the inner product $\ip{\cdot}{\cdot}_{\mH}=\ip{\cdot}{\mH \cdot}$ on $\fLL\oplus \fLJ$, we get for $g\in \fD_s$ (similar to calculation for $\drhb_s(\phi(g))$),
\begin{align*} 
\drhb_s(C(g))= \phi(\mH P_{\LL} P_{+}g)+ i\phi(\mH K P_{\LL} P_{-}g),
\end{align*}
Now Theorem \eqref{RegThm} (vii) gives that,
\[
[\phi(f), R(\lambda, g)]=i \smp{f}{g}{2}R(\lambda, g)^2\; on \;\cD.
\]
by which we calculate the identity 
\begin{equation}\label{eq:drbs3}
\drhb_s(R(\lambda, f)\otimes \one)= - R(\lambda, f)^2 \otimes (iC( K P_{\JJ}P_{+} f)+ C(P_{\JJ}P_{-}f))
\end{equation}
for $f\in (\cD_t\oplus \cD_s)$, similar to the calculation for $\drhb_s(\phi(g))$ above. Using that $C(f)=\fst(c(f)+c^*(Jf))$, the fact that $f \to c(f)$ is antilinear and $f\to c^*(f)$ is linear for all $f \in \fDL\oplus_{\mH}\fDJ$, and $J=P_{+}-P_{-}$ we get that
\begin{align*}
iC( K P_{\JJ}P_{+} f)+ C(P_{\JJ}P_{-}f)=&\;\fst(c( P_{\JJ}P_{+} f)+c(P_{\JJ}P_{-}f)-c^*( P_{\JJ}P_{+} f)+c^*(P_{\JJ}P_{-}f)),\\
=&\; \fst(c( P_{\JJ} f)-c^*( P_{\JJ}J f))
\end{align*}
hence from equation \eqref{eq:drbs3}:
\begin{equation*}
\drhb_s(R(\lambda, f)\otimes \one)= -\fst R(\lambda, f)^2 \otimes (c( P_{\JJ} f)-c^*( P_{\JJ}J f))
\end{equation*}
Since $\drhb_s$ agrees with $\drhb$ on the generators $\{\phi(f),R(\lambda,f),C(g)\,|\, f \in \fD_t \oplus \fD_S, \, g \in \fD_s\}$ for $\cE(\fD_s)$ this shows that $\drhb$ extends to a superderivation on each $\cE(\fD_s)$, hence to a superderivation on $\cE$. Since $\drhb(\pi(D_2(\drb))\subset \pi(D_2(\drb)$, this allows us to define a graded derivation $\drb: D_2(\drb) \to D(\drb)$ by 
\[
\drhb(\pi(\drb(A)))=:\pi(\drb(A)), \qquad \forall A \in D_2(\drb),
\]
where $\pi:=\pi_S \otimes \pi_0$. We check that $\drb$ agrees with the stated values on the generators
\[
\{R(\lambda, f), \zeta_1(g), \zeta_1(g)^*,\zeta_2(g), \zeta_2(g)^*\,|\,f \in (\fD_t \oplus \fD_s),\, g \in \fD_s \}
\]
of $D_2(\drb)$. We give the calculation for $\drhb_s(\zeta_1(g))$ where $g \in (\cD_s\cap \fDL)$, the calculation for the other generators being similar.
\begin{align*}
\drb(\zeta_1(g))=&\; \drb( R_{+}(g))R_{-}(Kg)C(g)+R_{+}(g)\drb(R_{-}(Kg))C(g)+R_{+}(g)R_{-}(Kg)\drb(C(g)),\\
=&\; -iR_{+}(g)^2C(KP_2P_{+}\mH g)R_{-}(Kg)C(g)+R_{+}(g)R_{-}(Kg)^2C(P_2P_{-}\mH Kg)C(g)\\
& +R_{+}(g)R_{-}(Kg)(\phi(\mH P_{+} g) + i \phi(K\mH P_{-} g)),\\
 =&\; -iR_{+}(g)C(KP_2P_{+}\mH g)R_{-}(Kg)C(g)+R_{+}(g)R_{-}(Kg)C(P_2P_{-}\mH Kg)]\zeta_1(g)\\
 &+(i-1)R_{+}(g)R_{-}(Kg)-iR_{+}(g)-R_{-}(Kg),
\end{align*}
using
\[
R_{+}\phi(\mH P_{+}g)= i R_{+}(g)-\one \qquad \text{and} \qquad R_{-}(Kg)\phi(\mH P_{-} Kg)= R_{-}(Kg)-\one
\]
by Theorem \eqref{RegThm} (vi).

\pfit(i): As already noted for $\cXs$ a finite dimensional subspace of $\fDL$, we have that $\bch_s^2=0$. Now for $A\in \cE(\fD_s)$ we have that $\drhb^2(A)= \drhb( \drhb_s(A))= \drhb_{s'}(\drhb_s(A))=\drhb_{s'}^2(A)=0$ for some finite dimensional $\cXsp \supset \cXs$. Therefore $\drhb^2=0$ and as $\pi_S \otimes \pi_0$ is faithful, we get that $\drb^2=0$ on $D_2(\drb)$.

\pfit(ii) and (iii): These follow from $\Qs^{\dag}:=J_{\omega}\Qs^*J_{\omega}=\Qs$ and calculations as in lemma \eqref{lm:propbsd} (ii) and (iii). Alternately, these can be easily verified directly for the generating tensors on of $D_2(\drb)$ and so extend to all $D_2(\drb)$
\end{proof}
\begin{rem} Some important points to note about the above superderivation:
\begin{itemize} 
\item[(i)] From the definition above we see that $\drb$ preserves $D_2(\drb)$ but \textit{not} $D_1(\drb)$. So $\drb^2$ makes sense on $D_2(\drb)$ but not on $D_1(\drb)$. The reason we define $D_1(\drb)$ is that it is a $*$-algebra whereas $D_2(\drb)$ is not.
\item[(ii)]  Theorem \eqref{th:calgdrb2} (iii) above states 
\begin{equation}\label{eq:delka2}
\drb(A)^*= - \alpha\circ\delta\circ\alpha \circ \gamma(A^*)
\end{equation}
for all $A \in D_2(\drb)$. By Remark \eqref{rm:DBRSII}(iv) we have that $\alpha(A^*) \in D_2(\drb)$ for all $A \in D_2(\drb)$ and so the $*$ on the RHS of equation \eqref{eq:delka2} does not give us domain problems. 
\end{itemize}
\end{rem}
\subsection{Strongly Regular States and Charge}\label{sbs:srs}
Now that the BRST structures are defined at the $C^*$-algebraic level using $\cR(\fD,\sigma_2)$, we will investigate what states give the structures as in Section \ref{sec:brstext}, and also calculate their BRST-physical subspace and BRST-physical algebra. We will define the BRST physical states as those which: 
\begin{itemize}
\item[(i)] Produce the natural Krein structures as related to the $J$-automorphisms in their GNS representation.
\item[(ii)] Produce ghost gradings as in Subsection \ref{sbs:ghgrad}.
\item[(iii)] Have associated cyclic vector being positive with respect to the Krein structure.
\item[(iv)] Have a (possibly unbounded) BRST charge $Q$ in the associated representation which generates $\drb$ on $D_2(\drb)$ in that representation, and has GNS-cyclic vector in $\ker Q$ (\ie $Q$ selects the vacuum).
\end{itemize}
Guided by the bounded $Q$ case we define (cf. Section \ref{sbs:bddcbrst}, in particular lemma \eqref{lm:bQkdel}):
\begin{definition}
Let $\fS_\drb$ be the set of states on $\cA$ of the form $\omega_1\otimes \omega_2$ such that $\omega_1\circ \beta=\omega_1\in \fS(\cR(\fD,\sigma_2))$, $\omega_2 \in \fS_g(\cA_g)$, and $(\omega_1\otimes \omega_2)(\delta(A))=0 \quad \forall A \in D(\drb)$.
\end{definition}
We easily see that the above definition of $\fS_{\drb}$ give representations that satisfy the first three criteria above:
\begin{itemize}
\item[(i)] Recall the Definition \eqref{df:csksrepterm}. For $\omega=\omega_1\otimes \omega_2\in \fS_{\drb}$ we have that $\omega\circ \alpha =\omega$ where $\alpha=\beta'\otimes \alpha'\in \aut(\cA)$ is the automorphism that encodes the Krein structure on $\cA$, \ie $A^{\dag}=\alpha(A^*)$ for all $A \in \cA$. 
From $\beta'^2=\iota$, $\alpha'^2=\iota$, $\alpha^2=(\alpha'\otimes \beta')^2=\iota$ and $\omega_1\circ \beta'=\omega_1$, $\omega_2\circ \alpha'=\omega_2$, $\omega\circ \alpha=\omega$ we see that $(\cH_{\omega_1},J_{\omega_1})$, $(\cH_{\omega_2},J_{\omega_2})$, $(\cH_{\omega},J_{\omega})$ are all Krein spaces as in Definition \eqref{df:csksrepterm}. 
\item[(ii)] By $\omega_2\in \fS_g(\cA_g)$ we have the correct ghost gradings on $\cA_g$ (cf. Proposition \eqref{pr:ghsp}).
\item[(iii)] Also $\omega \circ \alpha =\omega$ implies that $J_{\omega}\Omega_{\omega}=\Omega_{\omega}$ and so $\Omega_{\omega}$ is positive with respect to the inner product $\iip{\cdot}{\cdot}_{\omega}=\ip{\cdot}{J_{\omega}\cdot}_{\omega}$.
\end{itemize}
We need to examine the existence of a BRST charge. This will be done for a subset of states in $\fS_{\drb}$ for which we can construct a charge  as in Subsection \ref{sbsc:brch}. More general exitence criteria for the BRST charge will be examined in Section \ref{sec:ubbcsbrst} Theorem \eqref{th:abQ2}. 

The structers of Section \ref{sec:brstext} are constructed in the Fock representation using the fields $\phi(f)$ of the resolvents in $\cR(\fD,\sigma_2)$. So we start by investigating representations in which fields exist for all the resolvents, \ie the regular representations (cf. Definition \eqref{df:regrep}). We recall some useful facts and define terminology for these representations:
\begin{itemize}

\item Let the regular states of $\cR(\cD, \sigma_2)$ be denoted by $\Reg$ and let the strongly regular states be denoted $\SReg$ (cf. Definition \eqref{df:regrep} and Section \ref{StatesRep}). Recall that for $\omega_1\in  \SReg $ we have that $\phi_{\pi_{\omega_1}}(f)\in \op (\cH_{\omega_1})$ exists for all $f\in \fD$. In particular we can form the non-normed $*$-algebra
\[
\cE_{\omega_1}^{0}:= \salg{\phi_{\pi_{\omega_1}}(f), \, \pi_{\omega_1}(R(\lambda, f)) \, | \, f\in \fD, \, \lambda\in \mathbb{R} \backslash 0}
\]
\item For $\omega_1 \in \SReg$ we have that there exists a dense invariant domain $\cD_{\infty}$ for all $\phi_{\pi_S}(f)$, $f\in \fD$. Moreover the cyclic GNS-vector $\Omega_{\omega_1} \in \cD_{\omega_1}^0$, and so 
\[
\salg{\phi_{\pi_{\omega_1}}(f),| \, f\in \fD, \, \lambda\in \mathbb{R} \backslash 0} \Omega_{\omega_1} \subset \cD_{\infty}
\]
as $\cD_{\infty}$ is invariant for all  $\phi_{\pi_S}(f)$, $f\in \fD$. By Theorem \eqref{RegThm} (i) we get 
\begin{equation}\label{eq:dom1}
\cD_{\omega_1}^{0}:=\cE_{\omega_1}^{0}\Omega_{\omega_1}
\end{equation}
is an invariant domain for $\cE_{\omega_1}^{0}$, and as $\cD_{\omega_1}^{0} \supset \pi_{\omega_1}(\cR_0(\fD,\sigma_2))\Omega_{\omega_1}$ we get that $\cD_{\omega_1}^{0}$ is dense in $\cH_{\omega_1}$.   
\item For strongly regular states we can define the creators and annihilators of the fields. Let $\omega_1 \in \SReg(\cR(\fD, \sigma_2))$,then we define:
\[
a_{\pi_{\omega_1}}(f):=\fst(\phi_{\pi_{\omega_1}}(f)+i\phi_{\pi_{\omega_1}}(Kf)), \qquad f\in \fD
\]
where by Theorem \eqref{RegThm} (v) we have that $\cD^0_{\omega_1}$ is a dense invariant domain for $a_{\pi_{\omega_1}}(f)$. We have that
\[
a_{\pi_{\omega_1}}^*(f):=\fst(\phi_{\pi_{\omega_1}}(f)-i\phi_{\pi_{\omega_1}}(Kf)), \qquad f\in \fD
\]
which preserves $\cD^0_{\omega_1}$. For $f,g \in \fD$, 
\[
[a_{\pi_{\omega_1}}(f)\,,\,a_{\pi_{\omega_1}}(g)]\psi=[a_{\pi_{\omega_1}}^*(f)\,,\,a_{\pi_{\omega_1}}^*(g)]\psi=0,  
\]
\[
[a_{\pi_{\omega_1}}(f)\,,\,a_{\pi_{\omega_1}}^*(g)]\psi= \ip{f}{g}\psi,
\]
where $\psi \in \cD^0_{\omega_1}$. For consistency in notation we will denote $c_{\pi_{\omega_2}}(f):=\pi_{\omega}(\one \otimes c(f))$ for $c(f) \in \cA_g$ and $\omega_2 \in \fS_g(\cA_g)$ (cf. Section \ref{sec:ghost} for definition of $c(f)$ ).
\item Take $\omega=\omega_1\otimes \omega_2 \in \fS_\drb$ such that $\omega_1 \in \SReg$. Let $\cD_{\omega_2}^{0}:=\salg{C(g)\,|\, g\in \fDL \oplus_{\mH} \fDJ}\Omega_{\omega_2}$. Then  
\begin{equation}\label{eq:dfDomega}
\cD_{\omega}:= \cD_{\omega_1}^{0} \otimes \cD_{\omega_2}^{0}
\end{equation}
is a common dense invariant domain for $\pi_{\omega_1}(R(\lambda, f)) \otimes \one$, $\phi_{\pi_{\omega_1}}(f)\otimes \one$, and \\
$\pi_{\omega_1}(R(\lambda, f)) \otimes \one$ (cf equation \eqref{eq:dom1}). 
\item Let $\cXs$ be a  subspace of $\fDL$ and let $\fD_s:= \cXs \oplus J \cXs=\cXs \oplus_{\mH} J \cXs$. Let
\[
\cE(\cXs):=\salg{\phi_{\pi_{\omega_1}}(f)\otimes \one,\, \ \pi_{\omega_1}(R(\lambda, f)) \otimes \one, \,\one \otimes \pi_{\omega_2}(C(g))|\,f\in \fD_t\oplus \fD_s, \,g\in \fD_s, \, \lambda\in \mathbb{R} \backslash 0 },
\]
\item Define $\cE:=\cE(\fDL)$ and extend the ghost grading on $\cA$ to $\cE$ in the obvious way. 
\end{itemize}
\begin{rem}\label{rm:ED2} 
Note that for $\cXs \subset \fDL$ and $E \in \cE(\cXs)$, a monomial of generating tensors of $\cE(\cXs)$, we can choose a monomial of $R(\lambda_j, f_j)$'s, $\lambda_j \in \mathbb{R} \backslash 0 $, such that  ${\pi_{\omega}}(R(\lambda_1, f_1)\ldots R(\lambda_n, f_n)\otimes \one)E \in {\pi_{\omega}}(D_2(\drb))$. We can also take the $\lambda_j$'s to be arbitrarily large. 

To see this we first take $\cXs \subset \fDL$, $\lambda \in  \mathbb{R} \backslash 0 $,  $f \in \fD_t\oplus \fD_s$, $g \in \cXs$, and  $h \in J \cXs$. We show the above statement holds for $(\phi_{\pi_{\omega_1}}(f)\otimes \one)$ and ${\pi_{\omega}}(\one \otimes C(g))$.  
 
We have that ${\pi_{\omega_1}}(R(\lambda, f))(\phi_{\pi_{\omega_1}}(f)\otimes \one)= i \lambda {\pi_{\omega}}( R(\lambda, f))- \one \in {\pi_{\omega}}(D_2(\drb))$. Also, for $g \in \cXs$, we have that $ {\pi_{\omega}}(R(\lambda, \mH P_{+}g)R(\lambda,K\mH P_{-}g)\otimes \one ){\pi_{\omega}}(\one \otimes C(g))= (1/\lambda) {\pi_{\omega}}(\zeta_1 (g / \lambda)) \in {\pi_{\omega}}(D_2(\drb))$. 

Now $\pi_{\omega_1}(R(\lambda, f))\otimes \one,{\pi_{\omega}}(\one \otimes C(h))$ are elements of ${\pi_{\omega}}(D_2(\drb))$. Moreover, we can always multiply $E$ by the appropriate monomial $M$ of resolvents corresponding $(\phi_{\pi_{\omega_1}}(f)\otimes \one)$ and ${\pi_{\omega}}(\one \otimes C(g))$ in $E$ so that $ME \in {\pi_{\omega}}(D_2(\drb))$, and using the fact that $[\cR(\fD,\sigma_2)\otimes \one, \one \otimes \cA_g]=0$ and Theorem \eqref{RegThm} (vii) to pair the resolvents in $M$ with their corresponding terms in $E$.
\end{rem}
By the same construction as in the last subsection we now define a superderivation on $\cE$:
\begin{theorem}\label{th:srsdfch}
Let $\omega=\omega_1\otimes \omega_2 \in \fS_\drb$ is such that $\omega_1 \in \SReg$.

Define a map $\drhb:\cE \to \cE$ by:
\begin{align*}
\drhb(\phi_{\pi_{\omega_1}}(f)\otimes \one)&= \one \otimes \pi_{\omega_2} (iC( K P_{\JJ}P_{+} f)+ C(P_{\JJ}P_{-}f)),  \\
\drhb(\pi_{\omega_2}(C(h)))&= \phi_{\pi_{\omega_1}}(\mH P_{\LL} P_{+} h)+i \phi_{\pi_{\omega_1}}(\mH KP_{\LL} P_{-}h),\\
\drhb(\pi_{\omega}(A))&= \pi_{\omega}(\drb(A))
\end{align*}
for $f \in \fD$, $h \in \fDL \oplus \fDJ$, $A \in D_2(\drb)$.
 
Let $\cXs$ be a finite dimensional subspace  of $\fDL$, let $(f_j)$ be a finite $\fH$-orthonormal basis for $\cXs$ and define
\begin{align*}
\bch_s:= \sum_{j} & \Big[ (\phi_{\pi_{\omega_1}}(P_{+} f_j)+i \phi_{\pi_{\omega_1}}(KP_{-}f_j))\otimes \pi_{\omega_2}(C( Jf_j)) \\
&+ (\phi_{\pi_{\omega_1}}(KP_{+} f_j)-i\phi_{\pi_{\omega_1}}(P_{-}f_j)) \otimes \pi_{\omega_2}(C(K Jf_j)) \Big],
\end{align*}
where $D(\bch_s)=\cD_{\omega}$. Define the superderivation
\[
\drhb_s:\cE(\cXs) \to \cE(\cXs) \qquad \text{by} \qquad \drhb_s(E):=\sbr{\Qs}{E}, 
\]
on $ \cD_{\omega}$, graded with respect to the $\Z^2$-grading coming from the ghosts.
Then:
\begin{itemize}
\item[(i)] There exists a superderivation $\drhb:\cE \to \cE$ with respect to the $\Z^2$-grading coming from the ghosts that agrees with the map $\drhb$ defined on the generating tensors of $\cE$ above. Moreover $\drhb(E)=\drhb_s(E)$ for $E \in \cE(\cXs)$.
\item[(ii)] $\bch_s$ is 2-nilpotent, closable, Krein symmetric with respect to $\iip{\cdot}{\cdot}_{\omega}=\ip{\cdot}{J_{\omega}\cdot}_{\omega}$ and 
\[
\Qs=\sum_{j=1}^{n}[a_{\pi_{\omega_1}}^*( Jf_j)\otimes c_{\pi_{\omega_2}}( Jf_j) + a_{\pi_{\omega_1}}( f_j)\otimes c_{\pi_{\omega_2}}^*( f_j)]
\]
\item[(iii)] $D(\bch_s)\subset D(\bch_s^*)$ and 
\begin{align*}
\bch_s^*:= \sum_{j} &( (\phi_{\pi_{\omega_1}}(P_{+} f_j)-i \phi_{\pi_{\omega_1}}(KP_{-}f_j))\otimes \pi_{\omega_2}(C( f_j)) \\ &+ (\phi_{\pi_{\omega_1}}(KP_{+} f_j)+i\phi_{\pi_{\omega_1}}(P_{-}f_j)) \otimes \pi_{\omega_2}(C(K f_j)) ),
\end{align*}
on $\cD_{\omega}$.
\item[(iv)] For all $E \in \cE$,
\[
\ip{\Omega_{\omega}}{\drhb(E) \Omega_{\omega}}=0
\]
\item[(v)] We have 
\[
(a_{\pi_{\omega_1}}(f)\otimes \one)\Omega_{\omega}=0, \qquad(\one \otimes c_{\pi_{\omega_2}}(f))\Omega_{\omega}=0
\]
for all $f \in (\fDL\oplus \fDJ)$, hence $\Omega_{\omega}\in \ker \Qs \cap \ker \Qs^*$.
\end{itemize}
\end{theorem}
\begin{proof}
(i): We need only to check that the map $\drhb:\cE\to \cE$ extends to a superderivation on polynomials on the generating tensors in $\cE$. This was done in Theorem \eqref{th:calgdrb2}. 

\pfit(ii) and (iii): That $\bch_s$ is 2-nilpotent follows as in lemma \eqref{lm:2nQs} substituting $\phi_{\pi_{\omega_1}}(P_{+} f_j)+i \phi_{\pi_{\omega_1}}(KP_{-}f_j)$ for $A(f)$. 

For the remaining statements we note that $\cD_{\omega}=\cD_{\omega_1}^{0} \otimes \cH_{\omega_2}$,  that $\phi_{\pi_{\omega_1}}(f)$ is symmetric on $\cD_{\omega_1}^{0}$ and  $C(f)^*=C(Jf)$ for $f\in (\fDL\oplus \fDJ)$. Using these and that $\bch_s$ is a finite sum, we easily get that $D(\bch_s)\subset D(\bch_s^*)$ and $\Qs^*$ is given by the expression in (iii) on $D(\bch_s)$. Now  $\cD_{\omega}=D(\bch_s)\subset D(\bch_s^*)$ and so $ D(\bch_s^*)$ is dense in $\cH_{\omega}$. Hence by \cite{ReSi1972v1} Theorem VIII.1 p253  we get that $\Qs$ is closable. 

For Krein symmetry we use that $E^{\dag}=J_{\omega} E^* J_{\omega}$ for $E \in \cE$ and so
\begin{align*}
\phi_{\pi_{\omega_1}}(P_{+} f_j)^{\dag}&=\phi_{\pi_{\omega_1}}(JP_{+} f_j)^{*}=\phi_{\pi_{\omega_1}}(P_{+} f_j),\\ \phi_{\pi_{\omega_1}}(P_{-}f_j))^{\dag}&=\phi_{\pi_{\omega_1}}(JP_{-}f_j))=-\phi_{\pi_{\omega_1}}(JP_{-}f_j)),\\
\pi_{\omega_2}(C(Jf_j))^{\dag}&=\pi_{\omega_2}(C(J^2f_j)^*)=\phi_{\omega_2}(C(Jf_j)).
\end{align*}
using these on the terms in $\Qs$ gives that $\Qs$ is Krein-symmetric.

The expression for $\Qs$ in terms of creators and annihilators follows from the proof of lemma \eqref{lm:2nQs} (ii), since the algebraic input is the same.

\pfit(iv): Take a monomial $E$ of generating tensors of $\cE$. Now by Remark \eqref{rm:ED2} we can choose an appropriate monomial $M=\pi_{\omega}(R(\lambda_1, f_1)\ldots R(\lambda_n, f_n) \otimes \one)$ such that $ME \in \pi_{\omega}(D_2(\drb))$ , \ie $ME= \pi_{\omega}(A)$ for some $A \in D_2(\drb)$. Thus
\[
\ip{\Omega_{\omega}}{\drhb(ME)\Omega_{\omega}}=\omega(\drb(A))=0
\]
by the defining property of $\omega \in \fS_{\drb}$. As $ME$ is a monomial in a finite number of generating elements, there exists a finite dimensional subspace $\cXs$ of $\fDL$ such that $ME \in \cE(\cXs)$. Hence
\begin{align}
0=\ip{\Omega_{\omega}}{\drhb(ME)\Omega_{\omega}}&=\ip{\Omega_{\omega}}{\Qs ME\Omega_{\omega}}-\ip{\Omega_{\omega}}{ \gamma(ME)\Qs\Omega_{\omega}}, \notag\\
&=\ip{\Qs^*\Omega_{\omega}}{ ME\Omega_{\omega}}-\ip{\Omega_{\omega}}{ M\gamma(E)\Qs\Omega_{\omega}}\label{eq:Edom},
\end{align}
since $M$ is even so $\gamma(M)=M$ and by part (iii) that $\Omega_{\omega} \in \cD_{\omega} \subset D(\bch_s^*)$.
 
Now  $M$ is a function of the $\lambda_j$'s, and by Remark \eqref{rm:ED2} we have that $M(\lambda_1, \ldots, \lambda_n)E \in \pi_{\omega}(D_2(\drb))$, for all $\lambda_j \in \mathbb{R} \backslash 0$, and so equation \eqref{eq:Edom} holds for all $M(\lambda_1, \ldots, \lambda_n)$, $\lambda_j \in \mathbb{R} \backslash 0$.
 
By Theorem \eqref{RegThm} (ii) we have that, for all $\psi \in \cH_{\omega}$:
\[
\lim_{\lambda_1 \to \infty}\ldots \lim_{\lambda_n \to \infty} ((i)^n (\lambda_1 \ldots \lambda_n) \pi_{\omega}(M(\lambda_1, \ldots, \lambda_n) \psi)= \psi
\]
Take the limits  $\lambda_j \to \infty$ of $ (i)^n (\lambda_1 \ldots \lambda_n) \times (\text{equation }  \eqref{eq:Edom})$ to get that,
\begin{align*}
0=\ip{\Qs^*\Omega_{\omega}}{ E\Omega_{\omega}}-\ip{\Omega_{\omega}}{ \gamma(E)\Qs\Omega_{\omega}}=\ip{\Omega_{\omega}}{ \Qs E\Omega_{\omega}}-\ip{\Omega_{\omega}}{ \gamma(E)\Qs\Omega_{\omega}}=\ip{\Omega_{\omega}}{\drhb(E)\Omega_{\omega}}
\end{align*}
As any $E_2 \in \cE$ is a sum of monomials, this proves (iv).

\pfit(v): Take  $0 \neq f\in \fDL$, and let $\cXs= \{ \mathbb{C} f \}$. Therefore $\fD_s= [ f, Jf]$ and so as $J=P_{+}-P_{-}$ we have $P_{+}f, P_{-}f, KP_{+}f, KP_{-}f \in \fD_s$. Hence if we let
\begin{align*}
E= & \;(\phi_{\pi_{\omega_1}}(P_{+} f)-i \phi_{\pi_{\omega_1}}(KP_{-}f))\otimes \pi_{\omega_2}(C( f)) + (\phi_{\pi_{\omega_1}}(KP_{+} f)+i\phi_{\pi_{\omega_1}}(P_{-}f)) \otimes \pi_{\omega_2}(C(K f)),\\
&=a_{\pi_{\omega_1}}( Jf)\otimes c_{\pi_{\omega_2}}^*( Jf) + a_{\pi_{\omega_1}}^*( f)\otimes c_{\pi_{\omega_2}}( f),\\
=&\; \Qs^*
\end{align*}
then we have that $E \in \cE(\cXs)$. 

As $E \in \cE(\cXs)$, we have that $\drhb(E)=\drhb_s(E)$ and so following a calculation similar to equation \eqref{eq:kerD} we get
\begin{align*}
\drhb(E)=[\Qs, E\} = & \;( a_{\pi_{\omega_1}}^*( J f)a_{\pi_{\omega_1}}( J f)\{c_{\pi_{\omega_2}}( J f), c_{\pi_{\omega_2}}^*(  J f) \} + [a_{\pi_{\omega_1}}( J f), a_{\pi_{\omega_1}}^*(  J f)] c_{\pi_{\omega_2}}^*( J f) c_{\pi_{\omega_2}}( J f)) \notag \\
&+( a_{\pi_{\omega_1}}^*(  f)a_{\pi_{\omega_1}}(  f)\{c_{\pi_{\omega_2}}(  f), c_{\pi_{\omega_2}}^*(   f) \} + [a_{\pi_{\omega_1}}(  f), a_{\pi_{\omega_1}}^*(   f)] c_{\pi_{\omega_2}}^*(  f) c_{\pi_{\omega_2}}( J f)) \notag \\
= & \;( \norm{f}^2_{\fL}[a_{\pi_{\omega_1}}^*( J f)a_{\pi_{\omega_1}}( J f)+ a_{\pi_{\omega_1}}^*( f)a_{\pi_{\omega_1}}( f)] + [c_{\pi_{\omega_2}}^*(J f) c_{\pi_{\omega_2}}(J  f) + c_{\pi_{\omega_2}}^*( f) c_{\pi_{\omega_2}}( f)] )
\end{align*}
Now by part (iv) we have that
\begin{align*}
0= & \;\ip{\Omega_{\omega}}{\drhb(E)\Omega_{\omega}},\\
= & \;\ip{\Omega_{\omega}}{( \norm{f}^2_{\fL}[a_{\pi_{\omega_1}}^*( J f)a_{\pi_{\omega_1}}( J f)+ a_{\pi_{\omega_1}}^*( f)a_{\pi_{\omega_1}}( f)] + [c_{\pi_{\omega_2}}^*(J f) c_{\pi_{\omega_2}}(J  f) + c_{\pi_{\omega_2}}^*( f) c_{\pi_{\omega_2}}( f)]) \Omega_{\omega}}.
\end{align*}
As $f\in \fDL$ was arbitrary, and the RHS of $\drb(E)$ above is the sum of positive operators, we get that,
\[
(a_{\pi_{\omega_1}}(g)\otimes \one)\Omega_{\omega}=0, \qquad(\one \otimes c_{\pi_{\omega_2}}(g))\Omega_{\omega}=0,
\]
for all $g \in (\fDL \oplus \fDJ)$.  Hence, $\Omega_{\omega}\in \ker \Qs \cap \ker \Qs^*$ by the expression for $\Qs$ in (ii).
\end{proof}
Now Theorem \eqref{th:srsdfch} (v) says that strongly regular BRST physical states are Fock states when restricted to $\cR(\fDL \oplus \fDJ, \sigma_2)$. Using this we can construct a BRST charge $Q$ that generates $\drhb$ for these representations. We follow Subsection \ref{sbsc:brch}, in particular we have an almost identical version of lemma \eqref{lm:Qsind}.
\begin{lemma}\label{lm:Qsstuff}
Let $\omega=\omega_1\otimes \omega_2 \in \fS_\drb$ and  $\omega_1\in \fS_{sr}(\cR(\cD, \sigma_2))$. Let $\cXs \subset \cXsp$ be finite dimensional subspaces of $\fDL$. Then we have that,
\[
\bch_s \psi = \bch_{s'} \psi
\]
for $\psi \in \cE(\cXs) \Omega_{\omega}$.
\end{lemma}
\begin{proof}
Take $\psi \in \cE(\cXs) \Omega_{\omega}$, and suppose that $\dim(\cXs)=m< \dim(\cXsp)=n$. Now by lemma \eqref{lm:thbas}, take an $\fH$-orthonormal basis $\vL=(f_j)_{j=1}^{m}$ of $\cX_{s}$ such that it is also a $\fL$-orthogonal basis for $\cX_{s}$, and take an  $\fH$-orthonormal basis $\vL''=(f_j)_{j=m+1}^{n}$ of $\cXsp\ominus \cXs$ such that it is also a $\fL$-orthogonal basis for $\cXsp\ominus \cXs$.  Therefore $\vL'=(f_j)_{j=1}^{n}$ is an $\cH$-orthonormal basis and $\fL$-orthogonal basis for  $\cXsp$. Now,
\begin{align}\label{eq:Qbind2}
(\Qs-\Qsp )&=\sum_{j=m+1}^{n}[a^*( Jf_j)\otimes c( Jf_j) + a( f_j)\otimes c^*( f_j)].
\end{align}
By the way we chose $\vL'$ we have that $f_i \perp_{\fH} f_j$ and $f_i \perp_{\fL} f_j$ for $i\leq n$, $j>n+1$. So as $\psi=A \Omega_{\omega}$ where $A \in  \cE(\cXs)$, we see that we can (anti-)commute all the terms in the RHS of \eqref{eq:Qbind2} through the terms in $A$ to $\Omega_{\omega}$, which they annihilate by Theorem \eqref{th:srsdfch} (v). Therefore
\[
(\Qs-\Qsp )\psi=0
\]
\end{proof}
Now every $\psi \in \cD_{\omega}$ is also in $\cE(\cXs) \Omega_{\omega}$ for some finite dimensional subspace $\cXs \subset \fDL$ (cf. equation \eqref{eq:dfDomega}), and so using the above lemma we construct the BRST charge:
\begin{theorem} \label{th:Q}
Let $\omega=\omega_1\otimes \omega_2 \in \fS_\drb$ and  $\omega_1\in \fS_{sr}(\cR(\cD, \sigma_2))$. Let $\psi \in \cD_{\omega}$, let $\cXs$ be a finite dimensional subspace of $\fDL$ such that $\psi \in \cE(\cXs) \Omega_{\omega}$. Define 
\[
Q\psi:= \Qs \psi. 
\]
Then $Q$ extends to well defined operator on with domain $D(Q):=\cD_{\omega}$ such that
\[
\drhb(E)\psi=\sbr{Q}{E}\psi, 
\]
for all $\psi \in \cD_{\omega}$ and $E\in \cE$. Furthermore:
\begin{itemize}
\item[(i)] $Q$ preserves $D(Q)$, is Krein symmetric (hence closable by Proposition \eqref{pr:Krclos}) and $Q^2\psi=0$ for all $\psi \in D(Q)$. 
\item[(ii)] $\cQ$ preserves $D(\cQ)$ and $\cQ^2\psi=0$ for all $\psi \in D(\cQ)$. Hence $\cQ$ has \emph{dsp}-decomposition (cf. Theorem \eqref{th:Hdsp1})
\[
\cH_{\omega}^d\oplus \cH_{\omega}^s \oplus \cH_{\omega}^p,
\]
where $\cH_{\omega}^d= \overline {\ran \cQ}$, $\cH_{\omega}^s = \ker \cQ \cap \ker Q^* $, $\cH_{\omega}^p= \overline{\ran Q^*}$.
\item[(iii)] $(Q^*)^2\psi=0$ for all $\psi \in D(Q^*)$.
\end{itemize}
\end{theorem}
\begin{proof}
Let $\psi \in \cD_{\omega}$, then there exists a finite dimensional subspace $\cXs \subset \fDL$ such that $\psi \in \cE(\cXs) \Omega_{\omega}$. We check that $Q\psi=\Qs \psi$ is independent of the choice of this $\cX_s$. Let $\cXsp\subset \fDL$ be another finite dimensional subspace such that $\psi \in \cE(\cXsp) \Omega_{\omega}$ and suppose that $\cXs,\cXsp\subset\cXspp $ where $\cXspp$ is a finite dimensional subspace of $\fDL$. Then by lemma \eqref{lm:Qsstuff}
\[
\Qs\psi=\Qspp \psi=\Qsp \psi,
\]
and so $Q\psi=\Qs \psi$ is independent of the choice of this $\cX_s$. That $Q$ extends to a well defined linear operator on $D(Q)=\cD_{\omega}$ is now obvious. 

Let $E \in \cE$ and $\psi \in \cD_{\omega}$ and take a finite dimensional subspace $\cXs \subset \fDL$ such $E \in \cE(\cXs)$ and $\psi \in \cE(\cXs)\Omega_{\omega}$. Then $E\psi\in \cE(\cXs)$, and so $QE\psi=\Qs E\psi \in \cE(\cXs)\Omega_{\omega}$ and $EQ\psi =E\Qs \psi \in \cE(\cXs)\Omega_{\omega}$. Therefore by Theorem \eqref{th:srsdfch} (i) $\drhb(E)\psi=\sbr{\Qs}{E}\psi=\sbr{Q}{E}\psi$.

\pfit (i), (ii) and (iii): Let $\psi \in D(Q)$. Then $\psi \in \cE(\cXs) \Omega_{\omega}$ for some finite dimensional subspace of $\cXs \subset \fDL$.  So we have that $Q\psi=\Qs\psi \in \cE(\cXs) \Omega_{\omega}$, hence $Q D(Q)\subset D(Q)$. By the definition of $\Qs$ we have that $\Qs \cE(\cXs) \Omega_{\omega} \subset \cE(\cXs) \Omega_{\omega}$ and so we have that $Q^2\psi=\Qs^2\psi=0$. Furthermore $Q$ is $\dag$-symmetric as $\Qsp $ is for all finite dimensional subspaces $\cXsp\subset \fDL$. 

Now (i), (ii) and (iii) follow from lemma \eqref{lm:QessaCsa} and Theorem \eqref{th:Hdsp1}. 
\end{proof}
We want to calculate $\cH_{\omega}^s=\ker \cQ \cap \ker Q^*$ explicitly.
\begin{proposition}\label{pr:kersrRph}
Let $\omega=\omega_1\otimes \omega_2 \in \fS_\drb$ and  $\omega_1\in \fS_{sr}(\cR(\cD, \sigma_2))$, and let $Q$ and $\cH_{\omega}^s=\ker {\cQ} \cap \ker {\bch^*}$ be as in Theorem \eqref{th:Q} above. Let $\psi \in \cD_{\omega}$, then: 
\begin{itemize}
\item[(i)]$\psi \in \cH_{\omega}^s$ if and only if,
\[
(a_{\pi_{\omega_1}}(g)\otimes \one)\psi=0, \qquad(\one \otimes c_{\pi_{\omega_2}}(g))\psi=0,
\]
for all $g \in (\fDL \oplus \fDJ)$, hence $\Omega_{\omega}\in \cH_{\omega}^s$. This implies that for $\psi \in \cH_{\omega}^s$, the vector state $\omega_{\psi}$ is a Fock state when restricted to $\cR(\fDL \oplus \fDJ, \sigma_2)\otimes \one$.
\item[(ii)] Let  
\begin{gather*}
\cR^0_{u}:=\salg{ R(\lambda,f) \,|\, f \in \fDL \oplus \fDJ, \,\lambda \in \R \backslash \{0\}} \qquad \cR_{u}:=\overline{\cR^0_{u}}=\cR(\fDL\oplus \fDJ, \sigma_2)\\
\cF_{u}:=\salg{ \phi_{\pi_{\omega_1}}(f)\otimes \one \,|\, f \in \fDL\oplus \fDJ}.
\end{gather*}
Then
\[
\overline{(\pi_{\omega_1}(\cR_u)\otimes \one )\Omega_{\omega}}=\overline{\cF_{u} \Omega_{\omega}}
\]
\item[(iii)] Let 
\[
\cR^0_{ph}:=\salg{ R(\lambda,f) \,|\, f \in \fD_t, \,\lambda \in \R \backslash \{0\}} \qquad \cR_{ph}:=\overline{\cR^0_{t}}=\cR(\fD_t, \sigma_2).
\]
Then
\begin{align}
\cH_{\omega}=\overline{\pi_{\omega}(\cR(\fD,\sigma_2)\otimes \cA_g)\Omega_{\omega}}=&\;\overline{\pi_{\omega}(\cR_{ph}\cR_{u}\otimes\cA_g)\Omega_{\omega}},\notag \\
=&\; \overline{\pi_{\omega}(\cR_{ph}\otimes \one)\cF_{u}\pi_{\omega}(\one \otimes \cA_g)\Omega_{\omega}} \label{eq:Homstuff}
\end{align}
\item[(iv)] We have
 \begin{align}\label{eq:kerCQ}
\cH_{\omega}^s= \overline{\pi_{\omega}(\cR_{ph})\Omega_{\omega}}
\end{align}

\end{itemize}
\end{proposition}
\begin{proof}
(i): Take $\psi \in \cD_{\omega}$. Then there exists a finite dimensional subspace $\cXs \subset \fDL$ such that $\psi \in \cE(\cXs) \Omega_{\omega}$. Now as $\psi \in \cD_{\omega}$, we have $\cQ \psi= \bch_s \psi$ and $\bch^* \psi= \bch_s^* \psi$, and so $\psi \in (\ker \cQ \cap \ker \bch^*)$ if and only if $\psi \in \ker \{ \Qs, \Qs^* \}$. 

Using lemma \eqref{lm:thbas} we can choose a finite $\fH$-orthonormal basis, $(f_j)$, which is also $\fL$-orthogonal. So we  calculate as in equation \eqref{eq:kerD},
\begin{align*}
\{ \Qs, \Qs^* \} \psi =\left(\sum_{j=1}^{m}( \norm{f_j}^2_{\fL}[a^*( J f_j)a( J f_j)+ a^*( f_j)a( f_j)] + [c^*(J f_j) c(J  f_j) + c^*( f_j) c( f_j)] ) \right) \psi,
\end{align*}
Hence, as the RHS of the above equation is a sum of positive operators acting on $\psi$ we have that $\psi \in (\ker \cQ \cap \ker \bch^*)$ if and only if
\[
\psi \in  \{ \ker a_{\pi_{\omega_1}}(g)\otimes \one \cap \ker \one \otimes c_{\pi_{\omega_2}}(g)\,|\, g \in \cXs \oplus J\cXs\}.
\]
As $\psi \in \cE(\cXs) \Omega_{\omega}$, the (anti)commutation relations for the $a_{\pi_{\omega_1}}(\cdot)$'s and $c_{\pi_{\omega_2}}(\cdot)$'s give that this statement extends to $\psi \in (\ker \cQ \cap \ker \bch^*)$ if and only if
\[
\psi \in \{\ker a_{\pi_{\omega_1}}(g)\otimes \one \cap \ker \one \otimes c_{\pi_{\omega_2}}(g)\,|\, g\in \fDL \oplus J\fDL\}.
\]
By Theorem \eqref{th:srsdfch} (v) we have that $\Omega_{\omega} \in \cH_{\omega}^s$.

\pfit (ii): By (i) we have that $\omega$ is a Fock-state when restricted to $\cR_u=\cR(\fDL\oplus \fDJ)$ and so the result follows.

\pfit (iii): As $\omega_1$ is strongly regular and $\fD=\fD_t\oplus \fDL \oplus \fDJ$, Theorem \eqref{RegThm} (v) gives that 
$\overline{\pi_{\omega_1}(\cR(\fD,\sigma_2))\Omega_{\omega_1}}=\overline{\pi_{\omega_1}(\cR_{ph})\pi_{\omega_1}(\cR_{u})\Omega_{\omega_1}}$. 
Hence
\[
\cH_{\omega}=\overline{\pi_{\omega}(\cR(\fD,\sigma_2)\otimes \cA_g)\Omega_{\omega}}=\overline{\pi_{\omega}(\cR_{ph}\cR_{u}\otimes\cA_g)\Omega_{\omega}}= \overline{\pi_{\omega}(\cR_{ph}\otimes \one)\cF_{u}\pi_{\omega}(\one \otimes \cA_g)\Omega_{\omega}}
\]
where the last equality follows from (ii) and $[\cF_{u},\pi_{\omega}(\cR_{ph}\otimes \one)]=[\cF_{u},\pi_{\omega}(\one \otimes \cA_g)]=0$ on $\cD_{\omega}$.

\pfit (iv): Define 
\[
\cC_{u}:=\alg{a^*_{\pi_{\omega_1}}(f)\otimes \one, \one \otimes c^*_{\pi_{\omega_2}}(g)\,|\, f,g \in \fDL \oplus  \fDJ},
\]
acting on $\cD_{\omega}$. Equation \eqref{eq:Homstuff} and (i) give:
\begin{equation}\label{eq:Hwalstuff2}
\cH_{\omega}=[ A\Omega_{\omega},\, B\Omega_{\omega}\,|\, A \in \cR^{0}_{ph},\, B\in \cR^{0}_{ph}\cC_u].
\end{equation}
Let $\psi =A \Omega_{\omega}$ and $\xi=B\Omega_{\omega}$ where $A \in \cR^{0}_{ph}$ and $B\in \cR^{0}_{ph}\cC_u$. By definition, $B=TS$ where $T\in \cR^{0}_{ph}$ and $S \in \cC_u$. As $\fD_t \perp \fDL \oplus \fDJ$ with respect to $\ip{\cdot}{\cdot}$ on $\fD$ we have that the CCR's for the $\phi_{\pi_{\omega_1}}(f)$'s give that $[S,T]=[A,S^*]=0$ on $\cD_{\omega}$. Hence
\[
\ip{\psi}{\xi}_{\omega}=\ip{A\Omega_{\omega}}{TS\Omega_{\omega}}_{\omega}=\ip{AS^*\Omega_{\omega}}{T\Omega_{\omega}}_{\omega}=0
\]
where we used that $S^*$ is a polynomial of $a_{\pi_{\omega_1}}(f_j)\otimes \one$'s and $\one \otimes c_{\pi_{\omega_2}}(g_k)$'s where $ f_j,g_k \in \fDL \oplus  \fDJ$ and so $S^*\Omega_{\omega}=0$.

Therefore we have that $\cR^{0}_{ph}\Omega_{\omega}\perp \cR^{0}_{ph}\cC_u\Omega_{\omega}$ with respect to $\ip{\cdot}{\cdot}_{\omega}$ and so by equation \eqref{eq:Hwalstuff2} above we get that
\begin{equation}\label{eq:crdestblah}
\cH_{\omega}=\overline{\cR^{0}_{ph}\Omega_{\omega}}\oplus \overline{\cR^{0}_{ph}\cC_u\Omega_{\omega}} 
\end{equation}
By (i) $\cH_{\omega}^s=\cap \ker \{a_{\pi_{\omega_1}}(g)\otimes \one,\, \one \otimes c_{\pi_{\omega_2}}(g)\,|\, g\in \fDL \oplus \fDJ\}$ and $\Omega_{\omega}\in \cH_{\omega}^s$. Using the CCR's for the $a_{\pi_{\omega_1}}(g)$'s and the CAR's for the $c_{\pi_{\omega_2}}(g)$'s for $g \in \fDL \oplus \fDJ$ , we see that 
\[
\cF_{ph}\Omega_{\omega} \subset \cH_{\omega}^s, \qquad \text{and} \qquad (\cF_{ph}\cC_u\Omega_{\omega} \cap \cH_{\omega}^s)=\{0\},
\] 
where we recall that (i) implies that $\Omega_{\omega}$ is the Fock-vacuum for $\salg{\phi_{\pi_{\omega_1}}(f)\otimes \one\,\one \otimes \pi_{\omega_2}(C(f))|\, f \in \fDL \oplus \fDJ}$.

Hence by equation \eqref{eq:crdestblah} above we get that 
\[
\cH_{\omega}^s= \overline{\cR^0_{ph}\Omega_{\omega}}=\overline{\pi_{\omega}(\cR_{ph}) \Omega_{\omega}}.
\]

\end{proof}

\subsection{BRST-Algebra for QEM II}\label{sbs:BRSTIIalg}
Now that we have calculated $\cH_{\omega}^s$ for strongly regular states via Proposition \eqref{pr:kersrRph}, we want to define and investigate what the BRST-physical algebra will look like via constructions as in Section \ref{sbs:AbcuQ}. 
 
To use the results from Section \ref{sbs:AbcuQ} we first have to define in a given representation the superderivation we are using.  To make streamline notation, define $\cR_{\tp}:=\cR(\fD_t, \sigma_2)$ and $\cR_{\tilde{u}}:=\cR(\fDL \oplus \fDJ, \sigma_2)$. Note that $\fD_t$ and $\fDL \oplus \fDJ$ are $\sigma_2$-symplectic complements of each other. So by Theorem \eqref{TensorAlg}: 
\[
\cR_t:=C^*(\{ R(\lambda, f),\,R(\lambda, g)\,| f \in \fD_t,\, g \in \fDL \oplus \fDL \} )\cong \cR_{\tp}\otimes \cR_{\tilde{u}}.
 \]
Let $\vp$ be this $*$-isomorphism, and define $\cR_{ph}:= \vp^{-1}(\cR_{\tp} \otimes \one)$ and $\cR_{u}:= \vp^{-1}(\one \otimes \cR_{\tilde{u}})$. Recall that for a subspace $S\subset \fD$, we have that $\cR_0(S,\sigma_2)=\salg{ R(\lambda,f)\,|\, f \in S,\, \lambda \in \mathbb{R}\backslash \{0\}}$, and define:
\begin{gather}
\cR^0_{ph}:=\vp^{-1}(\cR_0(\fD_t,\sigma_2)), \qquad \cR^0_{u}:=\vp^{-1}(\cR_0(\fDL\oplus \fDJ,\sigma_2)), \label{eq:ders}\\
\cR^0_t:=\salg{ R(\lambda, f),\,R(\lambda, g)\,| f \in \fD_t,\, g \in \fDL \oplus \fDL  )} \notag
\end{gather}

\begin{definition}\label{df:BRS2domdel}
Define: $D(\drb_{\omega})=\pi_{\omega}(D_2(\drb))\cap (\cR^0_t\otimes \cA_g)) $
\begin{align*}
\drb_{\omega}&:D(\drb_{\omega}) \to D(\drb_{\omega}), \\ 
 \drb_{\omega}&(\pi_{\omega}(A))=\pi_{\omega}(\drb(A)),
\end{align*}
where $\drb:D_2(\drb)\to D_2(\drb)$ is as in Theorem \eqref{th:calgdrb2}.
\end{definition}
\begin{rem} The reason we do not define $D(\drb_{{\omega}})=\pi_{\omega}(D_2(\drb))$ is that the domain restriction makes calculations tractable. We do not lose information using the smaller domain and strongly regular states since for $\omega=\omega_1\otimes \omega_2 \in \fS_\drb$ and $\omega_1\in \SReg$, we have that $D(\drb_{{\omega}})=\pi_{\omega}(D_2(\drb)\cap (\cR^0_t\otimes \cA_g))$ is strongly dense in $\pi_{\omega}(\cA)$ by Theorem \eqref{RegThm} (v).
\end{rem}

From Theorem \eqref{th:Q} we have that if $\omega=\omega_1\otimes \omega_2 \in \fS_\drb$ where $\omega_1\in \SReg$, then $\drb_{\omega}$ has a generator $Q_{\omega}$, which has \emph{dsp}-decomposition $\cH_{\omega}=\cH^d_{\omega}\oplus \cH^s_{\omega}\oplus \cH^p_{\omega}$. Define
\begin{align*}
\Phi^{\omega}_s:D(\drb_{\omega})&\to D(\drb_{\omega})\qquad \text{by} \qquad \Phi^{\omega}_s(A)=P_s^{\omega}A P_s^{\omega},
\end{align*}
where $P^{\omega}_s$ is the projection onto $\cH^s_{\omega}$.
\begin{definition}\label{df:PQ2}
Let $\omega=\omega_1\otimes \omega_2 \in \fS_\drb$ and $\omega_1\in \SReg$. Let $Q_{\omega}$ be the generator of $\drb_{\omega}$ given in Theorem \eqref{th:Q}, and let,
\[
\cH_{\omega}=\cH^d_{\omega}\oplus \cH^s_{\omega}\oplus \cH^p_{\omega}
\]
be the \emph{dsp}-decomposition with respect to $\cQ_{\omega}$. We define
\[
\tilde{\cP}_{\omega}^{BRST}:=\ker \drb_{\omega}/(  \ker \drb_{\omega} \cap \ker \Phi^{\omega}_s) \cong \Phi^{\omega}_s(\ker \drb_{\omega}),
\]
where the above isomorphism is assumed to be algebraic, and we have assumed no topology on $\tilde{\cP}_{\omega}^{BRST}$.
\end{definition}
\begin{rem}
\begin{itemize} \item[(i)] By lemma \eqref{lm:alkerhom} we see that  $\Phi^{\omega}_s$ is an algebra homomorphism, hence $\tilde{\cP}_{\omega}^{BRST}$ is well defined as an algebra.
\item[(ii)] The above follows the alternative definition of the physical algebra (Subsection \eqref{sbs:alalg}). To connect with the $\ker \drb_{\omega}/ \ran \drb_{\omega}$ we need to extend $D(\drb_{{\omega}})$ to\\ $D(\drb_{{\omega}})^{ext}=\alg{D(\drb_{{\omega}}), P^{\omega}_s, P^{\omega}_p,P^{\omega}_d, Q, V}$ as in Section \eqref{sbs:AbcuQ}, which then gives  $ (\ker \drb_{\omega} \cap D(\drb_{{\omega}}))/ (\ran \drb_{\omega} \cap D(\drb_{{\omega}})) \cong \cP_{\omega}^{BRST}$ by Theorem \eqref{pr:krdel} (see Remark \eqref{rm:alus}).
\end{itemize}
\end{rem}

To calculate $\tilde{\cP_{\omega}}^{BRST}$, we first need Proposition \eqref{pr:kersrRph}  to show that BRST strongly regular states have a particularly simple form. 
\begin{lemma}\label{lm:fock}
\begin{itemize}
\item[(i)] Let $\omega \in \SReg(\cR(\fD, \sigma_2))$ be such that
\[
a_{\pi_{\omega}}(f)\Omega_{\omega}=0, \qquad \forall f \in (\fDL \oplus \fDJ),
\]
where $\Omega_{\omega}$ is the cyclic state in the GNS representation of $\omega$. Then we have that 
\[
\omega(AB)=\omega(A)\omega(B),
\]
for all $A \in \cR_{ph}$ and $B \in \cR_{u}$.
\item[(ii)] Let $\omega=\omega_1 \otimes \omega_2 \in \fS_\drb$ and $\omega_1 \in \SReg(\cR(\cD, \sigma_2))$. Then we have that for all $T \in (\cR^0_t\otimes\cA_g)$ there exists $S \in \cR_{ph}$ such that,
\[
\omega( R_1(T-S) R_2)=0, \qquad \forall R_1, R_2 \in \cR_{ph}
\]
\end{itemize}
\end{lemma}
\begin{proof}
(i): Let $\pi_{\omega}$ be the GNS representation for $\omega$. Define
\begin{align*}
\cG_{u}:=\salg{ \phi_{\pi_\omega}(f)\,|\, f \in (\fDL \oplus \fDJ) }
=\alg{a_{\pi_{\omega}}(f),\, a_{\pi_{\omega}}^*(f) \,|\, f \in (\fDL \oplus \fDJ) }
\end{align*}
acting of dense domain $\cD_{\omega}^0 \ni \Omega_{\omega}$ constructed as in \eqref{eq:dom1}. By the assumption that $a_{\pi_{\omega}}(f)\Omega_{\omega}=0$ for all $f\in (\fDL \oplus \fDJ)$, we have that $\omega$ is a Fock-state when restricted to $\cR_u=\cR(\fDL\oplus \fDJ,\sigma_2)$ which implies that 
\[
\overline{\cF_{u}\Omega_{\omega}}=\overline{\pi_{\omega}(\cR_u)\Omega_{\omega}}.
\]
If we take any element $T\in \cG_u$ then we can write $T$ in a normal ordering of $a_{\pi_{\omega}}(f_i)\text{'s},\, a_{\pi_{\omega}}^*(f_j)\text{'s}$ for $f_i, f_j \in \fDL \oplus \fDJ$. That is 
\[
T= a_0 \one + T_{1}
\]
where $a_0\in \mathbb{C}$ and $T_{1}$ is a polynomial with zero constant coefficient and  with all the $a_{\pi_{\omega}}^*(\cdot)$'s to the left of the $a_{\pi_{\omega}}(\cdot)$'s. Now using $a_{\pi_{\omega}}(f)\Omega_{\omega}=0$  for all $ f \in (\fDL \oplus \fDJ)$, we get that
\[
T\Omega_{\omega}= a_0\Omega_{\omega}+ T_{2}\Omega_{\omega}
\]
where $T_{2}$ is a polynomial with zero constant coefficient but now only in $a_{\pi_{\omega}}^*(\cdot)$'s. Thus:
\[
\ip{\Omega_{\omega}}{T\Omega_{\omega}}=\ip{\Omega_{\omega}}{a_0\Omega_{\omega}}+\ip{T_{2}^*\Omega_{\omega}}{\Omega_{\omega}}= a_0
\]
Now let $A\in \cR_{ph}$, then
\begin{align}
\ip{\Omega_{\omega}}{\pi_{\omega}(A)T\Omega_{\omega}}= & \;\ip{\Omega_{\omega}}{\pi_{\omega}(A)a_0\Omega_{\omega}}+\ip{\Omega_{\omega}}{\pi_{\omega}(A)T_{2}\Omega_{\omega}},\notag \\
= & \;\ip{\Omega_{\omega}}{\pi_{\omega}(A)a_0\Omega_{\omega}}+\ip{T_{2}^*\Omega_{\omega}}{\pi_{\omega}(A)\Omega_{\omega}},\notag \\
= & \; \omega(A)a_0,\notag\\
= & \; \omega(A)\ip{\Omega_{\omega}}{T\Omega_{\omega}} \label{eq:tenu}
\end{align}
since $[S,R]=0$ for all $S \in \cG_{u}$ and $R \in \pi_{\omega}(\cR_{ph})$ in the second equality. We have shown above that $\cF_u\Omega_{\omega}$ is dense in $\pi_{\omega}(\cR_{u})\Omega_{\omega}$, and so for $B \in  \cR_{u}$ there exists a sequence $(T_n)\subset \cG_u$ such that $\pi_{\omega}(B)\Omega_{\omega}=\lim_{n}T_n \Omega_{\omega}$ and so,
\begin{align*}
\omega(AB)=\lim_{n}\ip{\Omega_{\omega}}{\pi_{\omega}(A)T_n\Omega_{\omega}}=\omega(A)\lim_{n}\ip{\Omega_{\omega}}{T_n\Omega_{\omega}}
= \omega(A)\omega(B)
\end{align*}
where we used \eqref{eq:tenu} in the second equality. 

\pfit (ii): By Proposition \eqref{pr:kersrRph} (i) we get that $a_{\pi_{\omega_1}}(f)\Omega_{\omega_1}=0$ for all $f \in \fDL \oplus \fDJ$. Then part (i) gives 
\begin{equation}\label{eq:omfock}
\omega(AB\otimes C)= \omega_1(A)\omega_1(B)\omega_2(C)=\omega(A\otimes \one)\omega(B\otimes C)
\end{equation}
for all $A \in \cR_{ph}$, $B \in \cR_{u}$ and $C \in \cA_g$.

Let $T=\sum_{i=1}^{n}  A_i B_i \in \cR_{t}\otimes \cA_g$, $A_i \in (\cR_{ph}\otimes \one)$ $B_i \in (\cR_{u}\otimes \cA_g)$. Let $b_i= \omega(B_i)$ and let $S=\sum_{i=1}^{n}  b_i A_i$. Then we have that for all $R_1, R_2 \in \cR_{ph}$,
\begin{align*}
\omega(R_1(T-S)R_2)= & \;\sum_{i=1}^{n}\omega(R_1A_iR_2(B_i-b_i)),\\
= & \;\sum_{i=1}^{n}\omega(R_1A_iR_2)\omega(B_i-b_i),\\
= & \;0,
\end{align*}
where we have used that $[F,G]=0$ for all $F \in \cR_{u}$ and $G \in \cR_{ph}$ in the first line and equation \eqref{eq:omfock} in the second. As all elements in $\cR^0_{t} \otimes \cA_g$ are of the form of $T$ we are done.
\end{proof}
Using this we can calculate the algebra $\cP^{BRST}_{\omega}$. Before we proceed, we define 
\[
\tilde{\Phi}_s^{\omega}:\pi_{\omega}(\cA)\to \pi_{\omega}(\cA) \qquad \text{by} \qquad \tilde{\Phi}_s^{\omega}(A):=P^{\omega}_sAP^{\omega}
\]
\ie we have extended the domain of $\Phi^{\omega}_s$ to all of $\pi_{\omega}(\cA)$. Note that $\tilde{\Phi}_s^{\omega}$ is a linear transform \emph{not} a homomorphism on $\pi_{\omega}(\cA)$. As $\cH^s_{\omega}= =\overline{\pi_{\omega}(\cR_{ph}\otimes \one)\Omega_{\omega}}$ we have that $\tilde{\Phi}_s^{\omega}$ is a $*$-homomorphism on $\pi_{\omega}(\cR_{ph}\otimes \one)$. 
\begin{theorem}\label{th:brstpalgra2}
Let $\omega=\omega_1 \otimes \omega_2 \in \fS_\drb$ and $\omega_1 \in \SReg(\cR(\cD, \sigma_2))$, and define $\omega_{1,\cR_{ph}}:=\omega_1|_{\cR_{ph}}$. Then:
\begin{itemize}
\item[(i)] $(\Phi^{\omega}_s \circ \pi_{\omega})(\cR^{0}_{t} \otimes \cA_g) =(\Phi^{\omega}_s\circ\pi_{\omega})(\cR^{0}_{ph}\otimes \one)$
\item[(ii)] There exists an isometric isomorphism
\[
\vartheta_{\omega}:(\tilde{\Phi}^{\omega}_s\circ\pi_{\omega})(\cR_{ph}\otimes \one) \to \pi_{(\omega_1,\cR_{ph})}(\cR_{ph}),
\]
where $\pi_{(\omega_1,\cR_{ph})}$ is the GNS-representation for $\omega_{1,\cR_{ph}}$.
\item[(iii)] $\cP^{BRST}_{\omega}\cong \pi_{(\omega_1,\cR_{ph})}(\cR^0_{ph})$ as algebras with no topology.
\end{itemize}
\end{theorem}
\begin{proof}
(i): let $T=\sum_{i=1}^{n} a_i \pi_{\omega}(A_i)\pi_{\omega}( B_i) \in (\ker \drhb \cap \pi_{\omega}(\cR^{0}_t \otimes \cA_g) )$, where  $A_i \in (\cR^{0}_{ph}\otimes \one)$,  $B_i \in (\cR^{0}_{u}\otimes \cA_g)$. By lemma \eqref{lm:fock} (ii) we have that there exists $S \in \pi_{\omega}(\cR^{0}_{ph} \otimes \one) $ such that $\omega(R_1(T-S)R_2)=0$ for all $R_1, R_2 \in \cR_{ph}$. By Proposition \eqref{pr:kersrRph} (ii) we have $\cH^{\omega}_s= \overline{ \pi_{\omega}(\cR_{ph})\Omega_{\omega}}$, hence $P^{\omega}_s \pi_{\omega}((T-S))P^{\omega}_s =0$ and so
\[
P^{\omega}_s\pi_{\omega}(T)P^{\omega}_s=\Phi^{\omega}_s(\pi_{\omega}(T))=\Phi^{\omega}_s(\pi_{\omega}(S))
\]

\pfit(ii):  Let $\iota_{\cR_{ph}}:\cR_{ph} \to \cR_{ph}\otimes \one$ be the identity isomorphism. As $\cH^s_{\omega}=\overline{\pi_{\omega}(\cR_{ph}\otimes \one)\Omega_{\omega}}$ (Proposition \eqref{pr:kersrRph} (ii)), the representation $(\Phi^{\omega}_s\circ\pi_{\omega}\circ \iota_{\cR_{ph}}):\cR^{0}_{ph} \to B(\cH^s_{\omega})$ has an algebraic cyclic vector $\Omega_{\omega}$ (we take closures below. 

Define $\eta:=\tilde{\Phi}^{\omega}_s\circ\pi_{\omega}\circ \iota_{\cR_{ph}}$, then for all $A \in \cR_{ph}$
\begin{align*}
\ip{\Omega_{\omega}}{\eta(A)\Omega_{\omega}}_{\omega}=\ip{P^s_{\omega}\Omega_{\omega}}{\pi_{\omega}(A\otimes\one)P^s_{\omega}\Omega_{\omega}}_{\omega},
=&\;\omega(A\otimes\one),\\
=&\;\omega_{1,\cR_{ph}}(A),\\
=&\;\ip{\Omega_{(\omega_1,\cR_{ph})}}{\pi_{(\omega_1,\cR_{ph})}(A)\Omega_{(\omega_1,\cR_{ph})}}
\end{align*}
where  we used $P^s_{\omega}\Omega_{\omega}=\Omega_{\omega}$ in the second equality (Proposition \eqref{pr:kersrRph} (i)). Hence by \cite[Theorem 4.1.4, p143]{Mur90} we get that $\eta(\cR_{ph})$ unitarily equivalent to $\pi_{(\omega_1,\cR_{ph})}(\cR_{ph})$ and so there exists an isomorphism 
\[
\vartheta:\eta(\cR_{ph})\to\pi_{(\omega_1,\cR_{ph})}(\cR_{ph}).
\]
Note that as $\eta(\cR_{ph})$ and $\pi_{(\omega_1,\cR_{ph})}(\cR_{ph})$ are $C^*$-algebras, this isomorphism is isometric \cite[Theorem 3.1.5, p80]{Mur90}

Now $\eta(\cR_{ph})=  (\Phi^{\omega}_s\circ\pi_{\omega})(\cR_{ph}\otimes \one)$ so in fact $\vartheta:(\Phi^{\omega}_s\circ\pi_{\omega})(\cR_{ph}\otimes \one) \to \pi_{(\omega_1,\cR_{ph})}(\cR_{ph})$ is an isometric isomorphism.

\pfit (ii): By definition $\pi_{\omega}(\cR^0_{ph})$ contains only polynomials in resolvents with test functions in $\fD_t$. By the definition $\drb$ (cf. Theorem \eqref{th:calgdrb2}) 
\[
\drb_{\omega}(\pi_{\omega}(R(\lambda,f)))=-(1/\sqrt{2}) \pi_{\omega}(R(\lambda, f)^2(c(P_{\JJ}f)-c^*(P_{\JJ}Jf))),
\] 
so as $P_{\JJ}f=P_{\JJ}Jf=0$ for $f \in \fD_t$, we get that $\pi_{\omega}(\cR^0_{ph}) \in \ker \drb_{\omega}$. From this and part (i) we have,
\begin{equation*}
\Phi^{\omega}_s(\pi_{\omega}(\cR^{0}_{ph} \otimes \one)) \subset \Phi^{\omega}_s(\ker \drb_{\omega}) \subset \Phi^{\omega}_s(\pi_{\omega} (\cR^{0}_{t} \otimes \cA_g))= \Phi^{\omega}_s(\pi_{\omega}(\cR^{0}_{ph}\otimes \one)),
\end{equation*}
This and the definition of $\cP^{BRST}_{\omega}$ (cf. Definition \eqref{df:PQ2}) gives
\begin{equation}\label{eq:rphrkre}
\cP^{BRST}_{\omega}\cong\Phi^{\omega}_s(\ker \drb_{\omega})=(\Phi^{\omega}_s\circ\pi_{\omega})(\cR^{0}_{ph}\otimes \one)
\end{equation}
as algebras. 

\pfit (iii): Combining equation \eqref{eq:rphrkre} with the isomorphism $\vartheta$ above gives that 
\[
\cP^{BRST}_{\omega}\cong \pi_{(\omega_1,\cR_{ph})}(\cR^0_{ph}) 
\]
as algebras.
\end{proof}
From this we see that for each $\omega=\omega_1 \otimes \omega_2 \in \fS_\drb$ and $\omega_1 \in \SReg(\cR(\cD, \sigma_2))$, we have that $\cP^{BRST}_{\omega}$ is algebraically isomorphic to $\pi_{(\omega_1,\cR_{ph})}(\cR^0_{ph})$, which by Proposition \eqref{pr:physalgQEMii} is isomorphic to a dense subalgebra  $\cR^0_{ph}=\cR_0(\fD_t, \sigma_2)$ of the Dirac physical observables $\cP$. Hence we see that the BRST physical observables correspond closely to the Dirac physical observables from Subsection \ref{sec:cstbrstii}. We make this correspondence precise in the next section, after we have defined a the abstract BRST-physical algebra $\cP^{BRST}$ with suitable norm.
 
\subsection{QEM and Covariance II}\label{sbs:rescov}
A problematic feature of the auxiliary algebra $\cR(\fD,\sigma_2)$ is that the Poincar{\'e} automorphisms do not define naturally on it. This section is devoted to discussing this problem. 

Firstly we note that there is a `back door' solution to the problem. We found that in the covariant algebra case the Poincar{\'e} transformations define naturally and factor to the physical algebra (Subsection \ref{sbs:QEMcovI}). We have also found that the $T$-procedure for the covariant algebra and the auxiliary algebra produce the physical algebra, therefore we can use the factored Poincar{\'e} tranformations coming from the covariant algebra for the constrained auxiliary algebra also. 

This deals with relativistic covariance for the auxiliary algebra case but it is somewhat unsatisfactory from the standpoint that the covariant algebra and auxiliary algebra come from the same algebra of unbounded fields, and so either algebra should be able to contain the same information. Therefore we would expect that the Poincar{\'e} transformations should be able to be encoded as transformations on $\cR(\fD,\sigma_2)$ in some form. Such an encoding is not straightforward and the remainder of this subsection is devoted to a discussion of this issue.

Recall that we define the Poincar{\'e} transformations on the one-particle test function space as:
\begin{equation}\label{eq:oppoinact}
(V_gf)(p):=e^{ipa}\Lambda f(\Lambda^{-1}p) \quad \forall f \in \cS(\mathbb{R}^4,\mathbb{C}^4), \quad g=(\Lambda,a) \in \cP^{\uparrow}_{+}.
\end{equation}
where $\fD= \fX +i \fX$ as given in Subsection \eqref{sbs:testfunc}. Also recall as in Subsection \ref{sbs:QEMcovI}, that $V_g$ is $\sigma_1$-symplectic on $\fX$ and $\fD$ and so generates  $\alpha_g \in \Aut( \cR(\fD, \sigma_1))$ by Theorem \eqref{th:symautRA} (v). However note that we do not have that for all $g \in \cP^{\uparrow}_{+}$, $V_g$ is $\sigma_2$-symplectic on $\fX$ or $\fD$, for example take $g$ to be a Lorentz boost. Therefore we cannot use Theorem \eqref{th:symautRA} (v) to generate automorphisms for such $g$. 

At present the solution to directly encoding the Lorentz boosts for the auxiliary algebra is still a work in progress. Pursuing the issue here would be to great a digression, however strategies taken are discussed in Appendix \ref{ap:cvBRST II}. Briefly, these strategies are:
\begin{itemize}
\item We defined
 \[
\gamma_T:\sp(\fX, \sigma_2) \to\sp(\fX, \sigma_2)\qquad \text{defined by} \qquad M \to TMT, 
\]
where $T:\fX \to \fX$ is a \emph{real} linear operator such that $\gamma_T$ is a \emph{real} algebra isomorphism. Then $g \to \gamma_T(V_g) \in \sp(\fX,\sigma_2)$ is an real homomorphism of  $\cP^{\uparrow}_{+}\to \sp(\fX,\sigma_2)$. Unfortunately, as $\gamma_T$ is only real linear, it not extend to a complex homomorphism $\cP^{\uparrow}_{+}\to \sp(\fD,\sigma_2)$.
\item Now $\alpha_g:=Ad(\Gamma_{+}(V_g))$ is defined in the Fock representation in Subsection \ref{sbs:QEMcovI}. We see to what extent we define this on $\cR(\fD,\sigma_2)$. It turns out that $\Gamma_{+}(V_g)\in \op(\fF^{+}_0(\fD))$ is unbounded for general boosts, $\alpha_g$ does not preserve $\cR(\fD,\sigma_2)$ in the Fock representation and in fact maps bounded elements in $\cR(\fD,\sigma_2)$ to unbounded elements. A strategy is to encode  $\alpha_g$ in an infintesimal form (as a derivation) on parts $\cR(\fD,\sigma_2)$ from which we can recover a representation $g \to \alpha_g$ in certain representation. This strategy has difficulties discussed further in Appendix \ref{ap:cvBRST II}.
\end{itemize}

\section{General BRST for unbounded $\drb$}\label{sec:ubbcsbrst}
In this section we give an abstract definition of $C^*$-BRST for unbounded superderivations $\drb$, inspired by the bounded case and the BRST-QEM examples. We identify a set of states for which $\drb$ has a densely defined generating BRST charge $Q$, and define the abstract BRST-physical algebra $\cP^{BRST}$ with a suitable norm. Then we can establish the connection between the Dirac physical observables and BRST physical observables for the examples of abelian Hamiltonian BRST, and both versions BRST-QEM. 

The treatment for unbounded $\drb$ is very similar to the bounded case in Subsection \ref{sbs:bddcbrst}, but with differences due to the fact that we no longer have a BRST charge in the field algebra. Motivated by the BRST-QEM examples we modify the assumptions for bounded BRST in Subsection \ref{sbs:bddcbrst} as follows. 
\begin{itemize}
\item Let $\cA_0$ be a unital $C^*$-algebra and $\beta \in Aut(\cA_0)$ be an automorphism such that $\beta^2=\iota$ which encodes any Krein structure present in $\cA_0$. We call $\cA_0$ the \emph{Unextended Field Algebra} and we assume that it has a degeneracy, such as constraints. 
\item We tensor on a ghost algebra $\cA_g(\cH_2)$ where $\cH_2$ corresponds to the degrees of degeneracy, e.g. $\dim(\cH_2)=$ number of linear independent constraints in the Hamiltonian case (cf. Definition \ref{df:GA}), and let $\alpha'
\in \Aut(\cA_g)$ be the automorphism that corresponds to the Krein-Ghost stucture (cf. equation \eqref{eq:autghalg}), in particular $(\alpha')^2=\iota$. 
\end{itemize}
\begin{definition}\label{df:ufagbbrst}
The \emph{BRST-Field Algebra} is $\cA:=\cA_0\otimes \cA_g$ or $\cA:=\cA_0\otimes \rga$ and we let $\cA$ have $\mathbb{Z}_2$-grading with grading automorphism $\iota\otimes\gamma$ (cf. $\iota\otimes\gamma$ (cf. Definition \eqref{df:Z2grad}). The norm on $\cA$ is unique as $\cA_g$ is a CAR algebra and hence nuclear. Let $\alpha:=\beta\otimes \alpha' \in \Aut(\cA)$ and note that $\alpha^2=\iota$. Define the involution on $\cA$:
\[
A^{\dag}=\alpha(A^*), \qquad \forall A\in \cA.
\] 
Furthermore we assume that there exists a subalgebras $D_1(\drb)\subset D_2(\drb) \subset \cA$ that are not necessarily closed or norm dense in $\cA$, but that $\gamma(D_1(\drb))=D_1(\drb)$, $\gamma(D_2(\drb))=D_2(\drb)$ and $D_1(\drb)^*= D_1(\drb)$. We assume that there exists a (possibly unbounded) superderivation:
\begin{align*}
\drb:D_2(\drb) \to D_2(\drb),
\end{align*}
graded with respect to $\gamma$, and such that:
\begin{itemize}
\item[(i)] $\drb^2=0$ on $D_2(\drb)$,
\item[(ii)] $\gamma \circ \drb \circ \gamma= -\drb$
\item[(iii)] $\drb(A)^*=- \alpha \circ \drb \circ \alpha \circ \gamma(A^*)$
\end{itemize}
\end{definition}
The motivation for defining the domains $D_1(\drb)$ and $D_2(\drb)$ for unbounded $\drb$ is motivated by the BRST-QEM examples (cf. Definition \eqref{df:sddom2} and Definition \eqref{df:sddom1}). The assumptions (i), (ii), (iii) correspond to $Q^2=0$, $Q \in \cA^{-}$ and $Q^{\dag}=Q$ in the bounded case (cf. lemma \eqref{lm:propbsd}).

We can try to calculate the cohomogical version of the BRST physical algebra $\ker \drb / \ran \drb$ algebraically, however looking at examples such BRST-QEM (Theorem \eqref{th:calgdrb2}) we see that this is not always a tractable problem. The approach taken in the Hamiltonian BRST example and QEM-BRST II example was to use the projection $P_s$ from the \emph{dsp}-decomposition (Theorem \eqref{th:Hdsp}) for the charge $Q$ (with $P_s$ and $Q$ defined in the appropriate representation) to calculate the BRST-physical algebra $\cP^{BRST}$. This can then be connected to the cohomological BRST definition of the physical algebra via Theorem \eqref{pr:krdel}. For this strategy to work in general we need to select states such the projection $P_s$ exists in their corresponding representations. Hence we choose states $\omega \in \fS(\cA)$ such that $\drb$ will have a generator $Q_{\omega}$ in the GNS-representations associated to $\omega$ with the usual properties as in Section \ref{sec:dsprig}.

As for the bounded case we want states that give a space on which $\cA_0$ acts tensored to a ghost space, \ie:
\[
\fS_T:=\{ \omega \in \fS(\cA) \,|\, \omega=\omega_1 \otimes \omega_2, \; \omega_1 \in \fS(\cF), \omega_2\in \fS_g \}.
\]
By Definition \eqref{df:ghst} $\omega_2 \in \fS_g$ implies that $\omega_2 \circ \alpha' =\omega$ hence $\cH_{\omega_2}$ is a Krein space with $J_{\omega_2}$ implementing $\alpha$ and $J_{\omega_2}\Omega_{\omega_2}=\Omega_{\omega_2}$ hence $\Omega_{\omega_2}$ is positive in the Krein inner product with fundamental symmetry $J_{\omega_2}$ on $\cH_{\omega_2}$. If we want the same for $\cH_{\omega_1}$ then we assume that $\omega_1\circ \beta =\omega_1$. We want to choose a subset of $\fS_T$ such that we can construct an operator $Q$ that generates $\drb$.
\begin{definition}\label{df:ubsgsbdd}
Let $\fS_{\drb}$ be states of the form $\omega=\omega_1 \otimes \omega_2$, where $\omega_1 \circ \beta= \omega_1$, $\omega_2 \in \fS_{g}$, $\omega(\drb(A))=0$ $\forall A \in D_2(\drb)$ and $\overline{\pi_{\omega}(D_1(\drb))\Omega_{\omega}}=\overline{\pi_{\omega}(D_2(\drb))\Omega_{\omega}}=\cH_{\omega}$. Let $\iip{\cdot}{\cdot}_{\omega}$ and $J_{\omega}$ be the Krein inner product and fundamental symmetry as in Definition \eqref{df:csksrepterm}.
\end{definition}
The above definition is motivated by Defintion \eqref{df:bsgsbdd} and lemma \eqref{lm:bQkdel} (i). The following theorem justifies the above definition (the proof similar to \cite{Hendrik1991} p29).
\begin{theorem}\label{th:abQ2}
Let $\omega= \omega_1\otimes \omega_2 \in \fS_{\drb}$. Then there exists a $\iip{\cdot}{\cdot}_{\omega}$-symmetric, 2-nilpotent operator $Q_{\omega}$ on $\cH_{\omega}$ which preserves the dense domain $D(Q_{\omega}):=\pi_{\omega}(D_2(\drb)) \Omega_{\omega}$ and is such that 
\begin{align*}
\pi_{\omega}(\drb(A))\psi=\sbr{Q_{\omega}}{ \pi_{\omega}(A) } \psi, \qquad \forall \psi \in \pi_{\omega}(D_2(\drb))\Omega_{\omega},
\end{align*}
and for all $A \in D_2(\delta)$. Moreover $\Omega_{\omega}\in \ker Q_{\omega}$ and: 
\begin{itemize}
\item[(i)] $Q_{\omega}$ is closable.
\item[(ii)] $\ran \cQ_{\omega} \subset D(\cQ_{\omega})$ and $\cQ_{\omega}^2\psi =0$ for all $\psi \in D(\cQ)$. 
\item[(iii)] $(Q_{\omega}^*)^2\psi=0$ for all $\psi \in D(Q_{\omega}^*)$  
\end{itemize}
\end{theorem}
\begin{proof} As $\pi_{\omega}(D_2(\drb)) \Omega_{\omega}$ is dense in $\cH_{\omega}$ we define $Q_{\omega}$ on $\pi_{\omega}(D_2(\drb)) \Omega_{\omega}$. Take $\psi= \pi_{\omega}(A)\Omega_{\omega}$ where $A \in D_2(\drb)$. We define,
\[
Q_{\omega}\psi=Q_{\omega} \pi_{\omega}(A) \Omega_{\omega}:=\pi_{\omega}(\drb(A))\Omega_{\omega},
\]
and check that this is well defined. Linearity is obvious, so as $\pi_{\omega}(D_2(\drb))$ is dense in $\cH_{\omega}$, we just have to check that $\psi=\pi_{\omega}(A)\Omega_{\omega}=0$ implies that $Q_{\omega}\psi=\pi_{\omega}(\drb(A))\Omega_{\omega}=0$ for all $A \in D_2(\delta)$. For all $A \in D_2(\delta)$, $B \in D_1(\delta)$ we have
\begin{align*}
\omega(\drb(A)^*B)&=\ip{\pi_{\omega}(\drb(A))\Omega_{\omega}}{\pi_{\omega}(B)\Omega_{\omega}},\\
&=\overline{\omega(B^*\drb(A))},\\
 &= \overline{\omega(\drb(\gamma(B^*)A)-\drb(\gamma(B^*))A))},\\
 &=-\omega(A^*\drb(\gamma(B^*))^*) \qquad \text{(as $\omega(\drb(\cdot))=0$)}, \\
 &= -\ip{\pi_{\omega}(A)\Omega_{\omega}}{\pi_{\omega}(\drb(\gamma(B^*))^*)\Omega_{\omega}},
\end{align*}
and so $\pi_{\omega}(A)\Omega_{\omega}=0$ and $\pi_{\omega}(D_1(\drb))\Omega_{\omega}$ dense in $\cH_{\omega}$ imply that $\pi_{\omega}(\drb(A))\Omega_{\omega}=0$.

With $Q_{\omega}$ well defined we use the above calculation and $\drb(B)^*=- \alpha(\delta(\alpha(\gamma(B^*))))$ to calculate, 
\begin{align*}
\ip{Q_{\omega} \pi_{\omega}(A) \Omega_{\omega}}{\pi_{\omega}(B)\Omega_{\omega}}&=\omega(\drb(A)^*B), \\
&= -\omega(A^*\drb(\gamma(B^*))^*),\\
&= \omega(A^*\alpha(\delta(\alpha(B)))),\\
&=\omega(\alpha(A^*)\delta(\alpha(B))),\\
&=\ip{J_{\omega}\pi_{\omega}(A)\Omega_{\omega}}{\pi_{\omega}(\delta(\alpha(B)))\Omega_{\omega}},\\
&=\ip{J_{\omega}\pi_{\omega}(A)\Omega_{\omega}}{Q_{\omega}\pi_{\omega}(\alpha(B))\Omega_{\omega}},\\
&=\ip{\pi_{\omega}(A)\Omega_{\omega}}{(J_{\omega}Q_{\omega}J_{\omega}\pi_{\omega}(B))\Omega_{\omega}}.
\end{align*}
As $\pi_\omega(D_1(\drb))\Omega_{\omega}$ is dense in $\pi_\omega(D_2(\drb))\Omega_{\omega}$, and as $J_{\omega}^*=J_{\omega}$, we have that 
\[
Q_{\omega} \subset J_{\omega} Q_{\omega}^* J_{\omega},
\]
and so $Q_{\omega}$ is Krein symmetric.

Furthermore, let $\psi \in D(Q_{\omega})= \pi_{\omega}(D_2(\drb)\Omega_{\omega}$. Then $\psi =\pi_{\omega}(A) \Omega_{\omega}$ for some $A \in D_2(\drb)$. Therefore $Q_{\omega}\psi=\pi_{\omega}(\drb(A))\Omega_{\omega}  \in \pi_{\omega}(D_2(\drb)\Omega_{\omega}=D(Q_{\omega})$, and so $Q_{\omega}$ preserves $D(Q_{\omega})$. Furthermore $\drb^2(A)=0$, so we have that 
\[
Q_{\omega}^2 \psi= Q_{\omega} \pi_\omega(\drb(A))\Omega_{\omega} =\pi_\omega(\drb^2(A))\Omega_{\omega}=0. 
\]
and so $Q_{\omega}$ is 2-nilpotent. Letting $A=\one$ gives that $\Omega_{\omega} \in \ker Q_{\omega}$.

Now (i), (ii) and (iii) follow from lemma \eqref{lm:QessaCsa}
\end{proof}
An extra complication in the unbounded $\drb$ case as compared to the bounded $\drb$ case is that we would like to be able to restrict to onvenient subsets of $\fS_{\drb}$, e.g. the strongly regular states in the BRST-QEM case II (cf. Subsection \ref{sbs:srs}). We modify Definition \eqref{df:bdspdec} accordingly. 
\begin{definition}\label{df:ubdspdec} Let $\fS_{\WW} \subset \fS_{\drb}$,
\begin{itemize}
\item[(i)]Let the representation $\pi_{\WW}:\cA \to B(\cH_{\WW})$ be defined by,
\[
\cH_{\WW}:=\bigoplus\{ \cH_{\omega}\,|\, \omega \in \fS_{\WW}\}, \qquad \pi_{\WW}:\bigoplus \{ \pi_{\omega}\,|\, \omega \in \fS_{\WW}\}.  
\]
Denote the Hilbert inner product on $\cH_{\WW}$ by $\ip{\cdot}{\cdot}_{\WW}$, and let $P_{\omega}\in B(\cH_{\WW})$ denote the  projection onto $\cH_{\omega}$. 
\item[(ii)] For $\omega\in \fS_{\WW}$,  $\omega \circ \alpha =\omega$ hence $\alpha$ is unitarily implemented in each $\pi_{\omega}$ hence $\alpha$ is unitarily implemented in $\cH_{\WW}$. Denote the implementer for $\alpha$ in $\pi_{\WW}$ by $J^{\WW}$. As $J^{\WW}|_{\cH_{\omega}}=J_{\omega}$ it follows from $\alpha^2=\iota$ that $(J^{\WW})^2=\one$ and $J^{\WW*}=J^{\WW}$. By lemma \eqref{lm:JKsp} $\cH_{\WW}$ is a Krein space with fundamental symmetry $J^{\WW}$ and indefinite inner product $\iip{\cdot}{\cdot}_{\cH_{\WW}}:=\ip{\cdot}{J^{\WW}\cdot}_{\cH_{\WW}}$. 
\item[(iii)] For $\omega\in \fS_{W}$, let $Q_{\omega}$ be the BRST charge as in Theorem \eqref{th:abQ2}, with dense domain $D(Q_{\omega})$. Let 
\[
D(Q_\WW):=\left\lbrace \psi \in \cH_{\WW} \,\left| \, P_{\omega}\psi \in D(Q_{\omega}), \; \forall \omega\in \fS_{\WW}, \, \sum_{\omega  \in \fS_{\WW}}\right. \norm{Q_{\omega}P_{\omega}
\psi}_{\cH_{\omega}}^2 < \infty \right\rbrace 
\]
and 
\[
Q_{\WW} \psi:= \sum_{\omega \in \fS_{\WW}} Q_{\omega}P_{\omega}\psi, \qquad \psi  \in D(Q_{\WW}).
\]
\item[(iv)] We have $Q_{\WW}^2\psi=0$ and $Q_{\WW}$ is $\iip{\cdot}{\cdot}_{\WW}$-symmetric as $Q_{\omega}$ for each $\omega \in \fS_{\WW}$. Hence $Q_\WW$ is Krein symmetric and and $\cQ_\WW$ satisfies the conditions of the \emph{dsp}-decomposition (Theorem \eqref{th:Hdsp1}) via lemma \eqref{lm:QessaCsa}. Let $\cH_{\WW}=\cH^d_{\WW}\oplus \cH^s_{\WW}\oplus \cH^p_{\WW}$, $\cH_{\omega}=\cH^d_{\omega}\oplus \cH^s_{\omega}\oplus \cH^p_{\omega}$ be the \emph{dsp}-decompositions with respect to $\cQ_W$ and $\cQ_{\omega}$ where $\omega \in \fS_{\WW}$. Let $P^{k}_{j}$, $k=\omega, \WW$, $j=s,p,d$ be the corresponding projections on $\cH_{\WW},\cH_{\omega}$. 
\item[(v)] Let $\omega \in \fS_{\WW}$. For $j=\WW, \omega$ define,
\[
\cH^{BRST}_{phys,j}:=\ker Q_{j}/\cH^d_{j}.
\]
and let $\vp_{j}:\ker Q_{j} \to  \cH^{BRST}_{phys,j}$ be the factor map. Denote $\hat{\psi}:=\vp_{\WW}(\psi)$ for $\psi \in \cH_{\WW}$ and $\hat{\psi}{}^{\omega}=\vp_{\omega}(\psi)$ for $\psi \in \cH_{\omega}$. 
\end{itemize}
\end{definition}
\begin{rem}\label{rm:obvious2}
As $\cH_{\WW}=\oplus_{\omega \in \fS_{\WW}}\cH_{\omega}$ it is obvious that $\cH^{j}_{\WW}=\oplus_{\omega \in \fS_{\WW}}\cH^{j}_{\omega}$ for $j=d,s,p$, and hence $P^{\WW}_{j}=\oplus_{\omega \in \fS_{\WW}}P^{\omega}_{j}$, $j=d,s,p$.
\end{rem}
To get the spatial structures of Chapter \ref{ch:GenStruct}:
\begin{proposition}\label{pr:ucbiksp}
We have $\cH^{BRST}_{phys,\WW}$ in Definition \eqref{df:ubdspdec} has an indefinite inner product defined for all $\psi, \xi \in \ker \cQ_\WW$: 
\begin{equation*}
\iip{\hat{\psi} }{\hat{\xi}}_p:=\iip{\psi}{\xi}_{\cH_{\WW}}=\iip{P^{\WW}_s\psi}{P^{\WW}_s\xi }_{\cH_{\WW}}.
\end{equation*} 
Furthermore, if $J^{\WW}\cH^s_{\WW}=\cH^s_{\WW}$ then $(\cH^{BRST}_{phys,\WW}, \iip{\cdot}{\cdot}_{\cH_{\WW}})$ is Krein space with fundamental symmetry, $J^{\WW}_p\hat{\psi}:=\vp_{\WW}(J^{\WW}P^{\WW}_s\psi)$ and Hilbert inner product
\begin{equation}\label{eq:uBRSTpHilip2}
\ip{\hat{\psi}}{\hat{\xi}}_p:=\iip{\hat{\psi} }{J^{\WW}_p\hat{\xi}}_p= \ip{P^{\WW}_s\psi}{P^{\WW}_s\xi }_{\cH_{\WW}},
\end{equation} 
and norm $\norm{\hat{\psi}}_p:=\ip{\hat{\psi}}{\hat{\psi}}_p^{1/2}$. If
\begin{equation}\label{eq:uppossub}
J^{\WW}P^{\WW}_s=P^{\WW}_s,
\end{equation}
then $J^{\WW}_p=\one$ hence $\iip{\hat{\psi} }{\hat{\xi}}_p=\ip{\hat{\psi}}{J^{\WW}_p\hat{\xi}}_p$, \ie the Krein structure is the same as the Hilbert structure.
\end{proposition}
\begin{proof}
We have that $\cQ_{\WW}$ and $\cH_{\WW}$ satisfy the hypothesis of lemma \eqref{lm:ranQnll} hence  $\iip{\cdot}{\cdot}_{\cH_{\WW}})$ is well defined by the above formula on $\cH^{BRST}_{phys,\WW}$. The rest of the proposition follows from lemma \eqref{lm:phspksp}.
\end{proof}

To get the algebraic structures as in Chapter \ref{ch:GenStruct}:
\begin{definition}\label{df:ubddbrshom}
Define the linear map:
\begin{align*}
\Phi^{\WW}_s:\cA &\to B(\cH_{\WW}),\\
\Phi^{\WW}_s(A):&=P^{\WW}_s\pi_{\WW}(A)P^{\WW}_s,
\end{align*}
By lemma \eqref{lm:alkerhom} we have that $\Phi^{\WW}_s$ is a homomorphism on $\ker \drb$. We define the \emph{BRST-physical algebra} as,
\[
\cP_0^{BRST}:=\ker \drb / (\ker \drb \cap \ker \Phi^{\WW}_s)\cong \Phi^{\WW}_s(\ker \drb).
\]
where $\cong$ above denotes an algebra isomorphism. Let the factor map be $\tau: \ker \drb \to \cP^{BRST}$, and denote $\hat{A}:=\tau(A)\in \cP^{BRST}$ for $A \in \ker \drb$.  
\end{definition}
The above definition is motivated by the alternative definition of the BRST-physical algebra, cf. Subsection \ref{sbs:alalg} and Theorem \eqref{pr:krdel} for the connection to the usual cohomological definition of the BRST-physical algebra.  

Note that $\Phi^{\WW}_s$ is not necessarily a $*$-isomorphism as $\ker \drb$ is not necessarily a $*$-algebra (\eg for the BRST-QEM example in Subsection \ref{sbs:sdI} we have for $h \in \fDL$ that $\zeta_1(h)^* \in \ker \drb$, $\zeta_1(h)\notin \ker \drb$).
As in the bounded case, to get a natural norm for $\cP_0^{BRST}$ we use the norm on $\cH^{BRST}_{phys,\WW}$. 
\begin{proposition}\label{pr:urepban}
Define the representation,
\begin{align*}
\pi_{\WW, p}:\cP_0^{BRST}&\to B(\cH^{BRST}_{phys,\WW}),\\
\pi_{\WW, p}( \hat{A})\hat{\psi}:&= \widehat{\pi_{\WW}(A)\psi}, 
\end{align*}
and the seminorm on $\cP_0^{BRST}$ by: 
\[
\norm{\hat{A}}_p:=\norm{\pi_{\WW, p} ( \hat{A})}_{B(\cH^{BRST}_{phys,\WW})}.
\]
Then $\norm{\hat{A}}_p$ is a norm on $\cP_0^{BRST}$,
\begin{equation}\label{eq:unmfact}
\norm{\hat{A}}_p=\norm{\Phi^{\WW}_s(A)}_{\cH_{\WW}},
\end{equation}
and $\cP^{BRST}:=\overline{\cP_0^{BRST}}$ is a Banach algebra where closure is with respect to $\norm{\hat{A}}_p$. Furthermore ${\cP^{BRST}}\cong \overline{\Phi^{\WW}_s(\ker \drb)}^{\cB(\cH_{\WW})}$ where the isomorphism is isometric. Thus $\pi_{\WW,p}$ is a faithful representation of $\cP^{BRST}$ and all calculations can be done in this representation.
\end{proposition}
\begin{proof}
We adapt the proof of Proposition \eqref{pr:repban}. First,
\[
A \in \ker \drb \Rightarrow A \in D_2(\drb) \Rightarrow \pi_{\omega}(A)D(Q_{\omega}) \subset D(Q_{\omega}),\: \forall \omega \in \fS_{W} \Rightarrow \pi_{\WW}(A)D(Q_{\WW})\subset D(Q_{\WW}) .
\] 
Therefore $A \in \ker \drb \Rightarrow \pi_{\WW}(A)\ran Q_{\WW} \subset \ran Q_{\WW}$, and so as $\pi_{\WW}(A) \in B(\cH_{\WW})$ and $\ran Q_{\WW}$ is dense in $\cH^d_{\WW}$ we get that $\pi_{\WW}(A)\cH^d_\WW \subset \cH^d_\WW$. Therefore
\begin{equation}\label{eq:unAinker}
A \in \ker \drb \Rightarrow \pi_{\WW}(A)\cH^d_{\WW}\subset \cH^{d}_{\WW} \Rightarrow P^{\WW}_s \pi_{\WW}(A) \ker \cQ_{\WW}=P^{\WW}_s \pi_{\WW}(A)P^{\WW}_s \ker \cQ_{\WW}, 
\end{equation}
Using this we see that $\norm{\hat{A}}_p$ is a norm by the calculation:
\begin{align*}
\norm{\hat{A}}_p=&\; \sup \{ \norm{\pi_{\WW, p}( \hat{A})\hat{\psi}}_p\;|\; \psi \in \ker  \pi_{\WW}(Q),\,\norm{\hat{\psi}}_{p}\leq 1\},\\
=&\; \sup\{ \norm{ P^{\WW}_s \pi_{\WW}(A)\psi}_{\cH_{\WW}}\;|\; \psi \in \ker  \pi_{\WW}(Q),\,\norm{{P^{\WW}_s} \psi}_{\cH_{\WW}}=1\} \\
=&\;\norm{\Phi^{\WW}_s(A)}_{\cH_{\WW}},
\end{align*}
for $A \in \ker \drb$ where we have used in the second equality that  $\norm{\hat{\psi}}_p=\norm{P^{\WW}_s\psi}_{\cH_{\WW}}$ for $\psi \in \ker \cQ_\WW$ by equation \eqref{eq:uBRSTpHilip2}, and equation \eqref{eq:unAinker} in the third. Therefore $\norm{\hat{A}}_p=0$ iff $\norm{\Phi^{\WW}_s(A)}_{\cH_{\WW}}=0$ iff $\Phi^{\WW}_s(A)=0$  iff $\hat{A}=0$, and so $\norm{\cdot}_p$ is a norm. Let $A,B \in \ker \drb$. As $\norm{P^{\WW}_s}_{\cH_{\WW}}=1$ we have,
\[
\norm{\hat{A}\hat{B}}_p=\norm{\Phi^{\WW}_s(AB)}_{\cH_{\WW}}\leq\norm{\Phi^{\WW}_s(A)}_{\cH_{\WW}} \norm{\Phi^{\WW}_s(B)}_{\cH_{\WW}}=\norm{\hat{A}}_p\norm{\hat{B}}_p.
\]
This shows that $\overline{\cP^{BRST}}$ is a Banach algebra where closure is with respect to $\norm{\hat{A}}_p$.
\end{proof}
Note again that  $\ker \drb$ is not in general a $*$-algebra, and so in general ${\cP^{BRST}}\cong \overline{\Phi^{\WW}_s(\ker \drb)}^{\cB(\cH_{\WW})}$ is not a $*$-isomorphism. Now $\ker \drb$ is not a $*$-algebra in general so for $A \in \ker \drb$,  $\widehat{A^*}$ is not necessarily well defined and so we have to be careful of how we define a $*$-involution on  $\cP^{BRST}$. 
\begin{rem}
In order that the Krein involution $\dag$ on $\ker \drb$ factors to $\cP^{BRST}$ we assume below an extra condition that $J^{\WW}$ induces in indefinite innner product such that $\cH^{BRST}_{phys, \WW}$ is a Krein space in a natural way. This extra condition is sufficient but it may not be necessary. It is a reasonable condition to assume in this thesis as it holds for the all the examples considered. 
\end{rem}
\begin{proposition}\label{pr:ubddcshinv}
We have:
\begin{itemize}
\item[(i)] Assume $J^{\WW}\cH^{s}_{\WW}=\cH^{s}_{\WW}$. Then $(\ker \drb)^{\dag}=\ker \drb$ and $(\ker \drb \cap \ker \Phi^{\WW}_s)^{\dag}=(\ker \drb \cap \ker \Phi^{\WW}_s)$. Hence $\dag$ on $\ker \drb$ factors to the $\dag$-involution on $\cP^{BRST}$ which coincides with the $\dag$-involution with respect to the representation $\pi_{\WW,p}$, \ie
\[
\pi_{\WW,p}(\widehat{A^{\dag}})\hat{\psi}=\pi_{\WW,p}(\hat{A})^\dag\hat{\psi}
\]
for all $A \in \ker \drb$ where $\pi_{\drb,p}(\hat{A})^{\dag}$ is the the adjoint of $\pi_{\WW,p}(\hat{A})$ with respect to the inner product $\iip{\cdot}{\cdot}_p$. Furthermore,
\[
\cP^{BRST} \cong \overline{\Phi_s^{\WW}(\ker \drb)}
\]
where the above is an isometric $\dag$-isomorphism.
\item[(ii)] Let $\cM \in \ker \drb$ be a subalgebra such that $\Phi^{\WW}_s(\cM)=\Phi^{\WW}_s(\cM^*)$. Given $A \in \cM$, define
\begin{equation}\label{eq:udfhadpa}
\hat{A}^*:=\hat{B},
\end{equation} 
where $\Phi^{\WW}_s(A^*)=\Phi^{\WW}_s(B)$ for some $B \in \cM$. This defines an involution on $\cM/ (\cM \cap \ker \Phi^{\WW}_s)$ such that $\overline{\cM/ (\cM \cap \ker \Phi^{\WW}_s)}$ is a $C^*$-algebra where closure is with respect to $\norm{\cdot}_p$.
\item [(iii)] Let the physicality condition $J^{\WW}P^{\WW}_s=P^{\WW}_s$ hold (equation \eqref{eq:ppossub}). Then $\ker \drb$ satisfies the conditions on $\cM$ in (ii), hence $\cP^{BRST}$ is a $C^*$-algebra with respect to the norm $\norm{\cdot}_p$. Moreover, the $\dag$-involution from (i) and $*$-involution from (ii) coincide.
\end{itemize}
\end{proposition}
\begin{proof} We adapt the proof of Proposition \eqref{pr:bddcshinv} to this context. 
(i): $(\ker \drb)^{\dag}=\ker \drb$ follows from $\drb(A)^{\dag}=-  \drb \circ  \gamma(A^{\dag})$ (cf. Definition \eqref{df:ufagbbrst} (iii)). By the assumption $J^{\WW}\cH^{s}_{\WW}=\cH^{s}_{\WW}$ we have $[P^{\WW}_s, J^{\WW}]=0$ hence $(P^{\WW}_s)^{\dag}=J^{\WW}P^{\WW*}_sJ^{\WW}=P^{\WW}_s$ and so it follows that $(\ker \drb \cap \ker \Phi^{\WW}_s)$ is a $\dag$-subalgebra of  $\ker \drb$. Hence the involution $\dag$ factors to $\cP^{BRST}$.  Moreover, it coincides with the $\dag$-involution with respect to the representation $\pi_{\WW,p}$ which can be seen by this  calculation:
\[
\iip{\pi_{\WW,p}(\widehat{A^{\dag}})\hat{\psi}}{\hat{\xi}}_p=\iip{\widehat{\pi_{\WW}(A^{\dag}\psi)}}{\hat{\xi}}_p=\iip{\pi_{\WW}(A^{\dag})\psi}{\xi}=\iip{\psi}{\pi_{\WW}(A)\xi}=\iip{\hat{\psi}}{\pi_{\WW,p}(\hat{A})\hat{\xi}}_p,
\]
for all $\psi,\xi \in \cH^{BRST}_{phys, \WW}$ and all $A \in \ker \drb$, since by definition $\iip{\hat{\psi}}{\hat{\xi}}_p=\iip{\psi}{\xi}$ for all $\psi,\xi \in \ker \cQ_\WW$ (cf. Proposition \eqref{pr:ucbiksp}). 

Let $A \in \ker \drb$. Then using $[P^{\WW}_s, J^{\WW}]=0$,
\[
\Phi_s^{\WW}(A^{\dag})=\Phi_s^{\WW}(\alpha(A)^*)=P_s^{\WW}J_{\WW}\pi_{\WW}(A)^*J_{\WW}P_s^{\WW}=(J_{\WW}P_s^{\WW}\pi_{\WW}(A)P_s^{\WW}J_{\WW})^*=\Phi_s^{\WW}(A)^{\dag}
\]
where we used $J_{\WW}^*=J_{\WW}$ in the last line. Combining this with Proposition \eqref{pr:urepban} gives that there is a isometric $\dag$-isomorphism such that $\cP^{BRST} \cong \overline{\Phi_s^{\WW}(\ker \drb)}$.  

\smallskip
\noindent(ii): Let $A\in (\cM \cap \ker \drb)$ then by assumption there exists $B\in \cM$ such that $\Phi^{\WW}_s(A^*)=\Phi^{\WW}_s(B)$. Hence for all $\xi,\psi \in  \ker \cQ_{\WW}$ we have by equation \eqref{eq:uBRSTpHilip2} and equation \eqref{eq:unAinker} that,
\begin{align}
\ip{\hat{\xi}}{\pi_{\WW, p}(\hat{A}) \hat{\psi}}_p= & \;\ip{P^{\WW}_s\xi}{P^{\WW}_s \pi_{\WW}(A)\psi}_{\cH_{\WW}},\notag \\
= & \;\ip{P^{\WW}_s\xi}{P^{\WW}_s \pi_{\WW}(A)P^{\WW}_s\psi}_{\cH_{\WW}},\notag \\
= & \;\ip{\Phi^{\WW}_s(A^*)\xi}{\psi}_{\cH_{\WW}},\notag \\
= & \;\ip{\Phi^{\WW}_s(B)\xi}{\psi}_{\cH_{\WW}}, \notag \\
= & \;\ip{\pi_{\WW, p}(\hat{B})\hat{\xi}}{\hat{\psi}}_p.\label{eq:uipad}
\end{align}
Hence 
\[
\pi_{\WW, p}(\hat{A})^*=\pi_{\WW, p}(\hat{B})=\pi_{\WW, p}((\hat{A})^*)\in \cM/ (\cM \cap \ker \Phi^{\WW}_s)
\]
where we used equation \eqref{eq:udfhadpa} and that $B \in \cM$. This shows that $\cM/ (\cM \cap \ker \Phi^{\WW}_s)$ is a $*$-algebra. Furthermore Proposition \eqref{pr:urepban} gives that $\pi_{\WW, p}:\cM/ (\cM \cap \ker \Phi^{\WW}_s)\to B(\cH^{BRST}_{phys, \WW})$ is an isometric isomorphism, and so we have proved that $\pi_{\WW, p}$ is a $*$-isometric isomorphism. Hence as $B(\cH^{BRST}_{phys, \WW})$ is a $C^*$-algebra so is  $\overline{\cM/ (\cM \cap \ker \Phi^{\WW}_s)}$.

\pfit (iii): Now $A \in \ker \drb$ implies $A^{\dag} \in \ker \drb$, and by equation \eqref{eq:uppossub} 
\[
\Phi^{\WW}_s(A^{\dag})= P^{\WW}_s J^{\WW}\pi_{\WW}(A^*)J^{\WW} P^{\WW}_s=P^{\WW}_s \pi_{\WW}(A^*) P^{\WW}_s=\Phi^{\WW}_s(A^{*}).
\] 
Therefore we can apply (ii) with $\cM=\ker \drb$ to get that ${\cP^{BRST}}=\overline{\ker \drb/ (\ker \drb \cap \ker \Phi^{\WW}_s)}$ with norm $\norm{\cdot}_p$ is a $C^*$-algebra.  Moreover by (i) and the defining equation \eqref{eq:udfhadpa}, $\hat{A}^{\dag}=\widehat{A^{\dag}}=\hat{A}^*$.
\end{proof}
\begin{rem} \begin{itemize}
\item[(i)] Note that the above Proposition does \emph{not} assume that $\cM$ is a $C^*$-algebra. If is factoring out by $(\cM \cap \ker \Phi^{\drb}_s)$ and using the the Hilbert $*$-involution coming from the BRST physical space $\cH^{BRST}_{phys,\drb}$ that gives a $C^*$-algebra.
\item[(ii)] Proposition \eqref{eq:udfhadpa} (iii) shows that equation \eqref{eq:ppossub} is a good physicality condition as it ensures that the $\dag$-involution factors to a $C^*$-involution on the physical algebra $\cP^{BRST}$.
\end{itemize}
\end{rem}

Summarising we get:
\begin{theorem}\label{th:ubdbrststruc}
Let $\cA,\alpha,Q,\drb,\fS_{\WW}, \pi_{\WW}:\cA \to B(\cH_{\WW})$ be as in Definitions \eqref{df:ufagbbrst},\eqref{df:ubsgsbdd}, \eqref{df:ubdspdec}, and let 
\[
\Phi^{\WW}_s(A):=P^{\WW}_s\pi_{\WW}(A) P^{\WW}_s, \qquad A \in \cA.
\] 
as in Definition \eqref{df:ubddbrshom}. Then:
\begin{itemize}
\item[(i)]$\Phi^{\WW}_s$ is a homomorphism on $\ker \drb$ and on $(\ker \drb)^*$. 
\item[(ii)]Let  $\cP_0^{BRST}=\ker \drb / (\ker \drb \cap \ker \Phi^{\WW}_s)$ and $\pi_{\WW, p}:{\cP_0}^{BRST} \to B(\cH^{BRST}_{phys,\WW})$ be as in Proposition \eqref{pr:urepban}. Then $\cP^{BRST}$ has the norm 
\[
\norm{\hat{A}}_p:=\norm{\pi_{\WW, p} (\hat{A})}_{B(\cH^{BRST}_{phys,\WW})}=\norm{\Phi^{\WW}_s(A)}_{\cH_{\WW}}, \qquad \forall A \in \ker \drb
\]
with respect to which $\cP^{BRST}:=\overline{\cP_0^{BRST}}$ it is a Banach algebra and we have a isometric isomorphism such that ${\cP^{BRST}}\cong \overline{\Phi^{\WW}_s(\ker \drb)}^{\cB(\cH_{\WW})}$. Furthermore, if $J^{\WW}\cH_s^{\WW}=\cH_s^{\WW}$ then $\cP^{BRST}$ is a $\dag$-Banach algebra and is $\dag$-isometrically isomorphic to $ \overline{\Phi^{\WW}_s(\ker \drb)}^{\cB(\cH_{\WW})}$.
\item[(iii)] If $\Phi^{\WW}_s(\ker \drb)= \Phi^{\WW}_s((\ker \drb)^*)$, then ${\cP^{BRST}}$ is a $C^*$-algebra with norm $\norm{\cdot}_p$ and involution denoted by $*$ as defined in equation \eqref{eq:udfhadpa}. 

When the physcality condition $J^{\WW}P_s^{\WW}=P_s^{\WW}$ is satisfied, we have that $\Phi^{\WW}_s(\ker \drb)= \Phi^{\WW}_s((\ker \drb)^*)$, hence  $\cP^{BRST}$ is a $C^*$-algebra.
\end{itemize}
\end{theorem}
\begin{proof}
(i): Follows from lemma \eqref{lm:alkerhom} applied to  for $\cQ_{\WW}$ and $Q_{\WW}^*$. (ii) is Proposition \eqref{pr:urepban} and Proposition \eqref{pr:ubddcshinv} (i). (iii) follows from Proposition \eqref{pr:ubddcshinv}.
\end{proof}

\begin{rem}\label{rm:upscs}
\begin{itemize} 
\item Theorem \eqref{th:ubdbrststruc} (iii) gives a condition to check if ${\cP^{BRST}}$ is a $C^*$-algebra. For the case of bounded BRST charge $Q$, we saw in Example \eqref{ex:pscs} that by extending $\cA$ to $\tilde{\cA}=\osalg{P_s^\drb,\, \pi_u(\cA)}$ gives that this condition is satisfied hence ${\cP^{BRST}}$ defined using $\tilde{\cA}$ is a $C^*$-algebra. In the bounded case, extension to $\tilde{\cA}$ is straightforward as $Q \in \cA$ and $D(\drb)=\cA$. However, in the unbounded case a similar extension is not straightforward it would use the unbounded charge $Q_{\WW}$ to generate the extended structures, and would also require us to specify the extended domain of the unbounded $\drb$. For this reason, in the unbounded case we leave extensions similar to Example \eqref{ex:pscs} to a case by case basis.
\item An intrinsic characterization for $\cP_0^{BRST}$ without using $\Phi^{W}_s$ is possible, although more complicated than in the bounded case (cf.  Proposition \eqref{pr:bddinttriv}). We do not pursue it further however, as we do not need it to analyze the following examples.
\end{itemize}
\end{rem}

\subsection{Examples}\label{sbs:rigex}
We have now constructed the abstract structures associated to the BRST constraint process. Below we apply these to examples discussed so far: abelian Hamiltonian BRST for a finite number of constraints, and both versions of BRST-QEM. In particular we establish the connection between the physical algebra produced by the $T$-procedure (Dirac constraint procedure) and the BRST-physical algebra. 
\begin{itemize}
\item To describe the relations between the different algebras below, we will state explicitly the nature of the homomorphisms, \ie if they are algebra homomorphisms, $*$-homomorphisms, $\dag$-homomorphisms. The symbol `$\cong$' will denote isomorphism below but does not assume any adjointness property, \eg $*$-isomorphism. Any such property of the associated the isomorphism will explicitly stated along with the identity containing `$\cong$'.
\item  We will use the following basic result frequently (cf.\cite{Mur90}[Theorem 3.1.5 p80]:
\begin{theorem}
Let $\cA$ and $\cB$ be $C^*$-algebras. Let $\vp:\cA\to \cB$ be an $*$-homomophism with $\ker \vp=\{0\}$. Then $\vp$ is isomometric, \ie $\cA \cong \cB$.
\end{theorem}
\end{itemize}

\pfit \textbf{Example 1:} 

For bounded $Q$ let $D_1(\drb)=D_2(\drb)=\cA$, and let $\fS_W=\fS_{\drb}$. The definitions for bounded BRST in Section \ref{sbs:bddcbrst} agree with the unbounded case above. Recall abelian Hamiltonian BRST for a finite set of constraints in Example \eqref{ex:csfham}. For this case we get by lemma \eqref{lm:abcsst} that $\fS_{\drb}= \fS_D \otimes \fS_g$ and  $P_s=(\one - P)\otimes \one$ hence,
\[
\cP^{BRST}\cong \Phi^{\drb}_{s}(\cA_0\otimes \cA_g)=\big([(\one -P)\cA_0 (\one -P) ]\otimes \cA_g \big)\cap \ker \drb
\]
where the above isomorphism is a $\dag$-isomorphism as given by Proposition \eqref{pr:ubddcshinv} (i) but not in general a $*$-isomorphism as $\ker \drb$ is not a $*$-algebra in general. Theorem \eqref{th:ubdbrststruc} (iii) gives that $\cP^{BRST}$ is a $C^*$-algebra if $\Phi_s^{\WW}(\ker \drb)= \Phi_s^{\WW}((\ker \drb)^*)$. In Example \eqref{ex:pscs} we see that we can satisfy this condition by extending $\cA$ to $\tilde{A}=\osalg{ \{P_s\} \cup \cA}$. 

Alternatively we can restrict to the original algebra and get the traditional Dirac observables as in Proposition \eqref{pr:rabbdbral}, that is $\Phi_s^{\WW}(\ker \drb\cap (\cA_0 \otimes \one)) \cong (\cC'/(\cC'\cap \cD))\otimes \one$  where `$\cong$' denotes a $*$-isomorphism. 
 
\pfit \textbf{Example 2:} 

For BRST-QEM I in Section \ref{sbsc:cbrstv1} let $\cA=\cR(\fX,\sigma_1)\otimes \cA_g$ and $D_1(\drb), D_2(\drb)$ as in Definition \eqref{df:sddom1} and let $\fS_{\WW}=\fS_{\drb}$. We have by Proposition \ref{pr:rabrstd1} that $\fS_{\drb}=\fS_D\otimes \fS_g$,  where $\fS_D$ are the Dirac states of $\cR(\fD,\sigma_1)$ using $\fDL$ as a constraint test function space as in Subsection \ref{sbs:racons}. As discussed in Subsection \ref{sbs:problem} we have:
\begin{gather}
\pi_{\omega}\circ \drb=Q_{\omega}=0, \qquad P^{\omega}_s=\one, \qquad
\pi_{\omega}(D_2(\drb))=\pi_{\omega}(\cR_0(\fX_t, \sigma_1)\otimes \cA_g)\notag \\
\pi_{\omega}(D_2(\drb))=\pi_{\omega}(\cR_0(\fX_t, \sigma_1)\otimes \cA_g)\label{eq:qBR1stuff}
\end{gather}
for all $\omega \in \fS_{\drb}$ (cf. equation \eqref{eq:omg1} for equation \eqref{eq:qBR1stuff}). Using the facts we show that the BRST-physical algebra strictly contains the  Dirac physical observables and is contained in the Dirac physical observables tensored with the ghosts.
\begin{proposition}\label{pr:RA1phbalblah}
Let $\cP=\cO/\cD$ are the Dirac physical algebra from Subsection \ref{sbs:cnst1}. Then 
\[
\cP \cong \pi_{\WW}(\cR(\fX_t, \sigma_1)\otimes\one) \subset \overline{\Phi^{\WW}_s(\ker \drb)}^{B(\cH_\WW)} \subseteq \pi_{\WW}(\cR(\fX_t, \sigma_1)\otimes \cA_g) \cong \cP \otimes \cA_g
\]
where the containment $\pi_{\WW}(\cR(\fX_t, \sigma_1)\otimes\one) \subset \overline{\Phi^{\WW}_s(\ker \drb)}^{B(\cH_\WW)}$ is proper, and the tensor norm is unique as $\cA_g$ is nuclear. The above are $*$-isomorphisms.
\end{proposition}
\begin{proof}
As $P^{\omega}_s=\one$ for all $\omega \in \fS_{\WW}$ and $\cH_{\WW}=\oplus_{\omega\in \fS_{\WW}}$ we have that $\cH^{\WW}_s=\cH_{\WW}$ and so $P_s^{\WW}=\one \in B(\cH_{\WW})$. 

We show that $\overline{\Phi^{\WW}_s(\ker \drb)}^{B(\cH_\WW)} \subseteq \pi_{\WW}(\cR(\fX_t, \sigma_1)\otimes \cA_g) \cong \cP \otimes \cA_g$. From the discussion above $\pi_{\omega}\circ \drb=0$ for all $\omega \in \fS_{\WW}$, and as $\pi_{\WW}=\oplus_{\omega \in \fS_{\WW}}\pi_{\omega}$ we have
\[
\pi_{\WW}(\ker \drb )\subset \pi_{\WW}(D_2(\drb)).
\]
Also $P^{\WW}_s=\one$ and so
\begin{align}
P_s^{\WW} \pi_{\WW} (\ker \drb)P_s^{\WW}\subset \pi_{\WW}(D_2(\drb))= & \;\pi_{\WW}(\cR_0(\fX_t, \sigma_1)\otimes \cA_g),\notag \\
= & \;\pi_{\WW}(\cR_0(\fX_t, \sigma_1)\otimes \one)\pi_{\WW}(\one \otimes \cA_g) \label{eq:BRSIWhom}.
\end{align}
where we used equation \eqref{eq:qBR1stuff} in the second equality. Note also that the first containment need not be equality as $\ker \drb \subset D_2(\drb)$ where the containment is proper.

Recall from lemma \eqref{lm:resD2} the representation $\pi_D:\cR(\fX,\sigma_1) \to B(\cH_D)$ by: 
\[
\cH_{D}:=\bigoplus\{ \cH_{\omega}\,|\, \omega \in \fS_{D}\}, \qquad \pi_{D}:\bigoplus \{ \pi_{\omega}\,|\, \omega \in \fS_{D}\}.  
\]
By $\fS_{\WW}=\fS_{\drb}=\fS_D\otimes \fS_g$ we get,
\begin{align}
\overline{\pi_{\WW}(\cR_0(\fX_t, \sigma_1)\otimes \one)}^{B(\cH_\WW)}\cong & \;\overline{\pi_{D}(\cR_0(\fX_t, \sigma_1))}^{B(\cH_D)}\otimes \one, \notag \\
\cong& \cR(\fX_t, \sigma_1), \label{eq:homWP}\\
\cong& \cP \notag,
\end{align}
where the first isomorphism follows from lemma \eqref{lm:resD2} (ii) and the third from Proposition \eqref{pr:RA1phalg}. Also $\pi_{\WW}(\one \otimes \cA_g)\cong \cA_g$ as $\cA_g$ is a CAR algebra hence simple. Hence we have,
\begin{align*}
\overline{\Phi^{\WW}_s(\ker \drb)}^{B(\cH_{\WW})} =\overline{\pi_{\WW} (\ker \drb)}^{B(\cH_{\WW})} \cong &\;\overline{\pi_{\WW}(\cR_0(\fX_t, \sigma_1)\otimes \one)\pi_{\WW}(\one \otimes \cA_g)}^{B(\cH_{\WW})}\\
 \cong &\; \cR(\fX_t, \sigma_1)\otimes \cA_g\\
 \cong &\;\cP\otimes \cA_g
\end{align*}
where all the isomorphisms are $*$-isomorphisms, we used Theorem \eqref{th:ubdbrststruc} in the first isomorphism, $P^{\omega}_s=\one$ in the first equality and equation \eqref{eq:BRSIWhom}, equation \eqref{eq:homWP} and that the tensor norm is unique as $\cA_g$ is a nuclear.   

We now show that $\pi_{\WW}(\cR(\fX_t, \sigma_1)\otimes \one)\subset \overline{\Phi^{\WW}_s(\ker \drb)}^{B(\cH_{\WW})}$, where the containment is proper. By $P_2\fX_t=0$ and the definition  $\drb(R(\lambda,f))=\drb(R(\lambda, f))= i R(\lambda, f)^2 C( KP_{\JJ} f)$ (cf. Theorem \eqref{th:calgdrb}) we have that 
\[
\cR_0(\fX_t, \sigma_1)\otimes \one \subset \ker \drb
\]
hence $\pi_{\WW}(\cR(\fX_t, \sigma_1)\otimes \one)\subset \overline{\Phi^{\WW}_s(\ker \drb)}^{B(\cH_{\WW})}$ by equation \eqref{eq:homWP}. 

To see that the containment is proper consider $\zeta_1(h)^*:=R(-1, \mH h)\otimes C( Jh) \in D(\drb)$ for $h \in \fXL$. By the definition of $\drb$ in Theorem \eqref{th:calgdrb} we have that $\zeta_1(h)^* \in \ker \drb$. Now by $h\in \fXL$, lemma \eqref{lm:BRSTalgbas} (iii) and $\fS_{\WW}=\fS_{\drb}$, we have that $\pi_{\WW}(R(-1, \mH h))\otimes \one))=\pi_{\WW}(R(1, \mH h))\otimes \one))^* =i\one$. Hence 
\[
\pi_{\WW}(\zeta_1(h)^*)=\pi_{\WW}(R(-1, \mH h)\otimes \one )\pi_{\WW}(\one \otimes C( Jh))=i  \pi_{\WW}(\one \otimes C( Jh)) \notin \pi_{\WW}(\cR(\fX_t, \sigma_1),
\]
where we have used that $\pi_{\WW}(\one \otimes C( Jh))\neq 0$ as $\cA_g$ is simple. Hence the containment $\pi_{\WW}(\cR(\fX_t, \sigma_1)\otimes \one)\subset \overline{\Phi^{\WW}_s(\ker \drb)}^{B(\cH_{\WW})}$ is proper.
\end{proof}

Now by Proposition \eqref{pr:ubddcshinv},
\[
\cP^{BRST}\cong \overline{\Phi^{\WW}_s(\ker \drb)}^{B(\cH_\WW)}.
\]
by a $\dag$-isometric isomorphism, and so Propostion \eqref{pr:RA1phbalblah} therefore states that $\cP^{BRST}$ properly contains the Dirac observables obtained using the $T$-procedure $(\cR(\fX, \sigma_1), \cC_1)$ as in Subsection \ref{sbs:cnst1}, and is contained in the Dirac observables tensored with the ghost algebra. This shows rigorously that BRST using the Resolvent Algebra $\cR(\fX,\sigma_1)$ does \emph{not} give equivalent results to the Dirac algorithm ($T$-procedure) as it does not remove the ghosts.

By Proposition \eqref{pr:RA1keragcom} we can also easily encode Poincar{\'e} covariance in this picture.
\begin{proposition}\label{pr:RA1keragcom2}
Let $g\to \alpha_g\in \aut(\cA)$ be the representation of $\cP^{\uparrow}_{+}$ as in Proposition \eqref{pr:RA1keragcom}, \ie $\alpha_g$ has action
\begin{equation*}
\alpha_g(R(\lambda,f)\otimes C(h))=R(\lambda,V_gf)\otimes C(S_gh)
\end{equation*}
for $f\in \fX$, $h\in \fDL\oplus_{\mH}\fDJ $. Then
\begin{itemize}
\item[(i)] We have
\[
\{ \omega \circ \alpha_g \,|\, \omega\in \fS_{\drb},\,g \in \cP^{\uparrow}_{+}\}=\fS_{\drb}
\]
for all $g\in \cP^{\uparrow}_{+}$ there exists $\beta_g \in \Aut(\pi_{\drb}(\cA))$ such that 
\[
(\beta_g \circ \pi_{\drb})(A)=(\pi_{\drb}\circ \alpha_g)(A), \qquad \forall A \in \cA.
\]
Moreover, $g \to \beta_g$ is a representation $\cP^{\uparrow}_{+} \to \aut(\pi_{\drb}(\cA))$.
\item[(ii)] We have that $\alpha_g$ factors to $\hat{\alpha}_g \in \Aut(\cP_0^{BRST})$ such that
\begin{equation}\label{eq:facalg}
\hat{\alpha}_g\hat{A}=\widehat{\alpha_g(A)}
\end{equation}
for all $A \in \ker \drb$. This extends to a $\dag$-automorphism $\cP^{BRST}$, which we still denote $\hat{\alpha}_g$. Moreover
\[
g \to \alpha_g
\]
is a representation of $\cP^{\uparrow}_{+}$ in the $\dag$-automorphisms of $\cP^{BRST}$.
\end{itemize}
\end{proposition}
\begin{proof}
(i): Let $\omega \in \fS_{\drb}$ and $g \in \cP^{\uparrow}_{+}$. Then by Proposition \eqref{pr:RA1keragcom} we have $\alpha_{g^{-1}} \circ \drb \circ \alpha_g = \drb$ on $D_2(\drb)$ and $\alpha_g(D_2(\drb))=D_2(\drb)$, hence
\begin{equation}\label{eq:drbra1}
(\omega \circ \alpha_g) (\drb(A))=\omega (\drb ( \alpha_g (A)))=0
\end{equation}
for all $A \in D_2(\drb)$. Furthermore, by Proposition \eqref{pr:RA1keragcom} we have $\alpha_g\circ \alpha=\alpha \circ \alpha_g$ and hence 
\begin{equation}\label{eq:alpra1}
(\omega \circ \alpha_g) \circ \alpha = (\omega \circ \alpha) \circ \alpha_g =\omega \circ \alpha_g
\end{equation}
Obviouly $\omega\circ \alpha_g \in \fS_T=\{\omega_1 \otimes \omega_2\,|\, \omega_1 \in \fS(\cR(\fX,\sigma_1)),\, \omega_2 \in \fS(\cA_g)\}$, which combined with equation \eqref{eq:drbra1} and equation \eqref{eq:alpra1} shows that $\omega \circ \alpha_g \in \fS_{\drb}$ (cf. Definition \eqref{df:ubsgsbdd}). This shows 
\[
\{ \omega \circ \alpha_g \,|\, \omega\in \fS_{\drb},\,g \in \cP^{\uparrow}_{+}\}\subset\fS_{\drb}
\]
The reverse inclusion is obvious. Now $\pi_{\WW} =\oplus_{\omega \in \fS_{\WW}}\pi_{\omega}=\oplus_{\omega \in \fS_{\drb}}\pi_{\omega}$. As $\omega\circ \alpha_g \in \fS_{\drb}$ for all $\omega \in \fS_{\drb}$ we get that $\pi_{\WW}\circ \alpha_g$ is $\pi_{\omega}$ with the direct summands permuted. Such a direct summand can be done in $\cH_{\WW}$ by conjugation with a unitary, which we denote $\beta_g$. Conjugation by a unitary which preserves $\pi_{\WW}(\cA)$ is a $*$-automorphism we get that $\beta_g\in \Aut(\pi_{\WW}(\cA))$, and as $g \to \alpha_g$ and $\pi_{\WW}$ are representations, so is $g \to \beta_g$.

\pfit (ii): Let $g \in \cP^{\uparrow}_{+}$. By Proposition \eqref{pr:RA1keragcom} (i) we have $\alpha_g(\ker \drb)=\ker \drb$. As $P_s^{\WW}=\one$ and $\Phi_s^{\WW}(\cdot)=P_s^{\WW}(\cdot)P_s^{\WW}$ we have that $\ker \drb \cap \ker \Phi_s^{\WW}=\{0\}$. Hence $\alpha_g$ factors trivially to an automorphism $\hat{\alpha}_g$ on $\cP_0^{BRST}=\ker \drb/(\ker \drb \cap \ker \Phi_s^{\WW})$ with action given by equation \eqref{eq:facalg}.

Now by Proposition \eqref{eq:facalg} (iii), $\alpha_g(A^{\dag})=(\alpha_g \circ \alpha)(A^*)=(\alpha \circ \alpha_g)(A)^*=\alpha_g(A)^{\dag}$. Hence,
\[
\hat{\alpha}_g(\hat{A^{\dag}})=\widehat{\alpha_g(A^\dag)}=\widehat{\alpha_g(A)}^{\dag}
\]
where we have used that the $\dag$-invoution factors to a $\dag$-involution on $\cP_0^{BRST}$ as in Proposition \eqref{pr:ubddcshinv}. Hence $\hat{\alpha}_g$ is a $\dag$-involution.

We now show that $\hat{\alpha}_g$ extends to an $\dag$-automorhism on ${\cP^{BRST}}$. By (i), $\alpha_g (\ker \drb)=\ker \drb$ and the fact that $\Phi_s^{\WW}=\pi_{\WW}$ , we have that $\beta_g\in \aut(\pi_{\WW}(\cA))$ is such that 
\[
(\beta_g\circ\Phi_s^{\WW})(\ker \drb)=\Phi_s^{\WW}(\ker \drb)
\]
As $\beta_g$ is defined by a unitary conjugation, it is a $*$-automorphism on the $C^*$-algebra $\pi_{\WW}(\cA)$ and so isometric and so extends to an $\dag$-automorphism on $\overline{\Phi_s^{\WW}(\ker \drb)}^{\cH_{\WW}}$ but not a $*$-automorphism as $  \ker \drb$ is not a $*$-algebra. By definition $\cP^{BRST} \cong \overline{\Phi_s^{\WW}(\ker \drb)}^{\cH_{\WW}}$ and so if $\tau$ is this isomorphism, then $\beta_g \circ \tau^{-1}$ defines an automorphism on $\cP^{BRST}$. It is easy to check that this agrees with $\hat{\alpha_g}$ on $\cP_0^{BRST}$ and so $\beta_g\circ \tau^{-1}$ is the extension of $\hat{\alpha_g}$ to $\cP^{BRST}$. 

Now $g \to \hat{\alpha}_g$ is a representation of $\cP^{\uparrow}_{+} \to \Aut(\cP^{BRST})$ which follows as $g \to \beta_g$ and $\tau^{-1}$ are. 
\end{proof}
Proposition \eqref{pr:RA1keragcom2} (ii) encodes the Poincar{\'e} transformations on the BRST physical algebra $\cP^{BRST}$.

\pfit \textbf{Example 3:} 

For BRST-QEM II in Section \ref{sec:cstbrstii} we take $\cA=\cR(\fD,\sigma_2)\otimes \cA_g$ and $D_1(\drb), D_2(\drb)$ as in Definition \eqref{df:sddom2}. Recall from Definition \eqref{df:ubsgsbdd} that we needed $\pi_{\omega}(D_1(\drb))\Omega_{\omega}$ and $\pi_{\omega}(D_2(\drb))\Omega_{\omega}$ to be dense in $\cH_{\omega}$ for $\omega \in \fS_{\WW}$. We now show that if $\omega \in \fS_{T}(\cA)$ is regular on the resolvent part of $\cA$ it has this property.
\begin{lemma}\label{lm:domdel2}
Let $\omega= \omega_1\otimes \omega_2 \in \fS_{T}(\cA)$ and $\omega_1 \in \Reg(\cR(\fD, \sigma_2))$ (cf. Definition \eqref{df:regrep}). Then $\pi_\omega(D_1(\drb))\Omega_{\omega}$ is dense in $\cH_\omega$. 
\end{lemma}
\begin{proof} Recall from Definition \eqref{df:sddom2} that for all $g \in \fDL$,  $\zeta_1(g),\zeta_1(g)^* \in D_1(\drb)$ where
\[
\zeta_1(g)=R(1, \mH P_{+}g)R(1,K\mH P_{-}g)\otimes C(g) \qquad \text{and} \qquad \zeta_1(g)^*=R(-1, \mH P_{+}g)R(-1,K\mH P_{-}g)\otimes C(Jg),
\]
and furthermore recall the defining property $R(\lambda,f)=(1/ \lambda)R(1,f /\lambda)$ (cf. Definition \eqref{df:RA} \eqref{eq:R3}). By Theorem \eqref{RegThm} (ii), we have that for $\psi \in \pi_{\omega}(D_1(\drb))\psi$, $g \in \fDL$
\begin{align*}
\pi_{\omega}(C(g))\psi&=-\lim_{\lambda \to \infty} \lambda^2 \pi_{\omega}(R(\lambda, \mH P_{+} g)\pi_{\omega}(R(\lambda, K \mH P_{-} g)\pi_{\omega}(C(g))\psi,\\
& =\lim_{\lambda \to \infty} \lambda \pi_{\omega}(R(1, \mH P_{+} g/\lambda)\pi_{\omega}(R(1, K \mH P_{-} g/\lambda )\pi_{\omega}(C(g/\lambda ))\psi,\\
&=\lim_{\lambda \to \infty}\lambda \pi_{\omega}(\zeta_{1}(g/\lambda))\psi,
\end{align*}
and so $\pi_\omega(C(g))\psi \in \overline{ \pi_{\omega}(D_1(\drb))\psi}$. 

A similar argument using $\lambda \pi_{\omega}(\zeta_{1}(g/\lambda)^*)$, gives that for $g \in \fDL$, $\pi_\omega(C(J g))\psi \in \overline{ \pi_{\omega}(D_1(\drb))\Omega_{\omega}}$  and hence $\pi_\omega(C(Jg))\psi \in \overline{\pi_{\omega}(D_1(\drb))\Omega_{\omega}}$. 

In fact we can use the same argument to show that for any finite set $(g_i)_{i=1}^{n}\subset \fDL \cup \fDJ$, we have,
\[
\pi_{\omega}(C(f_1)\ldots C(f_n))\pi_{\omega}(D_1(\drb))\Omega_{\omega} \subset \overline{ \pi_{\omega}(D_1(\drb))\Omega_{\omega}},
\]
but 
\begin{align*}
[ \pi_{\omega}(C(f_1)\ldots C(f_n))&\pi_{\omega}(D_1(\drb))\,|\,  (f_i) \;\text{is a finite subset of $\fDL\cup \fDJ$}],\\
 &=[\pi_{\omega}(C(f_1)\ldots C(f_n)D_1(\drb))\,|\,  (f_i) \;\text{is a finite subset of $\fDL\cup \fDJ$}],\\
&=\pi_{\omega}(\cA_0),
\end{align*}
and so,
\[
\cH_{\omega}=\overline{\pi_{\omega}(\cA_0)\Omega_{\omega}}\subset \overline{ \pi_{\omega}(D_1(\drb))\Omega_{\omega}} \subset \cH_{\omega}.
\]
 \end{proof}
Hence we take for $\fS_{\WW}$ the set:
\[
\fS_{\WW}=\fS_{\WW r}:=\{ \omega=\omega_1 \otimes \omega_2\in \fS_{\drb} \,|\, \omega_1 \in \Reg(\cR(\fD, \sigma_2)) \} \subset \fS_{\drb},
\]
where $\fS_{\drb}$ is as in Definition \eqref{df:ubsgsbdd}. 

To calculate the BRST physical algebra explicitly as in Subsection \ref{sbs:srs} we  restrict to BRST ground states that are strongly regular states on $\cR(\fD,\sigma_2)$, ie,
\[
\fS_{\WW}=\fS_{\WW sr}:=\{ \omega=\omega_1 \otimes \omega_2\in \fS_{\drb} \,|\, \omega_1 \in \SReg(\cR(\fD, \sigma_2)) \}
\]
and restrict $\drb$ to the domains $\tilde{D_j}(\drb):=D_j(\drb)\cap (\cR^{0}_t\otimes \cA_g)$, for $j=1,2$ where $\cR^{0}_t\cong \cR_0(\fD_t,\sigma_2)\otimes \cR_0(\fDL\oplus \fDJ,\sigma_2)$ as in equation \eqref{eq:ders}. To use the domains $\tilde{D_1}(\drb)$ and $\tilde{D_2}(\drb)$ in the constructions above we prove.
\begin{lemma}
Let $\omega= \omega_1\otimes \omega_2 \in \fS(\cA)$ and $\omega_1 \in \SReg(\cR(\fDJ, \sigma_2))$. Then $\pi_\omega(\tilde{D_1}(\drb))\Omega_{\omega}$ is dense in $\cH_\omega$ and
\[
\pi_\omega(\tilde{D_1}(\drb))\Omega_{\omega}\subset\pi_\omega(\tilde{D_2}(\drb))\Omega_{\omega}\subset\pi_\omega(\cA)\Omega_{\omega}.
\]
\end{lemma}
\begin{proof}
The proof holds as for lemma \eqref{lm:domdel2}, replacing ${D_j}(\drb)$ by $\tilde{D_j}(\drb)$ for $j=1,2$.
\end{proof}
We can now calculate the physical algebra.
\begin{proposition} Let $\cA$, $\fS_{\WW sr}$ be as above for BRST-QEM II and let $\drb$ have domain $\tilde{D_2}(\drb)\supset \tilde{D_1}(\drb)$. Then
\begin{itemize}
\item[(i)] The physicality condition $J_{\WW sr}P_s^{\WW sr}=P_s^{\WW sr}$ is satisfied. Hence $\cP^{BRST}$ is a $C^*$-algebra and the $\dag$-involution coincides with the $*$-involution on $\cP^{BRST}$.

\item[(ii)] We have the $*$-isomorphisms:
\[
\cP^{BRST}\cong \cR(\fD_t,\sigma_2)\cong \cP,
\]
where $\cP=\cO/\cD$ is the Dirac physical algebra of observables for $(\cR(\fD,\sigma_2),\cC_2)$ as in Subsection \ref{sbs:qemdii}.
\end{itemize}
\end{proposition}
\begin{proof} Recall $\cR^0_{ph}= \salg{ R(\lambda,f)\,|\, f \in \fD_t}$ and $\cR_{ph}=\osalg{ R(\lambda,f)\,|\, f \in \fD_t}$. 

\pfit(i): Let $A \in \cR_{ph}\otimes \one$, then by $\alpha(R(\lambda,f)\otimes \one)=R(\lambda,Jf)\otimes \one$, and by $J_{\fD_t}=\one$ (Proposition, hence we have $\alpha(A)=A$. Now for $\omega \in \fS_{\WW sr}$ we have $\cH_s^{\omega}=\overline{\pi_{\omega}(\cR_{ph}\otimes \one)\Omega_{\omega}}$ (cf. Proposition \eqref{pr:kersrRph} (ii)),  hence 
\[
J_{\omega}\pi_{\omega}(A)\Omega_{\omega}=\pi_{\omega}(\alpha(A))\Omega_{\omega}=\pi_{\omega}(A)\Omega_{\omega},
\]
and so $J_{\omega}|_{\cH_s^{\omega}}=\one$. By Definition \eqref{df:ubdspdec} (i) and (iv) we have
\[
\cH^d_{\WW}\oplus \cH^s_{\WW}\oplus \cH^p_{\WW}=\cH_{\WW}=\bigoplus_{\omega \in \fS_{\WW sr}} \cH_{\omega}=\bigoplus_{\omega \in \fS_{\WW sr}} \cH^d_{\omega}\oplus \cH^s_{\omega}\oplus \cH^p_{\omega}
\]
from which it follows that $P_s^{\WW sr}|_{\cH_{\omega}}= P_s^{\omega}$ for all $\omega \in \fS_{\WW sr}$. Hence we have that $J_{\WW sr}P_s^{\WW sr}|_{\cH_{\omega}}=J_{\omega}|_{\cH_s^{\omega}}=\one$ hence $J_{\WW sr}P_s^{\WW sr}=\one$, \ie the physicality condition holds. The rest of the statements in (i) follow from Theorem \eqref{th:ubdbrststruc} (iii).

\pfit (ii): Theorem \eqref{th:brstpalgra2} (i) gives that 
\[
(\Phi^{\omega}_s \circ \pi_{\omega})(\cR^{0}_{t} \otimes \cA_g) =(\Phi^{\omega}_s\circ\pi_{\omega})(\cR^{0}_{ph}\otimes \one)
\]
for all $\omega \in \fS_{\WW sr}$.  Recall the definition
\[
\Phi^{\WW sr}_{s}(A)=P_s^{\WW sr} \pi_{\WW sr}(A)P_s^{\WW sr}, \qquad A \in \ker \drb.
\]
Combining $P_s^{\WW sr}|_{\cH_{\omega}}= P_s^{\omega}$ for all $\omega \in \fS_{\WW sr}$, $\pi_{\WW sr}=\oplus_{\omega \in \fS_{\WW sr}}\pi_{\omega}$ and $\tilde{D_2}(\drb)\subset \cR^{0}_{t} \otimes \cA_g$ give:
\begin{equation*}
\Phi^{\WW sr}_{s}(\cR^0_{ph} \otimes \one)\subset \Phi^{\WW sr}_{s}(\ker \drb \cap \tilde{D_2}(\drb)) \subset \Phi^{\WW sr}_{s}(\cR^{0}_{t} \otimes \cA_g)=\Phi^{\WW sr}_{s}(\cR^0_{ph} \otimes \one).
\end{equation*}
So by Theorem \eqref{th:ubdbrststruc} (ii), we get
\begin{equation}\label{eq:algisoBRS2}
\cP^{BRST}_0\cong \Phi^{\WW sr}_{s}(\ker \drb \cap \tilde{D_2}(\drb)) =\Phi^{\WW sr}_{s}(\cR^0_{ph} \otimes \one),
\end{equation}
where the above is a $*$-isomorphism as by (i) the $\dag$-involution coincides with the $*$-involution on $\cP^{BRST}_0$.

Now define
\[
\tilde{\Phi}^{\WW sr}_{s}:\cA\to \cA \qquad \text{by} \qquad \tilde{\Phi}^{\WW sr}_{s}(A):=P_s^{\WW sr}\pi_{\WW sr}(A)P_s^{\WW sr}
\]
\ie we have extended the domain of $\tilde{\Phi}^{\WW sr}$ to all of $\cA$. Note that $\tilde{\Phi}_s^{\WW sr}$ is a linear transform but \emph{not} a homomorphism on $\cA$, but that it is a $*$-homomorphism on $\pi_{\omega}(\cR_{ph}\otimes \one)$ as $\cH_{\omega}^s=\overline{\pi_{\omega}(\cR_{ph}\otimes \one)\Omega_{\omega}}$ for all $\omega \in \fS_{\WW sr}$. 

We want to show $\tilde{\Phi}^{\WW sr}(\cR_{ph} \otimes \one) \cong \cR_{ph}$. Let $\omega= \omega_1 \otimes \omega_2 \in \fS_{\WW sr}$, then by Theorem \eqref{th:brstpalgra2} (ii) we have a $*$-isometric isomorphism
\[
(\tilde{\Phi}^{\omega}_s\circ\pi_{\omega})(\cR_{ph}\otimes \one) \cong \pi_{(\omega_1,\cR_{ph})}(\cR_{ph}),
\]
where $\pi_{(\omega_1,\cR_{ph})}$ is the GNS-representation for $\omega_{1,\cR_{ph}}$. As $\omega_1|_{\cR_{ph}}$ is also strongly regular we have that $\pi_{(\omega_1,\cR_{ph})}$ is faithful. Hence 
\[
0=A \in \cR_{ph} \Leftrightarrow (\tilde{\Phi}^{\omega}_s\circ\pi_{\omega})(A\otimes \one)=0\: \forall \omega \in \fS_{\WW sr} \Leftrightarrow   \tilde{\Phi}^{\WW sr}(A)=0
\]
where the last equivalence is clear from $P_s^{\WW sr}|_{\cH_{\omega}}= P_s^{\omega}$ for all $\omega \in \fS_{\WW sr}$, $\pi_{\WW sr}=\oplus_{\omega \in \fS_{\WW sr}}\pi_{\omega}$. Hence $\tilde{\Phi}^{\WW sr}|_{\cR_{ph}\otimes \one}$ is faithful and so 
\begin{equation}\label{eq:algisoBRS22}
\cR_{ph} \cong (\cR_{ph}\otimes \one) \cong \Phi^{\WW sr}(\cR_{ph}\otimes \one),
\end{equation}
where the above isomorphisms are isometric. Combining equation \eqref{eq:algisoBRS22} and equation \eqref{eq:algisoBRS2} gives that
\[
\cR^0_{ph}\cong \cP^{BRST}_0
\]
where the above is a $*$-isomorphism between $C^*$-algebra and hence isometric. Taking closures gives,
\[
\cP^{BRST} \cong \cR_{ph}\cong \cP
\]
where the last isomorphism comes from Proposition \eqref{pr:physalgQEMii}.
\end{proof}

Poincar{\'e} covariance is a problem in this example. As in Subsection \ref{sbs:rescov}, we do not have that the Poincar{\'e} transformations define naturally on the auxiliary algebra $\cR(\fD,\sigma_2)$, so it is not clear how to define them on $\cA$ either. We do have the resolution to the problem that the final BRST-physical algebra is the same as that using the $T$-procedure where the Poincar{\'e} transformations are naturally defined, as discussed in Subsection \ref{sbs:rescov}. This is solution is in some sense the converse of a usual raison d'etre of BRST, that is that BRST is `manifestly covariant' while the Dirac method is not.  

Summarising, the BRST physical algebra is equivalent to the 
$T$-procedure  for \\
$(\cR(\fD, \sigma_1), \cC_1)$ or $(\cR(\fD, \sigma_2), \cC_2)$ in Subsection \ref{sbs:cnst1} Subsection \ref{sec:cstbrstii}, and so we can say that $C^*$-BRST for QEM using the auxiliary $\cR(\fD,\sigma_2)$ gives equivalent results to the Dirac method, but is not `manifestly covariant'.

\chapter{Results and Conclusion}

\section{Summary of Results}
A brief summary of the results in each chapter are as follows:

\pfit \textbf{Chapter \ref{ch:heu}:}  The heuristic BRST structures were described and the examples of BRST-QEM and Hamiltonian BRST with constraints that close were given. The Multiple Copies of the Physical Subspace (MCPS) problem was discussed for Hamiltonian BRST. In this introductory chapter there were no original results.

\pfit \textbf{Chapter \ref{ch:GenStruct}:} The heuristic structures in Chapter \ref{ch:heu} were made rigorous. The BRST charge $Q$ was analysed and the \emph{dsp}-decomposition was given for an unbounded $Q$ acting on a Krein space (as first proved in Horuzhy \cite{HoVo89}). The BRST physical subspace $(\cH^{BRST}_{phys}, \iip{\cdot}{\cdot}_p)$ was defined and shown to be a Krein space. The condition for physicality was given in equation \eqref{eq:spphcond} which implied that the Hilbert and Krein structure on $\cH^{BRST}_{phys}$ coincided.

The BRST superderivation was analysed. The basic spatial, algebraic and $\Z_2$-grading structures required to define $\drb$ was described at the level of algebras of unbounded operators acting on a dense invariant domain in a Krein space. Topological issues on the algebras were postponed till Chapter \ref{ch:CsBRST}. The BRST physical algebra $\cP^{BRST}$ was defined and its natural representation on $\cH^{BRST}_{phys}$ given. The connection between $\cP^{BRST}$ and the \emph{dsp}-decomposition was given in Theorem \eqref{pr:krdel}. This fundamental result can be found in the literature in \cite{Hen1989} p285, however the rigorous statement and proof for an infinite dimensional Hilbert space with unbounded BRST charge $Q$ is original. 

As a consequence of Theorem \eqref{pr:krdel} we showed that in a simple example with a single constraint, Hamiltonian BRST does \emph{not} remove the ghosts at a spatial or algebraic level in the final constrained theory, and that extra ghost number zero conditions do not fix the problem. This example and problem were communicated to me by Dr. Hendrik Grundling, however the proof via Theorem \eqref{pr:krdel} is original and shows the connection between the non-removal of the ghosts at the spatial and algebraic level.

\pfit \textbf{Chapter \ref{ch:BRSTQEM}:} We developed a rigorous model for the heuristic BRST-QEM model given in \cite{KuOj79,Schf2001}. We analysed the QEM test function space and defined the abstract \KOB test function space with analogous structures. We constructed the BRST superderivation $\drb$ and showed in the case of QEM that it gave the correct smeared version of the heuristic superderivation. We then constructed the BRST charge $Q$ for \KOB, calculated $\cH^{BRST}$ and $\cP^{BRST}$ using Theorem \eqref{pr:krdel}. We found that the \KOB model and Gupta-Bleuler model for QEM gave equivalent results, at both the spatial and algebraic level. The calculation of $\cP^{BRST}$ using Theorem \eqref{pr:krdel} is to the author's knowledge an original result. 

Using \KOB and Theorem \eqref{pr:krdel} we were also easily able to calculate the BRST physical subspace and algebra for the examples of BRST with a finite number of bosonic constraints, and BRST for massive abelian gauge theory. 

We answered in the affirmative the conjecture at the end of \cite{HoVo92} as to whether the BRST-physical state space can be used to calculate the BRST physical algebra more efficiently that by direct algebraic calculation (cf. Remarks \eqref{rm:fdkob}(i), \eqref{rm:ext}(ii)). 

Lastly, we synthesized Hamiltonian BRST with a finite set of commuting selfadjoint constraints with \KOB with a finite number of bosonic constraints, to get an abstract algorithm BRST for a finite set of commuting selfadjoint constraints that selects the correct physical subspace without the need for extra selection criteria. That is, it does not suffer the MCPS problem of the usual Hamiltonian algorithm and for simple examples, $\cP^{BRST}$ for this algorithm also coincided with the quantum Dirac physical observable algebra. This combined algorithm is an original construction.

\pfit \textbf{Chapter \ref{ch:CsBRST}:} We developed a $C^*$-algebraic framework for the structures in Chapters \ref{ch:GenStruct} and \ref{ch:BRSTQEM}. We did this by first investigating the case of bounded BRST charge $Q$. We used the selection criteria as in \cite{Hendrik1991} to give representations with all the structures of Chapter \ref{ch:GenStruct}, then defined $\cP^{BRST}$ as an abstract Banach algebra. We encoded \KOB in a $C^*$-algebraic form. 

To deal with issues related to the unboundedness of the fields we used the Resolvent Algebra to construct a `mollified' version of the BRST superderivation $\drb$. The first attempt at this gave a BRST-model on which we could encode the Poincar{\'e} transformations, but did \emph{not} give equivalent results to the quantum Dirac constraint procedure ($T$-procedure). The second attempt produced a BRST-model with BRST-physical algebra equivalent to that selected by the $T$-procedure, but did not admit a natural encoding of the Poincar{\'e} transformations. 

Finally, we developed an abstract $C^*$-algebraic framework for BRST for the case of an unbounded superderivation $\drb$ with non-norm dense domain. We show that the examples of Hamiltonian BRST with a finite number of constraints and $C^*$-\KOB Abelian (both versions) fit into this framework. The correspondence between the abstract BRST method and the $T$-procedure for the examples is: Hamiltonian BRST for a finite number of constraints suffers the MCPS problem and is not equivalent to the $T$-procedure; $C^*$-\KOB with covariant symplectic space $(\fD, \sigma_1)$ does not give equivalent results to the quantum Dirac method as the ghosts are not removed in $\cP^{BRST}$; $C^*$-\KOB with the auxiliary symplectic space $(\fD, \sigma_2)$ gives equivalent results to the $T$-procedure. 

The results regarding both versions of $C^*$-\KOB are original but with the construction of the mollified version of $\drb$ is heavily influenced by \cite{HendrikBuch2006}. The abstract construction of $\cP^{BRST}$  for general unbounded $\drb$, and its comparison  with the physical observable algebra selected by the $T$-procedure in the given examples is original.

\section{Conclusions}

We have now analysed the quantum BRST method of constraints in a well-defined mathematical framework. This was first done in the setting of a concrete Krein-space representation where the model independent structures of quantum BRST common to the standard BRST examples in the literature have been defined. The standard quantum BRST examples have then been developed rigorously in light of these frameworks, enabling a discussion of problematic issues related to BRST in the literature. 

The analysis the quantum BRST constraint method was then extended to a $C^*$-algebraic setting. The construction of this setting was one of the main goals of this thesis. This construction enabled a mathematically rigorous comparison of the results produced by the BRST method and the quantum Dirac constraint method ($T$-procedure), which was the other main aim of the thesis. We found that the results obtained from the different constraint methods were \emph{not} equivalent for the examples of Hamiltonian BRST with a finite number of constraints that close, and BRST-QEM using the covariant Resolvent Algebra $\cR(\fX,\sigma_1)$.

This is not to say that a $C^*$-algebraic framework is the only way to analyse quantum BRST rigorously, indeed there are many other approaches to rigorous quantum BRST as mentioned in the introduction. The $C^*$-algebraic viewpoint, however, was a methodology to draw together and analyse in a general and consistent mathematical framework examples coming from the different varieties of quantum BRST found in the literature. 

Even for basic examples, we found that a consistent treatment of quantum BRST was not straightforward. From the analysis, we draw several important conclusions:
\begin{itemize}
	\item The quantum BRST method and quantum Dirac method of constraints are \emph{not} equivalent in general. This has been seen in the rigorous examples given Subsection \ref{sbs:rigex}.  
	\item Quantum Hamiltonian BRST always suffers the MCPS problem and needs extra selection criteria to select the correct physical space. In simple examples, the MCPS problem also leads to the non-removal of the ghosts in the BRST physical algebra $\cP^{BRST}$. 
	\item Conversely, \KOB selects the correct physical algebra without extra selection criteria both at the level of unbounded operators acting on a Krein space, and at the $C^*$-algebraic level when using the auxiliary test function space $(\fD,\sigma_2)$. As \KOB is a rigorous example of Lagrangian BRST, we see that quantum Lagrangian  \KOB and quantum Hamiltonian BRST are \emph{not} equivalent constraint methods.
	\item Quantum BRST does not incorporate equivalent constraints well (cf. Remark \eqref{rm:wkcom}).
\end{itemize}

\section{Further Issues with Quantum BRST}

The inequivalence of quantum Hamiltonian BRST, Lagrangian BRST, and the Dirac method is an interesting topic in need of further investigation, particularly as all are equivalent at the classical level. Characterization of when the quantum equivalence of the different methods will hold has not been analysed in this thesis. This is primarily because we have found no general algorithm that encompasses all the varieties of methods labelled as `quantum BRST'.

The Lagrangian approach is used to model heuristically the important physical examples of quantum gauge theories, but is the least well understood mathematically. The $\drb$ and $Q$ are defined by `replacing the gauge parameter by ghost parameter'. As already stated, this is a vague concept and relies on gauge theory structures being present, rather than being an explicit algorithm such as the $T$-procedure which begins with a field algebra $\cA$ and a set of constraints $\cC \subset \cA$ that ultimately select the physical observable algebra. In fact in KO, the classical constraint equation $\partial^{\mu}A_{\mu}+\alpha_0B=0$ (cf. \cite{KuOj79} p14) is given but the quantum constraint set of the original unextended system are never explicitly stated. Natural questions to ask are: How do we do Lagrangian BRST-QEM with Gupta-Bleuler constraints? Or Coulomb constraints? Will they give equivalent results? To answer these questions, and the broader equivalence issues discussed above, we need to formulate rigorously the method of `replacing the gauge parameter by ghost parameter' in a way that produces the heuristic Lagrangian BRST structures in examples but is also valid when no gauge theory is present.

Related to the issue above is to what extent we can make other examples of heuristic Lagrangian BRST rigorous. We have done so for \KOB, but the main purpose for the use of BRST in \KO is to apply in the case of non-abelian quantum gauge theories (NAQGT). Issues to address with respect to Lagrangian BRST for NAQGT's are: will the physical subspace always be positive as in the abelian case; Will extra constraint conditions beyond $\drb$ and $Q$ be needed to select the correct BRST-physical space and algebra which were not needed in the abelian case, and if so what extra conditions? Resolving these issues, however, is very difficult as constructing a rigorous realistic NAQGT is a long standing open problem. 

We should note that there is a standard heuristic argument that BRST also gives the correct results for NAQGT based upon the `quartet mechanism" and asymptotic abelianess \cite{KuOj79} p46-47. Whether this can be made rigorous is difficult to answer given the lack of rigorous NAQGT's, but there is evidence that extra selection conditions beyond $\drb$ and $Q$ will be needed in the non-abelian case \cite{St1989, AbNak1996}. While very interesting and important to address, these issues are beyond the scope of this thesis.

A further issue requiring attention is the matter of equivalent quantum constraints and quantum BRST. It is easy to construct examples of different constraint sets of operators that select the same physical states, but have different commutants. This means that the corresponding traditional Dirac observables (\ie the commutant of the constraints) will be different although the physical states selected are the same. This issue is resolved by the $T$-procedure, which defines the observable as the abstract version of the \emph{weak commutant} of the constraints (cf. \cite{Hendrik2006} p100), and is a significant advantage which this generalized Dirac method of constraints has over the traditional Dirac method. In Remark \eqref{rm:wkcom} we see that quantum Hamiltonian BRST suffers the same problem with respect to equivalent constraints and poses the question of how to formulate a generalized version quantum Hamiltonian BRST. Whether a weak version of Lagrangian BRST is necessary is a difficult question to answer. Evidence against is that in both versions of \KOB for QEM, we found that the BRST method selected the same part of the Resolvent Algebra as the $T$-procedure. However, in the QEM case the $T$-procedure also selects the traditional Dirac observables. Hence to investigate the equivalent constraint issue for Lagrangian BRST we come back to need for a rigorous Lagrangian BRST algorithm.

\section{Directions for Further Analysis}
The most important task in finding a mathematically transperant understanding of the quantum BRST method is stating a well-defined algorithm that encodes Lagrangian BRST for a general quantum gauge theory as described by \KO \cite{KuOj79} . Until this is done, it is very difficult to resolve any of the issues of equivalence, extra selection conditions, positivity of the physical subspace, etc. as described above except on a case by case basis. But to complete this task, the author feels we must state what a quantum gauge theory is rigorously. As already stated, this is an open problem that is very difficult to solve, but should be first completed before we can feel confident that we fully understand the quantum BRST algorithm. 
A direction to take in doing this would be to further investigate the relationship between the structures developed in this thesis and the rigorous work related to PGI such as found in \cite{Schf2001,DuSh1999,DuFred1999,DuFred1998,Holl2008}, as this is developed in a Lagrangian BRST context and gives examples beyond BRST-QEM that would shed more light on the form of a general Lagrangian BRST constraint algorithm.

On a less grand scale, an issue that needs resolution is the connection between the Resolvent Algebra with covariant symplectic space $(\fD, \sigma_1)$ and the Resolvent Algebra with auxiliary test function space $(\fD,\sigma_2)$. As both symplectic forms have act on the same vector space $\fD$ and are related by the symplectic operator $\Sp(\fD,\sigma_1)\ni J \in \Sp(\fD,\sigma_2)$, it seems that they should encode the same information. However, the $C^*$-\KOB models constructed with the different Resolvent Algebras gave different results, \ie the $\cR(\fD,\sigma_1)$-model selected the wrong physical algebra but naturally admitted an encoding of the Poincar{\'e} transformations, while the $\cR(\fD,\sigma_2)$-model selected the correct physical algebra but did not admit an encoding of the Lorentz boosts. Hence understanding the correspondence is worth investigating. It will also have relevance in other areas such as supersymmetry models constructed as in \cite{HendrikBuch2006}, but for cases where the bose fields are Krein-symmetric gauge fields rather than Hilbert essentially-selfadjoint scalar fields.

With respect to quantum Hamiltonian BRST an interesting area to investigate is how to generalize the algorithm in a way that takes into account equivalent sets of constraints. It is not obvious to the author that the cohomological definition of the BRST physical algebra $\cP^{BRST}=\ker \drb/ \ran \drb$ can be generalized in the direction of using a `weak' version of $\ker \drb$ and $\ran \drb$ similar to the weak commutant of the constraints in the $T$-procedure.

A final direction to follow is to see to what extent we can extend the synthesized Hamiltonian BRST and \KOB algorithm in Section \ref{sec:Bsfinabhm}. We have constructed a general quantum BRST algorithm for the case corresponding to a finite set of selfadjoint commuting constraints such that $Q$ selects the correct physical subspace with no extra selection conditions (\ie it does not suffer the MCPS problem of Hamiltonian BRST). It would be extremely useful to extend this to the case where we have a finite set of constraints that close but do not commute. Or to a general set of non-commuting constraints. We would also like to extend the algorithm to infinite sets of constraints. This would be a valuable area to investigate as results would lead towards a general quantum BRST algorithm that removes the ghosts and selects the  correct physical objects with no extra selection conditions.

\chapter{Appendix}

\section{Superstuff}\label{ap:SS}
Let $\cA=\cA_{+}\oplus \cA_{-}$ be an algebra with $\Z_2$-grading where $\oplus$ denotes algebraic sum. We call $\cA_{+}$ the \emph{even} part of $\cA$ and $\cA_{-}$ the \emph{odd} part of $\cA$.
Define
\[
\epsilon_A:= 
\begin{cases}
0, \qquad \text{for $A\in \cA_{+}$}\\
1, \qquad \text{for $A\in \cA_{-}$}
\end{cases}
\]
Then we can define graded brackets on $\cA$ with the following properties:
\begin{itemize}
\item Super Bracket
\[
\sbr{A}{B}=AB-(-1)^{\epsilon_A \epsilon_B} BA
\]
Note that,
\begin{equation}\label{eq:sbrid1}
	\sbr{A}{\sbr{A}{A}}=0
\end{equation}
\item Superderivation
\begin{align*}
\sbr{A}{BC}= & \;\sbr{A}{B}C+(-1)^{\epsilon_A \epsilon_B}B \sbr{A}{C},\\
\sbr{BC}{A}= & \;B\sbr{C}{A}+(-1)^{\epsilon_A \epsilon_C} \sbr{B}{A}C,\\
\intertext{in particular, for $Q$ such that $\epsilon_Q=1$ and $\drb(A):=\sbr{Q}{A}$ we have}
\drb(BC)= & \;\drb(B)C+(-1)^{\epsilon_B}B\drb(C),
\end{align*}
\item Super Jacobi Identity
\[
(-1)^{\epsilon_A \epsilon_C}\sbr{A}{\sbr{B}{C}}+(-1)^{\epsilon_C \epsilon_B}\sbr{C}{\sbr{A}{B}}+(-1)^{\epsilon_B \epsilon_A}\sbr{B}{\sbr{C}{A}}=0,
\]
\end{itemize}

A preferred notation for $\Z_2$-graded algebras and graded brackets is the following. Define the grading automorphism on $\cA$ by
\[
\gamma(A_{+}+A_{-})=A_{+}-A_{-}
\]
for $A_{+}\in \cA_{+}$ and $A_{-}\in \cA_{-}$. Note that $\gamma^2 =\iota$. Then we can equivalently define the superbrackets
\[
\sbr{A}{B}=AB-\gamma(A)B, \qquad \forall A,B \in \cA.
\]
Using this we can restate all the properties of the superbrackets and superderivation above, such as the super Jacobi identity, using $\gamma$ notation instead on $\epsilon_A$ notation. In particular for $Q$ such that $\gamma(Q)=-Q$ we have
\[
\drb(AB)=\drb(A)B+\gamma(A)\drb(B) \qquad \forall A,B \in \cA.
\]

\section{Indefinite Inner Product Spaces}\label{ap:IIP}
This appendix gives some basic facts about indefinite inner product spaces relevent to the disussion in this thesis. For much more extensive developments see in particular \cite{Bog1974, AzIo1989}.

Let $\fD$ be a vector space with inner product $\iip{\cdot}{\cdot}$. A vector $f\in \fD$ such that $\iip{f}{f}$ is positive, negative, or zero, is called positive, negative or null. A subspace $\fL \subset \fD$ is called positive, negative, or null if the vectors in $\fL$ are positive negative or null. If $\fD$ is such that $\iip{f}{f}>0, \iip{f}{f}\geq 0, \iip{f}{f}<0,\iip{f}{f}\leq 0$ or is both positive and negative
then $\fD$ is called positive definite, positive semi-definite, negative, negative semi-definite or indefinite respectively.

Let $\fD_{+}=\{ f\in \fD \,|\, \iip{f}{f}>0 \}$, $\fD_{-}=\{ f\in \fD \,|\, \iip{f}{f}<0 \}$, $\fD_{n}=\{ f\in \fD \,|\, \iip{f}{f}=0 \}$, be the \emph{sets} of positive, negative and neutral vectors in $\fD$. If $\fD$ is indefinite then none of these sets are subspaces (\cite{Bog1974} Corollary 2.7 p7).

Let $\fD_0=\{ f \in \fD \,|\, \iip{g}{f}=0\, \forall g \in \fD\}$. Then $\fD_0$ is called the \emph{\textbf{isotropic}} subspace of $\fD$. If $\fD_0=\{0\}$ then $\fD$ is called \emph{\textbf{non-degenerate}}.

Now suppose that $\fDL,\fDJ \subset \fD$ are two subspaces such that $\iip{f_1}{f_2}=0$ for all $f_1 \in \fDL$, all $f_2 \in \fDJ$ and such that the sum of the two spaces is non-degenerate. Then it follows easily that the sum is a direct sum, which we denote by $\fDL [\oplus ] \fDJ$. 

We still have to discuss topology on $\fD$. We will only consider the special case when $\fD$ can be completed to a Krein space. More general cases are discussed in \cite{Bog1974}. Suppose that $\fD$ non-degenenerate and, 
\begin{equation}\label{eq:decks}
\fD=\fD_{1}[ \oplus ]\fDJ
\end{equation}
where $\fDL$ is a positive subspace and $\fDJ$ is a negative subspace. Let $P_{+}$, $P_{-}$ be the projections on $\fD_{1}$, $\fD_{2}$ respectively and let,
\[
J=P_{+}-P_{-}.
\]
We call $J$ the \emph{\textbf{fundamental symmetry}}  on $\fD$ corresponding to the decomposition \eqref{eq:decks} (\cite{Bog1974} p52), or the fundamental symmetry when the decomposition is understood.

Now we define,
\[
\ip{\cdot}{\cdot}:=\iip{\cdot}{J\cdot}.
\]
It is straightforward to check that $J^2=\one$, $\iip{Jf}{Jf}=\iip{f}{f}$, $\ip{Jf}{Jf}=\ip{f}{f}$.

Now $\ip{\cdot}{\cdot}$ is positive definite, and hence $\norm{f}:=\ip{f}{f}^{1/2}$ for $f \in \fD$ is a norm on $\fD$. Let $\cH$ be the completion of $\fD$ with respect to $\norm{\cdot}$. We call $\cH$ a \emph{\textbf{Krein Space}}. Note $J$ is isometric on $\fD$ with respect to $\norm{\cdot}$ and so can be extended to a unitary on $\cH$. Likewise the indefinite inner product $\iip{\cdot}{\cdot}$, and the projections $P_{+}$, $P_{-}$ can be extended to to $\cH$. This is not definition of a Krein space given in \cite{Bog1974} p100, but is equivalent for our purposes via \cite{Bog1974} Theorem 2.1 p102.

As we now have two inner products on $\cH$ we have two notions of orthogonality. \emph{Krein orthogonality} will be used when we refer to orthogonality with respect to the indefinite inner product, and we will use square brackets $[$ $ ]$ when writing relations with respect to Krein orthogonality. For example $f [\perp] g$ means that $\iip{f}{g}=0$ and $\fDL [\perp]\fDJ$ means that $f[\perp]g$ for all $f \in \fDL$, $g \in \fDJ$. \emph{Hilbert orthogonality} refers to the usual notion of orthogonality with respect to $\ip{\cdot}{\cdot}$.

As well as orthogonality, the two inner products on $\cH$ give rise to two different notions of adjoints of operators. Let $T \in \op(\cH)$ be a densely defined operator on $\cH$ and let,
\begin{align*}
D(T^{\dag})= & \;\{ f \in \cH\,|\, \iip{f}{Tg}=\iip{h}{g}, \, \forall g \in \cH\},\\
D(T^{*})= & \;\{ f \in \cH\,|\, \ip{f}{Tg}=\ip{h}{g}, \, \forall g \in \cH\}
\end{align*}
Then we define $T^{\dag}f=h$ is the operator such that $\iip{f}{Tg}=\iip{h}{g}$ for $f \in D(T^{\dag})$, and a similar for $T^*$. We refer to $T^{\dag}$ as the \emph{\textbf{Krein adjoint}} or \emph{\textbf{$\dag$-adjoint}} of $T$ and  $T^{*}$ as the \emph{\textbf{Hilbert adjoint}} or \emph{\textbf{$*$-adjoint}} of $T$.

Similarly as above we will also prefix properties of operators, such a self adjointness, symmetric, isometric, unitary with a Krein or Hilbert depending on whether it is in reference to the indefinite inner product $\iip{\cdot}{\cdot}$, or Hilbert inner product $\ip{\cdot}{\cdot}$.

We have the following relation between the two adjoints (\cite{Bog1974} lemma 2.1 p122)
\begin{lemma}\label{lm:khad}
Let $T$ be a densely defined operator on a Krein space $\cH$ with fundamental symmetry $J$. Then we have that,
\[
T^*=JT^{\dag}J
\]
\end{lemma}
We also have (\cite{AzIo1989} corollary 3.8 p105)
\begin{corollary}\label{cr:khsa}
Suppose that $T$ is $\dag$-self adjoint. Then $JT$ and $TJ$ are $*$-self adjoint.
\end{corollary}
We will also use,
\begin{lemma}\label{lm:JKsp}
Let $\cH$ be a Hilbert space with inner product $\ip{\cdot}{\cdot}$ giving the norm on $\cH$. Let $J\in B(\cH)$ be a unitary operator such that $J\neq \pm \one$, $J^*=J$, and define the indefinite inner product,
\[
\iip{\cdot}{\cdot}:=\ip{\cdot}{J\cdot},
\]
on $\cH$. Then $J^2=\one$ and $\cH$ is a Krein space with indefinite inner product $\iip{\cdot}{\cdot}$, and fundamental symmetry $J$.
\end{lemma}
\begin{proof} The proof follows \cite{Bog1974} Theorem IV.5.2 p89. 
First as $J$ is unitary and $J^*=J$ we have that $J^2=JJ^*=\one$.
Next as $J$ is unitary and self adjoint with respect to the Hilbert inner product $\ip{\cdot}{\cdot}$, it has spectrum $\sigma(J)\subset (\mathbb{R} \cap \mathbb{T})=\{1,-1\}$, where $\mathbb{T}$ is the unit circle in $\mathbb{C}$. As $J\neq \pm\one$ we have $\sigma(J)=\{1,-1\}$

Now by the spectral theorem for normal operators (\cite{Con1985} Theorem IX.2.2 (a) p263) we have that,
\[
J=\int_{\sigma} \lambda dE(\lambda)=E(1)-E(-1)
\]
where $E(\Delta)$ is the spectral measure for $J$ (\cite{Con1985} Definition 1.1 p256).  Now we define $\cH_{+}=E(1)\cH$ and $\cH_{-}=E(1)\cH$ and $P_{+}=E(1)$, $P_{-}=E(-1)$. Therefore $\cH=\cH_{+}\oplus \cH_{-}$ is a fundamental decomposition of $\cH$ with fundamental symmetry $J=P_{+}-P_{-}$, and so $\cH$ is a Krein space.
\end{proof}

A useful fact about Krein space operators is the following. Let $\cH$ is a Krein space with fundamental symmetry $J$. From \cite{Min1980} proposition2, p1843.
\begin{proposition}\label{pr:Krclos}
Every Krein symmetric operator is closable (in the $\ip{\cdot}{\cdot}$ topology).
\end{proposition}

\section{Symplectic Spaces}\label{ap:symp}

The following is a collection of basic facts about symplectic spaces as given in \cite{HendrikBuch2007}. In this section $\fX$ will be a real linear space with a 
nondegenerate symplectic form $\sigma:\fX\times \fX\to\R$, and for any
subspace $S\subset \fX$ its symplectic complement
will be denoted by $S^\perp:={\set f\in \fX,\sigma(f,S)=0.}$.
By $\fX=S_1\oplus S_2\oplus\cdots\oplus S_n$ we will mean that all
$S_i$ are nondegenerate and $S_i\subset S_j^\perp$ if $i\not=j,$
and each $f\in \fX$ has a unique decomposition
$f=f_1+f_2+\cdots+f_n$ such that $f_i\in S_i$ for all $i$.
\begin{lemma}
\label{SympFacts}
\begin{itemize}
\item[(i)] If $\fX$ is countably dimensional, then it has a symplectic basis, \ie 
a basis ${\big\{q_1,\,p_1;\,q_2,\,p_2;\ldots\big\}}$ such that
$\sigma(p_i,q_j)=\delta_{ij}$ and $0=\sigma(q_i,q_j)=\sigma(p_i,p_j)$ for all $i,\,j$.
\item[(ii)] {}For any symplectic space $\fX$ we have that if $S$ is a nondegenerate
finite--dimensional subspace, then $\fX=S\oplus S^\perp$
\item[(iii)] {}For any symplectic space $\fX$ and a finite linearly independent subset
${\big\{q_1,\,q_2,\ldots,\,q_k\big\}}\subset \fX$ such that $\sigma(q_i,q_j)=0$ 
for all $i,\,j,$ there is a set
${\big\{p_1,\,p_2,\ldots,\,p_k\big\}}\subset \fX$ such that 
$B:={\big\{q_1,\,p_1;\,q_2,\,p_2;\ldots;\,q_k,\,p_k\big\}}$
is a symplectic basis for ${\rm Span}(B)$.
\end{itemize}
\end{lemma}
\begin{beweis}
(i)  Let $(e_n)_{n \in \N}$ be a linear basis of $\fX$. 
We construct the basis elements $p_n, q_n$ inductively as follows. 
If $p_1,\ldots, p_k$ and $q_1, \ldots, q_k$ are already chosen, pick a minimal
$m$ with $e_m \not\in {\rm Span}\{p_1,\ldots, p_k, q_1,\ldots, q_k\}$ and put 
$$ p_{k+1} := e_m - \sum_{i=1}^k\big(\sigma(e_m, q_i)p_i +\sigma(p_i, e_m) q_i\big) $$
to ensure that this element is $\sigma$-orthogonal to all previous ones. 
Then pick $l$ minimal, such that $\sigma(p_{k+1}, e_l) \not=0$, put
$$ \tilde q_{k+1} := e_l - \sum_{i=1}^k\big( \sigma(e_l, q_i)p_i +\sigma(p_i, e_l) q_i\big) $$
and pick $q_{k+1} \in \R \tilde q_{k+1}$ with
$\sigma(p_{k+1}, q_{k+1}) = 1$.  
This process can be repeated \textit{ad infinitum} 
and produces the required 
basis of $\fX$ because for each $k$, the span of 
$\big\{p_1,\ldots, p_k,\, q_1,\ldots, q_k\big\}$ contains at least 
$\{e_1,\ldots, e_k\}$. \chop
(ii)
Since $S$ is finite dimensional and nondegenerate, we can choose 
by (i) a symplectic basis
${\big\{q_1,\,p_1;\,q_2,\,p_2;\ldots;\,q_k,\,p_k\big\}}$ for it.
Given any $v\in \fX$ then
\[
v_S:=\sum_{i=1}^k\big(\sigma(v,q_i)\,p_i +\sigma(p_i,v)\,q_i\big)
       \in S
\]
and $v-v_S\in S^\perp$, \ie $\sigma(v-v_S,S)=0$.
Thus $\fX={\rm Span}\{S\cup S^\perp\}$, and as $\sigma$
is nondegenerate $S\cap S^\perp=\{0\}$. Moreover, if
$0=v+w$ where $v\in S$ and $w\in S^\perp$, then 
$v=-w\in S\cap S^\perp=\{0\}$, and hence any decomposition
of an $x\in \fX$ as $x=x_1+x_2$ where $x_1\in S,$ $x_2\in S^\perp$
is unique. Thus $\fX=S\oplus S^\perp$. \chop
(iii) 
We first find via the method of part (i),
symplectic pairs ${\big\{\wt{q}_1,r_1;\ldots;\wt{q}_k,\,r_k\big\}}\subset \fX$
such that the nondegenerate subspaces
$S_j:= {{\rm Span}\big\{\wt{q}_1,r_1;\ldots;\wt{q}_j,\,r_j\big\}}\supset\{q_1,\ldots, q_j\}$
but $q_{j+1}\not\in S_j$.
We construct the basis elements $\wt{q}_i,\, r_i$ inductively as follows.
If $r_1,\ldots, r_j$ and $\wt{q}_1, \ldots, \wt{q}_j$ are already chosen,
put
$$ \wt{q}_{j+1} := q_{j+1} - \sum_{i=1}^k\big(\sigma(q_{j+1}, \wt{q}_i)r_i +
\sigma(r_i, q_{j+1}) \wt{q}_i\big) $$
to ensure that $\wt{q}_{j+1}\in S_j^\perp$. By (ii), $\fX=S_j\oplus S_j^\perp$
hence $S_j^\perp$ is nondegenerate, so there is an element $r_{j+1}\in S_j^\perp$
such that $\sigma(r_{j+1},\wt{q}_{j+1})=1$.
It follows that $q_{j+2}\not\in S_{j+1}$ and that
$\{q_1,\ldots, q_{j+1}\}\subset S_{j+1}$.
This process can be repeated to produce the required symplectic bases.
Next, we want to show that in $S_k$ we can choose 
${\big\{p_1,\,p_2,\ldots,\,p_k\big\}}$ such that 
${\big\{q_1,\,p_1;\,q_2,\,p_2;\ldots;\,q_k,\,p_k\big\}}$
is a symplectic basis for $S_k$.
Now $\{q_1,\ldots, q_k\}\subset\{q_1,\ldots, q_k\}^\perp$
where henceforth the symplectic complements are all taken in $S_k$.
We claim that the containment ${\{q_2,\ldots, q_k\}^\perp}\supset
{\{q_1,\,q_2,\ldots, q_k\}^\perp}$ is proper. The map
$\varphi:S_k\to S_k^*$ by $\varphi_x(y):=\sigma(x,y)$ is a linear isomorphism
by nondegeneracy of $\sigma$. Then for any set $R\subset S_k$ we have 
$\varphi\big(R^\perp\big)=R^0$
\ie the annihilator of $R$ in $S_k^*$, hence 
${\rm dim}(R^\perp)={\rm dim}(R^0)=2k-{\rm dim}\big({\rm Span}(R)\big)$.
Thus ${\rm dim}{\{q_1,\ldots, q_j\}^\perp}=2k-j$ from which the claim follows.
Thus there is an $r\in{\{q_2,\ldots, q_k\}^\perp}\backslash
{\{q_1,\,q_2,\ldots, q_k\}^\perp}$ such that $\sigma(r,q_1)\not=0$.
In particular, let $p_1$ be that multiple of $r$ such that
$\sigma(p_1,q_1)=1$. Let $T_1:={\rm Span}\{q_1,p_1\}$
then ${\{q_2,\ldots, q_k\}}\subset T_1^\perp,$ and by (ii)
we have $S_k=T_1\oplus T_1^\perp$ where $T_1^\perp$ is nondegenerate.
 Thus we can now
repeat this procedure in $T_1^\perp$ starting from $q_2$ to obtain
$p_2$. This procedure will exhaust $S_k$ to produce the 
desired symplectic 
basis ${\big\{q_1,\,p_1;\,q_2,\,p_2;\ldots;\,q_k,\,p_k\big\}}.$ 
\end{beweis}

\section{Quantum Dirac Constraints}\label{app:Tp}
The following account of Quantum Dirac Constraints in the algebraic context follows the  survey \cite{Hendrik2006} and references therin. 

A brief heuristic outline of the Quantum Dirac Constraints is as follows. Suppose that $\cH$ is a Hilbert space and that $\cC=\{G_j\,|\, j \in \Lambda \}$ is a set of operators that select a physical subspace, ie
\[
\cH_p:= \cap_{j \in \Lambda} \ker G_j
\]
Then the Dirac observables are traditionally taken as the commutant of the constraints $\cC'$, but can be enlarged to be the algebra generated by the self-adjoint operators that preserve $\cH_p$. The final constrained system is these observables restricted to $\cH_p$. This procedure can be problematic for reasons such as spectral issues (think of trying to impose momentum $p=0$ when $p$ has a canonical conjugate), hence we abstract the process to a $C^*$-algebraic setting that is independent of the defining representation, then look for representations where the problematic issues are no longer present.

We give a summary of the relevant aspects of this abstraction found in in the survey \cite{Hendrik2006}. Much more can be said about $C^*$-Dirac constraints than will be given here, see \cite{HendrikHu1985,Hendrik2006,HendrikHu1987,Hendrik1988, Hendrik2000} for more.

\begin{definition}
A \emph{quantum system with constraints} is a pair $(\cA, \cC)$ where the \emph{field algebra} $\cA$ is a unital $C^*$-algebra containing the \emph{constraint set} $\cC=\cC^*$. A \emph{constraint condtion} on $(\cA, \cC)$ consists of the selection of the physical state space by:
\[
\fS_D:=\{\omega \in \fS(\cA) \,|\, \pi_{\omega}(C)\Omega_{\omega}=0 \quad \forall C \in \cC\}
\]
where $\fS(\cA)$ denotes the state space of $\cA$, and $(\pi_{\omega}, \cH_{\omega}, \Omega_{\omega})$ denotes the \emph{GNS}-data of $\omega$. The elements of $\fS_D$ are called \emph{Dirac states}. The case of \emph{unitary constriants} means that $\cC=\cU-\one$ for a set of unitaries $\cU \subset \cA_u$, and for this we will also use the notation $(\cA,\cU)$. 
\end{definition}

Now observe that we have,
\begin{align*}
\fS_D= & \;\{\omega \in \fS(\cA)\,|\, \omega(C^*C)=0 \quad \forall C \in \cC \},\\
= & \;\{\omega \in \fS(\cA)\,|\, \cC \subset \cN_{\omega}\},\\
= & \; \cN^{\perp}\cap \fS(\cA).
\end{align*}
Here $\cN_{\omega}:=\{ A \in \cA \,|\, \omega(A^*A) \}$ is the left kernel of $\omega$ and $\cN:= \cap \{\cN_{\omega} \,|\, \omega \in \fS_{D} \}$, and the $\perp$ denotes the annihilator in the dual of $\cA$. 

We now have the equality $\cN=[\cA \cC]$. Since $\cC$ is self-adjoint and contained in $\cN$ we have $\cC \subset \osalg{\cC} \subset (\cN \cap \cN^*)= [\cA \cC]\cap [ \cC \cA]$. We can use these facts to get,
\begin{theorem}\label{th:DsD}
We have
\begin{itemize}
\item[(i)] $\fS_D=\{0\}$\, \text{iff}\, $\one \notin \osalg{\cC}$ \, \text{iff}\, $\one \notin \cD:= \cN\cap \cN^*$.
\item[(ii)] $\omega \in \fS_D$ \, \text{iff}\, $\pi_{\omega}(\cD)\Omega_{\omega}=0$.
\item[(iii)] An extreme Dirac state is pure.
\end{itemize}
\end{theorem}

A constraint set is \emph{first class} if $\one \notin \osalg{\cC}$ which by the above theorem is the assumption that the constraints are non-trivial.

Now we define the observable algebra as,
\[
\cO:=\{ A\in \cA \,|\, [A,D]\in \cD \quad \forall D \in \cD \}.
\]
We get
\begin{theorem}
We have
\begin{itemize}
\item[(i)] $\cD=\cN\cap \cN^*$ is the unique maximal \emph{$C^*$}-algebra in $\cap \{ \ker \omega \,|\, \omega \in \fS_D \}$. Furthermore $\cD$ is a hereditary \emph{$C^*$}-algebra of $\cA$.
\item[(ii)] $\cO=\cM_{\cA}(\cD):\{ A \in \cA \,|\, AD \in \cD \ni DA \quad \forall D \in \cD \}$, \ie it is the relative multiplier algebra of $\cD$ in $\cA$.
\item[(iii)] $\cO=\{A\in \cA\,|\, [A, \cC] \subset \cD \}$.
\item[(iv)] $\cD=[\cO\cC]=[\cC\cO]$.
\item[(v)] For the case of unitary constraints, \ie $\cC=\cU-\one$, we have $\cU \subset \cO$ and $\cC=\{A\in \cA \,|\, \alpha_U(A)-A\in \cD \quad \forall U \in \cU\}$ where $\alpha_{U}:=AdU$.
\end{itemize}
\end{theorem}
Therefore $\cD$ is a closed two-sided ideal in $\cO$ and the traditional observables $\cC'\subset \cO$ where $\cC'$ is the commutant of $\cC$ in $\cA$.

Define the \emph{maximal $C^*$-algebra of physical observables} as
\[
\cP:=\cO/\cD
\]
The factoring procedure is the step of imposing constraints. We call this method of imposing constraints the $\emph{T-procedure}$. We require that all physical information is contained in $(\cP, \fS(\cP))$. It is possible that $\cP$ is not simple and in the case we adjust as in \cite{Hendrik2006} p101. We have the following connection
\begin{theorem}
There exists a $w^*$-bijection between the Dirac states on $\cO$ and the states on $\cP$.
\end{theorem}

Although the $T$-procedure is an abstract $C^*$-algebra procedure not dependent on the orginal representation, it can be helpful to work in a representation to aid intuition and calculations. In fact as $\cD$ is a hereditary subalgebra of $\cA$, we can utilise results in \cite{Ped79} chapters 1-3 in relation to the universal enveloping von Neumann algebra. Denote the universal representation by $\pi_u$ on the universal Hilbert space $\cH_u$ and let $\cA''$ be the strong closure of $\pi_u(\cA)$ and make identification of $\cA$ with a subalgebra of $\cA''$, \ie we generally omit $\pi_u$ explicitly. Also if $\omega\in \fS(\cA)$ the we denote by $\omega$ the unique normal extension of $\omega$ from $\cA$ to $\cA''$.

From \cite{Ped79} we have that
\begin{definition}
A projection $P\in \cF''$ is called open if $\cL=\cA \cap (\cA''P)$ is a closed left ideal of $\cA$.
\end{definition}
From \cite{Ped79} Theorem 3.10.7, Proposition 3.11.9 and Remark 3.11.10 we have bijections of open projections with:
\begin{itemize}
\item[(i)] hereditary ideals of $\cA$ given by $P\to P\cA''P\cap \cA$.
\item [(ii)]closed left ideals of $\cA$ by $P \to \cA''P \cap \cA$.
\item[(iii)] weak $*$-closed faces containing $0$ of the quasi-state space $\mathcal(Q)(\cA)$ by
\[
P \to \{ \omega \in \mathcal{Q}(\cA) \,|\, \omega(P)=0 \}
\]
\end{itemize}
Using this we get the following results (\cite{Hendrik2006} Theorem 4, Theorem 5)

\begin{theorem}\label{th:vnT}
For the constraint system $(\cA, \cC)$ there exists an open projection, $P\in \cA''$, such that,
\begin{itemize}
\item [(i)] $\cN=\cA '' P \cap \cA$,
\item [(ii)] $\cD=P \cA'' P \cap \cA$,
\item [(iii)] $\fS_{D}= \{ \omega \in \fS(\cA)\, |\, \omega(P)=0 \}$
\item [(iv)] $\cO = \{ A \in \cA \, |\, PA(\one - P)=(\one - P)AP=0 \}= P'\cap \cA$,
\item [(v)] $\cP \cong (\one -P)(P' \cap \cA)= (\one -P)\cO (\one -P)$
\end{itemize}
\end{theorem}
Now in the universal representation we define the physical space to be $\cH^p_{u}:=\cap_{C\in \cC}\ker C$. Then from the above theorem we that $\cH_u=\ker (\one-P)$ where $P$ is projection in the statement of Theorem \eqref{th:vnT}. That is, $P$ is a condition that selects the same physical subspace as the constraints. Note however that $P$ is not in $\cA$ in general and so we would have to extend $\cA$ if we wanted to apply the $T$-procedure using $P$ as a constraint.

Still we have the decomposition of $\cH_u=P\cH_u\oplus (\one -P)\cH_u$, and with resepect to this decomposition using Theorem \eqref{th:vnT} (ii),(iii),(v) we may write
\begin{align*}
\cD= & \;\left\lbrace \left. A\in \cA \, \right| \, A=\left(
\begin{matrix}
D & 0 \\
0 & 0 \\
\end{matrix}
\right),
\quad D \in P\cA P 
\right\rbrace,\\
\cO= & \;\left\lbrace \left. A\in \cA \, \right| \, A=\left(
\begin{matrix}
A & 0 \\
0 & B \\
\end{matrix}
\right),
\quad A \in P\cA P,\,B\in(\one-P)\cA(\one-P)
\right\rbrace,\\
\intertext{and,}
\cP\cong&\left\lbrace \left. A\in \cA \, \right| \, A=\left(
\begin{matrix}
0 & 0 \\
0 & A \\
\end{matrix}
\right),
\quad A\in(\one-P)\cA(\one-P)
\right\rbrace,\\
\end{align*}
Now we define
\[
\Phi_p(A)=(\one-P)A(\one-P), \qquad A \in \cO
\]
From Theorem \eqref{th:vnT} (iii) we see that $\Phi_p$ is a homomorphism on $\cO$ and from Theorem \eqref{th:vnT} (v) we see that $\ker \Phi_p=\cD$. Therefore we have that, 
\begin{equation}\label{eq:Phkphys}
\cP\cong \cO/ \ker \Phi_p \cong \Phi_p(\cO)
\end{equation}

\begin{eje}\label{ex:DC}
\begin{itemize}
\item[(i)]Suppose that $\cH$ is a Hilbert space, $Q$ is a self-adjoint projection on $\cH$ and $(\cA,\cC)=(B(\cH),\{Q\})$. It is straightforward to see that $P=\pi_u(Q)$ where $P$ is the projection in Theorem \eqref{th:vnT}. Therefore $\cO = Q'$ and $\cP \cong (\one -Q)Q' (\one -Q)$. Hence in this case we have that the $T$-procedure gives that the observables are the same as the Dirac observables. 
\item[(ii)] Suppose that $\cH$ is a seperable Hilbert space, $K(\cH)$ are the compact operators on $\cH$ and that $(\cA, \cC)=(B(\cH),K(\cH))$. Then $\cC'=\mathbb{C}\one$ (see) but $\cO=B(\cH)$ as intersection of the kernel of all the finite rank operators is $\{0\}$.
\end{itemize}
\end{eje}
\section{Covariance for $C^*$-BRST II} \label{ap:cvBRST II}

The Fock-Krein representation of the auxiliary algebra in Section \ref{sec:FKCCR}, $\pi_F$, is faithful and the Poincar{\'e} transformations are defined there by $\alpha_g:= Ad(\Gamma(V_g^{\dag}))$. Therefore it is natural to see if we encode the Poincar{\'e} transformation on the auxiliary algebra via this representation, and so for the remainder of this section we will assume that we are working in the Fock representation and will identify $\cR(\fD, \sigma_2)$ with $\pi_F(\cR(\fD, \sigma_2))$. We will also use the notation $R_f=R(1,f)$.

An important point to note is that with respect to the IIP $\iip{\cdot}{\cdot}$ on $\fH$ (see subsection\eqref{sbs:testfunc}) $V_g$ is Krein-unitary, but not $\ip{\cdot}{\cdot}=\iip{\cdot}{J\cdot}$-unitary, \ie Hilbert unitary. Two important consequences of this are:
\begin{itemize}
\item[(i)] We have that $\alpha_g$ will be a $\dag$-automorphism on the unbounded algebra generated by the $A(f)$'s, but that in general $\alpha_g$ will not be a $*$-automorphism.
\item[(i)] $V_g$ is Krein unitary hence invertible and hence $\Gamma(V_g)$ will define Krein unitary operator on the finite particle space $\fF_0(\cH)$. However, for $g=(\Lambda,a)$ where $\norm{\Lambda}>1$ it follows from,
\[
\Gamma(V_g)(\psi_1 \otimes \ldots \otimes \psi_n)=V_g\psi_1 \otimes \ldots \otimes V_g\psi_n,
\]
that $V_g$ is unbounded and so does not extend to $\fF(H)$. Moreover, there will be problems defining $\alpha_g$ on all of $\cR(\fD, \sigma_2)$.
\end{itemize}
To understand the problem of $\alpha_g$ on $\cR(\fD, \sigma_2)$ more explicitly, we look at its action on $\phi(f)$.
\begin{align}
\alpha_g(\phi(f))&=\alpha_g(A(P_{+}f)+iA(iP_{-}f)),\notag \\
&=A(V_g^{\dag}P_{+}f)+iA(iV_g^{\dag}P_{-}f),\notag \\
&=\phi((P_{+}V_g^{\dag}P_{+}+P_{-}V_g^{\dag}P_{-})f)+i\phi(i(P_{+}V_g^{\dag}P_{-}+P_{-}V_g^{\dag}P_{+})f),\notag \\
&=(1/2)[\phi((V_g^{\dag}+V_g^{*})f)+i\phi(i(V_g^{\dag}-V_g^{*})f)],\label{eq:covgoop}\\
&=(1/\sqrt{2})[a(V_g^{\dag}f)+a^*(V_g^{*}f)]\notag
\end{align}
where the second last line follows from lemma \eqref{lm:plmnJ} and $V_g^*=JV_g^{\dag}J$. This shows explicitly that $\alpha_g$ will be a $*$-automorphism for all $f \in \fD$ iff $(V_g^{\dag}-V_g^{*})=0$, which is when $g$ is a rotation and/or translation.

Using this we can show that in fact that if $g$ is a Lorentz boost, then $\alpha_g$ can map resolvents to unbounded operators. First note that if we define coordinates $(x_0,x_1,x_2,x_3) \in \mathbb{R}^4$ and  $B_1$ is Lorentz boost in the $x_1$-direction, then 
\begin{equation}\label{eq:matB1}
B_1=B_1(t)=\left(
\begin{matrix}
\cosh(t) & \sinh(t) & 0 & 0\\
\sinh(t) & \cosh(t) & 0 & 0 \\
0 & 0 & 1 & 0 \\
0 & 0 & 0 & 1
\end{matrix}
\right),  \qquad t \in \mathbb{R}
\end{equation}
which is given in \cite{Na1964} p90-92, (Note that \cite{Na1964} uses the index $x_4$ for the time co-ordinate and the order $x=(x_1,x_2,x_3,x_4)$, while we use $x_0$ and  $x=(x_0,x_1,x_2,x_3)$, which leads to a different arrangement for the entries of $B_1$). 
 
Let $P_{-}=diag(-1,0,0,0)$ and $P_{+}=diag(0,1,1,1)$ be projections on $\mathbb{R}^4$. For $e_1=(0,1,0,0)\in \mathbb{R}^4$ then we see that $ P_{+}B_1P_{+}e_1 \neq 0$ and $ P_{1}B_1P_{+}e_1 \neq 0$. Hence by the definition of $V_{B_1}$ (equation \eqref{eq:oppoinact}), we see that we can choose $h \in P_{+}\fD$, such that $h_1=P_{+}V_{B_1}h\neq 0, \, h_2=iP_{-}V_{B_1}h \neq 0$.
Now $\phi(h)=A(P_{+}h)+iA(iP_{-}h)=A(h)$ so,
\[
\alpha_{B_1^{\dag}}(\phi(h))=A(V_{B_1}h)=\phi(P_{+}V_{B_1}h)+i\phi(iP_{-}V_{B_1}h)=\phi(h_1)+i\phi(ih_2).
\]
As $P_{+}\fD, P_{-}\fD$ $\sigma_2$-symplectically commute, we have $[\phi(h_1),\phi(h_2)]=[R_{h_1},R_{h_2}]=0$, and so we have joint spectral theory for $\phi(h_1),\phi(h_2)$. We can use this to calculate $\alpha_g(R_h)$. Take $\psi \in \fF_0(\cH)$, then,
\begin{align*}
\psi=(i-\phi(h))R_h\psi \Rightarrow \psi&=\alpha_g(i\one-\phi(h))\alpha_g(R_h)\psi,\\
 &=(i-\phi(h_1)-i\phi(h_2))\alpha_g(R_h) \psi,\\
 &=[\int_{\mathbb{R}^2}(i\one-\lambda-i\mu)dP(\lambda,\mu)]\alpha_g(R_h) \psi,\\
&=[\int_{\mathbb{R}^2}(i(1-\mu)-\lambda)dP(\lambda,\mu)]\alpha_g(R_h) \psi,\\
\Rightarrow \alpha_g(R_h)\psi&'='(\int_{\mathbb{R}^2}(i(1-\mu)-\lambda)^{-1}dP(\lambda,\mu))\psi,
\end{align*}
Now as $(i(1-\mu)-\lambda)^{-1}$ has a singularity at $\mu=1, \lambda=0$ then we have that $\alpha_g(R_h)$ is unbounded, and so we cannot define $\alpha_g$ abstractly on all of $\cR(\fD, \sigma_2)$. 

Although the Poincar{\'e} transformations as defined directly as in the Fock-Krein representation do not define on all of the auxiliary algebra, we can look for strategies to encode them in a way so that we can get correct transformations back on the physically interesting objects.

One strategy to use is to see if we can find a homomorphism from $\gamma:\Sp(\fD, \sigma_1) \to \Sp(\fD, \sigma_1)$, define the Poincar{\'e} transformations on the auxiliary algebra via Theorem \eqref{th:symautRA} (v) and check that these give the same transformations as in Subsection \ref{sbs:QEMcovI} when factored to the physical subspace.

A step in this direction begins with 
\begin{proposition} Define the \emph{real} linear operator
\[
T:\fX \to \fX \qquad\text{by} \qquad Tf:=P_{+}f + P_{-}Cf,
\]
where $C$ is as Proposition \eqref{pr:QM2tfdec}. Then
\[
\gamma_T:\sp(\fX, \sigma_2) \to\sp(\fX, \sigma_2)\qquad \text{defined by} \qquad M \to TMT, 
\]
is a \emph{real} algebra isomorphism.
\end{proposition}
\begin{proof}
By Proposition \eqref{pr:QM2tfdec} (ii) we have
\[
\smp{Cf}{Cg}{2}=-\smp{f}{g}{2}, \qquad f,g \in \fD.
\]
Hence for all $f,g \in \fD$ we have:
\begin{align*}
\smp{Tf}{Tg}{2}= & \;\smp{P_{+}f}{P_{+}g}{2}+\smp{CP_{-}f}{CP_{-}g}{2},\\
= & \;\smp{P_{+}^2f}{g}{2}-\smp{P_{-}^2f}{g}{2},\\
= & \;\smp{Jf}{g}{2},\\
= & \;\smp{f}{g}{1},
\end{align*}
Moreover, for $M\in \sp(\fX, \sigma_2)$ we have for all $f,g \in \fX$:
\[
\smp{TMTf}{TMTg}{1}=\smp{MTf}{MTg}{2}=\smp{Tf}{Tg}{2}=\smp{f}{g}{1}, 
\]
and as $C^2=\one$, it is obvious that $T^2=\one$ and so $\gamma_T$ is a \emph{real} algebra isomorphism from $\Sp(\fX, \sigma_1)\to \Sp(\fX, \sigma_2)$. 
\end{proof}
However $C$ and $K$ do not commute, hence $T$ does not extend to a linear or anti-linear operator on $\fD$. Hence we do \emph{not} have that $\gamma_T$ extends to an isomorphism or anti-isomorphism of $\Sp(\fD, \sigma_1)\to \Sp(\fD, \sigma_2)$. So although $\gamma_{T}$ is the first step in implementing the above strategy of mapping $\sp(\fD,\sigma_1)$ to $\sp(\fD,\sigma_2)$, it is not clear at present how to continue.

Another strategy towards implementing $\alpha_g$ is to reduce the problem to implementing Lorentz in the $x_1$-direction. There is no problem implementing rotations and translations as these are both $\Sp(\fD, \sigma_1)$ and $\Sp(\fD,\sigma_2)$ symplectic, and a general Lorentz transformations can always be written as $\Lambda=R_2B_1R_1$ where $R_2,R_1$ are rotations and $B_1$ is a boost in the $x_1$-direction (see \cite{Na1964} (III),p93).

In this direction we begin by breaking the transformations $V_g$ into more manageable pieces. Note that for $g=(\Lambda,a) \in \cP^{\uparrow}_{+}$, we have that $V_g= T_g U_{\Lambda}=U_{\Lambda}T_g$ where,
\[
(T_gf)(p):=e^{ipa} f(\Lambda^{-1}p), \qquad (U_{\Lambda}f)(p):=\Lambda f(p)\quad \forall f \in \cS(\mathbb{R}^4,\mathbb{C}^4). 
\]
now it is easy to check that $T_g$ is both $\sigma_1$-symplectic, and $\sigma_2$-symplectic, hence will generate automorphisms on  and $\cR(\fD, \sigma_2)$ (and $\cR(\fD, \sigma_1)$), and $\Gamma(T_g)$ implements these autormorphisms in the Fock representations. Therefore if we let $\cD=\fF_0(\cH)$, then,
\[
Ad(\Gamma(T_g))(\phi_{\pi}(f)) \cD = \phi_{\pi}(T_gf)\cD.
\]
As $[T_g, P_{+}]=[T_g, P_{-}]$, we use,
\[
A_{\pi}(f)\cD=(\phi_{\pi}(P_{+}f)+i\phi_{\pi}(P_{-}if))\cD,
\]  
to get,
\[
Ad(\Gamma(T_g))(A(f))\cD=A(T_gf)\cD.
\]
So we only need to consider how to encode of $U_{\Lambda}$ for $\Lambda \in \cL^{\uparrow}_{+}$, in $\cR(\fD, \sigma_2)$.

Now $\Lambda=R_2B_1R_1$ where $R_2,R_1$ are rotations and $B_1$ is a boost in the $x_1$-direction (see \cite{Na1964} (III),p93). It is straightforward to check that  $U_{R_1},U_{R_2}$ are both also $\sigma_2$ symplectic and $[U_{R_i}, P_{+}]=[U_{R_i}, P_{-}]$, $i=1,2$, hence as above,
\[
Ad(\Gamma(U_{R_i}))A(f)\cD=A(U_{R_i}f)\cD, \qquad  i=1,2.
\]
So all we have left to do is to construct $\alpha_{B_1}$.  $B_1$ is given by the matrix in equation \eqref{eq:matB1}, and by equation \eqref{eq:covgoop},
\begin{align*}
\alpha_{B_1(t)}(\phi(f))=(1/2)[\phi((B_1(t)^{\dag}+B_1(t)^{*})f)+i\phi(i(B_1(t)^{\dag}-B_1(t)^{*})f)].
\end{align*}
Furthermore, by \eqref{eq:matB1} $B_1^*=B_1$, so using $B_1^{\dag}=JB_1^*J$ and $J=diag(-1,1,1,1)$ we get,
\begin{equation*}
B_1(t)^{\dag}=\left(
\begin{matrix}
\cosh(t) & -\sinh(t) & 0 & 0\\
-\sinh(t) & \cosh(t) & 0 & 0 \\
0 & 0 & 1 & 0 \\
0 & 0 & 0 & 1
\end{matrix}
\right),  \qquad t \in \mathbb{R}.
\end{equation*}
Therefore,
\begin{equation*}
\left. \frac{\mathrm{d}}{\mathrm{d}t}\right|_{t=0}(B_1(t)^{\dag}+B_1(t)^{*})=2\left. \frac{\mathrm{d}}{\mathrm{d}t}\right|_{t=0}\left(
\begin{matrix}
\cosh(t) & 0 & 0 & 0\\
0 & \cosh(t) & 0 & 0 \\
0 & 0 & 1 & 0 \\
0 & 0 & 0 & 1
\end{matrix}
\right)
= 0,
\end{equation*}
and,
\begin{equation*}
\left. \frac{\mathrm{d}}{\mathrm{d}t}\right|_{t=0}(B_1(t)^{\dag}-B_1(t)^{*})=2\left. \frac{\mathrm{d}}{\mathrm{d}t}\right|_{t=0}\left(
\begin{matrix}
0 & \sinh(t) & 0 & 0\\
\sinh(t) & 0 & 0 & 0 \\
0 & 0 & 0 & 0 \\
0 & 0 & 0 & 0
\end{matrix}
\right)
= 
2\left(
\begin{matrix}
0 & 1 & 0 & 0\\
1 & 0 & 0 & 0 \\
0 & 0 & 0 & 0 \\
0 & 0 & 0 & 0
\end{matrix}
\right).
\end{equation*}
So we get that, 
\begin{align*}
\left. \frac{\mathrm{d}}{\mathrm{d}t}\right|_{t=0}\alpha_{B_1(t)}(\phi(f))&=(1/2)[\phi((B_1(t)^{\dag}+B_1(t)^{*})f)+i\phi(i(B_1(t)^{\dag}-B_1(t)^{*})f)],\\
&=i\phi(ib_1f)
\end{align*}
where,
\[
b_1:=\left(
\begin{matrix}
0 & 1 & 0 & 0\\
1 & 0 & 0 & 0 \\
0 & 0 & 0 & 0 \\
0 & 0 & 0 & 0
\end{matrix}
\right).
\]
We summarize in
\begin{lemma}
Let $B_1(t)$ be a boost in the $x_1$ direction. Then,
\begin{align*}
\left. \rd(\phi(f))\psi:=(\frac{\mathrm{d}}{\mathrm{d}t}\right|_{t=0}\alpha_{B_1(t)})(\phi(f))\psi=i\phi(ib_1f)\psi,
\end{align*}
where,
\[
b_1:=\left(
\begin{matrix}
0 & 1 & 0 & 0\\
1 & 0 & 0 & 0 \\
0 & 0 & 0 & 0 \\
0 & 0 & 0 & 0
\end{matrix}
\right),
\]
and $\psi \in \cD$.
\end{lemma}
From this lemma there are two directions to go:

\medskip
\noindent{\bf{Direction1: Mollify $\rd$}}

This approach is taken in a different context in \cite{HendrikBuch2006} p12. We can calculate $\rd$ on $R(\lambda,f)$ via,
\begin{align*}
R(\lambda,f)\psi= \alpha_{B_1(t)}(R(\lambda,f))\alpha_{B_1(t)}(i\lambda \one -\phi(f))R(\lambda,f)\psi,
\end{align*}
and so differentiating,
\begin{align*}
\rd(R(\lambda,f))(i\lambda \one -\phi(f))R(\lambda,f)\psi= & \;R(\lambda,f)\rd(i\lambda \one -\phi(f))R(\lambda,f)\psi,
\end{align*}
{which implies,}
\begin{align*}
\rd(R(\lambda,f))\psi= & \;iR(\lambda,f)\phi(ib_1f)R(\lambda,f)\psi,\\
= & \;i(\phi(ib_1f)+\smp{ib_1f}{f}{2}\one)R(\lambda,f)^2\psi
\end{align*}
where $R(\lambda,f)\psi \in \fF_0(\cH)$. Using the above we see that for any $A \in \cR(\fD, \sigma_2)_0$ we can find a monomial of resolvents $M_A$ to mollify $\rd(A)$. We have that $M_A\rd(A)\in \cR(\fD, \sigma_2)$ and so this encodes the Lorentz boosts on the auxiliary algebra in infintesimal form. With this done we can the aim is to recover $\rd$ on the auxiliary fields and then covariant fields in other representations by using tools such as Theorem \ref{RegThm}. A problem with this is that the mollifying monomial $M_A$ depends on the original $A \in  \cR(\fD, \sigma_2)_0$ and so the `mollified $\rd$' is a difficult object to analyse.

\medskip
\noindent{\bf{Direction2: Use $T \in \Sp(\fD,\sigma_2)$}}

The idea here is that $ ib_1$ is a self-adjoint operator with respect to $\ip{\cdot}{\cdot}$ on $\cH=\overline{\fD}$. Therefore $W_t:=\exp(itb_1)$ is unitary on $\cH=$ and so $\sigma_2$-symplectic on $\fD$ ($b_1$ preserves $\fD$ so so does $\exp(itb_1)$). Explicitly,
\begin{equation*}
W(t)=\left(
\begin{matrix}
\cos(t) & i\sin(t) & 0 & 0\\
i\sin(t) & \cos(t) & 0 & 0 \\
0 & 0 & 1 & 0 \\
0 & 0 & 0 & 1
\end{matrix}
\right),  \qquad t \in \mathbb{R}.
\end{equation*}
Therefore we can define the corresponding $\alpha_{W_t}\in \Aut(\cR(\fD,\sigma_2))$. In `nice' representations (such as the Fock) we have that $\alpha_{W_t}$ is unitarily implemented by $S_t$ and we can differentiate to get a self-adjoint generator $M$ with dense domain $D(M)$ for $S_t$. That is, for a `nice' representation, $\pi$, 
\[
\pi(\alpha_{W_t}(A))=S_t\pi(A)S_t, \qquad S_t=\exp(itM),
\]
and so when the appropriate fields exist, we get
\[
S_t \phi_{\pi}(f) S_t\psi = \phi_{\pi}(W_tf)\psi, \qquad \psi \in D(\phi_{\pi}(f))\cap D(\phi_{\pi}(W_tf)).
\]
Then we differentiate to get a $*$-derivation,
\[
\rd(\phi_{\pi}(f))\psi:=\left. \frac{\mathrm{d}}{\mathrm{d}t}\right|_{t=0}(\pi(\alpha_{W_t}(\phi_{\pi}(f))))\psi=\phi_{\pi}(ib_1f)\psi=[iM,\phi_{\pi}(f)]\psi.
\]
That is, under the appropriate regularity conditions, we we will have a well defined derivation that is generated by a self-adjoint operator, ie
\[
\rd(\phi_{\pi}(f))\psi= \phi_{\pi}(ib_1f)\psi=[iM,\phi_{\pi}(f)]\psi.
\]
Now if we let $L=iM$ then we have that,
\[
\rd_{tB_1}(A)\psi:=[A,tL]\psi,
\]
is a well defined derivation for $A$ that preserve the appropriate domain, and we have that,
\[
\rd_{tB_1}(\phi_{\pi}(f))\psi:=i\phi_{\pi}(itb_1f)\psi,
\]
the exact relation we want for the infintesimal version of $B_1$. Now we can check by direct computation that $[b_1,P_{+}]=[P_{-},b_1]$ and so we get that,
\begin{align*}
\rd_{tB_1}(A_{\pi}(f))\psi= & \;\rd_{tB_1}(\phi_{\pi}(P_{+}f)+i\phi_{\pi}(iP_{+}f))\psi,\\
= & \;i(\phi_{\pi}(itb_1P_{+}f)+i\phi_{\pi}(-tb_1P_{+}f))\psi,\\
= & \;i\phi_{\pi}(itP_{-}b_1f)+\phi_{\pi}(tP_{+}b_1f))\psi,\\
= & \;A(tb_1f)\psi
\end{align*}
Therefore $(\rd_{tB_1})^n(A_{\pi}(f))\psi=A_{\pi}((tb_1)^nf)\psi$ and so given strong convergence of the fields in their arguments, we get,
\[
\exp(\rd_{tB_1})(A_{\pi}(f))\psi=A_{\pi}(\exp(tb_1)f)\psi=A(B_1(t)f)\psi,
\]
and so we can reconstruct $\alpha_{B_1(t)}$. Also by the derivation property of $\rd_{tB_1}$ we can check that $\alpha_{B_1(t)}(AB)\psi=\alpha_{B_1(t)}(A)\alpha_{B_1(t)}(B)\psi$ for appropriate $A,B \in \op(\cH_{\pi})$ (such as the fields). Note that as $L$ is skew-self adjoint we will not have that $\alpha_{B_1(t)}$ is a $*$-automorphism or that it is even bounded, as in the case for the Fock representation.

To summarize for direction 2, we start with $W_t \in \Sp(\fD,\sigma_2)$ and look for restrictions on `nice' representations such that we can strongly differentiate and exponentiate and get convergence in the arguements for the fields, etc. Given these conditions we construct the derivation corresponding to the infintesimal version of $\alpha_{B_1}$ exponentiate this on the fields to get the $\alpha_{B_1}$. Once we have these it is straightforward to get any covariant transform in the nice representations. 

So we see that we have several strategies and directions for encoding the Poincar{\'e} transformations directly on the auxiliary algebra. The first strategy was mapping $\sigma_1$-symplectic transformations to $\sigma_2$-symplectic transformations via some isomorphism $\gamma$ and checking that the automorphisms generated are then match those in Subsection \ref{sbs:QEMcovI}  when factored to the physical algebra. At present such a $\gamma$ has only been constructed on $\Sp(\fD_C,\sigma_1)$ a real subalgebra of $\Sp(\fD,\sigma_1)$. The second strategy was to reduce the problem to encoding $\alpha_{U_{B_1}}$ where $B_1$ was a boost in the $x_1$-direction. This is still a difficult problem however we can instead encode the infintesimal version of $\alpha_{U_{B_1}}$ using $\rd$. There were to approaches to doing this: one was mollifying $\rd$ but this faced the problem that the `mollified $\rd$' depended on the $A \in \cR(\fD,\sigma_2)_0$ on which it was acting; The other was to use $W(t)\in \Sp(\fD, \sigma_2)$ and differentiate to get $\rd$, which faces the problem that we need to have nice representations where we can differentiate. Also, in both of the infintesimal strategies, we would need to restrict to representations where we could re-exponentiate $\rd$.

\bibliographystyle{plain}
\bibliography{Mybibfinal}

\begin{thebibliography}{100}

\bibitem{AbNak1988}
M.~Abe and N.~Nakanishi.
\newblock B{RS}-invariant {L}agrangian density in the new local supersymmetry
  of the vierbein formalism of {E}instein gravity.
\newblock {\em Progr. Theoret. Phys.}, 79(1):240--249, 1988.

\bibitem{AbNak1996}
M.~Abe and N.~Nakanishi.
\newblock Subtlety in the anomaly calculation of string theory in the harmonic
  gauge.
\newblock {\em Progr. Theoret. Phys.}, 96(6):1281--1290, 1996.

\bibitem{AbNak2002}
M.~Abe and N.~Nakanishi.
\newblock Exact solutions to the two-dimensional {$BF$} and {Y}ang-{M}ills
  theories in the light-cone gauge.
\newblock {\em Internat. J. Modern Phys. A}, 17(11):1491--1502, 2002.

\bibitem{AbNakOj1997}
M.~Abe, N.~Nakanishi, and I.~Ojima.
\newblock Resolution of the {BRS} singlet-pair problem in quantum {E}instein
  gravity.
\newblock {\em Nuclear Phys. B}, 486(1-2):466--478, 1997.

\bibitem{AzIo1989}
T.~Ya. Azizov and I.~S. Iokhvidov.
\newblock {\em Linear operators in spaces with an indefinite metric}.
\newblock Pure and Applied Mathematics (New York). John Wiley \& Sons Ltd.,
  Chichester, 1989.
\newblock Translated from the Russian by E. R. Dawson, A Wiley-Interscience
  Publication.

\bibitem{AzKh89}
T.~Ya. Azizov and S.~S. Khoruzhi{\u\i}.
\newblock Operators of ghost number and ghost conjugation in
  {BRST}-quantization formalism.
\newblock {\em Teoret. Mat. Fiz.}, 80(1):3--14, 1989.

\bibitem{BarBraHen00}
G.~Barnich, F.~Brandt, and M.~Henneaux.
\newblock Local {BRST} cohomology in gauge theories.
\newblock {\em Phys. Rep.}, 338(5):439--569, 2000.

\bibitem{Bat1987}
I.~A. Batalin.
\newblock The {F}radkin operator method.
\newblock In {\em Quantum field theory and quantum statistics, Vol.\ 1}, pages
  105--127. Hilger, Bristol, 1987.

\bibitem{BatFra1983}
I.~A. Batalin and E.~S. Fradkin.
\newblock A generalized canonical formalism and quantization of reducible gauge
  theories.
\newblock {\em Phys. Lett. B}, 122(2):157--164, 1983.

\bibitem{BRS1975}
C.~Becchi, A.~Rouet, and R.~Stora.
\newblock Renormalization of the abelian {H}iggs-{K}ibble model.
\newblock {\em Comm. Math. Phys.}, 42:127--162, 1975.

\bibitem{BRS1976}
C.~Becchi, A.~Rouet, and R.~Stora.
\newblock Renormalization of gauge theories.
\newblock {\em Ann. Physics}, 98(2):287--321, 1976.

\bibitem{Be1966}
F.~A. Berezin.
\newblock {\em The method of second quantization}.
\newblock Translated from the Russian by Nobumichi Mugibayashi and Alan
  Jeffrey. Pure and Applied Physics, Vol. 24. Academic Press, New York, 1966.

\bibitem{Bla2006}
B.~Blackadar.
\newblock {\em Operator algebras}, volume 122 of {\em Encyclopaedia of
  Mathematical Sciences}.
\newblock Springer-Verlag, Berlin, 2006.
\newblock Theory of $C *$-algebras and von Neumann algebras, Operator Algebras
  and Non-commutative Geometry, III.

\bibitem{Bog1974}
J.~Bogn{\'a}r.
\newblock {\em Indefinite inner product spaces}.
\newblock Springer-Verlag, New York, 1974.
\newblock Ergebnisse der Mathematik und ihrer Grenzgebiete, Band 78.

\bibitem{BorHerWal2000}
M.~Bordemann, H.~Herbig, and S.~Waldmann.
\newblock B{RST} cohomology and phase space reduction in deformation
  quantization.
\newblock {\em Comm. Math. Phys.}, 210(1):107--144, 2000.

\bibitem{BraRob21981}
O.~Bratteli and D.~W. Robinson.
\newblock {\em Operator algebras and quantum-statistical mechanics. {II}}.
\newblock Springer-Verlag, New York, 1981.
\newblock Equilibrium states. Models in quantum-statistical mechanics, Texts
  and Monographs in Physics.

\bibitem{HendrikBuch2006}
D.~Buchholz and H.~Grundling.
\newblock Algebraic supersymmetry: a case study.
\newblock {\em Comm. Math. Phys.}, 272(3):699--750, 2007.

\bibitem{HendrikBuch2007}
D.~Buchholz and H.~Grundling.
\newblock The resolvent algebra: A new approach to canonical quantum systems.
\newblock {\em J Funct Anal}, 254(11):2725--2779, 2008.

\bibitem{CapFer1981}
A.~Z. Capri and R.~Ferrari.
\newblock Schwinger model, chiral symmetry, anomaly and {$\theta $}-vacuums.
\newblock {\em Nuovo Cimento A (11)}, 62(4):273--294, 1981.

\bibitem{Con1985}
J.~B. Conway.
\newblock {\em A course in functional analysis}, volume~96 of {\em Graduate
  Texts in Mathematics}.
\newblock Springer-Verlag, New York, 1985.

\bibitem{DamGeo2004}
M.~Damak and V.~Georgescu.
\newblock Self-adjoint operators affiliated to $c^*$-algebras.
\newblock {\em Rev. Math. Phys.}, 16(2):257--280, 2004.

\bibitem{Dix77}
J.~Dixmier.
\newblock {\em {$C\sp*$}-algebras}.
\newblock North-Holland Publishing Co., Amsterdam, 1977.
\newblock Translated from the French by Francis Jellett, North-Holland
  Mathematical Library, Vol. 15.

\bibitem{DS2}
N.~Dunford and J.~T. Schwartz.
\newblock {\em Linear operators. {P}art {II}}.
\newblock Wiley Classics Library. John Wiley \& Sons Inc., New York, 1988.
\newblock Spectral theory. Selfadjoint operators in Hilbert space, With the
  assistance of William G. Bade and Robert G. Bartle, Reprint of the 1963
  original, A Wiley-Interscience Publication.

\bibitem{DuFred1998}
M.~D{\"u}tsch and K.~Fredenhagen.
\newblock Deformation stability of {BRST}-quantization.
\newblock In {\em Particles, fields, and gravitation (\L \'od\'z, 1998)},
  volume 453 of {\em AIP Conf. Proc.}, pages 324--333. Amer. Inst. Phys.,
  Woodbury, NY, 1998.

\bibitem{DuFred1999}
M.~D{\"u}tsch and K.~Fredenhagen.
\newblock A local (perturbative) construction of observables in gauge theories:
  the example of {QED}.
\newblock {\em Comm. Math. Phys.}, 203(1):71--105, 1999.

\bibitem{DuSh1999}
M.~D{\"u}tsch and G.~Scharf.
\newblock Perturbative gauge invariance: the electroweak theory.
\newblock {\em Ann. Phys. (8)}, 8(5):359--387, 1999.

\bibitem{DuvElTuy1990}
C.~Duval, J.~Elhadad, and G.~M. Tuynman.
\newblock The {BRS} method and geometric quantization: some examples.
\newblock {\em Comm. Math. Phys.}, 126(3):535--557, 1990.

\bibitem{EpGl1973}
H.~Epstein and V.~Glaser.
\newblock The role of locality in perturbation theory.
\newblock {\em Ann. Inst. H. Poincar\'e Sect. A (N.S.)}, 19:211--295 (1974),
  1973.

\bibitem{Fe1960}
J.~M.~G. Fell.
\newblock The dual spaces of $c^*$-algebras.
\newblock {\em Trans. Amer. Math. Soc.}, 94:365--403, 1960.

\bibitem{FiHenStTei1989}
J.~Fisch, M.~Henneaux, J.~Stasheff, and C.~Teitelboim.
\newblock Existence, uniqueness and cohomology of the classical {BRST} charge
  with ghosts of ghosts.
\newblock {\em Comm. Math. Phys.}, 120(3):379--407, 1989.

\bibitem{FraFradk1977}
E.~S. Fradkin and T.~E. Fradkina.
\newblock Quantization of relativistic systems with boson and fermion first-
  and second-class constraints.
\newblock {\em Phys. Lett. B}, 72(3):343--348, 1978.

\bibitem{FraVil1975}
E.~S. Fradkin and G.~A. Vilkovisky.
\newblock Quantization of relativistic systems with constraints.
\newblock {\em Phys. Lett. B}, 55(2):224--226, 1975.

\bibitem{vHo2006}
A.~Fuster and J.~W. van Holten.
\newblock A note on {BRST} quantization of {SU}(2) {Y}ang-{M}ills mechanics.
\newblock {\em J. Math. Phys.}, 46(10):102303, 12, 2005.

\bibitem{Gr2006}
D.R. Girgore.
\newblock Quantum strings and superstrings.
\newblock {\em E-print arXiv:hep-th/0506100v2.}

\bibitem{Got2000}
M.~J. Gotay.
\newblock Obstructions to quantization.
\newblock In {\em Mechanics: from theory to computation}, pages 171--216.
  Springer, New York, 2000.

\bibitem{GotHenTuy1996}
M.~J. Gotay, H.~B. Grundling, and G.~M. Tuynman.
\newblock Obstruction results in quantization theory.
\newblock {\em J. Nonlinear Sci.}, 6(5):469--498, 1996.

\bibitem{GrSchWit1987}
M.~B. Green, J.~H. Schwarz, and E.~Witten.
\newblock {\em Superstring theory. {V}ol. 1}.
\newblock Cambridge Monographs on Mathematical Physics. Cambridge University
  Press, Cambridge, 1987.
\newblock Introduction.

\bibitem{GrSc2003}
D.~R. Grigore and G.~Scharf.
\newblock The quantum supersymmetric vector multiplet and some problems in
  non-abelian supergauge theory.
\newblock {\em Ann. Phys. (8)}, 12(11-12):643--683, 2003.

\bibitem{Hendrik1991}
H.~Grundling.
\newblock {BRST}-quantum theory: A functional analytic approach.
\newblock {\em UNSW preprint 1991}.

\bibitem{Hendrik1988}
H.~Grundling.
\newblock Systems with outer constraints. {G}upta-{B}leuler electromagnetism as
  an algebraic field theory.
\newblock {\em Comm. Math. Phys.}, 114(1):69--91, 1988.

\bibitem{Hendrik2006}
H.~Grundling.
\newblock Quantum constraints.
\newblock {\em Rep. Math. Phys.}, 57(1):97--120, 2006.

\bibitem{HendrikHu1985}
H.~Grundling and C.~A. Hurst.
\newblock Algebraic quantization of systems with a gauge degeneracy.
\newblock {\em Comm. Math. Phys.}, 98(3):369--390, 1985.

\bibitem{HendrikHu1987}
H.~Grundling and C.~A. Hurst.
\newblock Algebraic structures of degenerate systems and the indefinite metric.
\newblock {\em J. Math. Phys.}, 28(3):559--572, 1987.

\bibitem{HenHur1988}
H.~Grundling and C.~A. Hurst.
\newblock The quantum theory of second class constraints: kinematics.
\newblock {\em Comm. Math. Phys.}, 119(1):75--93, 1988.

\bibitem{HendrikHur1993}
H.~Grundling and C.~A. Hurst.
\newblock The operator quantization of the open bosonic string: field algebra.
\newblock {\em Comm. Math. Phys.}, 156(3):473--525, 1993.

\bibitem{Hendrik2000}
H.~Grundling and F.~Lled{\'o}.
\newblock Local quantum constraints.
\newblock {\em Rev. Math. Phys.}, 12(9):1159--1218, 2000.

\bibitem{HaKa}
R.~Haag and D.~Kastler.
\newblock An algebraic approach to quantum field theory.
\newblock {\em J. Mathematical Phys.}, 5:848--861, 1964.

\bibitem{KuHa1979}
H.~Hata and T.~Kugo.
\newblock Subsidiary conditions and physical {$S$}-matrix unitarity in
  covariant canonical formulation of supergravity.
\newblock {\em Nuclear Phys. B}, 158(2-3):357--380, 1979.

\bibitem{Hen1985}
M.~Henneaux.
\newblock Hamiltonian form of the path integral for theories with a gauge
  freedom.
\newblock {\em Phys. Rep.}, 126(1):1--66, 1985.

\bibitem{Hen1988}
M.~Henneaux.
\newblock B{RST} symmetry in the classical and quantum theories of gauge
  systems.
\newblock In {\em Quantum mechanics of fundamental systems, 1 (Santiago,
  1985)}, Ser. Cent. Estud. Cient. Santiago, pages 117--144. Plenum, New York,
  1988.

\bibitem{Hen1989}
M.~Henneaux.
\newblock Duality theorems in {BRST} cohomology.
\newblock {\em Ann. Physics}, 194(2):281--302, 1989.

\bibitem{Hen1993II}
M.~Henneaux.
\newblock Errata: ``{R}emarks on the renormalization of gauge invariant
  operators in {Y}ang-{M}ills theory''.
\newblock {\em Phys. Lett. B}, 316(4):631, 1993.

\bibitem{Hen1993I}
M.~Henneaux.
\newblock Remarks on the renormalization of gauge invariant operators in
  {Y}ang-{M}ills theory.
\newblock {\em Phys. Lett. B}, 313(1-2):35--40, 1993.

\bibitem{HenTei1987}
M.~Henneaux and C.~Teitelboim.
\newblock B{RS} quantisation of generalised magnetic poles.
\newblock In {\em Quantum field theory and quantum statistics, Vol.\ 1}, pages
  165--180. Hilger, Bristol, 1987.

\bibitem{HenTei92}
M.~Henneaux and C.~Teitelboim.
\newblock {\em Quantization of gauge systems}.
\newblock Princeton University Press, Princeton, NJ, 1992.

\bibitem{Holl2008}
S.~Hollands.
\newblock Renormalized quantum yang-mills fields in curved spacetime.
\newblock {\em Rev. Math. Phys.}, 20:1033--1172, 2008.

\bibitem{HoVo89}
S.~S. Horuzhy and A.~V. Voronin.
\newblock Remarks on mathematical structure of {BRST} theories.
\newblock {\em Comm. Math. Phys.}, 123(4):677--685, 1989.

\bibitem{HoVo292}
S.~S. Horuzhy and A.~V. Voronin.
\newblock B{RST} quantization of the {S}chwinger model.
\newblock {\em J. Math. Phys.}, 33(8):2823--2841, 1992.

\bibitem{HorVo1992}
S.~S. Horuzhy and A.~V. Voronin.
\newblock B{RST} quantization of the {S}chwinger model.
\newblock {\em J. Math. Phys.}, 33(8):2823--2841, 1992.

\bibitem{HoVo92}
S.~S. Horuzhy and A.~V. Voronin.
\newblock A new approach to {BRST} operator cohomologies: exact results for the
  {BRST}-{F}ock theories.
\newblock {\em Teoret. Mat. Fiz.}, 93(2):342--353, 1992.

\bibitem{HoVo93}
S.~S. Horuzhy and A.~V. Voronin.
\newblock B{RST} and {$l(1,1)$}.
\newblock {\em Rev. Math. Phys.}, 5(1):191--208, 1993.

\bibitem{HoVo97}
S.~S. Horuzhy and A.~V. Voronin.
\newblock Representations of the {BRST} algebra and unsolvable algebraic
  problems.
\newblock {\em J. Math. Phys.}, 38(8):4301--4322, 1997.

\bibitem{Jak1985}
L.~Jak{\'o}bczyk.
\newblock Canonical quantization with indefinite inner product.
\newblock {\em Ann. Physics}, 161(2):314--336, 1985.

\bibitem{KaOg1982}
M.~Kato and K.~Ogawa.
\newblock Covariant quantization of string based on {BRS} invariance.
\newblock {\em Nuclear Phys. B}, 212:443--460, 1982.

\bibitem{KosSte1987}
B.~Kostant and S.~Sternberg.
\newblock Symplectic reduction, {BRS} cohomology, and infinite-dimensional
  {C}lifford algebras.
\newblock {\em Ann. Physics}, 176(1):49--113, 1987.

\bibitem{KuOj1978I}
T.~Kugo and I.~Ojima.
\newblock Manifestly covariant canonical formulation of the {Y}ang-{M}ills
  field theories. {I}. {G}eneral formalism.
\newblock {\em Progr. Theoret. Phys.}, 60(6):1869--1889, 1978.

\bibitem{KuOj1978II}
T.~Kugo and I.~Ojima.
\newblock Subsidiary conditions and physical {$S$}-matrix unitarity in
  indefinite-metric quantum gravitation theory.
\newblock {\em Nuclear Phys. B}, 144(1):234--252, 1978.

\bibitem{KuOj79}
T.~Kugo and I.~Ojima.
\newblock Local covariant operator formalism of nonabelian gauge theories and
  quark confinement problem.
\newblock {\em Progr. Theoret. Phys. Suppl.}, (66):130, 1979.

\bibitem{KuOJ1979II}
T.~Kugo and I.~Ojima.
\newblock Manifestly covariant canonical formulation of {Y}ang-{M}ills field
  theories. {II}. {${\rm SU}(2)$} {H}iggs-{K}ibble model with spontaneous
  symmetry breaking.
\newblock {\em Progr. Theoret. Phys.}, 61(1):294--314, 1979.

\bibitem{KuUe1980}
T.~Kugo and S.~Uehara.
\newblock The form of {BRS}-invariant operators.
\newblock {\em Progr. Theoret. Phys.}, 64(4):1395--1411, 1980.

\bibitem{LanLin1992}
N.~P. Landsman and N.~Linden.
\newblock Superselection rules from {D}irac and {BRST} quantisation of
  constrained systems.
\newblock {\em Nuclear Phys. B}, 371(1-2):415--433, 1992.

\bibitem{Mc1991}
D.~McMullan.
\newblock Gauge fixing and ghost variables.
\newblock {\em Nuclear Phys. B}, 363(2-3):451--485, 1991.

\bibitem{McPat1989I}
D.~McMullan and J.~Paterson.
\newblock Covariant factor ordering of gauge systems using ghost variables.
  {C}onstraint rescaling.
\newblock {\em J. Math. Phys.}, 30(2):477--486, 1989.

\bibitem{McPat1989II}
D.~McMullan and J.~Paterson.
\newblock Covariant factor ordering of gauge systems using ghost variables.
  {II}. {S}tates and observables.
\newblock {\em J. Math. Phys.}, 30(2):487--497, 1989.

\bibitem{Min1980}
M.~Mintchev.
\newblock Quantisation in indefinite metric.
\newblock {\em J. Phys. A}, 13(5):1841--1859, 1980.

\bibitem{Mur90}
G.~J. Murphy.
\newblock {\em $C^*$-algebras and operator theory}.
\newblock Academic Press Inc., Boston, MA, 1990.

\bibitem{Na1964}
M.~A. Naimark.
\newblock {\em Linear representations of the {L}orentz group}.
\newblock Translated by Ann Swinfen and O. J. Marstrand; translation edited by
  H. K. Farahat. A Pergamon Press Book. The Macmillan Co., New York, 1964.

\bibitem{Nak1979}
N.~Nakanishi.
\newblock On the general validity of the unitarity proof in the {K}ugo-{O}jima
  formalism of gauge theories.
\newblock {\em Progr. Theoret. Phys.}, 62(5):1396--1402, 1979.

\bibitem{NakOj1990}
N.~Nakanishi and I.~Ojima.
\newblock {\em Covariant operator formalism of gauge theories and quantum
  gravity}, volume~27 of {\em World Scientific Lecture Notes in Physics}.
\newblock World Scientific Publishing Co. Inc., Teaneck, NJ, 1990.

\bibitem{Nish1984}
K.~Nishijima.
\newblock Representations of {BRS} algebra.
\newblock {\em Nuclear Phys. B}, 238(3):601--620, 1984.

\bibitem{Nish1996}
K.~Nishijima.
\newblock B{RS}-invariance, asymptotic freedom and color confinement (a
  review).
\newblock {\em Czechoslovak J. Phys.}, 46(1):1--40, 1996.

\bibitem{Ota1984}
S.~{\=O}ta.
\newblock Closed linear operators with domain containing their range.
\newblock {\em Proc. Edinburgh Math. Soc. (2)}, 27(2):229--233, 1984.

\bibitem{Ped79}
G.~K. Pedersen.
\newblock {\em $C^*$-algebras and their automorphism groups}, volume~14 of {\em
  London Mathematical Society Monographs}.
\newblock Academic Press Inc. [Harcourt Brace Jovanovich Publishers], London,
  1979.

\bibitem{Pin1998}
M.~A. Pinsky.
\newblock {\em Partial differential equations and boundary value problems with
  applications}.
\newblock International Series in Pure and Applied Mathematics. McGraw-Hill
  Inc., New York, third edition, 1998.

\bibitem{Zw1998}
J~Portegies~Zwart.
\newblock Brst reduction and quantization of constrained hamiltonian systems.
\newblock Master's thesis, The University of Amsterdam, 1998.

\bibitem{RazRyb1990}
A.~V. Razumov and G.~N. Rybkin.
\newblock State space in {BRST}-quantization of gauge-invariant systems.
\newblock {\em Nuclear Phys. B}, 332(1):209--223, 1990.

\bibitem{ReSi1972v1}
M.~Reed and B.~Simon.
\newblock {\em Methods of modern mathematical physics. {I}. {F}unctional
  analysis}.
\newblock Academic Press, New York, 1972.

\bibitem{ReSi1975v2}
M.~Reed and B.~Simon.
\newblock {\em Methods of modern mathematical physics. {II}. {F}ourier
  analysis, self-adjointness}.
\newblock Academic Press [Harcourt Brace Jovanovich Publishers], New York,
  1975.

\bibitem{Rob1999}
P.~L. Robinson.
\newblock The {B}erezin calculus.
\newblock {\em Publ. Res. Inst. Math. Sci.}, 35(2):123--194, 1999.

\bibitem{Rud1987}
W.~Rudin.
\newblock {\em Real and complex analysis}.
\newblock McGraw-Hill Book Co., New York, third edition, 1987.

\bibitem{RuRuiAl2004}
H.~Ruegg and M.~Ruiz-Altaba.
\newblock The {S}tueckelberg field.
\newblock {\em Internat. J. Modern Phys. A}, 19(20):3265--3347, 2004.

\bibitem{Schf2001}
G.~Scharf.
\newblock {\em Quantum gauge theories}.
\newblock Wiley-Interscience [John Wiley \& Sons], New York, 2001.
\newblock A true ghost story.

\bibitem{Sch1994}
R.~Schmid.
\newblock Local cohomology in gauge theories, {BRST} transformations and
  anomalies.
\newblock {\em Differential Geom. Appl.}, 4(2):107--116, 1994.

\bibitem{Sch64}
S.~Schwebber.
\newblock {\em An Introduction to Relativistic Quantum Field Theory}.
\newblock Harper \& Row, New York, 1964.

\bibitem{Sl1989}
A.~A. Slavnov.
\newblock A unitarity condition in covariant quantum field theory with an
  indefinite metric.
\newblock {\em Teoret. Mat. Fiz.}, 79(3):347--358, 1989.

\bibitem{Sta1998}
J.~Stasheff.
\newblock The (secret?) homological algebra of the {B}atalin-{V}ilkovisky
  approach.
\newblock In {\em Secondary calculus and cohomological physics (Moscow, 1997)},
  volume 219 of {\em Contemp. Math.}, pages 195--210. Amer. Math. Soc.,
  Providence, RI, 1998.

\bibitem{St1989}
O.~Steinmann.
\newblock On the characterization of physical states in gauge theories.
\newblock {\em Ann. Inst. H. Poincar\'e Phys. Th\'eor.}, 51(3):299--321, 1989.

\bibitem{Str1967}
F.~Strocchi.
\newblock Gauge problem in quantum field theory.
\newblock {\em Phys. Rev. (2)}, 162:1429--1438, 1967.

\bibitem{Tak2002}
M.~Takesaki.
\newblock {\em Theory of operator algebras. {I}}, volume 124 of {\em
  Encyclopaedia of Mathematical Sciences}.
\newblock Springer-Verlag, Berlin, 2002.
\newblock Reprint of the first (1979) edition, Operator Algebras and
  Non-commutative Geometry, 5.

\bibitem{Ty1975}
I.V. Tyutin.
\newblock Unpublished.
\newblock {\em Lebedev preprint FIAN}, 39, 1975.

\bibitem{vHo1990}
J.~W. van Holten.
\newblock The {BRST} complex and the cohomology of compact {L}ie algebras.
\newblock {\em Nuclear Phys. B}, 339(1):158--176, 1990.

\bibitem{vHol2005}
J.~W. van Holten.
\newblock Aspects of {BRST} quantization.
\newblock In {\em Topology and geometry in physics}, volume 659 of {\em Lecture
  Notes in Phys.}, pages 99--166. Springer, Berlin, 2005.

\bibitem{VoKh2000}
A.~V. Voronin and S.~S. Khoruzhi{\u\i}.
\newblock Conformal theories, the {BRST} approach, and representations of {L}ie
  superalgebras.
\newblock {\em Tr. Mat. Inst. Steklova}, 228(Probl. Sovrem. Mat.
  Fiz.):155--167, 2000.

\bibitem{Wein2005II}
S.~Weinberg.
\newblock {\em The quantum theory of fields. {V}ol. {II}}.
\newblock Cambridge University Press, Cambridge, 2005.
\newblock Modern applications.

\bibitem{Yos1980}
K.~Yosida.
\newblock {\em Functional analysis}, volume 123 of {\em Grundlehren der
  Mathematischen Wissenschaften [Fundamental Principles of Mathematical
  Sciences]}.
\newblock Springer-Verlag, Berlin, sixth edition, 1980.

\bibitem{ZJ1993}
J.~Zinn-Justin.
\newblock {\em Quantum field theory and critical phenomena}, volume~85 of {\em
  International Series of Monographs on Physics}.
\newblock The Clarendon Press Oxford University Press, New York, second
  edition, 1993.
\newblock Oxford Science Publications.

\end{thebibliography}

\end{document}